\documentclass[spanningrule]{cambridge7A} 

\usepackage{amsthm}
\usepackage[numbers]{natbib}
\usepackage{rotating}
\usepackage{floatpag}
\rotfloatpagestyle{empty}
\usepackage{graphicx}
\usepackage{multind}

\usepackage{amsmath,amsfonts,amssymb,mathrsfs}

\input{PaulTaylorDiagramsV393}

\theoremstyle{plain} 

\newtheorem{global-theorem}{Theorem}
\newtheorem{theorem}{Theorem}[section]
\newtheorem{lemma}[theorem]{Lemma}
\newtheorem{remark}[theorem]{Remark}
\newtheorem{scholium}[theorem]{Scholium}
\newtheorem{corollary}[theorem]{Corollary}
\newtheorem{condition}[theorem]{Condition}

\newtheorem{conjecture}[theorem]{Conjecture}
\newtheorem{definition}[theorem]{Definition}
\newtheorem{proposition}[theorem]{Proposition}
\newtheorem{prop-def}[theorem]{Proposition-Definition}
\newtheorem{lemma-def}[theorem]{Lemma-Definition}

\diagramstyle[l>=1em]

\newarrow{Equalarr}=====
\newarrowtail{<-}\leftarrow\rightarrow\uparrow\downarrow

\newarrow{Botharr}{<-}---{->}
\newarrow{Botharrow}{<-}---{->}
\newarrow{Surjarrow}----{->>}
\newarrow{TeXto}----{->}
\newarrow{Hookarrow}C---{->}

\newcommand{\rightarr}{\rTeXto}
\newcommand{\rightarrt}{\rTeXto}
\newcommand{\leftarr}{\lTeXto}
\newcommand{\downarr}{\dTeXto}
\newcommand{\uparr}{\uTeXto}

\newcommand{\sqb}[1]{[ #1 ]}

\newcommand{\leftbrack}{[}

\newcommand{\nocom}{}

\newcommand{\eop}{\ \hfill $\Box$}

\numberwithin{equation}{section}

\newcommand{\mynolabel}[1]{}

\newcommand{\nn}{{\mathbb N}}

\newcommand{\rr}{{\mathbb R}}

\newcommand{\zz}{{\mathbb Z}}

\newcommand{\aaa}{{\mathbb A}}

\newcommand{\pP}{{\bf P}}

\newcommand{\pA}{{\mathcal A}}
\newcommand{\pB}{{\mathcal B}}
\newcommand{\pC}{{\mathcal C}}
\newcommand{\pD}{{\mathcal D}}
\newcommand{\pE}{{\mathcal E}}
\newcommand{\pF}{{\mathcal F}}
\newcommand{\pG}{{\mathcal G}}

\newcommand{\pQ}{{\mathcal Q}}
\newcommand{\pR}{{\mathcal R}}

\newcommand{\pT}{{\mathcal T}}
\newcommand{\pU}{{\mathcal U}}
\newcommand{\pV}{{\mathcal V}}

\newcommand{\pX}{{\mathcal X}}
\newcommand{\pY}{{\mathcal Y}}

\newcommand{\mA}{{\mathscr A}}
\newcommand{\mB}{{\mathscr B}}
\newcommand{\mC}{{\mathscr C}}

\newcommand{\mE}{{\mathscr E}}
\newcommand{\mF}{{\mathscr F}}
\newcommand{\mG}{{\mathscr G}}

\newcommand{\mK}{{\mathscr K}}

\newcommand{\mM}{{\mathscr M}}
\newcommand{\mN}{{\mathscr N}}

\newcommand{\mP}{{\mathscr P}}
\newcommand{\mQ}{{\mathscr Q}}

\newcommand{\mV}{{\mathscr V}}
\newcommand{\mW}{{\mathscr W}}
\newcommand{\mX}{{\mathscr X}}
\newcommand{\mY}{{\mathscr Y}}

\newcommand{\Ee}{{\mathscr E}}

\newcommand{\Ff}{{\mathscr F}}

\newcommand{\Rr}{{\mathscr R}}
\newcommand{\Mm}{{\mathscr M}}
\newcommand{\Vv}{{\mathscr V}}

\newcommand{\Ww}{{\mathscr W}}

\newcommand{\Cc}{{\mathscr C}}
\newcommand{\Kk}{{\mathscr K}}

\newcommand{\Hom}{\mbox{\sc Hom}}
\newcommand{\uHom}{\underline{\Hom}}

\newcommand{\Sets}{\mbox{\sc Set}}
\newcommand{\Top}{\mbox{\sc Top}}

\newcommand{\inj}{{\bf inj}}

\newcommand{\cof}{{\bf cof}}

\newcommand{\cell}{{\bf cell}}

\newcommand{\Cell}{{\mbox{{\sc Cell}}}}

\newcommand{\Ob}{{\rm Ob}}
\newcommand{\ob}{{\rm Ob}}

\newcommand{\Mor}{{\rm Mor}}

\newcommand{\uMor}{{\underline{\Mor }}}

\newcommand{\glob}{{\bf G}}

\newcommand{\colim}{{\rm colim}}
\newcommand{\mylim}{{\rm lim}}

\newcommand{\disc}{{\bf disc}}
\newcommand{\codisc}{{\bf codsc}}

\newcommand{\univa}{{\mathbb U}}
\newcommand{\univb}{{\mathbb V}}

\newcommand{\diag}{\mbox{\sc Func}}

\newcommand{\func}{\mbox{\sc Func}}

\newcommand{\precat}{{\bf PC}}

\newcommand{\presh}{{\rm Presh}}

\newcommand{\scone}{{\bf C}}
\newcommand{\sconeplus}{{\bf C}^+}
\newcommand{\sconemap}{{\bf c}}

\newcommand{\src}{{\rm src}}
\newcommand{\targ}{{\rm targ}}

\newcommand{\Alg}{\mbox{\sc Alg}}

\newcommand{\Arr}{\mbox{\sc Arr}}

\newcommand{\Ho}{\mbox{\rm ho}}

\newcommand{\gen}{{\bf gen}}
\newcommand{\fix}{{\bf fix}}

\newcommand{\Enr}{{\bf Enr}}

\newcommand{\hec}{{\bf eh}}

\newcommand{\Gen}{\mbox{\sc Gen}}

\newcommand{\Fix}{\mbox{\sc Fix}}

\newcommand{\Int}{\mbox{\sc Int}}

\newcommand{\Cat}{\mbox{\sc Cat}}

\newcommand{\Seg}{\mbox{\sc Seg}}

\newcommand{\Iso}{{\rm Iso}}

\newcommand{\resol}{{\bf R}}

\newcommand{\Upstild}{\widetilde{\Upsilon}}

\newcommand{\bfi}{{\bf i}}

\newcommand{\sk}{\mbox{sk}}
\newcommand{\csk}{\mbox{csk}}

\newcommand{\dgt}{{\bf d}}

\newcommand{\underlying}{{\mathcal U}}

\newcommand{\mylabel}[1]{\label{#1}}

\newcommand{\Latch}{{\rm Latch}}
\newcommand{\Match}{{\rm Match}}

\newcommand{\latch}{{\rm latch}}
\newcommand{\match}{{\rm match}}

\newcommand{\Opetopes}{{\bf Otp}}

\newcommand{\Glob}{{\mathbb G}}
\newcommand{\globe}{{\bf g}}

\newcommand{\gpasting}{\pG\pP}

\begin{document}

\vspace*{4cm}

\begin{center}
{\bf \LARGE
Homotopy Theory \vspace*{2mm}\\
of  \vspace*{4mm}\\ 
Higher Categories}
\end{center}

\thispagestyle{empty}

\vspace*{3cm}

\begin{center}
Carlos T. Simpson
\vspace*{.3cm}
\\
CNRS---INSMI
\\
Laboratoire J.\ A.\ Dieudonn\'e
\\
Universit\'e de Nice-Sophia Antipolis
\\
\verb+carlos@unice.fr+
\end{center}

\vfill

Ce papier a b\'en\'efici\'e d'une aide de l'Agence Nationale de la Recherche
portant la r\'ef\'erence ANR-09-BLAN-0151-02 (HODAG). 

\vspace*{3mm}

This is draft material from a forthcoming book to be published by 
Cambridge University Press in the New Mathematical Monographs series. This 
publication is in copyright.  \copyright Carlos T. Simpson 2010.

\newpage

\vspace*{4cm}

\begin{abstract}
This is the first draft of a book about higher categories approached by
iterating Segal's method, as in Tamsamani's definition of $n$-nerve and Pelissier's thesis.
If $M$ is a tractable left proper cartesian model category,
we construct a tractable left proper cartesian model structure on the category 
of $M$-precategories. The procedure can then be iterated, leading to model categories 
of $(\infty , n)$-categories. 
\end{abstract}

\pagenumbering{roman}
\setcounter{page}{5}

\tableofcontents

\cleardoublepage

\pagenumbering{arabic}
\setcounter{page}{1}

\copyrightline{This is draft material from a forthcoming book to be 
published by Cambridge University Press in the New Mathematical Monographs 
series. This publication is in copyright. \copyright Carlos T. Simpson 2010.}


\chapter*{Preface}
\label{intro}

The theory of $n$-categories is currently under active consideration by a number of different
research groups around the world. The history of the subject goes back a long way, on
separate but interrelated tracks in algebraic topology, algebraic geometry, and category theory. 
For a long time, the crucial definition of {\em weakly associative higher category} remained
elusive, but now on the contrary we have a plethora of different possibilities available. 
One of the next major problems in the subject will be 
to achieve a global comparison between these
different approaches. Some work is starting to come out in this direction, but in the current state of
the theory the various different approaches remain distinct. After the comparison is achieved, they
will be seen as representing different facets of the theory, so it is important to continue working
in all of these different directions. 

The purpose of the present book is to concentrate on one of the methods of defining and working
with higher categories, very
closely based on the work of Graeme Segal in algebraic topology many years earlier.
The notion of ``Segal category'', which is a kind of category weakly enriched over simplicial sets, was considered by Vogt and Dwyer, Kan and Smith. The application of this method
to $n$-categories was introduced by Zouhair Tamsamani. And then 
put into a strictly iterative
form, with a general model category as input,
by Regis Pelissier following a suggestion of Andr\'e Hirschowitz. Our treatment will integrate important ideas contributed by Julie Bergner, Clark Barwick, Jacob Lurie and others. 

The guiding principle is to use
the category of simplices $\Delta$ as the basis for all the higher coherency conditions which come in 
when we allow weak associativity. The objects of $\Delta$ are nonempty finite ordinals  
$$
\begin{array}{lcc}
{[ 0 ]} &=& \{ \upsilon _0\}\\
{[ {1} ]} &= &\{ \upsilon _0,\upsilon _1 \} \\
{[ {2} ]} &=& \{ \upsilon _0,\upsilon _1,\upsilon _2 \} \\
& \ldots & 
\end{array}
$$
whereas the morphisms are nondecreasing maps between them. Kan had already introduced this category
into algebraic topology, considering {\em simplicial sets} which are functors $\Delta ^o\rightarrow \Sets$.
These model the homotopy types of CW-complexes. 

One of the big problems in algebraic topology in the 1960's
was to define notions of {\em delooping machines}. Segal's way was to consider {\em simplicial spaces}, or functors $\pA :\Delta ^o \rightarrow
\Top$, such that the first space is just a point $\pA _0= \pA ([0])=\ast$. In $\Delta$ there are three nonconstant maps 
$$
f_{01},f_{12},f_{02}: [1]\rightarrow [2]
$$
where $f_{ij}$ denotes the map sending $\upsilon _0$ to $\upsilon _i$ and $\upsilon _1$ to $\upsilon _j$.
In a simplicial space $\pA $ which is a contravariant functor on $\Delta$, we get three maps
$$
f_{01}^{\ast},f_{12}^{\ast},f_{02}^{\ast}: \pA _2\rightarrow \pA _1.
$$
Organize the first two as a map into a product, giving a diagram of the form
$$
\begin{diagram}
\pA _2 & \rightarr^{\sigma _2} & \pA _1\times \pA _1 \\
\downarr^{ f_{02}^{\ast}} & & \\
\pA _1 & &
\end{diagram} .
$$
If we require the {\em second Segal map} $\sigma _2:= (f_{01}^{\ast},f_{12}^{\ast})$ to be an isomorphism between $\pA _2$ and $\pA _1\times \pA _1$,
then $f_{02}^{\ast}$ gives a product on the space $\pA _1$. The basic idea of Segal's delooping machine is that if
we only require 
$\sigma _2$ to be a weak homotopy equivalence of spaces, then 
$f_{02}^{\ast}$ gives what should be considered as a ``product defined up to homotopy''. One salient aspect of this approach is that
no map $\pA _1\times \pA _1\rightarrow \pA _1$ is specified, and indeed if the spaces involved have bad properties there might exist no
section of $\sigma _2$ at all.

The
term ``delooping machine'' refers to any of several kinds of further
mathematical structure on the loop space $\Omega X$, enhancing
the basic composition of loops up to homotopy, which should allow one to reconstruct
a space $X$
up to homotopy. In tandem with 
Segal's machine in which $\Omega X=\pA_1$, the notion of {\em operad} introduced by Peter May
underlies the best known and most studied family of delooping machines. There were
also other techniques such as PROP's
which are starting to receive renewed interest. 
The various kinds of delooping machines are sources for the various different
approaches to higher categories, after a multiplying effect whereby each
delooping technique leads to several different definitions of higher categories. 
Our technical work in
later parts of the book will concentrate on the particular direction of iterating
Segal's approach while maintaining a discrete set of objects, but we will survey some of the
many other approaches in later chapters of Part I. 

The relationship between categories and simplicial objects was noticed early on with the {\em nerve construction}. Given a category
$\pC $ its {\em nerve} is the simplicial set $N\pC :\Delta ^o\rightarrow \Sets$ 
such that $(N\pC )_m$ is the set of composable 
sequences of $m$ arrows 
$$
x_0\stackrel{g_1}{\rightarrow}x_1 \stackrel{g_2}{\rightarrow} \cdots x_{m-1}\stackrel{g_m}{\rightarrow}x_m
$$
in $\pC $. The operations of functoriality for maps $[m]\rightarrow [p]$ are obtained using the composition law and the identities of $\pC $. 
The first piece is just the set $(N\pC )_0=\Ob (\pC )$ of objects of $\pC $, and the nerve satisfies a relative version of the Segal condition:
$$
\sigma _m: (N\pC )_m \stackrel{\cong}{\longrightarrow} (N\pC )_1\times _{(N\pC )_0}(N\pC )_1 \times _{(N\pC )_0}\cdots \times _{(N\pC )_0}(N\pC )_1.
$$
Conversely, any simplicial set $\Delta ^o\rightarrow \Sets$ satisfying these conditions comes from a unique category
and these constructions are inverses. In other words, categories may be considered as simplicial sets satisfying 
the Segal conditions. 

In comparing this with Segal's situation, recall that he required $\pA _0=\ast$, which is like looking at a category with
a single object. 

An obvious way of putting all of these things together is to consider simplicial spaces $\pA : \Delta ^o \rightarrow \Top$
such that $\pA _0$ is a discrete set---to be thought of as the ``set of objects''---but considered as a space, and such that the Segal maps 
$$
\sigma _m : \pA _m \rightarrow \pA _1\times _{\pA _0}\pA _1\times _{\pA _0}\cdots \times _{\pA _0}\pA _1
$$
are weak homotopy equivalences for all $m\geq 2$. Functors of this kind are {\em Segal categories}.
We use the same terminology when $\Top$ is replaced by the category of simplicial sets $\Kk := \Sets ^{\Delta ^o} = \diag (\Delta ^o, \Sets )$.
The possibility of making this generalization was clearly evident at the time of Segal's papers \cite{SegalHspaces} \cite{Segal}, but was
made explicit only later by Vogt \cite{Vogt} and Dwyer, Kan and Smith \cite{DKS}. 

Segal categories provide a good way of considering categories enriched over spaces. However, a more elementary approach is
available, by looking at categories {\em strictly} enriched over spaces, i.e. simplicial categories. A simplicial category
could be viewed as a Segal category where the Segal maps $\sigma _m$ are isomorphisms. More classically it can be considered as a
category enriched in $\Top$ or $\Kk$, using the definitions of enriched category theory. In a simplicial category,
the product operation is well-defined and strictly associative. 

Dwyer, Kan and Smith showed that
we don't lose any generality at this level by requiring strict associativity: every Segal category is equivalent to
a simplicial category \cite{DKS}. Unfortunately, we cannot just iterate the construction
by continuing to look at categories strictly 
enriched over the category of simplicial categories
and so forth. Such an iteration leads to higher categories 
with strict associativity and strict
units in the middle levels. One way of seeing why this isn't good 
enough\footnote{Paoli has shown \cite{PaoliAdvances} 
that $n$-groupoids can be {\em semistrictified} in any
single degree, however one cannot get strictness in many different degrees as will be seen in Chapter \ref{nonstrict1}.}
is to note
that the Bergner model structure on strict simplicial categories, is not cartesian:
products of cofibrant objects are no longer cofibrant. This suggests the need for a
different construction which preserves the cartesian condition, and Segal's method works. 

The iteration then says:
a Segal $(n+1)$-category is a functor from $\Delta ^o$ to the category of Segal $n$-categories,
whose first element is a discrete set, and such that the Segal maps are equivalences. 
The notion of equivalence needs to be defined in the inductive process \cite{Tamsamani}. 
This iterative
point of view towards higher categories is the topic of our book. 

We emphasize an algebraic approach within the world of homotopy theory, 
using Quillen's homotopical
algebra \cite{Quillen}, but also paying particular attention to 
the process of creating a higher category from 
generators and relations. For me this goes back to Massey's book \cite{Massey} which was
one of my first references both for for algebraic topology, and for the notion of a group 
presented by generators and relations.

One of the main inspirations for the recent interest in higher categories came from Grothendieck's manuscript 
{\em Pursuing stacks}. He set out a wide vision of the possible developments and applications of the theory
of $n$-categories going up to $n=\omega$. Many of his remarks continue to provide important research directions, 
and many others remain untouched.

The other main source of interest stems from the {\em Baez-Dolan conjectures}.
These extend, to higher categories in all degrees, the relationships explored by
many researchers  between various categorical structures and
phenomena of knot invariants and
quantum field theory. Hopkins and Lurie have recently
proven a major part of these conjectures. These motivations incite us to search for a good understanding of the algebra of
higher categories, and I hope that the present book can contribute in a small way. 

\vspace*{3mm}

The mathematical discussion of the contents of the main part of the book
will be continued in more detail in Chapter \ref{cartenr1} at the end of
Part I. The intervening chapters of Part I serve to introduce the problem by 
giving some motivation for why higher categories are needed, by explaining why 
strict $n$-categories aren't enough, and by considering some of the many other approaches
which are currently being developed. 

In Part II we collect our main tools from the theory of categories, including locally
presentable categories and  closed model categories. A small number of these items,
such as the discussion of enriched categories,
could be useful for reading Part I. The last chapter of Part II concerns
``direct left Bousfield localization'', which is a special case of left Bousfield localization
in which the model structure can be described more explicitly. This Chapter \ref{direct1}, together with
the general discussion of cell complexes and 
the small object argument in Chapters \ref{cattheor1}
and \ref{modcat1}, are intended
to provide some ``black boxes'' which can then be used in the rest of the book without
having to go into details of cardinality arguments and the like. It is hoped that this
will make a good part of the book accessible to readers wishing to avoid too many
technicalities of the theory of model categories, although some familiarity is obviously
necessary since our main goal is to construct a cartesian model structure. 

In Part III starts the main work of looking at weakly $\mM$-enriched (pre)categories.
This part is entitled ``Generators and relations'' because the process of starting
with an $\mM$-precategory and passing to the associated weakly $\mM$-enriched
category by enforcing the
Segal condition, should be viewed as a higher or weakly enriched analogue of the
classical process of describing an algebraic object by generators and relations. 
We develop several aspects of this point of view, including a detailed discussion
of the example of categories weakly enriched over the model category $\mK$
of simplicial sets. We see in this case how to follow along the calculus
of generators and relations, taking as example the calculation of the loop space of $S^2$.

Part IV contains the construction of the cartesian model category, after the two steps
treating specific elements of our categorical situation: products, including 
the proof of the cartesian condition, and intervals. 

Part V, not yet present in the current version,
will discuss various directions going towards basic techniques in the theory 
of higher categories,
using the formalism developed in Parts II-IV. 

The first few chapters of Part V should
contain discussions of inversion of morphisms, limits and colimits, and adjunctions,
based to a great extent on my preprint \cite{limits}. 

For the case of $(\infty ,1)$-categories, these topics are treated in
Lurie's recent book \cite{LurieTopos} about the analogue of Grothendieck's 
theory of topoi, using 
quasicategories. 

Other topics from higher category theory which will only be discussed very briefly
are the theory of higher stacks, and the Poincar\'e $n$-groupoid and van Kampen theorems.
For the theory of higher
stacks the reader can consult \cite{descente} which starts from 
the model categories
constructed here.

For a discussion of the Poincar\'e $n$-groupoid, the reader can
consult Tamsamani's original paper \cite{Tamsamani} as well as Paoli's discussion of
this topic in the context of special $Cat^n$-groups \cite{PaoliAdvances}. 
 
We will spend a chapter looking at the {\em Breen-Baez-Dolan stabilization hypothesis},
following the preprint \cite{BBDSH}.
This is one of the first parts of the famous {\em Baez-Dolan conjectures}. These
conjectures have strongly motivated the development of higher category theory.
Hopkins and Lurie have recently proven important pieces of the main conjectures.
The stabilization conjecture is a preliminary statement about the behavior of the
notion of $k$-connected $n$-category, understandable with the basic techniques we
have developed here. We hope that this will serve as an introduction to an exciting
current area of research.


\chapter*{Acknowledgements}

I would first like to thank Zouhair Tamasamani, whose original work on this 
question led to
all the rest. His techniques for gaining access to a theory of $n$-categories
using Segal's delooping machine, set out the basic contours of the theory, 
and continue to inform and guide our understanding.   
I would like to thank Andr\'e Hirschowitz for much encouragement and many
interesting conversations in the course of our work on descent for $n$-stacks, 
$n$-stacks of complexes, and higher Brill-Noether. 
And to thank Andr\'e's thesis student Regis
Pelissier who took the argument to a next stage of abstraction, braving the multiple
difficulties not the least of which were the cloudy reasoning and several 
important errors in one of my preprints. We are following quite closely the main idea
of Pelissier's thesis, which is to iterate a construction whereby a good model category
$\mM$ 
serves as input, and we try to get out a model category of $\mM$-enriched precategories.
Clark Barwick then added a further crucial insight, which was that the argument could be
broken down into pieces, starting with a fairly classical left Bousfield localization.
Here again, Clark's idea serves as groundwork for our approach. Jacob Lurie continued with
many contributions, on different 
levels most of which are beyond our immediate grasp; but still
including some quite understandable innovations such as the idea of
introducing the category $\Delta _X$ of finite ordered sets decorated with elements of 
the set $X$ ``of objects''. This leads to a significant lightening of the hypotheses
needed of $\mM$. His approach to cardinality
questions in the small object argument is groundbreaking, and we give here
an alternate treatment which is certainly less streamlined but might help the reader to
situate what is going on. These items are of course subordinated to the use of Smith's recognition principle and Dugger's notion of combinatorial model category, on which our constructions are based. 
Julie Bergner has gained important information about a whole range of model structures
starting with her consolidation of the Dwyer-Kan structure on the category of
simplicial categories. Her characterization of fibrant objects in the model structures
for Segal categories, carries over easily to our case and provides the basis for
important parts of the statements of our main results. 

Bertrand Toen was 
largely responsible for teaching me
about model categories. His philosophy that they are a good down-to-earth yet powerful
approach to higher categorical questions, is suffused throughout this work. I would like
to thank Joseph Tapia and Constantin Teleman for their encouragement in this direction too,
and Georges Maltsiniotis and Alain Brugui\`eres who were able to explain the 
higher categorical meaning of the Eckmann-Hilton
argument in an understandable way.  
I would similarly like to thank Clemens Berger, Ronnie Brown, Eugenia Cheng,
Denis-Charles Cisinski, Delphine Dupont, Joachim Kock, Peter May, and
Simona Paoli,
for many interesting and informative conversations about
various aspects of this subject; and to thank 
my current 
doctoral students for continuing discussions in directions extending the present work, 
which motivate the completion of
this project.  I would also like to thank my co-workers on related
topics, things which if they are not directly present here, have still contributed a lot
to the motivation for the study of higher categories. 

I would specially like to thank Diana Gillooly of Cambridge University Press, and Burt Totaro, for setting this project in motion. 

Paul Taylor's diagram package is used for the commutative diagrams and even
for the arrows in the text. 

For the title, we have chosen something almost the same as the title of a special
semester in Barcelona a while back; we apologize for this overlap.

This research is partially supported by the Agence Nationale de la Recherche, grant
ANR-09-BLAN-0151-02 (HODAG). I would like to thank the Institut de Math\'ematiques de Jussieu for their hospitality during the completion of this work. 




\part{Higher categories}


\chapter{History and motivation}
\label{why1}

The most basic motivation for introducing higher categories is the observation that $\Cat _{\univa}$, the category of
$\univa$-small categories, naturally has a structure of $2$-category: the objects are categories, the morphisms are functors,
and the $2$-morphisms are natural transformations between functors. If we denote this $2$-category by $\Cat ^{2{\rm cat}}$ then 
its truncation $\tau _{\leq 1}\Cat ^{2{\rm cat}}$ to a $1$-category would have, as morphisms, the equivalence classes of functors up to natural equivalence. 
While it is often necessary to consider two naturally equivalent functors as being ``the same'', identifying them formally leads to a loss of information. 

Topologists are confronted with a similar situation when looking at the category of spaces. In homotopy theory one thinks of two homotopic maps
between spaces as being ``the same''; however the {\em homotopy category} $\Ho (\Top )$ obtained after dividing by this equivalence relation, doesn't retain
enough information. This loss of information is illustrated by the question of diagrams. Suppose $\Psi$ is a small category. A {\em diagram of spaces}
is a functor $T:\Psi \rightarrt \Top$, that is a space $T(x)$ for each object $x\in \Psi$ and a map $T(a):T(x)\rightarrt T(y)$ for each arrow $a\in \Psi (x,y)$,
satisfying strict compatibility with identities and compositions. The category of diagrams $\diag (\Psi , \Top )$ has a natural subclass of morphisms: a morphism
$f:S\rightarrt T$ of diagrams is a {\em levelwise weak equivalence} if each $f(x):S(x)\rightarrt T(x)$ is a weak equivalence. Letting $\mW = 
\mW _{\diag (\Psi , \Top )}$ 
denote this subclass,
the homotopy category of diagrams $\Ho ( \diag (\Psi , \Top ))$ is defined to be the Gabriel-Zisman localization
$\mW ^{-1}\diag (\Psi , \Top )$. There is a natural functor
$$
\Ho ( \diag (\Psi , \Top ))\rightarrt \diag (\Psi , \Ho (\Top )),
$$
which is {\em not} in general an equivalence of categories. In other words $\Ho (\Top )$ doesn't retain enough information to recover 
$\Ho ( \diag (\Psi , \Top ))$. Thus, we need to consider some kind of extra structure
beyond just the homotopy category. 

This phenomenon occurs in a number of different places. 
Starting in the 1950's and 1960's, the notion of {\em derived category}, an abelianized version of $\Ho (\;\; )$, became crucial to 
a number of areas in modern homological algebra and particularly for
algebraic geometry. The notion of localization of a category seems to
have been proposed in this context
by Serre, and appears in Grothendieck's Tohoku paper \cite{GrothendieckTohoku}. 
A systematic
treatment is the subject of Gabriel-Zisman's book \cite{GabrielZisman}. 

As the example of diagrams illustrates, in many derived-categorical 
situations one must first make some intermediate constructions on
underlying categorical data, then pass to the derived category. 
A fundamental example of this kind of reasoning was Deligne's approach to the Hodge theory of
simplicial schemes using the notion of ``mixed Hodge complex'' \cite{HodgeIII}.

In the nonabelian or homotopical-algebra case, Quillen's notion of closed model category formulates a good collection of requirements that can be made
on the intermediate categorical data. Quillen in \cite{Quillen} asked for a general structure which would encapsulate all of the higher homotopical data.
In one way of looking at it, the answer lies in the notion of {\em higher category}. Quillen had already provided this answer with his definition of ``simplicial
model category'', wherein the  simplicial subcategory of cofibrant and fibrant objects
provides a homotopy invariant higher categorical structure. 
As later became clear with the work of Dwyer and Kan, this simplicial
category contains exactly the right information.  
The notion of Quillen model category is still one of the best ways of approaching the problem of calculation with homotopical objects, so much so that we
adopt it as a basic language for dealing with notions of higher categories. 

Bondal and Kapranov introduced the idea of {\em enhanced derived categories} 
\cite{KapranovInv} \cite{BondalKapranov}, whereby
the usual derived category, which is the Gabriel-Zisman localization of the category of
complexes, is replaced by a {\em differential graded (dg) category} containing the required
higher homotopy information. The notion of dg-category actually appears near the
end of Gabriel-Zisman's book \cite{GabrielZisman} (where it is compared with the notion of $2$-category),
and it was one of the motivations for Kelly's theory
of enriched categories \cite{Kelly}. The notion of dg-category, now further developped by Keller \cite{Keller}, Tabuada \cite{Tabuada}, Stanculescu \cite{Stanculescu}, Batanin \cite{BataninAinfty}, Moriya \cite{Moriya} and others, is one possible answer to the search for higher
categorical structure in the $k$-linear case, pretty much analogous to the notion
of strict simplicial category. The corresponding weak notion is that of {\em $A_{\infty}$-category} used for example by Fukaya \cite{Fukaya}
and Kontsevich \cite{Kontsevich}. This 
definition is based on Stasheff's notion of $A_{\infty}$-algebra  \cite{Stasheff}, 
which is an example of the
passage from delooping machinery to higher categorical theories. 

In the far future one could imagine starting directly with a notion of higher category and bypassing the model-category step entirely, but for now this
raises difficult questions of bootstrapping. Lurie has taken this kind of program a long way in \cite{LurieTopos} \cite{LurieAlgebra}, using the notion of {\em quasicategory}
as his basic higher-categorical object. But even there, the underlying model category theory remains important. The reader is invited to reflect on this
interesting problem. 

The original example of the $2$-category of categories, suggests using $2$-categories and their eventual iterative generalizations, as higher
categorical structures. This point of view occured as early as Gabriel-Zisman's book, where they introduce a $2$-category enhancing the structure of
$\Ho (\Top )$ as well as its analogue for the category of complexes, and proceed to use it to treat questions about homotopy groups. 

Benabou's monograph
\cite{Benabou} introduced the notion of {\em weak $2$-category}, as well as various notions
of weak functor. These are also related to Grothendieck's notion of {\em fibered 
category} in that a fibered category may be viewed as some kind of weak functor from
the base category to the $2$-category of categories.

Starting with Benabou's book, it has been clear that there would be two types of generalization from $2$-categories to $n$-categories. The {\em strict $n$-categories}
are defined recurrently as categories enriched over the category of strict $n-1$-categories. By the Eckmann-Hilton argument, these don't contain enough objects,
as we shall discuss in Chapter \ref{nonstrict1}. For this reason, these are not our main objects of study and we will use the terminology {\em strict $n$-category}. The relative
ease of defining strict $n$-categories nonetheless makes them attractive for learning some
of the basic outlines of the theory, the starting point of Chapter \ref{strict1}. 

The other generalization would be to consider {\em weak $n$-categories} also called {\em lax $n$-categories}, and which we usually call just ``$n$-categories'', 
in which the composition would be associative
only up to a natural equivalence, and similarly for all other operations. The requirement that all equalities between sequences of operations
be replaced by natural equivalences at one
level higher, leads to a combinatorial explosion because the natural equivalences themselves are to be considered as operations. For this reason, the theory 
of weak $3$-categories developped by Gordon, Power, Street \cite{GordonPowerStreet}
following the
path set out by Benabou for $2$-categories in \cite{Benabou}, was already very complicated; for $n=4$ it became next to impossible (see however \cite{Trimble4})
and development of this line stopped there.

In fact, the problem of defining and studying the higher operations which are needed
in a weakly associative category, had been considered rather early on by
the topologists who noticed that
the notion of ``$H$-space'', that is to say a space with an operation which provides
a group object in the homotopy category, was insufficient to capture the data contained
in a loop space. 
One needs to specify, for example, a homotopy of associativity 
between $(x,y,z)\mapsto x(yz)$ and $(x,y,z)\mapsto (xy)z$. This ``associator'' 
should itself be subject to some kind of higher associativity laws, called 
{\em coherence relations},
involving composition of four or more elements.

One of the first discussions of the resulting higher coherence structures
was Stasheff's notion of $A_{\infty}$-algebra \cite{Stasheff}. This was placed in the realm
of differential graded algebra, but not long thereafter the notion of ``delooping machine''
came out, including MacLane's notion of PROP, then May's operadic and Segal's simplicial 
delooping machines. 

In Segal's case, the higher coherence relations come about by requiring
not only that $\sigma _2$ be a weak equivalence, but that all of the ``Segal maps''
$$
\sigma _m : \pA _m \rightarrt \pA _1\times \ldots \times \pA _1
$$
given by $\sigma _m = (f_{01}^{\ast}, f_{12}^{\ast}, \ldots , f_{m-1,m}^{\ast})$ should be weak homotopy equivalences. This was iterated by Dunn \cite{Dunn}. 
In the operadic viewpoint, the coherence relations
come from contractibility of the spaces of $n$-ary operations. 

By the late 1960's and early 1970's, the topologists had their delooping machines
well in hand. A main theme of the present work is that these delooping machines can 
generally lead to definitions of higher categories, but that doesn't seem to have
been done explicitly at the time. A related notion also appeared in the book of
Boardman and Vogt \cite{BoardmanVogt}, 
that of {\em restricted Kan complex}. These objects are now known as
``quasicategories'' thanks to Joyal's work \cite{JoyalQC}.  At that time, in algebraic geometry, an elaborate theory of
derived categories was being developed, but it relied only on $1$-categories which
were the $\tau _{\leq 1}$ of the relevant higher categories. This difficulty was worked
around at all places, by techniques of working with explicit resolutions and
complexes. Illusie gave the definition of {\em weak equivalence of simplicial 
presheaves} which, in retrospect, leads later to the idea of {\em higher stack} via the
model categories of Jardine and Joyal. Somewhere in these works is the idea, which seems
to have been communicated to Illusie by Deligne, of looking at the derived category
of diagrams as a functor of the base category; this was later taken up by Grothendieck and
Cisinski under the name ``derivator'' \cite{GrothendieckDerivateurs} \cite{CisinskiDerivateurs}.  

In 1980, Dwyer and Kan came out with their theory of simplicial localization,
allowing the association of a simplicial category to any pair $(\mM , \mW )$
and giving the higher categorical version of Gabriel-Zisman's theory. 
They developped an extensive theory of simplicial categories, including
several different
constructions of the simplicial localization which inverts
the morphisms of $\mW$ in a homotopical sense. This construction
provides the door passing from the world of categories to the world of higher categories, because even if we start with
a regular $1$-category, then localize by inverting a collection of morphisms, the simplicial localization is in general 
a simplicial category which is not a $1$-category. The simplicial localization maps to the usual or Gabriel-Zisman localization
but the latter is only the $1$-truncation. So, if we want to invert a collection of morphisms in a ``homotopically correct'' way,
we are forced to introduce some kind of higher categorical structure, at the very least the notion of simplicial category. 
Unfortunately, the importance
of the Dwyer-Kan construction doesn't seem
to have been generally noticed at the time. 

During this period, the category-theorists and particularly the Australian school,
were working on fully understanding the theory of strictly associative $n$-categories
and $\infty$-categories. In a somewhat different direction,
Loday introduced the notion of $cat ^n$-group
which was obtained by iterating the internal category construction in a different way,
allowing categories of objects as well as of morphisms.  
Ronnie Brown worked on various aspects of the problem of
relating these structures to homotopy theory: the strictly associative $n$-categories
don't model all homotopy types (Brown-Higgins), whereas the $cat^n$-groups do (Brown-Loday).

A major turning point in the history of higher categories was Alexandre Grothendieck's famous manuscript
{\em Pursuing Stacks}, which started out as a collection of letters to different colleagues with many parts crossed out and rewritten, the whole circulated in mimeographed form. I was lucky to be able to consult a copy in the back room of the Princeton math library, and later
to obtain a copy from Jean Malgoire; a published version edited by Georges Maltsiniotis should appear soon \cite{Grothendieck}. 
Grothendieck introduces the problem of defining a notion of
weakly associative $n$-category, and points out that many 
areas of mathematics could benefit from such a theory, 
explaining in particular how a theory of higher stacks should
provide the right kind of coefficient system for higher nonabelian cohomology.

Grothendieck made important progress in investigating the topology and category theory 
behind this question. He introduced the notion of {\em $n$-groupoid}, an $n$-category
in which all arrows are invertible up to equivalences at the next higher level. 
He conjectured the existence of a  {\em Poincar\'e $n$-groupoid construction}
$$
\Pi _n : \Top \rightarrt n \mbox{-}{\rm Gpd} \subset n \mbox{-}{\rm Cat}
$$
where $n \mbox{-}{\rm Gpd}$ is the collection of weakly associative
$n$-groupoids. He postulated that this
functor should provide an equivalence of homotopy theories between 
$n$-truncated spaces\footnote{A space 
$T$ is {\em $n$-truncated} if $\pi _i(T,t)=0$ for all $i>n$ and all basepoints $t\in T$. The $n$-truncated spaces are the
objects which appear in the Postnikov tower of fibrations, and one can define the truncation $T\rightarrt \tau _{\leq n}(T)$ for any space $T$,
by adding cells of dimension $\geq n+2$ to kill off the higher homotopy 
groups.}
and $n$-groupoids.

In his search for algebraic models for homotopy types, Grothendieck was inspired by one of the  pioneering works in this direction, the notion of {\em $Cat ^n$-groups} of Brown and 
Loday. This
is what is now known as the ``cubical'' approach where the set
of objects can itself have a structure for example of $n-1$-category, so it
isn't
quite the same as the approach we are looking for, commonly called the
``globular'' case.\footnote{Paoli has recently defined a notion of {\em special $Cat^n$-group}
\cite{PaoliAdvances} that imposes the globularity condition weakly.}

Much of ``Pursuing stacks'' is devoted to the more general question of modeling 
homotopy types by algebraic objects such as presheaves on a fixed small category,
developing a theory of ``test categories'' which has now blossomed into a distinct
subject in its own right thanks to the further work of Maltsiniotis \cite{MaltsiniotisAsterisque} and Cisinski \cite{CisinskiAsterisque}. 
One of the questions their theory aims to address is, which presheaf categories provide good
models for homotopy theory. One could ask a similar question with respect to Segal's
utilisation of $\Delta$, namely whether other categories could be used instead. We don't currently have any good information about this. As a start,
throughout the book we will try to point
out in discussion and counterexamples the main places where special properties of $\Delta$ are used. 

In the parts of ``Pursuing stacks'' about $n$-categories, the following
theme emerges: the notion of $n$-category with strictly
associative composition, is not sufficient. This is seen from the fact that
strictly associative $n$-categories satisfying a weak groupoid condition, do not serve
to model homotopy $n$-types as would be expected. Fundamentally due to Godement and the Eckmann-Hilton
argument, this observation was refined over time by Brown and Higgins \cite{BrownHiggins}
and Berger \cite{BergerEckmannHilton}. 
We discuss it in some detail in Chapter \ref{nonstrict1}.

Since strict $n$-categories aren't enough, it leads to the question of defining a notion of weak $n$-category, which is the
main subject of our book.
Thanks to a careful reading by Georges 
Maltsiniotis, we now know that Grothendieck's manuscript
in fact contained a definition of weakly associative
$n$-groupoid \cite{MaltsiniotisGroGpd}, and that his 
definition is very similar to Batanin's definition of 
$n$-category \cite{MaltsiniotisGroBat}. Grothendieck enunciated the deceptively simple
rule \cite{Grothendieck}:
\label{grothendieckgroupoidquote}
\begin{quote}
Intuitively, it means that whenever we have two ways of associating to a finite
family $(u_i)_{i\in I}$ of objects of an $\infty$-groupoid, $u_i\in F_{n(i)}$, subjected
to a standard set of relations on the $u_i$'s, an element of some $F_n$, in
terms of the $\infty$-groupoid structure only, then we have automatically a ``homotopy''
between these built in in the very structure of the $\infty$-groupoid, provided
it makes sense to ask for one \ldots 
\end{quote}
The structure of this as a definition was not immediately evident upon any initial reading,
all the more so when one takes into account the directionality of arrows,
so ``Pursuing stacks'' left open the problem of giving a good definition of 
weakly associative $n$-category.

Given the idea that an equivalence $\Pi _n$ between homotopy $n$-types and $n$-groupoids should exist, it becomes possible to think of replacing the notion of $n$-groupoid by
the notion of $n$-truncated space. This motivated Joyal to define a closed model
structure on the category of simplicial sheaves, and Jardine to extend this to
simplicial presheaves. These theories give an approach to the notion of 
$\infty$-stack, and were used by Thomason, Voevodsky, Morel and others in $K$-theory. 

Also explicitly mentioned in ``Pursuing stacks'' was the limiting case $n=\omega$,
involving $i$-morphisms of all degrees $0\leq i <\infty$. Again, an $\omega$-groupoid
should correspond, via the inverse of a Poincar\'e construction $\Pi _{\omega}$, to a full homotopy type up to weak equivalence. 

We can now get back to the discussion of simplicial categories. These are categories enriched over spaces, and applying $\Pi _{\omega}$ (which is supposed to be compatible
with products) to the morphism spaces, we can think of simplicial categories
as being categories enriched over $\omega$-groupoids. Such a thing is itself an $\omega$-category
$\pA $, with the property that the morphism $\omega$-categories $\pA (x,y)$ are groupoids. In other words, the $i$-morphisms are invertible for $i\geq 2$,
but not necessarily for $i=1$. Jacob Lurie introduced the terminology {\em $(\infty , 1)$-categories} for these things, where more generally
an $(\infty ,n)$-category would be an $\omega$-category such that the $i$-morphisms are invertible up to equivalence, for $i>n$. 
The  point of this 
discussion---of notions which have not yet been defined---is to say
that the notion of simplicial category is a perfectly good substitute for the notion of $(\infty , 1)$-category even if we don't know what
an $\omega$-category is in general. 

This replacement no longer works if we want to look at $n$-categories with noninvertible morphisms at levels $\geq 2$, or somewhat similarly,
$(\infty ,n)$-categories for $n\geq 2$. Grothendieck doesn't seem to have
been aware of Dwyer and Kan's work, just prior to ``Pursuing stacks'', on simplicial categories;\footnote{Paradoxically, Grothendieck's unpublished manuscript is responsible
in large part for the regain of interest in Dwyer and Kan's published papers!} 
however he was well aware that the notions of $n$-category for small values of $n$ had been extensively investigated earlier in
Benabou's book about $2$-categories \cite{Benabou}, and Gordon, Powers and Street on $3$-categories \cite{GordonPowerStreet}. 
The combinatorial explosion inherent in these  explicit theories
was why Grothendieck
asked for a different form of definition which could work in general.

As he forsaw  in a vivid passage \cite[First letter, p. 16]{Grothendieck},
there are currently many different definitions of $n$-category. This started with Street's proposal
in \cite{Street}, of a definition of weak $n$-category
as a simplicial set satisfying a certain variant of the Kan condition where one
takes into account the directions of arrows, 
and also using the idea of
``thinness''. His suggestion, in retrospect
undoubtedly somewhat similar to
Joyal's iteration of the notion of quasicategory, wasn't worked out at the time,
but has recieved renewed interest, see Verity \cite{Verity} for example. 

The Segal-style approach to weak topological categories
was introduced by Dwyer, Kan, Smith \cite{DKS} and Schw\"{a}nzl, Vogt \cite{Vogt},
but the fact that they immediately
proved a rectification result relating Segal categories back to
strict simplicial categories, seems to have slowed down their further consideration of this
idea. Applying Segal's idea seems to have been the topic of a letter from Breen to
Grothendieck in 1975, see page \pageref{grobreenquote} below.

Kapranov and Voevodsky in \cite{KapranovVoevodsky}
considered a notion of ``Poincar\'e $\infty$-groupoid'' which
is a strictly associative $\infty$-groupoid but where the arrows are
invertible only up to equivalence.  It now appears likely that their constructions 
should best be
interpreted using some kind of weak unit condition \cite{Kock}.

At around the same time in the mid-1990's, three distinct approaches to defining weak $n$-categories came out: Baez and Dolan's approach used {\em opetopes} \cite{BaezDolanLetter} \cite{BaezDolanIII}, 
Tamsamani's approach used iteration of the Segal delooping machine \cite{TamsamaniThesis}
\cite{Tamsamani}, and 
Batanin's approach used
{\em globular operads} \cite{Batanin} \cite{Batanin2}. 
The Baez-Dolan and Batanin approaches will be discussed in Chapter  \ref{operadic1}. 

The work of Baez and Dolan was motivated by a far-reaching program of conjectures
on the relationship between $n$-categories and physics \cite{BaezDolan} \cite{catBD},
which has led to important developments most notably the recent proof by Hopkins and Lurie.

In relationship with Grothendieck's manu\-script, as we pointed out above,
Batanin's approach is the one which most closely resembles what Grothendieck was asking for,
indeed Maltsiniotis generalized the
definition of $n$-groupoid which he found in ``Pursuing stacks'', to a definition of $n$-category which is similar to Batanin's one \cite{MaltsiniotisGroBat}. 

In the subsequent period, a number of other definitions have appeared, and people have begun working more seriously on the approach which had been suggested by
Street. 
Batanin, in mentioning the letter from Baez and Dolan to Street \cite{BaezDolanLetter}, also 
points out that Hermida, Makkai and Power have worked on the opetopic ideas, leading to \cite{HermidaMakkaiPower}.
M. Rosellen suggested in 1996 to give a  version of the 
Segal-style definition, using the theory of operads.
He didn't concretize this but Trimble gave a definition along these lines, now playing
an important role in work of Cheng \cite{ChengComparison}. Further ideas include those
of Penon, Leinster's multicategories, and others.
Tom Leinster has collected together ten different definitions in the useful compendium \cite{Leinster}. The somewhat mysterious \cite{Kondratieff} could also be pointed out.
In the simplicial direction Rezk's complete Segal spaces 
\cite{Rezk} can be iterated as suggested by Barwick \cite{RezkCartesian},
and Joyal proposes an iteration of the method of quasicategories \cite{JoyalTheta}.

We shall discuss the simplicial definitions in Chapter \ref{simplicial1} and 
the operadic definitions in Chapter
\ref{operadic1}. One of the main tasks in the future will be to understand the relationships between all of these approaches. Our goal here is more down-to-earth: we would
like to develop the tools necessary for working with Tamsamani's $n$-categories. We hope that similar tools can be developped for the other approaches,
making an eventual comparison theory into a powerful method whereby the particular advantages of each definition could all be put in play at the same time. 

Tamsamani defined the Poincar\'e $n$-groupoid functor for his notion of $n$-category,
and showed Grothendieck's conjectured equivalence with the theory of homotopy $n$-types \cite{Tamsamani}. The same has also been done for Batanin's theory,
by Berger in \cite{BergerCellularNerve}.

It is interesting to note that the two main ingredients in Tamsamani's approach, the 
multisimplicial nerve construction and Segal's delooping machine, are both mentioned in ``Pursuing stacks''.
In particular,
Grothendieck reproduces a letter from himself to Breen dated July 1975, in which Grothendieck acknowledges
having recieved a proposed definition of non-strict $n$-category from Breen, a
definition which according to {\em loc. cit} ``...has certainly the merit of
existing...''. It is not clear whether this proposed construction was ever
worked out. Quite apparently, Breen's suggestion for using Segal's delooping
machine must have gone along the lines of what we are doing here. 
Rather than taking up
this direction, Grothendieck elaborated a general {\em ansatz} whereby $n$-categories would
have various different composition operations, and natural equivalences between any two 
natural compositions with the same source and target, an idea now fully developed in
the context of Batanin's and related definitions.

Once one or more points of view for defining $n$-categories are in hand, 
the main problem which needs to be considered
is to obtain---hopefully within the same point of view---an $n+1$-category $nCAT$
parametrizing the $n$-categories of that point of view.  This problem,
already clearly posed in ``Pursuing stacks'', is one of our main goals in
the more technical central part of the book, for one model.

It turns out that Quillen's technique of model categories, subsequently
deepened by several generations of mathematicians, is a great way of attacking this problem.
It is by now well-known that closed model categories provide an excellent environment
for studying homotopy theory, as became apparent from the work of Bousfield,
Dwyer and Kan on closed model categories of diagrams, and the generalization
of these ideas by Joyal, Jardine, 
Thomason and Voevodsky who used model categories to study simplicial presheaves under
Illusie's condition of weak equivalence. In the Segal-style paradigm of weak enrichment,
we look at functors $\Delta ^o \rightarrt (n-1){\rm Cat}$, so we are certainly also studying diagrams and it is reasonable to expect the notion of model category to bring
some of the same benefits as for the above-mentioned theories. 

To be more precise about this motivation, recall from ``Pursuing stacks'' 
that 
$nCAT$ should be an $n+1$-category whose objects are in one-to-one correspondence 
with the $n$-categories of a given universe. The structure of $n+1$-category 
therefore consists of specifying the morphism 
objects $\Hom _{nCAT}(\pA ,\pB )$ which should themselves be $n$-categories parametrizing
``functors'' (in an appropriate sense) from $\pA$ to $\pB$.

In the explicit theories for $n=2,3,4$ this is one of the places where a combinatorial
explosion takes place: the functors from $\pA $ to $\pB$ have to be taken in a weak sense,
that is to say we need a natural equivalence between the image of a composition and
the composition of the images, together with the appropriated coherence data at all levels. 

The following simple example shows that, even if we were to consider only
strict $n$-categories, the strict morphisms are not
enough. Suppose $G$ is a group and $V$ an abelian group and we set $\pA $
equal to the category with one object and group of automorphisms $G$, and
$B$ equal to the strict $n$-category with only one $i$-morphism for $i<n$ and
group $V$ of $n$-automorphisms of the unique $n-1$-morphism; then for $n=1$
the equivalence classes of strict morphisms from $\pA$ to $\pB$ are the elements of
$H^1(G,V)$ so we would expect to get $H^n(G,V)$ in general, but for $n>1$ there
are no nontrivial strict morphisms from $\pA $ to $\pB$.
So some kind of weak notion of functor is needed. 

Here is where the notion of 
model category comes in: one can view this situation as being similar to the problem 
that usual maps between simplicial sets are generally too rigid and don't reflect
the homotopical maps between spaces. Kan's fibrancy condition and Quillen's formalization
of this in the notion of model category, provide the solution: we should require the
target object to be fibrant and the source object to be cofibrant in an appropriate
model category structure. Quillen's axioms serve to guarantee that the notions of
cofibrancy and fibrancy go together in the right way. So, in the application
to $n$-categories we would like to define a model structure and then 
say that the usual maps $\pA \rightarrt \pB$
strictly respecting the structure, are the right ones provided that $\pA $ is cofibrant
and $\pB$ fibrant. 

To obtain $nCAT$ a further property is needed, indeed we are not just looking to find
the right maps from $\pA $ to $\pB$ but to get a morphism object $\Hom _{nCAT}(\pA ,\pB)$
which should itself be an $n$-category. It is natural to apply the idea of ``internal $\uHom$'', that is to put
$$
\Hom _{nCAT}(\pA ,\pB ):= \uHom (\pA ,\pB )
$$
using an internal $\uHom$ in our model category. For our purposes, it is sufficient to
consider $\uHom$ adjoint to the direct product operation, in other words a map 
$$
\pE \rightarrt \uHom (\pA ,\pB )
$$
should be the same thing as a map $\pE\times \pA\rightarrt \pB$. This obviously implies
imposing further axioms on the closed model structure, in particular
compatibility between $\times$ and cofibrancy since the direct product is used
on the source side of the map. It turns out that the required
axioms are already well-known in the notion of {\em monoidal model category} \cite{Hovey},
which is a model category provided with an additional operation $\otimes$, and 
certain axioms of compatibility with the cofibrant objects. In our case, the operation
is already given as the direct product $\otimes = \times$ of the model category,
and a model category which is monoidal for the direct product operation will be called
{\em cartesian} (Chapter \ref{cartmod1}). 

In the present book, we are concentrating on Tamsamani's approach to $n$-categories,
which in \cite{descente} was modified to `` Segal $n$-categories''  
in the course of discussions with Andr\'e Hirschowitz.
In Tamsamani's theory an $n$-category
is viewed as a category enriched over $n-1$-categories, using Segal's machine to
deal with the enrichment in a homotopically weak way. 

In Regis Pelissier's thesis, following a question posed by Hirschowitz, 
this idea was pushed to a next level: to study weak Segal-style
enrichment over a more general model category, with the aim of making the iteration
formal. A small link was missing in this process at the end of \cite{Pelissier},
essentially because of an error in \cite{svk} which Pelissier discovered. He provided
the correction when the iterative procedure is applied to the model category of
simplicial sets. But in fact, his patch applies much more generally if we just 
consider the operation of functoriality under change of model categories.  

This is what we will be doing here. But instead of following Pelissier's argument too closely,
some aspects will be set into a more general discussion of certain kinds of 
left Bousfield localizations. The idea of breaking down the construction into 
several steps including a main step of left Bousfield localization, is due to
Clark Barwick.  

The Segal $1$-categories are, as was originally proven in \cite{DKS}, equivalent to strict
simplicial categories. Bergner has shown that this equivalence takes the form of a
Quillen equivalence between model categories \cite{BergnerThreeModels}. However, the
model category of simplicial categories is not appropriate for the considerations
described above: it is not cartesian, indeed the product of two cofibrant simplicial
categories will not be cofibrant.\footnote{This remark also applies to the projective model structure for weakly enriched Segal-style categories, whereas on the other hand the projective structure is much more practical for calculating maps.}  
It is interesting to imagine several possible ways
around this problem: one could try to systematically apply the cofibrant replacement
operation; this would seem to lead to a theory very similar to the consideration of
Gray tensor products of Leroy \cite{Leroy} and Crans \cite{CransTensor}; or one could
hope for a general construction replacing a model category by a cartesian one
(or perhaps, given
a model category with monoidal structure incompatible with cofibrations,
construct a monoidal model category in some sense equivalent to it). 

As Bergner pointed out,
the theories of simplicial categories and Segal categories are also equivalent to
Charles Rezk's theory of {\em complete Segal spaces}. As we shall discuss further
in Chapter \ref{simplicial1}, Rezk requires that the Segal maps
be weak equivalences, but rather than having $\pA_0$ be a discrete simplicial set
corresponding to the set of objects, he asks that $\pA_0$ be a simplicial set weakly
equivalent to the ``interior'' Segal groupoid of $\pA$. Barwick pointed out that Rezk's theory could also be iterated, and Rezk's recent preprint \cite{RezkCartesian} 
shows that the resulting model category is cartesian. So, this route also leads to
a construction of $nCAT$ and can serve as an alternative to what we are doing here. 
It should be possible to extend Bergner's comparison result to obtain equivalences
between the iterates of Rezk's theory and the iterates we consider here. If our current
theory is perhaps simpler in its treatment of the set of objects, Rezk's theory has
a better behavior with respect to homotopy limits.

As more different points of
view on higher categories are up and running, the comparison problem will be posed: to find an
appropriate way to compare different points of view on $n$-categories and (one
hopes) to say that the various points of view are equivalent and in particular
that the various $n+1$-categories $nCAT$ are equivalent via these comparisons.
Grothendieck gave a vivid description of this problem (with remarkable foresight,
it would seem
\cite{Leinster}) in the first letter of \cite{Grothendieck}. 
He pointed out that 
it is not actually clear  what type of general setup one should use for
such a comparison theory. Various possibilities would include the model category
formalism, or the formalism of $(\infty ,1)$-categories starting with
Dwyer-Kan localization and moving through Lurie's theory. 

Within the domain of simplicial theories, we have mentioned
Bergner's comparison
between three different approaches to $(\infty , 1)$-categories \cite{BergnerThreeModels}. 
A further comparison of these theories with
quasicategories is to be found in Lurie \cite{LurieTopos}. 

A recent result due to Cheng \cite{ChengComparison} gives a comparison between
Trimble's definition and 
Batanin's definition (with some modifications on both sides due to Cheng and Leinster). 
Batanin's approach used operads
more as a way of encoding general algebraic structures, and is the closest to 
Grothendieck's original philosophy.
While also
operadic, Trimble's definition is much closer to the philosophy we are developing
in the present book, whereby one goes from topologists' delooping machinery (in his
case, operads) to an iterative theory of $n$-categories.  It is to be hoped that Cheng's result
can be expanded in various directions to obtain comparisons between a wide range
of theories, maybe using Trimble's definition as a bridge towards the 
simplicial theories. This should clearly be pursued in the near future, but 
it would go beyond the scope of the present work.

We now turn to the question of potential applications. 
Having a good theory of $n$-categories should open up the possibility to pursue
any of the several programs such as that outlined by Grothendieck
\cite{Grothendieck}, the generalization to $n$-stacks and $n$-gerbs of the work
of Breen \cite{Breen}, or the program of Baez and Dolan in
topological quantum field theory \cite{BaezDolan}. Once the theory of
$n$-stacks is off the ground this will give an algebraic approach to the
``geometric $n$-stacks'' considered in \cite{geometricN}.

As the title indicates, Grothendieck's manuscript was intended to develop a foundational
framework for the theory of higher stacks. In turn, higher stacks should be the natural
coefficients for nonabelian cohomology, the idea being to generalize Giraud's \cite{Giraud}
to $n\geq 3$. 

The example of diagrams of spaces translates, via the construction $\Pi _n$, to a notion of {\em diagram of $n$-groupoids}. This is a strict version of the notion 
of {\em $n$-prestack in groupoids} which would be a weak functor from the base category $\Psi$ to the $n+1$-category $GPD^{n}$ of $n$-groupoids. Grothendieck 
introduced the notion of {\em $n$-stack} which generalizes to $n$-categories the classical notion of stack. A full discussion of this theory would go
beyond the scope of the present work: we are just trying to set up the $n$-categorical foundations first. The notion of $n$-stack, maybe with $n=\infty$,
has applications in many areas as predicted by Grothendieck.

Going backwards along $\Pi _n$, it turns out that diagrams of spaces or equivalently simplicial presheaves, serve as a very adequate replacement \cite{JoyalLetter} \cite{Jardine} \cite{ThomasonENS}
\cite{Voevodsky} \cite{MorelVoevodsky}. 
So, the notion of $n$-categories as a prerequisite for higher stacks has proven
somewhat illusory.  
And in fact, the model category theory developed for simplicial presheaves 
has been useful for attacking the theory
of $n$-categories as we do here, and also for going from a theory of $n$-categories to
a theory of $n$-stacks, as Hollander has done for $1$-stacks \cite{Hollander} and 
as Hirschowitz and I did for $n$-stacks in \cite{descente}.

An $n$-stack on 
a site $\mX$ will be a morphism $\mX \rightarrt nCAT$. This requires a construction for
the $n+1$-category $nCAT$, together with the appropriate notion of morphism between
$n+1$-categories. The latter is almost equivalent to knowing how to construct the
$n+2$-category $(n+1)CAT$ of $n+1$-categories. From this discussion the need for an
iterative approach to the theory of $n$-categories becomes clear.

My own favorite application of stacks is that they lead in turn to
a notion of {\em nonabelian cohomology}. Grothendieck 
says \cite{Grothendieck}:
\begin{quote}
Thus $n$-stacks, relativized over a topos to ``$n$-stacks over $X$'', are viewed primarily
as the natural ``coefficients'' in order to do (co)homological algebra of dimension 
$\leq n$ over $X$.
\end{quote}

The idea of using higher categories for nonabelian cohomology goes back to Giraud \cite{Giraud}, and had been extended to the cases $n=2,3$ by Breen somewhat more
recently \cite{BreenAsterisque}. Breen's book motivated us to proceed to the case of
$n$-categories at the beginning of Tamsamani's thesis work. 

Another utilisation of the notion of $n$-category is to model homotopy types. For
this to be useful one would like to have as simple and compact a definition as possible,
but also one which lends itself to calculation. The simplicial approach developped here
is direct, but it is possible that the operadic approaches which will be mentioned in
Chapter \ref{operadic1} could be  more amenable to topological computations. 
An iteration of the classical Segal delooping machine has been considered by Dunn \cite{Dunn}.

The Poincar\'e $n$-groupoid of a space is a generalization of the Poincar\'e groupid $\Pi _1(X)$, a basepoint-free version of the fundamental group $\pi _1(X)$ popularized by
Ronnie Brown \cite{RBrown}. Van Kampen's theorem allows for computations
of fundamental groups, and as Brown has often pointed out, it takes a particularly 
nice form when written in terms of the Poincar\'e groupoid: it says that if 
a space $X$ is written as a pushout $X=U\cup ^WV$ then the Poincar\'e groupoid is a
pushout in the $2$-category of groupoids:
$$
\Pi _1(X)=\Pi _1(U)\cup ^{\Pi _1(W)}\Pi _1(V).
$$
This says that $\Pi _1$ commutes with colimits. 

Extending this theory to the case of Poincar\'e $n$-groupoids is one of the motivations for
introducing colimits and indeed the whole model-categoric machinery for $n$-categories.
We will then be able to write, in case of a pushout of spaces $X=U\cup ^WV$,
$$
\Pi _n(X)=\Pi _n(U)\cup ^{\Pi _n(W)}\Pi _n(V).
$$
Of course the pushout diagram of spaces should satisfy some excision
condition as in the original
Van Kampen theorem, and this may be abstracted by refering to simplicial sets instead. 

The homotopy theory and nonabelian cohomology motivations may be combined by looking for a higher-categorical
theory of {\em shape}. For a space $X$ we can define the {\em nonabelian cohomology
$n$-category}
$H(X,A)$ with coefficients in an $n$-stack $\Ff$ over $X$. This applies in
particular to the constant stack $A_X$ associated to an $n$-category $A$. 
The functor 
$$
A\mapsto H(X,A_X)
$$
is co-represented by the universal element 
$$
\eta _X\in H(X,\Pi _n(X)_X),
$$
giving a way of characterizing $\Pi _n(X)$ by universal property. 
This essentially tautological observation paves the way for more nontrivial shape theories,
consisting of an $n$-category $COEFF$ and a functor ${\rm Shape}(X):COEFF \rightarrt COEFF$.
A particularly useful version is when $COEFF$ is the $n$-category of certain 
$n$-stacks over a site $\mY$, and ${\rm Shape}(X)(\Ff )= \uHom (\Pi _n(X)_{\mY} , \Ff )$
where $\Pi _n(X)_{\mY}$ denotes the constant stack on $\mY$ with values
equal to $\Pi _n(X)$. This leads to subjects generalizing Malcev completions and 
rational homotopy
theory \cite{Hain} \cite{HainMalcev} \cite{HainRelativeMalcev}, 
such as the schematization of homotopy
types \cite{Toen} \cite{KaPaTo} \cite{Pridham} \cite{Moriya}, de Rham shapes and nonabelian Hodge theory. 

One of the main advantages to a theory of higher categories, is that the notions
of {\em homotopy limit} and {\em homotopy colimit}, by now classical in algebraic
topology, become {\em internal} notions in a higher category. Indeed, they become 
direct analogues of the notions of limit and colimit in a usual $1$-category,
with corresponding universal properties and so on. This has an interesting application to
the ``abelian'' case: the structure of triangulated category is automatic once we
know the $(\infty ,1)$-categorical structure. This was pointed out by Bondal and Kapranov
in the dg setting \cite{BondalKapranov}: their {\em enhanced triangulated categories} are just dg-categories satisfying some further axioms; the structure of triangles comes from the dg structure. 
Historically one can trace this observation back to the end of Gabriel-Zisman's book \cite{GabrielZisman}, although nobody seems to have noticed it until rediscovered by
Bondal and Kapranov. 

I first learned of the notion of ``2-limit'' from the paper of Deligne and
Mumford \cite{DeligneMumford}, where it appears at the beginning with very
little explanation (their paper should also be added to the list of 
motivations for developing the theory
of higher stacks). Several authors have since considered $2$-limits and $2$-topoi,
originating with Bourn
\cite{BournDitopos} and continuing recently with Weber \cite{Weber} for example.

The notions of homotopy limits and colimits internalized in an $(\infty ,1)$-category
have now recieved an important foundational formulation with Lurie's work on $\infty$-topoi
\cite{LurieTopos}. 

Power has given an extensive discussion of the motivations for higher categories
stemming from logic and computer science, in \cite{Power}. He points out the role played by weak limits. Recently, Gaucher, Grandis and others have
used higher categorical notions to study directed and concurrent processes 
\cite{Gaucher} \cite{Grandis}. It would be
interesting to see how these theories interact with the notion of $\infty$-topos.

Recall that gerbs played an important role in descent theory and non-neutral tannakian 
categories \cite{lnm900}. 
Current developments where the notion of higher category is more or less essential
on a foundational level, include ``derived algebraic geometry'' and higher tannakian theory. 
It would go beyond our present scope to discuss these here but the reader may search for
numerous references. 

Stacks and particular gerbes of higher groupoids 
have found many interesting applications in the
mathematical physics literature, starting with explicit considerations for $1$- and $2$-gerbes. Unfortunately it would go beyond our scope to list all of these. However, one of the main contributions from mathematical physics has been to 
highlight the utility of higher categories which are not groupoids, in which there can
be non-invertible morphisms. Explicit first cases come about when we consider {\em monoidal categories}: they may be considered as $2$-categories with a single object. And then {\em braided monoidal categories} may be considered as $3$-categories with a single object
and a single $1$-morphism, where the braiding isomorphism comes from the Eckmann-Hilton
argument. These entered into the vast program of research on combinatorial quantum field theories and knot invariants---again the reader is left to fill in the references here. 

John Baez and Jim Dolan provided a major impetus to the theory of higher categories, by
formulating a series of conjectures about how the known relationships between
low-dimensional field theories and $n$-categories for small values of $n$, should generalize
in all dimensions \cite{BaezDolan} \cite{BaezIntroduction} \cite{BaezDolanIII} \cite{catBD}.
On the topological or field-theoretical side, they conjecture the 
existence of a  $k$-fold monoidal $n$-category (or equivalently, a $k$-connected $n+k$-category) representing $k$-dimensional
manifolds up to cobordism, where the higher
morphisms should correspond to manifolds with corners. On the $n$-categorical side,
they propose a notion of {\em $n$-category with duals} in which all morphisms should have
internal adjoints. Then, their main conjecture relating these two
sides is that the cobordism $n+k$-category should be the universal $n+k$-category with duals
generated by a single morphism in degree $k$. The specification of a field theory is 
a functor from this $n+k$-category to some other one, and it suffices to specify the image
of the single generating morphism. They furthermore go on to investigate possible candidates
for the target categories of such functors, looking at higher Hilbert spaces and other
such things. We will include some discussion (based on \cite{BBDSH}) of one of Baez and Dolan's preliminary conjectures,
the {\em stabilization hypothesis}, in an ulterior version of the present manuscript.

The Baez-Dolan conjectures step outside of the realm of $n$-groupoids, so they
really require an approach which can take into account non-invertible morphisms.
In their ``$n$-categories with duals'', they generalize the fact that 
the notion of adjoint functor can be expressed in $2$-categorical
terms within the $2$-category $1CAT$. Mackaay describes the application of internal adjoints
to $4$-manifold invariants in \cite{Mackaay}. 
The notion of adjoint generalizes within an $n$-category
to the notion of dual of any $i$-morphism for $0<i<n$. At the top level of $n$-morphisms,
the dual operation should  either be: ignored; imposed as additional structure; or
pushed to $\infty$ by considering directly the theory of $\infty$-categories. 
Of course, a morphism which is really
invertible is automatically dualizable and its dual is the same as its inverse, so 
the interesting $n$-categories with duals have to be ones which are not $n$-groupoids.

As 
Cheng has pointed out \cite{ChengInfty}, 
in the last case one obtains a structure which looks algebraically
like an $\infty$-groupoid, so the distinction between invertible and dualizable
morphisms should probably be considered as
an additional more analytic structure in itself. 
We don't yet have the tools to 
fully investigate the theory of $\infty$-categories. Further comments on these issues
will be made in Section \ref{sec-towards}. 

In a very recent development, Hopkins and Lurie have announced a proof of a major part
of the Baez-Dolan conjectures, saying that the category of manifolds with appropriate
corners, and cobordisms as $i$-morphisms, is the universal $n$-category with
duals generated by a single element. This universal property allows one to define
a functor from the cobordism $n$-category to any other $n$-category with duals,
by simply specifying a single object. I hope that some of the techniques presented here
can help in understanding this fascinating subject.


\chapter{Strict $n$-categories}
\label{strict1}

Classically, the first and easiest notion of higher category was that of 
{\em strict $n$-category}. We review here some basic definitions, as they introduce
important notions for weak $n$-categories. In Chapter \ref{nonstrict1} we will point out
why the strict theory is generally considered not to be sufficient. 

In the current chapter only, {\em all $n$-categories are meant
to be strict $n$-categories}. For this reason we try to put in the adjective
``strict'' as much as possible when $n>1$; but in any case, the very few
times that we speak of weak $n$-categories, this will be explicitly stated.
We mostly restrict our attention to $n\leq 3$.

In case that isn't already clear, it should be stressed that everything we do
in this section (as well as most of the next and even the subsequent one as
well) is very well known and classical, so much so that I don't know what are
the original references.

To start with, a {\em
strict $2$-category} $\pA $ is a collection of objects $\pA _0$ plus, for each pair of
objects $x,y\in \pA _0$ a category ${\pA} (x,y)$ together with a morphism
$$
\pA(x,y)\times {\pA}(y,z)\rightarrt \pA(x,z)
$$
which is strictly associative in the obvious way; and such that a unit exists,
that is an element $1_x\in \Ob ( {\pA}(x,x))$ with the property that
multiplication by $1_x$ acts trivially on objects of ${\pA}(x,y)$ or
${\pA}(y,x)$ and multiplication by $1_{1_x}$ acts trivially on morphisms of
these categories.

A {\em strict $3$-category} $\pC $ is the same as above but where $\pC (x,y)$ are
supposed to be strict $2$-categories. There is an obvious notion of direct
product of
strict $2$-categories, so the above definition applies {\em mutatis mutandis}.

For general $n$, the well-known definition is most
easily presented by induction on $n$. We assume known the definition of strict
$n-1$-category for $n-1$, and we assume known that the category of strict
$n-1$-categories is closed under direct product. A {\em strict $n$-category}
$\pC $ is then a category enriched \cite{Kelly} over the category of strict
$n-1$-categories. This means that $\pC $ is composed of a {\em set of objects}
$\Ob (\pC )$ together with, for each pair $x,y\in \Ob (\pC )$, a {\em morphism-object}
$\pC (x,y)$ which is a strict $n-1$-category; together with a strictly
associative composition law
$$
\pC (x,y)\times \pC (y,z) \rightarrt \pC (x,z)
$$
and a morphism $1_x: \ast \rightarrt \pC (x,x)$ (where $\ast$ denotes the
final object cf below) acting as the identity for the composition law. The {\em
category of strict $n$-categories} denoted $nStrCat$ is the category whose
objects are as above and whose morphisms are the transformations strictly
perserving all of the structures. Note that $nStrCat$ admits a direct product:
if $\pC $ and $\pC '$ are two strict $n$-categories then $\pC \times \pC '$ is the strict
$n$-category with
$$
\Ob (\pC \times \pC '):= \Ob (\pC ) \times \Ob (\pC ')
$$
and for $(x,x'), \; (y,y') \in \Ob (\pC \times \pC ')$,
$$
(\pC \times \pC ')((x,x'), (y,y')):= \pC (x,y)\times {\pC '}(x',y')
$$
where the direct product on the right is that of $(n-1)StrCat$. Note that the
final object of $nStrCat$ is the strict $n$-category $\ast$ with exactly one
object $x$ and with ${\ast}(x,x)= \ast$ being the final object of
$(n-1)StrCat$.

The induction inherent in this definition may be worked out explicitly to give
the definition as it is presented in \cite{KapranovVoevodsky} for example. In doing this one
finds that underlying a strict $n$-category $\pC $ are the sets $\Mor ^i(\pC )$ of {\em
$i$-morphisms} or {\em $i$-arrows}, for $0\leq i\leq n$. The set of $0$-morphisms
is by
definition the set of objects $\Mor ^0(\pC ):= \Ob (\pC )$, and $\Mor ^i(\pC )$ is the disjoint union over all
pairs $x,y$ of the $\Mor ^{i-1}(\pC (x,y))$. 
These fit together in a diagram called a {\em (reflexive) globular set}:
$$
\begin{diagram}
\cdots \Mor ^{i+1}(\pC ) &\pile{\rightarr^s \\ \leftarr \\ \rightarr_t} & 
\Mor ^{i}(\pC ) &\pile{\rightarr^s \\ \leftarr \\ \rightarr_t} & \Mor ^{i-1}(\pC )\cdots 
\Mor ^{1}(\pC )&\pile{\rightarr^s \\ \leftarr \\ \rightarr_t} & \Mor ^{0}(\pC )
\end{diagram}
$$
where the rightward maps are the {\em source} and {\em target} maps
and
the leftward maps are the {\em identity} maps\footnote{The adjective ``reflexive'' refers to the inclusion of these leftward ``identity'' maps; a {\em non-reflexive globular set}
would have only the $s$ and $t$.}. These may be defined inductively using
the definition we have given of $\Mor ^i(\pC )$. 
The structure of strict $n$-category on this underlying globular set is determined
by further {\em composition laws} at each stage: the $i$-morphisms
may be composed with respect to the $j$-morphisms for any $0\leq j < i$,
operations denoted in \cite{KapranovVoevodsky} by $\ast _j$.
These are partially defined depending on iterations of the source and target maps. 
For a more detailed explanation, refer to the standard references including
\cite{BrownHiggins} \cite{Street} \cite{KapranovVoevodsky} \cite{Benabou}
\cite{GabrielZisman}.

\section{Godement, Interchange or the Eckmann-Hilton argument}
\label{sec-eckhilt}

One of the most important of the axioms satisfied by the various compositions
in a strict $n$-category is variously known under the name of ``Eckmann-Hilton
argument'', ``Godement relations'', ``interchange rules'' etc.
This property comes from the fact that the
composition law
$$
\pC (x,y)\times \pC (y,z)\rightarrt \pC (x,z)
$$
is a morphism with domain the direct product of the two morphism
$n-1$-categories from $x$ to $y$ and from $y$ to $z$. In a direct product,
compositions in the two factors by definition are independent (commute).
Thus, for $1$-morphisms in $\pC (x,y)\times \pC (y,z)$ (where the
composition $\ast _0$ for these $n-1$-categories is actually the composition
$\ast _1$ for $\pC $ and we adopt the latter notation),
we have
$$
(a,b) \ast _1 (c,d) = (a\ast _1c, b \ast _1d).
$$
This leads to the formula
$$
(a\ast _0b) \ast _1 (c\ast _0d) = (a\ast _1c) \ast _0 (b \ast _1d).
$$
This seemingly innocuous formula takes on a special meaning when we start
inserting identity maps. Suppose $x=y=z$ and let $1_x$ be the identity of $x$
which may be thought of as an object of $\pC (x,x)$. Let $e$ denote the
$2$-morphism of $\pC $, identity of $1_x$; which may be thought of as a
$1$-morphism of $\pC (x,x)$. It acts as the identity for both compositions
$\ast _0$ and $\ast _1$ (the reader may check that this follows from the part
of the axioms for an $n$-category saying that the morphism $1_x: \ast
\rightarrt \pC (x,x)$ is an identity for the composition).

If $a, b$ are also
endomorphisms of $1_x$, then the above rule specializes to:
$$
a\ast _1b =
(a\ast
_0e) \ast _1 (e\ast _0b) = (a\ast _1e) \ast _0 (e \ast _1b)  = a \ast _0b.
$$
Thus in this case the compositions $\ast _0$ and $\ast _1$ are the same.
A different ordering gives the formula
$$
a\ast _1b =
(e\ast
_0a) \ast _1 (b\ast _0e) = (e\ast _1b) \ast _0 (a \ast _1e)  = b \ast _0a.
$$
Therefore we have
$$
a\ast _1b= b\ast _1a = a \ast _0b = b\ast _0a.
$$
This argument says, then, that $\Ob  ({\pC (x,x)}(1_x, 1_x))$
is a commutative monoid and the two natural multiplications are the same.

The same argument extends to the whole monoid structure on the $n-2$-category
${\pC (x,x)}(1_x, 1_x)$:

\begin{lemma}
\label{godement}
The two composition laws on the strict $n-2$-category
${\pC (x,x)}(1_x, 1_x)$ are equal, and this law is commutative.
In other words, ${\pC (x,x)}(1_x, 1_x)$ is an abelian monoid-object in
the category $(n-2)StrCat$.
\end{lemma}
\eop

There is a partial converse to the above observation: if the only object is $x$
and the only $1$-morphism is $1_x$ then nothing else can happen and we get the
following equivalence of categories.

\begin{lemma}
\label{scholium}
Suppose $G$ is an abelian
monoid-object in the category $(n-2)StrCat$. Then there is a unique strict
$n$-category $\pC $ such that
$$
\Ob (\pC )= \{ x\}\;\;\; \mbox{and}\;\;\; Mor ^1(\pC )=\Ob (\pC (x,x))=\{
1_x\}
$$
and such that ${\pC (x,x)}(1_x, 1_x)=G$ as an abelian
monoid-object. This construction establishes an equivalence between the
categories of abelian monoid-objects in $(n-2)StrCat$, and the strict
$n$-categories having only one object and one $1$-morphism.
\end{lemma}
{\em Proof:}
Define the strict $n-1$-category $\pU$ with $\Ob (\pU )= \{ u\}$ and $\pU (u,u)=G$
with its monoid structure as composition law. The fact that the composition law
is commutative allows it to be used to define an associative and commutative
multiplication
$$
\pU \times \pU  \rightarrt \pU .
$$
Now let $\pC $ be the strict $n$-category
with $\Ob (\pC )=\{ x\}$ and $\pC (x,x)=\pU $ with the above multiplication. It is
clear that this construction is inverse to the previous one.
\eop

It is clear from the construction (the fact that the multiplication on $\pU$ is
again commutative) that the construction can be iterated any number of times.
We obtain the following ``delooping'' corollary.

\begin{corollary}
\label{iterate}
Suppose $\pC $ is a strict $n$-category with only one object and only one
$1$-morphism. Then there exists a strict $n+1$-category $\pB $ with only one
object $b$ and with $\pB (b,b)\cong \pC $.
\end{corollary}
{\em Proof:}
By the previous lemmas, $\pC $ corresponds to an abelian monoid-object $G$ in
$(n-2)StrCat$. Construct $\pU $ as in the proof of \ref{scholium}, and note that
$\pU $ is an abelian monoid-object in $(n-1)StrCat$. Now apply the result of
\ref{scholium} directly to $\pU $ to obtain $\pB \in (n+1)StrCat$, which will have
the desired property.
\eop

\section{Strict $n$-groupoids}
\label{sec-strgpd}

Recall that a {\em groupoid} is a category where all morphisms are invertible.
This definition generalizes to strict $n$-categories in the following way,
as was pointed out by Kapranov and Voevodsky
\cite{KapranovVoevodsky}. We give a theorem stating that three versions of this definition are
equivalent (one of these equivalences was left as an exercise
in \cite{KapranovVoevodsky}).

Note that, following \cite{KapranovVoevodsky}, we {\em do not} require strict invertibility of
morphisms, thus the notion of strict $n$-groupoid is more general than the
notion employed by Brown and Higgins \cite{BrownHiggins}.

Our discussion is in many ways parallel to the treatment of the groupoid
condition for weak $n$-categories \cite{Tamsamani} to be discussed in the next
chapter, and our treatment in this
section comes in large part from discussions with Tamsamani about this.

We state a few definitions and results which will then be proven all together by induction
on $n$. 

\begin{theorem}
\label{gpdI}
Fix an integer $n\geq 1$. 
Suppose $\pA $ is a strict $n$-category. The following three conditions
are equivalent (and in this case we say that $\pA $ is a {\em strict
$n$-groupoid}). \newline
(1)\, $\pA $ is an $n$-groupoid in the sense of 
\cite{KapranovVoevodsky};
\newline
(2)\, for all $x,y\in \pA $, $\pA (x,y)$ is a strict $n-1$-groupoid, and for any
$1$-morphism $f:x\rightarrt y$ in $\pA $, the two morphisms of composition with
$f$
$$
\pA (y,z)\rightarrt \pA (x,z),\;\;\;\;
\pA (w,x)\rightarrt \pA (w,y)
$$
are equivalences of strict $n-1$-groupoids (see below);
\newline
(3)\, for all $x,y\in \Ob (\pA ) $, $\pA (x,y)$ is a strict $n-1$-groupoid, and
the truncation 
$\tau _{\leq
1}\pA $ (defined in the next proposition) is a $1$-groupoid.
\end{theorem}

\begin{proposition}
\label{gpdII}
If $\pA $ is a strict $n$-groupoid, then define $\tau _{\leq k}\pA $ to be the
strict $k$-category whose $i$-morphisms are those of $\pA $ for $i<k$ and whose
$k$-morphisms are the equivalence classes of $k$-morphisms of $\pA $ under the
equivalence relation that two are equivalent if there is a $k+1$-morphism
joining them. The fact that this is an equivalence relation is a statement about
$n-k$-groupoids. The set $\tau _{\leq 0}\pA $ will also be denoted $\pi _0\pA $.
The truncation is again a $k$-groupoid, and for $n$-groupoids $\pA $ the truncation
coincides with the operation defined in \cite{KapranovVoevodsky}.
\end{proposition}

If $\pA$ is a strict $n$-groupoid, define $\pi _0(\pA ):= \tau _{\leq 0}(\pA )$.
For $x\in \Ob (\pA )$ define $\pi _1(\pA , x):= (\tau _{\leq 1}(\pA ))(x,x)$,
which is a group since $\tau _{\leq 1}(\pA )$ is a groupoid by the previous proposition.
For $2\leq i\leq n$ define by induction $\pi _i(\pA , x):= \pi _{i-1}(\pA (x,x), 1_x)$.
The interchange property allows to show that this is an abelian group. These classical
definitions are recalled in \cite{KapranovVoevodsky}. 

\begin{definition}
\label{gpdIII}
A morphism $f:\pA \rightarrt \pB $ of strict $n$-groupoids is said to be an {\em
equivalence} if the following equivalent conditions are satisfied:
\newline
(a)\, $f$ induces an isomorphism $\pi _0\pA \rightarrt
\pi _0\pB $, and for every object $a\in \Ob (\pA )$ $f$ induces an isomorphism
$\pi _i(\pA ,a)\rightarrt^{\cong} \pi _i(\pB , f(a))$;
\newline
(b)\, $f$ induces a surjection $\pi _0\pA \rightarrt \pi _0\pB $ and, for every pair
of objects $x,y\in \Ob (\pA )$, an equivalence of $n-1$-groupoids
$\pA (x,y)\rightarrt \pB (f(x), f(y))$;
\newline
(c)\, if $u,v$ are $i$-morphisms in $\pA $ sharing the same source and target, and
if $r$ is an $i+1$-morphism in $\pB $ going from $f(u)$ to $f(v)$ then there exists
an $i+1$-morphism $t$ in $\pA $ going from $u$ to $v$ and an $i+2$-morphism in $\pB $
going from $f(t)$ to $r$ (this includes the limiting cases $i=-1$ where $u$ and
$v$ are not specified, and $i=n-1, n$ where ``$n+1$-morphisms'' mean equalities
between $n$-morphisms and ``$n+2$-morphisms'' are not specified).
\end{definition}

\begin{lemma}
\label{gpdLemmas}
If $f: \pA \rightarrt \pB $ and $g: \pB \rightarrt \pC $ are
morphisms of strict $n$-groupoids and if any two of $f$, $g$ and $gf$ are
equivalences, then so is the third.
If
$$
\pA \rightarrt ^f\pB \rightarrt^g\pC \rightarrt^h \pD 
$$
are morphisms of strict $n$-groupoids and if $hg$ and $gf$ are equivalences,
then $g$ is an equivalence.
\end{lemma}

It is clear for $n=0$, so we assume $n\geq 1$ and
proceed by induction on $n$: we assume that Theorem \ref{gpdI} and the
subsequent Proposition \ref{gpdII}, Definition \ref{gpdIII}, 
as well as Lemma \ref{gpdLemmas},
are known for  strict $n-1$-categories.

We first discuss the existence of truncation (Proposition \ref{gpdII}), for $k\geq 1$. Note that
in this case $\tau _{\leq k}\pA $ may be defined as the strict $k$-category
with the
same objects as $\pA $ and with
$$
(\tau _{\leq k}\pA )(x,y):= \tau _{\leq k-1}\pA (x,y).
$$
Thus the fact that the relation in question is an equivalence relation, is a
statement about $n-1$-categories and known by induction. Note that the
truncation operation clearly preserves any one of the three groupoid conditions
(1), (2), (3). Thus we may affirm in a strong sense that
$\tau _{\leq k}(\pA )$ is a
$k$-groupoid without knowing the equivalence of the
conditions (1)-(3).

Note also that the truncation operation for $n$-groupoids is the same as that
defined in \cite{KapranovVoevodsky} (they define truncation for general strict
$n$-categories but for $n$-categories which are not groupoids, their definition
is different from that of \cite{Tamsamani} and not all that useful).

For $0\leq k\leq k'\leq n$ we have
$$
\tau _{\leq k}(\tau _{\leq k'}(\pA )) = \tau _{\leq k}(\pA ).
$$
To see this note that the equivalence relation used to define the $k$-arrows of
$\tau _{\leq k}(\pA )$ is the same if taken in $\pA $ or in $\tau _{\leq
k+1}(\pA )$---the existence of a $k+1$-arrow going between two $k$-arrows
is equivalent to the existence of an equivalence class of $k+1$-arrows going
between the two $k$-arrows.

Finally using the above remark we obtain the existence of the truncation $\tau
_{\leq 0}(\pA )$: the relation is the same as for the truncation $\tau _{\leq
0}(\tau _{\leq 1}(\pA ))$, and $\tau _{\leq 1}(\pA )$ is a strict $1$-groupoid in
the usual sense so the arrows are invertible, which shows that the relation
used to define the $0$-arrows (i.e. objects) in $\tau _{\leq 0}(\pA )$ is
in fact an equivalence relation.

We complete our discussion of truncation
by noting that there is a natural morphism
of strict $n$-categories
$\pA \rightarrt \tau _{\leq k}(\pA )$, where the right hand side ({\em a priori} a
strict $k$-category) is considered as a strict $n$-category in the obvious way.

\bigskip

We turn next to the notion of equivalence (Definition \ref{gpdIII}), and prove that conditions
(a) and (b) are equivalent. This notion for $n$-groupoids will not enter
into the
subsequent treatment of Theorem \ref{gpdI}---what does enter is the notion of equivalence
for $n-1$-groupoids, which is known by induction---so we may assume the
equivalence of definitions (1)-(3) for our discussion of Definition \ref{gpdIII}.

Recall first of all the definition of the homotopy groups. Let $1^i_a$ denote
the $i$-fold iterated identity of an object $a$; it is an $i$-morphism, the
identity of $1^{i-1}_a$ (starting with $1^0_a=a$). Then
$$
\pi _i(\pA ,a):= {\tau _{\leq i}(\pA )}(1^{i-1}_a, 1^{i-1}_a).
$$
This definition is completed by setting $\pi _0(\pA ):= \tau _{\leq 0}(\pA )$.
These definitions are the same as in \cite{KapranovVoevodsky}. Note directly from the
definition that for $i\leq k$  the truncation morphism induces isomorphisms
$$
\pi _i(\pA ,a)\rightarrt^{\cong} \pi _i(\tau _{\leq k}(\pA ), a).
$$
Also for $i\geq 1$ we have
$$
\pi _i(\pA ,a)= \pi _{i-1}(\pA (a,a), 1_a).
$$
One shows that the $\pi _i$ are abelian for $i\geq 2$.  This is a consequence of the
Eckmann-Hilton argument discussed in the previous section. 

Suppose $f:\pA \rightarrt \pB $ is a morphism of strict $n$-groupoids satisfying
condition (b). From the immediately preceding formula and the inductive
statement for $n-1$-groupoids, we get that $f$ induces isomorphisms
on the $\pi _i$ for $i\geq 1$. On the other hand, the truncation
$\tau _{\leq 1}(f)$ satisfies condition (b) for a morphism of $1$-groupoids,
and this is readily seen to imply that $\pi _0(f)$ is an isomorphism. Thus $f$
satisfies condition (a).

Suppose on the other hand that $f:\pA \rightarrt \pB $
is a morphism of strict $n$-groupoids satisfying
condition (a). Then of course $\pi _0(f)$ is surjective. Consider two objects
$x,y\in \pA $ and look at the induced morphism
$$
f^{x,y}: \pA (x,y)\rightarrt \pB (f(x), f(y)).
$$
We claim that $f^{x,y}$ satisfies condition (a) for a morphism of
$n-1$-groupoids. For this, consider a $1$-morphism from $x$ to $y$, i.e. an
object $r\in \pA (x,y)$. By version (2) of the groupoid condition
for $\pA $, multiplication by $r$ induces an equivalence of $n-1$-groupoids
$$
m(r): \pA (x,x)\rightarrt \pA (x,y),
$$
and furthermore $m(r)(1_x)=r$. The same is true in $\pB $: multiplication by
$f(r)$ induces an equivalence
$$
m(f(r)): \pB (f(x), f(x))\rightarrt \pB (f(x), f(y)).
$$
The fact that $f$ is a morphism implies that these fit into a commutative
square
$$
\begin{diagram}
\pA (x,x)&\rightarr &\pA (x,y)\\
\downarr && \downarr \\
\pB (f(x), f(x))&\rightarr &\pB (f(x), f(y)).
\end{diagram}
$$
The equivalence condition (a) for $f$ implies that the left vertical morphism
induces isomorphisms
$$
\pi _i(\pA (x,x), 1_x)\rightarrt^{\cong}
\pi _i(\pB (f(x), f(x)), 1_{f(x)}).
$$
Therefore the right vertical morphism (i.e. $f_{x,y}$) induces isomorphisms
$$
\pi _i(\pA (x,y), r)\rightarrt^{\cong}
\pi _i(\pB (f(x), f(y)), f(r)),
$$
this for all $i\geq 1$.
We have now verified these isomorphisms for any base-object $r$. A similar
argument implies that $f^{x,y}$ induces an injection on $\pi _0$. On the other
hand, the fact that $f$ induces an isomorphism on $\pi _0$ implies that
$f^{x,y}$ induces a surjection on $\pi _0$ (note that these last two statements
are reduced to statements about $1$-groupoids by applying $\tau _{\leq 1}$ so we
don't give further details). All of these statements taken together  imply that
$f^{x,y}$ satisfies condition (a), and by the inductive statement of
the theorem for $n-1$-groupoids this implies that $f^{x,y}$ is an equivalence.
Thus $f$ satisfies condition (b).

We now remark that condition (b) is equivalent to condition (c) for a morphism
$f:\pA \rightarrt \pB $.
Indeed, the part of condition (c) for $i=-1$ is, by the definition of $\pi _0$,
identical to the condition  that $f$ induces a surjection $\pi _0(\pA )\rightarrt
\pi _0(\pB )$. And the remaining conditions for $i=0,\ldots , n+1$ are identical to
the  conditions of (c) corresponding to $j=i-1=-1,\ldots , (n-1)+1$ for all the
morphisms of $n-1$-groupoids  $\pA (x,y)\rightarrt \pB (f(x), f(y))$.
(In terms of $u$ and $v$ appearing in the condition in question, take $x$ to
be the source of the source of the source \ldots , and take $y$ to be the
target of the target of the target \ldots ).   Thus by induction on $n$
(i.e. by the equivalence $(b)\Leftrightarrow (c)$ for $n-1$-groupoids),
the conditions (c) for $f$ for $i=0,\ldots , n+1$, are equivalent to the
conditions that  $\pA (x,y)\rightarrt \pB (f(x), f(y))$
be equivalences of $n-1$-groupoids. Thus condition (c) for $f$
is equivalent to condition (b) for $f$, which completes the proof of the equivalence
of the different parts of Definition \ref{gpdIII}. 

\bigskip

We now proceed with the proof of Theorem \ref{gpdI}. Note first of
all that the implications $(1)\Rightarrow (2)$ and $(2)\Rightarrow (3)$ are
easy. We give a short discussion of $(1)\Rightarrow (3)$ anyway, and then
we prove $(3)\Rightarrow (2)$ and $(2)\Rightarrow (1)$.

Note also that
the equivalence $(1)\Leftrightarrow (2)$ is the content of Proposition 1.6 of
\cite{KapranovVoevodsky}; we give a proof here because the proof of
Proposition 1.6 was ``left to the reader'' in \cite{KapranovVoevodsky}.

\medskip

\noindent
{\bf $(1)\Rightarrow (3)$:}\,\, Suppose $\pA $ is a strict $n$-category
satisfying condition $(1)$. This condition (from \cite{KapranovVoevodsky}) is compatible with
truncation, so $\tau _{\leq 1} (\pA )$ satisfies condition $(1)$ for
$1$-categories; which in turn is equivalent to the standard condition of being a
$1$-groupoid, so we get that $\tau _{\leq 1}(\pA )$ is a $1$-groupoid. On the other
hand, the conditions $(1)$ from \cite{KapranovVoevodsky} for $i$-arrows, $1\leq i \leq n$,
include the same conditions for the $i-1$-arrows  of $\pA (x,y)$ for any
$x,y\in \Ob (\pA )$ (the reader has to verify this by looking at the definition in
\cite{KapranovVoevodsky}). Thus by the inductive statement of the present theorem for strict
$n-1$-categories, $\pA (x,y)$ is a strict $n-1$-groupoid. This shows that  $\pA $
satisfies condition $(3)$.

\medskip

\noindent
{\bf $(3)\Rightarrow (2)$:}\,\, Suppose $\pA $ is a strict $n$-category satisfying
condition $(3)$. It already satisfies the first part of condition $(2)$, by
hypothesis. Thus we have to show the second part, for example that for $f:
x\rightarrt y$ in $\Ob (\pA (x,y))$, composition with $f$ induces an
equivalence
$$
\pA (y,z)\rightarrt \pA (x,z)
$$
(the other part is dual and has the same proof which we won't repeat here).

In order to prove this, we need to make a digression about the effect of
composition with $2$-morphisms. Suppose $f,g\in \Ob  (\pA (x,y))$ and
suppose that $u$ is a $2$-morphism from $f$ to $g$---this last supposition
may be
rewritten
$$
u\in \Ob ({\pA (x,y)}(f,g)).
$$
{\em Claim:}
Suppose $z$ is another object; we claim that if composition with $f$
induces an equivalence $\pA (y,z)\rightarrt \pA (x,z)$, then composition
with $g$ also induces an equivalence $\pA (y,z)\rightarrt \pA (x,z)$.

To prove the claim, suppose that  $h,k$ are two $1$-morphisms from $y$ to $z$.
We now obtain a diagram
$$
\begin{diagram}
Hom _{\pA (y,z)}(h,k) & \rightarr & Hom _{\pA (x,z)}(hf, kf)\\
\downarr && \downarr \\
Hom _{\pA (x,z)}(hg, kg) & \rightarr & Hom _{\pA (x,z)}(hf, kg),
\end{diagram}
$$
where the top arrow is given by composition $\ast _0$ with $1_f$; the left arrow
by composition $\ast _0$ with $1_g$; the bottom arrow by composition $\ast _1$
with the $2$-morphism $h\ast _0u$; and the right morphism is given by
composition
with $k\ast _0u$. This diagram commutes (that is the ``Godement rule'' or
``interchange rule'' cf \cite{KapranovVoevodsky} p. 32). By the inductive statement of the
present theorem (version (2) of the groupoid condition) for the $n-1$-groupoid
$\pA (x,z)$, the morphisms on the bottom and on the right in the above
diagram are equivalences. The hypothesis in the claim that $f$ is an
equivalence means that the morphism along the top of the diagram is an
equivalence; thus by the first part of Lemma \ref{gpdLemmas} applied
to the $n-2$-groupoids in the diagram, we get that the morphism on the left of
the diagram is an equivalence. This provides the second half of the criterion
(b) of Definition \ref{gpdIII} for showing that the morphism of composition with $g$,
$\pA (y,z)\rightarrt \pA (x,z)$, is an equivalence of $n-1$-groupoids.

To finish the proof of the claim, we now verify the first half of criterion (b)
for the morphism of composition with $g$ (in this part we use directly the
condition (3) for $\pA $ and don't use either $f$ or $u$).  Note that $\tau _{\leq
1}(\pA )$ is a $1$-groupoid, by the condition (3) which we are assuming. Note also
that (by definition)
$$
\pi _0\pA (y,z)= {\tau _{\leq 1}\pA }(y,z) \;\;\; \mbox{and}\;\;\;
\pi _0\pA (x,z)= {\tau _{\leq 1}\pA }(x,z),
$$
and the morphism in question here is just the morphism of composition by the
image of $g$ in $\tau _{\leq 1}(\pA )$. Invertibility of this morphism in
$\tau _{\leq 1}(\pA )$ implies that the composition morphism
$$
(\tau _{\leq 1}\pA )(y,z)\rightarrt (\tau _{\leq 1}\pA )(x,z)
$$
is an isomorphism. This completes verification of the first half of criterion
(b), so we get that composition with
$g$ is an equivalence. This completes the proof of the claim.

We now return to the proof of the composition condition for (2). The fact that
$\tau _{\leq 1}(\pA )$ is a $1$-groupoid implies that given $f$ there is another
morphism $h$ from $y$ to $x$ such that the class  of $fh$ is equal to the
class of
$1_y$ in $\pi _0\pA (y,y)$, and the class of $hf$ is equal to
the class of $1_x$ in $\pi _0\pA (x,x)$. This means that there exist
$2$-morphisms $u$ from $1_y$ to $fh$, and $v$ from $1_x$ to $hf$. By the above
claim (and the fact that the compositions with $1_x$ and $1_y$ act as the
identity and in particular are equivalences), we get that composition
with $fh$ is an equivalence
$$
\{ fh\} \times \pA (y,z) \rightarrt \pA (y,z),
$$
and that composition with $hf$ is an equivalence
$$
\{ hf\} \times \pA (x,z)\rightarrt \pA (x,z).
$$
Let
$$
\psi _f: \pA (y,z)\rightarrt \pA (x,z)
$$
be the morphism of composition with $f$, and let
$$
\psi _h: \pA (x,z)\rightarrt \pA (y,z)
$$
be the morphism of composition with $h$. We have seen that $\psi _h\psi _f$ and
$\psi _f\psi _h$ are equivalences. By the second
statement of Lemma \ref{gpdLemmas} applied to $n-1$-groupoids, these imply that $\psi _f$ is an
equivalence.

The proof for composition in the other direction is the same; thus we have
obtained condition (2) for $\pA $.

\medskip

\noindent
{\bf $(2)\Rightarrow (1)$:}\,\, Look at the condition (1) by refering to
\cite{KapranovVoevodsky}: in question are the conditions $GR'_{i,k}$ and $GR''_{i,k}$ ($i<k\leq
n$) of Definition 1.1, p. 33 of \cite{KapranovVoevodsky}. By the inductive version of the
present equivalence for $n-1$-groupoids and by the part of condition (2) which
says that the $\pA (x,y)$ are $n-1$-groupoids, we obtain the conditions
$GR'_{i,k}$ and $GR''_{i,k}$ for $i\geq 1$. Thus we may now restrict our
attention to the condition $GR'_{0,k}$ and  $GR''_{0,k}$.  For a $1$-morphism
$a$ from $x$ to $y$, the conditions $GR'_{0,k}$ for all $k$ with respect to $a$,
are the same as the condition that for all $w$, the morphism of
pre-multiplication by $a$ $$
\pA (w,x)\times \{ a\} \rightarrt \pA (w,y)
$$
is an equivalence according to the version (c) of the notion of equivalence (Definition 
\ref{gpdIII}). Thus, condition $GR'_{0,k}$ follows from the
second part of condition (2) (for pre-multiplication). Similarly condition
$GR''_{0,k}$ follows from the second part of condition (2) for
post-multiplication by every $1$-morphism $a$. Thus condition (2) implies
condition (1). This completes the proof of Part (I) of the theorem.

\bigskip

For proof of the first part of Lemma \ref{gpdLemmas}, using the fact that isomorphisms
of sets satisfy the same ``three for two'' property, and using the
characterization of equivalences in terms of homotopy groups (condition (a))
we immediately get two of the three statements: that if $f$ and $g$ are
equivalences then $gf$ is an equivalence; and that if $gf$ and $g$ are
equivalences then $f$ is an equivalence. Suppose now that $gf$ and $f$ are
equivalences; we would like to show that $g$ is an equivalence. First of all
it is clear that if $x\in \Ob (\pA )$ then $g$ induces an isomorphism
$\pi _i(\pB , f(x))\cong \pi _0(\pC , gf(x))$ (resp. $\pi _0(\pB )\cong \pi _0(\pC )$).
Suppose now that $y\in \Ob (\pB )$, and choose a $1$-morphism $u$ going from $y$ to
$f(x)$ for some $x\in \Ob (\pA )$ (this is possible because $f$ is surjective on $\pi
_0$).  By condition (2) for being a groupoid, composition with
$u$ induces equivalences along the top
row of the diagram
$$
\begin{diagram}
\pB (y,y) &\rightarr &\pB (y, f(x))&\leftarr &\pB (f(x), f(x))\\
\downarr && \downarr && \downarr \\
\pC (g(y),g(y)) &\rightarr &\pC (g(y), gf(x))&\leftarr &
\pC (gf(x), gf(x)).
\end{diagram}
$$
Similarly composition with $g(u)$ induces equivalences along the bottom row.
The sub-lemma for $n-1$-groupoids applied to the sequence
$$
\pA (x,x)\rightarrt \pB (f(x), f(x))\rightarrt \pC (gf(x), gf(x))
$$
as well as the hypothesis that $f$ is an equivalence, imply that the rightmost
vertical arrow in the above diagram is an equivalence. Again applying the
sub-lemma to these $n-1$-groupoids yields that the leftmost vertical arrow is
an equivalence. In particular $g$ induces isomorphisms
$$
\pi _i(\pB ,y) = \pi _{i-1}(\pB (y,y), 1_y) \rightarrt^{\cong}
\pi _{i-1}(\pC (g(y), g(y)),1_{g(y)}) = \pi _i(\pC , g(y)).
$$
This completes the verification of condition (a) for the morphism $g$,
completing the proof of part (IV) of the theorem.

Finally we prove the second part of Lemma \ref{gpdLemmas} (from which we
now adopt the notations $\pA ,\pB ,\pC ,\pD ,f,g,h$). Note first of all that
applying $\pi _0$ gives the same situation for maps of sets, so $\pi _0(g)$ is
an isomorphism. Next, suppose $x\in \Ob (\pA )$. Then we obtain a sequence
$$
\pi _i(\pA ,x)\rightarrt \pi _i(\pB ,f(x))
\rightarrt \pi _i(\pC , gf(x))\rightarrt \pi _i(\pD , hgf(x)),
$$
such that the composition of the first pair and also of the last pair are
isomorphisms; thus $g$ induces an isomorphism
$\pi _i(\pB  , f(x))\cong \pi _i(\pC , gf(x))$. Now, by the same argument as for Part
(IV) above, (using the hypothesis that $f$ induces a surjection $\pi
_0(\pA )\rightarrt \pi _0(\pB )$) we get that for any object $y\in \Ob (\pB )$, $g$
induces an isomorphism $\pi _i(\pB  , y)\cong \pi _i(\pC , g(y))$. By Definition \ref{gpdIII} (a) we have now shown that $g$ is an equivalence.
This completes the proof of the statements in question.
\eop

Let $nStrGpd$ be the category of strict $n$-groupoids.
In Chapter \ref{nonstrict1} we shall see that for any realization functor $\Re : nStrGpd\rightarrt \Top$ preserving homotopy groups, topological spaces with nontrivial
Whitehead products cannot be weakly equivalent to any $\Re (\pA )$. This was 
Grothendieck's motivation for proposing to look for a definition of weak 
$n$-category. Our first look at the case of strict $n$-categories serves nevertheless as a
guide to the outlines of any general theory of weak $n$-categories.


\chapter{Fundamental elements of $n$-categories}
\label{ncats1}

The observation that the theory of strict $n$-groupoids 
is not enough to give a good model for homotopy $n$-types (detailed in the next Chapter \ref{nonstrict1}), 
led Grothendieck to ask for a theory of $n$-categories with {\em weakly associative composition}.
This will be the main subject of our book, in particular we use the terminology {\em $n$-category} to mean some kind of object in a possible theory with weak associativity,
or even ill-defined composition, 
or perhaps some other type of weakening (as will be briefly discussed in Chapter \ref{operadic1}). 

There are a certain number of basic elements expected of any theory of $n$-categories,
and which can be explained without refering to a full definition. It will be useful to start
by considering these. Our discussion follows Tamsamani's paper \cite{Tamsamani}, but really 
sums up the general expectations for a theory of $n$-categories which were developped over many years starting
with Benabou and continuing through the theory of strict $n$-categories and Grothendieck's manuscript. 

For this chapter, we will use the terminology ``$n$-category'' to mean
any object in a generic theory of $n$-categories. We will sometimes use the idea that our generic theory
should admit direct products and disjoint sums. 

\section{A globular theory}

We saw that a strict $n$-category has, in particular, an underlying globular set.
This basic structure should be conserved, in some form, in any weak theory. 

\noindent
(OB)---An $n$-category $\pA $ should have an {\em underlying set of objects} denoted $\ob (\pA )$. If $i=0$ then
the structure $\pA $ is identified with just this set $\ob (\pA )$, that is to say a $0$-category is just a set.

\noindent
(MOR)---If $i\geq 1$ then for any two elements $x,y\in \ob (\pA )$, there should be an {\em $n-1$-category of morphisms from $x$ to $y$}
denoted $\uMor _{\pA} (x,y)$.  From these two things, we obtain by induction a whole family of sets called the {\em sets of $i$-morphisms of $\pA $}
for $0\leq i\leq n$. 

\noindent
(PS)---With respect to direct products and disjoint sums, we should have $\ob (\pA \times \pB ) = \ob (\pA )\times \ob (\pB )$
and $\ob (\pA \sqcup \pB ) = \ob (\pA ) \sqcup \ob (\pB )$. 

The set of $i$-morphisms of $\pA $ can be defined inductively by the following procedure. 
Put 
$$
\uMor [\pA ] := \coprod _{x,y\in \ob (\pA )}\uMor _{\pA}(x,y); 
$$
this is the {\em $n-1$-category of morphisms of $\pA $}.

By induction we obtain the $n-i$-category of $i$-morphisms of $\pA $, denoted by
$$
\uMor ^i[\pA ]:= \uMor [\uMor [\cdots [\pA ]\cdots ]].
$$
This is defined whenever $0\leq i \leq n$, with $\uMor ^0[\pA ]:=\pA $ and 
$\uMor ^n[\pA ]$ being a set. 

Define 
$$
\Mor ^i[\pA ]:= \ob (\uMor ^i[\pA ]).
$$
This is a set, called the {\em set of $i$-morphisms of $\pA $}.

From the above definitions we can write
$$
\uMor ^i[\pA ] = \coprod _{x,y\in \Mor ^{i-1}[\pA ]}\uMor _{\uMor ^{i-1}[\pA ]}(x,y),
$$
and by compatibility of objects with coproducts, 
$$
\Mor ^i[\pA ] = \coprod _{x,y\in \Mor ^{i-1}[\pA ]}\ob (\uMor _{\uMor ^{i-1}[\pA ]}(x,y)).
$$
In particular, we have maps $s_i$ and $t_i$ from $\Mor ^i[\pA ]$ to $\Mor ^{i-1}[\pA ]$ 
taking an element $f\in \Mor ^i[\pA ]$ lying in the piece of the coproduct indexed by $(x,y)$, to $s_i(f):= x$ or $t_i(f):=y$ 
respectively. These maps are called {\em source} and {\em target} and if no confusion arises, the index $i$ may be dropped. 

If $u,v\in \Mor ^i[\pA ]$, let $\Mor ^{i+1}_{\pA}(u,v)$ denote the preimage of the pair $(u,v)$ by the map $(s_{i+1}, t_{i+1})$.
It is nonempty only if $s_i(u)=s_i(v)$ and $t_i(u)=t_i(v)$ and when using the notation $\Mor ^{i+1}_{\pA}(u,v)$ we generally mean
to say that these conditions are supposed to hold. Similarly, we get $n-i-1$-categories denoted $\uMor ^{i+1}_{\pA}(u,v)$. 

In this way, starting just from the principles (OB) and (MOR) together with the compatibility with sums in (PS),
we obtain from an $n$-category a collection of sets 
$$
\Mor ^0[\pA ]= \ob (\pA ); \,\, \Mor ^1[\pA ], \ldots , \Mor ^n[\pA ]
$$
together with pairs of maps
$$
s_i,t_i : \Mor ^i[\pA ]\rightarrt \Mor ^{i-1}[\pA ].
$$
They satisfy 
$$
s_is_{i+1} = s_it_{i+1},\;\; 
t_is_{i+1} = t_it_{i+1}.
$$
These elements make our theory of $n$-categories into a {\em globular theory}.

Among other things, starting from this structure we can draw pictures in a way which is usual for the theory of $n$-categories. 
These pictures explain why the theory is called ``globular''.
A $0$-morphism is just a point, and a $1$-morphism is pictured as a usual arrow
$$
{\setlength{\unitlength}{.5mm}
\begin{picture}(60,20)
\put(10,10){\circle*{2}}
\put(50,10){\circle*{2}}
\put(12,10){\line(1,0){36}}
\put(48,10){\vector(1,0){0}}
\end{picture}
} 
$$
A $2$-morphism is pictured as
$$
{\setlength{\unitlength}{.5mm}
\begin{picture}(60,40)
\put(10,20){\circle*{2}}
\put(50,20){\circle*{2}}
\qbezier(11,21)(30,40)(47,23)
\qbezier(11,19)(30,0)(47,17)
\put(49,21){\vector(1,-1){0}}
\put(49,19){\vector(1,1){0}}
\put(29,18){\ensuremath{\Downarrow}}
\end{picture}
} 
$$
whereas a $3$-morphism should be thought of as a sort of ``pillow''
which might be pictured as
$$
{\setlength{\unitlength}{.5mm}
\begin{picture}(80,80)
\put(10,20){\circle*{3}}
\put(70,60){\circle*{3}}
\qbezier(10,20)(55,5)(71,58)
\qbezier(10,20)(25,70)(69,62)

\qbezier(28,57)(36,66)(46,58)
\qbezier(22,56)(28,65)(40,56)
\put(37,56){\line(1,0){7}}
\qbezier(44,56)(46,58)(48,60)

\qbezier(35,35)(41,30)(47,28)
\qbezier(39,37)(45,32)(51,32)
\qbezier(51,28)(46,27)(41,26)
\qbezier(51,28)(52,32)(53,36)

\put(36,52){\line(0,-1){11}}
\put(38,52){\line(0,-1){11}}
\put(40,52){\line(0,-1){11}}

\qbezier(38,40)(40,41)(43,42)
\qbezier(38,40)(36,41)(33,42)

\linethickness{.3mm}
\qbezier(11,20)(50,25)(70,58)
\qbezier(10,21)(30,55)(68,60)

\qbezier(70,58)(68,57)(66,56)
\qbezier(70,58)(70,56)(70,55)

\qbezier(68,60)(67,59)(65,58)
\qbezier(68,60)(67,61)(65,62)
\end{picture}
} 
$$

\section{Identities}
\label{sec-identities}

For each $x\in\ob (\pA )$ there should be a natural element $1_x\in \Mor _{\pA}(x,x)$, called the {\em identity of $x$}.
One can envision theories in which the identity is not well-defined but exists only up to homotopy, see Kock and Joyal \cite{Kock} \cite{JoyalKock}. 
However, the theory considered here will have canonical identities. 

Following the same inductive procedure as in the previous section, we get morphisms for any $0\leq i<n$,
$$
e_i:\Mor ^i[\pA ]\rightarrt \Mor ^{i+1}[\pA ]
$$
such that $s_{i+1}e_i(u)=u$ and $t_{i+1}e_i(u)=u$. We call $e_i(u)$ the {\em identity $i+1$-morphism of the $i$-morphism $u$}.

Some authors introduce a {\em category of globules} $\Glob^n$ having objects $M_i$ for $0\leq i\leq n$, with 
generating morphisms $s_i,t_i: M_i\rightarrt M_{i-1}$ and $e_i:M_i\rightarrt M_{i+1}$ subject to the relations
$$
s_is_{i+1} = s_it_{i+1},\;\; 
t_is_{i+1} = t_it_{i+1}, \;\; s_{i+1}e_i = 1_{M_i}, \;\; t_{i+1}e_i = 1_{M_i}.
$$
An {\em $n$-globular set} is a functor $\Glob ^n\rightarrt \Sets$; with the
identities this should be called ``reflexive''. Any $n$-category $\pA $ induces an {\em underlying globular set}
constructed as above. Other authors (such as Batanin) use a category of globules
which doesn't have the identity arrows $e_i$, leading to non-reflexive globular
sets, indeed we shall use that notation
in Section \ref{sec-batanin}.

The first and basic idea for defining a theory of $n$-categories is that an $n$-category should consist of an underlying globular set (with or without identities),
plus additional structural morphisms satisfying certain properties. Whereas the Batanin-type theories \cite{Batanin} \cite{Leinster} \cite{MaltsiniotisGroBat}
are closest to this ideal, the Segal-type theories we consider in the present book will add additional structural sets to the basic globular set of $\pA $.

\section{Composition, equivalence and truncation}
\label{sec-compequivtrunc}

For objects $x,y,z\in \ob (\pA )$ there should be some kind of morphism of $n-1$-categories
\begin{equation}
\label{eq-comp}
\uMor _{\pA}(x,y)\times \uMor _{\pA}(y,z)\rightarrt \uMor _{\pA}(x,z)
\end{equation}
corresponding to {\em composition}. 
In the Segal-type theories considered in this book, the composition morphism is not well defined and
may not even exist, rather existing only in some homotopic sense. 

Nonetheless, in order best to motivate the following discussion,
assume for the moment that we know what composition means, particularly how to define $g\circ f\in \Mor^1 _{\pA}(x,z)$ for
$f\in \Mor ^1_{\pA}(x,y)$ and $g\in \Mor ^1_{\pA}(x,y)$. 

We can then inductively define a notion of {\em equivalence}.
Tamsamani calls this {\em inner equivalence} \cite{Tamsamani} to emphasize that we are speaking of arrows in our $n$-category $\pA $ which are equivalences internally in $\pA $. 
To be more precise, we will define what it means for $f\in \Mor ^1_{\pA}(x,y)$ to be an {\em inner equivalence between $x$ and $y$}. 
If such an $f$ exists, we say that {\em $x$ and $y$ are equivalent} and write $x\sim y$.

Inductively we suppose known what this means for $n-1$-categories, and in particular within the $n-1$-categories $\uMor _{\pA}(x,x)$ or $\uMor _{\pA}(y,y)$.

The definition then proceeds by saying that $f\in \Mor ^1_{\pA}(x,y)$ is an inner equivalence between $x$ and $y$, if there exists $g\in \Mor ^1_{\pA}(y,x)$
such that $g\circ f$ is equivalent to $1_x$ in $\uMor _{\pA}(x,x)$ and $f\circ g$ is equivalent to $1_y$ in $\uMor _{\pA}(y,y)$. 

This notion should be transitive in the sense that if $f$ is an equivalence from $x$ to $y$ and $g$ is an equivalence from $y$ to $z$, then 
$g\circ f$ should be an equivalence from $x$ to $z$. The relation ``$x\sim y$'' is therefore a transitive equivalence relation on the set $\ob (\pA )$.

Define the {\em truncation $\tau _{\leq 0}(\pA )$} to be the quotient set $\ob (\pA )/\sim $.

We can go further and define the $1$-categorical truncation $\tau _{\leq 1}(\pA )$, a $1$-category, as follows:
$$
\ob (\tau _{\leq 1}(\pA )):=\ob (\pA ),
$$
$$
\Mor ^1_{\tau _{\leq 1}(\pA )}(x,y):= \tau _{\leq 0}(\uMor _{\pA}(x,y)).
$$
In other words, the objects of $\tau _{\leq 1}(\pA )$ are the same as the objects of $\pA $, but the morphisms of $\tau _{\leq 1}(\pA )$
are the equivalence classes of $1$-morphisms of $\pA $, under the equivalence relation on the objects of the $n-1$-category $\uMor _{\pA}(x,y)$. 

Composition of morphisms in  $\tau _{\leq 1}(\pA )$ should be defined by composing representatives of the equivalence classes. One of the main requirements for
our theory of $n$-categories is that this composition in  $\tau _{\leq 1}(\pA )$ should be well-defined, independent of the choice of representatives and
indeed independent of the choice of notion of composition morphism introduced at the start of this section. 

Denote also by $\sim$ the equivalence relation obtained in the same way on the objects of the $n-i$-categories $\uMor ^i[\pA ]$. Noting that it is compatible
with the source and target maps, we get an equivalence relation $\sim$ on $\Mor ^i_{\pA}(u,v)$ for any $i-1$-morphisms $u$ and $v$. 

The above discussion presupposed the existence of some kind of composition operation, but in the Segal-style theory we consider in this book, such a
composition morphism is not canonically defined. Thus, we restart the discussion without assuming existence of a composition morphism of $n-1$-categories. 
The first fundamental structure to be considered is thus:

\noindent
(EQUIV)---on each set $\Mor ^i[\pA ]$ we have an equivalence relation $\sim$ compatible with the source and target maps, giving the {\em set of $i$-morphisms 
up to equivalence} $\Mor ^i[\pA ]/\sim $. For $i=n$ this equivalence relation should be trivial. The induced relation on $\Mor ^i_{\pA}(u,v)$ is also denoted $\sim$. 

We can then consider the structure of composition which is well-defined up to equivalence, in other words it is given by a map on quotient sets. 

\noindent
(COMP)---for any $0<i\leq n$ and any three $i-1$-morphisms $u,v,w$ sharing the same sources and  the same targets, we have a well-defined composition map
$$
\left( \Mor^i_{\pA}(u,v)/\sim \right)\times  \left( \Mor^i_{\pA}(v,w)/\sim \right) \rightarrt \Mor^i_{\pA}(u,w)/\sim
$$
which is associative and has the classes of identity morphisms as left and right units. 

These two structures are compatible in the sense that composition is defined after passing to the quotient by $\sim$. As a matter of simplifying
notation, given $f\in \Mor^i_{\pA}(u,v)$ and $g\in \Mor^i_{\pA}(v,w)$ then denote by $g\circ f$ 
any representative in $\Mor^i_{\pA}(u,w)$ for the composition of the class of $g$ with the class of $f$. This is well-defined up to equivalence and
by construction independent, up to equivalence, of the choices of representatives $f$ and $g$ for their equivalence classes. 

Equivalence and composition also satisfy the following further compatibility condition, expressing the notion of equivalence in terms which
closely resemble the classical definition of equivalence of categories. 

\noindent
(EQC)---for any $0\leq i < n$ and $u,v\in \Mor ^i_{\pA}$ sharing the same source and target (i.e. $s_i(u)=s_i(v)$ and $t_i(u)=t_i(v)$ in case $i>0$),
then $u \sim v$ if and only if there exist $f\in \Mor ^{i+1}_{\pA}(u,v)$ and $g\in \Mor ^{i+1}_{\pA}(v,u)$ such that $g\circ f \sim 1_u$ and $f\circ g\sim 1_v$. 

With these structures, we can define the $1$-categories $\tau _{\leq 1}\uMor  ^i_{\pA}(u,v)$, having objects the elements of $\Mor ^i_{\pA}(u,v)$ and
as morphisms between $w,z\in \Mor ^i_{\pA}(u,v)$ the equivalence classes $\Mor ^{i+1}_{\pA}(w,z)/\sim$. The composition of (COMP) gives this a structure of $1$-category,
and $w\sim z$ if and only if $w$ and $z$ are isomorphic objects of $\tau _{\leq 1}\uMor  ^i_{\pA}(u,v)$. At the bottom level we obtain a $1$-category denoted
$\tau _{\leq 1}(\pA )$ and called the {\em $1$-truncation of $\pA $},
whose set of objects is $\ob (\pA )$ and whose set of morphisms is $\Mor ^1[\pA ]/\sim$. These constructions are compatible with the 
induction in the sense that $\tau _{\leq 1}\uMor  ^i_{\pA }(u,v)$ is indeed the $1$-truncation of the $n-i$-category $\uMor  ^i_{\pA}(u,v)$. 

Suppose $x,y\in \Mor ^{i-1}[\pA ]$ and $u,v\in \Mor ^i_{\pA}(x,y)$. An element $f\in \Mor ^{i+1}_{\pA}(u,v)$ 
is said to be an {\em internal equivalence} between $u$ and $v$, if its class is an isomorphism in  $\tau _{\leq 1}\uMor ^{i}_{\pA}(x,y)$.
This is equivalent to requiring the existence of $g\in \Mor ^{i+1}_{\pA}(v,u)$ such that $g\circ f \sim 1_u$ and $f\circ g\sim 1_v$.

\section{Enriched categories}
\label{sec-enrichedcats}

The natural first approach to the notion of $n$-category is to ask for $n-1$-categories of morphisms $\uMor _{\pA}(x,y)$, with composition operations
\eqref{eq-comp} which are strictly associative and have the $1_x$ as strict left and right units. This gives a structure of {\em category enriched over 
$n-1$-categories}.  In an intuitive sense the reader should think of an $n$-category in this way. However, if the definition is applied inductively
over $n$, that is to say that the $n-1$-categories $\uMor _{\pA}(x,y)$ are themselves enriched over $n-2$-categories and so forth, one gets to the notion of
strict $n$-category considered in the previous chapter. But, as 
we shall discuss in the next Chapter \ref{nonstrict1} below, the strict $n$-categories
are not sufficient to capture all of the homotopical behavior we want for $n\geq 3$. 

Paoli has shown \cite{PaoliAdvances} that homotopy $n$-types can be modelled by 
{\em semistrict} $n$-groupoids, in other words $n$-categories which are strictly enriched
over weak $n-1$-categories. Bergner showed a corresponding strictification theorem for Segal categories, and the analogous strictification from $A_{\infty}$-categories
to dg categories has been known to the experts for some time. Lurie's technique \cite{LurieTopos} for constructing the model category structure we consider here, gives
additionally the strictification theorem generalizing Bergner's result. So, as we shall discuss briefly in Section \ref{strictMcats}, 
the objects of our Segal-type theory of $n$-categories
can always be assumed equivalent to semistrict ones, that is to categories strictly enriched over the model category for $n-1$-precategories. This doesn't mean
that we can go inductively towards strict $n$-categories because the strictification operation is not compatible with direct product, so if applied to the
enriched morphism objects it destroys the strict enrichment structure. As Paoli notes in \cite{PaoliAdvances}, semistrictification at one level is as far as we can go.

\section{The $(n+1)$-category of $n$-categories}

One of the main goals of a theory of $n$-categories is to provide a structure of $n+1$-category on the collection of all $n$-categories.
Of course, some discussion of universes is needed here: the collection of all $n$-categories in a universe $\univa$ should form an $n+1$-category
in a bigger containing universe $\univb \supset \univa$. This precision will be dropped from most of our discussions below. 

Recall that the notion of $2$-category was originally introduced because of the familiar observation that
``the set of all categories is actually a $2$-category'', a $2$-category
to be denoted $1CAT$. Its objects are the $1$-categories (in the smaller universe); the $1$-morphisms of $1CAT$ are the functors between categories,
and the $2$-morphisms between functors are the natural transformations. 

In general, we hope and expect to obtain an $n+1$-category denoted $nCAT$, whose objects are the $n$-categories. The $1$-morphisms of $nCAT$ are the 
``true'' functors between $n$-categories, and we obtain all of the $\Mor ^i[nCAT]$ for $0\leq i\leq n+1$ which are higher analogues of natural transformations
and so on. 

The notion of internal equivalence within $nCAT$ itself, yields the notion of {\em external equivalence} between $n$-categories: a functor $f:\pA \rightarrt \pB $
of $n$-categories, by which we mean in the most general sense an element of $\Mor ^1_{nCAT}(\pA ,\pB )$, is said to be an {\em external equivalence}
if it is an internal equivalence considered as a $1$-morphism in $nCAT$. 

In practice, a theory of $n$-categories will usually involve defining some kind of mathematically structured set or collection of sets,
which naturally generates a {\em usual $1$-category} of $n$-categories, which we can denote by $nCat$. We expect then that $\ob (nCat ) = \ob (nCAT)$
but that there is a natural projection $\Mor ^1[nCat]\rightarrt \Mor ^1[nCAT]$ compatible with composition, indeed it should come from a 
morphism of $n+1$-categories $nCat\rightarrt nCAT$ (which is to say, a $1$-morphism in the $n+2$-category $(n+1)CAT$!).

However, the notion of external equivalence in $nCat$ will not generally speaking have the same characterization as in $nCAT$: if
$f\in \Mor ^1[nCat]$ projects to an equivalence in $nCAT$ it means that there is $g\in \Mor ^1[nCAT]$ such that $fg$ and $gf$ are equivalent
to identities; however the essential inverse $g$ will not necessarily come from a morphism in $nCat$. 
For precisely this reason, one of the main tasks needed to get a theory of $n$-categories off the ground, is to give a different definition of
when a usual morphism $f\in \Mor ^1_{nCat}(\pA ,\pB )$ is an external equivalence. Of course it is to be expected and---one hopes---later proven that
this standalone definition of external equivalence, should become equivalent to the above definition once we have $nCAT$ in hand. 

A morphism of $n$-categories $f:\pA \rightarrt \pB $ should induce a morphism of sets $\ob (\pA )\rightarrt \ob(\pB )$ (usually denoted just by $x\mapsto f(x)$), and
for any $x,y\in \ob (\pA )$ it should induce a morphism of $n-1$-categories $\uMor _{\pA}(x,y)\rightarrt \uMor _{\pB}(f(x),f(y))$. Just as was the case for composition,
the morphism part of $f$ needn't necessarily be very well defined, but it should be well-defined up to an appropriate kind of equivalence.  

It is now possible to state, by induction on $n$, the second or ``standalone'' definition of external equivalence. 
A morphism $f:\pA \rightarrt \pB $ is said to be {\em fully faithful}
if for every $x,y\in \ob (\pA )$ the morphism  $\uMor _{\pA}(x,y)\rightarrt \uMor _{\pB}(f(x),f(y))$ is an external equivalence 
between $n-1$-categories. And $f$ is said to be {\em essentially surjective} if it induces a surjection $\ob (\pA ) \twoheadrightarrow \ob (\pB )/\sim$.
Then $f$ is said to be an {\em external equivalence} if it is fully faithful and essentially surjective. 

We can state the required compatibility between the two notions:

\noindent
(EXEQ)---a morphism $f:\pA \rightarrt \pB $, an element of $\Mor ^1_{nCAT}(\pA ,\pB )$, is an external equivalence (fully faithful and essentially surjective),
if and only if it is an inner equivalence in $nCAT$ (i.e. has an essential inverse $g$ such that $fg$ and $gf$ are equivalent to the identities). 

Note that the fully faithful condition implies (by an inductive consideration and comparison
with the truncation operation) that if $f:\pA \rightarrt \pB $ is an
equivalence in either of the two equivalent senses, and if $a,b$ are $i$-morphisms
of $\pA $ with the same source and target, then the set of inner equivalence classes of
$i+1$-morphisms from $a$ to $b$ in $\pA $, is isomorphic via $f$ to the
set of inner equivalence classes of $i+1$-morphisms from $f(a)$ to $f(b)$ in $\pB $. 

When developing a theory of $n$-categories, we therefore expect to be in the following situation: having first obtained a $1$-category of $n$-categorical
structures denoted $nCat$, we obtain a notion of when a morphism $f$ in this category, or first kind of functor, is an external equivalence. On the other hand,
in the full $n+1$-category $f$ should have an essential inverse $g$. So, one of the main steps towards construction of $nCAT$ is to formally invert the
external equivalences. This is a typical localization problem. Furthermore $nCAT$ should be closed under limits and colimits, so it is very natural
to use Quillen's theory of model categories, and all of the localization machinery that is now known to go along with it, as our main technical tool
for going towards the construction of $nCAT$. 

To finish this section, note one of the interesting and important features of $nCAT$: it is, in a certain sense, {\em enriched over itself}. In other words,
for two objects $\pA ,\pB \in \ob (nCAT)$, we get an $n$-category of morphisms $\uMor _{nCAT}(\pA ,\pB )$ which, since it is an $n$-category (and furthermore in the
same universe level as $\pA $ and $\pB $), is itself an object of $nCAT$:
$$
\uMor _{nCAT}(\pA ,\pB )\in \ob (nCAT).
$$
This is the motivation for using the theory of cartesian model categories with internal $\uHom$ as a preliminary model for $nCAT$. 

\section{Poincar\'e $n$-groupoids}

An $n$-category is said to be an {\em $n$-groupoid} if all $i$-morphisms are inner equivalences. More generally, we say that $\pA $ is {\em $k$-groupic}
if all $i$-morphisms are inner equivalences for $i>k$. 
Lurie introduces the notation {\em $(n,k)$-category} for a $k$-groupic $n$-category, mostly used in the limiting case $n=\infty$. 

Fundamental to Grothendieck's vision in ``Pursuing stacks'' was the {\em Poincar\'e $n$-groupoid of a space}. If $X$ is a topological space,
this is to be an $n$-groupoid denoted by $\Pi _n(X)$, with the following properties:
$$
\ob (\Pi _n(X))  = X,
$$
for $0\leq i<n$
$$
\Mor ^i[\Pi _n(X)]= C^0_{\rm glob}([0,1]^i,X)
$$
where the right hand side is the subset of maps of the $i$-cube into $X$ satisfying certain {\em globularity conditions} (explained below),
and at $i=n$ we have 
$$
\Mor ^n[\Pi _n(X)]= C^0_{\rm glob}([0,1]^n,X)/\sim
$$
where $\sim$ is an equivalence relation similar to the one considered above (and indeed, it is the same in the context of $\Pi _k(X)$ for $k>n$).

The globularity condition is automatic for $i=1$, thus the $1$-morphisms in $\Pi _n(X)$ are continuous paths $p:[0,1]\rightarrt X$ with source 
$s_1(p):=p(0)$ and target $t_1(p):= p(1)$. In the limiting case $n=1$, the $1$-morphisms in $\Pi _1(X)$ are homotopy classes of paths
with homotopies fixing the source and target, and $\Pi _1(X)$ is just the classical Poincar\'e groupoid of the space $X$. 

The {\em globularity condition} is most easily understood in the case $i=2$: a $2$-morphism in $\Pi _n(X)$ should be a homotopy between 
paths, that is to say it should be a map $\psi :[0,1]^2\rightarrt X$ such that $\psi (0,t)$ and $\psi (1,t)$ are independent of $t$. 

For $2\leq i\leq n$, the globularity condition on a map $\psi : [0,1]^i\rightarrt X$ says that for any $0\leq k <i$ and any $z_1,\ldots , z_{k-1}\in [0,1]$,
the functions
$$
(z_{k+1},\ldots , z_i)\mapsto \psi (z_1,\ldots , z_{k-1},0,z_{k+1},\ldots , z_i) 
$$
and
$$
(z_{k+1},\ldots , z_i)\mapsto \psi (z_1,\ldots , z_{k-1},1,z_{k+1},\ldots , z_i)
$$
are constant in $(z_{k+1},\ldots , z_n)$. 
The {\em source} and {\em target} of $\psi$ are defined by 
$$
s_i\psi : (z_1,\ldots , z_{i-1})\mapsto \psi (z_1,\ldots , z_{i-1},0),
$$
$$
t_i\psi : (z_1,\ldots , z_{i-1})\mapsto \psi (z_1,\ldots , z_{i-1},1).
$$

Grothendieck's fundamental prediction was that this globular set should have a natural structure of $n$-groupoid denoted $\Pi _n(X)$; that there should be a {\em realization construction}
taking an $n$-groupoid $G$ to a topological space $|G|$; and that these two constructions should set up an equivalence of homotopy theories between
$n$-truncated spaces (i.e. spaces with $\pi _i(X)=0$ for $i>n$) and $n$-groupoids. This would generalize the classical correspondence between $1$-groupoids
and their classifying spaces which are disjoint unions of Eilenberg-MacLane $K(\pi , 1)$-spaces. 

\section{Interiors}

A useful notion which should exist in a theory of $n$-categories is the notion of {\em $k$-groupic interior} denoted $\Int ^k(\pA )$. This is the largest sub-$n$-category of
$\pA $ which is $k$-groupic, and we should have 
$$
\Mor ^i[\Int ^k(\pA )]= \Mor ^i[\pA ], \;\;\; i\leq k
$$
whereas for $i>k$ the $i$-morphisms of the interior $\Mor ^i[\Int ^k(\pA )]$ should consist only of those $i$-morphisms of $\pA $ which are equivalences,
i.e. invertible up to equivalence. Specifying exactly the structure of $\Int ^k(\pA )$ will depend on the particular theory of $n$-categories, but in
any case it should be a $k$-groupic $n$-category i.e. an $(n,k)$-category. 

The usual case is for $k=1$, and sometimes the index $1$ will then be forgotten. Thus, if $\pA $ is an $n$-category then we get a $1$-groupic $n$-category
$\Int (\pA )$. As will be discussed below in the next section and more extensively in  Chapter \ref{simplicial1}, the notion of $1$-groupic $n$-category
is well modeled by the notions of simplicial category, Segal category, Rezk complete Segal space, or quasicategory. Simplicial categories typically
arise as Dwyer-Kan localizations $L(\;\; )$ of model categories, and one feature of the model category $\precat ^n(\Sets )$ constructed in this book is that 
$$
\Int (nCAT) = L(\precat ^n(\Sets )).
$$
This adds to the motivation for why it is interesting and important to use  model categories as a substrate for the theory of $n$-categories: we get
a calculatory model for an important piece of $nCAT$ namely its interior. 

\section{The case $n=\infty$}

Constructing a theory of $\infty$-categories in general, represents a new level of difficulty and to do this in detail would go beyond the scope of the present
book. We include here a few comments about this problem, largely following Cheng's observation \cite{ChengInfty}. The example considered in her paper shows that
in an algebraic sense, any $\infty$-category whose $i$-morphisms have duals at all levels, looks like an $\infty$-groupoid. However, it is clear that we don't want
to identify such $\infty$-categories with duals, and $\infty$-groupoids. Indeed, they occupy almost dual positions in the general theory as predicted by Baez and Dolan.
From this paradox one can conclude that the notion of ``equivalence'' in an $\infty$-category is not merely an algebraic one. One way around this problem would be
to include the notion of equivalence in the initial structure of an $\infty$-category: it would be a globular set with additional structure similar to the
structure used for the case of $n$-categories; but also with the information of a subset of the set of $i$-morphisms which are to be designated as ``equivalences''.
These subsets would be required to satisfy a compatibility condition similar to (EXEQ) above. Then, in Cheng's example \cite{ChengInfty} there would be two
compatible choices: either to designate everybody as an equivalence, in which case we get an $\infty$-groupoid; or to designate only some (or potentially none other
than the identities) as equivalences, yielding an ``$\infty$-category with duals'' more like what Baez and Dolan are looking for. 

In view of these problems, it is tempting to take a shortcut towards consideration of certain types of $\infty$-categories. The shortcut is motivated by
Grothendieck's Poincar\'e $n$-groupoid correspondence, which he says should also extend to an equivalence between $\infty$-groupoids and the homotopy theory
of all CW-complexes. Turning this on its head, we can use that idea to {\em define} the notion of $\infty$-groupoid as simply being a homotopy type of a space. 
The iterative enrichment procedure yields the notion of $n$-categories when started from $0$-categories being sets. If instead we start with $\infty$-groupoids
being spaces, then iterating gives a notion of $(\infty , n)$-categories which are $n$-groupic $\infty$-categories, i.e. ones whose $i$-morphisms are supposed and declared
to be invertible for all $i>n$. 

At the first stage, an $(\infty , 1)$-category is therefore a category enriched over spaces. In Part II we will consider in detail many of the different current
approaches to this theory. At the $n$-th stage, in the Segal-type theory pursued here, we obtain the notion of {\em Segal $n$-category}. The possibility of
doing an iterative definition in various different cases, motivates our presentation here of a general iterative construction of the theory of categories
weakly enriched over a cartesian model category. The cartesian condition corresponds to the fact that the composition morphism goes out from a product,
so the model structure should have a good compatibility property with respect to direct products.  
When the iteration starts from the model category of sets, we get the theory of $n$-categories; and starting from
the model category of spaces leads to the theory of Segal $n$-categories which is one approach to $(\infty , n)$-categories. 

Looking only at $(\infty , n)$-categories rather than all $\infty$-categories is compatible with the notion of $nCAT$, in the sense that $(\infty , n)CAT$, 
the collection of all $(\infty , n)$-categories, is expected to have a natural structure of $(\infty , n+1)$-category. In the theory presented here,
this will be achieved by constructing a cartesian model category for Segal $n$-categories. The cartesian property thus shows up at the output side of the iteration,
and at the input side because we need to handle products in order to talk about weak composition morphisms. Thus, one of our main goals is to 
construct an iteration step starting with a cartesian model category $\mM$ and yielding a cartesian model category $\precat (\mM )$ representing the homotopy
theory of weakly $\mM$-enriched categories.


\chapter{The need for weak composition}
\label{nonstrict1}

In this chapter, we take a break from the general theory to consider the
problem of realization of homotopy $3$-types. It is here that the phenomenon
of ``weak composition'' shows up first, in that strict $3$-groupoids are not sufficient
to model all $3$-truncated homotopy types. The notations here refer and continue those
of Chapter \ref{strict1}.

The classical Eckmann-Hilton argument, originally used to show that the homotopy groups $\pi _i$ are abelian for $i\geq 2$,
applies in the context of strict $n$-categories to give a vanishing of certain homotopy operations. Indeed, not only the $\pi _i$ but
the $i$-th loop spaces are abelian objects, and this forces the Whitehead products to vanish. This observation,
which I learned from G. Maltsiniotis and A. Brugui\`eres, had been used by many people
to argue that strict $n$-categories do not contain sufficient information to model homotopy $n$-types, as soon as $n\geq 3$.
See for example Brown \cite{RBrown1}, with Gilbert
\cite{BrownGilbert} and with Higgins \cite{BrownHiggins} \cite{BrownHiggins2}; 
Grothendieck's discussion of this in various places in
\cite{Grothendieck}, and the paper of Berger \cite{BergerEckmannHilton}.

In R.
Brown's terminology, strict $n$-groupoids correspond to {\em crossed complexes}. While
a nontrivial  action of $\pi _1$ on the $\pi _i$ can occur in a crossed complex,
the higher Whitehead operations such as $\pi _2\otimes \pi _2\rightarrt
\pi _3$
must vanish.  The Eckmann-Hilton argument for strict $n$-categories
is also known as the ``interchange rule'' or
``Godement relation''. This effect occurs when
one takes two $2$-morphisms $a$ and $b$ both with source and target a
$1$-identity  $1_x$. There are various ways of composing $a$ and $b$ in this
situation, and comparison of these compositions leads to the conclusion that all
of the compositions are commutative. In a weak $n$-category, this commutativity
would only hold up to higher homotopy, which leads to the notion of
``braiding''; and in fact it is exactly the braiding which leads to the
Whitehead operation. However, in a strict $n$-groupoid, the commutativity is
strict and applies to all higher arrows, so the Whitehead operation is trivial.

The same may be said in the setting of $3$-categories not
necessarily groupoids: there are some examples (which G. Maltsiniotis
pointed out to me) in Gordon-Power-Street \cite{GordonPowerStreet} 
of weak $3$-categories not equivalent to strict ones.  This
in turn is related to the difference between braided monoidal categories and
symmetric monoidal categories, see for example the nice discussion in Baez-Dolan
\cite{BaezDolan}.

This chapter gives a variant on these observations; it is a modified version 
of the preprint \cite{hty3types}. 
We will show that one cannot obtain all homotopy $3$-types by any reasonable
realization functor
from strict $3$-groupoids (i.e. groupoids in the sense of
\cite{KapranovVoevodsky}) to spaces.  More precisely we show that one does not obtain
the $3$-type of $S^2$. 
This constitutes a small generalization of Berger's theorem \cite{BergerEckmannHilton},
which concerned the standard realization functor.  We define the notion of 
possible ``reasonable realization functor'' in Definition
\ref{realizationdef} to be any functor $\Re$ from the category of strict
$n$-groupoids to $\Top$, provided with a natural transformation $r$ from the set
of objects of $\pG $ to the points of $\Re (\pG )$, and natural isomorphisms $\pi
_0(\pG )\cong \pi _0(\Re (\pG ))$ and $\pi _i(\pG ,x) \cong \pi _i(\Re (\pG ), r(x))$.
This axiom is fundamental to the question of whether one can realize homotopy
types by strict $n$-groupoids, because one wants to read off the homotopy
groups of the space from the strict $n$-groupoid. The standard
realization functors satisfy other properties beyond this minimal one.

In order to apply Definition
\ref{realizationdef}, the interchange argument is written in a particular way. 
We get a picture of strict
$3$-groupoids having only one object and one $1$-morphism, as being equivalent
to abelian monoidal objects $(\pG ,+)$ in the category of groupoids, such that
$(\pi _0(\pG ),+)$ is a group. In the case in question, this group will be $\pi
_2(S^2)=\zz$. Then comes the main part of the argument. We show that, up to
inverting a few equivalences, such an object has a morphism giving a splitting
of the Postnikov tower (Proposition \ref{diagramme}). It follows that for any
realization functor respecting homotopy groups, the Postnikov tower of the
realization (which has two stages corresponding to $\pi _2$ and $\pi _3$)
splits. This implies that the $3$-type of $S^2$ cannot occur as a realization,
Theorem \ref{noS2}.

\section{Realization functors}

Recall that $nStrGpd$ is the category of strict $n$-groupoids as defined in Chapter \ref{strict1}. Let $\Top$ be the category of topological spaces. The following
definition encodes the minimum of what one would expect for a reasonable
realization functor from strict $n$-groupoids to spaces.

\begin{definition}
\label{realizationdef}
A {\em realization functor for strict $n$-groupoids} is a functor
$$
\Re : nStrGpd \rightarrt \Top
$$
together with the following natural transformations:
$$
r:\Ob  (\pA ) \rightarrt \Re (\pA );
$$
$$
\zeta _i(\pA ,x): \pi _i (\pA , x) \rightarrt \pi _i (\Re (\pA ), r(x)),
$$
the latter including
$\zeta _0(\pA ): \pi _0(\pA )\rightarrt \pi _0(\Re (\pA ))$; such that
the $\zeta _i(\pA ,x)$ and $\zeta _0(\pA )$ are isomorphisms for $0\leq i \leq n$,
such that $\zeta _0$ takes the isomorphism class of $x$ to the connected component of $r(x)$,
and such that the $\pi _i(\Re (\pA ), y)$ vanish for $i>n$.
\end{definition}

\begin{theorem}
\label{realization}
There exists a realization functor $\Re$ for strict  $n$-groupoids.
\end{theorem}

Kapranov and Voevodsky \cite{KapranovVoevodsky} construct such a functor. Their construction
proceeds by first defining a notion of ``diagrammatic set''; they define a
realization functor from $n$-groupoids to diagrammatic sets (denoted $Nerv$),
and then define the topological realization of a diagrammatic set (denoted $|
\cdot |$). The composition of these two constructions gives a realization
functor  $$
\pG  \mapsto \Re _{KV}(\pG ):= | Nerv(\pG )|
$$
from
strict $n$-groupoids to spaces. Note that this functor $\Re_{KV}$ satisfies the
axioms of \ref{realizationdef} as a consequence of Propositions 2.7 and 3.5 of
\cite{KapranovVoevodsky}.

One obtains a different construction by considering strict $n$-groupoids as weak
$n$-groupoids in the sense of \cite{Tamsamani} (multisimplicial sets) and then
taking the realization of \cite{Tamsamani}.  This construction is actually
due to the Australian school many years beforehand---see \cite{BergerEckmannHilton}---and we
call it the {\em standard realization} $\Re _{\rm std}$. The properties of
\ref{realizationdef} can be extracted from \cite{Tamsamani} (although again
they are classical results).

We don't claim here that any two realization functors must be the same. 
This is why we shall work, in what follows, with an arbitrary
realization functor satisfying the axioms of \ref{realizationdef}.

\begin{proposition}
\label{realizationconsequences}
If $\pC \rightarrt \pC '$ is a
morphism of strict $n$-groupoids inducing isomorphisms on the $\pi _i$ then $\Re
(\pC )\rightarrt \Re (\pC ')$ is a weak homotopy equivalence.
Conversely if $f:\pC \rightarrt \pC '$ is a morphism of strict $n$-groupoids which
induces a weak equivalence of realizations then $f$ was an equivalence.
\end{proposition}
\begin{proof}
Apply version (a) of the equivalent conditions in Definition \ref{gpdIII},
together with the property of Definition \ref{realizationdef}. 
\end{proof}

\section{$n$-groupoids with one object}

Let $\pC $ be a strict
$n$-category with only one object $x$. Then $\pC $ is an $n$-groupoid if and only
if $\pC (x,x)$ is an $n-1$-groupoid and
$\pi _0\pC (x,x)$ (which has a structure of monoid) is a group.
This is version (3) of the definition of groupoid in Theorem \ref{gpdI}. Iterating
this remark one more time we get the following statement.

\begin{lemma}
\label{scholiumgpd}
The construction of \ref{scholium} establishes an equivalence of categories
between the strict $n$-groupoids having only one object and only one
$1$-morphism, and the abelian monoid-objects $\pG $ in $(n-2)StrGpd$ such that the
monoid $\pi _0(\pG )$ is a group.
\end{lemma}
\begin{proof}
Lemma \ref{scholium} gives an equivalence between the
categories of abelian monoid-objects in $(n-2)StrCat$, and the strict
$n$-categories having only one object and one $1$-morphism. The groupoid
condition for the $n$-category is equivalent to saying that $\pG$ is
a groupoid, and that $\pi _0(\pG )$ is a group.
\end{proof}

\begin{corollary}
\label{iterategpd}
Suppose $\pC $ is a strict $n$-category  having only one object and only one
$1$-morphism, and let $\pB $ be the strict $n+1$-category of \ref{iterate}
with one object $b$ and $\pB (b,b)= \pC $. Then $\pB $ is a strict $n+1$-groupoid if
and only if $\pC $ is a strict $n$-groupoid.
\end{corollary}
{\em Proof:}
Keep the notation $\pU$ of the proof of \ref{iterate}.
If $\pC $ is a groupoid this means that $\pG $ satisfies the condition that $\pi
_0(\pG )$ be a group, which in turn implies that $\pU $ is a groupoid. Note that $\pi
_0(\pU )=\ast$ is automatically a group; so applying the observation
\ref{scholiumgpd} once again, we get that $\pB $ is a groupoid. In the other
direction, if $\pB $ is a groupoid then $\pC =\pB (b,b)$ is a groupoid by versions
(2) and (3) of the definition of groupoid.
\eop

\section{The case of the standard realization}

Before getting to our main result which concerns an arbitrary realization
functor satisfying \ref{realizationdef}, we take note of an easier argument
which shows that the standard realization functor cannot give rise to arbitrary
homotopy types.

\begin{definition}
\label{compatiblelooping}
A collection of realization functors $\Re ^n$ for $n$-groupoids ($0\leq n <
\infty$) satisfying \ref{realizationdef} is said to be {\em compatible with
looping} if there exist transformations natural in an $n$-groupoid $\pA $ and an
object $x\in \Ob (\pA )$,
$$
\varphi (\pA , x): \Re ^{n-1}(\pA (x,x))\rightarrt \Omega ^{r(x)}\Re ^n(\pA )
$$
(where $\Omega ^{r(x)}$ means the space of loops based at $r(x)$), such that
for $i\geq 1$ the following diagram commutes:
$$
\begin{diagram}
\pi _i(\pA , x) & = \pi _{i-1}(\pA (x,x), 1_x) \rightarr &
\pi _{i-1}(\Re ^{n-1}(\pA (x,x)), r(1_x))\\
\downarr &&\downarr \\
\pi _i(\Re ^n(\pA ), r(x)) & \leftarr & \pi _{i-1}(  \Omega ^{r(x)}\Re ^n(\pA ),
cst(r(x)))
\end{diagram}
$$
where the top arrow is $\zeta _{i-1}(\pA (x,x), 1_x)$, the left arrow is
$\zeta _{i}(\pA ,x)$, the right arrow is induced by $\varphi (\pA , x)$, and the
bottom arrow is the canonical arrow from topology. (When $i=1$, suppress the
basepoints in the $\pi _{i-1}$ in the diagram.)
\end{definition}

{\em Remark:} The arrows on the top, the bottom and the left are isomorphisms
in the above diagram, so the arrow on the right is an isomorphism and we obtain
as a corollary of the definition that the $\varphi (\pA ,x)$ are actually weak
equivalences.

{\em Remark:} The collection of standard realizations $\Re ^n_{\rm std}$ for
$n$-groupoids, is compatible with looping. We leave this as an exercise.

Recall the statements of \ref{iterate} and \ref{iterategpd}:
if $\pA $ is a strict  $n$-category with only one object $x$ and only one
$1$-morphism $1_x$, then there exists a strict $n+1$-category $\pB $ with
one object $y$, and with $\pB (y,y)=\pA $; and  $\pA $ is a strict  $n$-groupoid if
and only if $\pB $ is a strict $n+1$-groupoid.

\begin{corollary}
\label{forstandard}
Suppose $\{ \Re ^n\}$ is a collection of realization  functors
\ref{realizationdef} compatible with looping \ref{compatiblelooping}.
Then if $\pA $ is a $1$-connected strict $n$-groupoid (i.e. $\pi _0(\pA )=\ast$ and
$\pi _1(\pA ,x)=\{ 1\}$), the space $\Re ^n(\pA )$ is weak-equivalent to a loop
space.
\end{corollary}
{\em Proof:}
Let $\pA '\subset \pA $ be the sub-$n$-category having one object $x$ and one
$1$-morphism $1_x$. For $i\geq 2$ the inclusion induces isomorphisms
$$
\pi _i(\pA ', x) \cong \pi _i(\pA ,x),
$$
and in view of the $1$-connectedness of $\pA $ this means (according to the
Definition \ref{gpdIII} (a)) that the morphism $\pA '\rightarrt \pA $ is an
equivalence. It follows (by definition \ref{realizationdef}) that $\Re
^n(\pA ')\rightarrt \Re ^n(\pA )$ is a weak equivalence.  Now $\pA '$ satisfies the
hypothesis of \ref{iterate}, \ref{iterategpd} as recalled above, so there is an
$n+1$-groupoid $\pB $ having one object $y$ such that $\pA '=\pB (y,y)$. By the
definition of ``compatible with looping'' and the subsequent remark that the
morphism $\varphi (\pB ,y)$ is a weak equivalence, we get that $\varphi (\pB ,y)$
induces a weak equivalence
$$
\Re ^n(\pA ') \rightarrt \Omega ^{r(y)}\Re ^{n+1}(\pB ).
$$
Thus $\Re ^n(\pA )$ is weak-equivalent to the loop-space of $\Re ^{n+1}(\pB )$.
\eop

The following corollary is due to C. Berger
\cite{BergerEckmannHilton} (although the same statement appears without proof in Grothendieck
\cite{Grothendieck}). See also R. Brown and coauthors \cite{RBrown1}
\cite{BrownGilbert} \cite{BrownHiggins} \cite{BrownHiggins2}.

\begin{corollary}
\label{berger}
{\rm (C. Berger \cite{BergerEckmannHilton})} There is no strict $3$-groupoid $\pA $ such that
the standard realization $\Re _{\rm std} (\pA )$ is weak-equivalent to the
$3$-type of $S^2$.
\end{corollary}
{\em Proof:}
The $3$-type of $S^2$ is not a loop-space. By the previous corollary
(and the fact that the standard realizations are compatible with looping, which
we have above left as an exercise for the reader), it is impossible for
$\Re _{\rm std}(\pA )$ to be the $3$-type of $S^2$.
\eop

\section{Nonexistence of strict $3$-groupoids giving rise to the $3$-type
of $S^2$}

The present discussion aims to extend Berger's negative result 
to {\em any} realization functor satisfying the minimal
definition \ref{realizationdef}.

The first step is to prove the following statement (which
contains the main part of the argument). It basically says that the Postnikov
tower of a simply connected strict $3$-groupoid $\pC $, splits.
The intermediate $\pB$ is not really necessary for the statement but corresponds to
the technique of proof. 

\begin{proposition}
\label{diagramme}
Suppose $\pC $ is a strict $3$-groupoid with an object $c$ such that  $\pi
_0(\pC )=\ast$, $\pi _1(\pC ,c)=\{ 1\}$, $\pi _2(\pC ,c)
 = \zz$ and $\pi _3(\pC ,c)=H$ for an abelian group $H$. Then there exists a
diagram of strict $3$-groupoids
$$
\pC  \leftarr^{g}\pB  \leftarr^{f} \pA 
\rightarrt^{h} \pD 
$$
with objects $b\in \Ob (\pB )$, $a\in \Ob (\pA )$, $d\in \Ob (\pD )$ such that
$f(a)=b$, $g(b)=c$, $h(a)=d$. The diagram is
such that $g$ and $f$ are equivalences of strict $3$-groupoids, and such that
$\pi _0(\pD )=\ast$, $\pi _1(\pD ,d)=\{ 1\}$, $\pi _2(\pD ,d)=\{ 0\}$, and
such that $h$ induces an isomorphism
$$
\pi _3(h): \pi _3(\pA ,a)=H \rightarrt^{\cong} \pi _3(\pD ,d).
$$
\end{proposition}
\begin{proof}
Start with a strict groupoid $\pC $ and
object $c$, satisfying the hypotheses of \ref{diagramme}.

The first step is to construct $(\pB ,b)$. We let $\pB \subset \pC $ be the
sub-$3$-category having only one object $b=c$, and only one $1$-morphism
$1_b=1_c$. We set
$$
\Hom _{\pB (b,b)}(1_b, 1_b):=\Hom _{\pC (c,c)}(1_c, 1_c) ,
$$
with the same composition law.
The map $g: \pB \rightarrt \pC $ is the inclusion.

Note first of all that $\pB $ is a strict
$3$-groupoid. This is easily seen using version (1) of the definition
in Theorem \ref{gpdI} (but one has to look at the conditions in \cite{KapranovVoevodsky}).
We can also verify it using condition (3). Of course $\tau _{\leq 1}(\pB )$ is the
$1$-category with only one object and only one morphism, so it is a groupoid.
We have to verify that $\pB (b,b)$ is a strict $2$-groupoid. For this, we
again apply condition (3) of \ref{gpdI}. Here we note that
$$
\pB (b,b)\subset \pC (c,c)
$$
is the full sub-$2$-category with only one object $1_b=1_c$. Therefore, in view
of the definition of $\tau _{\leq 1}$, we have that
$$
\tau _{\leq 1}(\pB (b,b))\subset \tau _{\leq 1}(\pC (c,c))
$$
is a full subcategory. A full subcategory of a $1$-groupoid is again a
$1$-groupoid, so $\tau _{\leq 1}(\pB (b,b))$ is a $1$-groupoid. Finally,
$\Hom _{\pB (b,b)}(1_b, 1_b)$ is a $1$-groupoid since by construction it
is the same as  $\Hom _{\pC (c,c)}(1_c, 1_c)$ (which is a groupoid by condition
(3) applied to the strict $2$-groupoid $\pC (c,c)$). This shows that
$\pB (b,b)$ is a strict $2$-groupoid an hence that $\pB $ is a strict
$3$-groupoid.

Next, note that $\pi _0(\pB )=\ast$ and $\pi _1(\pB ,b)=\{ 1\}$. On the other hand,
for $i=2,3$ we have
$$
\pi _i(\pB ,b)= \pi _{i-2}(\Hom _{\pB (b,b)}(1_b, 1_b), 1^2_b)
$$
and similarly
$$
\pi _i(\pC ,c)= \pi _{i-2}(\Hom _{\pC (c,c)}(1_c, 1_c), 1^2_c),
$$
so the inclusion $g$ induces an equality $\pi _i(\pB ,b) \rightarrt^{=}
\pi _i(\pC ,c)$. Therefore, by definition (a) of equivalence \ref{gpdIII}, $g$ is
an equivalence of strict $3$-groupoids. This completes the construction and
verification for $\pB $ and $g$.

Before getting to the construction of $\pA $ and $f$, we analyze the strict
$3$-groupoid $\pB $ in terms of the discussion of \ref{scholium} and
\ref{scholiumgpd}. Let
$$
\pG := \Hom _{\pB (b,b)}(1_b, 1_b).
$$
It is an abelian monoid-object in the category of $1$-groupoids, with abelian
operation denoted by $+: \pG \times \pG \rightarrt \pG $ and unit element denoted $0\in
\pG $ which is the same as $1_b$. The operation $+$ corresponds to both of the
compositions $\ast _0$ and $\ast _1$ in $\pB $.

The hypotheses on the homotopy
groups of $\pC $ also hold for $\pB $ (since $g$ was an equivalence). These translate
to the statements that $(\pi _0(\pG ), +) = \zz$ and $\pG (0,0)=H$.

We now construct $\pA $ and $f$ via \ref{scholium} and
\ref{scholiumgpd}, by constructing a morphism $(\pG ',+)\rightarrt (\pG ,+)$ of
abelian monoid-objects in the category of $1$-groupoids. We do this by a type
of ``base-change'' on the monoid of objects, i.e. we will first define a
morphism $\Ob (\pG ')\rightarrt \Ob (\pG )$ and then define $\pG '$ to be the groupoid
with object set $\Ob (\pG ')$ but with morphisms corresponding to those of $\pG $.

To accomplish the ``base-change'', start with the following construction. If
$S$ is a set, let $\codisc (S)$ denote the groupoid with $S$ as set of
objects, and with exactly one morphism between each pair of objects. If $S$ has
an abelian monoid structure then $\codisc (S)$ is an abelian monoid object in
the category of groupoids.

Note that for any groupoid $\pU $ there is a morphism of groupoids
$$
\pU \longrightarrow \codisc (\Ob (\pU )),
$$
and by ``base change'' we mean the following operation: take a set $S$ with a
map $p:S\rightarrt \Ob (\pU )$ and look at
$$
\pV := \codisc (S)\times _{\codisc (\Ob (\pU ))}\pU .
$$
This is a groupoid with $S$ as set of objects, and with
$$
\pV (s,t)= \pU (p(s), p(t)).
$$
A similar construction will be used later in Chapter \ref{precat1}
under the notation $\pV = p^{\ast}(\pU )$. 
For the present purposes, note that
if $\pU $ is an abelian monoid object in the category of groupoids, if $S$ is an
abelian monoid and if $p$ is a map of monoids then $\pV $ is again an abelian
monoid object in the category of groupoids.

Apply this as follows. Starting with $(\pG ,+)$ corresponding to $\pB $
via \ref{scholium} and
\ref{scholiumgpd} as above,
choose objects $a,b \in \Ob (\pG )$ such that the image of $a$ in $\pi
_0(\pG )\cong \zz$
corresponds to $1\in \zz$, and such that the image of $b$ in $\pi _0(\pG )$
corresponds to $-1\in \zz$. Let $N$ denote the abelian monoid, product of two
copies of the natural numbers, with objects denoted $(m,n)$ for nonnegative
integers $m,n$. Define a map  of abelian monoids
$$
p:N \rightarrt \Ob (\pG )
$$
by
$$
p(m,n):= m\cdot a + n\cdot b := a+a+\ldots +a \, + \, b+b+\ldots +b.
$$
Note that this induces the surjection $N\rightarrt \pi _0(\pG )=\zz$
given by $(m,n)\mapsto m-n$.

Define $(\pG ',+)$ as the base-change
$$
\pG ':=  \codisc (N) \times _{\codisc (\Ob (\pG ))} \pG ,
$$
with its induced abelian monoid operation $+$. We have
$$
\Ob  (\pG ')= N,
$$
and the second projection $p_2: \pG '\rightarrt \pG $ (which induces $p$ on object
sets) is fully faithful i.e.
$$
{\pG '}((m,n), (m',n'))= \pG (p(m,n), p(m',n')).
$$
Note that $\pi _0(\pG ')=\zz$ via the map induced by $p$ or equivalently $p_2$.
To prove this, say that: (i) $N$ surjects onto $\zz$ so the map induced by
$p$ is
surjective; and (ii) the fact that $p_2$ is fully faithful implies that the
induced map $\pi _0(\pG ')\rightarrt \pi _0(\pG )=\zz$ is injective.

We let $\pA $ be the strict $3$-groupoid corresponding to $(\pG ',+)$
via \ref{scholium}, and let $f: \pA \rightarrt \pB $ be the
map corresponding to $p_2: \pG '\rightarrt \pG $ again via \ref{scholium}.
Let $a$ be the unique object of $\pA $ (it is mapped by $f$ to the unique object
$b\in \Ob (\pB )$).

The fact that $(\pi _0(\pG '),+)=\zz$ is a group implies that $\pA $ is a
strict $3$-groupoid (\ref{scholiumgpd}). We have $\pi _0(\pA )=\ast$ and $\pi
_1(\pA ,a)=\{ 1\}$. Also,
$$
\pi _2(\pA ,a)= (\pi _0(\pG '), +) = \zz
$$
and $f$ induces an isomorphism from here to $\pi _2(\pB ,b)=(\pi _0(\pG ), +)=\zz$.
Finally (using the notation $(0,0)$ for the unit object of $(N,+)$ and the
notation $0$ for the unit object of $\Ob (\pG )$),
$$
\pi _3(\pA ,a)= {\pG '}((0,0),(0,0)),
$$
and similarly
$$
\pi _3(\pB ,b)=\pG (0,0)=H;
$$
the map $\pi _3(f): \pi _3(\pA ,a)\rightarrt \pi _3(\pB ,b)$
is an isomorphism because it is the same as the map
$$
{\pG '}((0,0),(0,0))\rightarrt \pG (0,0)
$$
induced by $p_2: \pG '\rightarrt \pG $, and $p_2$ is fully faithful.
We have now completed the verification that $f$ induces isomorphisms on the
homotopy groups, so by  version (a) of the definition of equivalence
\ref{gpdIII}, $f$ is an equivalence of strict $3$-groupoids.

We now construct $\pD $ and define the map $h$ by an explicit calculation in
$(\pG ',+)$. First of all, let $[H]$ denote the $1$-groupoid with one object
denoted $0$, and with $H$ as group of endomorphisms:
$$
{[H]}(0,0):= H.
$$
This has a structure of abelian monoid-object in the category
of groupoids,  denoted $([H], +)$, because $H$ is an abelian group.
Let $\pD $ be the strict $3$-groupoid corresponding to $([H], +)$ via
\ref{scholium} and \ref{scholiumgpd}. We will construct a morphism
$h: \pA \rightarrt \pD $ via \ref{scholium} by constructing a morphism of abelian
monoid objects in the category of groupoids,
$$
h:(\pG ', +)\rightarrt ([H], +).
$$
We will construct this morphism so that it induces the identity  morphism
$$
{\pG '}((0,0), (0,0))=H \rightarrt \leftbrack H](0,0)=H.
$$
This will insure that the morphism $h$ has the property required for
\ref{diagramme}.

The object $(1,1)\in N$ goes to $0\in \pi _0(\pG ')\cong \zz$. Thus we may choose
an isomorphism $\varphi : (0,0)\cong (1,1)$ in $\pG '$.  For any $k$ let $k\varphi$
denote the isomorphism $\varphi + \ldots +\varphi$ ($k$ times) going from
$(0,0)$ to $(k,k)$.  On the other hand, $H$ is the automorphism group of
$(0,0)$ in $\pG '$. The operations $+$ and composition coincide on $H$. Finally,
for any $(m,n)\in N$ let $1_{m,n}$ denote the identity automorphism of the
object $(m,n)$.  Then any arrow $\alpha$ in $\pG $ may be uniquely written in the
form
$$
\alpha = 1_{m,n} + k\varphi + u
$$
with $(m,n)$ the source of $\alpha$, the target being $(m+k, n+k)$, and where
$u\in H$.

We have the following formulae for the composition $\circ$ of arrows in $\pG '$.
They all come from the basic rule
$$
(\alpha \circ \beta ) + (\alpha ' \circ \beta ')=
(\alpha + \alpha ') \circ (\beta + \beta ')
$$
which in turn comes simply from the fact that $+$ is a morphism of groupoids
$\pG '\times \pG '\rightarrt \pG '$ defined on the cartesian product of two copies of
$\pG $. Note in a similar vein that $1_{0,0}$ acts as the identity for the
operation
$+$ on arrows, and also that
$$
1_{m,n} + 1_{m',n'} = 1_{m+m', n+n'}.
$$

Our first equation is
$$
(1_{l,l} +k\varphi )\circ l\varphi = (k+l)\varphi .
$$
To prove this note that $l\varphi + 1_{0,0}= l\varphi$ and our basic formula
says
$$
(1_{l,l}\circ l_{\varphi} ) + (k\varphi \circ 1_{0,0})
=
(1_{l,l} +k\varphi )\circ (l\varphi + 1_{0,0} )
$$
but the left side is just $l\varphi + k\varphi = (k+l)\varphi$.

Now our basic formula, for a composition starting with
$(m,n)$, going first to $(m+l,n+l)$, then going to $(m+l+k, n+l+k)$, gives
$$
(1_{m+l,n+l} + k\varphi + u)\circ (1_{m,n} + l\varphi + v)
$$
$$
= (1_{m,n} + 1_{l,l} + k\varphi + u)\circ (1_{m,n} + l\varphi + v)
$$
$$
= 1_{m,n}\circ 1_{m,n}  + (1_{l,l} +k\varphi )\circ l\varphi
+ u\circ v
$$
$$
= 1_{m,n} + (k+l)\varphi + (u\circ v)
$$
where of course $u\circ v=u+v$.

This formula shows that the morphism $h$ from arrows of $\pG '$ to the group $H$,
defined by
$$
h(1_{m,n} + k\varphi + u):= u
$$
is compatible with composition.  This implies that it provides a morphism of
groupoids $h:\pG \rightarrt \leftbrack H]$ (recall from above that $[H]$ is defined to
be the
groupoid with one object whose automorphism group is $H$).  Furthermore the
morphism $h$ is obviously compatible with the operation $+$ since
$$
(1_{m,n} + k\varphi + u)+ (1_{m',n'} + k'\varphi + u')=
$$
$$
(1_{m+m',n+n'} + (k+k')\varphi + (u+u'))
$$
and once again $u+u'=u\circ u'$ (the operation $+$ on $[H]$ being given by the
commutative operation $\circ$ on $H$).

This completes the construction of a morphism
$h: (\pG , +)\rightarrt ([H], +)$ which induces the identity on $\Hom (0,0)$.
This corresponds to a morphism of strict $3$-groupoids $h: \pA \rightarrt \pD $
as required to complete the proof of Proposition \ref{diagramme}.
\end{proof}

We can now give the nonrealization statement. 

\begin{theorem}
\label{noS2}
Let $\Re$ be any realization functor satisfying the properties of Definition
\ref{realizationdef}. Then there does not exist a strict $3$-groupoid $\pC $ such
that $\Re (\pC )$ is weak-equivalent to the $3$-truncation of the homotopy type of
$S^2$.
\end{theorem}
\begin{proof}
Suppose for the moment that we know Proposition \ref{diagramme}; with this we
will prove \ref{noS2}. Fix a realization functor $\Re$ for strict
$3$-groupoids satisfying the axioms \ref{realizationdef}, and assume that $\pC $ is
a strict $3$-groupoid such that $\Re (\pC )$ is weak homotopy-equivalent to the
$3$-type of $S^2$. We shall derive a contradiction.

Apply Proposition \ref{diagramme} to $\pC $. Choose an object $c\in \Ob (\pC )$. Note
that, because of the isomorphisms between homotopy sets or groups
\ref{realizationdef}, we have $\pi _0(\pC )=\ast$,
$\pi _1(\pC ,c)=\{ 1\}$, $\pi _2(\pC ,c)
 = \zz$ and $\pi _3(\pC ,c)=\zz $, so \ref{diagramme} applies with $H=\zz$.
We obtain a sequence of strict $3$-groupoids
$$
\pC  \leftarr^{g} \pB  \leftarr^{f} \pA 
\rightarrt^{h} \pD .
$$
This gives the diagram of spaces
$$
\Re (\pC ) \leftarr^{\Re (g)} \Re (\pB ) \leftarr^{\Re
(f)} \Re
(\pA )  \rightarrt^{\Re (h)} \Re (\pD ).
$$
The axioms \ref{realizationdef} for $\Re$ imply that $\Re$ transforms
equivalences of strict $3$-groupoids into weak homotopy equivalences of spaces.
Thus $\Re (f)$ and $\Re (g)$ are weak homotopy equivalences and we get that
$\Re (\pA )$ is weak homotopy equivalent to the $3$-type of $S^2$.

On the other hand, again by the
axioms \ref{realizationdef}, we have that $\Re (\pD )$ is $2$-connected, and $\pi
_3(\Re (\pD ), r(d))=H$ (via the isomorphism $\pi _3(\pD ,d)\cong H$ induced by $h$,
$f$ and $g$).
By the Hurewicz theorem there is a class $\eta \in H^3(\Re (\pD ),
H)$ which induces an isomorphism
$$
{\bf Hur}(\eta ): \pi _3(\Re (\pD ), r(d))\rightarrt^{\cong} H.
$$
Here
$$
{\bf Hur} : H^3(X , H)\rightarrt Hom (\pi _3(X,x), H)
$$
is the Hurewicz map for any pointed space $(X,x)$; and the cohomology is
singular cohomology (in particular it only depends on the weak homotopy type of
the space).

Now look at the pullback of this class
$$
\Re (h)^{\ast}(\eta )\in H^3(\Re (\pA ), H).
$$
The hypothesis that $\Re (u)$ induces an isomorphism on $\pi _3$ implies that
$$
{\bf Hur}(\Re (h)^{\ast}(\eta )): \pi _3(\Re (\pA ),
r(a))\rightarrt^{\cong} H.
$$
In particular, ${\bf Hur}(\Re (h)^{\ast}(\eta ))$ is nonzero so
$\Re (h)^{\ast}(\eta )$ is nonzero in $H^3(\Re (\pA ), H)$. This is a
contradiction because $\Re (\pA )$ is weak homotopy-equivalent to the $3$-type
of $S^2$, and $H=\zz$, but $H^3(S^2 , \zz )=\{ 0 \}$.

This contradiction completes the proof of the theorem.
\end{proof}

As was discussed in \cite{Grothendieck}, this result motivates the search for
a notion of higher category weaker than the notion of strict $n$-category. 
Following the yoga described by Lewis \cite{Lewis}, it
appears to be sufficient to weaken any single particular aspect.


\chapter{Simplicial approaches}
\label{simplicial1}

There are a number of approaches to weak higher categories based on the simplicial category $\Delta$, including the
Segal approach and its iterations which are the main subject of our book. We also discuss several other approaches,
which concern first and foremost the theory of $(\infty ,1)$-categories. 

\section{Strict simplicial categories}

A {\em simplicial category} is a $\mK$-enriched category. It has a set of objects $\Ob (\pA )$, and for each pair
$x,y\in \Ob (\pA )$ a simplicial set $\pA (x,y)$ thought of as the ``space of morphisms'' from $x$ to $y$. The composition maps
are morphisms of simplicial sets $\pA (x,y)\times \pA (y,z)\rightarrt \pA (x,z)$ satisfying the associativity condition strictly,
that is for any $x,y,z,w$ the diagram of simplicial sets 
$$
\begin{diagram}
\pA (x,y)\times \pA (y,z)\times \pA (z,w) & \rightarr & \pA (x,y)\times \pA (y,w) \\
\downarr & & \downarr \\
\pA (x,z)\times \pA (z,w)& \rightarr &  \pA (x,w) 
\end{diagram}
$$
commutes. The identities of $\pA $ are points (i.e. vertices) of the simplicial sets $\pA (x,x)$ satisfying left and right identity conditions
which are equalities of maps $\pA (x,y)\rightarrt \pA (x,y)$. 

A {\em functor} of simplicial categories $f:\pA \rightarrt \pB $ consists of a map $f:\Ob (\pA )\rightarrt \Ob (\pB )$ and for each $x,y\in\Ob (\pA )$,
a  map of simplicial sets $f_{x,y}:\pA (x,y)\rightarrt \pB (f(x),f(y))$, compatible with the composition maps and identities in an obvious way.
In keeping with our general notation for enriched categories, the category of simplicial categories is denoted $\Cat (\mK )$. 

Given a simplicial category $\pA $, we define its {\em truncation} $\tau _{\leq 1}(\pA )$ to be the category whose set of objects is the same
as $\Ob (\pA )$, but for any $x,y\in \Ob (\pA )$
$$
\tau _{\leq 1}(\pA )(x,y):= \pi _0(\pA (x,y)).
$$
The composition maps and identities for $\pA $ define composition maps and identities for $\tau _{\leq 1}(\pA )$, and we obtain a functor
$$
\tau _{\leq 1}: \Cat (\mK )\rightarrt \Cat .
$$

A functor $f:\pA \rightarrt \pB $ between simplicial categories is said to be {\em fully faithful} if for every $x,y\in \Ob (\pA )$ the map
$f_{x,y}: \pA (x,y)\rightarrt \pB (f(x),f(y))$ is a weak equivalence of simplicial sets, in other words a weak equivalence in the standard
model structure of $\mK$. A functor $f$ is said to be {\em essentially surjective} if the functor $\tau _{\leq 1}(f)$ between usual $1$-categories
is essentially surjective, in other words it induces a surjection on sets of isomorphism classes $\Iso \tau _{\leq 1} (\pA ) \twoheadrightarrow 
\Iso \tau _{\leq 1} (\pB )$. 
A functor $f:\pA \rightarrt \pB $ is said to be a {\em Dwyer-Kan equivalence} between simplicial categories, if it is fully faithful and essentially surjective.
In this case, $\tau _{\leq 1}(f)$ is also an equivalence of categories, in particular it is bijective on sets of isomorphism classes. 

Given a simplicial category $\pA $, its {\em underlying category} is the category with objects $\Ob (\pA )$, but the morphisms from $x$ to $y$ are
the set of points or vertices of the simplicial set $\pA (x,y)$. This is not to be confused with $\tau _{\leq 1}(\pA )$, but there is a natural projection functor from the
underlying category to the truncated category. An ``arrow'' in $\pA $ from $x$ to $y$
means a map in this underlying category; such an arrow is said to be an {\em internal equivalence} if it projects to an isomorphism in $\tau _{\leq 1}(\pA )$. 
In these terms, the essential surjectivity condition for a functor $f:\pA \rightarrt \pB $ may be rephrased as saying that any object of $\pB $ is internally
equivalent to the image of an object of $\pA $. 

Dwyer, Kan and Bergner have constructed a model category structure on $\Cat (\mK )$ such that the weak equivalences are the Dwyer-Kan equivalences,
and the fibrations are the functors $f:\pA \rightarrt \pB $ such that each $f_{x,y}$ is a fibration of simplicial sets, and furthermore $f$ satisfies
an additional lifting condition which basically says that an internal equivalence in $\pB $ should lift to $\pA $ if one of its endpoints lifts.

It is interesting to note that Dwyer and Kan started first by constructing a model structure on $\Cat (X,\mK )$, the category of simplicial categories
with a fixed set of objects $X$. Refer to their papers \cite{DK1} \cite{DK2} \cite{DK3}.
We will also adopt this route, following a suggestion by Clark Barwick. 

Simplicial categories appear in an important way in homotopy theory. Quillen defined the notion of {\em simplicial model category},
and of $\mN$ is a simplicial model category then we obtain a simplicial category $\mN ^{\rm spl}_{cf}$ of fibrant and cofibrant objects, such that its truncation
$\tau _{\leq 1}(\mN ^{\rm spl}_{cf})\cong \Ho (\mN )$ is the homotopy category of $\mN$. Dwyer and Kan then developped the theory of {\em simplicial 
localization} which gives a good simplicial category even when $\mN$ doesn't have a simplicial model structure. If $\mC$ is any category, and 
if we are given a subcollection of arrows $\Sigma \subset \Arr (\mC )$, then Dwyer and Kan define a simplicial category $L(\mC , \Sigma )$
whose truncation is the classical Gabriel-Zisman localization: $\tau _{\leq 1}(L(\mC , \Sigma ))\cong \Sigma ^{-1}(\mC )$.
In the case where $\mN$ is a simplicial model category, then the two options $L(\mN , \mN_w )$ (where $\mN_w$ denotes the class of weak equivalences) and
$\mN ^{\rm spl}_{cf}$, are Dwyer-Kan equivalent as simplicial categories \cite{DK2} \cite{DK3}. 

Even though simplicial categories have strictly associative composition, they
are weaker than strict $n$-categories in the sense that the higher categorical structure is
encoded by the simplicial morphism sets rather than by strict $n-1$-categories. 
Hence, the need for weak composition described in the previous chapter, is not
contradicatory with the fact that strict simplicial categories model all 
$(\infty ,1)$-categories. For the weaker versions to be discussed next, 
one can rectify back
to a strict simplicial category, as was originally shown by Dwyer-Kan-Smith \cite{DKS} and 
Schw\"{a}nzl-Vogt \cite{Vogt}, then extended to Quillen equivalences between
the corresponding model structures by Bergner \cite{BergnerThreeModels}. 

\section{Segal's delooping machine}

The best-known version of Segal's theory is his notion of infinite delooping machine
or $\Gamma$-space. Grothendieck mentioned some correspondence from Larry Breen 
in 1975 concerning this idea: 
\label{grobreenquote}
\begin{quote}
Dear Larry, 

\ldots 
The construction which you propose for the notion of a non-strict $n$-category, and of the nerve functor, has certainly the merit of existing, and of being a first precise approach
\ldots 

Otherwise, not having understood the idea of Segal in your last letter\ldots 
\end{quote} 

In the first letter of 1983, Grothendieck also 
mentioned the notion of multisimplicial 
nerve of a (strict) $n$-category. So it would seem that the idea of applying
Segal's delooping machine much as was done in \cite{Tamsamani}, was present in some sense
at the time. 

The starting point is Segal's $1$-delooping machine. 
Recall from topology that for a pointed space $(X,x)$ the {\em loop space} $\Omega X$ is the space
of pointed loops $(S^1,0)\rightarrt (X,x)$. These can be composed by reparametrizing the loops in a well-known way,
although the resulting composition is only associative and unital up to homotopy. It is possible to replace $\Omega X$ by
a topological group, for example with Quillen's realization of homotopy types as coming from simplicial groups. 
A classifying space construction usually denoted $B(\cdot )$ allows one to get back to the original space:
$$
B(\Omega X) \sim X.
$$
A popular question in topology in the 1960's was how to define various types of structure on spaces homotopy equivalent to 
$\Omega X$, which would be weaker than the strong structure of topological group, but which would include sufficiently
much homotopical data to let us get back to $X$ by a classifying space construction $B(\cdot )$. Such a kind of structure was known
as a ``delooping machine''. There were a number of examples including $A_{\infty}$-algebras (in the linearized case),
PROP's, operads, and the one which we will be
considering: Segal's simplicial delooping machine. 

Any of these delooping machines should lead to one or several notions of higher category, indeed this has been the case as we shall discuss elsewhere.

Let $\Delta$ denote the simplicial category whose
objects are denoted $m$ for positive integers $m$, and where the morphisms
$p\rightarrt m$ are the (not-necessarily strictly) order-preserving maps
$$
\{ 0,1,\ldots , p\} \rightarrt \{ 0,1,\ldots , m\} .
$$
A morphism $1\rightarrt m$ sending $0$ to $i-1$ and $1$ to $i$ is called a
{\em principal edge} of $m$.  A morphism which is not injective is called a
{\em degeneracy}.

A simplicial set $\pA : \Delta ^o\rightarrt \Sets$ such that $\pA _0=\ast$ and such that 
the {\em Segal maps} obtained by the principal edges 
$(01), (12), \ldots , ((n-1)n)\subset (0123\cdots n)= [n]$
$$
\pA _m \rightarrt \pA _1 \times \cdots \times \pA _1
$$
are isomorphisms of sets, corresponds to a structure of monoid on the set $\pA _1$. Indeed, the diagram 
$$
\begin{diagram}
\pA _2 & \rightarr^{\cong} & \pA _1\times \pA _1 \\
\downarr & & \\
\pA _1 & & 
\end{diagram}
$$
where the horizontal map is the Segal map and the vertical map is given by the third edge $(02)\subset (012)$,
provides a composition law $\pA _1\times \pA _1\rightarrt \pA _1$. The degeneracy map $\pA _0\rightarrt \pA _1$ provides a 
unit---proved using the degeneracy maps for $\pA _2$---and consideration of $\pA _3$ gives the proof of associativity.

In Segal's $1$-delooping theory, this characterization of monoids is weakened by replacing the condition of isomorphism by
the condition of weak homotopy equivalences (i.e. maps inducing isomorphisms on the $\pi _i$). Thus, 
a loop space is defined to be a simplicial space 
$$
\pA _{\cdot} : \Delta ^o \rightarrt \Top
$$
such that $\pA _0$ is a single point, such that the  Segal maps, again using the principal edges 
$(01), (12), \ldots , ((n-1)n)\subset (0123\cdots n)= [n]$
$$
\pA _m \rightarrt \pA _1 \times \cdots \times \pA _1
$$
are weak homotopy equivalences, and which is {\em grouplike} in that the monoid which results when we compose 
$$
\pi _0 \circ \pA _{\cdot} : \Delta ^o \rightarrt \Top \rightarrt \Sets 
$$
should be a group. Segal explains in \cite{Segal} how to deloop such an object: if $\Top$ is replaced by the
category of simplicial sets then the structure $\pA _{\cdot}$ is a bisimplicial set, and its delooping $B(\pA {\cdot})$ is just the diagonal realization. 

As was well known at the time, the characterization of monoids generalizes to give a characterization of the nerve of a category in terms
of Segal maps being isomorphisms. Indeed, a monoid can be viewed as a category with a single object, and a small change in the definition
makes it apply to the case of categories with an arbitrary set of objects: a simplicial set
$$
\pA _{\cdot} : \Delta ^o \rightarrt \Sets
$$
is the nerve of a $1$-category if and only if the Segal maps made using fiber products are isomorphisms
$$
\pA _m \rightarrt^{\cong} \pA _1 \times _{\pA _0}\cdots \times _{\pA _0}\pA _1.
$$
Here the fiber products are taken over the two maps $\pA _1\rightarrt \pA _0$ corresponding to $(0)\subset (01)$ and $(1)\subset (01)$,
alternatingly and starting with $(1)\subset (01)$. These correspond to the inclusions of the intersections of adjacent principal edges.

\section{Segal categories}

A {\em Segal precategory} is a bisimplicial set
$$
\pA = \{ \pA _{p,k},
\;\; p,k\in \Delta \}
$$
in other words a  functor $\pA : \Delta ^o\times \Delta ^o\rightarrt \Sets$
satisfying the {\em globular condition} that the simplicial set $k\mapsto
\pA _{0,k}$ is constant equal to a set which we denote by $\pA _0$ (called the set of
{\em objects}).

If $\pA $ is a Segal precategory then for $p\geq 1$ we obtain a simplicial set
$$
k\mapsto \pA _{p,k}
$$
which we denote by $\pA _{p/}$.  This yields a simplicial collection of simplicial
sets, or a functor $\Delta ^o\rightarrt \mK$ to the Kan-Quillen model category
$\mK$ of simplicial sets. 
One could instead look at {\em simplicial spaces}, that is functors $\Delta ^o\rightarrt \Top$ such that $\pA _0$ is a discrete space thought of as a set.  This gives an equivalent
theory, although there are degeneracy problems which apparently need to be
treated in an appendix in that case \cite{MayThomason} \cite{Thomason}. We often
use the ``simplicial space'' point of view when speaking informally, as it is more
intuitively compelling; however, we don't want to get into details of defining a model
structure on $\Top$, and instead use $\mK$ for technical statements. 

For each $m\geq 2$ there is a morphism of simplicial sets whose components are
given by the principal edges of $m$, which we call the {\em Segal map}:
$$
\pA _{m/}\rightarrt \pA _{1/} \times _{\pA _0} \ldots \times _{\pA _0} \pA _{1/}.
$$
The morphisms in the fiber product $\pA _{1/}\rightarrt \pA _0$ are alternatively the
inclusions $0\rightarrt 1$ sending $0$ to the object $1$, or to the object $0$.

We would like to think of the inverse image $\pA _{1/}(x,y)$ of a pair
$(x,y)\in \pA _0\times \pA _0$ by the two maps $\pA _{1/}\rightarrt \pA _0$ referred to
above, as the {\em simplicial set of maps from $x$ to $y$}.

We say that a Segal precategory $\pA $ is  a {\em Segal category} if for all $m\geq
2$ the Segal maps
$$
\pA _{m/}\rightarrt \pA _{1/} \times _{\pA _0} \ldots \times _{\pA _0} \pA _{1/}.
$$
are weak equivalences of simplicial sets. This notion was introduced by
Dwyer, Kan and Smith \cite{DKS} and Schw\"{a}nzl and Vogt \cite{Vogt}. 

Given a strict simplicial category $\pA$, we obtain a corresponding Segal precategory
by setting 
$$
\pA _{n/}:= \coprod _{(x_0,\ldots , x_n)\in \Ob (\pA )^{n+1}} \pA (x_0,x_1)\times
\cdots \times \pA (x_{n-1},x_n).
$$
This is a Segal category, because the Segal maps are {\em isomorphisms}. In the other
direction, a Segal category such that the Segal maps are isomorphisms comes from a unique
strict simplicial category.  

The ``generators and relations'' operation introduced in Chapter \ref{genrel1} is a way of starting with a Segal precategory
and enforcing the condition of becoming a Segal category, by forcing the
condition of weak equivalence on the Segal maps. As a general matter we will
call operations of this type $\pA \mapsto \Seg (\pA )$. 

Suppose $\pA $ is a Segal category. Then the simplicial set $p\mapsto \pi
_0(\pA _{p/})$ is the nerve of a category which we call $\tau _{\leq 1} \pA $.
We say that $\pA $ is a {\em Segal groupoid} if $\tau _{\leq 1} \pA $ is a groupoid.
This means that the $1$-morphisms of $\pA $ are invertible up to equivalence.

In fact we can make the same definition even for a Segal precategory $\pA $: we define
$\tau _{\leq 1} \pA $ to be the simplicial set $p\mapsto \pi _0(\pA _{p/})$.

We can now describe exactly the situation envisaged in \cite{Adams}
\cite{Segal}: a Segal category $\pA $ with only one object, $\pA _0 = \ast$.
We call this a {\em Segal monoid}. If $\pA $ is a groupoid then the homotopy
theorists' terminology is to say that it is {\em grouplike}.

\subsection{Equivalences of Segal categories}

The basic intuition is to think of Segal categories as the natural weak version
of the notion of topological (i.e. $\Top$-enriched) 
category. One of the main concepts in category
theory is that of a functor which is an ``equivalence of categories''.
This may be generalized to Segal categories. 
For simplicial (i.e. $\mK$-enriched) categories, this notion is due to Dwyer and
Kan, and is often called {\em DK-equivalence}. 
The same thing in the context
of $n$-categories is well-known (see Kapranov-Voevodsky \cite{KapranovVoevodsky} for
example); in the weak case it is described in Tamsamani's paper \cite{Tamsamani}.

We say
that a morphism $f:\pA \rightarrt \pB $ of Segal categories is an {\em equivalence}
if it is {\em fully faithful}, meaning that for $x,y\in \pA _0$ the map
$$
\pA _{1/}(x,y) \rightarrt \pB _{1/}(f(x),f(y))
$$
is a weak equivalence of simplicial sets; and  {\em essentially
surjective}, meaning that the induced functor of categories
$$
\tau _{\leq 1} (\pA )\rightarrt^{\tau _{\leq 1}f} \tau _{\leq 1}(\pB )
$$
is surjective on isomorphism classes of objects. Note that this induced
functor $\tau _{\leq 1}f$ will be an equivalence of categories as a consequence of the
fully faithful condition.

The homotopy theory that we are interested in is that of the category of Segal
categories modulo the above notion of equivalence. In particular, when we
search for the ``right answer'' to a question, it is only up to the above type
of equivalence. Of course when dealing with Segal categories having only one
object (as will actually be the case in what follows) then the essentially
surjective condition is vacuous and the fully faithful condition just amounts
to equivalence on the level of the``underlying space'' $\pA _{1/}$.

In order to have an appropriately reasonable point of view on the homotopy
theory of Segal categories one should look at the closed model structure
(which is one of our main goals, specialized to the model category
$\mK$, see Chapter \ref{secat1}):  the right notion of weak morphism
from $\pA $ to $\pB $ is that of a morphism from $\pA $ to $\pB '$ where $\pB \hookrightarrow
\pB '$ is a fibrant replacement of $\pB $.  

\subsection{Segal's theorem}

We define the {\em realization} of a Segal category $\pA $ to be the space $|\pA |$
which is the realization of the bisimplicial set $\pA $.  Suppose $\pA _0 = \ast$.
Then we have a morphism
$$
|\pA _{1/} | \times {[}0,1] \rightarrt |\pA |
$$
giving a morphism
$$
|\pA _{1/}| \rightarrt \Omega |\pA |.
$$
The notation $|\pA _{1/}|$ means the realization of the simplicial set $\pA _{1/}$
and $\Omega |\pA |$ is the loop space based at the basepoint $\ast = \pA _0$.

\begin{theorem}
\mylabel{segal}
{\rm (G. Segal \cite{Segal}, Proposition 1.5)}
Suppose $\pA $ is a Segal groupoid with one object. Then the morphism
$$
|\pA _{1/}| \rightarrt \Omega |\pA |.
$$
is a weak equivalence of spaces.
\end{theorem}

Refer to Segal's paper, or also May (\cite{MayFibs} 8.7), for a proof.
Tamsamani noted that the same works in the case of many objects, and indeed this was
a key step in his proof of the topological realization theorem for $n$-categories.

\begin{corollary}
Suppose $\pA$ is a Segal groupoid. Then the morphism
$$
|\pA _{1/}| \rightarrt \Omega |\pA |.
$$
is a weak equivalence in $\mK$.
\end{corollary}
\begin{proof}
\cite{Tamsamani}.
\end{proof}

In order to do these things inside the world of simplicial spaces,
the additional cofibrancy conditions in $\Top$ would necessitate a discussion of ``whiskering'' as is standard in delooping and classifying space
constructions (cf \cite{Segal} \cite{MayLoops}, \cite{MayThomason},
\cite{Thomason}). This is why we have replaced ``spaces''
by ``simplicial sets'' in the above discussion, and corresponds also to our use of Reedy
model structures in the main chapters. 

\subsection{$(\infty , 1)$-categories}

Simplicial categories and Segal categories are two models for what Lurie calls the notion of {\em $(\infty ,1)$-category}, meaning $\infty$-categories where
the $i$-morphisms are invertible (analogous to being an inner equivalence) for $i\geq 2$. 
The Dwyer-Kan simplicial localization may be viewed as the localization in the $(\infty , 2)$-category of $(\infty ,1)$-categories. 
Part of our goal in this book is to develop algebraic formalism useful for looking at these situations, as well as their iterative counterparts for $(\infty ,n)$-categories. 

\subsection{Iteration}

A subtle point is that simplicial categories don't behave well under direct products: the Dwyer-Kan-Bergner model structure on $\Cat (\mK )$
is not cartesian for the direct product because a product of two cofibrant objects is no longer cofibrant. Thus, if we try to continue by looking at 
$\Cat (\Cat (\mK ))$ the resulting theory doesn't have the right properties. 
This hooks up with what we have seen in Chapter \ref{nonstrict1}, that iterating the strict enriched category construction, doesn't lead to enough objects.

If we use Segal's method, on the other hand, one can iterate the construction with a
better effect. 
This leads to Tamsamani's iterative definition of $n$-categories. See Chapter \ref{cartenr1}.
It is also related to Dunn's theory of iterated $n$-fold Segal delooping machines
\cite{Dunn}, and it will undoubtedly be profitable to compare \cite{Tamsamani} and 
\cite{Dunn}.

The next two sections will be devoted to brief descriptions of two other major points of view on $(\infty ,1)$-categories. After that, we discuss the comparison between
the various theories. 

\subsection{Strictification and Bergner's comparison result}

The various models of $(\infty ,1)$-categories discussed above all furnish essentially the same homotopy theory.  Such a rectification result was known very early for homotopy
monoid structures. For Segal categories, the first result of this kind was due to Dwyer, Kan and Smith who showed how to rectify a Segal category into a strict simplicial category
in \cite{DKS}. Similarly, Scw\"{a}nzl and Vogt showed the same thing in their paper
introducing Segal categories \cite{Vogt}. The full homotopy equivalence result,
stating that the rectification operation is a Quillen equivalence between model
structures,  was shown by Bergner \cite{BergnerThreeModels}
at the same time as she constructed the
model structures in question. The model structure for Segal categories
is the special case $\mM = \mK$ of the global construction we are doing in the present book.

With respect to the models we are going to discuss next,
Bergner also gave a Quillen equivalence with Rezk's model category of complete Segal
spaces, and the comparison can be extended to quasicategories too, as shown by Joyal and Tierney in \cite{JoyalTierneyQCSC}.

\section{Rezk categories}

Rezk has given a different way of using the Segal maps to specify an $(\infty ,1)$-categorical structure. 
Barwick showed how to iterate this construction, and this iteration has now also been taken up by Lurie and
Rezk. Their iteration is philosophically similar to what we are doing in the main part of this book. In the present section
we discuss Rezk's original case, which he called ``complete Segal spaces''. These objects enter into Bergner's three-way comparison
\cite{BergnerThreeModels}. 

It will be convenient to start our discussion by refering to a generic notion of
$(\infty ,1)$-category, which could be concretized by simplicial categories, or Segal categories.
Recall that an $(\infty ,1)$-category $\pA$ has an $(\infty ,0)$-category or $\infty$-groupoid,
as its {\em interior} denoted $\pA ^{\rm int}$. The interior is the universal $\infty$-groupoid mapping to $\pA$. For any $x,y\in \Ob (\pA )$, the mapping space
is the subspace
$\pA ^{\rm int}(x,y)\subset \pA (x,y)$ union of all the connected components corresponding
to maps which are invertible up to equivalence. By Segal's theorem (which will be discussed
further in Chapter \ref{secat1}), this corresponds to a space which we can denote by 
$|\pA ^{\rm int}|$. It is the ``moduli space of objects of $\pA$ up to equivalence'':
there is a separate connected component for each equivalence class of objects. 
The vertices coming from the $0$-simplices correspond to the original objects $\Ob (\pA )$,
and
within a connected component the space of paths from one vertex to another, is
the space $\pA ^{\rm int}(x,y)$
of equivalences between the corresponding objects.

In Rezk's theory, our $(\infty ,1)$-category is represented by a simplicial space $\pA ^R$
with $\pA ^R_0= |\pA ^{\rm int}|$ in degree zero. The {\em homotopy fiber} of
the map 
$$
\pA ^R_1\rightarrt \pA ^R_0\times \pA ^R_0
$$
over a point $(x,y)$ is (canonically equivalent to) the space of morphisms $\pA (x,y)$.
The categorical structure is defined by imposing a Segal condition on homotopy
fiber products: for any $n$, there is a version of the Segal map going to the
homotopy fiber product 
$$
\pA ^R_n\rightarrt \pA ^R_1\times^h_{\pA ^R_0}\pA ^R_1\times^h_{\pA ^R_0}\cdots 
\times^h_{\pA ^R_0}\pA ^R_1
$$
and this is required to be a weak equivalence. A {\em complete Segal space} is a
simplicial space satisfying these Segal conditions, and also the {\em completeness}
condition which corresponds to the requirement $\pA ^R_0= |\pA ^{\rm int}|$.
We had formulated that requirement by first considering a generic theory of $(\infty ,1)$-categories. 
Internally to Rezk's theory, the completeness condition says that $\pA ^{R,{\rm int}}_1\rightarrt
\pA ^R_0\times \pA ^R_0$ should be equivalent to the path space fibration, where
$\pA ^{R,{\rm int}}_1\subset \pA ^R_1$ denotes the union of connected components
corresponding to morphisms which are invertible up to equivalence. This condition
is shown to be equivalent to a more abstract condition useful for manipulating
the model structure, see \cite[6.4]{Rezk} \cite[3.7]{BergnerThreeModels}
\cite[Section 4]{JoyalTierneyQCSC}. 

Rezk's theory is a little bit more complicated in its initial stages than the theory of 
Segal categories. The Segal maps go to a homotopy fiber product, which nevertheless
can be assumed
to be a regular fiber product by imposing a Reedy fibrant condition on the
simplicial space, for example. Since the set of objects is not really too well-defined,
the kind of reasoning which we are considering here (and which was also followed
by Dwyer and Kan in their series of papers), breaking up the problem into 
first a problem for higher categories with a fixed set of objects, then varying
the set of objects, is less available. 

On the other hand, Rezk's theory has the advantage that $\pA^R$
is a canonical model for $\pA$ up to levelwise homotopy equivalence 
in the category of diagrams $\Delta ^o\rightarrt \Top$. Thus, a map
$\pA ^R\rightarrt \pB ^R$ of complete Segal spaces, is an equivalence if and only if
each $\pA ^R_n\rightarrt \pB ^R_n$ is a weak equivalence of spaces (and it suffices
to check $n=0$ and $n=1$ because of the Segal conditions).
This contrasts with the
case of Segal categories, where the set of objects 
$\pA _0=\Ob (\pA )$ is not invariant under equivalences
of categories. 

As Bergner has pointed out \cite{BergnerHoFiProd}, the canonical nature of the spaces
involved makes Rezk's theory particularly amenable to calculating limits.
For example, if 
$$
\pA ^R\rightarrt \pB^R \leftarr \pC ^R
$$
are two arrows between complete Segal spaces, then the levelwise homotopy fiber product
$$
\pU ^R_n:= \pA ^R_n\times _{\pB ^R_n}^h \pC ^R_n
$$
is again a complete Segal space, and it is the right homotopy fiber product in the world
of $(\infty , 1)$-categories. 

This again contrasts with the case of Segal categories, or indeed even usual $1$-categories.
For example, letting $\pE$ denote the $1$-category with two isomorphic objects $\upsilon _0$ and $\upsilon _1$, the inclusion maps 
$$
\{ \upsilon _0 \} \rightarrt \pE \leftarr \{ \upsilon _1\}
$$
are equivalences of categories, so the homotopy fiber product in any reasonable model structure for $1$-categories, should also be equivalent to a discrete singleton category.
However, the  fiber product of categories, or of simplicial sets (the nerves) is empty.
In Rezk's theory, the degree $0$ space $\pE ^R_0$ will again be contractible,
since there is only a single equivalence class of objects of $\pE$ and the have no 
nontrivial automorphisms. 

As usual, for treating the technical aspects of the theory it is better to look at
bisimplicial sets rather than simplicial spaces. Rezk constructs a model structure on
the category of bisimplicial sets, such that the fibrant objects are complete Segal 
spaces which are Reedy fibrant as $\Delta ^o$-diagrams and levelwise fibrant \cite{Rezk}.
Bergner considers further this theory and shows the equivalence with simplicial categories
and Segal categories \cite{BergnerThreeModels}. 

Barwick has suggested to iterate this construction to  a Rezk-style theory of 
$(\infty , n)$ categories for all $n$, and Rezk has taken this up in \cite{RezkCartesian}.
He shows that the resulting model categories are cartesian, in particular this gives
a construction of the $(\infty , n+1)$-category of $(\infty , n)$-categories.

\section{Quasicategories}

Joyal and Lurie have developped extensively the theory of {\em quasicategories}. These
first appeared in the book of Boardman and Vogt \cite{BoardmanVogt} under the name
``restricted Kan complexes''. An important example appeared in work of Cordier and Porter \cite{CordierPorter}.

A good place to start is to recall Kan's original horn-filling conditions for 
the category of simplicial sets $\mK$. As $\mK$ is a category of diagrams $\Delta ^o\rightarrt \Sets$, we have in particular the {\em representable diagrams} 
which we shall denote $R(n)$, defined by $R(n)_m:= \Delta ([m],[n])$. 
This is the ``standard $n$-simplex'', classically denoted by $R(n)=\Delta [n]$.
For our purposes this classical
notation would seem to risk some confusion with too many symbols $\Delta$ around,
so we call it $R(n)$ instead. Now, $R(n)$ has a standard simplicial subset denoted 
$\partial R(n)$, which is the ``boundary''. It can be defined as the $n-1$-skeleton
of $R(n)$, or as the union of the $n-1$-dimensional faces of $R(n)$. The
faces are indexed by $0\leq k\leq n$; in terms of linearly ordered sets, the
$k$-th face corresponds to the linearly ordered subset of $[n]$ obtained by
crossing out the $k$-th element. Now, the $k$-th {\em horn} $\Lambda (n,k)$ is
the subset of $\partial R(n)$ which is the union of all the $n-1$-dimensional
faces except the $k$-th one. 

If $\pX\in \mK$ is a simplicial set, then the universal property of the representable
$R(n)$ says that $\pX_n=\Hom _{\mK}(R(n),\pX)$. Kan's {\em horn-filling condition} says
that any map $\Lambda (n,k)\rightarrt \pX$ extends to a map $R(n)\rightarrt \pX$. 
The simplicial sets $\pX$ satisfying this condition are the fibrant objects of 
the model structure on $\mK$. 

Boardman and Vogt introduced the {\em restricted Kan condition}, satisfied by a simplicial
set $X$ whenever any map $\Lambda (n,k)\rightarrt \pX$ extends to $R(n)\rightarrt \pX$,
for each $0<k<n$. In other words, they consider only the horns obtained by taking out
any except for the first and last faces. 

This condition corresponds to keeping a directionality of the $1$-cells in $\pX$.
This may be seen most clearly by looking at the case $n=2$. A $2$-cell may be drawn as
$$
{\setlength{\unitlength}{.5mm}
\begin{picture}(100,40)

\put(-25,17){\ensuremath{R(2):}}

\put(10,10){\circle*{2}}
\put(50,30){\circle*{2}}
\put(90,10){\circle*{2}}

\qbezier(10,10)(30,20)(50,30)
\qbezier(50,30)(70,20)(90,10)
\put(10,10){\line(1,0){80}}

\put(29,23){\ensuremath{h}}
\put(69,23){\ensuremath{g}}
\put(49,2){\ensuremath{f}}
\end{picture}
} 
$$
where $h$, $g$ and $f$ are the $1$-cells corresponding to
edges $(01)$, $(12)$ and $(02)$ respectively.
Such a $2$-cell is thought of as the relation $f=gh$. 
In the usual Kan condition, there are three horns which need to be filled:
$$
{\setlength{\unitlength}{.5mm}
\begin{picture}(100,40)

\put(-25,17){\ensuremath{\Lambda (2,0):}}

\put(10,10){\circle*{2}}
\put(50,30){\circle*{2}}
\put(90,10){\circle*{2}}

\qbezier(10,10)(30,20)(50,30)

\put(10,10){\line(1,0){80}}

\put(29,23){\ensuremath{h}}
\put(65,20){\ensuremath{?}}
\put(49,2){\ensuremath{f}}
\end{picture}
} 
$$
$$
{\setlength{\unitlength}{.5mm}
\begin{picture}(100,40)

\put(-25,17){\ensuremath{\Lambda (2,1):}}

\put(10,10){\circle*{2}}
\put(50,30){\circle*{2}}
\put(90,10){\circle*{2}}

\qbezier(10,10)(30,20)(50,30)
\qbezier(50,30)(70,20)(90,10)

\put(29,23){\ensuremath{h}}
\put(69,23){\ensuremath{g}}
\put(49,9){\ensuremath{?}}
\end{picture}
} 
$$
$$
{\setlength{\unitlength}{.5mm}
\begin{picture}(100,40)

\put(-25,17){\ensuremath{\Lambda (2,2):}}

\put(10,10){\circle*{2}}
\put(50,30){\circle*{2}}
\put(90,10){\circle*{2}}

\qbezier(50,30)(70,20)(90,10)
\put(10,10){\line(1,0){80}}

\put(33,20){\ensuremath{?}}
\put(69,23){\ensuremath{g}}
\put(49,2){\ensuremath{f}}
\end{picture}
} 
$$
However, in the restricted Kan condition, only the middle horn $\Lambda (2,1)$
is required to be filled. This corresponds to saying that for any composable
arrows $g$ and $h$, there is a composition $f=gh$.
On the other hand, filling the horn $\Lambda (2,0)$ would correspond to saying that
given $f$ and $h$, there is $g$ such that $f=gh$, which essentially means we look
for $g=fh^{-1}$; and filling $\Lambda (2,2)$ would correspond to saying that
given $f$ and $g$ there is $h$ such that $f=gh$, that is $h=g^{-1}f$.

When we look at things in this way, it is clear that the full Kan condition corresponds
to imposing, in addition to the categorical composition of arrows, some kind of groupoid
condition of existence of inverses. It isn't surprising, then, that Kan complexes 
correspond to $\infty$-groupoids. 

Following through this philosophy has led Joyal to the theory of quasicategories,
which are simplicial sets satisfying the restricted Kan condition, but viewed as
$(\infty ,1)$-categories with arrows which are not necessarily invertible. 

Making the translation from restricted Kan 
simplicial sets to $(\infty ,1)$-categories is not altogether
trivial, most notably for any two vertices $x,y$ of a quasicategory $\pX$ we need
to define the {\em simplicial mapping space} $\pX (x,y)$; one possibility
is to say that it is the Kan simplicial set
$$
k\mapsto \Hom ^{x,y}(R(1)\times R(k), \pX )
$$
where the superscript indicates maps sending $0\times R(k)$ to $x$ and $1\times R(k)$ to
$y$.  There is also a way of describing directly a simplicial category which is the
rectification of the corresponding $(\infty , 1)$-category; see \cite{Riehl}
for a detailed discussion. 

Joyal constructs a model category structure whose underlying category
is that of simplicial sets, in for which the fibrant objects are exactly those satisfying
the restricted Kan condition. The passage from a general simplicial set to its
fibrant replacement, done by enforcing the restricted Kan horn filling conditions using
the small object argument, is a version of the ``calculus of generators and relations''
very similar to what we will be discussing in Chapters \ref{genrel1} and \ref{secat1} 
for the
case of Segal categories. 

In the three basic kinds of simplicial objects which we now have representing
$(\infty ,1)$-categories with weak composition, we can see a trade-off between 
information content and simplicity. The simplest model is that of quasicategories,
which are just simplicial sets satisfying a very classical horn-filling condition;
but in this case it isn't easy to get back some of the main pieces of information in
an $(\infty ,1)$-category such as the simplicial mapping sets. At the other end,
in Rezk's complete Segal spaces, the full information of the $\infty$-groupoid
interior is contained within the object, to the extent that the homotopy type of the
$\Delta ^o$-diagram is an invariant of the $(\infty ,1)$-category up to equivalence;
on the other hand, the initial steps of the theory are more complicated. The theory of
Segal categories fits in between: a Segal category has more information readily
at hand than a quasicategory, but less than a complete Segal space; and the initial
theory is more complicated than for quasicategories but less than for complete Segal spaces.

\section{Going between Segal categories and $n$-categories}
\label{sec-goingbetween}

We mention briefly the relationship between the notions of Segal category and
$n$-category. 
Tamsamani's definition of $n$-category is recursive. The basic idea is to use
the same definition as above for Segal category, but where the $\pA _{p/}$ are
themselves $n-1$-categories. The appropriate condition on the Segal maps is
the condition of equivalence of $n-1$-categories, which in turn is defined
(inductively) in the same way as the notion of equivalence of Segal categories
explained above.

Tamsamani shows
that the homotopy category of $n$-groupoids is
the same as that of $n$-truncated spaces. The two relevant functors are the
realization and Poincar\'e $n$-groupoid $\Pi _n$ functors. Applying this to
the $n-1$-categories $\pA _{p/}$ we obtain the following relationship. An
$n$-category $\pA $ is said to be {\em $1$-groupic} (notation introduced in
\cite{limits}) if the $\pA _{p/}$ are  $n-1$-groupoids. In this case, replacing
the $\pA _{p/}$ by their realizations $|\pA _{p/}|$ we obtain a simplicial space
which satisfies the Segal condition. Conversely if $\pA _{p/}$ are spaces or
simplicial sets then replacing them by their $\Pi _{n-1}(\pA _{p/})$ we obtain
a simplicial collection of $n-1$-categories, again satisfying the Segal
condition.  These constructions are not quite inverses because
$$
| \Pi _{n-1}(\pA _{p/})| = \tau _{\leq n-1} (\pA _{p/})
$$
is the Postnikov truncation. If we think (heuristically) of setting
$n=\infty$ then we get inverse constructions. Thus---in a sense which I will
not currently make more precise than the above discussion---one can say that
Segal categories are the same thing as $1$-groupic $\infty$-categories.

The passage from simplicial sets to Segal categories is the same as the
inductive passage from $n-1$-categories to $n$-categories. In \cite{svk}
was introduced the notion of {\em $n$-precat}, the analogue of the above
Segal precat.  Noticing that the results and arguments in \cite{svk} are
basically organized into one gigantic inductive step passing from $n-1$-precats
to $n$-precats, the same step applied only once works to give the analogous
results in the passage from simplicial sets to Segal precats.

The notion of Segal category thus presents, from a technical point of view,
an aspect of a ``baby'' version of the notion of $n$-category.  On the other
hand, it allows a first introduction of homotopy going all the way up to
$\infty$ (i.e. it allows us to avoid the $n$-truncation inherent in the
notion of $n$-category).

One can easily imagine combining the two into a notion of ``Segal
$n$-category'' which would be an $n$-simplicial simplicial set satisfying the
globular condition at each stage. It is interesting and historically
important to note that the notion of Segal $n$-category with only one
$i$-morphism for each $i\leq n$, is the same thing as the notion of {\em
$n$-fold delooping machine}.
This translation comes out of Dunn \cite{Dunn}, which apparently dates
essentially back to 1984.  In retrospect it is not too hard to see
how to go from Dunn's notion of $E_n$-machine, to Tamsamani's notion of
$n$-category, simply by relaxing the conditions of having only one object.
Metaphorically, $n$-fold delooping machines correspond to the
Whitehead tower, whereas $n$-groupoids correspond to the Postnikov tower.

There are other proposals for simplicial models for $n$-categories which we haven't been able to discuss.
For example, 
Street proposed a model based on simplicial sets with certain distinguished
simplicial subsets which he calls ``thin subcomplexes''.

\section{Towards weak $\infty$-categories}
\label{sec-towards}

We mention here some ideas for going towards a theory of $\infty$-categories. 
As the iterative approach makes clear, there is no direct generalization of our theory to the case $n=\infty$
(which means the first infinite ordinal). The notion of equivalence of Segal categories, crucial to everything,
is defined by a top-down induction, so by its nature it is related to some kind of $n$-categories.
As our procedure makes clear, and as came out in \cite{descente} and \cite{Pelissier}, this iteration
can start with any model category such as $\mK$, allowing us to define Segal $n$-categories which in Lurie's notation
correspond to $(\infty , n)$-categories i.e. $\infty$-categories where all morphisms are invertible starting from degree $n+1$.

Cheng's argument \cite{ChengInfty} shows that if $\pA $ is to be an $\infty$-category with duals at all levels
then, in an algebraic sense, all morphisms look invertible. However, it is clear from Baez-Dolan's predictions about the theory,
that we don't want to identify $\infty$-categories with duals and $\infty$-groupoids, indeed they represent somehow complementary
points of view, the first being related to quantum field theory and the second to topology. 
Hermida, Makkai and Power have discussed these issues in \cite{HermidaMakkaiPower}.

From these observations and thinking about specific kinds of examples, the following idea emerges: in a true $\infty$-category $\pA $,
the information about which $i$-morphisms should be considered as invertible, should be thought of as an {\em additional structure}
beyond the algebraic structure of some kind of weak multiplication operations. Going back to the strict case, one can 
well imagine a strictly
associative $\infty$-category $\pA $ in which any $i$-morphism $u$ has a morphism going in the other direction $v$
and $i+1$-morphisms $uv\rightarrt 1$ and $vu\rightarrt 1$. Then, we could either declare all morphisms to be invertible,
in which case $v$ would be the inverse of $u$ since the $i+1$-morphisms going to the identities would be invertible; or
alternatively we could declare that no morphisms (other than the identities) are invertible. 

Both choices would be reasonable, and would lead to {\em different} $\infty$-categories sharing the same underlying algebraic
structure $\pA $. 

So, if we imagine a theory of weak $\infty$-categories in which the information of which morphisms are invertible is somehow present
then it becomes reasonable to define the {\em truncation operations} $\tau _{\leq n}$ as in \cite{Tamsamani} but going 
from weak $\infty$-categories to weak $n$-categories. Thus an $\infty$-category $\pA $ would lead to a compatible system of $n$-categories 
$\tau _{\leq n}(\pA )$ for all $n$. 

This suggests a definition: 
it appears that one should get the right theory by taking a homotopy inverse limit of the theories of $n$-categories. Jacob Lurie had mentioned something like this in correspondence
some time ago. Given the compatibility of the Rezk-Barwick theory with homotopy
limits \cite{BergnerHoFiProd}, that might be specially adapted to this task. 

One might alternatively be able to view the theory of $\infty$-categories as
some kind of first ``fixed point'' of the operation $\mM \rightarrt \precat (\mM )$ which we will construct in the
main chapters.

We will leave these considerations on a speculative level for now, hoping only that the
techniques to be developed in the main part of the book will be useful in attacking
the problem of $\infty$-categories later.


\chapter{Operadic approaches}
\label{operadic1}

Apart from the simplicial approaches, the other main direction is comprised of a number of {\em operadic approches}, definitions
of higher categories based on Peter May's notion of ``operad''. This dichotomy is not surprising, given that operads and simplicial objects
are the two main ways of doing delooping machines in algebraic topology. The operadic approaches are not the main subject of this book, so
our presentation will be more succinct designed to inform the reader of what is out there. 

Tom Leinster's book \cite{LeinsterBook} has a very complete
discussion of the relationship between operads and higher categories. His paper 
\cite{Leinster} gives a brief but detailed exposition of numerous different definitions
of higher categories, including several in the operadic direction, and that was the first
appearence in print of some definitions such as
Trimble's for example. 

\section{May's delooping machine}

We start by recalling Peter May's delooping machine. An {\em operad} is a collection of sets $O(n)$, thought of as the ``set of $n$-ary operations'',
together with some maps which are thought of as the result of substitutions:
$$
\psi : O(m) \times O(k_1) \times \cdots \times O(k_m) \rightarrt O(k_1+\cdots + k_m)
$$
for any uplets of integers $(m; k_1,\ldots , k_m)$. If we think of an element $u\in O(n)$ as representing a function $u(x_1,\ldots , x_n)$ then
$\psi (u; v_1,\ldots , v_m)$ is the function of $k_1+\ldots + k_m$ variables
$$
{\scriptstyle
u(v_1(x_1,\ldots , x_{k_1}),v_2(x_{k_1+1},\ldots , x_{k_{k_1+k_2}}), \ldots , v_m(x_{k_1+\cdots + k_{m-1}+1},\ldots , x_{k_1+\cdots + k_m})).
}
$$
The substitution operation $\psi$ is required to satisfy the appropriate axioms \cite{MayLoops}. 

A {\em topological operad} is the same, but where the $O(n)$ are topological spaces; it is said to be {\em contractible} if each $O(n)$ is contractible. More generally,
we can consider the notion of operad in any category admitting finite products. 

There is a notion of {\em action of an operad $O$ on a set $X$}, which means an association to each $u\in O(n)$ of an actual $n$-ary function 
$X\times \cdots \times X\rightarrt X$ such that the substitution functions $\psi$ map to function substitution as described above.

If $X$ is a space and $O$ is a topological operad, we can require that the action consist of continuous functions 
$O(n)\times X^n\rightarrt X$, or more generally if $O$ is an operad in $\mM$ then
an action of $O$ on $X\in \mM$ is a collection of morphisms $O(n)\times X^n\rightarrt X$
satisfying the appropriate compatibility conditions. 

A {\em delooping structure} on a space $X$, is an action of a contractible operad on it. 

The typical example is that of the {\em little intervals operad}: here $O(n)$ is the space of inclusions of $n$ consecutive intervals into the
given interval $[0,1]$, and substitution is given by pasting in. 
A variant is used in Section \ref{sec-trimble} below.

\section{Baez-Dolan's definition}

In Baez and Dolan's approach, the notion of operad is first and foremost used to
determine the shapes of higher-dimensional cells. They introduce a category of {\em opetopes}
and a notion of {\em opetopic set} which, like the case of simplicial sets, just
means a presheaf on the category of opetopes. They then impose filler conditions.
Their scheme of filler conditions is inductive on the dimension of the opetopes, but
is rather intricate. We describe the category of opetopes 
by drawing some of the standard pictures, and then give an informal
discussion of the filler conditions. In addition to  the original
papers \cite{BaezDolan}
\cite{BaezDolanIII},
Leinster's \cite{Leinster} was one of
our main sources. Readers may also consult \cite{ChengOpetopes} \cite{KockJoyalBataninMascari} 
for other approaches
to defining and calculating with opetopes. 

An $n$-dimensional 
opetope should be thought of as a roughly
globular $n$-dimensional object, with an output face which is an 
$n-1$-dimensional opetope, and an input face which is a {\em pasting diagram}
of $n-1$-dimensional opetopes. To paste opetopes together, match up the output
faces with the different pieces of the input faces. 

The only $0$-dimensional opetope is a point. The only $1$-dimensional  opetope
has as input and output a single point, so it is a single arrow 
$$
{\setlength{\unitlength}{.5mm}
\begin{picture}(60,20)
\put(10,10){\circle*{2}}
\put(50,10){\circle*{2}}
\put(12,10){\line(1,0){36}}
\put(48,10){\vector(1,0){0}}
\end{picture}
} 
$$
A pasting diagram of $1$-dimensional opetopes can therefore be composed
of several arrows joined head to tail:
$$
{\setlength{\unitlength}{.5mm}
\begin{picture}(110,30)
\put(10,10){\circle*{2}}
\put(40,25){\circle*{2}}
\put(70,25){\circle*{2}}
\put(100,10){\circle*{2}}
\put(12,11){\line(2,1){26}}
\put(42,25){\line(1,0){26}}
\put(72,24){\line(2,-1){26}}
\put(38,24){\vector(2,1){0}}
\put(68,25){\vector(1,-0){0}}
\put(98,11){\vector(2,-1){0}}
\end{picture}
} 
$$
A $2$-dimensional opetope can then have such a pasting diagram as input,
with a single $1$-dimensional opetope or arrow as output:
$$
{\setlength{\unitlength}{.5mm}
\begin{picture}(110,50)
\put(10,10){\circle*{2}}
\put(40,40){\circle*{2}}
\put(70,40){\circle*{2}}
\put(100,10){\circle*{2}}
\put(53,23){\ensuremath{\Downarrow}}

\qbezier(11,12)(25,32)(38,39)
\qbezier(42,41)(55,45)(68,41)
\qbezier(72,39)(85,32)(99,12)
\qbezier(12,9)(55,4)(98,9)
\put(38,39){\vector(3,2){0}}
\put(68,41){\vector(4,-1){0}}
\put(99,12){\vector(2,-3){0}}
\put(98,9){\vector(4,1){0}}
\end{picture}
} 
$$
A pasting diagram of $2$-dimensional opetopes would arise if we add
on some other opetopes, with the output edges attached to the 
input edges. This can be done recursively. A picture where we add
three more opetopes on the three inputs, with $2$, $1$ and $3$ input
edges respectively; then a fourth one on the 
second input edge of the first new one, would look like this:
$$
{\setlength{\unitlength}{.5mm}
\begin{picture}(110,50)
\put(10,10){\circle*{2}}

\put(12,40){\circle*{2}}
\put(12,39){\vector(1,4){0}}

\put(25,50){\circle*{2}}
\put(24,50){\vector(3,2){0}}

\put(40,40){\circle*{2}}
\put(38,39){\vector(3,2){0}}
\put(39,40){\vector(1,0){0}}
\put(39,42){\vector(3,-2){0}}

\put(70,40){\circle*{2}}
\put(68,41){\vector(1,-1){0}}
\put(68,40){\vector(1,0){0}}

\put(88,41){\circle*{2}}
\put(87,41){\vector(4,-1){0}}

\put(100,27){\circle*{2}}
\put(100,28){\vector(1,-2){0}}

\put(100,10){\circle*{2}}
\put(101,11){\vector(-1,-4){0}}
\put(99,12){\vector(2,-3){0}}
\put(98,9){\vector(4,1){0}}

\qbezier(10,12)(9,25)(12,39)
\qbezier(13,40)(27,42)(39,40)
\qbezier(13,41)(19,48)(24,50)
\qbezier(26,50)(32,48)(39,41)

\qbezier(41,41)(55,55)(68,41)

\qbezier(70,40)(77,43)(88,41)
\qbezier(88,41)(97,35)(100,27)
\qbezier(100,27)(103,18)(100,10)

\put(53,23){\ensuremath{\Downarrow}}
\put(53,42){\ensuremath{\Downarrow}}
\put(23,42){\ensuremath{\Downarrow}}

\put(18,30){\ensuremath{\Downarrow}}
\put(87,29){\ensuremath{\Downarrow}}

\qbezier(11,12)(25,28)(38,39)
\qbezier(41,40)(55,38)(68,40)
\qbezier(72,39)(85,28)(99,12)
\qbezier(12,9)(55,4)(98,9)

\end{picture}
} 
$$

Now, a $3$-dimensional opetope could have the above pasting diagram as input;
in that case, the output opetope is supposed to have the shape of the boundary
of the pasting diagram:
$$
{\setlength{\unitlength}{.5mm}
\begin{picture}(110,50)
\put(10,10){\circle*{2}}

\put(12,40){\circle*{2}}
\put(12,39){\vector(1,4){0}}

\put(25,50){\circle*{2}}
\put(24,50){\vector(3,2){0}}

\put(40,40){\circle*{2}}
\put(39,42){\vector(3,-2){0}}

\put(70,40){\circle*{2}}
\put(68,41){\vector(1,-1){0}}

\put(88,41){\circle*{2}}
\put(87,41){\vector(4,-1){0}}

\put(100,27){\circle*{2}}
\put(100,28){\vector(1,-2){0}}

\put(100,10){\circle*{2}}
\put(101,11){\vector(-1,-4){0}}

\qbezier(10,12)(9,25)(12,39)
\qbezier(13,41)(19,48)(24,50)
\qbezier(26,50)(32,48)(39,41)

\qbezier(41,41)(55,55)(68,41)

\qbezier(70,40)(77,43)(88,41)
\qbezier(88,41)(97,35)(100,27)
\qbezier(100,27)(103,18)(100,10)

\put(50,24){\ensuremath{\Downarrow}}

\qbezier(12,9)(55,4)(98,9)

\end{picture}
} 
$$
That is to say, it is a $2$-dimensional opetope with $7$ input edges. 
We don't draw the $3$-dimensional opetope; it should look like a ``cushion''. 

Baez and Dolan define in this way the {\em category of opetopes} $\Opetopes $. 
The reader is refered to the main references \cite{BaezDolan}
\cite{BaezDolanIII} \cite{HermidaMakkaiPower} \cite{ChengOpetopes} \cite{KockJoyalBataninMascari} 
\cite{Leinster} for the precise definitions. 
A {\em Baez-Dolan $n$-category} is an {\em opetopic set}, that is a functor 
$$
\pA : \Opetopes \rightarrt \Sets ,
$$
which is required to satisfy some conditions. Refer again to 
\cite{BaezDolan}
\cite{BaezDolanIII} \cite{HermidaMakkaiPower} \cite{ChengOpetopes} \cite{KockJoyalBataninMascari} 
\cite{Leinster} for the precise statements of these conditions. 
see also \cite{Zawadowski} for one of the most recent treatments. 

Included
in the conditions are things having to do with truncation so that the cells of dimension $>n+1$ don't matter, when speaking of an $n$-category. 

Beyond this truncation, 
one of the main features of the Baez-Dolan viewpoint is that an $n$-dimensional opetope
represents an $n$-morphism which is not necessarily invertible, from the ``composition''
of the opetopes in the input face, to that of the output face. In particular, the output
face of an opetope is not necessarily the same as the composition of the input faces; this
distinguishes their setup from the Segal-style picture we are mostly
considering in the present book, in which a simplex represents a composition with
the outer edge being equivalent to the composition of the principal edges. 

The idea that opetopes represent arbitrary 
morphisms from the composition of the inputs, to the output,
makes it somewhat complicated to collect together and write down all of the appropriate
conditions which an opetopic set should satisfy in order to be an $n$-category. 
The main step is to designate certain cells as ``universal'', meaning that the output
edge is equivalent to the composition of the input edges, which is basically a homotopic
initiality property. The necessary uniqueness has to be treated up to equivalence,
whence the need for a definition of equivalence in the induction. All of this leads to 
Baez and Dolan's notions of {\em niches}, {\em balanced cells} and so forth. 

As a rough approximation, we can say that these conditions are a form of horn-filler conditions, generalizing the restricted Kan condition
used in the definition of quasicategory \cite{BoardmanVogt} \cite{JoyalQC}
but adapted to the opetopic context.

To close out this section, here are a few thoughts on the possible relationship between 
this theory and the many other theories of $n$-categories. For one thing, the fact that 
the opetopic cells represent explicitly the $n$-morphisms which are not necessarily invertible, would seem to render this theory particularly well adapted to looking at things like lax functors (what Benabou calls ``morphisms'' as opposed to ``homomorphisms'' \cite{Benabou}). 

On the other hand, for a comparison with other theories, it would
be interesting to investigate functors ${\bf F}:\Opetopes \rightarrt \mP$ where
$\mP$ is a closed model category of ``$n$-categories'' (for example, the $\mP = \precat ^n(\Sets )$ which we are going to be constructing in the rest of the book). 
Given such a functor, if $\pB \in \mP$ is a fibrant object, we would obtain an opetopic set 
${\bf F}^{\ast}(\pB )$, and conversely given an opetopic set $\pA$ we could construct
its realization ${\bf F}_! (\pA )$ in $\mP$. Under the right
hypotheses, one hopes, 
these should set up a Quillen adjunction between a model category of opetopic
sets, and the other model category $\mP$. 

If both of the above remarks could be realized, it would lead to the introduction of
powerful new techniques for treating lax functors
in any of the other theories of $n$-categories.

\section{Batanin's definition}
\label{sec-batanin}

Batanin's definition is certainly the closest to Grothendieck's original vision.
Recall the passage that we have quoted from ``Pursuing Stacks'' on p. \pageref{grothendieckgroupoidquote} in Chapter \ref{why1} above, saying that whenever two
morphisms which are naturally obtained as some kind of composition, have the
same source and target, there should be a homotopy between them at one level up. 
Batanin's definition puts this into place, by carefully studying the notion of possible
composition of arrows in a higher category. It was recently pointed out by Maltsiniotis
that Grothendieck had in fact given a definition of higher groupoid \cite{MaltsiniotisGroGpd}
and that a small modification of that approach yields a definition of higher category
which is very close to Batanin's \cite{MaltsiniotisGroBat}. 

Tom Leinster and Eugenia Cheng have refined Batanin's original work, and our discussion
will be informed by their expositions \cite{Leinster} \cite{LeinsterBook} \cite{ChengComparison},
to which the reader should refer for more details. Leinster has also introduced 
some related definitions of weak $n$-category based on the notion of multicategory,
again for this we refer to \cite{LeinsterBook}. 

One of Batanin's innovations was to introduce a notion of operad adapted to higher
categories, based on the notion of {\em globular set}, a 
presheaf on the category $\Glob$ which has objects $\globe _i$ for each $i$,
and maps $s^{!},t^{!}:\globe _i\rightarrt \globe _{i+1}$. The objects $\globe _i$ are
supposed to represent globular pictures of $i$-morphisms, for example 
$\globe _2$ may be pictured as 
$$
{\setlength{\unitlength}{.5mm}
\begin{picture}(60,40)
\put(10,20){\circle*{2}}
\put(50,20){\circle*{2}}
\qbezier(11,21)(30,40)(47,23)
\qbezier(11,19)(30,0)(47,17)
\put(49,21){\vector(1,-1){0}}
\put(49,19){\vector(1,1){0}}
\put(29,18){\ensuremath{\Downarrow}}
\end{picture}
} 
$$
and the maps $s^{!},t^{!}$ are viewed as inclusions of smaller-dimensional globules on
the boundary. Thus, the inclusion maps are subject to the relations 
$$
s^{!}s^{!} = t^{!}s^{!},\;\;\; t^{!}s^{!} = t^{!}t^{!}.
$$
Dually, a globular set is a collection of sets $\pA _i:= \pA (\glob _i)$ 
with source and target maps
$$
\begin{diagram}
\pA _i &\pile{\rightarr^s \\ \rightarr_t} & \pA _{i-1} \cdots \pA _1 &
\pile{\rightarr^s \\ \rightarr_t} & \pA _0
\end{diagram}
$$
subject to the relations $s\circ s = s\circ t$ and $s\circ t= t\circ t$. 

This definition differs slightly from the one which was suggested in Chapter \ref{strict1};
there we looked at what should probably be called a {\em unital globular set} having
identity maps going back in the other direction $i:\pA _i\rightarrt \pA _{i+1}$. 
For the present purposes, we consider globular sets without identity maps. 

Batanin's idea is to use the notion of globular set to generate the appropriate kind of
``collection'' for a notion of higher operad. One can think of an operad as specifying
a collection of operations of each possible arity; and an arity is a possible configuration of the collection
of inputs. For Batanin's globular operads, an input configuration is represented by a 
{\em globular pasting diagram} $P$, and the family of operations of arity $P$ should
itself form a globular set. A globular pasting diagram in dimension $n$ is an $n$-cell
in the free strict $\infty$-category generated by a single non-identity cell in each dimension.

Batanin gives an explicit description of the $n$-cells in this free $\infty$-category,
using planar trees. The nodes of the trees come from the source and target maps
between stages in a globular set, which is not to be confused with the occurence of trees
in parametrizations of elements of a free nonassociative algebra (where the nodes
correspond to parenthetizations)---for the free $\infty$-category involved here, it
is strictly associative and going up a level in the tree corresponds rather to 
going from $i$-morphisms to $i+1$-morphisms. 

As Cheng describes in a detailed example \cite{ChengComparison}, a globular pasting diagram 
is really just a picture of how one might compose together various $i$-morphisms.
Given any globular set $\pA$, and a globular pasting diagram $P$, we get a 
set denoted $\pA _P$ consisting of all the ways of filling in $P$ with labels from
the globular set $\pA$, consistent with sources and targets. 

Now, a globular operad $\mB$ consists of specifying,
for each globular pasting diagram $P$, 
a globular set $\mB (P)$ of operations of arity $P$. This globular set consists
in particular of sets $\mB (P)_n$ for each $n$ calle the {\em set $P$-ary operations 
of level $n$}, together with sources and targets
$$
\begin{diagram}
\mB (P)_{n+1} &\pile{\rightarr ^s \\ \rightarr _t} & \mB (P)_{n}
\end{diagram}
$$
satisfying the globularity relations. An {\em action} of $\mB$ on a globular set $\pA$
consists of specifying for each globular pasting diagram $P$ of degree $n$, a map 
of sets
$$
\mB (P)_n\times \pA _P\rightarrt \pA _n
$$
such that the source diagram 
$$
\begin{diagram}
\mB (P)_{n}\times \pA _P& \rightarr & \pA _{n} \\
\downarr _s & & \downarr ^s \\
\mB (sP)_{n-1}\times \pA _{sP}& \rightarr & \pA _{n-1}
\end{diagram}
$$
commutes and similarly for the target diagram. Here $sP$ and $tP$ are the source
and target of the pasting diagram $P$, which are pasting diagrams of degree $n-1$.
Of course we need to describe an 
additional operadic structure on $\mB$  and the action should be compatible
with this too, but it is easier
to first consider what data an  algebra should have.

To describe what an operadic structure should mean, notice that there is
an operation of substitution of globular pasting diagrams, indeed the free
$\infty$-category on one cell in each dimension can be denoted $\gpasting$,
and if $P\in \gpasting _m$ is a globular pasting diagram in degree $m$
then we get a map of globular sets 
$$
\gpasting _P \rightarrt \gpasting ,
$$
that is to say that given a labeling of the cells of $P$ where the labels $L_j$ 
are themselves
globular pasting diagrams, we can substitute the labels into the cells and obtain
a big resulting globular pasting diagram. 

Now the operadic structure, which is in addition to the structure of $\mB$ described
above, should say that for any element of $\gpasting _P$ consisting of labels
denoted $L_j\in \gpasting _{n_j}$ for the cells of $P$ and yielding an output pasting 
diagram $S$, if we are given
elements of the $\mB (L_j)$ plus an element of $\mB (P)$ then there should be
a big output element of $\mB (S)$.  This needs to be compatible with source and target
operations, as well as compatible with iteration in the style of usual operads. 

Batanin describes explicitly the combinatorics of this using the identification
between pasting diagrams and trees, whereas Leinster takes a more abstract approach 
using monads and multicategories. We refer the reader to the references for 
further details. 

The main point is that this discussion establishes a language in which to say that
the system of coherencies should satisfy a globular contractibility condition. 
Observe that contractibility is a very easily defined property of a globular set, and
it doesn't depend on any kind of composition law: a globular set $\pA$ is contractible
if it is nonempty, and if for any two $f,g\in \pA _m$ with $s(f)=s(g)$ and $t(f)=t(g)$,
there exists an $h\in \pA _{m+1}$ with $s(h)=f$ and $t(h)=g$. 

Note that this definition
is really only appropriate if we are working with globular sets which are truncated
at some level (i.e. trivial above a certain degree $n$), or ones which are supposed to represent
$(\infty , n)$-categories, that is ones in which all $i$-morphisms are declared
invertible for $i>n$. For the purposes of Batanin's definition of $n$-categories,
this is the case. 

Now, a Batanin $n$-category is a globular set $\pA$ provided with an action of a
globular operad $\mB$, such that $\mB$ is contractible. Batanin constructs a {\em universal}
globular operad, but it can also be convenient to work with other contractible
globular operads, as in Cheng's comparison with Trimble's work which we discuss next. 

The elements of $\mB (P)_n$ are the ``natural operations taking a collection of 
morphisms of various degrees, and combining together to get an $n$-morphism''.
Contractibility of $\mB$ really puts into effect Grothendieck's dictum that, given
two natural operations $f$ and $g$ with the same source and target, there should
be an operation $h$ one level higher whose source is $f$ and whose target is $g$. 
So, Batanin's definition is the closest to what Grothendieck was asking for. 

\section{Trimble's definition and Cheng's comparison}
\label{sec-trimble}

Trimble's definition of higher category has acquired a central role because of 
Cheng's recent work comparing it to Batanin's definition \cite{ChengComparison}.
This has also been taken up by Batanin, Cisinski and Weber in \cite{BataninCisinskiWeber}.
Trimble's framework has the advantage that it is iterative. In the future 
it should be possible to establish a comparison with the iterative Segal approach
we are discussing in the rest of the book. Such a comparison result would be very interesting,
but we don't discuss it here. Instead we just give the basic outlines of Trimble's
approach and state Cheng's comparison theorem. We are following very closely her article 
\cite{ChengComparison}.

Consider the following operad $O^{T}$ in $\Top$: 
$$
O^{T}(n) \subset C^0([0,1], [0,n])
$$
is the subset of endpoint-preserving continuous maps. 
The operad structure $\psi ^T$ is given by 
$$
\psi (f; g_1,\ldots , g_m) (t):= g_j(f(t)-(j-1)), \;\;\;\; f(t)\in [j-1,j] \subset [0,m].
$$
The spaces $O^{T}(n)$ are contractible. 
This topological  operad is particularly adapted to loop spaces and path spaces. 
If $Z\in \Top$, and $x,y\in Z$ let ${\rm Path}^{x,y}(Z)$ denote the space of paths
$\gamma : [0,1]\rightarrt Z$ with $\gamma (0)=x$ and $\gamma (1)=y$. 
For any sequence of points $x_0,\ldots  , x_m\in X$ we have a ``substitution'' map
$$
O^{T}(m) \times {\rm Path}^{x_0,x_1}(Z)\times \cdots \times {\rm Path}^{x_{n-1},x_n}(Z),
$$
and these are compatible with the operad structure. In particular when the points are all
the same, this gives an action of $O^{T}$ on the loop space $\Omega ^x(Z)$. 

A {\em Trimble topological category} consists of a set $X$ ``of objects'', together
with a collection of spaces $\pA (x,y)$ for any $x,y\in X$, and collection of maps
$$
\phi _{x_{\cdot}}: O^{T}(n)\times \pA (x_0,x_1)\times \cdots \times \pA (x_{n-1},x_n)
\rightarrt \pA (x_0,x_n)
$$
for any sequence $(x_0,\ldots , x_n)$ in $X$. These should satisfy a compatibility condition
with the operad structure $\psi^T$ for $O^{T}$: given a sequence $x_0,\ldots , x_m$
and sequences $y^i_0,\ldots , y^i_{k_i}$ with $y^i_0=x_{i-1}$ and $y^i_{k_i}=x_i$,
$$
\phi (\psi ^T(f;g_1,\ldots , g_m); u^1_1,\ldots , u^1_{k_1},\ldots , u^m_1,\ldots , u^m_{k_m})
$$
$$
=
\phi (f; \phi (g_1;u^1_1,\ldots , u^1_{k_1}),\ldots , \phi (g_1m;u^m_1,\ldots , u^m_{k_m}).
$$

Similarly, for any category $\mM$ admitting finite products, and any operad
$(O,\psi )$ in $\mM$, we can define a notion of $(\mM , O,\psi )$-category; this
consists of a set $X$ of objects, together with a collection of $\pA (x,y)\in \mM$
for any $x,y\in X$, and a collection of maps $\phi$ as above satisfying the same compatibility condition. 

Suppose we are given a category $\mM$ with finite products, and a functor $\Pi :\Top \rightarrt
\mM$, then we obtain an operad $\Pi (O^{T})$ in $\mM$. We get the notion of
$(\mM ,\Pi (O^T),\Pi (\psi ^T))$-category. If contractibility makes sense in $\mM$ and
$\Pi (O^T(n))$ is contractible, then this is a generalization due to Cheng \cite{ChengComparison}, of Trimble's notion of higher category enriched in
$\mM$. 

Trimble's original definition, which first appeared publicly in \cite{Leinster}, 
included an inductive construction of the Poincar\'e $n$-groupoid functor $\Pi _n$. He defines inductively a sequence of categories which 
Cheng denotes by $\mV _n$, starting with $\mV _0=\Sets$; together with product-compatible
Poincar\'e $n$-groupoid functors $\Pi _n:\Top \rightarrt \mV _n$. The inductive
definition is that $\mV _{n+1}$ is the category of $(\mV _n,\Pi _n(O^T),\Pi _n(\psi ^T))$-categories, and 
$$
\Pi _{n+1}(Z)= (X, \pA , \phi )
$$
where $X:= Z^{\disc}$ is the discrete set of points of $Z$, where $\pA$ is defined by
$\pA (x,y):= \Pi _n({\rm Path}^{x,y}(Z))$, and $\phi$ is defined using the action described above
of $O^T$ on the 
path spaces. See \cite{ChengComparison}  for further details,
as well as for the generalization to the case where an arbitrary contractible operad $P_n$
replaces $\Pi _n(O^T)$. 

Cheng goes on to compare this family of definitions, with Batanin's definition: she shows how
to combine the $P_0,P_1,\ldots , P_{n-1}$ together to form a contractible
globular operad $Q^{(n)}$ such that the category of globular algebras of $Q^{(n)}$ is 
$\mV _n$. This expresses $\mV _n$ as a category of Batanin $n$-categories for this particular
choice of contractible globular operad. The reader is refered to \cite{ChengComparison}, as
well as to \cite{BergerCellularNerve}, \cite{BataninCisinskiWeber} and
\cite{CisinskiBatanin} for related aspects of this kind of comparison.

\section{Weak units}

In the course of investigating the nonrealization of homotopy $3$-types by strict 
$3$-groupoids (\cite{hty3types}, see Chapter \ref{nonstrict1}), the main obstruction seemed
to be the strict unit condition in a strict $n$-category. This is one of the main aspects
which allows the Eckmann-Hilton argument to work. That was explained to me by
Georges Maltsiniotis and Alain Brugui\`eres, but has of course been well-known for a long time.  This led to the conjecture that maybe it would be sufficient to keep the strict associativity of composition, but to weaken the unit condition. 

Recall from homotopy theory (cf \cite{Lewis}) the yoga that 
it suffices to weaken any one of the
principal structures involved.  Most weak notions of $n$-category involve a
weakening of the associativity, or eventually of the Godement interchange
conditions.

O. Leroy \cite{Leroy} and apparently, independantly, Joyal and
Tierney \cite{JoyalTierney} were the first to do this in the context of $3$-types. See also
Gordon, Power, Street \cite{GordonPowerStreet} and Berger \cite{BergerEckmannHilton} for
weak $3$-categories and $3$-types. Baues \cite{Baues} showed that
$3$-types correspond to {\em quadratic modules} (a generalization of
the notion of crossed complex) \cite{Baues}. Then come the models for weak 
higher categories which we are considering in the rest of the book. 

It seems likely that the arguments of \cite{KapranovVoevodsky} would show that
one could instead weaken the condition of being {\em unital}, that is  having
identities, and keep associativity and Godement.  We give
a proposed definition of what this would mean and then state two conjectures.

This can be motivated by looking at the {\em
Moore loop space} $\Omega ^x_M(X)$ of a space $X$ based at $x\in X$, cited in  \cite{KapranovVoevodsky} as a motivation for their
construction. Recall that $\Omega ^x_M(X)$ is the space of {\em pairs} $(r,
\gamma )$ where  $r$ is a real number $r\geq 0$ and $\gamma = [0,r]\rightarrt
X$ is a path starting and ending at $x$. This has the advantage of being a
strictly associative monoid. On the other side of the coin, the ``length''
function
$$
\ell : \Omega ^x_M(X)\rightarrt \leftbrack 0,\infty )\subset \rr
$$
has a special behavoir over $r=0$. Note that over the open half-line
$(0,\infty )$ the length function $\ell$ is a fibration (even a fiber-space)
with fiber homeomorphic to the usual loop space. However, the fiber over $r=0$
consists of a single  point, the constant path $[0,0]\rightarrt X$ based at
$x$. This additional point (which is the unit element of the monoid
$\Omega ^x_M(X)$) doesn't affect the topology of $\Omega ^x_M$ (at least if $X$
is locally contractible at $x$) because it is glued in as a limit of paths which
are more and more concentrated in a neighborhood of $x$. However, the map
$\ell$ is no longer a fibration over a neighborhood of $r=0$.  This is a bit of
a problem because $\Omega ^x_M$ is not compatible with direct products of the
space $X$; in order to obtain a compatibility one has to take the fiber product
over $\rr$ via the length function:
$$
\Omega ^{(x,y)}_M(X\times Y)= \Omega ^x_M(X) \times _{\rr} \Omega ^y_M(Y),
$$
and the fact that $\ell$ is not a fibration could end up causing a problem in
an attempt to iteratively apply a construction like the Moore loop-space.

Things seem to get better if we restrict to
$$
\Omega ^x_{M'}(X):=\ell ^{-1}((0,\infty ))\subset \Omega ^x_M(X) ,
$$
but this associative  monoid no longer has a strict unit. Even so, the constant
path of any positive length gives a weak unit.

A motivation coming from a different direction was an observation made by
Tamsamani early in the course of his thesis work. He was trying to define a
strict $3$-category $2Cat$ whose objects would be the strict $2$-categories and
whose morphisms would be the weak $2$-functors between $2$-categories (plus
notions of weak natural transformations and $2$-natural transformations).
At some point he came to the conclusion that one could adequately define
$2Cat$ as a strict $3$-category except that he couldn't get strict identities.
Because of this problem we abandonned the idea and looked toward weakly
associative $n$-categories. In retrospect it would be interesting to pursue
Tamsamani's construction of a strict $2Cat$ but with only weak identities.

In \cite{hty3types} was introduced a preliminary
definition of {\em weakly unital strict $n$-category} (called ``snucategory'' there),
including a notion of direct product. The proposed definition went as follows. 
Suppose we know what these are for $n-1$. Then a weakly unital strict $n$-category $\pC$ consists
of a set $\Ob (\pC )$ of objects together with, for every pair of objects $x,y\in \Ob (\pC )$
a weakly unital strict $n-1$-category $ \pC (x,y)$ and composition morphisms
$$
\pC (x,y)\times \pC (y,z) \rightarrt \pC (x,z)
$$
which are strictly associative,  such that a {\em weak unital condition} holds.  
We now explain this condition. An element $e_x\in \pC (x,x)$ is
called a {\em weak identity} if:
\newline
---composition with $e$ induces  equivalences of weakly unital strict $n-1$-categories
$$
\pC (x,y)\rightarrt \pC (x,y) , \;\;\;
\pC (y,x)\rightarrt \pC (y,x);
$$
---and if $e\cdot e$ is equivalent to $e$. It would be best to complete this last condition
to the fuller collection of coherence conditions introduced by Kock \cite{Kock}. 

In order to complete the recursive definition we must define the notion of when
a morphism of weakly unital strict $n$-categories 
is an equivalence, and we must define what it
means for two objects to be equivalent. A morphism is said to be an equivalence
if the induced morphisms on the $\pC (x,y)$ are equivalences of 
weakly unital strict $n-1$-categories and
if it is essentially surjective on objects: each object in the target is
equivalent to the image of an object. It thus remains just to be seen what internal
equivalence of objects means. For this we introduce the {\em truncations}
$\tau _{\leq i}\pC $ of a weakly unital strict $n$-category $\pC $. Again this is done in the same way
as usual: $\tau _{\leq i}\pC $ is the weakly unital strict
$i$-category with the same objects as $\pC $
and whose morphism $i-1$-categories  are the truncations
$$
Hom_{\tau _{\leq i}\pC }(x,y):=\tau _{\leq i-1}\pC (x,y).
$$
This works for $i\geq 1$ by recurrence, and for $i=0$ we define the truncation
to be the set of isomorphism classes in $\tau _{\leq 1}\pC $.  Note that
truncation is compatible with direct product (direct products are defined in
the obvious way) and takes equivalences to equivalences. These statements used
recursively allow us to show that the truncations themselves satisfy the weak
unary condition. Finally, we say that two objects are equivalent if they map
to the same thing in $\tau _{\leq 0}C$.

Proceeding in the same way as in Chapter \ref{strict1}, we can define the
notion of weakly unital strict $n$-groupoid.

\begin{conjecture}
\label{wusagpd}
There are functors $\Pi _n$ and $\Re$ between the categories of
weakly unital strict $n$-groupoids and $n$-truncated spaces (going in the usual directions)
together with adjunction morphisms inducing an equivalence between the
localization of weakly unital strict $n$-groupoids 
by equivalences, and $n$-truncated spaces by
weak equivalences.
\end{conjecture}

Joachim Kock has developed the right definition of
weakly unitary strict $n$-category which takes into account the full collection
of higher coherence relations \cite{Kock} \cite{Kock2} rather than just asking that
$e\sim e\cdot e$; we refer the reader there for
his definition which supersedes the preliminary version described above. 

Joyal and Kock have proven Conjecture \ref{wusagpd} for the case $n=3$ in \cite{JoyalKock}.
For general $n$, one could hope to apply the argument of \cite{KapranovVoevodsky}. 

These results concern the case of groupoids, however we might also expect that
weakly unital strict $n$-categories serve to model all weak $n$-categories: 

\begin{conjecture}
\label{wusacat}
The localization of the category of weakly unital strict $n$-categories by equivalences,
is equivalent to the localizations of the categories of weak $n$-categories
of Tamsamani and/or Baez-Dolan and/or Batanin by equivalences.
\end{conjecture}

While we're discussing the subject of unitality conditions, the following remark 
is in order. 
The role of strict unitality conditions in the interchange or Eckmann-Hilton relations,
and the consequent nonrealization of homotopy types with nontrivial Whitehead bracket,
suggests that we need to take some care about this point in the general argument
which will be developped in Parts III and IV. It turns out that, in order to insure
a good cartesian property, our Segal-style weakly enriched categories should nonetheless
be endowed with strict units in a certain sense. These correspond to the degeneracies 
in the simplicial category $\Delta$, and are important for the Eilenberg-Zilber argument
which yields the cartesian property. They don't correspond to full strict units
in the maximal possible way, because the composition operation will not even be well-defined;
that is why we will be able to  impose the unitality condition in Part III, 
without running up against the problems
identified in Chapter \ref{nonstrict1}.

\section{Other notions}

The theory of $n$-categories is an essentially {\em globular} theory: an $i$-morphism
has a single $i-1$-morphism as source, and a single one as target. This basic shape
can be relaxed in many ways. For example, Leinster and others have investigated
notions of {\em multicategory} where the input is a collection of objects rather
than just one object. This is somewhat related to the opetopic shapes introduced
by Baez and Dolan. 

Another way of relaxing the globular shape is to iterate the internal category construction.
Brown and Loday constructed the first algebraic representation for homotopy $n$-types,
with the notion of {\em $Cat^n$-group}. Let $Cat ^1(Gp)$ denote the category
of internal categories in the category $Gp$, then $Cat ^{n+1}(Gp)$ is the
category of internal categories in $Cat ^{n}(Gp)$. There is a natural realization
functor from $Cat ^{n}(Gp)$ to homotopy $n$-types, and Brown and Loday prove that
all homotopy $n$-types are realized. 

Paoli has recently refined this model to go back in the globular direction, 
by introducing a notion of {\em special $Cat^n$-group} \cite{PaoliAdvances}. If we think of an internal category as being a pair of objects
connected by several morphisms, then an internal $n$-fold category may be seen as a
collection of objects arranged at the vertices of an $n$-cube. The speciality condition requires that certain faces of the cube be contractible. 

The speciality condition is a sort of weak globularity condition. An internal $n$-fold
category is not in and of itself a globular object, because the object of objects may
be nontrivial. A strict globular condition would have the object of objects be a discrete
set; the speciality condition requires only that it be a disjoint union of contractible
$Cat ^{n-1}$-groups. 

Paoli shows that the special $Cat^n$-groups model homotopy $n$-types, and she relates this
model to Tamsamani's model. This provides a semistrictification result saying that 
we can have a strictly associative composition at any  one stage of Tamsamani's model.
An alternative proof of this result for semistrictification at the last stage, may be
obtained using the fact that Segal categories are equivalent to strict simplicial categories.

Paoli's model relates to Lewis's principle cited above \cite{Lewis} in an interesting
way: in a $Cat^n$-group all the structures---associativity, units, inverses, interchange---are strict; the special $Cat^n$-groups weaken instead the globularity condition itself.

Penon's definition \cite{Penon} is completely algebraic in the sense that a weak $n$-category
is an algebra over a monad. For a rapid description of the monad, one can refer
to \cite[pp 14--17]{Leinster}, see also \cite{BataninPenon} \cite{ChengLauda} 
\cite{ChengMakkai} \cite{Futia}. 
Penon introduces a category $\mQ$
consisting of arrows of ``$\omega$-magmas'' $M \rightarrt S$
where $S$ is a strict $\omega$-category, together with a contractive
structure on $\pi$. The monad is adjoint to the functor ``underlying globular set of $\pA ^{\sharp}$''. If $\pA $ is a globular set then it goes to an element of $\mQ$
with $S$ being the free strict $\omega$-category generated by $\pA $.
Note that this free category is the set of globular pasting diagrams as in
Section \ref{sec-batanin}, and $M$ may
be viewed as some sort of family of elements over $S$ with a contractible structure.
In this sense Penon's definition uses objects of the same sort as Batanin's definition,
indeed Batanin has made
a more precise comparison in \cite{BataninPenon}. Penon's definition uses
globular sets with identities (``reflexive globular sets'')
whereas Batanin's were without them, so \cite{BataninPenon}
proposes a modified version of Penon's definition with non-reflexive globular sets. 
Cheng and Makkai have pointed out that it is better to use the non-reflexive version,
since the reflexive version doesn't lead to all the objects one would want \cite{ChengMakkai}, essentially because of the Eckmann-Hilton argument.
Futia proposes a generalized family of Penon-style definitions in \cite{Futia}.  
In these definitions, 
one could say that globally the goal is to be able to parametrize higher compositions
sorted according to their shapes which are globular pasting diagrams.


\chapter{Weak enrichment over a cartesian model category: an introduction}
\label{cartenr1}

To close out the first part of the book, we describe in this chapter the
basic outlines of the theory which will occupy the rest of the work.
The basic idea, already considered in Pelissier's thesis, is to abstract 
Tamsamani's iteration process to obtain a theory of $\mM$-enriched categories,
weak in Segal's sense,
for a model category $\mM$. 

\section{Simplicial objects in $\mM$}

The original definitions of Segal category, Tamsamani $n$-category, and Pelissier's
enriched categories, took as basic object a functor $\pA :\Delta ^o\rightarrt \mM $.
The first condition is that the image $\pA _0$ of $[0]\in\Delta $ should be a ``discrete
object'', that is the image of a set under the natural inclusion $\Sets \rightarrt \mM$
which sends a set $X$ to the colimit of $\ast$ over the discrete category corresponding
to $X$. This version of the theory therefore requires, at least, some axioms saying that
the functor $\Sets \rightarrt \mM$ is fully faithful and compatible with disjoint
unions. Thus $\pA _0$ may be viewed as
a set and the expression $x\in \pA _0$ means that $x$ is an element of the corresponding
set, equivalently $x:\ast \rightarrt \pA _0$. The higher elements of the simplicial
object will be denoted $\pA _{m/}$. 

Then, the Segal category condition says that the Segal maps 
$$
\pA _{m/} \rightarrt \pA _{1/}\times _{\pA _0}\cdots \times _{\pA _0} A_{1/}
$$
are supposed to be weak equivalences.

The pair of structural maps 
$(\partial _0,\partial _1):\pA _{1/}\rightarrt \pA _0\times \pA _0$
serves to decompose 
$$
\pA _{1/}= \coprod _{x,y\in \pA_0}\pA (x,y)
$$
where $\pA (x,y)$ is the inverse image of $(x,y)\in \pA _0\times \pA _0$, 
or more precisely 
$$
\pA (x,y):= \pA _{1/}\times _{\pA _0\times \pA _0}\ast
$$
with the right map of the fiber product being given by $(x,y):\ast \rightarrt \pA _0\times \pA _0$.
We can similarly decompose
$$
\pA _{m/}= \coprod _{(x_0,\ldots , x_m)\in \pA _0^{m+1}} \pA (x_0,\ldots , x_m)
$$
and the Segal condition may be expressed equivalently as saying that
$$
\pA (x_0,\ldots , x_m)\rightarrt \pA (x_0,x_1)\times \cdots \times \pA (x_{m-1},x_m)
$$
is a weak equivalence in $\mM$. 

\section{Diagrams over $\Delta _X$}

Upon closer inspection, most of the arguments about $\mM$-Segal categories
can really be phrased in terms of the objects $\pA (x_0,\ldots , x_m)$;
and in these terms, the Segal condition involves only a product rather than
a fiber product. So, it is natural and useful to consider the objects $\pA (x_0,\ldots , x_m)$
as the primary objects of study rather than the $\pA _{m/}$. This economizes some hypotheses
and arguments about discrete objects and fiber products. 

This point of view has been introduced by Lurie \cite{LurieGC}. 
For any set $X$, define the category
$\Delta _X$ whose objects are finite linearly ordered sets decorated by elements of $X$,
that is to say an object of $\Delta _X$ is an ordered set $[m]\in \Delta$ plus a
map of sets $x_{\cdot}:[m]\rightarrt X$. This pair will be denoted $(x_0,\ldots , x_m)$,
that is it is an $m+1$-tuple of elements of $X$. The morphisms in the category
$\Delta _X$ are the morphisms of $\Delta$ which respect the decoration,
so for example the three standard morphisms $[1]\rightarrt \nocom [2]$ yield morphisms
of the form 
$$
(x_0,x_1)\rightarrt (x_0,x_1,x_2),\;\;\;
(x_1,x_2)\rightarrt (x_0,x_1,x_2),\;\;\;
(x_0,x_2)\rightarrt (x_0,x_1,x_2).
$$
Now, an $\mM$-Segal category will be a pair $(X,\pA )$ where $X$ is a set, called the 
{\em set of objects}, and $\pA :\Delta ^o_X\rightarrt \mM$ is a functor denoted by
$$
([m],x_{\cdot})=(x_0,\ldots ,x_m)\mapsto \pA (x_0,\ldots , x_m)
$$
or just $(x_0,\ldots ,x_m)\mapsto \pA (x_0,\ldots , x_m)$ if there is no danger of
confusion,
such that the Segal maps 
$$
\pA (x_0,\ldots , x_m)\rightarrt \pA (x_0,x_1)\times \cdots \times \pA (x_{m-1},x_m)
$$
are weak equivalences in $\mM$.
At $m=0$ the Segal condition says that $\pA _0(x_0)\rightarrt \ast$ is a weak equivalence.
This is a sort of weak unitality condition, but for our purposes it is generally
speaking better to impose the {\em strict unitality condition} that $\pA _0(x_0)=\ast$
for any $x_0$. This condition becomes essential when we consider direct products. 

At $m=1$, the morphism space between two elements $x,y\in X$ is $\pA (x,y)$. 
At $m=2$ the usual diagram using the three standard morphisms, serves to define
the composition operation in a weak sense:
\begin{equation}
\label{compdiagram}
\pA (x,y)\times \pA (y,z)\leftarr \pA _2(x,y,z)\rightarrt \pA (x,z)
\end{equation}
with the leftward arrow being a weak equivalence in $\mM$. For higher values of $m$
we get the higher homotopy coherence conditions starting with associativity at $m=3$. 

\section{Hypotheses on $\mM$}

In the weak enrichment, the composition operation is given by a diagram
of the form \eqref{compdiagram} above, using
the usual direct product $\times$ in $\mM$ and where the leftward arrow is
the Segal map which is required to be a weak equivalence. 

Therefore, the main condition which we need to impose upon $\mM$ is that it be {\em cartesian}, that is to say a monoidal model category whose monoidal operation is the
direct product. This insures that direct product is compatible with the
cofibrations and weak equivalences. 
Monoidal model categories have been considered by many authors,
see Hovey \cite{Hovey} for example, and the cartesian theory is a special
case. This condition will be discussed in Chapter \ref{cartmod1}.

For convenience we also impose the conditions that $\mM$ be {\em left proper},
and {\em tractable}. Tractability is Barwick's slight modification
of J. Smith's notion of {\em combinatorial model category}. Recall that a combinatorial
model category is a cofibrantly generated
one whose underlying category is locally presentable---locally presentable categories are
the most appropriate environment for using the small object argument, one of our staples.

Barwick's tractability
adds the condition that the domains of the generating cofibrations and trivial 
cofibrations, be themselves cofibrant objects. This is useful at some technical
places in the small object argument. Our discussion of these topics is
put together in Chapters \ref{cattheor1} and \ref{modcat1}. 

In Section \ref{sec-interpretations} we consider some additional hypotheses on $\mM$
saying that disjoint unions behave like we think they do; if $\mM$ satisfies
these hypotheses then the discrete-set objects in $\mM$ work
well, and we can use the notation $\pA _{m/}$ for the disjoint
union of the $\pA (x_0,\ldots , x_m)$. This reduces to consideration of simplicial
objects in $\mM$ rather than functors from $\Delta _X^o$. If $\mM$ is a category
of presheaves over a connected category $\Phi$ then it satisfies the additional
hypotheses, and the category of $\mM$-precategories discussed next will again
be a category of presheaves (over a quotient of $\Delta \times \Phi$). The 
fact that iteratively we stay within the world of presheaf categories, is convenient
if one wants to think of the small object argument in a simplified way.

\section{Precategories}

A tractable left proper cartesian
model category $\mM$ is fixed. For the original case of $n$-categories, 
$\mM$ would
be the model category for $(n-1)$-categories constructed according to the
inductive hypothesis. In order for the induction to work, the main goal is to construct
from $\mM$ a new model category, whose objects represent up to homotopy the $\mM$-enriched
Segal categories, and which satisfies the same hypotheses of tractability, left properness,
and the cartesian condition.

If $\mM$ satisfies the additional hypotheses on disjoint unions, then
an $\mM$-enriched category is a functor $\pA :\Delta ^o\rightarrt \mM$ such that 
$\pA _0$ is a discrete set also called $\Ob (\pA )$, and such that the Segal maps
$$
\pA _n\rightarrt \pA _{1/}\times _{\pA _0}\cdots \times _{\pA _0}\pA _{1/}
$$
are weak equivalences in $\mM$.

However, looking at the category of all such functors,
the Segal condition is not preserved by limits or colimits of diagrams.
It would be preserved by homotopy limits, but not even by homotopy colimits, and indeed
the problem of taking a homotopy colimit of diagrams and then imposing the Segal condition
is our main technical difficulty. 

So, in order to obtain a model category structure, we have to relax the Segal condition.
This leads to the basic notion of {\em $\mM$-precategory}.
Our utilisation of the word ``precategory'' is similar to
but not the same as that of \cite{Janelidze}. The reader may refer to the introduction
of \cite{svk} for a discussion of this notion in the original $n$-categorical
context. 

If $\mM$
satisfies the additional condition about disjoint unions, then an $\mM$-precategory
may be defined as a functor $\pA : \Delta ^o\rightarrt \mM$ such that $\pA _0$
is a discrete set, that is to say a disjoint sum of copies of the coinitial 
object $\ast \in \mM$.

In the more general case, an $\mM$-enriched precategory is a pair $(X,\pA )$ where $X$
is a set (often denoted by $\Ob (\pA )$), and $\pA :\Delta ^o_X\rightarrt \mM$
is a functor satisfying the {\em unitality condition} that $\pA (x_0)=\ast$ for
any sequence of length zero (i.e. having only a single element).

In either of these situations, the category $\precat (\mM )$ of such diagrams is
closed under limits and colimits, and furthermore if $\mM$ is locally presentable then
$\precat (\mM )$ is locally presentable too. The category $\precat (\mM )$ will
thus serve as a suitable substrate for our model structure. 
The fibrant objects of the model structure should additionally satisfy the Segal conditions.

An additional benefit of the notation $\Delta _X$ is that it allows us to break down
the argument into two pieces, a suggestion of Clark Barwick \cite{BarwickThesis}. Indeed, 
we obtain two different categories, $\precat (X;\mM )$ and $\precat (\mM )$.
The first consists of all $\mM$-precategories with a fixed set of objects $X$. It is
just the full subcategory of the diagram category $\func (\Delta ^o_X,\mM )$ consisting
of diagrams satisfying the strict unitality condition $\pA _0(x_0)=\ast$. The study of 
$\precat (X;\mM )$, considered first, 
is therefore almost the same as the study of the category of $\mM$-valued diagrams on a fixed category $\Delta ^o_X$.

The category $\precat (\mM )$ is obtained by letting $X$ vary, with a natural definition
of morphism $(X,\pA )\rightarrt (Y,B)$. Once everything is well under way and the
objects of $\precat (\mM )$ become our main objects of study, then we will drop the
set $X$ from the notation: an object of $\precat (\mM )$ will be denoted $\pA $ and
its set of objects by $\Ob (\pA )$, but with the same letter for the functor
$\pA :\Delta ^o_{\Ob (\pA )}\rightarrt \mM$, in other words $\pA $ denotes $(\Ob (\pA ),\pA )$.

\section{Unitality}

The strict unitality condition says that $\pA (x_0)=\ast$. The reason for imposing
this condition, aside from its convenience, is that it is needed to obtain the
cartesian condition on the model category of $\mM$-precategories. Indeed, if we don't
impose the unitality condition, then the precategories must be allowed to have
$\pA (x_0)$ arbitrary, even those would be forced to be contractible by the
Segal condition of length $0$. Product with a non-unital precategory
such that $\pB (x_0,\ldots ,x_m)=\emptyset$ for all sequences of objects,
is not compatible with weak equivalences (as will be discussed in Section
\ref{sec-unitalnecessary}).

Given the fact that the Eckmann-Hilton argument rules out a number of different
approaches to higher categories, as we have seen in Chapter \ref{nonstrict1}
but also as in Cheng and Makkai's remark in Penon's original definition \cite{ChengMakkai},
we should justify why our version of the unitality condition which says that
$\pA (x_0)=\ast$ doesn't also lead to an Eckmann-Hilton argument.

The point is that in the Segal-style definitions, the composition is not a well-defined
operation. So, even if there exist cells which are supposed to be the ``identities'',
there is not a single well-defined composition with the identity. The degeneracies
provide $2$-cells which say that, for an $i$-morphism $f$, some possible composition
of the form $1_{t(f)}\circ f$  or $f\circ 1_{s(f)}$ will be equal to $f$,
but these choices (which we call respectively the left and right degeneracies)
are not the only possible ones. The main step of the Eckmann-Hilton
argument going from 
$$
{\setlength{\unitlength}{.5mm}
\begin{picture}(120,60)
\put(10,30){\circle*{2}}
\put(60,30){\circle*{2}}
\put(110,30){\circle*{2}}
\put(10,30){\line(1,0){100}}
\qbezier(10,30)(35,58)(60,30)
\qbezier(10,30)(35,2)(60,30)
\qbezier(60,30)(85,58)(110,30)
\qbezier(60,30)(85,2)(110,30)

\put(30,35){\ensuremath{f \Downarrow}}
\put(30,22){\ensuremath{1 \Downarrow}}
\put(80,35){\ensuremath{1 \Downarrow}}
\put(80,22){\ensuremath{g \Downarrow}}
\end{picture}
} 
$$
to
$$
{\setlength{\unitlength}{.5mm}
\begin{picture}(90,60)
\put(10,30){\circle*{2}}
\put(80,30){\circle*{2}}
\put(10,30){\line(1,0){70}}
\qbezier(10,30)(45,60)(80,30)
\qbezier(10,30)(45,0)(80,30)

\put(30,35){\ensuremath{1\circ f = f \Downarrow}}
\put(30,22){\ensuremath{g\circ 1 = g \Downarrow}}
\end{picture}
} 
$$
involves glueing the left degeneracy for $f$ on top to the right degeneracy for $g$
on the bottom,
generating a coproduct of cells which doesn't fit into any canonical 
global composition operation for the
four $2$-morphisms at once. And similarly for the step involving vertical compositions.
The information on composition with units which comes from the unitality condition and
the degeneracies of $\Delta$, is luckily not enough to make the Eckmann-Hilton argument 
work. Because we're close to the borderline here,
it is clear that some care should be taken to verify everything
related to the unitality condition in the technical parts of our construction. 

Unitality will therefore be considered n the context of more general up-to-homotopy
finite product theories, in Chapter \ref{algtheor1}.

\section{Rectification of $\Delta _X$-diagrams}

The reader should now be asking the following question: wouldn't it be better
to consider $\mM$ as some kind of higher category, and to look at {\em weak}
functors $\Delta ^o_X\rightarrt \mM$? This would certainly seem like the most
natural thing to do. Unfortunately, this idea leads to ``bootstrapping'' 
problems both philosophical as well as practical.  
On the philosophical level, the really good version of $\mM$ as
a higher category, is to think of $\mM$ as being enriched over itself. We exploit this
point of view starting from Section \ref{sec-ncatcat1} where $\mM$ considered as an
$\mM$-enriched category is called $\Enr (\mM )$. However, if we are looking to define
a notion of $\mM$-enriched category, then we shouldn't start with something which is
itself an $\mM$-enriched category. One can imagine getting around this problem by noting that
$\mM$, considered as a category enriched over itself, is actually strictly 
associative; however for looking at functors to $\mM$ we need to go to a 
weaker model, and we end up basically having at least to pass through the notion
of strict functors $\Delta _X\rightarrt \mM$. One could alternatively
say that instead of requiring that $\mM$ be considered as an $\mM$-enriched category,
we could look at a slightly easier structure such as the Dwyer-Kan simplicial 
localization associated to $\mM$. In this case, we would need a theory of
weak functors from $\Delta _X^o$ to a simplicial category. This theory has already
been done by Bergner \cite{BergnerModel} \cite{BergnerThreeModels}, so it would be possible
to go that route. However, it would seem to lead to many notational and mathematical
difficulties. 

Luckily, we don't need to worry about this issue. It is well-known that any kind of
weak functors from a usual $1$-category, to a higher category such as comes from a model
category, can be {\em rectified} (or ``strictified'') to actual $1$-functors. 
This became apparent as early as \cite{SGA1} where Grothendieck pointed out that
fibered categories are equivalent to strictly cartesian fibered categories. Since then
it has been well-known to homotopy theorists working on diagram categories, and indeed
the various model structures on the category of diagrams $\func (\Delta ^o_X,\mM )$
serve to provide model categories whose corresponding higher categories, in whatever
sense one would like, are equivalent to the higher category of weak functors. In the
context of diagrams towards Segal categories, an argument is given in \cite{descente}. 

Due to the philosophical bootstrapping
problem mentioned above, I don't see any way of making the argument given in the previous
paragraph into anything other than the heuristic consideration that it is. But, taking 
it as a basic principle, we shall stick to the notion of a usual functor $\Delta _X^o\rightarrt \mM$ as being the underlying object of study.

\section{Enforcing the Segal condition}

We relaxed the Segal condition in order to get a good locally presentable category
$\precat (\mM )$ of $\mM$-precategories. The Segal condition 
should then be built into the model structure, 
for example it is supposed to be satisfied by the fibrant objects.
This guides our construction of the model structure:
a fibrant
replacement should impose or ``force'' the Segal condition, and 
such a process in turn tells us how
to define the notion of weak equivalence.

To understand this,
one should view an $\mM$-enriched precategory as being a prescription for constructing
an $\mM$-enriched category by a collection of ``generators and relations''. The notion of
precategory was made necessary by the need for colimits, so one should think of a 
precategory as being a colimit of smaller pieces. The associated $\mM$-enriched category
should then be seen as the homotopy
colimit of the same pieces, in the model category we are looking for. That is to say
it is an object specified by generators and relations.  
This is explained in some detail
for the case of $1$-categories in Section \ref{genrel1cat}.

The {\em calculus of generators and relations} is 
the process whereby $\pA$ may be replaced with $\Seg (\pA )$ which is, in a homotopical
sense, the minimal object satisfying the Segal conditions with a map $\pA \rightarrt \Seg (\pA )$. Another way of putting it is that we enforce the Segal conditions using the
small object argument.  
In order to
find the model category, we should define and investigate closely this process of
generating an $\mM$-enriched category. 

The construction $\Seg (\pA )$ doesn't in itself change the set of objects of $\pA $,
so we can look at it in the smaller category $\precat (X,\mM )$. There, it can be considered
as a case of left Bousfield localization. This way of breaking up the procedure
was suggested by Barwick. Luckily, 
the left Bousfield localization which occurs here has a particular form which
we call ``direct'', in which the the weak equivalences may be characterized
explicitly, and we develop that theory with general
notations in Chapter \ref{direct1}. Going to a more general situation helps to clarify
and simplify notations at each stage; it isn't clear that these discussions would have
significant other applications although that cannot be ruled out. Continuing in this way,
we discuss in Chapter \ref{algtheor1} the application of direct left Bousfield localization
to algebraic theories in diagram categories. This formalizes the idea of requiring certain
direct product maps to be weak equivalences, with the objective of applying it to the
Segal maps. Here we refer implicitly to the theories of sketches and algebraic theories.

Then, in Chapters \ref{weakenr1} and \ref{genrel1} we apply the preceding general
discussions to the case of $\mM$-enriched precategories, and define the
operation $\pA \mapsto \Seg (\pA )$ which to an $\mM$ enriched precategory $\pA $ associates
an $\mM$-enriched category, i.e. a precategory satisfying the Segal condition. 
As a rough approximation the idea is to ``force'' the Segal condition in a minimal way,
an operation that can be accomplished using a series of pushouts along standard
cofibrations. 

The passage from precategories to Segal $\mM$-categories 
is inspired by the workings of the theory of simplicial presheaves
as developped by Joyal and Jardine \cite{JoyalLetter} \cite{Jardine}. 
Whereas their ultimate objects of interest were simplicial
presheaves satisfying a descent condition, it was most convenient to consider all
simplicial presheaves and impose a model structure such that the fibrant objects will
satisfy descent. The Segal condition is very close to a descent condition as has
been remarked by Berger \cite{BergerCellularNerve}. 

As in \cite{Jardine}, we are tempted to use the {\em injective model structure}
for diagrams, defining the cofibrations to be all maps of diagrams $\pA \rightarrt \pB$
which induce cofibrations
at each stage $\pA _n\rightarrt \pB_n$. It turns out that a slightly better alternative
is to use a Reedy definition of cofibration, see Chapter \ref{cofib1}. If $\mM$ is itself an injective model
category then they coincide. It can also be helpful to maintain a parallel {\em projective
model structure} where the cofibrations are generated by elementary cofibrations
as originally done by Bousfield. However, the projective structure is not helpful
at the iteration step: it will not generally give back a cartesian model category. 

Once we have a construction $\pA \mapsto \Seg (\pA )$ 
which enforces the Segal condition, a map $\pA \rightarrt \pB$ is said to be a
``weak equivalence'' if $\Seg (\pA )\rightarrt \Seg (\pB )$ satisfies the usual conditions
for being an equivalence of enriched categories, essential surjectivity and 
full faithfulness.  A map of $\mM$-enriched
categories $f:\pA \rightarrt \pB$ is {\em fully faithful} if, for any two objects $x,y\in \Ob (\pA )$
the map $\pA (x,y)\rightarrt \pB (f(x),f(y))$ is a weak equivalence in $\mM$.
Taking a homotopy class projection $\pi _0:\mM \rightarrt \Sets$ gives a 
truncation operation $\tau _{\leq 1}$ from $\mM$-enriched categories to $1$-categories,
and we say that $f:\pA \rightarrt \pB$ is {\em essentially surjective} if
$\tau _{\leq 1}(f):\tau _{\leq 1}(\pA )\rightarrt \tau _{\leq 1}(\pB)$ is an essentially
surjective map of $1$-categories, i.e. it is surjective on isomorphism classes.
The isomorphism classes of $\tau _{\leq 1}(\pA )$ should be thought of as the
``equivalence classes'' of objects of $\pA $. Putting these together, we say that a
map $f:\pA \rightarrt \pB$ is an {\em equivalence of $\mM$-enriched
categories} if it is fully faithful and essentially surjective. 

Now, we say that a map $f:\pA \rightarrt \pB$ of $\mM$-enriched precategories, is a
{\em weak equivalence} if the corresponding map between the $\mM$-enriched categories
obtained by generators and relations $\Seg (f):\Seg (\pA )\rightarrt \Seg (\pB)$
is an equivalence in the above sense. With this definition and any one of the classes
of cofibrations briefly referred to above  and
considered in detail in Chapter \ref{cofib1}, the specification of the model structure
is completed by defining the fibrations to be the morphisms satisfying right lifting
with respect to trivial cofibrations i.e. cofibrations which are weak equivalences.

\section{Products, intervals and the model structure}

The introduction of $\mM$-precategories together with the operation $\Seg$
allows us to define {\em pushouts of weakly $\mM$-enriched categories}: if $\pA \rightarrt
\pB$ and $\pA \rightarrt \pC$ are morphisms of weak $\mM$-enriched categories, then the
pushout of diagrams $\Delta ^o\rightarrt \mM$ gives an $\mM$-precategory
$\pB\cup ^{\pA }\pC$. The associated pushout in the world of weakly $\mM$-enriched categories
is supposed to be $\Seg (\pB\cup ^{\pA }\pC)$. Proving that the collections of maps we have
defined above, really do define a closed model category, may be viewed as
showing that this pushout operation behaves well. As came out pretty clearly in
Jardine's construction \cite{Jardine} but was formalized in 
{\em Smith's recognition principle} \cite{Beke} \cite{DuggerCombinatorial} \cite{Barwick}, 
the key step is to prove that 
pushout by a trivial cofibration is again a trivial cofibration.

Before getting to the proof of this property, 
one has to calculate something somewhere, which is what is done in the Chapters
\ref{freecat1} and \ref{product1} leading up to the theorem that the calculus of
generators and relations is compatible with direct products:
$$
\Seg (\pA )\times \Seg (\pB )\rightarrt \Seg (\pA \times \pB)
$$
is a weak equivalence. 
Our proof of this compatibility 
really starts in Chapter \ref{freecat1} about free ordered $\mM$-enriched categories.
These may be used as basic building blocks for the generators defining the model 
structure, so it suffices to check the product condition on them, which is then done 
in Chapter \ref{product1}.

The compatibility between $\Seg$ and direct products leads to what will be 
the main part of the cartesian property for the model category which is
being constructed. This is a categorical analogue of the Eilenberg-Zilber theorem
for simplicial sets. It wouldn't be true if we hadn't kept the degeneracy maps in
$\Delta$, and the strict unitality condition seems to be essential too. 

From this result on direct products, 
a trick lets us conclude the main result for constructing the model
structure via a Smith-type recognition theorem: 
that trivial cofibrations are preserved by
pushout.
For that trick, 
one requires also a good notion of interval, which was the subject of Pelissier's
correction \cite{Pelissier} to an errof in \cite{svk}. Although Pelissier discussed
only the case of Segal categories enriched over the model category $\mK$ of
simplicial sets, his construction transfers
to $\precat (\mM )$ by functoriality using a functor $\mK \rightarrt \mM$. 
A somewhat similar correction was made by Bergner in her construction
of the model category structure
for simplicial categories originally suggested by Dwyer and Kan \cite{BergnerModel}.

The construction of
a natural ``interval category'' is described in Chapter \ref{interval1}.
It is a sort of versal replacement for the 
simple category $0\rBotharrow 1$ with two isomorphic objects. This is the point
where Pelissier's correction \cite{Pelissier} of \cite{svk} comes in, and
in order to make the process fully iterative we just have to point out that
an interval for the case of the standard model category $\mM = \mK$ of simplicial
sets, leads by functoriality to an interval for any other $\mM$. The
good version of the versal interval constructed by Pelissier \cite{Pelissier}
is similar to the interval object for dg-categories subsequently introduced by
Drinfeld \cite{DrinfeldIntervalDG}. We modify slightly Pelissier's construction,
but one could use his original one too.

Once all of these ingredients are in place, we can construct the model structure
in Chapter \ref{mproof1}. We obtain a model category structure on $\precat (\mM )$
which again satisfies all of the hypotheses which were required of $\mM$,
so the process can be iterated. 

Starting with the trivial model category
structure on $\Sets$, the $n$-th iterate $\precat ^n(\Sets )$ is the model category
structure for $n$-precategories as considered in \cite{svk}. If instead we start with
the Kan-Quillen model category $\mK$ of simplicial sets, then $\precat ^n(\mK )$
is the model category of Segal $n$-precategories which was used in the work \cite{descente}
about $n$-stacks. 
We discuss these iterations for weak $n$-categories
(which are Tamsamani's {\em $n$-nerves}) and Segal $n$-categories,
together with a few variants
where the initailzing category of sets is replaced by a category of
graphs or other things, in Chapter \ref{chap-iterate}. 

The internal $\uHom$ operation then 
leads to a category enriched over our new model category, which in
the iterative scheme for $n$-categories gives a construction of the $n+1$-category
$nCAT$.

The last part of the book, not yet included in the present
version, will be dedicated to considering how to write in 
this language some basic elements of higher
category theory, such as inverting morphisms, and limits and colimits. 
We also hope to discuss the Breen-Baez-Dolan stabilization hypothesis,
about the behavior of the theories of $n$-categories for different values of $n$.



\part{Categorical preliminaries}


\chapter{Some category theory}
\label{cattheor1}

In this chapter, we regroup various things which can be said in the context of abstract category theory.
Our discussion is based in large part on the book of Adamek and Rosicky \cite{AdamekRosicky} about
locally presentable and accessible categories. Refer there for historical
remarks about these notions. The applicability of this theory to model categories came
out with J. Smith's notion of {\em combinatorial model category} \cite{SmithCombi}, slightly modified
by Barwick with his notion of {\em tractable model category} \cite{Barwick}. 

One of our main goals is to provide a fairly general discussion of the notion of
cell complex in a locally presentable category. Hirschhorn has formalized the use of cell complexes
for the small object argument and left Bousfield localization, in \cite{Hirschhorn}. However, he used an additional
assumption of a monomorphism property of elements of the generating set of arrows $I\subset \Arr (\mM )$,
encoded in his notion of {\em cellular model category}. We would like to avoid this hypothesis. Indeed, 
one of the main examples which we can use to start out our induction is the
model category on $\Sets$, but as Hirschhorn pointed out this is not cellular.
It turns out that a somewhat more abstract approach to cell complexes works pretty well.

Our discussion covers much the same materiel as Lurie in the appendix to \cite{LurieTopos}.
Lurie introduces a notion of ``tree'' generalizing the standard transfinite cell-addition
process. The basic idea is that once we have attached a certain number of cells, the
next cell is attached along a $\kappa$-presentable subcomplex, but this information is lost
under the usual indexation by an ordinal. In our discussion, just to be different, we'll
stick to the standard ordinal presentation, but we introduce a category of ``inclusions of cell
complexes'', and show that the category of $\kappa$-small inclusions of cell complexes into a given one,
is $\kappa$-filtered. Roughly speaking, an inclusion of cell complexes corresponds to a downward-closed subset of
a tree. We sketch a proof of Lurie's theorem \cite[Proposition A.1.5.12]{LurieTopos} that cofibrations are
cell complexes over $\kappa$-small cofibrations, rather than just retracts of such. 

The main application of this result is to construct the generating set for injective cofibrations.
Again we give a brief account of a proof of the main technical result
in the present chapter, although the reader can also refer to \cite{LurieTopos} and \cite{Barwick}.

This discussion prepares the way for the ``recognition principle'' introduced in Chapter
\ref{modcat1},
based on Smith's recognition principle as reported by  
Barwick \cite{Barwick}. Our addition is
to give a statement which encodes the accessibility argument. The advantage is that the notion of 
accessibility no longer appears in the statement, so we can then use that in later chapters to
construct model categories without needing to discuss the notion of accessibility anymore.  

So, in a certain sense what we are doing here is to evacuate some of the more technical details 
in the theory of model categories, towards these first two chapters. We hope that this will be helpful
to the reader who wishes to avoid this kind of discussion: if willing to take for granted the recognition
principle which will be stated as Theorem \ref{recog} in the next chapter, the reader may largely skip over
the most technical parts of these first two chapters.

In order to avoid repetitive language, we often apply the following conventions about universes.
We assume given at least two universes $\univa \in \univb$. Recall that these are sets which themselves provide models
for ZFC set theory. A {\em category} will mean a category object in $\univb$. 
An example is the category $\Sets_{\univa}$ of sets in $\univa$. 
A {\em small category} will be a category
object in $\univa$, which is also one in $\univb$. Often a category $\Cc$ will have {\em small morphism sets},
that is for any $x,y\in \ob (\Cc )$ the set $Hom _{\Cc}(x,y)\in \univb$ is isomorphic to a set in $\univa$. 

Depending on context, the word ``category'' can sometimes mean ``small category'', or sometimes ``category with small morphism sets''. 
However, when we need to consider categories outside of $\univb$ this will be explicitly mentioned.

Recall that an {\em ordinal} is a set $a$, such that if $x\in y$ and $y\in a$ then $x\in a$; 
and such that $a$ is well-ordered by the strict relation 
$x<y\Leftrightarrow x\in y$ for $x,y\in a$.  For the corresponding non-strict order relation we then have $x\leq y \Leftrightarrow x\subset y$, and
for any $x\in a$ the successor of $x$ is $x\cup \{ x\}$. 

A {\em cardinal} is an ordinal $a$ with the property
that for any $b\in a$, $b$ is not isomorphic to $a$. Any set $x$ has a unique {\em cardinality} $|x|$ which
is a cardinal such that $x\cong |x|$. 

An {\em ordinal (resp. cardinal) of $\univa$} is an ordinal (resp. cardinal) which is an element of $\univa$. These are the ordinals 
(resp. cardinals) for the
model of set theory given by $\univa$. 
In particular, for any $x\in \univa$ we have $|x|\in \univa$. 

We say that an ordinal $\alpha$ is {\em approached by a sequence of cardinality $\lambda$} if there
is a subset $x\subset \alpha$ with $|x|=\lambda$, such that $\alpha$ is the least upper bound of $x$. 
A cardinal $\kappa$ is {\em regular} if it is not approached by any sequence of cardinality $<\kappa$.

\section{Locally presentable categories}

Fix a regular cardinal $\kappa$.
A category $\Phi$ is said to be {\em $\kappa$-filtered} 
if for any collection of $< \kappa$ objects $X_i\in \Cc$,
there exists an object $Y$ and morphisms $X_i\rightarrt Y$; and for any pair of objects $X$ and $Y$,
and any collection of $< \kappa$ morphisms $f_i:X\rightarrt Y$
there exists a morphism $g:Y\rightarrt Z$ such that all the $gf_i$ are equal. Note that
taking an empty set of objects in the first condition implies that $\Cc$ is nonempty. 
A {\em $\kappa$-filtered colimit} is a colimit over a $\kappa$-filtered index category.
See \cite[Remark 1.21]{AdamekRosicky}.

Let $\Cc$ be a category. We assume that $\Cc$ admits $\kappa$-filtered colimits. Then, 
say that an object $X\in \Cc$ is {\em $\kappa$-presentable} if, for any $\kappa$-filtered colimit
$\colim _{i\in \Phi}Y_i = Z$, the map
\begin{diagram}
\colim _{i\in \Phi}\Hom _{\Cc}(X, Y_i) &\rightarr &\Hom _{\Cc}(X,Z)
\end{diagram}
is an isomorphism of sets. 
An object $Z\in \Cc$ is said to be {\em $\kappa$-accessible} if it can be expressed as a 
$\kappa$-filtered colimit of $\kappa$-presentable objects. 

\begin{definition}
Let $\kappa$ be a regular cardinal. 
A category $\Cc$ with small morphism sets is called {\em $\kappa$-accessible} if: 
\newline
(1)---$\Cc$ admits $\kappa$-filtered colimits;
\newline
(2)---the full subcategory of $\kappa$-presentable objects is equivalent to a small category;
\newline
(3)---every object of $\Cc$ is $\kappa$-accessible.
\newline
Furthermore, if $\Cc$ admits all small colimits then we say that $\Cc$ is {\em locally $\kappa$-presentable}.
A category which is locally $\kappa$-presentable for some regular cardinal $\kappa$ is called {\em locally presentable}.
\end{definition}

For the countable cardinal $\kappa = \omega$, the terminology ``locally finitely presentable'' is 
interchangeable with ``locally $\omega$-presentable''. 
 
\begin{theorem}
Suppose $\Cc$ is locally presentable. Then it is complete, i.e. it admits small limits too.
Each object has only a small set of subobjects up to isomorphism. All $\kappa$-filtered
colimits commute with $\kappa$-small limits (i.e. limits over categories
cardinality $<\kappa$). For any $X\in \Cc$, the subcategory $\Cc _{\kappa}/X$ of $\kappa$-presentable objects
of $\Cc$, is $\kappa$-filtered and $X$ is canonically the colimit of the forgetful
functor on $\Cc _{\kappa}/X$. 

If $\Cc$ is locally $\kappa$-presentable then for any regular cardinal $\kappa  >\kappa$
it is locally $\kappa '$-presentable too.
\end{theorem}
\begin{proof}
See \cite{MakkaiPare}, or \cite{AdamekRosicky}, Proposition 1.22,
Corollary 1.28, Remark 1.56 and Proposition 1.59.
For the last sentence see \cite[page 22]{AdamekRosicky}. 
\end{proof}

\begin{lemma}
\label{diagpres}
Suppose $\Psi$ is a small category, and $\Cc$ is locally $\kappa$-presentable. 
Then the category $\diag (\Psi , \Cc )$ of diagrams from $\Psi$ to $\Cc$, is locally $\kappa$-presentable.
The $\kappa$-presentable diagrams in $\diag (\Psi , \Cc )$ are exactly the functors $F:\Phi \rightarrt \Cc$
such that $F(a)$ is $\kappa$-presentable in $\Cc$ for every $a\in \Phi$. 
\end{lemma}
\begin{proof}
See Makkai and Pare \cite{MakkaiPare}, or Adamek and Rosicky \cite[Corollary 1.54]{AdamekRosicky}.
\end{proof}

\begin{corollary}
\label{preshpres}
Suppose $\Psi$ is a small category, then the category $\presh (\Psi )$ of presheaves of sets on $\Psi$, is locally finitely presentable.
\end{corollary}
\begin{proof}
The category of sets is locally finitely presentable. 
\end{proof}

If $\Cc$ is a category, let $\Arr (\Cc )$ be the category of arrows of $\Cc$,
whose objects are the diagrams of shape $X\rightarrt^f Y$ in $\Cc$. The 
morphisms  in $\Arr (\Cc )$ from $X\rightarrt^f Y$
to $X'\rightarrt ^{f'} Y'$ are the commutative squares 
$$
\begin{diagram}
X & \rightarr & X' \\
\downarr^{f} & & \downarr_{f'}\\
Y & \rightarr & Y' 
\end{diagram}.
$$

\begin{corollary}
Suppose $\Cc$ is a locally $\kappa$-presentable category. Then  $\Arr (\Cc )$
is locally $\kappa$-presentable. 
\end{corollary}
\begin{proof}
Indeed, $\Arr (\Cc )= \diag (\Ee , \Cc )$ where $\Ee$ is the category with two objects $0,1$
and a single morphism $0\rightarrt 1$ besides the identities. Therefore \ref{diagpres} applies.
\end{proof}

In a similar way, any category of commutative diagrams of a given shape in a locally presentable
category, will again be locally presentable.

A functor $q:\alpha \rightarrt \beta$ is {\em cofinal} if:
\newline
---for any object $i\in \beta$ there exists $j\in \alpha$ and an arrow $i\rightarrt q(j)$;
\newline
---for any pair of arrows $i\rightarrt q(j)$ and $i\rightarrt q(j')$ in $\beta$ there
are arrows $j\rightarrt j''$ and $j'\rightarrt j''$ in $\alpha$ such that the 
diagram 
$$
\begin{diagram}
i&\rightarr & q(j)\\
\downarr & & \downarr \\
j'&\rightarr & q(j'')
\end{diagram}
$$
commutes (see \cite[0.11]{AdamekRosicky}).

Recall that if $q:\alpha\rightarrt \beta$ is a cofinal functor, then it induces an equivalence between
the theory of colimits indexed by $\beta$ and the theory of colimits indexed by $\beta$, see \cite[page 4]{AdamekRosicky}.

A basic and motivating  example of a locally presentable category is when $\Phi^o$ is a site, and $\mM\subset \Sets_{\univa}^{\Phi}$
is the subcategory of sheaves. In this case $H^{\ast}$ just denotes the identity
inclusion, whereas $T= H_!$ is the sheafification functor. This motivates the following characterization.

\begin{proposition}
A $\univb$-category $\Cc$ with $\univa$-small morphism sets is locally presentable if it has $\univa$-small limits and colimits, and if
there exists a $\univa$-small category $\Phi$ and  adjunction
$$
H_! : \Sets_{\univa}^{\Phi} \rBotharrow \Cc : H^{\ast}
$$
and a regular  cardinal $\kappa \in \univa$
such that $H_!H^{\ast}$ is the identity, and the composition $T:= H^{\ast} H_! : \Sets_{\univa}^{\Phi}{\rightarrt} \Sets_{\univa}^{\Phi}$
commutes with $\beta$-directed colimits. 
Or equivalently,  that $H^{\ast}$ itself preserves $\kappa$-directed colimits. 
\end{proposition}
\begin{proof}
This rephrases \cite{AdamekRosicky}, Theorem 1.46, using an adjoint pair of functors
such that one composition is the identity, instead of a full reflective subcategory. 
\end{proof}

A locally $\kappa$-presentable category $\Cc$ will in general have the property that
$|\Hom (X,Y)|>\kappa $ for two $\kappa$-presentable objects. For example, 
$\Sets _{\univa}^{\Phi}$ is locally $\kappa$-presentable whatever the size of $\Phi$,
but if $|\Phi |>\kappa$ then there can be $>\kappa$ morphisms between $\kappa$-presentable
objects. This is rectified by taking $\kappa$ big enough, but the bound has
to be exponential in $\kappa$ because we look at maps from objects of size $<\kappa$.  

\begin{corollary}
\label{lambdamuexp}
Suppose $\Cc$ is a locally presentable category. There is 
a regular cardinal $\kappa$ such that $\Cc$ is locally
$\kappa$-presentable, such that for any regular cardinals $\lambda , \mu >\kappa$
and objects $X,Y$ such that $X$ is $\lambda$-presentable and $Y$ is
$\mu$-presentable, the set of morphisms $\Hom _{\Cc}(X,Y)$ has size $<\mu^{\lambda}$. 
\end{corollary}
\begin{proof}
Suppose $\Cc$ is $\kappa _0$-presentable to begin with. 
The total cardinality of a presheaf $A\in \Sets _{\univa}^{\Phi}$ is the 
sum of the cardinalities of the values $A(x)$ for $x\in \Phi$. 
Use the characterisation of the previous proposition, and choose a new regular cardinal
$\kappa_1 >\sup (\kappa _0,|\Phi |)$. For any $\lambda \geq \kappa_1$  
the $\lambda$-presentable objects of $\Cc$ are those of the form $H_!(A)$ for
presheaves $A:\Phi \rightarrt \Sets _{\univa}$ of total cardinality $<\lambda$
(see [Example 1.31]{AdamekRosicky}). Choose a regular cardinal $\kappa \geq \kappa _1$ so that
the total cardinality of $H^{\ast}H_!(B)$ is $<\kappa$ for any 
presheaf $B$ of total cardinality $<\kappa _1$.

Now suppose $\lambda ,\mu>\kappa$. 
By \cite[Remark 1.30(2)]{AdamekRosicky}, the $\lambda$-presentable objects of $\Cc$
are $\kappa_1$-filtered colimits of size $<\lambda$, of $\kappa_1$-presentable objects
(and the same for $\mu$).
Suppose $X$ and $Y$ are $\lambda$-presentable and $\mu$-presentable objects respectively. Write 
$X=\colim _{i\in I}H_!(A_i)$ (resp. $Y=\colim _{j\in J}H_!(B_j)$) where $I$ (resp. $J$) is a
$\kappa_1$-filtered category of size $<\lambda$ (resp. of size $<\mu$), 
and $A_i$ and $B_j$ are presheaves of total cardinality $<\kappa _1$. 
Then $H^{\ast}H_!(B_j)$ has total cardinality $<\kappa$. 
Now 
$$
\Hom _{\Cc}(X,Y) =\lim _{i\in I}\left( \colim _{j\in J}\Hom _{\Cc}(A_i,H^{\ast}H_!(B_j))\right) .
$$
This has size $<\mu ^{\lambda }$. 

One should also be able to prove this using the characterization of locally presentable
categories as categories of models of limit theories \cite[Theorem 5.30]{AdamekRosicky}. 
\end{proof}

The following lemma will be useful in dealing with unitality conditions in Chapter \ref{algtheor1}.

\begin{lemma}
\label{connstar}
Suppose $\mM$ is a category with coinitial object $\ast$, and suppose $\alpha$ is a nonempty connected small category
(that is, a category whose nerve is a connected simplicial set). Then the colimit of the constant functor $C_{\cdot}:\alpha \rightarrt \mM$
defined by $C_i= \ast$ for all $i\in \alpha$, exists and is equal to $\ast$. 
\end{lemma}
\begin{proof}
There is a unique compatible system of morphisms $\phi _i: C_i=\ast \rightarrt \ast$. We claim that this makes $\ast$ into 
a colimit of $C_{\cdot}$. Suppose $U\in \mM$ and $\psi _i: C_i\rightarrt U$ is a compatible system of morphisms.
Pick $i_0\in \alpha$ and use $f:= \psi _{i_0}: C_{i_0}=\ast \rightarrt U$. We claim that for any $j\in \alpha$ the 
composition $f\phi _j$ is equal to $\psi _j$. Let $\alpha '$ be the subset of objects of $\alpha$ for which this is true,
nonempty since it contains $i_0$.
If $j'\in \alpha '$ and if $g:j\rightarrt j'$ is an arrow of $\alpha$ then 
$$
f\phi _j = f\phi _{j'}C_g = \psi _{j'}C_g = \psi _j
$$
so $j\in \alpha '$. Suppose $j'\in \alpha '$ and if $g:j'\rightarrt j$ is an arrow of $\alpha$.
Then 
$$
f \phi _jC_g = f\phi _{j'} = \psi _{j'} = \psi _jC_g
$$
but $C_g$ is an isomorphism so $f\phi _j=\psi _j$ and $j\in \alpha '$. These two steps imply inductively, using connectedness of $\alpha$,
that $\alpha '=\Ob (\alpha )$. Thus $f\phi _{\cdot} = \psi _{\cdot}$. Clearly $f$ is unique, so we get the required universal property.
\end{proof}

\section{Monadic projection}
\label{sec-monadic}

Suppose $\Cc$ is a category, and $\Rr \subset \Cc$ a full subcategory. 
We assume that $\Rr$ is stable under isomorphisms. 

A {\em monadic projection} from $\Cc$ to $\Rr$ is a functor $F: \Cc \rightarrt \Cc$
together with a natural transformation $\eta _X: X\rightarrt F(X)$,
such that:
\newline
(Pr1)---$F(X)\in \Rr$ for all $X\in \Cc$; 
\newline
(Pr2)---for any $X\in \Rr$, $\eta _X$ is an isomorphism; and
\newline
(Pr3)---for any $X\in \Cc$, the map $F(\eta _X): F(X)\rightarrt F(F(X))$ is an isomorphism. 

\begin{lemma}
\label{monadicequal}
Suppose $(F,\eta )$ is a monadic projection. Then the two isomorphisms
$F(\eta _X)$ and $\eta _{F(X)}$ from $ F(X)$ to  $F(F(X))$ are equal. 
\end{lemma}
\begin{proof}
Naturality of $\eta$ with respect to the morphism $\eta _X$ gives the commutative diagram
$$
\begin{diagram}
X & \rightarr^{\eta _X} & F(X) \\
\downarr^{\eta _X} & & \downarr_{\eta _{F(X)}}\\
F(X) & \rightarr^{F(\eta _X)} & F(F(X))
\end{diagram} .
$$
For any $X\in \Rr$, $\eta _X$ is an isomorphism so composing with its inverse
we get $F(\eta _X)= \eta _{F(X)}$.

On the other hand, we can also apply $F$ to the above diagram. By (Pr3), for any $X\in \Cc$
we have that $F(\eta _X)$ is an isomorphism so again composing with its inverse
we conclude that $F(F(\eta _X)) = F(\eta _{F(X)})$
for any $X\in \Cc$. 

For arbitrary $X\in \Cc$, consider the diagram
$$
\begin{diagram}
F(X) & \rightarr^{F(\eta _X)} & F(F(X)) \\
\downarr^{\eta _{F(X)}} & & \downarr_{\eta _{F(F(X))}}\\
F(F(X)) & \rightarr^{F(F(\eta _X))} & F(F(F(X)))
\end{diagram} .
$$
It commutes by naturality of $\eta$ with respect to the morphism $F(\eta _X)$. On the other hand,
by the first statement we proved above, and noting that $F(X)\in \Rr$ by (Pr1), we have
$\eta _{F(F(X))} = F(\eta _{F(X)})$. On the other hand, by the second statement we proved above,
$F(F(\eta _X)) = F(\eta _{F(X)})$. We have now shown that both second maps in the two equal
compositions of this diagram, are the same isomorphism. It follows that the two first maps
along the top and the left side, are the same. This proves the lemma.
\end{proof}

\begin{proposition}
\label{monadicunique}
Suppose $\Rr \subset \Cc$ are as above, and $(F,\eta )$ and $(G,\varphi )$ are two monadic  projections
from $\Cc$ to $\Rr$. Then, for any $X\in \mM$ the maps $F(\varphi _X): F(X)\rightarrt F(G(X))$
and $G(\eta _X): G(X)\rightarrt G(F(X))$ are isomorphisms. The diagram of isomorphisms
$$
\begin{diagram}
F(X) & \rightarr^{F(\varphi _X)} & F(G(X)) \\
\downarr^{\varphi _{F(X)}} & & \uparr_{ \eta _{G(X)}}\\
G(F(X)) & \leftarr^{G(\eta _X)} & G(X)
\end{diagram}
$$
commutes. 
\end{proposition}
\begin{proof}
Define the functor $H(X):= G(F(X))$, with a natural transformation $\psi _X: X \rightarrt H(X)$
defined as the composition
$$
X  \rightarrt ^{\varphi _X} G(X) \rightarrt^{G(\eta _X)} G(F(X)).
$$

Naturality of the transformation $\varphi$ with respect to the morphism $\eta _X$ gives a commutative diagram 
$$
\begin{diagram}
X & \rightarr^{\varphi _X}& G(X) \\
\downarr^{\eta _X} & & \downarr_{G(\eta _X)}\\
F(X) & \rightarr^{\varphi _{F(X)}} & G(F(X))
\end{diagram}
$$
giving the the alternative expression for $\psi$, 
\begin{equation}
\label{alternative}
\psi _X = G(\eta _X) \varphi _X = \varphi _{F(X)}\eta _X.
\end{equation}

We claim that $(H,\psi )$ is again a monadic projection from $\Cc$ to $\Rr$. The first (Pr1) 
is a direct consequence of the same conditions for $F$ and $G$. 

For the second condition, suppose $X\in \Rr$. Then $\eta _X$ is an isomorphism by (Pr2) for $(F,\eta )$,
and $\varphi _{F(X)}$ is an isomorphism by (Pr2) for $(G,\varphi )$ plus (Pr1) for $F$. 
By the expression \eqref{alternative} we get that $\psi _X = \varphi _{F(X)}\eta _X$ is an isomorphism which is (Pr2) for 
$(H, \psi )$. One could instead
use the expression $\psi _X = G(\eta _X) \varphi _X$ and the fact that a functor $G$ preserves isomorphisms. 

For the third condition, note from \eqref{alternative} again,
that $H(\psi _X)$ is obtained by applying $G$ to the composed map
\begin{equation}
\label{composed}
F(X)  \rightarrt^{F(\eta _X)} F(F(X)) \rightarrt^{F(\varphi _{F(X)})} F(G(F(X))).
\end{equation}
The first arrow $F(\eta _X)$ is an isomorphism by (Pr3) for $(F,\eta )$. Consider the diagram
$$
\begin{diagram}
F(X) & \rightarr^{\eta _{F(X)}}& F(F(X)) \\
\downarr^{\varphi _{F(X)}} & & \downarr_{F(\varphi _{F(X)})}\\
G(F(X)) & \rightarr^{\eta _{G(F(X))}} & F(G(F(X)))
\end{diagram} .
$$
It commutes by naturality of $\eta$ with respect to the morphism $\varphi _{F(X)}$. 
The right vertical arrow is the the second arrow in the composition \eqref{composed}.
The top arrow is an isomorphism by (Pr1) and (Pr2) for $(F,\eta )$. The
left arrow is an isomorphism by (Pr1) for $F$ and (Pr2) for $(G,\varphi )$. The bottom arrow is
an isomorphism by (Pr1) for $G$ and (Pr2) for $(F,\eta )$. It follows that the
right vertical arrow $F(\varphi _{F(X)})$ is an isomorphism. 

Now apply the above conclusions together with the fact that $G$ preserves isomorphisms, in the expression 
$$
H(\psi _X) = G\left( F(\varphi _{F(X)}) \circ  F(\eta _X) \right) 
$$
to conclude that $H(\psi _X)$ is an isomorphism. This completes the proof of (Pr3) to show that $(H,\psi )$
is a monadic projection.

Continue now the proof of the proposition. We have a morphism 
$$
G(X)\rightarrt^{G(\eta _X) } G(F(X)) = H(X) .
$$
Consider the diagram
$$
\begin{diagram}
G(X) & \rightarr^{G (\varphi _X) }& G(G(X)) & &  \\
\downarr^{G(\eta _X)} & & \downarr_{G(\eta _{G(X)})} & & \\
H(X) & \rightarr^{H(\varphi _X)} & H(G(X)) &  \rightarr^{H(G(\eta _X))} H(H(X))
\end{diagram} .
$$
The square commutes because it is obtained by applying $G$ to the naturality diagram for $\eta$
with respect to the morphism $\varphi _X$. The top arrow is an isomorphism by (Pr3) for $(G,\varphi )$.
The middle vertical arrow is an isomorphism by applying $G$ to (Pr2) for $(F,\eta )$ and using (Pr1) for $G$.
This shows that the composition in the square is an isomorphism, so we deduce that the composition of the
left bottom arrow with the left vertical arrow is an isomorphism. 

The composition along the bottom is equal to $H(\psi _X)$, indeed $\psi _X = G(\eta _X)\varphi _X$ by definition. 
Hence, by (Pr3) for $(H,\psi )$ which was proven above, the composition along the bottom is an isomorphism. 

The morphism $H(\varphi _X)$ at the bottom of the square, now has maps which compose on the left and the right to isomorphisms.
It follows that this map is an isomorphism,
and in turn that the left vertical arrow $G(\eta _X)$ is an isomorphism. 
This is one of the statements to be  proven in the proposition.

The other statement, that $F(\varphi _X)$ is an isomorphism, is obtained by symmetry with the roles of $F$ and $G$ reversed. 

To finish the proof, we have to show that the square diagram of isomorphisms commutes (which means that
it commutes as a usually shaped diagram when the inverses of the isomorphisms are included). 

Apply $FG$ to the diagram in question, and add on another square to get the diagram of isomorphisms
$$
\begin{diagram}
FGF(X) & \rightarr^{FGF(\varphi _X)} & FGFG(X) \\
\downarr^{FG(\varphi _{F(X)})} & & \uparr^{ FG(\eta _{G(X)})}\\
FGGF(X) & \leftarr^{FGG(\eta _X)} & FGG(X) \\
\uparr^{ FG(\varphi _{F(X)})} & & \uparr^{ FG(\varphi _X)}\\
FGF(X)) & \leftarr^{FG(\eta _X)} & FG(X) .
\end{diagram}
$$
The bottom square commutes by $FG$ applied to the naturality square for $\varphi$ with respect to $\eta _X$. 
On the left side, we have the same isomorphism going in both directions. Consider the outer square
$$
\begin{diagram}
FGF(X) & \rightarr^{FGF(\varphi _X)} & FGFG(X)) \\
\uEqualarr & & \uparr^{ FG(\eta _{G(X)}\varphi _X)}\\
FGF(X)) & \leftarr^{FG(\eta _X)} & FG(X) .
\end{diagram}
$$
It is obtained by applying $FG$ to the square
$$
\begin{diagram}
X & \rightarr^{\eta _X} & F(X) \\
\downarr^{\varphi _X} & & \downarr_{F(\varphi _X)}\\
G(X) & \rightarr^{\eta _{G(X)}} & FG(X) 
\end{diagram}
$$
which commutes by naturality of $\eta$ with respect to $\varphi _X$. Since the outer square, and the bottom square
of the above diagram of isomorphisms commute, it follows that the upper square commutes. The upper square was
obtained by applying $FG$ to the square of isomorphisms in $\Rr$ in question, and $FG$ is a functor which is naturally isomorphic to
the identity functor on $\Rr$. Therefore 
the diagram of the proposition, which is a diagram of isomorphisms in $\Rr$, commutes. 
\end{proof}

\section{Miscellany about limits and colimits}
\label{sec-miscellany}

Here is an elementary observation about limits. 

\begin{lemma}
\label{colimmap}
Suppose $F:I\rightarrt J$ is a functor between small categories, suppose $\mM$ is cocomplete, 
and suppose $B: J\rightarrt \mM$ is a diagram. Pullback along $F$ induces the diagram $B\circ F$ which is noted $F^{\ast}B$. 
There is an induced map in $\mM$
$$
\colim _IF^{\ast}B \rightarrt \colim _J B
$$
\end{lemma}
\begin{proof}
Let $C_I$ and $C_J$ be the ``constant diagram'' constructions. One way of defining colimits is by adjunction with $C_{\cdot}$. 
We have a universal map in $\mM ^J$
$$
B\rightarrt C_J(\colim _J B ), 
$$
and applying restriction along $F$ we get a map
$$
F^{\ast}B \rightarrt F^{\ast}(C_J(\colim _J B )) = C_I(\colim _J B ).
$$
By universality of the map $B\circ F \rightarrt C_I(\colim _I F^{\ast}B )$ we get a factorization
$$
F^{\ast}B \rightarrt C_I(\colim _I (B\circ F) ) \rightarrt C_I(\colim _J B )
$$
for a unique map $\colim _I F^{\ast}B \rightarrt \colim _J B$ which is the map in question for the lemma.
\end{proof}

Here is another fact useful in the construction of adjoints. 

\begin{lemma} 
\label{colimsforadj}
Suppose $p: \alpha \rightarrt \beta$ is a 
functor between small categories, and suppose $A:\beta \rightarrt \mM$ is a diagram with values in a cocomplete category $\mM$.
Then there is a natural map $\colim _{\alpha} p^{\ast}(A)\rightarrt \colim _{\beta}A$,
satisfying compatibility conditions in case of compositions of functors.
\end{lemma}
\begin{proof}
Indeed, there is a 
tautological natural transformation of
$\beta$-diagrams from $A$ to the constant diagram with values $\colim _{\beta}A$, and the pullback of this natural transformation to $\alpha$
is a natural transformation of $\alpha$-diagrams from $p^{\ast}(A)$ to the constant diagram with values $\colim _{\beta}A$, which
gives the map $\colim _{\alpha} p^{\ast}(A)\rightarrt \colim _{\beta}A$ in question. If $q:\delta \rightarrt \alpha$
is another functor then the composition of the maps for $p$ and for $q$
$$
\colim _{\alpha} q^{\ast}(p^{\ast}(A))\rightarrt \colim _{\alpha} p^{\ast}(A)\rightarrt \colim _{\beta}A
$$
is the natural map for $pq$. Similarly, if $p$ is the identity functor then the associated natural map is the identity. 
\end{proof}

\section{Diagram categories}
\label{sec-diagramcats}

Suppose $\Phi$ is a small category and $\mM$ a category. Consider the {\em diagram category} $\diag (\Phi , \mM )$ of functors $\Phi \rightarrt \mM $.
If $\mM$ is complete (resp. cocomplete) then so is $\diag (\Phi ,\mM )$ and limits (resp. colimits) of diagrams are computed objectwise,
that is over each object of $\Phi$. 

Suppose $f:\Phi \rightarrt \Psi$ is a functor between small categories. Given any diagram $A:\Psi \rightarrt \mM$ then the composition $A\circ f$
is a diagram $\Phi \rightarrt \mM$ which will also be denoted $f^{\ast}(A)$. This gives a functor $f^{\ast}: \diag (\Psi , \mM )\rightarrt \diag (\Phi , \mM )$.
If $\mM$ is complete (resp. cocomplete) then the functor $f^{\ast}$ preserves limits (resp. colimits) since they are computed objectwise in
both $\diag (\Phi , \mM )$ and $\diag (\Psi , \mM )$. 

We consider the left and right adjoints of $f^{\ast}$. 
This parallels the similar discussion worked out with A. Hirschowitz in 
in \cite[Chapter 4]{descente} but of course these things are of a nature to have
been well-known much earlier.  See also \cite[A.2.8.7]{LurieTopos}.

If $\mM$ is cocomplete then we can construct a left adjoint 
$$
f_!: \diag (\Phi , \mM )\rightarrt \diag (\Psi , \mM )
$$ 
as follows.
For any object $y\in \Psi$ consider the category $f/y$ of pairs $(x,a)$ where $x\in \Phi$ and $a:f(x)\rightarrt y$ is an arrow in $\Psi$.
There is a forgetful functor $r_{f,y}:f/y\rightarrt \Phi$ sending $(x,a)$ to $x$. 

Suppose $A\in \diag (\Phi , \mM )$. Put $f_!(A)(y):= \colim _{f/y} r_{f,y}^{\ast}(A)$. Suppose $g:y\rightarrt y'$ is an arrow. Then
we obtain a functor $c_g: f/y\rightarrt f/y'$ sending $(x,a)$ to $(x,ga)$. Furthermore this commutes with the forgetful functors in
the sense that $r_{f,y'}\circ c_g = r_{f,y}$. Thus $r_{f,y}^{\ast} (A) = c_g^{\ast}(r_{f,y'}^{\ast}(A))$. Applying the above remark
about colimits, we get a natural map 
$$
f_!(A)(y):= \colim _{f/y} c_g^{\ast}(r_{f,y'}^{\ast}(A)) \rightarrt \colim _{f/y'} r_{f,y'}^{\ast}(A) =: f_!(A)(y').
$$
Using the last part of the paragraph about colimits above, we see that this collection of maps turns $f_!(A)$ into a functor
from $\Psi$ to $\mM$, that is an object in $\diag (\Psi , \mM )$. The construction is functorial in $A$ so it defines a functor
$$
f_!: \diag (\Phi , \mM )\rightarrt \diag (\Psi , \mM ).
$$
The structural maps for the colimit defining $f_!(A)(y)$ are maps $A(x)\rightarrt f_!(A)(y)$ for any $a:f(x)\rightarrt y$. 
In particular when $y=f(x)$ and $a$ is the identity we get maps $A(x)\rightarrt f_!(A)(f(x)) = f^{\ast}f_!(A)(x)$. 
This is a natural transformation from the identity on $\diag (\Phi , \mM )$ to $f^{\ast}f_!$. 

On the other hand, suppose $A= f^{\ast}(B)$ for $B\in \diag (\Psi , \mM )$. Then, for any $(x,a)\in f/y$ we get a
map $r_{f,y}^{\ast}(f^{\ast}(B))(x,a) = B(f(x)) \rightarrt^{B(a)} B(y)$. This gives a map from
$r_{f,y}^{\ast}(f^{\ast}(B))$ to the constant diagram with values $B(y)$, hence a map on the colimit
$$
f_!f^{\ast}(B)(y) \colim _{f/y}r_{f,y}^{\ast}(f^{\ast}(B)) \rightarrt B(y).
$$
It is functorial in $y$ and $B$ so it gives a natural transformation from $f_!f^{\ast}$ to the identity. 

\begin{lemma}
\label{fshriek}
Supposing that $\mM$ is cocomplete and with these natural transformations, $f_!$ becomes left adjoint to 
$f^{\ast}$ and one has the formula
$$
f_!(A)(y)= \colim _{f/y} r_{f,y}^{\ast}(A).
$$
\end{lemma}
\begin{proof}
Suppose $A\in \diag (\Phi , \mM )$ and $B\in \diag (\Psi , \mM )$. A map $A\rightarrt f^{\ast}(B)$
consists of giving, for each $x\in \Phi$, a map $A(x)\rightarrt B(f(x))$. This is equivalent to giving,
for each $y\in \Psi$ and $(x,a)\in f/y$, a map $A(x)\rightarrt B(y)$ subject to some naturality constraints
as $x,a,y$ vary. This in turn is the same as giving a map of $f/y$-diagrams from $r_{f,y}^{\ast}(A)$ to
the constant diagram with values $B(y)$, which in turn is the same as giving a map from
$f_!(A) = \colim _{f/y}r_{f,y}^{\ast}(A)$ to $B$. It is left to the reader to verify that these identifications
are the same as the ones given by the above-defined adjunction maps. 
\end{proof}

\begin{lemma}
\label{frightadj}
Suppose that $\mM$ is complete. Then $f^{\ast}$ has a right adjoint denoted
$f_{\ast}$ given by the formula
$$
f_{\ast}(A)(y)= \mylim _{f\backslash y} s_{f,y}^{\ast}(A)
$$
where $f\backslash y$ is the category of pairs $(z,u)$ where 
$z\in \Phi$ and $f(z)\rightarrt^u y$ is an arrow in $\Psi$,
and $s_{f,y}:f\backslash y\rightarrt \Phi$ is the forgetful functor.
\end{lemma}
\begin{proof}
Apply the previous lemma to the functor $f^o:\Phi ^o\rightarrt \Psi ^o$ for
diagrams in the opposite category $\mM ^o$. 
\end{proof}

\section{Enriched categories}
\label{sec-enriched}

Enriched categories have been familiar objects for quite a while, see Kelly's book \cite{Kelly}. Our overall goal is to
discuss a homotopical analogue susceptible of being iterated, so it is worthwhile to recall the classical theory. 
Furthermore, this will provide an important intermediate step of our argument: a weakly enriched category over a model category $\mM$
gives rise to a $\Ho (\mM )$-enriched category in the classical sense, and this construction is conservative for weak equivalences.
This is used notably for Proposition \ref{globalretract32}. 

Suppose $\mE$ is a category admitting finite direct products. This includes existence of the coinitial object $\ast$ which is the
empty direct product.

If $X$ is a set, then a {\em $\mE$-enriched category on object set $X$}
is a collection of objects $A(x,y)\in \mE$ for $x,y\in X$, together with morphisms $A(x,y)\times A(y,z)\rightarrt A(x,z)$
and $\ast \rightarrt A(x,x)$, such that for any $x,y$ the composed map
$$
\ast \times A(x,y)\rightarrt A(x,x)\times A(x,y) \rightarrt A(x,y)
$$
is the identity; the composed map 
$$
A(x,y)\times \ast \rightarrt A(x,y)\times A(y,y) \rightarrt A(x,y)
$$
is the identity; and for any $x,y,z,w$ the diagram 
$$
\begin{diagram}
A(x,y)\times A(y,z)\times A(z,w) & \rightarr & A(x,y)\times A(y,w) \\
\downarr && \downarr \\
A(x,z)\times A(z,w) & \rightarr & A(x,w)
\end{diagram}
$$
commutes. 

An {\em $\mE$-enriched category} is a pair $(X,A)$ as above, but often this will be denoted just by $A$ with $X=\Ob (A)$.
A functor between $\mE$-enriched categories $f:A\rightarrt B$ consists of a map of sets $f:\Ob (A)\rightarrt \Ob (B)$,
and for each $x,y\in \Ob (A)$ a morphism $f_{x,y}:A(x,y)\rightarrt B(f(x),f(y))$ in $\mE$,
such that the diagrams 
$$
\begin{diagram}
\ast & \rightarr & A(x,x) \\
\downarr & & \downarr \\
\ast & \rightarr & B(f(x),f(x))
\end{diagram}
$$
and 
$$
\begin{diagram}
A(x,y)\times A(y,z) & \rightarr & A(x,z) \\
\downarr & & \downarr \\
B(f(x),f(y)) \times B(f(y),f(z))& \rightarr & B (f(x),f(z))
\end{diagram}
$$
commute. 

Let $\Cat (\mE )$ denote the category of $\mE$-enriched categories. It admits direct products too: if $A$ and $B$ are $\mE$-enriched
categories the 
$$
\Ob (A\times B) = \Ob (A)\times \Ob (B),
$$
and for $(x,x')$ and $(y,y')$ in $\Ob (A)\times \Ob (B)$ we have
$$
(A\times B) ( (x,x'), (y,y')) = A (x,y)\times B(x',y').
$$
Hence this construction can be iterated and we can obtain $\Cat ^n(\mE )$, the category of strict $n$-categories enriched in $\mE$ at the top level.  

Starting with $\mE = \Sets$ yields $\Cat (\Sets )= \Cat$, and $\Cat ^n(\Sets )$ is the category of strictly associative and strictly unital $n$-categories.
These objects have been studied a great deal. However, as we have seen in 
Chapter \ref{nonstrict1}, these objects do not have a sufficiently rich homotopy theory, in particular the groupoid objects therein do not model 
homotopy types in any reasonable way. This obervation is the motivation for considering weakly associative objects as in the remainder of this work.

If $\varphi : \mE \rightarrt \mE '$ is a functor compatible with direct products, applying it to the morphism objects of an $\mE$-enriched 
category $A$ gives an $\mE '$-enriched category $\Cat (\varphi )(A)$ with the same set of objects, and morphism objects defined by 
$$
\Cat (\varphi ) (A)(x,y):= \varphi (A(x,y)).
$$
We apply this in particular to the functor $\tau _{\leq 0}: \mE \rightarrt \Sets$ defined by $\tau _{\leq 0}(E):= \Hom _{\mE } (\ast , E)$.
Define $\tau _{\leq 1}:= \Cat (\tau _{\leq 0})$, i.e. 
$$
\tau _{\leq 1}A(x,y)= \Hom _{\mE } (\ast , A(x,y)).
$$
With this, $\tau _{\leq 1}A \in \Cat (\Sets )$ is a usual category. Let $\Iso \tau _{\leq 1}A$ denote its set of isomorphism classes. 
A functor $A\rightarrt B$ of $\mE$-enriched categories is said to be {\em essentially surjective}
if the induced map 
$$
\Iso \tau _{\leq 1}A \rightarrt \Iso \tau _{\leq 1}B
$$
is surjective. A functor $f:A\rightarrt B$ of $\mE$-enriched categories is said to be {\em fully faithful} if, for each $x,y\in \Ob (A)$
the morphism 
$$
f_{x,y}:A(x,y)\rightarrt B(f(x),f(y))
$$
is an isomorphism in $\mE$. A functor is an {\em equivalence of $\mE$-enriched categories} if it is essentially surjective and fully faithful.

These definitions are useful because, as we shall see in Chapter \ref{weakenr1}, a morphism of $\mM$-enriched precategories is a global weak
equivalence if and only if the associated morphism of $\Ho (\mM )$-enriched categories is an equivalence in the present sense. 
So, we can already obtain versions of some of the main closure properties: closure under retracts and 3 for 2.

\begin{lemma}
\label{transferequiv}
Suppose $\varphi : \mE \rightarrt \mE '$ is a functor commuting with direct products. 
Suppose $f:A\rightarrt B$ is an equivalence of $\mE$-enriched precategories. Then 
$$
\Cat (\varphi )(f): \Cat (\varphi )(A)\rightarrt \Cat (\varphi )(B)
$$
is an equivalence of $\mE '$-enriched precategories. This applies in particular to the functor $\varphi = \tau _{\leq 0}$,
to conclude that $\tau _{\leq 1}(f): \tau _{\leq 1}(A)\rightarrt \tau _{\leq 1}(B)$ is an equivalence of categories,
hence $f$ induces an isomorphism of sets $\Iso \tau _{\leq 1}A \cong \Iso \tau _{\leq 1}B$.
\end{lemma}
\begin{proof}
The functor $\varphi$ sends $\ast _{\mE}$ to $\ast _{\mE '}$ so it induces a map
$$
\Hom _{\mE}(\ast , A)\rightarrt \Hom _{\mE '}(\ast , \varphi (A)).
$$
This natural 
transformation 
induces a natural transformation of functors
$$
\tau _{\leq 1,\varphi }: \tau _{\leq 1,\mE}\rightarrt \tau _{\leq 1,\mE '}\circ \Cat (\varphi )
$$
from $\Cat (\mE )$ to $\Cat $, which is the identity on underlying sets of objects. Therefore
the resulting natural transformation 
$$
\Iso \tau _{\leq 1,\varphi }(A): \Iso\tau _{\leq 1,\mE} (A)\rightarrt \Iso\tau _{\leq 1,\mE '}(\Cat (\varphi )A)
$$
is surjective. It follows that if $f:A\rightarrt B$ is essentially surjective in $\Cat (\mE )$ then
$\Cat (\varphi )f : \Cat (\varphi )A\rightarrt \Cat (\varphi )B$ is also essentially surjective. 
On the other hand, $\varphi$ takes isomorphisms to isomorphisms, so
if $f$ is fully faithful then for any $x,y\in \Ob (\Cat (\varphi )A)= \Ob (A)$,
the map
$$
\Cat (\varphi )(f)_{x,y}= \varphi (f_{x,y}): \varphi (A(x,y))\rightarrt \varphi (B(f(x),f(y)))
$$
is an isomorphism in $\mE '$. This shows that if $f$ is an equivalence in $\Cat (\mE )$ then $\Cat (\varphi )(f)$ is
an equivalence in $\Cat (\varphi )(\mE ')$. 

For the second part of the statement, apply this to $\varphi := \tau _{\leq 0}$ which preserves products. 
\end{proof}

\begin{lemma}
\label{Eretract32}
In any category $\mE$, the class of isomorphisms is closed under retracts and satisfies 3 for 2.
\end{lemma}
\begin{proof}
It is easy to see for the category of sets. 
A retract of objects in $\mE$ gives a retract of representable functors to sets; so if we have a retract of an isomorphism then
the resulting retracted natural transformation
is a natural isomorphism, and a morphism in $\mE$ which induces an isomorphism between representable functors, is an isomorphism. 

The 3 for 2 property is easy to see using the inverses of the isomorphisms in question. 
\end{proof}

\begin{theorem}
\label{enrichedretract32}
The notion of equivalence of $\mE$-enriched categories is closed under retracts and satisfies 3 for 2. If $A\rightarrt ^fB$ and $B\rightarrt ^gA$ are functors between $\mE$-enriched
categories such that $fg$ and $gf$ are equivalences, then $f$ and $g$ are equivalences.  
\end{theorem}
\begin{proof}
Suppose $f: A\rightarrt B$ is a retract of an equivalence of $\mE$-enriched precategories, by
a commutative diagram
$$
\begin{diagram}
A & \rightarr & U & \rightarr & A \\
\downarr^{f} & & \downarr & & \downarr_{f}\\
B & \rightarr & V & \rightarr & B
\end{diagram}
$$
such that the horizontal compositions are the identity and the 
middle vertical arrow is an equivalence. We get a corresponding diagram of sets 
$$
\begin{diagram}
\Iso \tau _{\leq 1}(A) & \rightarr & \Iso \tau _{\leq 1}(U) & \rightarr & \Iso \tau _{\leq 1}(A) \\
\downarr & & \downarr & & \downarr \\
\Iso \tau _{\leq 1}(B) & \rightarr & \Iso \tau _{\leq 1}(V) & \rightarr & \Iso \tau _{\leq 1}(B)
\end{diagram}
$$
where the middle vertical arrow is an isomorphism by Lemma \ref{transferequiv}. 
It follows that $\Iso \tau _{\leq 1}(A)\rightarrt \Iso \tau _{\leq 1}(B)$ is an isomorphism, 
in particular $f$ is essentially surjective.

If $x_0,x_1\in  \Ob (A)$ is a pair of objects,
denote the image objects in $U$, $V$ and $B$ respectively by $u_i$, $v_i$ and $y_i$. Then 
we get a commutative diagram 
$$
\begin{diagram}
A (x_0, x_1)& \rightarr & U (u_0, u_1)& \rightarr & A (x_0, x_1)\\
\downarr & & \downarr & & \downarr \\
B(y_0, y_1) & \rightarr & V(v_0, v_1) & \rightarr & B (y_0, y_1) 
\end{diagram}
$$ 
in which the horizontal compositions are the identity and the middle vertical map is an equivalence. 
The class of isomorphisms in $\mE$ is closed under retracts, by the previous lemma. It follows that 
$A (x_0,x_1)\rightarrt B(y_0,y_1)$ is an isomorphism in $\mE$. This proves that $f$ is fully faithful,
completing the proof that the class of equivalences between $\mE$-enriched categories is closed under retracts.

Turn to the proof of the 3 for 2 property which says that if 
$$
A\rightarr^{f}B \rightarr^{g} C
$$
is a composable pair of morphisms and if any two of $f$, $g$ and $gf$ are  equivalences, then
the third one is too. 

Suppose that some two of $f$, $g$ and $gf$ are equivalences of $\mE$-enriched categories.
Applying the truncation functor gives a composable pair of morphisms of sets
$$
\Iso \tau _{\leq 1}(A) \rightarrt \Iso \tau _{\leq 1}(B) \rightarrt  \Iso \tau _{\leq 1}(C) 
$$
and the corresponding two of the maps are isomorphisms by Lemma \ref{transferequiv}.
From 3 for 2 for isomorphisms of sets, it follows that the third map is also an isomorphism, hence the third map among the $f$, $g$ and $gf$ 
is essentially surjective. 

For the fully faithful condition, consider first the two easy cases. If $f$ and $g$ are global weak equivalences then
for any pair of objects $x_0,x_1\in \Ob (A)$ we have a factorization
$$
A(x_0,x_1) \rightarrt B(f(x_0), f(x_1))\rightarrt C(gf(x_0), gf(x_1))
$$
where both maps are weak equivalences in $\mE$. The previous lemma gives 3 for 2 for isomorphisms in $\mE$, so the composed map is an isomorphism, which shows that $gf$ is
fully faithful. 

Similarly, if we assume known that $g$ and $gf$ are equivalences, then in the same factorization we
know that the composed map and the second map are isomorphisms in $\mE$, so again by Lemma \ref{Eretract32} it follows that the first
map is a weak equivalence, showing that $f$ is fully faithful. 

More work is needed for the third case: with the assumption that $f$ and $gf$ are equivalences, to show that $g$ is an
equivalence. Applying Lemma \ref{transferequiv} we get that $\Iso \tau _{\leq 1}(f)$ and 
$\Iso \tau _{\leq 1}(g)\circ \Iso \tau _{\leq 1}(f)$ are isomorphisms of sets, so $\Iso \tau _{\leq 1}(g)$ is an isomorphism.
The problem is to show the fully faithful condition. Suppose $x,y\in \Ob (B)$. Choose $x',y'\in \Ob (A)$ and
isomorphisms in $\tau _{\leq 1}(B)$,  $u\in \tau _{\leq 1}(B)(f(x'),x)$ and $v\in \tau _{\leq 1}(B)(f(y'),y)$ plus their inverses denoted $u^{-1}$ and
$v^{-1}$.
These are really maps $u:\ast \rightarrt B(f(x'),x)$ and $v:\ast \rightarrt B(f(y'),y)$ and similarly for $u^{-1}$ and
$v^{-1}$. The composition map 
$$
B(f(x'),x)\times B(x,y)\times B(y,f(y'))\rightarrt B(f(x'),f(y'))
$$
composed with 
$$
u\times 1\times v^{-1}: \ast \times B(x,y)\times \ast \rightarrt B(f(x'),x)\times B(x,y)\times B(y,f(y'))
$$
gives 
$$
B(x,y)\rightarrt B(f(x'),f(y')).
$$
Using $u^{-1}$ and $v$ gives $B(f(x'),f(y'))\rightarrt B(x,y)$, and the associativity axiom, definition of inverses and unit axioms
combine to say that these are inverse isomorphisms of objects in $\mE$. The action of $g$ on morphism objects respects all of these operations so
we get a diagram 
$$
\begin{diagram}
B(x,y)&\rightarr & C(g(x),g(y)) \\
\downarr & & \downarr \\
B(f(x'),f(y'))& \rightarr & C(gf(x'),gf(y')).
\end{diagram}
$$
The first vertical map is an isomorphism as described above. The second vertical map is an isomorphism for the same reason applied to $C$
and noting that $g(u)$ and $g(v)$ are isomorphisms in $\tau _{\leq 1}(C)$. The bottom map is an isomorphism by Lemma \ref{Eretract32} and because of the
hypotheses that $f$ and $gf$ are fully faithful. Therefore the top map is an isomorphism, showing that $g$ is fully faithful. 

For the last part, suppose given 
functors between $\mE$-enriched
categories $A\rightarrt ^fB$ and $B\rightarrt ^gA$ such that $fg$ and $gf$ are equivalences.
On the level of sets $\tau _{\leq 0}$ we get a pair of maps whose compositions in
both directions are
isomorphisms, it follows that $\tau _{\leq 0}(f)$ and $\tau _{\leq 0}(g)$ are
isomorphisms. It also easily follows that for any objects $x,y\in \Ob (A)$,
the map 
$$
B(f(x),f(y))\rightarrt A(gf(x),gf(y))
$$
has both a left and right inverse, so it is invertible. But since any object of $B$ is
isomorphic to some $f(x)$, an argument similar to the previous one shows that
$B(u,v)\rightarrt A(g(u),g(v))$ is an isomorphism for any $u,v\in \Ob (B)$. 
Doing the same in the other direction we see that both $f$ and $g$ are equivalences. 
\end{proof}

\subsection{Interpretation of enriched categories as functors  $\Delta _X^o\rightarrow \mE$}

Lurie \cite{LurieTopos}
has used an important variant on the notion of nerve of a category (undoubtedly well-known in the $1$-categorical context).
This makes the usual nerve construction apply to an enriched category, without needing to assume anything further about $\mE$. 
In this point of view, the set of objects is singled out as a set while the morphism objects of an $\mE$-enriched category are considered as objects of
$\mE$. The nerve is then neither a functor from $\Delta ^o$ to $\Sets$, nor a functor to $\mE$, but rather a mixture of the two. We will adopt this
point of view when defining $\mM$-enriched precategories later on. It has the advantage of allowing us to avoid consideration of disjoint unions,
sidestepping some of the difficulties of \cite{Pelissier}. 

If $X$ is a set, define the category $\Delta _X$ to consist of all sequences of elements of $X$ denoted $(x_0,\ldots , x_n)$ with $n\geq 0$. 
One should think of such a sequence as a decoration of the basic object $[n]\in \Delta$. The morphisms of $\Delta _X$ are defined in an obvious way
generalizing the morphisms of $\Delta$, by just requiring compatibility of the decorations on the source and target. This will be discussed further in
Chapter \ref{precat1}. 

If $A$ is an $\mE$-enriched category, let $X:= \Ob (A)$. The {\em nerve} of $A$ is the functor $\Delta _X^o \rightarrt \mE$, denoted also by $A$,
defined by 
$$
A(x_0,\ldots , x_n):= A(x_0,x_1)\times \cdots \times A(x_{n-1},x_n).
$$
Here the notations $A(x,y)$ used on the right side are those of the enriched category, but after having made the definition they are seen to be the
same as the notations for the nerve, so there is no contradiction in this notational shortcut. The transition maps for the functor 
$A:\Delta _X^o \rightarrt \mE$ are obtained using the composition maps of $A$, for example 
$$
A(x_0,x_1,x_2)= A(x_0,x_1)\times A(x_1,x_2)\rightarrt A(x_0,x_2)
$$
in the main case that was discussed in the introduction. 

\begin{theorem}
\label{enriched-interp}
The category $\Cat (\mE )$ of $\mE$-enriched categories becomes equivalent, via the above construction, to the category of pairs $(X,A)$ where
$X$ is a set and $A:\Delta ^o_X\rightarrt \mE $ is a functor satisfying the {\em Segal condition} that for any sequence of elements $x_0,\ldots , x_n$
the morphism 
$$
A(x_0,\ldots , x_n)\rightarrt  A(x_0,x_1)\times \cdots \times A(x_{n-1},x_n)
$$
is an isomorphism in $\mE$. 
\end{theorem}

Understanding this theorem is crucial to understanding the weakly enriched version which is the object of this book. It is left as an exercise for the reader.
Note that the definition of morphisms between pairs $(X,A)\rightarrt (Y,B)$ is done in an obvious way, but the reader may consult the corresponding discussion
in Chapter \ref{precat1} below.

\section{Internal $Hom$}
\label{cat-internalHom}

An important aspect of the theory is the {\em cartesian} condition on the model categories involved.
In this section, we explain how compatibility with products induces an internal $\uHom$.

Suppose $\mM$ is a locally presentable category. Say that
{\em direct product distributes over colimits} if for any small diagram $A:\eta \rightarrt \mM$ and any object $B\in \mM$
the natural map 
$$
\colim _{i\in \eta}(A_i \times B)\rightarrt (\colim _{i\in \eta}A_i )\times B
$$
is an isomorphism. If this is the case, then for any pair of diagrams $A:\eta \rightarrt \mM$
and $B:\zeta \rightarrt \mM$, the natural map 
$$
\colim _{(i,j)\in \eta \times \zeta}(A_i \times B_j)\rightarrt (\colim _{i\in \eta}A_i )\times (\colim _{j\in \zeta}B_j )
$$
is an isomorphism, so these two ways of stating the condition are equivalent. 

We say that $\mM$ {\em admits an internal $\uHom$} if, for any $A,B\in \mM$ the functor $E\mapsto \Hom _{\mM}(A\times E,B)$
is representable by an object $\uHom (A,B)$ contravariantly functorial in $A$ and covariantly functorial in $B$
together with a natural transformation $\uHom (A,B)\times A\rightarrt B$. That is to say that a map $E\rightarrt \uHom (A,B)$ is
the same thing as a map $A\times E\rightarrt B$. 

\begin{proposition}
\label{prop-ihom}
Suppose $\mM$ is a locally presentable category such that direct product distributes over colimits. Then $\mM$ admits 
an internal $\uHom$.
\end{proposition}
\begin{proof}
See \cite{JoyOfCats}, Theorem 27.4, applying the fact that locally presentable
categories are co-wellpowered (\cite[Remark 1.56(3)]{AdamekRosicky}). 
\end{proof}

\begin{corollary}
Suppose $\Phi$ is a small category. Then the category of presheaves of sets $\presh (\Phi )=\diag (\Phi ^p, \Sets )$
admits an internal $\uHom$. This may be calculated as follows: if $A,B$ are presheaves of sets on $\Phi$
then $\uHom (A,B)$ is the presheaf which to $x\in \Phi$ associates the set $\Hom _{\presh (\Phi /x)}(A|_{\Phi /x},A|_{\Phi /x})$.
\end{corollary}
\begin{proof}
Note that $\presh (\Phi )$ is locally presentable. 
Direct products and colimits are calculated objectwise, so direct product distributes over colimits since the same is true for 
the category $\Sets$. The explicit description of $\uHom (A,B)$ is classical. 
\end{proof}

\begin{corollary}
Suppose $\mM$ is a locally presentable category such that direct product distributes over colimits, and 
suppose $\Phi$ is a small category. Then direct product distributes over colimits in $\diag (\Phi , \mM )$
and this category admits an internal $\uHom$.
\end{corollary}
\begin{proof}
Again, direct product and colimits are calculated objectwise in $\diag (\Phi , \mM )$. 
\end{proof}

\section{Cell complexes}

For the basic treatment we follow Hirschhorn \cite{Hirschhorn}. Fix a locally presentable category $\mM$.

If $\alpha$ is an ordinal, denote by $[\alpha ]$ the set $\alpha +1$, that is the set of all ordinals
$j\leq \alpha$. This notation extends the usual notation $[n]$ used in designating objects of
the simplicial category $\Delta$. Note that by definition $[\alpha ]$ is again an ordinal. 
We can write interchangeably $i\leq \alpha$ or $i\in [\alpha ]$. 

A {\em sequence} is a pair $(\beta , X_{\cdot})$ where $\beta$ is an ordinal, and $X:[\beta ] \rightarrt \mM$ is a functor.
We usually denote this by the collection of objects $X_i$ for $i\leq \beta$, with morphisms $\phi _{ij}:X_i\rightarrt X_j$ whenever $i<j\leq \beta$. 
A sequence is {\em continuous} if for any $j\leq  \beta$ such that $j$ is
a limit ordinal, the map
$$
\colim _{j<j}X_i \rightarrt X_j
$$
is an isomorphism. 
A sequence yields in particular a morphism $X_0\rightarrt X_{\beta}$. 

Suppose we are given a set of arrows $I\subset \Arr (\mM )$. An {\em $I$-cell complex}
is a continuous sequence $X:[\beta ]\rightarrt \mM$ together with the data, for each $i< \beta$,
of $f_i\in I$ and a diagram
\begin{equation}
\label{attaching}
\begin{diagram}
U_{i} & \rightarr^{f_i} & V_{i} \\
\downarr^{u_i} & & \downarr_{v_i}\\
X_i & \rightarr & X_{i+1}
\end{diagram}
\end{equation}
inducing an isomorphism $X_i \cup ^{U_{i}}V_{i} \cong X_{i+1}$. In general, the choice of 
data $\{ (f_i, u_i, v_i)\} _{i\leq \beta}$ is not uniquely determined by the choice of sequence 
$X_{\cdot}$ although by abuse of notation we usually say just that $\{ X_i \} _{i\leq \beta}$ is a cell complex. 

Denote by $\cell (I)$ the class of arrows in $\mM$ which are of the form $X_0\rightarrt X_{\beta}$ for a cell complex indexed by $[\beta ]$.

\subsection{Cell complexes in presheaf categories}

We would like to develop the idea of a cell complex as corresponding to a sequence of additions of cells,
whereby one could notably envision to change the order of attachment of the cells. 
It is useful to consider first the special case when $\mM$ is a presheaf category and the elements of $I$ are
nontrivial injections of presheaves. In this case, an objects of $\mM$ has its own ``underlying set'' and the cells
in a cell complex can be identified with subsets. This case is much easier to understand, and it covers many if not
almost all of the examples we want to consider. We describe the notion of inclusion of cell complexes in this special case.
It is easier to understand, but on the other hand this facility obscures what is required for treating the general case,
so this discussion should be considered as optional. 

Suppose $\Psi$ is a small category, and $\mM = \diag (\Psi , \Sets )$. For $A\in \mM$ define its {\em underlying set} to be
$$
\underlying (A):= \coprod _{x\in \Psi} A(x).
$$
This construction gives a functor $\underlying : \mM \rightarrt \Sets$ which is faithful and compatible with colimits.
A morphism in $\mM$ is a monomorphism if and only if it goes to an injection of underlying sets. 

Suppose given a set of morphisms $I$ in $\mM$, which we suppose to be 
monomorphisms but not isomorphisms. Then, given a cell complex $\{ X_i \} _{i\leq \beta}$,
the set $\underlying (X_{\beta})-\underlying (X_0)$ is partitioned into nonempty subsets which we call the {\em cells}, and which are indexed by
the successor ordinals 
$i+1\leq \beta$. The cell $C_{i+1}$ is by definition the complement $\underlying (X_{i+1})-\underlying (X_i)$.
Because of the assumption that elements of $I$ are monomorphisms, all of the transition maps in the cell complex are monomorphisms.

Associated to a cell $C_{i+1}$ is its {\em attaching diagram} \eqref{attaching} as above. With those notations,
in our case the map $u_i:U_i\rightarrt X_i$ is uniquely determined by the composed map $U_i\rightarrt X_{\beta}$.

If $\{ Y_j\} _{j\leq \alpha}$ is another cell complex, an {\em inclusion of cell complexes}
consists of morphisms $Y_0\rightarrt X_0$ and $Y_{\alpha}\rightarrt X_{\beta}$ such that the morphism 
$$
\underlying (Y_{\alpha})-\underlying (Y_0) \rightarrt \underlying (X_{\beta})-\underlying (X_0)
$$
is injective, respects the partitions, induces an isomorphism from each cell for $Y_{\cdot}$ to one of the cells for $X_{\cdot}$,
and such that the attaching maps for these corresponding cells are the same, where the attaching maps are considered
as maps into $Y_{\alpha}$ and $X_{\beta}$ respectively.

A presheaf category will satisfy the hypotheses used by Hirschhorn \cite{Hirschhorn}, so in this case
we can refer there. The reader may note that starting with a presheaf
category $\mM$, our construction
of the model category $\precat (X,\mM )$ keeps us within the realm of presheaf categories,
and this remark can be applied iteratively. So, in a certain sense for the main examples
to be envisioned, the case of cell complexes in a presheaf category is sufficient.

\subsection{Inclusions of cell complexes}
\label{inclusioncell}

In spite of the previous remark, it seems like a good idea to search for the highest
convenient level of generality. Thus, 
we turn now to cell complexes in the general situation where $\mM$ is a locally presentable category and $I\subset \Arr (\mM )$ is a 
subset of morphisms. 

In order to develop a theory allowing us to interchange the order of cell attachment, we first
define a notion of inclusion of cell complexes. An intermediate approach was
developped by Hirschhorn in \cite{Hirschhorn} with his notion of {\em cellular model category}. The general case has been treated by Lurie in one of the first appendix
to \cite{LurieTopos}. We give a somewhat different discussion which is certainly less
streamlined than Lurie's, but we hope it will help to gain an intuitive picture of what is
going on. 

Consider a strictly increasing map of ordinals $q:[\alpha ]{\rightarrt} \nocom [\beta ]$. Define the map
$q^- : [\alpha ]{\rightarrt} \nocom [\beta ]$ by
$$
q^- (i):= \inf \{ k, \;\; \forall j<i, q(j)<k\} .
$$ 
Thus $q^- (0) = 0$ and for $i>0$, $q^-(i) = \sup _{j<i}(q(j) + 1)$. 
For successor ordinals, $q^-(i+1) = q(i) + 1$, whereas if $i$ is a limit ordinal then 
$q^-(i)$ is the limit of $q(j)$ for $j<i$. It follows that $q^-$ is strictly increasing, and also continuous, which is to say that if $i$ is a limit ordinal 
then $q^-(i)$ is the limit of $q^- (j)$ for $j<i$.  Furthermore $q^-(i) \leq q(i)$. 
One can see from the above characterization that the intervals 
$$
[q^-(i),q(i)]:= \{ j\in [\beta ],\;\; q^-(i)\leq j \leq q(i)\}
$$
are disjoint and cover $[\beta]$, which is another way of understanding why to introduce
$q^-$. 

The basic idea of the following definition is that a cell complex is a 
transfinite sequence of pushouts along specified morphisms in $I$,
with these specifications forming part of the data of the cell complex (in particular,
a cell complex consists of more than just the resulting
morphism $X_0\rightarrt X_{\beta}$). An ``inclusion of cell complexes'' should  mean
a map compatible with the sequences of pushouts, via a map of ordinals
which should make the labels in $I$ correspond. 

Here is the precise definition. 
Suppose we are given two  cell complexes $(\alpha ,  X_{\cdot}) = (\alpha ,  X_i , \phi _{ij}, f_i, u_i, v_i )$ and 
$(\beta ,  Y_{\cdot})= (\beta ,  Y_i , \psi _{ij}, g_i, r_i, s_i )$. Here $g_i : R_i\rightarrt S_i$ and
$r_i:R_i \rightarrt Y_i$ with $S_i\rightarrt^{\psi _{i,i+1}\cup s_i}Y_{i+1}$. 

An {\em inclusion of cell complexes} from $(\alpha ,  X_{\cdot})$ to $(\beta ,  Y_{\cdot})$
is a pair consisting of a strictly increasing  map $q:\alpha \rightarrt \beta $, and a collection of morphisms
$\xi _i : X_i\rightarrt Y_{q^-(i)}$, subject to some conditions. Before explaining the
conditions, extend the map to $q:[\alpha ]{\rightarrt} \nocom [\beta ]$
by setting $q(\alpha ):=\beta$ and 
notice that the maps $\xi _i$ induce maps denoted $\xi _{i,+}:X_i \rightarrt Y_{q(i)}$,
defined by $\xi _{i,+}:= \psi _{q^-(i)q(i)}\circ \xi _i$ because $q^-(i) \leq q(i)$. 
This includes $\xi _{\alpha , +}: X_{\alpha}\rightarrt Y_{\beta}$. 
In view of the relation $q^-(i+1)= q(i)+1$, we have $\xi _{i+1}: X_{i+1}\rightarrt Y_{q(i)+1}$.  
The conditions are as follows: 
\newline
---first of all that $g_i = f_{q(i)}$ as elements of the set $I$; 
\newline
---second, if $i$ is a limit ordinal then $\xi _i$ is the colimit of the $\xi _j$ for $j<i$,
going from $X_i = \colim _{j<i}X_j$ to $Y_{q^-(i)} = \colim _{j<i}Y_{q^-(j)}$; 
\newline
---and third, for any $i<\alpha$ the composition
$$
\begin{diagram}
R_i &\rightarr^{r_i} & X_i & \rightarr^{\xi _{i,+}} & Y_{q(i)} \\
\downarr^{g_i} & & \downarr & & \downarr \\
S_i &\rightarr^{s_i} & X_{i+1} & \rightarr^{\xi _{i+1}} & Y_{q(i)+1}
\end{diagram}
$$
is equal to the diagram 
$$
\begin{diagram}
U_{q(i)} &\rightarr^{u_{q(i)}} & Y_{q(i)} \\
\downarr^{f_{q(i)}} & & \downarr  \\
V_{q(i)} &\rightarr^{v_{q(i)}}  & Y_{q(i)+1}
\end{diagram}
$$
for the complex $Y_{\cdot}$.

Suppose $(\gamma , Z_{\cdot})$ is a third cell complex and $(p,\zeta _{\cdot})$ is a map of cell complexes
from $(\beta , Y_{\cdot})$ to $(\gamma , Z_{\cdot})$. We define the {\em composition} $(p, \zeta _{\cdot})\circ (q,\xi _{\cdot})$
to be the map of cell complexes $(h, \eta_{\cdot})$ given as follows. First of all, $h:= p\circ q : [\alpha ]{\rightarrt} \nocom [\gamma ]$.
Notice that $h^-(i) = p^-(q^-(i))$, as can be seen from the definitions
of $p^-$ and $q^-$. Thus it makes sense to define $\eta _i$ to be the
composition
$$
X_i \rightarrt^{\xi _i} Y_{q^-(i)} \rightarrt^{\zeta _{q^-(i)}} Z_{h^-(i)}= Z_{p^-(q^-(i))}.
$$
This collection of compositions satisfies the three conditions for being a morphism of cell complexes. 

Let $\Cell (\mM ; I)$ denote the category whose objects are $I$-cell complexes in $\mM$ and whose morphisms are the inclusions of cell complexes.
This is provided with functors ${\bf s}$ and ${\bf t}$ to $\mM$ defined by ${\bf s}(\alpha , X_{\cdot} ):= X_0$ and ${\bf t}(\alpha , X_{\cdot}):= X_{\alpha}$.
If $(q,\xi _{\cdot})$ is an inclusion from $(\alpha , X_{\cdot} )$ to $(\beta , Y_{\cdot} )$ then $\xi _{\alpha , +}: X_{\alpha}\rightarrt Y_{q(\alpha )}$
and we can compose with the transition map for $Y_{\cdot}$ to get a map to $Y_{\beta}$. This defines the functoriality maps for ${\bf t}$.
The collection of maps $X_0\rightarrt X_{\alpha}$ provides a natural transformation from ${\bf s}$ to ${\bf t}$ or equivalently a functor
${\bf a}: \Cell (\mM ; I)\rightarrt \Arr (\mM )$ whose compositions with the source and target functors are 
${\bf s}$ and ${\bf t}$ respectively. 

Recall that $\cell (I)$ denotes the class of arrows in the essential image of ${\bf a}$, in other words it
is the class of morphisms $f:X\rightarrt Y$ in $\mM$ such that there exists $(\alpha , X_{\cdot})\in \Cell (\mM ; I)$ with
$X_0\rightarrt X_{\alpha}$ isomorphic to $f$ in $\Arr (\mM )$.

\begin{lemma}
\label{pushoutcell}
If $(\alpha , X_{\cdot})$ is a cell complex and $X_0\rightarrt Y_0$ is a morphism, then define $Y_i:= X_i \cup ^{X_0}Y_0$.
Using the induced maps for attaching data, we get a cell complex $(\alpha , Y_{\cdot})$. In particular, 
if 
$$
\begin{diagram}
X & \rightarr & Y \\
\downarr & & \downarr \\
Z & \rightarr & W
\end{diagram}
$$
is a cocartesian diagram in $\mM$ such that the top arrow is in $\cell (I)$, then the bottom arrow is also in $\cell (I)$.
There is a tautological inclusion of cell complexes $(t,\phi ): (\alpha , X_{\cdot})\rightarrt (\alpha , Y_{\cdot})$
where $t:[\alpha ]{\rightarrt} \nocom [\alpha ]$ is the identity and $\phi _i:X_i\rightarrt Y_i$ is the tautological inclusion. 
\end{lemma}
\begin{proof}
With the notations used in the definition of cell complex, given that $X_i \cup ^{U_{i}}V_{i} \cong X_{i+1}$
we get 
$$
Y_{i+1}= Y_0\cup ^{X_0}X_{i+1} = Y_0\cup ^{X_0}(X_i \cup ^{U_{i}}V_{i}) = (Y_0\cup ^{X_0}X_i ) \cup ^{U_{i}}V_{i} = Y_i \cup ^{U_{i}}V_{i}.
$$
Similarly $Y_{\cdot}$ satisfies the continuity condition at limit ordinals.
\end{proof}

\begin{lemma}
\label{compositioncell}
Given a continuous sequence indexed by an ordinal such that the transition maps
$X_i\rightarrt X_{i+1}$ are in $\cell (I)$, the composition $X_0\rightarrt
\colim _iX_i$ is again in $\cell (I)$. 
\end{lemma}
\begin{proof}
Just combine together the sequences. 
\end{proof}

Inclusions of cell complexes enjoy some rigidity properties because of the fact that
the map of labels is included as part of the data.

\begin{lemma}
\label{cellinclusionrigid1}
Suppose $q:[\alpha ]{\rightarrt} \nocom [\beta ]$ is a strictly increasing map, and suppose that we
are given two cell complexes 
$(\alpha , X_{\cdot} )$ and $(\beta ,Y_{\cdot})$. For a given initial map $\xi _0: X_0\rightarrt Y_0$,
there is at most one inclusion of cell complexes
$$
(q,\xi _{\cdot}) : (\alpha , X_{\cdot} )\rightarrt (\beta ,Y_{\cdot}),
$$
based on $q$ and starting with $\xi _0$. If $q$ and $\xi_0$ are isomorphisms, and if $(q,\xi _{\cdot})$ exists
then it is an isomorphism. 
\end{lemma}
\begin{proof}
Suppose $(q,\xi '_{\cdot})$ is another inclusion of cell complexes with $\xi '_0=\xi _0$. 
We prove by induction that $\xi '_j = \xi _j$ for all $j\leq \alpha$. Suppose it is known for all $i<j$. If $j$ is a limit ordinal
then the universal property of the colimit $X_j = \colim _{i<j}X_i$ means that the map 
$$
\xi _j: X_j \rightarrt Y_{q^-(j)}
$$
is determined by the $\xi _i$ for $i<j$, and the same is true of $\xi '_j$. The inductive hypothesis that $\xi '_i=\xi _i$ for
$i<j$ therefore implies $\xi '_j=\xi _j$. This treats the case of limit ordinals; the other possible case is
assume $j=i+1$. Then the map $\xi _j= \xi _{i+1}$ fits into the diagram
$$
\begin{diagram}
R_i &\rightarr^{r_i} & X_i & \rightarr^{\xi _{i,+}} & Y_{q(i)} \\
\downarr^{g_i} & & \downarr & & \downarr \\
S_i &\rightarr^{s_i} & X_{i+1} & \rightarr^{\xi _{i+1}} & Y_{q(i)+1}
\end{diagram}
$$
where the horizontal compositions are given as the attaching maps for $(\beta , Y_{\cdot})$.
On the other hand, $X_{i+1}$ is a pushout in the left square, that is the left square is cocartesian.
Thus $\xi _{i+1}$ is uniquely determined by $\xi _{i,+}$ and the attaching map 
$S_i = V_{q(i)} \rightarrt^{v_{q(i)}}   Y_{q(i)+1}$ which identifies with $\xi _{i+1}\circ s_i$.
The same determination holds for $\xi '_{i+1}$, but
$\xi _{i,+}$ is determined by $\xi _i$ so $\xi '_{i,+}=\xi _{i,+}$, therefore $\xi '_{i+1} = \xi _{i+1}$.
This completes the induction step. 

Suppose $q$ and $\xi _0$ are isomorphisms. Then $\alpha = \beta$ and $q^-(i)=q(i) = i$ for all $i\leq \alpha$,
and a similar induction shows that the $\xi _i$ are all isomorphisms. 
\end{proof}

An application of the rigidity property will help with the uniqueness part of the main accessibility result. 

\begin{corollary}
\label{cellinclusionrigid2}
Suppose $(\alpha , X_{\cdot})$ is a cell complex, and suppose
$$
(p,\eta _{\cdot}) : (\beta , Y_{\cdot} )\rightarrt (\gamma ,Z_{\cdot})
$$
is an inclusion of cell complexes. 
Suppose given a set $A$ and a family of inclusions of cell complexes
$$
(q(a),\xi (a)): (\alpha , X_{\cdot})\rightarrt (\beta , Y_{\cdot} )
$$
indexed by $a\in A$, such that the compositions $(p,\eta _{\cdot})\circ (q(a),\xi (a))$
are all the same. Suppose furthermore that the maps $X_0\rightarrt^{\xi _0(a)} Y_0$
are all the same. Then the $(q(a),\xi (a))$ are all the same. 
\end{corollary}
\begin{proof}
Since $p$ is strictly increasing, it is injective. Thus, the condition that the $p\circ q(a)$ are all the
same implies that the $q(a)$ are all the same; the previous lemma immediately says that the $(q(a),\xi (a))$ are all the same.
\end{proof}

It is useful to have a different representation of an inclusion of cell complexes
$(q,\xi _{\cdot}): (\alpha , X_{\cdot})\rightarrt (\beta , Y_{\cdot} )$.
Define a family of objects denoted $X^{j}$ and indexed by $j\in [\beta]$ as follows: let 
$$
X^{j}:= X_i, \;\;\; q^-(i)\leq j \leq q(i).
$$
Note that for any $j\in [\beta ]$ there exists a unique $i\in [\alpha ]$ such that 
$q^-(i)\leq j \leq q(i)$ (see one of the first properties of $q^-$ above). 
We have maps $X^{j}\rightarrt Y_j$, and compatible
transition maps $X^{j}\rightarrt X^{k}$ for any $j\leq k$. 
For exponents in the same sub-interval $q^-(i)\leq j\leq k \leq q(i)$,
the transition maps are the identity. 

In these terms,
the cell attaching data consists of a subset $c(X^{\cdot})\subset \beta$, defined as the set of all $q(j)$ for
$j\in \alpha$, 
which we think of as the {\em subset of cells in $X^{\cdot}$}; plus, for each $j\in c(X^{\cdot})$
two maps $u^j(X^{\cdot})$ and $v^j(X^{\cdot})$ fitting into a diagram
$$
\begin{diagram}
U_j&\rightarr^{u^j(X^{\cdot})} &X^{j} &\rightarr & Y_j \\
\downarr^{g_j} && \downarr && \downarr \\
V_j&\rightarr^{v^j(X^{\cdot})} &  X^{j+1} &\rightarr & Y_{j+1} 
\end{diagram}
$$
where $g_j:U_j\rightarrt V_j$ is the cell attached 
to $Y_{\cdot}$ at $j\in \beta$ and the horizontal compositions are the
attaching maps
$u_j:U_j\rightarrt Y_j$ and $v_j:V_j\rightarrt Y_j$ for $Y_{\cdot}$. For any $j\in c(X^{\cdot})$ the above diagram 
has a cocartesian left square. The composed outer
square is cocartesian, being is the attaching diagram for $Y_j\rightarrt Y_{j+1}$, 
which implies that the right square is cocartesian also. On the other hand, for $j\not\in c(X)$ the map
$X^j\rightarrt X^{j+1}$ is an isomorphism. In the new notation note that the last element is $X^{\beta}= X_{\alpha}$. The maps $u^j(X^{\cdot})$ and $v^j(X^{\cdot})$ are 
contained in the data of the cell complex $(\alpha , X_{\cdot})$ with appropriate
renumbering via $c(X^{\cdot})\cong \alpha$. 

The collection of data described above is the same as the inclusion of cell complexes. 
Furthermore, compositions of inclusions of cell complexes can be understood in this notation:
if $X^{\cdot}$ is a subcomplex of $Y_{\cdot}$ then the indexing ordinal $\alpha  $
for $X$ is isomorphic to the subset $c(X^{\cdot})\subset \beta$;
a subcomplex $Z^{\cdot}$ of $X^{\cdot}$ then consists of a subset 
$$
c(Z^{\cdot})\subset \alpha \cong c(X^{\cdot})\subset \beta ,
$$
and the composed inclusion of cell complexes from $Z^{\cdot}$ to $Y_{\cdot}$
corresponds to the subset $c(Z^{\cdot})\subset \beta$ obtained using transport of
structure along the isomorphism in the middle. 

The category of ordinals with strictly increasing maps, isn't closed under sequential transfinite colimits. However, the
category of ordinals mapping to a given fixed ordinal, does admit sequential colimits. Therefore the same is true of our notion of
inclusion of cell complexes: whereas $\Cell (\mM ; I)$ is not itself closed under sequential colimits, on the other hand
if we look at cell complexes included into a given $(\beta , Z_{\cdot})$, then we can take sequential colimits.  

\begin{proposition}
\label{cellfiltcolim}
Suppose $(\beta , Z_{\cdot})\in \Cell (\mM ; I)$ is a cell complex. Then the category 
$\Cell (\mM ; I)/(\beta , Z_{\cdot})$  is closed under filtered colimits: if $(\alpha (k) , X_{\cdot}(k))$ is a family of
cell complexes indexed by a filtered category $k\in \zeta$ such that the transition maps are inclusions of cell complexes,
all provided with compatible inclusions $(\alpha (k) , X_{\cdot}(k))\rightarrt (\beta , Z_{\cdot})$,
then there is a colimit cell complex $(\alpha (\zeta ), Y_{\cdot})$ again mapping to
$(\beta , Z_{\cdot})$ and such that $Y_0 = \colim _{k\in \zeta} X_0(k)$
and $Y_{\alpha (\zeta )} = \colim _{k\in \zeta} X_{\beta (k)}(k)$. 
\end{proposition}
\begin{proof}
Use the alternative description of an inclusion of cell complexes to $(\beta ,Z_{\cdot})$. We obtain a family
$X^j(k)$ doubly indexed by $j\in [\beta]$ and $k\in \zeta$. Put
$$
Y^j := \colim _{k\in \zeta} X^j(k).
$$
The maps $Y^j\rightarrt Z^j$ are given by the universal property of the colimit. 
The subset of cells is defined as $c(Y^{\cdot}):= \bigcup _{k\in \zeta} c(X^{\cdot}(k))$. 
If $j\in c(Y^{\cdot})$, choose $k_j\in \zeta$ such that $j\in c(X_{\cdot}(k_j))$. Then let 
$k_j\backslash \zeta $ be the category of arrows $k_j \rightarrt m$ in $\zeta$. 
The functor $k_j\backslash \zeta _j\rightarrt \zeta$ is cofinal, so 
$$
Y^j = \colim _{m\in k_j\backslash \zeta} X^j(m).
$$
For each $k_j\rightarrt m$ we have the attaching maps
$$
U_j\rightarrt X^j(k_j)\rightarrt X^j(m),\;\;\;
V_j\rightarrt X^{j+1}(k_j)\rightarrt X^{j+1}(m),
$$
which yield attaching maps 
$$
U_j\rightarrt Y^j,\;\;\;
V_j\rightarrt Y^{j+1}.
$$
If $j\not \in c(Y^{\cdot})$ then for any $k\in \zeta$ the map $X^j(k)\rightarrt X^{j+1}(k)$ is
an isomorphism, so the map $Y^j\rightarrt Y^{j+1}$ is an isomorphism. If $j$ is a limit ordinal then
$$
X^j(k)= \colim _{i<j}X^i(k)
$$
so passing to the double colimit, we have $Y^j= \colim _{i<j}Y^i$. This gives $Y^{\cdot}$ the
structure of cell complex with cell inclusion to $Z^{\cdot}$ according to our second description above. 

This provides a colimit in the category $\Cell (\mM ; I)/(\beta , Z_{\cdot})$,
indeed if $(\beta ', Z'_{\cdot})$ is an object of the category provided with 
inclusions from all the $(X^{\cdot}(k))$, then we can understand the above construction
as corresponding to the same construction in  $\Cell (\mM ; I)/(\beta ', Z'_{\cdot})$,
which gives the inclusion of cell complexes from $Y^{\cdot}$ to $(\beta ', Z'_{\cdot})$.
\end{proof}

\subsection{Cutoffs}
\label{sec-cutoffs}

Suppose $(\beta , Z_{\cdot})$ is a cell complex, and $\beta '\leq \beta$. Define the {\em cutoff}
$C_{\beta '}(\beta , Z_{\cdot})$ to be the cell complex consisting of ordinal $\beta '$ and
family of $Z_i$ for $i\leq \beta ' \leq \beta$ with the same structural data restricted to 
the subset of values of $i$. 
If
$$
(q, \xi _{\cdot}): (\alpha , X_{\cdot})\hookrightarrow (\beta , Z_{\cdot})
$$
is an inclusion of cell complexes,  define its {\em relative cutoff} 
$$
C_{\beta '}(q, \xi _{\cdot}) : C_{\alpha '}(\alpha , X_{\cdot})\rightarrt C_{\beta '}(\beta , Z_{\cdot})
$$
as follows. Using the strictly increasing map 
$q:\alpha \rightarrt \nocom {\beta }$, 
put $\alpha ':= \sup \{ i\leq \alpha , q(i)\leq \beta '\}$.
Define a new map $q'$ equal to $q$ on $\alpha$, but extended to
$[\alpha ]$ by setting $q'(\alpha ):= \beta '$. In general this is different from
$q(\alpha )$. 

Define the relative cutoff inclusion of cell complexes 
to be given by the restricted map $q':[\alpha ']\rightarrt \nocom {[\beta ']}$ and the
family of maps $\xi _{i}$ and associated cell identification data for $X_{\cdot}$ 
restricted to $i\leq \alpha '$. 

In terms of the alternative notation for inclusions of cell complexes, we denote 
the cutoff of a subcomplex $X^{\cdot}$
by just
$C_{\beta '}(X^{\cdot})$. In these terms, the subset of cells is just the intersection
$$
c(C_{\beta '}(X^{\cdot}))=c(X^{\cdot})\cap \beta '
$$
and the attaching data are the same as those of $X^{\cdot}$. 

\subsection{The filtered property for subcomplexes}
\label{filteredsubcomplex}

In the context of our 
next main accessibility result for cell complexes, we  would like to
consider the {\em join} of a
set of inclusions of cell complexes. Suppose $(\beta , Z_{\cdot})= (\beta , Z_{\cdot}, f_{\cdot} , u_{\cdot}, v_{\cdot})$ is a cell
complex, $J$ is a set, and $(q^j, \xi _{\cdot}(j)): (\alpha ^j, X_{\cdot}(j))\hookrightarrow (\beta , Z_{\cdot})$
is a family of inclusions of cell complexes. The ``join'' should be  
a cell complex sitting in the middle 
$$
(\alpha ^j, X_{\cdot}(j))\hookrightarrow (\varphi , Y_{\cdot})\hookrightarrow 
(\beta , Z_{\cdot}).
$$
The set of cells $c(Y^{\cdot})$ should be the union of the sets of cells 
$c(X_{\cdot}(j))\subset \beta$. In order to start off the process, we 
need to make a choice of the place to start
$Y_0$ which should fit into factorizations 
$$
X_0^j\rightarrt Y_0\rightarrt Z_0
$$
for all $j$. From there, one can proceed to attach the required cells, without problem
if $\mM$ is a {\em cellular model category} \cite{Hirschhorn}. This condition basically
means that the arrows in $I$ are monomorphisms enjoying good properties; for example
a presheaf category in which the cofibrations are contained in the injections, is cellular. 

In the general case, we run into the problem that the cell attaching maps needed to
construct $Y^{\cdot}$ are not uniquely determined by those of $Z_{\cdot}$. So, we proceed
differently. The general situation has been treated by Lurie in the appendix of \cite{LurieTopos}. Between the cellular case treated in \cite{Hirschhorn} and the
general treatment in \cite{LurieTopos} our present discussion is undoubtedly superfluous,
included for completeness. 

Suppose $\mM$ is locally $\kappa$-presentable, that the elements of $I$ are arrows whose source and
target are $\kappa$-presentable, and that $I$ has $<\kappa$ elements. 
Let $\Cell (\mM ; I)_{\kappa}$ denote the full subcategory of cell complexes  
$(\alpha , X_{\cdot})$ such that 
$|\alpha | < \kappa$ and $X_0$ is $\kappa$-presentable. 
The following theorem lets us replace the general notion of cofibration, by
a cell complex. This was done by Lurie in \cite{LurieTopos}, A.1.5.12 so the reader
could refer there instead. 

\begin{theorem}
\label{cellover}
With the above hypotheses, 
suppose given a cell complex $(\beta , Z_{\cdot})\in \Cell (\mM ; I)$. Then the category 
$\Cell (\mM ; I)_{\kappa}/(\beta , Z_{\cdot})$ is $\kappa$-filtered. Furthermore
$(\beta , Z_{\cdot})$ is the colimit in $\Cell (\mM ; I)$ of the tautological functor
$$
\Cell (\mM ; I)_{\kappa}/(\beta , Z_{\cdot})\rightarrt \Cell (\mM ; I)
$$
and the arrow $Z_0\rightarrt Z_{\beta}$ is the colimit in $\Arr (\mM )$ of
the composition of the tautological functor with $\Cell (\mM ; I)\rightarrt \Arr (\mM )$. 
\end{theorem}
\begin{proof}
The objects of $\Cell (\mM ; I)_{\kappa}/(\beta , Z_{\cdot})$ are
inclusions of cell complexes $(\alpha , X_{\cdot})\hookrightarrow (\beta , Z_{\cdot})$ such that 
$|\alpha | < \kappa$ and such that $X_0$ is $\kappa$-presentable.

We first show the uniqueness half of the $\kappa$-filtered property. 
Suppose 
$$
(p,\eta _{\cdot} ): (\alpha , X_{\cdot})\rightarrt (\beta , Z_{\cdot})
$$
and 
$$
(p',\eta ' _{\cdot}): (\alpha ', X'_{\cdot})\rightarrt (\beta , Z_{\cdot})
$$
are two objects of $\Cell (\mM ; I)_{\kappa}/(\beta , Z_{\cdot})$,
and suppose given a set $A$ of cardinality $|A|<\kappa$
indexing a family of morphisms from one to the other in 
$\Cell (\mM ; I)_{\kappa}/(\beta , Z_{\cdot})$, that is to say a family of inclusions of cell complexes 
$$
(q(a),\xi _{\cdot}(a)): (\alpha , X_{\cdot})\rightarrt (\alpha ', X'_{\cdot})
$$
such that $(p',\eta ' _{\cdot})\circ (q(a),\xi _{\cdot}(a))= (p,\eta _{\cdot} )$.
Injectivity of $p'$ implies that the $q(a)$ are all the same. 
On the other hand we get the family of maps 
$$
\xi _0(a):X_0\rightarrt X'_0
$$
such that $\eta '_0\circ \xi _0(a)= \eta _0$. Recall that $X_0$ and $X'_0$ are required 
to be $\kappa$-presentable. The same uniqueness part of the property that
$\mM _{\kappa}/Z_0$ is $\kappa$-filtered,
says that the map $X'_0\rightarrt Z_0$
factors as $X'_0\rightarrt^{\phi _0} Y_0\rightarrt Z_0$
through a $\kappa$-presentable object $Y_0$, such that all of the maps 
$$
\phi _0\circ \xi _0(a): X_0\rightarrt Y_0
$$ 
are the
same. As in Lemma \ref{pushoutcell}, define a new cell complex $(\alpha ', Y_{\cdot})$ by setting 
$$
Y_i:= Y_0\cup ^{X'_0}X'_i,
$$
and keeping the same attaching data as for $X'_{\cdot}$, and with the tautological inclusion of cell complexes $(t, \phi _{\cdot})$
from $(\alpha ', X'_{\cdot})$ to $(\alpha ', Y_{\cdot})$. 
The condition on the choice of $Y_0$
together with the fact that the $q(a)$ are all the same, provide the hypotheses required in order to apply Corollary \ref{cellinclusionrigid2}
to conclude that all the maps of cell complexes 
$$
(t, \phi _{\cdot})\circ (q(a),\xi _{\cdot}(a)) : (\alpha , X_{\cdot})\rightarrt (\alpha ', Y_{\cdot})
$$
are the same. 

The elements $Y_i$ are $\kappa$-presentable, so $(\alpha ', Y_{\cdot})\in \Cell (\mM ; I)_{\kappa}$. 
The map of cell complexes $(p',\eta ' _{\cdot})$ extends, using the pushout expressions for $Y_i$ and the maps $Y_0\rightarrt Z_i$,
to a map of cell complexes
$$
(p',\zeta _{\cdot}): (\alpha ', Y_{\cdot})\rightarrt (\beta , Z_{\cdot}).
$$
Via this map, we may consider $(\alpha ', Y_{\cdot})$ as an element of $\Cell (\mM ; I)_{\kappa}/(\beta , Z_{\cdot})$.

We have 
$$
(p',\zeta _{\cdot})\circ (t, \phi _{\cdot}) = (p',\eta ' _{\cdot}),
$$
so $(t, \phi _{\cdot})$ may be viewed as a map in $\Cell (\mM ; I)_{\kappa}/(\beta , Z_{\cdot})$. This gives exactly 
a map there whose compositions with the $(q(a),\xi _{\cdot}(a))$ are all the same, serving to prove the uniqueness half of the
$\kappa$-filtered property.

For the remainder of the theorem, we prove the full statement 
by induction on the length $\beta$ of
the cell complex $Z$. If $\beta = 0$ there is nothing to prove. Hence, we may assume that the statement of
the theorem is known for all cell complexes $(\beta ', Z'_{\cdot})$ with $\beta '<\beta$. 

The main step is to prove that 
$\Cell (\mM ; I)_{\kappa}/(\beta , Z_{\cdot})$ is $\kappa$-filtered. Suppose that $X^{\cdot}(a)$ is
a collection of subcomplexes of $(\beta , Z_{\cdot})$ which are $\kappa$-small, that is with $X^0$ being $\kappa$-presentable and
$| c(X^{\cdot})|<\kappa$, and indexed by a set $a\in A$ with $|A|<\kappa$. 
We have to show that they all map to a single $\kappa$-small subcomplex $Y^{\cdot}$. 

For any $\beta '< \beta$, consider the collection of cutoff complexes (using the
alternate notation at the end of Section \ref{sec-cutoffs}) $C_{\beta '}(X^{\cdot}(a))$ mapping by inclusions
of cell complexes to $C_{\beta '}(\beta , Z_{\cdot})$. 
By the induction hypothesis, there is a single $\kappa$-small subcomplex $Y^{\cdot}(\beta ' )$ of $C_{\beta '}(\beta , Z_{\cdot})$,
such that all of the $C_{\beta '}(X^{\cdot}(a))$ map to $Y^{\cdot}(\beta ' )$. 

Assume that $\beta$ is a limit ordinal approached by a sequence of $\beta '<\beta$ of cardinality $<\kappa$. Then, by transfinite induction over this sequence
we may assume that the choice of $Y^{\cdot}(\beta ' )$ is provided with transition inclusion maps from $Y^{\cdot}(\beta ' )$ to 
$Y^{\cdot}(\beta '')$ whenever $\beta ' \leq \beta ''$ are two members of the sequence. Furthermore, we may take the colimit of all of the 
$Y^0(\beta ')$ and use this as starting element (it is a colimit of length $<\kappa$ so it remains $\kappa$-presentable).
Thus we may assume that the $Y^0(\beta ')$ are all the same, so that Lemma \ref{cellfiltcolim} applies.
Set $Y^{\cdot}:= \colim _{\beta '}Y^{\cdot}(\beta ')$. There are maps from $C_{\beta '}(X^{\cdot}(a))$ to $Y^{\cdot}$.
The different maps from the $X^0(a)$ to $Y^0$, depending on $\beta '$, compose to the same map into $Z^0$. Since 
$X^0(a)$ is $\kappa$-presentable, and the $\kappa$-presentable objects mapping to $Z^0$ form a $\kappa$-filtered category \cite[Proposition 1.22]{AdamekRosicky},
and also
since both $A$ and the sequence of $\beta '$ have size $<\kappa$, we may choose $W$ fitting into $Y^0\rightarrt W^0\rightarrt Z^0$
but with $W$ still being $\kappa$-presentable. Then replace $Y^j$ by $W\cup ^{Y^0}Y^j$. This way, all of the maps 
$X^0(a)\rightarrt Y^0$ will be the same independent of $\beta '$. Then a transfinite induction argument on the cutoff level $\beta '$
shows that the maps $C_{\beta '}(X^{\cdot}(a))\rightarrt Y^{\cdot}$ fit together in the colimit as $\beta '\rightarrt \beta$,
to a collection of maps $X^{\cdot}(a)\rightarrt Y^{\cdot}$. This completes the proof of the filtered property when $\beta$ is
a limit ordinal approached by a sequence of cardinality $<\kappa$.

If on the other hand $\beta$ is a limit ordinal which is not approached by any sequence of cardinality $<\kappa$, then
there is some $\beta '<\beta$ such that all of the cells in $c(X^{\cdot}(a))$ are at $j<\beta '$ for all $a\in A$.
Therefore $C_{\beta '}(X^{\cdot}(a)) = X^{\cdot}(a)$, and these are cell complexes included in 
$C_{\beta '}(\beta , Z_{\cdot})$. The inductive hypothesis says that they all map to a same subcomplex $Y^{\cdot}$ of
$C_{\beta '}(\beta , Z_{\cdot})$, and this $Y^{\cdot}$ serves to show the filtered property for the family 
of cell complex inclusions $X^{\cdot}(a)\rightarrt (\beta , Z_{\cdot})$.

This leaves us with the case when $\beta$ is a successor ordinal: $\beta = \eta + 1$. It is in some sense the main case where we need some work.
This extra work is occasionned by the goal of working without a monomorphism axiom for the elements of $I$, such as was used
by Hirschhorn \cite{Hirschhorn} with his notion of ``cellular model category''. The basic problem is that when we try to add in the last
cell, the attaching maps may be ill-defined because of various different collections of cell attachments up until then.
Hence, we need to backtrack and add some more cells so as to stabilise the collection of attaching maps. For this we use the full
inductive hypothesis about the colimit over our filtered category, which applies to the cell complex $C_{\eta}(\beta , Z_{\cdot})$
of length one less. This says that $Z_{\eta}$ is a $\kappa$-filtered colimit of things obtained by inclusions of $\kappa$-small 
cell complexes. 

To set things up more carefully, write $A = A'\cup A''$ where $A'$ consists of those
$a$ such that $X^{\cdot}(a)$ involves the last cell,
that is to say $\eta \in c(X^{\cdot}(a))$; and $A''$ consists of those $a$ for which 
$\eta \not\in c(X^{\cdot}(a))$. One could suppose that $A''$ consists of a single element,
indeed the $X^{\cdot}(a)$ which don't involve the last cell, are all
subcomplexes of a single subcomplex of $C_{\eta}(\beta , Z_{\cdot})$ by the inductive hypothesis, and we could take this one 
as the single subcomplex indexed by $A''$. 

Now for each $a\in A'$ we have $X^{\eta}(a)\rightarrt Z_{\eta}$. Furthermore we have attaching maps
$$
u_{\eta}(X^{\cdot}(a)): U_{\eta}\rightarrt X^{\eta}(a), 
\;\;\;
v_{\eta}(X^{\cdot}(a)): V_{\eta}\rightarrt X^{\beta}(a)\cong X^{\eta}(a)\cup ^{U_{\eta}}V_{\eta}.
$$
These lift the attaching maps $u_{\eta}:U_{\eta}\rightarrt Z_{\eta}$ and $v_{\eta}:V_{\eta}\rightarrt Z_{\beta}$.
By the inductive hypothesis, there is a single $\kappa$-small cell complex $Y^{\cdot}$ mapping by a cell complex inclusion to $C_{\eta}(\beta , Z_{\cdot})$,
such that all of the $C_{\eta}(X^{\cdot}(a))$ for $a\in A'$, and all of the $X^{\cdot}(a)$ for $a\in A''$,
map to $Y^{\cdot}$. Choose such maps for each $a$, denoted $\psi ^{\cdot}(a)$. Then for any $a\in A'$ we get a composed attaching map
$$
u_{\eta}(a):=  \psi ^{\eta }(a)\circ u_{\eta}(X^{\cdot}(a)) : U_{\eta}\rightarrt Y^{\eta}.
$$
The composition of these maps with $Y^{\eta}\rightarrt Z_{\eta}$ are all the same, they are the attaching map for 
the original complex $Z_{\cdot}$. 

The inductive hypothesis tells us that $Z_{\eta}$ is a $\kappa$-filtered colimit of $W^{\eta}$ as
$W^{\cdot}$ runs over the $\kappa$-filtered category of $\kappa$-small cell complexes with inclusions to 
$C_{\eta}(\beta , Z_{\cdot})$. 
The category of objects $W^{\cdot}$ under $Y^{\cdot}$ is cofinal in here, so we can view this colimit
as being over factorizations $Y^{\cdot}\rightarrt W^{\cdot}\rightarrt C_{\eta}(\beta , Z_{\cdot})$.

As, on the other hand, the category of $\kappa$-presentable objects mapping to 
$Z_{\eta}$ is $\kappa$-filtered, there is a factorization $Y^{\eta}\rightarrt T\rightarrt Z_{\eta}$
with $T$ being $\kappa$-presentable,
such that all of the above maps $u_{\eta}(a)$ compose into the same map $U_{\eta}\rightarrt T$. 
Because $Z_{\eta}$ is a $\kappa$-filtered colimit of the $W^{\eta}$, the $\kappa$-presentable
property of $T$ tells us that there is a factorization $T\rightarrt W^{\eta}\rightarrt Z_{\eta}$ 
for some $Y^{\cdot}\rightarrt W^{\cdot}\rightarrt C_{\eta}(\beta , Z_{\cdot})$.
Thus, all of the composed maps
$$
U_{\eta}\rightarrt^{u_{\eta}(a)} Y^{\eta}\rightarrt W^{\eta}
$$
are the same. We may now use this unique map as attaching map to add on the last cell,
to create a cell complex $\tilde{W}^{\cdot}$
with an inclusion to $(\beta , Z_{\cdot})$, whose cutoff at $\eta$ is 
$C_{\eta}(\tilde{W}^{\cdot}) = W^{\cdot}$. The identity of all of the composed attaching maps provides
us with inclusions of cell complexes $X^{\cdot}(a)\rightarrt \tilde{W}^{\cdot}$ for all $a\in A$
(the cell attachment provides these inclusions for $a\in A'$ and they are automatic for $a\in A''$). 
This finishes the proof of the $\kappa$-filtered property for 
$\Cell (\mM ; I)_{\kappa}/(\beta , Z_{\cdot})$. 

To complete the proof of the inductive step we just have to note that the colimit of the $\kappa$-small cell complex
inclusions, which we now know is $\kappa$-filtered and hence exists by Lemma \ref{cellfiltcolim}, is equal to the full $(\beta , Z_{\cdot})$.
For this, note first of all that the colimit of all of the $X^0$ will be $Z^0$ by the locally $\kappa$-presentable property of $\mM$. 
To conclude, it suffices to note that every $j\in \beta$ is a cell in at least one of the 
$\kappa$-small subcomplexes $X^{\cdot}$. If $\beta$ is a limit ordinal, apply the inductive hypothesis
saying that we know the statement of the theorem for any $\beta '<\beta$, and note that for any $j<\beta$ we can choose
$j<\beta '<\beta$. If $\beta = \eta +1$ is a successor ordinal, the inductive hypothesis gives the required statement for
any $j<\eta$ so we may assume $j=\eta$. Then, argue as above. We know that $C_{\eta}(\beta , Z_{\cdot})$ is
a $\kappa$-filtered colimit of its $\kappa$-small subcomplexes, and as before it follows that there will be a single
$\kappa$-small subcomplex $W^{\cdot}$ such that the attaching map $U_{\eta}\rightarrt Z_{\eta}$ factors 
through $U_{\eta}\rightarrt W^{\eta}$. The complex $\tilde{W}^{\cdot}$ obtained by attaching the last cell
$U_{\eta}\rightarrt V_{\eta}$ onto $W^{\cdot}$, satisfies the current requirement that $j=\eta$ be an element of
$c(\tilde{W}^{\cdot})$.  

This completes the inductive proof of the theorem. 
\end{proof}

We now restate the result of the theorem in its simplified form which will be used later
for the pseudo-generating set construction of a model category structure in Chapter \ref{modcat1}. 

\begin{corollary}
\label{AfiltcolimCor}
Suppose $f:X\rightarrt Y$ is a morphism in $\cell (I)$. Then we can express $f$ as a $\kappa$-filtered colimit
of arrows $f_i:X_i\rightarrt Y_i$, in particular $X= \colim _iX_i$ and $Y= \colim _iY_i$, such that
$X_i,Y_i$ are $\kappa$-presentable and $f_i\in \cell (I)$ are cell complexes of length $<\kappa$.
\end{corollary}
\begin{proof}
By definition of $\cell (I)$, there is a cell complex $(\beta , Z_{\cdot})\in \Cell (\mM ; I)$ such that
$f$ is equal to the  map $Z_0\rightarrt Z_{\beta}$. The theorem then says exactly that $f$ is a $\kappa$-filtered
colimit of morphisms coming from cell complexes in $\Cell (\mM ; I)_{\kappa}$.
\end{proof}

We now describe how to fine-tune a map between colimits so that it comes from
a levelwise map of diagrams. 
Let $\mM$ be a locally $\kappa$-presentable category. 
Suppose $\alpha$ is a $\kappa$-filtered category, and $F:\alpha \rightarrt \mM$ is a diagram such that each $F_i$ is $\kappa$-presentable.
Suppose $\beta$ is another $\kappa$-filtered category and $G:\beta \rightarrt \mM$ is another diagram. Suppose given a map
$$
f:\colim _{i\in \alpha} F_i \rightarrt \colim _{j\in \beta} G_j.
$$
Assume that the $F_i$ are $\kappa$-presentable objects.

\begin{lemma}
\label{colimap}
In the above situation, there are a $\kappa$-filtered category $\psi$ and cofinal functors
$p:\psi \rightarrt \alpha$ and $q:\psi \rightarrt \beta$ 
and a natural collection of maps $f_k:F_{p(k)}\rightarrt G_{q(k)}$ depending on $k\in \psi$,
i.e. a natural transformation $p^{\ast}F\rightarrt q^{\ast}(G)$,
such that the composition
$$
\colim _{i\in \alpha} F_i = \colim _{k\in \psi}F_{p(k)}\rightarrt^{\colim f_k}  \colim _{k\in \psi} G_{q(k)} = \colim _{j\in \beta}G_j
$$
is equal to $f$. 
\end{lemma}
\begin{proof}
Let $\psi$ be the category of triples $(i,j,u)$ where $i\in \alpha$, $j\in \beta$ and $u:F_i\rightarrt G_j$ is a morphism
such that the diagram
$$
\begin{diagram}
F_i & \rightarr & \colim _{i\in \alpha} F_i \\
\downarr & & \downarr \\
G_j & \rightarr & \colim _{j\in \beta} G_j
\end{diagram}
$$
commutes. Since the $F_i$ are $\kappa$-presentable, any $i$ is part of a triple $(i,j,u)\in \psi$. Any $j$ sufficiently far out in $\beta$ is also part of a triple, and 
indeed $\psi$ is $\kappa$-filtered, and the forgetful functors $\psi \rightarrt \alpha$ and
$\psi \rightarrt \beta$ are cofinal. The third variable $u$ provides the desired natural transformation $\{ f_k\}$ such that the diagram
$$
\begin{diagram}
\colim _{k\in \psi}F_{p(k)} & \rightarr^{\cong} & \colim _{i\in \alpha} F_i \\
\downarr^{\colim f_k} & & \downarr_{f}\\
\colim _{k\in \psi} G_{q(k)} & \rightarr^{\cong} & \colim _{j\in \beta} G_j
\end{diagram}
$$
commutes.
\end{proof}

The following proposition represents the idea that given an inclusion of cell complexes, we can rearrange things so that the
subcomplex comes first and then the rest of the complex is added on later. 

\begin{proposition}
\label{cellmove}
If $(\alpha , X_{\cdot})\hookrightarrow (\beta , Y_{\cdot})$ is an inclusion of cell complexes, then there is another inclusion of cell complexes
$(\delta , Z_{\cdot})\hookrightarrow (\beta , Y_{\cdot})$ 
with $Z_0 = X_{\alpha} \cup ^{X_0}Y_0$ and $Z_{\delta} = Y_{\beta}$, such that $\beta$
is the disjoint union of $c(X^{\cdot})$ and $c(Y^{\cdot})$. 
\end{proposition}
\begin{proof}
Left to the reader.
\end{proof}

\section{The small object argument}
\label{sec-soa}

Suppose $A\subset \Arr (\mM )$ is a class of morphisms. We say that a morphism $f:X\rightarrt Y$ 
{\em satisfies the right lifting property with respect to $A$} if, for any commutative diagram
$$
\begin{diagram}
U & \rightarr^{a} & X \\
\downarr^{u} && \downarr_{f} \\
V & \rightarr^{b} & Y
\end{diagram}
$$
such that $u\in A$, there exists a lifting $s:V\rightarrt X$ such that $fs=b$ and $su=a$. 
We say that $f$ {\em satisfies the left lifting property with respect to $A$} if,
for any commutative diagram
$$
\begin{diagram}
X & \rightarr^{a} & U \\
\downarr^{f} && \downarr_{v} \\
Y & \rightarr^{b} & V
\end{diagram}
$$
such that $v\in A$, there exists a lifting $s:Y\rightarrt U$ such that $vs=b$ and $sf=a$.

If $I\subset \Arr (\mM )$ is a  subset of morphisms, we have defined above $\cell (I) \subset \Arr (\mM )$ to
be the class of morphisms $f:X\rightarrt Y$ such that there exists an  $I$-cell complex $(\beta , X_{\cdot})$ with
$X_0=X$ and $X_{\beta} = Y$, with $f$ being the transition map from $X_0$ to $X_{\beta}$. 
The class $\cell (I)$ defined this way is clearly closed under
compositions. 

Let $\inj (I)\subset \Arr (\mM )$ be the class of maps which satisfy the right lifting property with respect to $\cell (I)$,
and let $\cof (I)\subset \Arr (\mM )$ be the class of maps which satisfy the left lifting property with respect to $\inj (I)$.
Note that $\cell (I)\subset \cof (I)$. 
The famous {\em small object argument} as it applies to locally presentable categories, can be summed up in the following theorem.

\begin{theorem}
\label{smallobject}
Suppose $\mM$ is a locally presentable category, and $I\subset \Arr (\mM )$ is a small subset of morphisms. 
Then:
\newline 
---any morphism $f:X\rightarrt Y$ admits a factorization $f=pg$ where $X\rightarrt^{g}Z\rightarrt^{p}Y$,
such that $g\in \cell (I)$ and $p\in \inj (I)$;
\newline
---one may choose a factorization which is functorial in $f$;
\newline
---the class $\cof (I)$ is closed under retracts and is equal to the class of morphisms
$f:X\rightarrt Y$ such that $Y$ is a retract of some $g:X\rightarrt Z$ in $\cell (I)$, in the category $X \backslash \mM$ of objects
under $\mM$; 
\newline
---the class $\inj (I)$ is also equal to the class of morphisms which satisfy the right lifting property with respect to $\cof (I)$.
\end{theorem}
\begin{proof}
Refer to the numerous discussions in the literature. 
\end{proof}

We record here some basic facts about lifting properties. 

\begin{lemma}
\label{rlpretract}
Suppose $\mM$ is a category and $\mA$ is a class of arrows in $\mM$. Then the class $\mF$ of morphisms
in $\mM$ which satisfy the right lifting property with respect to all morphisms of $\mA$, is closed under
retracts. Similarly, if  $\mB$ is a class of arrows then the class $\mG$ of morphisms
which satisfy the left lifting property with respect to all morphisms of $\mB$, is closed under retracts. 
\end{lemma}
\begin{proof}
Suppose $f:X\rightarrt Y$ is in $\mF$, and $g:A\rightarrt B$ is
a retract of $f$ with a retract diagram 
$$
\begin{diagram}
A & \rightarr^{i} & X & \rightarr^{r} & A \\
\downarr & & \downarr & & \downarr \\
B & \rightarr^{j} & Y & \rightarr^{s} & B .
\end{diagram}
$$
Suppose 
$$
\begin{diagram}
U & \rightarr^{a} & A  \\
\downarr^{u} & & \downarr_{g}  \\
V & \rightarr^{b} & B 
\end{diagram}
$$
is a diagram with $u\in \mA$. Compose with the left square of the retract diagram to get
$$
\begin{diagram}
U & \rightarr^{ia} & X  \\
\downarr^{u} & & \downarr_{f}  \\
V & \rightarr^{jb} & Y .
\end{diagram}
$$
The lifting property for $f$ says that there exists $t:V\rightarrt X$ with $tu=ia$ and $ft=jb$.
Composing with $r$ gives a lifting $rt:V\rightarrt A$ such that $rtu=ria=a$ and 
$grt=sft = sjb = b$, so $rt$ is a lifting for the original diagram to $g$. 

The proof for $\mG$ is similar. 
\end{proof}

\begin{corollary}
\label{injretracts}
Suppose $I$ is a set of morphisms in a locally presentable category $\mM$. Then the class $\inj (I)$ is 
closed under retracts. 
\end{corollary}
\begin{proof}
The class $\inj (I)$ is defined by the right lifting property with respect to $\cell (I)$.
\end{proof}

\begin{lemma}
\label{rlplims}
Suppose $\mM$ is a complete and cocomplete category, and $\mA$ is a class of arrows in $\mM$. Then the class $\mF$ of morphisms
in $\mM$ which satisfy the right lifting property with respect to all morphisms of $\mA$, is closed under
fiber products and sequential inverse limits. Similarly, if  $\mB$ is a class of arrows then the class $\mG$ of morphisms
which satisfy the left lifting property with respect to all morphisms of $\mB$, is closed under pushouts and transfinite composition. 
\end{lemma}
\begin{proof}
Left to the reader. 
\end{proof}

\section{Injective cofibrations in diagram categories}

We turn now to two of the main results from the appendix to Lurie's \cite{LurieTopos}; these are notably refered to
by Barwick \cite{Barwick} who gives a general discussion of the injective model structure on section categories of
left Quillen presheaves. We need these results in order to construct injective model categories in what follows,
and we need to understand something about the proof in order to apply it to the case of unital diagram categories
where we require that some values are equal to $\ast$. The novice reader is invited to skip this section,
refering to \cite{LurieTopos} for the results, and imagining their extension to the unital case. Alternatively,
the arguments can be made more concrete in the case (due originally to Heller \cite{Heller})
where $\mM$ is a presheaf category and the cofibrations are monomorphisms. 

Fix a locally $\kappa$-presentable  category $\mM$ and a subset $I$ of morphisms. Let $\cof (I)$ denote the class of $I$-cofibrations,
that is morphisms which are retracts of morphisms in $\cell (I)$. Let $\cof (I)_{\kappa}$ denote a set of representatives for 
the isomorphism classes of morphisms 
$A\rightarrt B$ in $\cof (I)$ such that $A$ and $B$ are $\kappa$-presentable. Lurie's first theorem \cite[A.1.5.12]{LurieTopos} is:

\begin{theorem}
\label{lurie-thm1}
In the above situation, $\cell (\cof (I)_{\kappa}) = \cof (\cof (I)_{\kappa}) = \cof (I)$, that is to say
any $I$-cofibration can be expressed as a cell complex whose attaching maps are taken from the set $\cof (I)_{\kappa}$.
\end{theorem}
\begin{proof}
The inclusions $\cell (\cof (I)_{\kappa}) \subset \cof (\cof (I)_{\kappa}) \subset \cof (I)$
are immediate, we need to show that $\cof (I)\subset \cell (\cof (I)_{\kappa})$.
Suppose $w:A\rightarrt V$ is in $\cof (I)$. Using the small object argument, we can choose a diagram
$$
\begin{diagram}
& & B \\
& \ruTeXto(2,2)^f & \downarr^p \uparr_s \\
A & \rightarr & V
\end{diagram}
$$
with a map $A\rightarrt^f B$ in $\cell (I)$, a projection $B\rightarrt^p V$ in $\inj (I)$,
and a section $V\rightarrt^s B$ compatible with the map from $A$, such that $ps = 1_V$. 
Let $\pi _B:= s p : B\rightarrt B$ be the idempotent $\pi _B^2 = \pi_B$, compatible with the identity map of $A$.

For the transfinite induction step,
it will be convenient to consider more generally the case when $\pi_B$ is compatible with
$\pi_A$ which can be a nontrivial idempotent on $A$.

The idea is to find $\kappa$-presentable cell complexes $f_i:A_i\rightarrt B_i$ included into $f$, together with definitions of $\pi_{B,i}$ and $\pi _{A,i}$
on the source and target. Note that, given an idempotent $\pi _{B,i}:B_i\rightarrt B_i$ we can let $V_i$ be the colimit of the
diagram 
$$
B_i\rightarrt^{\pi _{B,i}} B_i\rightarrt^{\pi _{B,i}} B_i\rightarrt^{\pi _{B,i}} \ldots ,
$$
then we have a projection $p_i:B_i\rightarrt V_i$ and the system of maps $\pi _i$ on the colimit gives $s_i:V_i\rightarrt B_i$
with $p_is_i=1_{U_i}$. Similarly let $U_i$ be the image of $\pi _{A,i}$ defined in the
same way, so if $A_i\rightarrt B_i$ are in $\cell (I)$ then their retracts $U_i\rightarrt V_i$ are in $\cof (I)$ since $\cof (I)$
is closed under retracts.

In order to construct the $f_i$, proceed as follows. Let 
$$
\xi := \Cell (\mM ; I)_{\kappa}/f
$$
denote the category of inclusions of cell complexes over $f$. By Theorem \ref{cellover}, $\xi$ is $\kappa$-filtered.
For $g\in \xi$ let $h(g): X(g)\rightarrt Y(g)$ denote the corresponding cell complex, mapping to $f$ by an inclusion of
cell complexes. We have $A=\colim _{g\in \xi} X(g)$, $B=\colim _{g\in \xi} Y(g)$, and $f$ is the colimit of the $g\in \xi$.
Let $x(g):X(g)\rightarrt A$ and $y(g):Y(g)\rightarrt B$ denote the maps to the colimits.
On the other hand, $X(g)$ and $Y(g)$ are $\kappa$-presentable. In particular, the map 
$$
\pi_B y(g): Y(g)\rightarrt B
$$
has to factor through a map $\pi_B (g): Y(g)\rightarrt Y(n(g))$
for some $g\rightarrt n(g)\in \xi$ which is a function of $g$. Furthermore for any $v:g\rightarrt h$ in $\xi$ there is
$n(g), n(h)\rightarrt n'(v)$ such that the projections of $\pi_B (g)$ and $\pi (h)$ into $Y(n'(v))$ are the compatible along $Y(v)$,
and finally there is $n(n(g))\rightarrt n''(g)$ such that $\pi_B (n(g))\pi (g) = \pi_B (g)$ after projecting into $Y(n''(g))$. 
Choose similarly $\pi _A(g)$ and assume that the same choices work for $\pi _A(g)$, and furthermore that $\pi_B (g)h(g)=\pi _A(g)$ after projection
to $Y(n(g))$. 

A {\em good filtered subcategory} $\xi _i\subset \xi$ is a full subcategory with $<\kappa$ objects, filtered
(note however that it would be too small to be $\kappa$-filtered), and such that for any $g\in \xi_i$ the elements $n(g)$, $n''(g)$
and their arrows are also in $\xi _i$; and for an arrow $v$ in $\xi _i$ the object $n'(v)$ with its arrows is also in $\xi _i$.
If $\xi _i$ is a good filtered subcategory, then 
$$
B_i:= \colim _{g\in \xi _i}Y(g)
$$
has an endomorphism $\pi _{B,i}:B_i\rightarrt B_i$ defined by using the $\pi_B (g)$.
Similarly $\pi _{A,i}:A_i\rightarrt A_i$ is defined using the $\pi _A(g)$; these are compatible 
and $\pi _i^2= \pi _i$. By Proposition \ref{cellfiltcolim}, 
$f_i$ is in $\cell (I)$, so the images $U_i$ (resp. $V_i$) of the idempotents $\pi _{B,i}$ (resp. $\pi _{A,i}$) as defined above,
have an $I$-cofibration $w_i: U_i\rightarrt V_i$. Note that $w_i\in \cof (I)_{\kappa}$.

Recall that by a {\em cell} of $f$ we mean an element of the ordinal indexing the cell attachements.
An inclusion of cell complexes generates a corresponding subset of cells, although in the general case
we are currently considering, specification of the subset of cells is not sufficient to specify the
inclusion of cell complexes, because our attaching maps are not necessarily monomorphisms. 

For any given cell of $f$, we can choose a good filtered subcategory such that the inclusion of cell complexes from $B_i$ to $B$,
contains the given cell. Hence $A\cup ^{A_i}B_i$ is a pushout of $A$ along an element $w_i\in \cof (I)_{\kappa}$ containing the given cell.

Start now with our given situation $w:A\rightarrt V$ being a retract of $f:A\rightarrt B$. 
We construct by transfinite induction, a sequence of inclusions of cell complexes $A\rightarrt C_j\rightarrt B$
together with definitions of the idempotent $\pi _{C,j}$ on $C_j$, compatible with $1_A$ and $\pi _B$. Start with $C_0=1$. Assume $C_j$ is chosen for $j<j_0$.
If $j_0$ is a limit ordinal then let $C_j:= \colim _{j<j_0}C_j$. If not, $j_0=k+1$ and we are in the general situation envisioned above:
$C_k$ has its idempotent $\pi _{C,k}$ compatible with $\pi _B$ via the cell complex $b_k:C_k\rightarrt B$. Fix the smallest cell of $B$ not
contained in the subset of cells of $C_k$, and choose according to the previous procedure applied to the map $b_k$, a $\kappa$-presentable
cell complex $w_j:A_j\rightarrt B_j$ with 
an inclusion of cell complexes from $w_i$ to $(b_k:C_k\rightarrt B)$, and with
compatible idempotents $\pi _{A,j}$ and $\pi _{B,j}$. Set 
$$
C_j= C_{k+1}:= C_k\cup ^{A_j}B_j,
$$
with its induced idempotent which will be called $\pi _{C,j}$. 

The process stops when there are no longer any cells of $B$ but not in $C_k$, that is to say $C_k=B$. 
Now let $W_0:= A$ and let $W_j$ be the image of the idempotent $\pi _{C,j}$ defined as 
$$
W_j:= \colim \left( C_j\rightarrt^{\pi _{C,j}} C_j\rightarrt^{\pi _{C,j}} \ldots \right) .
$$
If $U_j$ and $V_j$ denote the images of $\pi _{A,j}$ and $\pi _{B,j}$ respectively, 
we have $U_j\rightarrt V_j \in \cof (I)_{\kappa}$. Note that
$$
W_{j} = W_k \cup ^{U_{j}}V_{j}
$$
by commutation of colimits.  Similarly if $j$ is a limit ordinal then $W_j=\colim _{i<j}W_i$. 
The sequence $\{ W_j\}$ is an expresssion of $f:A\rightarrt V$ as an element of $\cell (\cof (I)_{\kappa})$. 
\end{proof}

Lurie's second theorem is the application to diagram categories. Keeping the previous notations,
suppose $\Phi$ is a $\kappa$-small category (i.e. its object and morphism sets
have cardinality $<\kappa$). Within the category $\diag (\Phi , \mM )$ say that a map $A\rightarrt B$ is an {\em injective cofibration}
if $A(x)\rightarrt B(x)$ is in $\cof (I)$ for any $x\in\Ob (\Phi )$. This class of morphisms is denoted $\diag (\Phi ,  \cof (I))$. Let
$\diag (\Phi ,  \cof (I))_{\kappa}$ denote a set of representatives for isomorphism classes of injective cofibrations between 
$\kappa$-presentable objects, which is also the set of morphisms $A\rightarrt B$ such that each $A(x)\rightarrt B(x)$ is in
$\cof (I)_{\kappa}$. The following statement is \cite[A.2.8.3]{LurieTopos}; the argument is similar to the previous one.

\begin{theorem}
\label{lurie-thm2}
In the above situation, 
$$
\diag (\Phi ,  \cof (I)) = \cof (\diag (\Phi ,  \cof (I))_{\kappa}),
$$
i.e. the set $\diag (\Phi ,  \cof (I))_{\kappa}$ generates the class of injective cofibrations. 
\end{theorem}
\begin{proof}
In light of the previous theorem, we may assume that $\cof (I)=\cell (I)$. Suppose 
$f:A\rightarrt B$ is in $\diag (\Phi ,  \cof (I))$, then each $f(x)$ can be
given a structure of $I$-cell complex. Make such a choice for each $x\in \Ob (\Phi )$, although these choices are not compatible
with the structure of diagram. 

Let $\xi (x)$ be the $\kappa$-filtered category of $\kappa$-presentable inclusions of cell complexes into $f(x)$
(see Theorem \ref{cellover}). 
For $i\in \xi (x)$ we have a $\kappa$-presentable cell complex $f(x,i):A(x,i)\rightarrt B(x,i)$
with an inclusion to $f$, and $B(x)=\colim _{i\in \xi (x)}B(x,i)$. If $i\in \xi (x)$ and
$\varphi : x\rightarrt y$ is
an arrow in $\Phi$ then there is $j\in \xi (y)$ plus a lifting to a map $B(x,i)\rightarrt B(y,j)$ compatible with $B(\varphi )$.
We can assume that the same $j$ works for any $\varphi$; for $x,y,z$ and $i\in \xi (x)$ we can choose $k\in \xi (z)$ such that
the compositions of any transition maps from $B(x,i)$ to $B(y,j)$ then to $B(z,k)$ satisfy the product rule. Continuing
in this way, we obtain a collection of operations defined on $\bigcup _{x}\xi (x)$ together with lifting data, such 
that if $\xi '(x)\subset \xi (x)$ is a collection of filtered subcategories preserved by these operations, then
the collection of $B'(x):= \colim _{i\in\xi '(x)}B_i$ is given a structure of diagram $B'$.  Similarly for $A'$
and the map $A'\rightarrt B'$ is an injective cofibration. 
We can choose $\xi '(x)$ with cardinality $<\kappa$, so $A'\rightarrt B'$ is a
$\kappa$-presentable injective cofibration, and in such a way that $B'(x)$ contains any given cell of some $B(x_0)$. 

Now applying the same inductive argument as at the end of the previous theorem, which we don't repeat, gives an expression of $A\rightarrt B$ as
a transfinite composition of pushouts along such $A'\rightarrt B'$.
\end{proof}

The version we really need later is a variant in which the diagram category is replaced by the {\em unital diagram category}
denoted $\diag (\Phi /\Phi _0, \mM )$, discussed in more detail in Chapter \ref{algtheor1}. 
Here $\Phi _0$ is a subset of objects of $\Phi$ and $A\in \diag (\Phi /\Phi _0, \mM )$
means that $A:\Phi \rightarrt \mM$ is a diagram such that $A(x)\cong \ast$ is the coinitial object of $\mM$, for all $x\in \Phi _0$.
As before if $I$ is a fixed set of maps in $\mM$, a map $A\rightarrt B$ in this category is an {\em injective cofibration}
if each $A(x)\rightarrt B(x)$ is in $\cof (I)$. Note that isomorphisms are always contained in $\cof (I)$ so this condition
holds automatically whenever $x\in \Phi _0$. 

\begin{theorem}
\label{lurie-unital}
The class $\diag (\Phi /\Phi _0 ,  \cof (I))$ of injective cofibrations in the unital diagram category,
equals $\cof (\diag (\Phi /\Phi _0,  \cof (I))_{\kappa})$, i.e. it is
generated by the set $\diag (\Phi ,  \cof (I))_{\kappa}$ of representatives for isomorphism classes
of $\kappa$-presentable injective cofibrations. 
\end{theorem}
\begin{proof}
Follow the same procedure as in the previous theorem, noting that all coproducts involved are taken over connected
categories, and by Lemma \ref{connstar} the unitality condition is preserved by such coproducts. The $\kappa$-presentable
objects of $\diag (\Phi /\Phi _0 ,\mM )$ are the unital diagrams $A$ which are $\kappa$-presentable as diagrams, or
equivalently such that each $A(x)$ is $\kappa$-presentable in $\mM$.
So, the process described in the proof of the previous theorem leads to an expression of an arbitrary
injective cofibration as a transfinite composition of pushouts along $\kappa$-presentable ones, and all of the
objects intervening here remain unital.
\end{proof}


\chapter{Model categories}
\label{modcat1}

In this chapter, we recall some of the basic elements of Quillen's theory of model categories
and modern variants. We consider model structures on diagram categories, Quillen functors, and
left Bousfield localization. One of the goals is to understand a quick version of left Bousfield localization
which is easier to understand than the general version, but which needs a restrictive collection of hypotheses.
It turns out that these hypotheses will  hold in the cases we need. 

Along the way, we also discuss some abstract notions in category theory which are useful for formulating the 
small object argument. 

As pointed out in the previous chapter, the motivation for introducing Quillen's theory into the world of $n$-categories
is the fact that there is a localization problem at the heart of the route towards construction of the $n+1$-category $nCAT$ of
all $n$-categories. Grothendieck was undoubtedly aware of this connection at least on an intuitive level, because ``Pursuing Stacks'' started
out as a series of letters to Quillen. 

\section{Quillen model categories}

Start by recalling Quillen's definition of a closed model category. This is a category $\mM$ provided with
three classes of morphisms, the ``weak equivalences'', the ``cofibrations'', and the ``fibrations''.
The intersection of the classes of cofibrations and weak equivalences is called the class of ``trivial cofibrations'', 
and similarly 
the intersection of the classes of fibrations and weak equivalences is called the class of ``trivial fibrations''.
Quillen asks that these should satisfy the following axioms.  
\newline
(CM1)---The category $\mM$ should be closed under finite limits and colimits. Following modern tradition, we really require that it be closed
under all small limits and colimits. 
\newline
(CM2)---The class of weak equivalences satisfies {\em 3 for 2}: given a composable sequence of arrows in $\mM$
$$
X \rightarrt^{f} Y\rightarrt^{g} Z
$$
if any two of $f$, $g$ and $gf$ are weak equivalences, then so is the third one. 
\newline
(CM3)---The classes of cofibrations, fibrations, and weak equivalences should be closed under retracts: given a diagram
$$
\begin{diagram}
A & \rightarr^{i} & X & \rightarr^{r} & A \\
\downarr & & \downarr & & \downarr \\
B & \rightarr^{j} & Y & \rightarr^{s} & B
\end{diagram}
$$
such that $ri=1_A$ and $sj=1_B$, such that the two outer downward arrows are the same, if the middle
arrow $g:X\rightarrt Y$ is a cofibration (resp. fibration, weak equivalence) then the outer arrow $f:A\rightarrt B$ is also
a cofibration (resp. fibration, weak equivalence). 
\newline
(CM4)---Cofibrations satisfy the left lifting property with respect to trivial fibrations, and
trivial cofibrations satisfy the left lifting property with respect to fibrations. 
\newline
(CM5)---If $f:X\rightarrt Y$ is any morphism, then there exist factorizations (i) and (ii)
of $f$ as the composition $X\rightarrt^{g} Z \rightarrt^{p} Y$ such that
\newline
(i) \,\, $g$ is a cofibration and $p$ is a trivial fibration;  
\newline
(ii) \,\, $g$ is a trivial cofibration and $p$ is a fibration. 

It follows from these axioms that any two of the classes of cofibrations, fibrations
and weak equivalences, determine the third. 
For example a morphism is a fibration (resp. trivial fibration)
if and only if it satisfies right lifting
with respect to any trivial cofibration (resp. cofibration), and dually. 
A morphism is a weak equivalence if and only if it factors as a composition of
a trivial fibration and a trivial cofibration, these classes being determined from the
classes of cofibrations and fibrations respectively by the lifting property. 
Each of the three classes contains any isomorphism.

We point out, as was done in \cite{svk}, that the diagram included in \cite{Quillen} for the definition of ``retract'' is visibly
wrong, so his notion of ``retract'' is not well defined. There could be two reasonable interpretations of
this condition. For condition (CM2) we have adopted the weak interpretation. The stronger interpretation would have the arrows
on the bottom row going in the opposite direction. If $f$ is a strong retract of $g$ then it is also a weak retract of $g$.
Hence, closure under retracts as we require in (CM2) also implies closure under strong retracts. This choice coincides with what was said in Dwyer-Spalinski \cite[2.6]{DwyerSpalinski}. 
Similarly, 
Hinich \cite{HinichHAHA} uses the retract condition stated as we have done above.

\begin{lemma}
\label{determine}
If $\mM$ is a closed model category i.e. satisfies (CM1)--(CM5), then the classes
of trivial cofibrations and fibrations determine each other by the lifting property; and similarly
the classes of cofibrations and trivial fibrations determine each other. In other words, a morphism
is a fibration (resp. trivial fibration) if and only if it satisfies the right lifting property
with respect to all trivial cofibrations (resp. all cofibrations). And a morphism is a 
cofibration (resp. trivial cofibration) if and only if it satisfies the left lifting property
with respect to all trivial fibrations (resp. all fibrations).
\end{lemma}
\begin{proof}
There are four things to prove. Consider for example the fact that a morphism is a fibration
if and only if it satisfies the right lifting property with respect to trivial cofibrations; the
other three arguments are identical. If $f$ is a fibration then by (CM4) it satisfies the
right lifting property with respect to all trivial cofibrations. Suppose on the other hand that $f$ is
a morphism which satisfies the
right lifting property with respect to all trivial cofibrations. Using (CM5), factor $f=pg$,
$$
X \rightarrt^{g} Z\rightarrt^{p} Y
$$ 
where
$g$ is a trivial cofibration and $p$ is a fibration. Apply the right lifting property being assumed for $f$,
to the diagram 
$$
\begin{diagram}
X & \rightarr^{=} & X \\
\downarr^{g} & & \downarr_{f} \\
Z & \rightarr^{p} & Y 
\end{diagram}
$$
to get a morphism $s:Z\rightarrt X$ such that $sg=1_X$ and $fs=p$. Putting $s$ into the diagram
$$
\begin{diagram}
X & \rightarr^{g} & Z & \rightarr^{s} & X \\
\downarr & & \downarr & & \downarr \\
Y & \rightarr^{=} & Y & \rightarr^{=} & Y
\end{diagram}
$$
now says that $f$ is a retract of the fibration $p:Z\rightarrt Y$, so by (CM3) $f$ is a fibration.
\end{proof}

\section{Cofibrantly generated model categories}
\label{sec-cofibgen}

There is an important class of model categories which are particularly easy to work with, and 
which contains many if not most of the examples currently considered as important. Useful examples
of model categories which are not cofibrantly generated, will often be Quillen equivalent to cofibrantly generated ones.
On the other hand, this notion is very helpful for a number of the operations we need, such as taking the model category
of $\mM$-diagrams, and left Bousfield localization. The reader is referred to Hirschhorn's excellent book 
\cite{Hirschhorn} for a full explanation of everything concerning cofibrantly generated model categories. 

Suppose $J\subset \Mor (\mM )$ is a small subset of morphisms, 
then we define classes of morphisms  ${\cell}(J)$, ${\inj}(J)$ and ${\cof}(J)$,
see Sections \ref{inclusioncell} and \ref{sec-soa}. 

Recall that ${\cof}(J)$ consists of arrows $X\rightarrt Y$ which are retracts, in the category of objects under $X$,
of elements of ${\cell}(J)$. This condition means that there should exist $X\rightarrt Z$ in ${\cell}(J)$
and maps $Y\rightarrt Z\rightarrt Y$ compatible with the maps from $X$ and composing to the identity of $Y$.

Recall that a $\univa$ model category $\mM$ is {\em cofibrantly generated} if it is closed under 
$\univa$-small limits and colimits, and if there are $\univa$-small subsets of arrows $I,J\subset {\bf Arr}(\mM )$
satisfying the following properties: 
\newline
(CG1) \, the sources and targets of arrows in $I$ and $J$ are small in $\mM$;
\newline
(CG2a) \, the arrows in $I$ are cofibrations;
\newline
(CG2b) \, the trivial fibrations of $\mM$ are the morphisms satisfying lifting with respect to arrows in $I$;
\newline
(CG3a) \, the arrows in $J$ are trivial cofibrations;
\newline
(CG3b) \, the fibrations of $\mM$ are the morphisms satisfying lifting with respect to arrows in $J$.

As a matter of notation, we say that {\em $(\mM, I, J)$ is a cofibrantly generated model category}, if 
$\mM$ is a model category and $I$ and $J$ are subsets of arrows satisfying the above axioms. 

This notation
is convenient in that $I$ and $J$ determine the model category structure, indeed $I$ determines the class of 
trivial fibrations by (CG2b) and $J$ determines the class of fibrations by (CG3b). These in turn determine
respectively the classes of cofibrations and trivial cofibrations by the saturated lifting properties.
Then the class of weak equivalences is determined by either one of the factorization properties.

Given a triple $(\mM , I, J)$ consisting of a category admitting $\univa$-small limits and colimits, and
two subsets of arrows, we can try to define a model structure following the recipe of the previous paragraph. 
If these classes of (trivial) cofibrations, (trivial) fibrations and weak equivalences do in fact form a closed model structure,
and if (CG1) is satisfied, then $(\mM, I, J)$ is a cofibrantly generated model category.

See \cite{Hirschhorn} \cite{GoerssJardine} and other references, for discussions of various recognition properties
telling when a triple $(\mM , I, J)$ yields a cofibrantly generated model category. 
For example, the following statement will be useful. 

\begin{proposition}
\label{cofibgenrec}
If $\mM$ is a locally presentable category and we start with the three classes of cofibrations, fibrations and
weak equivalences; then define trivial cofibrations as the intersection of cofibrations
and weak equivalences and similarly for trivial fibrations; if $I$ and $J$ satisfy
properties (CG1)--(CG3b); and if furthermore we know that weak equivalences are
closed under retracts and satisfy 3 for 2, then it is a cofibrantly generated model category.
\end{proposition}
\begin{proof}
Indeed, (CM5) comes from the small object argument, (CM4) comes from the hypotheses 
(CG2b) and (CG3b), (CM3) is supposed for weak equivalences and 
follows for $\inj (J)$ and $\cof (I)$, and (CM1) and (CM2) are supposed. 
\end{proof}


\section{Combinatorial and tractable model categories}

It is most convenient to combine the notion of cofibrantly generated model category with a good category-theoretical condition on
$\mM$ which guarantees the small object argument. As was observed by J. Smith and
reported by D. Dugger in several papers for example
\cite{DuggerCombinatorial}, the appropriate condition is to require that
$\mM$ be locally presentable. 

A {\em combinatorial model category} 
is a cofibrantly generated model category which is also locally presentable.
In this case all objects are small in $\mM$ so (CG1) is automatic.

Barwick then refined this by requiring that the domains of
generating cofibrations and trivial cofibrations also be cofibrant:
a combinatorial model category is {\em tractable} if in addition for the given generating sets $I$ and $J$:
\newline
(TR)---the sources of arrows in $I$ and $J$ are cofibrant.

As we shall discuss below, diagrams with values in a combinatorial model category
have injective and projective model structures. If the target category is tractable,
then the projective model structure on diagrams is also tractable; it isn't clear
whether this is known in the case of the injective structure.

\section{Homotopy liftings and extensions}
\label{sec-htyliftext}

For $X\in \mM$ a {\em cylinder object} is a diagram 
$$
X\cup ^{\emptyset} X \rightarrt^{i_0\cup i_1} C\rightarrt ^p X
$$
such that $i_0\cup i_1$ is a cofibration and $p$ is a trivial fibration.
The existence is guaranteed by Axiom (CM5)(i). Similarly, if $X\rightarrt Y$
is a cofibration there exists a {\em relative cylinder object}
which is a diagram 
$$
Y\cup ^{X} Y \rightarrt^{i_0\cup i_1} C\rightarrt ^p  Y
$$
again with $i_0\cup i_1$ a cofibration and $p$ a trivial fibration.
Quillen shows that
if $X$ is cofibrant and $A$ fibrant, then two maps $f_0,f_1:X\rightarrt A$ are
homotopic, that is project to the same map in $\Ho (\mM )$, if and only if
there is a cylinder object and a map
(or for any cylinder object there is a map) $C\rightarrt A$
inducing the two given maps. Given a cofibration $X\rightarrt Y$ and
two maps $f_0,f_1:Y\rightarrt A$ which agree when restricted to $X$,
we say they are {\em homotopic relative to $X$} if there is a relative cylinder
object and a map $C\rightarrt A$ inducing the two maps.

Weak equivalences between fibrant objects can be characterized by a homotopy
lifting property, being careful to look at homotopies relative to the 
subobject of the cofibration. This is a classical fact from the theory of model categories
but, as it is a technical step needed in
Chapter \ref{algtheor1}, the proof is included here for completeness.

\begin{lemma}
\label{homotopylifting}
Suppose $g:X\rightarrt Y$ is a map between fibrant objects in a tractable
model category $\mM$, such that for any diagram 
$$
\begin{diagram}
U& \rightarr ^u & X \\
\downarr ^f & & \downarr _g\\
V & \rightarr ^v & Y
\end{diagram}
$$
where $f$ is a generating cofibration, there exists 
a lifting $r:V\rightarrt X$ such that $rf=u$, together with
a relative cylinder object 
$$
V\cup ^{U} V\rightarrt^{j_0\cup j_1} IV\rightarrt ^q V
$$
with $j_0\cup j_1$ a cofibration and $q$ a weak equivalence,
and a map $h:IV\rightarrt Y$ restricting to $gu=vf$ on $U$,
such that $hj_0=gr$ and $hj_1=v$.
Then $g$ is a weak equivalence. 
\end{lemma}
\begin{proof}
Factor $X\rightarrt ^{k}X'\rightarrt^{g'}Y$ where $k$ is a trivial cofibration and
$g'$ is a fibration. Since $X$ is assumed fibrant
there is a retraction $s$ from $X'$ to $X$ with $sk=1_X$. 
We show that $g'$ is a trivial fibration by the  lifting property.
Suppose given a diagram
$$
\begin{diagram}
U& \rightarr ^{u'} & X' \\
\downarr ^f & & \downarr _{g'}\\
V & \rightarr ^{v'} & Y
\end{diagram}
$$
with $f$ a generating cofibration. 
Put $u:=su':U\rightarrt X$; there is a homotopy between $ku$ and $u'$,
given by a cylinder object
$$
U\cup ^{\emptyset} U\rightarrt^{i_0\cup i_1} IU\rightarrt ^pU
$$
with a map $z:IU\rightarrt X'$ with $zi_0=ku$ and $zi_1=u'$.
Choose a compatible cylinder object for $V$,
fitting into a diagram
$$
\begin{diagram}
U\cup ^{\emptyset} U&\rightarrt^{i_0\cup i_1}& IU&\rightarrt ^p&U\\
\downarr^f & & \downarr^{If} & & \downarr_f \\
V\cup ^{\emptyset} V&\rightarrt^{l_0\cup l_1} &I'V&\rightarrt ^q& V .
\end{diagram}
$$
Extend the map $IU\cup ^{i_1,U,f}V\rightarrt Y$ given by $g'z\cup v'$,
to a map $t:I'V\rightarrt Y$ (this is an extension along a trivial cofibration, for
maps to the fibrant object $Y$). Restricting along $l_0$ now gives a map $v:V\rightarrt Y$
such that $vf=g'ku=gu$. We obtain a diagram 
$$
\begin{diagram}
U& \rightarr ^{u} & X \\
\downarr ^f & & \downarr _{g}\\
V & \rightarr ^{v} & Y.
\end{diagram}
$$
By hypothesis there is a lifting $r:V\rightarrt X$ such that $rf=u$ and
$gr$ is homotopic to $v$ relative to $U$. This gives a lifting $kr:V\rightarrt X'$
such that $krf=ku$ and $g'kr$ is homotopic to $v$ relative to $U$. 
Putting this back together with the previous homotopy, we get the top map 
in the diagram
$$
\begin{diagram}
IU\cup ^{i_0,U,f}V & \rightarr ^{z\cup kr} & X'\\
\downarr ^{If\cup l_0} & & \downarr _{g'}\\
I'V & \rightarr ^t & Y 
\end{diagram}
$$
where the left map is a trivial cofibration and $g'$ is a fibration. 
This diagram doesn't commute, however by the hypothesis of the lemma 
(and using the notations from there)
there is a homotopy relative to $IU$ making
it commute. Adding this on and shifting the map from the
component $V$ on the upper left, to the other side of the new cylinder
object $IV$, gives the commutative diagram 
$$
\begin{diagram}
IU\cup ^{i_0,U,f}V & \rightarr ^{z\cup kr} & X'\\
\downarr ^{If\cup j_1} & & \downarr _{g'}\\
I'V \cup ^{l_0,V,j_0}IV& \rightarr ^{t\cup h} & Y 
\end{diagram}
$$
where the left map remains a trivial cofibration.
Hence
there exists a lifting $I'V \cup ^{l_0,V,j_0}IV\rightarrt X'$ which, when restricted along $l_1$ gives 
the desired lifting to show that $g'$ is a trivial fibration. This completes the proof.
\end{proof}

A similar consideration holds for extensions. 

\begin{lemma}
\label{htyext}
Suppose given two cofibrations $X\rightarrt ^a Y\rightarrt ^b Z$ and 
maps $Y\rightarrt^f A$ and $Z\rightarrt ^gA$ to a fibrant object $A$,
such that $fa = gba$. Suppose that $gb$ is homotopic to $f$ relative to $X$.
Then there exists a map $Z\rightarrt ^{g'}A$ such that $g'b=f$. 
\end{lemma}
\begin{proof}
Let $C$ be a relative cylinder object for $X\rightarrt Y$.
Choose a factorization 
$$
Z \cup ^{Y,i_0} C \cup ^{Y,i_1} Z \rightarrt^d D \rightarrt^q Z
$$
of a trivial fibration composed with a cofibration. Then $D$
is also a relative cylinder object for $X\rightarrt Z$ and we have
the compatibility diagram 
$$
\begin{diagram}
Y\cup ^{X} Y &\rightarr^{i_0\cup i_1} &C&\rightarr ^p & Y \\
\downarr^{b\cup b} & & \downarr_c & & \downarr_b \\
Z\cup ^{X} Z &\rightarr^{j_0\cup j_1} &D&\rightarr ^q & Z.
\end{diagram}
$$
Choose a map $C\rightarrt^hA$ such that $hi_0=f$ and $hi_1=gb$. 
We get a map 
$$
C \cup ^{Y,i_1}Z \rightarrt ^{h\cup g} A
$$
but the cofibration $Y\rightarrt ^{i_1}C$ is a weak equivalence, so
the pushout along $i_1$ 
is a trivial cofibration $Z\rightarrt C \cup ^{Y,i_1}Z$. Using the
fact that $d$ was chosen to be a cofibration at the start, it follows that 
$$
C \cup ^{Y,i_1}Z\rightarrt^{d'} D
$$
is a trivial cofibration. The fibrant condition for $A$ now allows us to extend the
previous map $h\cup g$ along $d'$, so we get a map $D\rightarrt^{h'} A$ extending $h$.
The restriction $g':= h'j_0: Z\rightarrt A$ provides the required map. 
\end{proof}

\section{Left properness}

Recall that a model category $\mM$ is {\em left proper} if, in any pushout square
$$
\begin{diagram}
X & \rightarr & Y \\
\downarr  & & \downarr \\
Z & \rightarr & W
\end{diagram}
$$
such that $X\rightarrt Y$ is a cofibration and $X\rightarrt Z$ is a weak equivalence, then $Y\rightarrt W$ is also a weak equivalence. 

\begin{lemma}
\label{leftpropinvariance}
Suppose $\mM$ is a left proper model category, and suppose we are given a diagram 
$$
\begin{diagram}
X & \rightarr & Y \\
\downarr  & \rdTeXto & \downarr \\
Z &  & V
\end{diagram}
$$
such that  $Y\rightarrt V$ is a weak equivalence. 
Suppose either that $X\rightarrt Z$ is a cofibration, or that 
both maps $X\rightarrt Y$ and $X\rightarrt V$ are cofibrations. Then
the map 
$$
Z\cup ^X Y\rightarrt Z\cup ^X V
$$
is a weak equivalence.
\end{lemma}
\begin{proof}
In the case where $X\rightarrt Z$ is a cofibration, is a straightforward
application of the definition of left properness. Suppose that $X\rightarrt Y$ and
$X\rightarrt V$ are cofibrations. 
Choose a factorization of the map $X\rightarrt Z$ into $X\rightarrt^{c} W \rightarrt^{f} Z$
such that $c$ is a cofibration and $f$ is a trivial fibration. The map
$$
W\cup ^X Y \rightarrt (W \cup ^X Y)\cup ^Y V = W\cup ^X V
$$
is the pushout of the weak equivalence $Y\rightarrt V$ along the cofibration $Y\rightarrt W\cup ^X Y$, so
by the left properness condition, it is a weak equivalence. On the other hand, the morphisms 
$W\cup ^X Y\rightarrt Z \cup ^X Y$ and $W\cup ^X V\rightarrt Z \cup ^X V$ are pushouts of the weak equivalence
$W\rightarrt Z$ along the cofibrations $W\rightarrt W\cup ^X Y$ and $W\rightarrt W\cup ^X V$ respectively.
Again by left properness, these maps are weak equivalences. Hence, in the diagram
$$
\begin{diagram}
W\cup ^X Y & \rightarr & W\cup ^X V \\
\downarr  & \rdTeXto & \downarr \\
Z \cup ^X Y &  & Z \cup ^X V
\end{diagram}
$$
the top arrow and both vertical arrows are weak equivalences; hence the bottom arrow is a weak equivalence as desired.
\end{proof}

\begin{corollary}
\label{pushoutequivcor}
Suppose $\mM$ is a left proper model category and suppose given a
diagram 
$$
\begin{diagram}
X & \leftarr & Y & \rightarr & Z \\
\downarr & & \downarr & & \downarr \\
A & \leftarr & B & \rightarr & C 
\end{diagram}
$$
such that the vertical arrows are weak equivalences, and the left horizontal
arrows are cofibrations. Then the map 
$$
X\cup ^YZ \rightarr A \cup ^BC
$$
is a weak equivalence.
\end{corollary}
\begin{proof}
Applying the previous lemma on both sides gives the required statement 
in the case $Y=B$. In general let $A':= X\cup ^YB$. The map $X\rightarrt A'$
is a weak equivalence by left properness, so $A'\rightarrt A$
is a weak equivalence by 3 for 2. The map $B\rightarrt A'$ is a cofibration. 
Applying the case where the middle map is an isomorphism gives the
statement that 
$$
A'\cup ^BC\rightarr A\cup ^BC
$$
is a weak equivalence. That case, or really the previous lemma, also
implies that the map 
$$
X\cup ^YZ\rightarrt X\cup ^YC = (X\cup ^YB)\cup ^BC = A'\cup ^BC
$$
is a weak equivalence. Putting these together gives the required statement. 
\end{proof}

We sometimes need an analogous invariance property for transfinite compositions. 

\begin{proposition}
\label{prop-transfiniteleftproper}
Suppose $\mM$ is a tractable left proper model category. Then 
transfinite compositions of cofibrations are invariant under homotopy; that is, if
$\{ X_n \} _{n< \beta }$ and $\{ Y_n \} _{n< \beta}$ are continuous
sequences indexed by an ordinal $\beta$,
with cofibrant transition maps,
and if we have a map of sequences given by a compatible collection of weak equivalences $X_n\rightarrt^{g_n} Y_n$,
the induced map $\colim _{n<\beta}X_n \rightarrt ^{g_{\beta}} \colim _{n<\beta}Y_n$ is a weak equivalence.
\end{proposition}
\begin{proof}
The proof is by induction on the ordinal $\beta$; we may assume it is known for
all sequences indexed by ordinals $\alpha <\beta$. 
Suppose given continuous 
sequences $\{ X_n \} _{n< \beta }$ and $\{ Y_n \} _{n< \beta}$ with cofibrant
transition maps, and $X_n\rightarrt^{g_n} Y_n$ compatible with the transition maps of the
sequences. Put $Z_n:= X_n\cup ^{X_0}Y_0$. By left properness, 
the maps $X_n\rightarrt Z_n$ are weak equivalences, hence by 3 for 2 
the same holds for the
maps $Z_n\rightarrt Y_n$ induced by the universal property of the pushout. 
Furthermore, the map
$X_0\rightarrt \colim _{n<\beta }X_n$ is a cofibration, so again 
by left properness the map
$$
\colim _{n<\beta }X_n\rightarrt
\colim _{n<\beta }Z_n = (\colim _{n<\beta }X_n)\cup ^{X_0}Y_0
$$
is a weak equivalence. 

Choose inductively $W_n$ fitting into a diagram 
$$
\begin{diagram}
Z_0 & \rightarr & Z_1 & \rightarr & \cdots & \rightarr & Z_n & \rightarr & \cdots \\
\dEqualarr & & \downarr_{i_1} &  & & & \downarr_{i_n} & & \\
W_0 & \rightarr & W_1 & \rightarr & \cdots & \rightarr & W_n & \rightarr & \cdots \\
\dEqualarr & & \downarr_{p_1} &  & & & \downarr_{p_n} & & \\
Y_0 & \rightarr & Y_1 & \rightarr & \cdots & \rightarr & Y_n & \rightarr & \cdots 
\end{diagram}
$$
such that $i_n$ are trivial cofibrations and $p_n$ are trivial fibrations, for any
non-limiting ordinal $n<\beta$. At a limit ordinal assume $W_{\cdot}$ is continuous,
that is $W_n= \lim _{j<n}W_j$. Note that the map $Z_n\rightarrt W_n$ will still
satisfy the left lifting property against any fibration, so $Z_n\rightarrt W_n$ is
still a trivial cofibration when $n$ is a limit ordinal, and similarly
once we are done the map $\colim _{n<\beta}Z_n\rightarrt \colim _{n<\beta}W_n$
will be a trivial cofibration. On the other hand, for a limiting ordinal $n<\beta$
the induction hypothesis implies that the map $p_n$ is still a weak equivalence. 
We may also choose $W_{\cdot}$ to have cofibrant transition maps, indeed 
given the choice up to $n$ we apply the factorization property to the
map $W_n\cup ^{Z_n}Z_{n+1}\rightarrt Y_n$, this map is a weak equivalence by
left properness, so it factors as a trivial cofibration followed by a trivial 
fibration. 

The trivial fibration property at successor ordinals, and continuity at limit
ordinals, plus the fact that the transition maps of $Y_{\cdot}$ are cofibrations,
allow us to choose a sequence of maps $Y_n\rightarrt ^{s_n}W_n$
which are sections of $p_n$ and commute with the transition  maps. 
The map $p_{\beta}: \colim _{n<\beta}W_{n}\rightarrt \colim _{n<\beta} Y_{n}$
therefore admits a section $s_{\beta}$
given by the colimit of the $s_n$. This shows that $p_{\beta}$ admits a right inverse.

We have now proven that our original map $g_{\beta}$ admits a right inverse in
the homotopy category. However, $s_{\beta}$ comes from a system satisfying the
same hypotheses, so $s_{\beta}$ also admits a right inverse in the homotopy category; 
but it has by construction
a left inverse, so $s_{\beta}$ is a weak equivalence; then its left inverse $p_{\beta}$
is a weak equivalence which implies that $g_{\beta}$ is one too. 
\end{proof}

\section{Quillen adjunctions}

Suppose $\mM$ and $\mN$ are model categories. A {\em Quillen adjunction} from $\mM$ to $\mN$ is an adjoint pair of functors $L:\mM \rightarrt \mN$
and $R: \mN \rightarrt \mM$, with $L$ left adjoint and $R$ right adjoint, such that:
\newline
(QA1)---$L$ sends cofibrations to cofibrations and trivial cofibrations to trivial cofibrations; 
\newline
(QA2)---$R$ sends fibrations to fibrations and trivial fibrations to trivial fibrations.

Either one of these conditions implies the other, by the adjunction formula. 

If $L:\mM \rightarrt \mN$ is a functor admitting a right adjoint $R$ such that $(L,R)$ is a Quillen adjunction, we say that $L$ is a {\em left 
Quillen functor}, similarly the right adjoint $R$ is a {\em right Quillen functor}. We say that $L$ or $R$ or $(L,R)$ is a {\em Quillen equivalence}
if a map $L(x)\rightarrt y$ is a weak equivalence if and only if the corresponding
map $x\rightarrt R(y)$ is a weak equivalence.

\section{The Kan-Quillen model category of simplicial sets}

The most important model category is the category of simplicial sets with the structure of Quillen model category originally defined by Kan.
We denote this model category by $\mK$. As a category, it is the category of functors $\Delta ^o\rightarrt \Sets _{\univa}$ where $\univa$ denotes
our main chosen universe. The cofibrations are just the monomorphisms. The fibrations are the maps satisfying the Kan lifting condition for
all horns over standard simplices. A modern proof that these classes of maps provide a model structure, is given in \cite{GoerssJardine}.

For $\mK$ there are particularly good generating sets. Let $h(n)$ denote the standard simplex of dimension
$n$, which as a diagram $\Delta ^o\rightarrt \Sets _{\univa}$ is just the functor represented by 
$[n]$. Let $\partial h(n)$ denote the boundary, defined by the condition that $(\partial h(n))_k$ is
the subset of simplices in $h(n)_k$ which factor through some principal face $h(n-1)\subset h(n)$.
The cofibrations are generated by the inclusions of boundaries $\partial h(n)\hookrightarrow h(n)$. 
The {\em $k$-th horn} $\partial \langle k\rangle h(n)$ is the subobject spanned by all principal faces 
except the $k$-th one. The trivial cofibrations are generated by the inclusions 
$\partial \langle k\rangle h(n)\hookrightarrow h(n)$.

\section{Model structures on diagram categories}

Suppose $\Phi$ is a small category, and $\mM$ is a combinatorial model category. 
Recall that $\diag (\Phi , \mM )$ is the category of functors $\Phi \rightarrt \mM$.
There are two main types of model structure on $\diag (\Phi ,\mM )$ usually known as the {\em projective} and 
{\em injective} model structures. In both model structures, the weak equivalences are defined to be the {\em objectwise} ones,
in other words $A\rightarrt B$ is a weak equivalence if and only if $A(x)\rightarrt B(x)$ is a weak equivalence for each $x\in \Phi$. 

In the projective structure, the fibrations are defined as the objectwise ones. The cofibrations, on the other hand, 
will be generated by
the $i_{x,!}(f)$ where $i_x:\{ x\} \rightarrt \Phi$ is inclusion of a single object,
and $f$ runs through a generating set of cofibrations for $\mM$. 

Dually, in the injective structure the cofibrations are defined as the objectwise ones, but the fibrations don't have
an easy description. 

\begin{theorem}
\label{injectiveprojective}
Suppose $\Phi$ is a small category, and $\mM$ is a combinatorial model category. Then $\diag (\Phi ,\mM )$ 
has two structures of combinatorial model category, the projective and the injective ones, with classes of morphisms described above. If $\mM$ is left proper then
so are the injective and projective diagram structures. 
\end{theorem}

See \cite{Barwick}, \cite[Proposition A.2.8.2]{LurieTopos} for the general case. 
Left properness may be verified levelwise. The main difficulty is in the construction of the injective
model structure, to get a generating set for the cofibrations. This is done using Lurie's theorem which was recalled as \ref{lurie-thm2}
in the preceding chapter. Of course, these model structures have a long history going back to Bousfield-Kan  and others, more recently considered by Hirschhorn \cite{Hirschhorn},
Smith (see Beke \cite{Beke}, Dugger \cite{DuggerCombinatorial}), Blander \cite{Blander}, Barwick \cite{Barwick}. 
For most special cases including many of the cases of interest to us, the reader may refer to any of a number of these references.

Barwick generalises this result to the relative case \cite{Barwick}, using a generalization of
Lurie's theorem (which we have stated above as Theorem \ref{lurie-thm2}). 
Recall that a {\em left Quillen presheaf} is
a presheaf of model categories over a base category, such that the transition functors are left Quillen functors.
The {\em category of sections} is the category of sections of the associated fibered category.

\begin{theorem}
\label{barwickinjective}
The category of sections of a left Quillen presheaf whose values are combinatorial model categories,
has an ``injective'' and a ``projective'' combinatorial model structure.
\end{theorem}

In the case we need which corresponds to diagrams in a constant combinatorial model category this was given by
Lurie in \cite{LurieTopos}. 
See \cite{Barwick} for the proof in general, which for the injective structure uses Lurie's techniques as discussed in Theorems \ref{lurie-thm1} and
\ref{lurie-thm2} above. The notion of left Quillen presheaf will not be used below,
the above statement was given for informational purposes. 

In Chapter \ref{algtheor1} 
below, we consider a variant of diagrams called {\em unital diagram categories}.
Given $\Phi _0\subset \Ob (\Phi )$ and a single model category $\mM$, we could define a presheaf of model categories
by setting $\mM _x:= \mM$ if $x\not \in \Phi _0$, with $\mM _x:=\{ \ast \}$ for $x\in \Phi _0$. The unital diagram category
$\diag (\Phi /\Phi _0; \mM )$ is the category of sections of the associated fibered category. This presheaf of categories is usually not a left Quillen presheaf
so Propositions \ref{unitalinjective} and \ref{unitalprojective}
 can't be viewed as a corollary of \ref{barwickinjective}, however the techniques of proof are the same. 

Another important type of model structure on certain diagram categories is the {\em Reedy} model structure, which is defined
when $\Phi$ has a structure of ``Reedy category''. This lies in between the projective and the injective structures, and indeed a closely related variant
will be provide an important class of cofibrations in Chapter \ref{cofib1} below. 
Some references for this discussion are \cite{Reedy} \cite{Bousfield-Kan}
\cite{Hirschhorn} \cite{DwyerHirschhornKan} \cite{DK3} \cite{GoerssJardine} \cite{BarwickReedy}
and \cite[Chapter 17]{descente}. 

A {\em Reedy category} is a category $\Phi$ provided with two subcategories on the same 
object set, called the {\em direct subcategory} $\Phi ^d$ and the {\em inverse subcategory}
$\Phi ^i$
and a function {\em degree} from the set of objects to an ordinal (usually $\omega$),
such that the non-identity direct maps strictly increase the degree, the
non-identity inverse maps strictly decrease the degree, and any morphism
$f$ factors uniquely as $f^df^i$ where $f^d$ is direct and $f^i$ is inverse. 
The degree of the middle object in this factorization is $\leq$ the degrees of
the source and target of $f$, and in case of equality $f$ is either direct or inverse
(or both in which case it is the identity). 

For $y\in \Phi$ define ${\Latch}(y)$ to be $\Phi ^d/y$ minus $\{y\}$
and ${\Match}(y)$ to be $y/\Phi ^i$ minus $\{ y\}$. 

Suppose $A:\Phi \rightarrt \mM$ is a diagram. The {\em latching and matching objects} at $y\in \Phi$
are defined to be
$$
{\latch}(A,y):= \colim A|_{{\Latch}(y)},
$$
$$
{\match}(A,y):= \mylim A|_{{\Match}(y)}.
$$
A diagram $A\in \diag (\Phi , \mM )$ is {\em Reedy cofibrant} if the morphisms 
$$
{\latch}(A,y)\rightarrt A(y)
$$
are cofibrations in $\mM$, and {\em Reedy fibrant} if 
the morphisms 
$$
A(y)\rightarrt \nocom {\rm match}(A,y)
$$
are fibrations in $\mM$. 
A morphism  $A\rightarrt^fB$ in $\diag (\Phi , \mM )$ is said to be {\em Reedy 
cofibrant} if the maps
$$
{\latch}(f,y):= {\latch}(B,y)\cup ^{{\latch}(A,y)}A(y)\rightarrt B(y)
$$
are cofibrations, and {\em Reedy fibrant} if the maps 
$$
A(y)\rightarrt \nocom {\match}(A,y)\times _{{\match}(B,y)}B(y) =: {\match}(f,y)
$$
are fibrations. 

\begin{proposition}
\label{reedystructure}
The category of diagrams $\diag (\Phi , \mM )$ provided with the levelwise
weak equivalences, and the above classes of cofibrations and fibrations, is a
closed model category, fitting in the middle of a 
sequence of left Quillen functors
$$
\diag _{\rm proj}(\Phi , \mM )\rightarrt \diag {\rm Reedy}(\Phi , \mM )\rightarrt 
\diag _{\rm inj}(\Phi , \mM ).
$$
The Reedy structure is combinatorial (resp. tractable, left proper)
whenever $\mM$ is. 
\end{proposition}
\begin{proof}
See the references \cite{Reedy} \cite{Bousfield-Kan}
\cite{Hirschhorn} \cite{DwyerHirschhornKan} \cite{DK3} \cite{GoerssJardine}
\cite{BarwickReedy}. 
In particular, inheritance of the combinatorial or tractable properties
is shown by Barwick in Lemmas 3.10, 3.11 of \cite{BarwickReedy}. 
The left properness condition may be verified levelwise, since Reedy cofibrations
are injective ones \cite[Lemma 3.1]{BarwickReedy}.  
\end{proof}

\section{Pseudo-generating sets}
\label{sec-pseudogen}

The cofibrant generating sets for a cofibrantly generated model category do not often have as simple and geometric a meaning as for the 
original case of simplicial sets $\mK$. This problem tends to occur particularly for the generating set for trivial cofibrations, which is often obtained
by an abstract accessibility argument leading to a set containing all trivial cofibrations up to a given cardinality. In this section,
we explore a way of defining a cofibrantly generated model category using sets $I$ and $K$. The first will be the generating set for
cofibrations; but the set $K$ will only generate the weak equivalences in a roundabout way, thus the terminology ``pseudo-generating sets''. 
The construction of a cofibrantly generated model structure in Theorem \ref{recog} uses exactly the argument where we throw in everything
up to a given cardinality. Once we have the statement of  Theorem \ref{recog} we will be able to apply it in later chapters without
having to come back to this cardinality argument. The reader hoping to avoid too much theory of model categories could therefore
skip this section and just take Theorem \ref{recog} as a ``black box'' for constructing model structures. The pseudo-generating sets
$I$ and $K$ used later will have some geometric meaning and hence be motivated outside of the technicalities of model category theory. 

Suppose we are given a locally presentable category $\mM$, and two sets of morphisms
$I,K\subset \Arr (\mM )$. Here both $I$ and $K$ are assumed to be small sets. 
Say that 
a morphism $f:A\rightarrt B$ is a {\em weak equivalence} if and only if (PG) there exists a diagram
$$
\begin{diagram}
A & \rightarr & A' \\
\downarr^{f} && \downarr _{f'}\\
B & \rightarr & B' 
\end{diagram}
$$
such that the horizontal morphisms are in $\cell (K)$, and the morphism $f'$ is in $\inj (I)$.
Note in particular that any morphism in $\inj (I)$ is a weak equivalence. 
 
Define the class of cofibrations to be $\cof (I)$, the trivial cofibrations to be the cofibrations which are weak equivalences, and
the fibrations to be the morphisms satisfying right lifting with respect to trivial cofibrations. As usual a cofibrant object
means an object $X$ such that the morphism $\emptyset \rightarrt X$ is a cofibration. 

Suppose the following axioms:
\newline
(PGM1)---record here the hypotheses that $\mM$ is locally presentable, and $I$ and $K$ are small sets of morphisms;
\newline
(PGM2)---the domains of arrows in $I$ and $K$ are cofibrant, and $K\subset \cof (I)$; 
\newline
(PGM3)---the class of weak equivalences is closed under retracts;
\newline
(PGM4)---the class of weak equivalences satisfies 3 for 2;
\newline
(PGM5)---the class of trivial cofibrations is closed under pushouts;
\newline
(PGM6)---the class of trivial cofibrations is closed under transfinite composition.

Using these axioms, we would like to show that these classes define a cofibrantly generated, and indeed tractable, model category structure on $\mM$.
Note that the class of trivial fibrations
is defined as the intersection of the fibrations and the weak equivalences; we don't know {\em a priori} that this is the same
as $\inj (I)$, that will have to be proven as a consequence of the axioms. One can say, however, that $\inj (I)$ is
contained in the class of trivial fibrations, indeed an element of $\inj (I)$ satisfies right lifting with respect to
$\cof (I)$ so it is a fibration, and from the definition (PG) it is a weak equivalence. 

One important preliminary result is Corollary \ref{AfiltcolimCor} from the previous chapter, saying that any morphism in $\cell (I)$ is
a $\kappa$-filtered colimit of $\kappa$-presentable cell complexes. The other main ingredient is the following
observation which is a sort of accessibility property for morphisms in $\inj (I)$. 

\begin{lemma}
\label{Afilt}
We suppose given three regular cardinals $\mu < \lambda < \kappa$ such that $\mM$ is locally $\lambda$-presentable
(hence also locally $\kappa$-presentable), such that $|I| < \kappa$, and such that the sources and targets of arrows in $I$ are
$\mu$-presentable. We assume that $2 ^{\mu}<\kappa$. 

Suppose $f:X\rightarrt Y$ is in $\inj (I)$, and suppose it can be expressed $f=\colim _{i\in \alpha}f_i$
where $f_i: X_i\rightarrt Y_i$ are arrows between $\kappa$-presentable objects, and $\alpha$ is $\kappa$-filtered. 
Then there is a collection of $\lambda$-filtered categories $\beta _j$ of size $|\beta _j| < \kappa $,
together with functors $q_j:\beta _j\subset \alpha$,
all indexed by a $\kappa$-filtered poset $j\in \psi $,
such that for any $j$, $f(\beta _j):= \colim _{i\in \beta _j}f_i$ is in $\inj (I)$,
and $f = \colim _{j\in \psi}f(\beta _j)$. 
\end{lemma}
\begin{proof}
The first step is to say the following. 
For each $i\in \alpha$ there is $t(i)\in \alpha$ with an arrow $i\rightarrt t(i)$, such that
for any diagram
$$
\begin{diagram}
U & \rightarr & X_i \\
\downarr^u && \downarr \\
V & \rightarr & Y_i
\end{diagram}
$$
with $u\in I$, there exists a lifting $V\rightarrt X_{t(i)}$
so that the two triangles in 
$$
\begin{diagram}
U & \rightarr & X_{t(i)} \\
\downarr & \ruTeXto & \downarr \\
V & \rightarr & Y_{t(i)}
\end{diagram}
$$
commute. Here the horizontal maps are the compositions of the previous ones, with the transition maps
$X_i\rightarrt X_{t(i)}$ and $Y_i\rightarrt Y_{t(i)}$.

To prove this, note that there exist liftings $V\rightarrt X$ for any diagram. However,
the $X_i$ and $Y_i$ are $\kappa$-presentable, the arrows in $I$ are $\mu$-presentable, and $|I|< \kappa$.
Thus, the cardinality of the set of diagrams we need to consider is $<\kappa^{\mu} \leq \kappa$, by Corollary \ref{lambdamuexp}. Since again $U$ is $\kappa$-presentable,
there is some $t\in \alpha$ which works for each diagram; but since $\alpha$ is $\kappa$-filtered we can choose a single
$t(i)$ for all the diagrams at a given value of $i$. This completes the proof of the first step. 

Now we can exhaust $\alpha$ by a family of $\lambda$-filtered subcategories $\beta _j$ with $|\beta _j|< \kappa$, 
indexed by $j\in \psi$ where $\psi$
is a $\kappa$-filtered partially ordered set. We can do it in such a way that for any $i\in \beta _j$ we also have $t(i)\in \beta _j$. 
The exhaustion condition means that for any $i\in \alpha$ there is
$j\in \psi$ and $k\in \beta _j$ with $i\rightarrt k$, and it implies that $f = \colim _{j\in \psi}f(\beta _j)$.

Now put 
$$
\begin{diagram}
X(\beta _j)&:=& \colim _{i\in \beta _j}X_i\\
& & \downarr_{f(\beta _j):= \colim _{i\in \beta _j}f_i} \\
Y(\beta _j)&:=& \colim _{i\in \beta _j}Y_i .
\end{diagram}
$$
We claim that $f(\beta _j)\in \inj (I)$. If 
$$
\begin{diagram}
U & \rightarr & X(\beta _j) \\
\downarr^u && \downarr \\
V & \rightarr & Y(\beta _j)
\end{diagram}
$$
is a diagram with $u\in I$, then since $\beta _j$ is $\lambda$-filtered and $u$ is $\lambda$-small,
there exists a lifting to a diagram of the form
$$
\begin{diagram}
U & \rightarr & X_i \\
\downarr^u && \downarr \\
V & \rightarr & Y_i 
\end{diagram}
$$
for some $i\in \beta _j$.
By our choice of $t(i)$ plus the hypothesis that $t(i)\in \beta _j$ whenever $i\in \beta _j$, we get a
lifting $V\rightarrt X_{t(i)}\rightarrt X(\beta _j)$ which makes the triangles in 
$$
\begin{diagram}
U & \rightarr & X(\beta _j) \\
\downarr^u & \ruTeXto & \downarr \\
V & \rightarr & Y(\beta _j)
\end{diagram}
$$
commute. This shows that $f(\beta _j)\in \inj (I)$. 
\end{proof}

A {\em cofibrant replacement} of an  object $X\in \mM$ is a morphism $p:X'\rightarrt X$
such that $p\in \inj (I)$ and $\emptyset \rightarrt X'$ is in $\cell (I)$. This exists by the small object argument
for $I$, and by the definition (PG) of weak equivalence, it follows that
$p\in \inj (I)$
is a weak equivalence.

Construct as follows a set $J$ of trivial cofibrations with cofibrant domains. 
Recall that the category $\mM$ is locally $\kappa$-presentable. 
Choose a small set $N_1$ of representatives for the isomorphism classes of 
arrows $f:X\rightarrt Y$ which are weak equivalences between $\kappa$-presentable objects.
Since the isomorphism classes of $\kappa$-presentable objects form a small set (axiom (2) for the locally $\kappa$-presentable category $\mM$), 
we can choose $N_1$ as a
small set.  For each $f\in N_1$, choose a cofibrant replacement $p:X'\rightarrt X$, let $f':X'\rightarrt Y$ be the
composition $f'=fp$, and choose a factorization $f'=gh$ where $h:X'\rightarrt Z$ is in $\cell (I)$ and
$g:Z\rightarrt Y$ is in $\inj (I)$, in particular $g$ is a trivial fibration. Let $J$ be the set of all cofibrations $h$
obtained in this way. Using the facts that $g$ and $p$ are weak equivalences, and the hypothesis that $f$ is a weak equivalence,
we get that $f'$ and then $h$ are weak equivalences by (PGM4). Thus the elements of $J$ are trivial cofibrations with
cofibrant domains. We have to show that a map in $\inj (J)$ is a new fibration, that is that it satisfies lifting with respect to 
any new trivial cofibration. Or equivalently, to show that any new trivial cofibration is in $\cof (J)$.

\begin{lemma}
\label{cellJ}
The elements of $\cof (J)$ are trivial cofibrations.
\end{lemma}
\begin{proof}
Trivial cofibrations are closed under pushouts (PGM5) and transfinite compositions (PGM6), hence the elements
of $\cell (J)$ are trivial cofibrations. The class $\cof (J)$ is the closure of $\cell (J)$ under retracts;
by definition the class $\cof (I)$ of cofibrations is closed under retracts, and by (PGM3) the class of weak equivalences
is closed under retracts, so elements of $\cof (J)$ are trivial cofibrations. 
\end{proof}

The following proposition says that the class of weak equivalences is accessible. In this argument, we are following Barwick \cite{Barwick}. 

\begin{proposition}
\label{nweaccessible}
Suppose $f:X\rightarrt Y$ is a weak equivalence. Then $f$ can be expressed as a $\kappa$-filtered colimit of
arrows $f_i$ which are weak equivalences between $\kappa$-presentable objects.
\end{proposition}
\begin{proof}
The category $\Arr (\mM )$ is locally $\kappa$-presentable, so we can write the arrow $f:X\rightarrt Y$ as
a $\kappa$-filtered
colimit of arrows $\colim _{i\in \alpha} f_i$, such that $f_i:X_i\rightarrt Y_i$ is $\kappa$-presentable in $\Arr (\mM )$. In particular,
$X = \colim _{i\in \alpha} X_i$ and $Y = \colim _{i\in \alpha}Y_i$ are expressed 
as $\kappa$-filtered colimits of $\kappa$-presentable  objects in $\mM$ (see Lemma \ref{diagpres}).

By the definition of weak equivalence depending on $I$ and $K$, there exists a diagram
$$
\begin{diagram}
X & \rightarr^{a} & A \\
\downarr ^f& & \downarr_{ g} \\
Y & \rightarr^{b} & B
\end{diagram}
$$
such that $a$ and $b$ are in $\cell (K)$, and such that $g\in \inj (I)$. 

Apply Corollary \ref{AfiltcolimCor}
to $a$ and $b$. Thus
$$
A= \colim _{i\in\alpha} A_i,\;\;\; B= \colim _{i\in\alpha} B_i
$$
such that $a$ and $b$ are the colimits of systems $a_i:X_i\rightarrt A_i$ and $b_i:Y_i\rightarrt B_i$,
with $A_i$ and $B_i$ being $\kappa$-presentable, and $a_i,b_i\in \cell (K)$. 
For each $i$ we have a diagram
$$
\begin{diagram}
X_i & \rightarr^{a_i} & A_i \\
\downarr^{ f_i} & & \downarr _{g_i} \\
Y_i & \rightarr^{b_i} & B = \colim _{j\in \alpha}B_j .
\end{diagram}
$$
This gives a collection of maps, natural in $i\in \alpha$,
$$
Y_i\cup ^{X_i}A_i \rightarrt B
$$
with the sources being $\kappa$-presentable.  
By Lemma \ref{colimap} of the previous chapter, there is a functor $q:\alpha \rightarrt \alpha$
together with factorizations
$$
Y_i\cup ^{X_i}A_i\rightarrt B_{q(i)}\rightarrt B
$$
natural in $i$.

This gives a functor  of
diagrams depending on $i\in \alpha$
$$
\begin{diagram}
X_i & \rightarr^{a_i} & A_i \\
\downarr^{f_i} & & \downarr _{g_i} \\
Y_{q(i)} & \rightarr^{b_{q(i)}} & B_{q(i)}
\end{diagram}
$$
such that $g=\colim g_i$, $a=\colim a_i$, $b=\colim b_{q(i)}$ and $a_i,b_{q(i)}\in \cell (K)$.

Apply Lemma \ref{Afilt} to the map $g\in \inj (I)$; we conclude that the original diagram may be seen as a 
$\kappa$-filtered colimit over $j\in \psi$ of diagrams of the form 
$$
\begin{diagram}
X(\beta _j) & \rightarr^{a(\beta _j)} & A(\beta _j) \\
\downarr^{ f(\beta _j)} & & \downarr _{g(\beta _j)} \\
Y (\beta _j)& \rightarr^{b(\beta _j)} & B(\beta _j)
\end{diagram}
$$
where $X(\beta _j)= \colim _{i\in \beta _j}X_i$ etc.,
and such that the maps $g(\beta _j)$ are in $\inj (I)$. The maps $a(\beta _j)$ and $b(\beta _j)$ are still in $\cell (K)$,
so this shows that $f(\beta _j):X(\beta _j)\rightarrt Y(\beta _j)$ are weak equivalences. The objects 
$X(\beta _j)$ and $Y(\beta _j)$ are $\kappa$-presentable, and $f=\colim _{j\in \psi}f(\beta _j)$, which completes the proof of
the proposition. 
\end{proof}

\begin{theorem}
\label{access}
Our set of arrows $J$ generates the trivial cofibrations, in other words
the class $\cof (J)$ equals the class of all trivial cofibrations. 
In particular, $\inj (J)$ equals the class of fibrations. 
\end{theorem}
\begin{proof}
Suppose we are given a  trivial cofibration $f:X\rightarrt Y$. 
For now suppose also that $f$ is in $\cell (I)$. 

Choose a cellular expression for $f$, that is
a sequence $\{X_n\} _{n\leq \beta}$ indexed by an ordinal $\beta$ with $X_0=X$, $X_{\beta}=Y$, $X_m= \colim _{n<m}X_m$ when
$m$ is a limit ordinal, and $X_n\rightarrt X_{n+1}$ obtained by pushout along an element of $I$. 

We will choose a sequence of diagrams $X_n\rightarrt Z_n\rightarrt Y$, where $\{ Z_n\}_{n\leq \beta}$ is a transfinite sequence, 
such that $X\rightarrt Z_n$ is in $\cell (J)$ in particular a new trivial cofibration;
$Z_n\rightarrt Z_{n+1}$ is pushout along an element of $J$, and again $Z_m=\colim _{n<m}Z_n$ when $m$ is a limit ordinal.
The transition maps in the sequence $\{ Z_n\}_{n\leq \beta}$ are supposed to be compatible
with those of the sequence  $\{X_n\} _{n\leq \beta}$ and with the maps to $Y$.

For the induction step at a limit ordinal $m$, just let $Z_m$ be the colimit of the $Z_{n}$ for $n<m$. 

Suppose $n$ is a successor  ordinal with the $Z_i$ for $i\leq n-1$ given, and we want to choose $Z_{n}$. 
Consider the map $w_n:U_n\rightarrt V_n$ in $I$ such that $X_{n-1}\rightarrt X_n$ is pushout along $w_n$;  
we have a commutative diagram
$$
\begin{diagram}
U_n & \rightarr & X_{n-1} & \rightarr & Z_{n-1} \\
\downarr & & \downarr & & \downarr \\
V_n & \rightarr & X_n & \rightarr & Y
\end{diagram}
$$
where the first square is a pushout. The map $X\rightarrt Z_{n-1}$ is a weak equivalence by the inductive hypothesis,
and $X\rightarrt Y$ is a weak equivalence by assumption, so by (PGM4), the map $Z_{n-1}\rightarrt Y$ is
a weak equivalence. On the other hand, $w_n: U_n\rightarrt V_n$ is in $I$, so by the choice of $\kappa$ it is $\kappa$-presentable in $\Arr (\mM )$. 
By Proposition \ref{nweaccessible} the map $Z_{n-1}\rightarrt Y$ is a $\kappa$-filtered colimit of $\kappa$-presentable arrows in $\Arr (\mM )$ which
are new weak equivalences---that is, elements of $N_1$. Since $w_n$ is $\kappa$-presentable it factors through one of them. Therefore, there exists a factorization
$$
\begin{diagram}
U_n & \rightarr^{a} & A & \rightarr^{z} & Z_{n-1} \\
\downarr & & \downarr & & \downarr \\
V_n & \rightarr^{b} & B & \rightarr^{y} & Y
\end{diagram}
$$
with the middle vertical arrow in $N_1$. Consider the choice of 
cofibrant replacement $A'\rightarrt A$ used for the definition of $J$ and, since $U_n$ is cofibrant (PGM2),
a lifting $U_n\rightarrt A'$. In particular we may replace $A$ by $A'$ in the previous diagram. Consider the choice of
factorization $A'\rightarrt^{u} C \rightarrt^{v} B$ with 
$u\in J$ and $v\in \inj (I)$, and choose a lifting $c:V_n\rightarrt C$ with $vc = b$ and $c$ restricting to the given map on $U_n$. 
This gives a factorization
$$
\begin{diagram}
U_n & \rightarr^{a} & A' & \rightarr & Z_{n-1} \\
\downarr & & \downarr & & \downarr \\
V_n & \rightarr^{c} & C & \rightarr^{yv} & Y
\end{diagram}
$$
where the middle vertical arrow is in $J$. Let $Z_n:= Z_{n-1}\cup ^{A'}C$ be the pushout along this arrow. 
Then $X\rightarrt Z_n$ is again in $\cell (J)$, in particular it is a new trivial cofibration (Lemma \ref{cellJ}).
From the second square in the above diagram we get a map $Z_n\rightarrt Y$ compatible with the map on $Z_{n-1}$.
But, given that $X_n= X_{n-1}\cup ^{U_n}V_n$ we get a map $X_n\rightarrt Z_n$ compatible with given map on $X_{n-1}$.
This completes the inductive step, giving the construction of the sequence $\{ Z_n\}_{n\leq \beta}$ as desired. 

Letting $Z:=Z_{\beta}$ we get a map $g:X\rightarrt Z$ in $\cell (J)$ with projection 
$p:Z\rightarrt Y$ and a splitting $s:Y\rightarrt Z$, $ps=1$ so that 
$f$ is a retract of $g$ in objects under $X$. This  shows $f\in \cof (J)$. 

We have shown that a map $f$ in $\cell (I)$ which is a new weak equivalence, is in $\cof (J)$.
Suppose $f\in \cof (I)$ is a weak equivalence. Choose a factorization $f=gh$ with 
$X\rightarrt^{h} V\rightarrt^{g}Y$ where $h\in \cell (I)$ and $g\in\inj (I)$.
Then $g$ is a weak equivalence from the definition (PG), so 
$h$ is a weak equivalence (PGM4). We have shown above that $h\in \cof (J)$. 
On the other hand, since $f\in \cof (I)$ it satisfies lifting with respect to $h$ so 
there is a section $s:Y\rightarrt Z$ making $f$ into a retract of $h$ in objects under $X$. 
In general the cofibrations for a set
of arrows forms a class closed under retracts, and $f$ is a retract of $h\in \cof (J)$ so $f\in \cof (J)$.  

This shows that $\cof (J)$ equals the class of trivial cofibrations. Given a map $g\in \inj (J)$ it satisfies lifting with respect to $\cell (J)$, and any trivial cofibration is a retract of something in $\cell (J)$. It follows that $g$ satisfies lifting with respect to any
trivial cofibration, so $g$ is a fibration. Conversely a fibration is in $\inj (J)$
as follows from Lemma \ref{cellJ}.
\end{proof}

\begin{corollary}
\label{newfactor}
Any map $f:X\rightarrt Y$ in $\mM$ factors as $f=gh$ with 
$X\rightarrt^{h} Z \rightarrt^{g} Y$
where $h$ is a trivial cofibration (which can be assumed even in $\cell (J)$) and $g$ is a fibration. 
\end{corollary}
\begin{proof}
This is just the factorization into $h\in \cell (J)$ and $g\in \inj (J)$ for the set $J$,
given by the small object argument. Note that $h$ is a
new trivial cofibration by Lemma  \ref{cellJ}, and $g$ is a new fibration by Theorem \ref{access}.
\end{proof}

\begin{lemma}
\label{trivfib}
The class of trivial fibrations is equal to $\inj (I)$. 
\end{lemma}
\begin{proof}
An element of $\inj (I)$ satisfies lifting with respect to all cofibrations, in particular with respect to
trivial ones, so $\inj (I)$ is contained in the class of fibrations. It is contained in the class of
weak equivalences by the definition (PG), which shows that it is contained in the class of trivial fibrations. 

Suppose $f:X\rightarrt Y$ is a trivial fibration. Applying Theorem \ref{access},
this means that $f\in \inj (J)$ and $f$ is a weak equivalence. Use the small object argument to
choose a factorization
of $f$ as the composition of
$$
X\rightarrt^{g} Z\rightarrt^{p} Y,
$$
such that $g\in \cell (I)$ and $p$ is in $\inj (I)$. 
It follows from the definition (PG) that $p$ is a weak equivalence.
By 3 for 2, the map $g$ is 
also a weak equivalence, so it is a trivial cofibration. Thus, 
the condition that $f$ be a fibration implies that it satisfies right
lifting with respect to $g$. Use this lifting property on the 
square with the identity $1_X$ along the top, $f$ on the right,
$g$ on the left and $p$ on the bottom: 
hence there is a map $t:Z\rightarrt X$ such that $tg=1_X$ and $ft=p$. 
This presents $f$ as a retract of $p$, but $p\in \inj (I)$ and $\inj (I)$ is
closed under retracts (see Lemma \ref{rlpretract} of the previous chapter),
so $f\in \inj (I)$.
\end{proof}

We can now put together everything above to obtain a model category structure on $\mM$.
This follows Smith's recognition theorem, as exposed in \cite{Beke}, \cite{Barwick}.
For the reader's convenience 
we repeat the definition of weak equivalence (PG) and the axioms equivalent to (PGM1)--(PGM6).
The present statement sums up the main result from the first two chapters which will be used later on.

\begin{theorem}
\label{recog}
Suppose $\mM$ is a locally presentable category, and  $I\subset \Arr (\mM )$ and $K\subset \cof (I)$ are sets of morphisms.
Say that a morphism $f:A\rightarrt B$ is a {\em weak equivalence} if and only if there exists a diagram
$$
\begin{diagram}
A & \rightarr & A' \\
\downarr ^f&& \downarr _{f'}\\
B & \rightarr & B' 
\end{diagram}
$$
such that the horizontal morphisms are in $\cell (K)$, and $f'$ is in $\inj (I)$.
Define the class of cofibrations to be $\cof (I)$. Define the trivial cofibrations to be the cofibrations which are weak equivalences, and
the fibrations to be the morphisms satisfying right lifting with respect to trivial cofibrations. Suppose:
\newline
---the domains of arrows in $I$ and $K$ are cofibrant; 
\newline
---the class of weak equivalences contains $\inj (I)$, is closed under retracts, and satisfies 3 for 2;
\newline
---the class of trivial cofibrations is closed under pushouts and transfinite composition. 
\newline
Then the class of trivial fibrations is exactly $\inj (I)$, and
$\mM$ with the given classes is a cofibrantly generated and indeed tractable model category.
\end{theorem}
\begin{proof}
Note that the class of weak equivalences is the one defined by condition (PG) above,
and the hypotheses stated in the theorem are equivalent to the system of axioms (PGM1)--(PGM6). 
Let $J$ be the set of morphisms given for Theorem \ref{access}, which says that $\inj (J)$ is the
class of fibrations and $\cof (J)$ is the class of trivial cofibrations. 
By Lemma \ref{trivfib}, $\inj (I)$ is the class of trivial fibrations. 

Verify first the axioms 
for a closed model category. 

Axiom (CM1) comes from the fact that $\mM$ is locally presentable.

Axiom (CM2) is a hypothesis.

Axiom (CM3) for weak equivalences is a hypothesis. In a  locally presentable category, for a given subset of 
morphisms $I$, the class $\cof (I)$ is closed under retracts. This gives (CM3) for the cofibrations.
The class of fibrations is defined to be the class of morphisms which satisfy the right lifting property with
respect to the trivial cofibrations. By Lemma \ref{rlpretract} of the previous chapter, the class of fibrations is closed under
retracts to give (CM3).

The fibrations are defined to be the maps satisfying the right lifting property with respect to trivial cofibrations.
Therefore, the trivial cofibrations satisfy the left lifting property with respect to any fibrations. This is one half of
(CM4). For the other half, Lemma \ref{trivfib} says that the class of trivial fibrations is equal to $\inj (I)$,
but the class of cofibrations is equal to $\cof (I)$ so the cofibrations satisfy the left lifting property with
respect to trivial fibrations. 

For (CM5), suppose $f:X\rightarrt Y$ is any morphism.
It can be factored by the small object argument in two ways, as
$$
X\rightarrt^{g} Z \rightarrt^{p} Y
$$
with $g\in \cell (J)\subset \cof (J)$ and $p\in \inj (J)$, or as
$$
X\rightarrt^{g'} Z' \rightarrt^{p'} Y
$$
with $g'\in \cell (I)\subset \cof (I)$ and $p'\in \inj (I)$. 
In the first, $g$ is a trivial cofibration and $p$ is a fibration by Theorem \ref{access};  
in the second $g'$ is a cofibration in view of the definition of cofibrations,
and $p'$ is a  trivial fibration by Lemma \ref{trivfib}. 
This proves the two parts of (CM5) and completes the proof of the Quillen model structure.

Next we note that the model structure is cofibrantly generated, indeed we have exhibited the
required generating sets $I$ and $J$. Axiom (CG1) is automatic since $\mM$ is locally presentable
(all objects are small); axiom (CG2a) is the definition of cofibrations and (CG2b) is Lemma \ref{trivfib};
and axioms (CG3a) and (CG3b) are given by Theorem \ref{access}.

As $\mM$ is locally presentable, the model structure is combinatorial. Furthermore it is
tractable: indeed we have required in (PGM2) that the domains of arrows in $I$ are cofibrant,
and the arrows in $J$ have cofibrant domains by construction. 
\end{proof}

It should perhaps be stressed that the  main work was done in the previous chapter in the discussion of
cell complexes and the proof of Theorems \ref{cellover} and \ref{lurie-thm1}, then used in the proof of Theorem \ref{access}
above. The advantage of the statement of Theorem \ref{recog} is that
it makes no reference to accessibility (other than the hypothesis that $\mM$ be locally presentable),
so when we apply it later on we don't need to do any more work with cardinals and presentability.

We can get some further information on trivial cofibrations, fibrant objects,
and fibrations between fibrant objects. This result will, as it is applied
successively in later chapters, eventually turn into
the version for constant object set
of Bergner's result characterizing fibrant Segal categories \cite{BergnerSegal}. 

\begin{proposition}
\label{furtherinfo}
In the situation of Theorem \ref{recog},
suppose $X\rightarrt ^fY$ is a cofibration (i.e. in $\cof (I)$). 
Then $f$ is a trivial cofibration if and only if there exists a diagram 
$$ 
\begin{diagram}
X & \rightarr ^a & A \\
\downarr ^f & & \uparr^s \downarr_g \\
Y & \rightarr^b & B
\end{diagram}
$$
such that $ga=bf$, $sbf = a$, $gs = 1_B$,
and $a,b\in \cell (K)$. 

An object $U\in \mM $ is fibrant if and only if $U\in \inj (K)$. 
If this is the case, then a morphism $W\rightarr^p U$ is a fibration if and only if
$p\in \inj (K)$.
\end{proposition}
\begin{proof}
If $f$ is a new trivial cofibration, then it fits into a diagram with 
a two maps $a,b\in \cell (K)$ and 
a morphism $g\in \inj (I)$ by the definition of weak equivalences
(PG). Since $bf$ is
cofibrant, there exists a lifting $s$ with $gs = 1_B$ and $sbf = a$.
This is a diagram of the required form. 

Suppose $f$ is a cofibration such that there exists a diagram as above
(but in this case $g$ is no longer assumed in $\inj (I)$). Then $bf$ is a retract of 
$a$ which is a weak equivalence. By closure of
weak equivalences under retracts and then 3 for 2, it follows that $f$ is a new weak equivalence hence a  trivial cofibration. This proves the first statement.

If $U$ is fibrant then it is clearly $K$-injective since the elements of $K$ are
trivial cofibrations. Suppose $U\in \inj (K)$.
If $X\rightarrt^f Y$ is a trivial cofibration, there exists a diagram as in the first
part.  For any map $X\rightarrt U$ we can extend it to a map $A\rightarrt U$
and composing with $sb$ gives the required extension to $Y\rightarrt A$. This
shows that $U$ is fibrant. 

Suppose now that $U$ is fibrant and $W\rightarrt ^p U$ is a map. If $p$ is a fibration
then clearly it is in $\inj (K)$. On the other hand suppose $p\in \inj (K)$, and  
$$
\begin{diagram}
X & \rightarr^w & W \\
\downarr ^f & & \downarr _p\\
Y & \rightarr^u & U
\end{diagram}
$$
is a diagram with $f$ being a trivial cofibration. Choose a diagram as in the first
part of the proposition for $f$. Since $U$ is assumed fibrant, we can extend the bottom
map to a map $B\rightarrt^t U$ with $tb=u$ and
which composes with $g$ to give a diagram 
$$
\begin{diagram}
X & \rightarr^w & W \\
\downarr ^a & & \downarr _p\\
A & \rightarr^{tg} & U.
\end{diagram}
$$
Now $a\in \cell (K)$ and $p\in \inj (K)$ so there is a lifting $A\rightarrt^r W$
with $ra=w$ and $pr=tg$. This gives a map $Y\rightarrt ^{rsb}W$ such that $rsbf = ra = w$
and $prsb = tgsb = tb = u$, which is the required right lifting property showing
that $p$ is a fibration. 
\end{proof}


\chapter{Cartesian model categories}
\label{cartmod1}
\label{sec-cartmodcat}

The notion of {\em cartesian model category} plays two important roles. First of all,
the Segal conditions involve direct products, so it is important that they behave
well homotopically. The second application is the fact that a cartesian model category $\mP$
leads to a $\mP$-enriched category obtained by looking at its fibrant and cofibrant objects
and using the internal $\uHom$ to define morphism objects. 
So, the cartesian condition is one of the main hypotheses but also one of the main properties which we would like to prove for our construction of a model
category $\precat (\mM )$ of $\mM$-enriched precategories. This compatibility with 
products will be shown in 
Chapter \ref{product1}, furthermore it will provide a useful trick to help establishing
the model structure.

By ``cartesian model category'', we mean a symmetric monoidal model category in the sense of Hovey \cite{Hovey} for the monoidal operation given by direct product. We add a condition about commutation of direct products
with colimits, in order to be able to get an internal $\uHom$. 

Recall that $\emptyset$ denotes the initial object and $\ast$ the coinitial object of $\mM$.

\begin{definition}
\label{def-cartesian}
Say that a combinatorial model category $\mM$ is {\em cartesian} if:
\newline
(DCL)---the direct product preserves colimits: if $\{ A_i\} _{i\in \alpha}$ and $\{ B_j\} _{j\in \beta}$
are diagrams, then 
$$
\colim _{\alpha \times \beta} A_i\times B_j = (\colim _{\alpha}A_i)\times (\colim _{\beta}B_j);
$$
\newline
(AST)---the map $\emptyset \rightarrt \ast$ is cofibrant;
\newline
(PROD)---for any cofibrations $A\rightarrt B$ and $C\rightarrt D$
the map
\begin{equation}
\label{pushprod}
A\times D \cup ^{A\times C}B\times C \rightarrt B\times D
\end{equation}
is a cofibration; if in addition at least one of $A\rightarrt B$ or $C\rightarrt D$ is a trivial cofibration then 
\eqref{pushprod} is a trivial cofibration. 
\end{definition}

See the definition of monoidal model category given for example in \cite{HoveyArxiv99} \cite{ShipleySchwede02}. Our first axiom says that the unit object 
for the direct product is cofibrant, which is stronger than the unit axiom of \cite{ShipleySchwede02}. 
Note that since $\mM$ is combinatorial, the factorizations can be chosen functorially so as to approach Hovey's definition. 

We record now some first consequences of this definition. 

\begin{lemma}
\label{emptyempty}
Condition (DCL) implies that for any object $X$ the natural map  $\emptyset\rightarrt X\times \emptyset$ is an isomorphism. 
In turn this implies that the object $\emptyset$ is empty: if $X\rightarrt \emptyset$ is any morphism, then $X\cong \emptyset$.
\end{lemma}
\begin{proof}
Note that $\emptyset$ is the colimit of the empty diagram, that is $\emptyset = \colim _{\emptyset{\rm cat}} F$ where
$\emptyset{\rm cat}$  is the empty category and 
$F:\emptyset{\rm cat} \rightarrt \mM$ is the unique functor. Compatibility of products and colimits (DCL)
therefore implies that for any $X = \colim _{1}X$ (here $1$ is the one-object category) we have
$$
X\times \emptyset = \colim _{1}X\times \colim _{\emptyset{\rm cat}} F = \colim _{1\times\emptyset{\rm cat}} (X\boxtimes F) = \emptyset 
$$
the last equality coming from the fact that $1\times \emptyset{\rm cat} = \emptyset{\rm cat}$ is again the empty category. 

Suppose now that $f:X\rightarrt \emptyset$ is a morphism. Letting $e:\emptyset \rightarrt X$ denote the unique morphism
we get $fe = 1_{\emptyset}$. On the other hand, the morphism $(1_X,f): X\rightarrt X\times \emptyset$ factors through
$(e,1_{\emptyset}):\emptyset \rightarrt^{\cong} X\times \emptyset$. Projecting to the second factor shows that the
factorization map is $f:X\rightarrt \emptyset$ and projecting to the first factor then shows that $1_X=ef$. Thus $f$ is an isomorphism inverse to $e$.  
\end{proof}

\begin{lemma}
\label{prodconsequence}
If $\mM$ satisfies (PROD) and (DCL) then, for any pair of weak equivalences 
$A\rightarrt^f B$ and $C\rightarrt^g D$, the resulting product map $A\times C\rightarrt^{(f,g)} B\times D$ is
a weak equivalence. 
\end{lemma}
\begin{proof}
Suppose first that $f$ and $g$ are trivial cofibrations between cofibrant objects.
Apply (PROD) to $f$ and $\emptyset \rightarrt C$, 
then to $\emptyset \rightarrt B$ and $g$, and use (DCL) via
the result of Lemma \ref{emptyempty}. 
These give that
both maps  
$$
A\times C\rightarrt B\times C \rightarrt B\times D
$$ 
are trivial cofibrations; hence their composition is a weak equivalence. 

Suppose $f$ is a cofibration between cofibrant objects, and $C$ is any object.
Choose a replacement $C'\rightarrt^p C$ where $C'$ is cofibrant and $p$ is a trivial fibration.
By the first paragraph, $A\times C'\rightarrt B\times C'$ is a weak equivalence.
On the other hand, $A\times C'\rightarrt A\times C$ and $B\times C'\rightarrt B\times C$
are trivial fibrations as can be seen directly from the lifting property. Writing a commutative
square and using 3 for 2 we get that the map $A\times C\rightarrt B\times C$ is a
weak equivalence. 

Suppose now that $f$ is an arbitrary  weak equivalence and $C$ is any object. 
Choose a cofibrant replacement
$A'\rightarrt A$, then complete to a square
$$
\begin{diagram}
A' & \rightarr & B' \\
\downarr && \downarr \\
A &\rightarr & B
\end{diagram}
$$
such that $B'\rightarrt B$ is a cofibrant replacement, and $A'\rightarrt B'$ is
a trivial cofibration. The vertical maps are trivial fibrations so their products with
$C$ remain trivial fibrations, and by the previous paragraph the map 
$A'\times C\rightarr  B'\times C$ is a weak equivalence. By 3 for 2, the map 
$A\times C\rightarrt B\times C$ is a weak equivalence. 

Similarly $B\times C\rightarrt B\times D$ is a weak equivalence, and composing
we get the statement of the lemma. 
\end{proof}

\begin{lemma}
\label{prodproperties}
Suppose $\mM$ is a cartesian model category. Then if $A$ and $B$ are cofibrant, so is $A\times B$. If $A\rightarrt B$ and 
$C\rightarrt D$ are cofibrations (resp. trivial cofibrations)
between cofibrant objects then $A\times C\rightarrt B\times D$ is a cofibration
(resp. trivial cofibration). 
\end{lemma}
\begin{proof}
Follow the argument in the first paragraph of the proof of the previous lemma.  
\end{proof}

Condition (DCL) holds for a wide variety of categories,  notably
any presheaf category. 

\begin{lemma}
If $\mM = \presh (\Phi )$ is the category of presheaves of sets over a small category $\Phi$, then
it satisfies condition (DCL).
\end{lemma}
\begin{proof}
Products and colimits are computed objectwise, and these properties hold in $\Sets$. 
\end{proof}

\section{Internal $Hom$}
\label{internalhom}

uppose $\mP$ is a tractable left proper cartesian model category.
In practice, $\mP$ will be the model category $\precat (\mM )$ 
of $\mM$-enriched precategories which we
are going to construct, starting with a tractable left cartesian model
category $\mM$. This is why we use the notation $\mP$ rather than $\mM$ in
the present discussion. 

The underlying category of $\mP$ 
is locally presentable, so commutation of direct product with colimits
yields an internal $\uHom$ (Proposition \ref{prop-ihom}).
For any $A,B$ this is an object  $\uHom (A,B)$ together with a map
$\uHom (A,B)\times A\rightarrt B$ such that for any  
$E\in \mM$, to give a map $E\rightarrt \uHom( A,B)$ is the same
as to give a map $E\times A\rightarrt B$.

The cartesian condition is designed exactly so that the internal $\uHom$ will be compatible with the model structure. 

\begin{theorem}
\label{internalmodel}
Suppose $\mM$ is a tractable cartesian combinatorial model category. Then:
\newline
(a)---the internal $\uHom (A,B)$ exists for any $A,B\in \mM$, and takes pushouts
in the first variable or fiber products in the second variable, to fiber products. For
example, given morphisms $A\rightarrt B$ and $A\rightarrt C$ and any $D$, we have
$$
\underline{Hom} (B\cup ^AC, D) = \underline{Hom} (B, D)
\times _{\underline{Hom}(A,D)}\underline{Hom}(C,D).
$$
\newline
(b)---if $A$ is cofibrant and $B$ is fibrant, then $\uHom (A,B)$ is fibrant;
\newline
(c)---if $A'\rightarrt A$ is a cofibration (resp. trivial cofibration) 
and $B$ is fibrant, 
then the induced map $\uHom (A,B)\rightarrt \uHom (A',B)$ is a fibration (resp. trivial
fibration);
\newline
(d)---if $A$ is cofibrant and $B\rightarrt B'$ is a 
fibration (resp. trivial fibration),
then the induced map $\uHom (A,B)\rightarrt \uHom (A,B')$ is a fibration
(resp. trivial fibration);
\newline
(e)---if $A'\rightarrt A$ (resp. $B\rightarrt B'$) is a weak equivalence between
cofibrant (resp. fibrant) objects then the induced map 
$\uHom (A,B)\rightarrt \uHom (A',B')$ is a weak equivalence.
\end{theorem}
\begin{proof}
Part (a) comes from Proposition \ref{prop-ihom}, and the fiber product formulae
come from the 
adjunction definition of $\uHom$ and the
fact that direct product preserves colimits. 

For (c), suppose $A\rightarrt ^f A'$ is a trivial cofibration and $B$ is fibrant. Then
for any cofibration $E\rightarrt ^hF$, the map 
$$
A\times F \cup ^{A\times E} A'\times E \rightarrt A'\times F
$$
is a trivial cofibration. Hence $B$ satisfies the right lifting property with respect to it.
This translates, by the adjunction property of $\uHom$, to the right lifting property
for 
$$
f^{\ast}: \uHom (A',B)\rightarrt \uHom (A,B)
$$
along $h$.  This shows that $f^{\ast}$ is a trivial fibration. Similarly if $f$ was
a cofibration then $f^{\ast}$ is a fibration. 

For (d), suppose $A$ is cofibrant and $B\rightarrt^f B'$ is a fibration; 
then if $E\rightarrt^h F$ is a trivial
cofibration, the product $A\times E\rightarrt^{1_A\times h} A\times F$ is also a 
trivial cofibration.
Since $f$ is fibrant, it satisfies right lifting with respect to this product map
$1_A\times h$,
which is equivalent to the right lifting property for $\uHom (A,f)$ with respect to
$h$ by the adjunction property of $\uHom$. Thus $\uHom (A,B)\rightarrt \uHom (A,B')$ is fibrant. 
Applied to $B'=\ast$ this says that if $B$ is fibrant and $A$ cofibrant then $\uHom (A,B)$ is fibrant, giving (b). For the other part of (d), note that similarly a trivial fibration is tranformed to
a trivial fibration.

For (e), any weak equivalence between cofibrant objects $A\rightarrt ^f A'$
can be decomposed as $f=pi$ where  $i$ is a trivial
cofibration from $A$ to some $A'$ and $p$ is a trivial fibration from $A''$ to $A'$; 
in turn $p$ admits a section $s:A'\rightarrt A''$ with $ps=1_{A'}$. Now $\uHom$ transforms
(contravariantly) $i$ to $i^{\ast}$ which is a trivial fibration, in particular
invertible in the homotopy category. Thus, in $\Ho (\mP )$ we have
$$
\Ho (s^{\ast})\Ho (i^{\ast})^{-1} \Ho (f^{\ast}) = 
\Ho (s^{\ast})\Ho (i^{\ast})^{-1} \Ho (i^{\ast})\Ho (p^{\ast}) = 
\Ho (s^{\ast}) \Ho (p^{\ast}) = \Ho ((ps)^{\ast}) = 1
$$
which says that $\Ho (f^{\ast})$ admits a left inverse. On the other hand, 
$s$ is also a weak equivalence between cofibrant objects so $\Ho (s^{\ast})$
also admits a left inverse, whereas from the above formula it admits a 
right inverse too. Thus $\Ho (s^{\ast})$ is invertible, which then implies that
$\Ho (f^{\ast})$ is invertible and $f^{\ast}$ is a weak equivalence. 
This was just a contravariant version of the standard argument which appears elsewhere.

If $A$ is cofibrant and $E\rightarrt ^hF$ is a cofibration then
$A\times E\rightarrt A\times F$ is again a cofibration; it follows as usual
from the adjunction property that
a trivial fibration $B\rightarrt B'$ induces a
trivial fibration $\uHom (A,B)\rightarrt \uHom (A,B')$, and then repeating 
an argument similar to the previous one (but covariantly this time) yields that
if $B\rightarrt B'$ is a weak equivalence between fibrant objects, then $\uHom (A,B)\rightarrt \uHom (A,B')$ is a weak equivalence. 
\end{proof}

\section{The enriched category associated to a cartesian model category}
\label{sec-enr-cart}

Keep the a tractable left proper
cartesian model category $\mP$. Then we obtain a $\mP$-enriched category of cofibrant and fibrant objects, denoted
$\Enr (\mP )$ defined as follows. It is in the next higher universe level.
The object class of $\Enr (\mP )$ is defined to be 
$\Ob (\mP _{cf})$, the class of
cofibrant and fibrant objects of $\mP$. For any two such objects $A$ and $B$, put 
$$
\Enr (\mP )(A,B):= \uHom (A,B)\in \mP .
$$
The structural morphisms for the enriched category structure come from 
the standard morphisms for the internal $\uHom$.
We get a structure of $\mP$-enriched category.

Although $\Enr (\mP )$ is strictly associative as a $\mP$-enriched category,
we can also consider it as a weakly $\mP$-enriched category, that is
$$
\Enr (\mP )\in \precat (\mP )
$$
in the notation of Chapter \ref{precat1}. 

An early example using an internal $\uHom$ to obtain a higher categorical structure was
the notion of ``enhanced triangulated category'' of Bondal and Kapranov \cite{BondalKapranov}.


\chapter{Direct left Bousfield localization}
\label{direct1}

Suppose $\mM$ is a model category and $K$ a subset of morphisms. A {\em left Bousfield localization} is a 
left Quillen functor $\mM \rightarrt \mN$ sending elements of $K$ to weak equivalences, universal for this property, and
furthermore which induces an isomorphism of underlying categories and an isomorphism of classes of cofibrations. 
If it exists, it is unique up to isomorphism. 

It is pretty well-known that the left Bousfield localization exists whenever $\mM$ is a left proper combinatorial model category.
We refer to the references and particularly Hirschhorn \cite{Hirschhorn}
and Barwick \cite{Barwick} (see also Rosicky-Tholen \cite{RosickyTholen}) 
for this existence theorem and for some of the main characterizations, statements, details and proofs concerning this notion in general.
Recall that the general definition of $K$-local objects and the localization functor depend on the notion of simplicial mapping spaces.
Of course, simplicial mapping spaces are exactly the kind of thing we are looking at in the present book, but to start with these as basic
building blocks would stretch the notion
of `bootstrapping' pretty far. Therefore, in the
present chapter, we consider a special case of left Bousfield localization in which everything is much  more explicit.

\section{Projection to a subcategory of local objects}

Start with a left proper tractable model category $(\mM , I, J)$, that is a left proper cofibrantly generated model category
such that $\mM$ is locally $\kappa$-presentable for some regular cardinal $\kappa$, and the domains of arrows in $I$ and $J$
are cofibrant.

Suppose we are given a 
subclass of objects considered as a full subcategory $\Rr \subset \mM$,
and a subset of morphisms $K \subset \Mor (\mM )$.
We assume that: 
\newline
(A1)---$K$ is a small subset;
\newline
(A2)---$J\subset K$;
\newline
(A3)---$K\subset \cof (I)$ 	and the domains of arrows in $K$ are cofibrant;
\newline
(A4)---if $X\in \Rr$ and $X\cong Y$ in $\Ho (\mM )$ then $Y\in \Rr$; and
\newline
(A5)---$\inj (K)\subset \Rr$. 
\newline
Say that $(\Rr , K)$ is {\em directly localizing} if in addition to the above conditions:
\newline
(A6)---for all $X\in \Rr$ such that $X$ is fibrant (i.e. $J$-injective), and for any 
$X\rightarrt Y$ which is a pushout by an element of $K$,  there exists $Y\rightarrt Z$ in $\cell (K)$
such that $X\rightarrt Z$ is a weak equivalence.

\begin{lemma}
\label{remove}
Under the above hypotheses, we can remove the requirement that $X$ be fibrant in the direct localizing condition: 
if $X\in \Rr$ and $X\rightarrt Y$ is a pushout by an element of $K$ then there exists $Y\rightarrt Z$ in $\cell (K)$
such that $X\rightarrt Z$ is a weak equivalence. 
\end{lemma}
\begin{proof}
Choose a map $X\rightarrt X'$ in $\cell (J)$, in particular a trivial cofibration,
such that $X'$ is fibrant. 
By invariance of $\Rr$ under weak equivalences (A3), $X'\in\Rr$. Let $Y':= Y\cup ^X X'$, then $X'\rightarrt Y'$ is
again a pushout by the same element of $K$. The above condition now applies: there is a map 
$Y'\rightarrt Z$ in $\cell (K)$ such that $X'\rightarrt Z$ is a weak equivalence. Note that $X\rightarrt Z$ is then
also a weak equivalence. 
On the other hand, $Y\rightarrt Y'$ is in $\cell (J)$ hence also in $\cell (K)$ because $J\subset K$ by (A2). Thus the composition $Y\rightarrt Z$
is in $\cell (K)$. We get the desired properties. 
\end{proof}

In our main examples, $\Rr$ will be the subcategory of precategories which satisfy the Segal conditions. Pelissier called this the subcategory of ``regal objects'' \cite{Pelissier}. 
The main step is to use $K$ to construct a monadic projection from $\mM$ to $\Rr$ up to
homotopy. 

Let $G:\mM \rightarrt \mM$ with $\eta _X:X\rightarrt G(X)$ denote a $K$-injective replacement functor,
such that $\eta _X\in \cell (K)$ for all $X$. This exists by small object argument 
for the locally presentable category $\mM$, Theorem \ref{smallobject}. 

For any $X\in \ob (\mM )$, since  $G(X)\in \inj (K)$, condition (A5) implies that $G(X)\in \Rr$. 

Our first step will be to augment the direct localizing property (A6), from morphisms in $K$ to morphisms in $\cell (K)$.

\begin{proposition}
\label{forward}
Under the above assumptions, if $X\in \Rr$ and $X\rightarrt Y$ is in $\cell (K)$ then 
there exists $Y\rightarrt Z$ also in $\cell (K)$ such that $X\rightarrt Z$ is a weak equivalence.
\end{proposition}
\begin{proof}
Write $Y= \colim _{ i} X_i$ with $X_0=X$, and the colimit ranges over an ordinal $\beta$.
Suppose this is a standard presentation of a cell complex, that is $X_i\rightarrt X_{i+1}$ is
a pushout by an element of $K$, and for a limit ordinal $i$ we have $X_i = \colim _{j<i}X_j$.

We construct a system of morphisms $v_i:X_i\rightarrt Z_i$ (for $i\in \beta$)
together with morphisms $g_{ji}:Z_j\rightarrt Z_i$ for $j\leq i$ forming a transitive system,
and starting with $Z_0=X_0$ (and $v_0$ is the identity), such that 
\newline
---each $v_i$ is an element of $\cell (K)$;
\newline
---for a limit ordinal $i$ we have $Z_i = \colim _{ j<i}Z_j$; 
\newline
---each $Z_i\in \Rr$;
\newline
---each $Z_i\rightarrt Z_{i+1}$ is a weak equivalence, and an element of $\cell (K)$.

For any ordinal $\alpha \leq \beta$ let $\Vv (\alpha )$ denote the set of such compatible collections for $i<\alpha$. 
We have restriction maps $\Vv (\alpha ')\rightarrt \Vv (\alpha )$ for $\alpha \leq \alpha '$. 
If $\alpha$ is a limit ordinal then $\Vv (\alpha ) = \mylim _{ \eta < \alpha} \Vv (\eta )$. 
We will show later that $\Vv (\alpha + 1)\rightarrt \Vv (\alpha )$ is surjective.

Assuming this for now, we show by transfinite induction on $\alpha \leq \beta $ that for any $\eta \leq \alpha$ the map 
$\Vv (\alpha )\rightarrt \Vv (\eta )$ is surjective. Assuming the contrary,
let $\alpha _0$ be the smallest ordinal $\leq \beta$ such that this is not true. Choose an $\eta$ (which we may assume $<\alpha _0$)
such that $\Vv (\alpha _0)\rightarrt \Vv (\eta )$ is not surjective. Either $\alpha _0$ is a limit ordinal, or 
it is a successor ordinal. If $\alpha _0$ is a limit ordinal
then the limit in the above expression can be restricted to 
ordinals bigger than $\eta$:
$$
\Vv (\alpha _0) = \mylim _{\eta \leq \phi  < \alpha_0} \Vv (\phi ).
$$
Each of the maps $\Vv (\phi )\rightarrt \Vv (\eta )$ is surjective, and the
transition maps $\Vv (\phi )\rightarrt \Vv (\phi ')$ in
the inverse system are surjective, so the map 
$$
\mylim _{\eta \leq \phi  < \alpha_0} \Vv (\phi )\rightarrt \Vv (\eta )
$$
is surjective, a contradiction. Suppose $\alpha _0=\alpha +1$ is a successor ordinal. 
Then by the above claim (which we will show next) the map $\Vv (\alpha _0)\rightarrt \Vv (\alpha )$
is surjective; hence by the induction hypothesis for any $\eta \leq \alpha $ the map 
$\Vv (\alpha _0)\rightarrt \Vv (\eta )$
is surjective, and of course it also works for $\eta = \alpha _0$, so this again gives a contradiction.
We have completed the proof of surjectivity of $\Vv (\alpha )\rightarrt \Vv (\eta )$ for all $\eta \leq \alpha \leq \beta$.

Since $\Vv (0)$ is nonempty, this shows that $\Vv (\beta )$ is nonempty, i.e. there exists a system of the required kind. 

We now have to show the claim: that $\Vv (\alpha + 1)\rightarrt \Vv (\alpha )$ is surjective.
We are given data $v_i$ for all $i<\alpha $ and need to construct $v_{\alpha }: X_{\alpha }\rightarrt Z_{\alpha}$.

If $\alpha$ is a limit ordinal, put 
$$
Z_{\alpha} := \colim _{i<\alpha} Z_i.
$$
Let $v_{\alpha}$ be the natural map from $X_{\alpha}\cong \colim _{ i<\alpha} X_i$ to $Z_{\alpha}$.
We claim that this is again in $\cell (K)$, indeed a directed colimit of maps in $\cell (K)$ is again in $\cell (K)$.
The second condition holds automatically by construction. For the fourth condition, note that the maps
$Z_i\rightarrt Z_{\alpha}$ are in $\cell (K)$ because $Z_{\alpha}$ is by construction a transfinite composition of the previous maps
which were in $\cell (K)$, see Lemma \ref{compositioncell}.
Furthermore, the  maps in the system $\{ Z_i\} _{i<\alpha}$ are trivial cofibrations, and a directed limit of trivial cofibrations
satisfies the required lifting property against fibrations, so it is again a trivial cofibration. Thus $Z_i\rightarrt Z_{\alpha}$
are trivial cofibrations, in particular they are weak equivalences. This finishes the fourth condition, and it also gives the
third condition: since $Z_i\in \Rr$ we have $Z_{\alpha}\in \Rr$ by invariance of $\Rr$ under weak equivalences. 

We now treat the other case of the claim: where $\alpha$ is a successor ordinal, say $\alpha = \eta +1$. 
We need to construct $v_{\eta + 1}: X_{\eta + 1}\rightarrt Z_{\eta +1}$ starting from $v_{\eta}$. 
Note that $X_{\eta}\rightarrt X_{\eta + 1}$ is a pushout along an element of $K$. 
Let $W$ be the pushout in the square
$$
\begin{diagram}
X_{\eta} & \rightarr & X_{\eta + 1} \\
\downarr & & \downarr \\
Z_{\eta} & \rightarr & W .
\end{diagram}
$$
The map along the bottom is again a pushout by the same element of $K$ as the map along the top. 
By our inductive hypothesis, $Z_{\eta} \in \Rr$. By the assumed properties (improved in Lemma \ref{remove}) there
is a new $Z_{\eta +1}$ and a map $W\rightarrt Z_{\eta +1}$ in $\cell (K)$ such that $Z_{\eta}\rightarrt Z_{\eta + 1}$ is 
a weak equivalence. Let $v_{\eta +1}$ be the composed map $X_{\eta +1}\rightarrt Z_{\eta + 1}$. This map is in $\cell (K)$ because
the map $X_{\eta +1}\rightarrt W$ is in $\cell (K)$ (since it is a pushout of the map $v_{\eta}$ which was in $\cell (K)$);
and the map $W\rightarrt Z_{\eta +1}$ is in $\cell (K)$ by construction. This gives the first required property. The
second required property doesn't say anything because $\eta +1$ is not a limit ordinal; the
fourth property comes from the above construction, and as usual the third property follows from the condition $Z_{\eta}\in \Rr$
and stability of $\Rr$ under weak equivalences. 

This finishes the proof of the proposition. 
\end{proof}

\begin{corollary}
\label{maincor}
Suppose $f:X\rightarrt Y$ is a morphism in $\cell (K)$, such that $X,Y\in \Rr$. Then $f$ is a weak equivalence. 
\end{corollary}
\begin{proof}
By the proposition, there exists a morphism $\varphi : Y\rightarrt Z$ in $\cell (K)$ such that $\varphi \circ f : X\rightarrt Z$ is a weak equivalence. 
Applying the proposition again, there exists a morphism $\xi : Z\rightarrt W$ in $\cell (K)$ such that $\xi \circ \varphi : Y\rightarrt W$ is
a weak equivalence. Looking at the images of everything in the homotopy category $\Ho (\mM )$ we find that the image of $\varphi$ is
a map with both left and right inverses. It follows that $\varphi$ goes to an isomorphism in the homotopy category, hence $\varphi$ is
a weak equivalence. By $3$ for $2$ we get that $f$ is a weak equivalence. 
\end{proof}

\begin{corollary}
\label{determined}
Under the above assumptions, for an object $X\in \mM$ the following are equivalent:
\newline
(1)---$X\in \Rr$;
\newline
(2)---$\eta _X:X\rightarrt G(X)$ is a weak equivalence;
\newline
(3)---for any map $f:X\rightarrt Y$ in $\cell (K)$ such that $Y\in \inj (K)$, $f$ is a weak equivalence;
\newline
(4)---there exists a map $f:X\rightarrt Y$ in $\cell (K)$ such that $Y\in \inj (K)$ and $f$ is a weak equivalence.
\end{corollary}
\begin{proof}
Suppose $\eta _X$ is a weak equivalence. Note that $G(X)$ is $K$-injective, so by our hypothesis on $\Rr$ we have
$G(X)\in \Rr$. Then since $\Rr$ is supposed to be stable under weak equivalences, we get $X\in \Rr$. 

Suppose $X\in \Rr$. Apply the previous Corollary \ref{maincor} to the map $\eta _X:X\rightarrt G(X)$ which is in $\cell (K)$,
between elements of $\Rr$. It says that $\eta _X$ is a weak equivalence. 

We have now shown $(1)\Leftrightarrow (2)$. It is clear that $(3)\rightarrt (2)\rightarrt (4)$. Suppose (4) with
$f:X\rightarrt Y$ a weak equivalence and $Y\in \inj (K)$. Then $Y\in \Rr$ by (A5) and $X\in \Rr$ by (A4)n which is (1). 
\end{proof}

This corollary says that $\Rr$ is determined by $K$. Thus, we can say that $K$ is {\em directly  localizing}
if there exists a class of objects $\Rr$ satisfying properties (A1)--(A6).  The class $\Rr$ can be
assumed to be defined by conditions (2), (3) or (4) of the corollary.

\begin{lemma}
\label{triangle}
Suppose given a diagram 
$$
\begin{diagram}
X & & \\
\downarr^{ f}  & \rdTeXto  & \\
Y & \rightarr^{g} & Z ,
\end{diagram}
$$
with $Y,Z\in \Rr$ such that $f:X\rightarrt Y$ and $gf : X\rightarrt Z$ are in $\cell (K)$.
Then $g$ is a weak equivalence. 
\end{lemma}
\begin{proof}
Let $U:= Y\cup ^XG(X)$, denote the two morphisms by $a:Y\rightarrt U$ and
$b: G(X)\rightarrt U$ and compose with the map $\eta _U : U\rightarrt G(U)$ to get a diagram
$$
\begin{diagram}
X &\rightarr^{f}& Y \\
\downarr^{\eta _X} & & \downarr_{ \eta _Ua}\\
G(X) & \rightarr^{\eta _Ub} & G(U) .
\end{diagram}
$$
Note here that $\eta _Ub$ is not necessarily equal to $G(b\eta _X)$ (this kind of problem is the difficulty of the present proof). 
By hypothesis $f\in \cell (K)$ so all of the maps in this square are in $\cell (K)$. 

Next, put $V:= Z\cup ^YG(U)$ and let $c:Z\rightarrt V$ and $d:G(U)\rightarrt V$ denote the morphisms. Compose again with $\eta _V$
to get the diagram 
$$
\begin{diagram}
Y &\rightarr^{g}& Z \\
\downarr^{\eta _U a} & & \downarr _{ \eta _Vc}\\
G(U) & \rightarr^{\eta _Vd} & G(V) .
\end{diagram}
$$
The map $c$ comes by pushout from the map $\eta _Ua : Y\rightarrt G(U)$ which is in $\cell (K)$, so $c\in \cell (K)$
and $\eta _Vc\in \cell (K)$. However, we don't know that $g$, $d$ or $\eta _Vd$ are in $\cell (K)$. 

Put these together into a big diagram of the form 
\begin{equation}
\label{rectangle}
\begin{diagram}
X & \rightarr^{f} & Y &\rightarr^{g} & Z \\
\downarr^{ \eta _X} & & \downarr_{\eta _Ua}& & \downarr _{\eta _Vc}\\
G(X) & \rightarr^{\eta _Ub} & G(U) &\rightarr^{\eta_V d} & G(V).
\end{diagram}
\end{equation}
All of the vertical maps are in $\cell (K)$. The horizontal maps in the square on the left are in $\cell (K)$. 
The composition $gf$ along the top is in $\cell (K)$; we would like to show the same for the composition along
the bottom. For this, note that we have a diagram 
$$
\begin{diagram}
X & \rightarr^{f} & Y &\rightarr^{g} & Z \\
\downarr^{ \eta _X} & & \downarr _{ a}& & \downarr \\
G(X) & \rightarr^{b} & U  &\rightarr & Z\cup ^YU = Z\cup ^X G(X) \\
& & \downarr _{ \eta _U}& & \downarr _{ \eta _V}\\
& & G(U)  &\rightarr & G(U)\cup ^U (Z\cup ^YU) = V \\
& & & & \downarr _{ \eta _V} \\
& & & & G(V) .
\end{diagram}
$$
The horizontal map $G(X)\rightarrt Z\cup ^XG(X)$ is a pushout of the morphism $gf\in \cell (K)$ by
the morphism $\eta _X$, so it is in $\cell (K)$. The vertical map from 
$Z\cup ^XG(X)$ to $V$ is a pushout along $\eta _U$, so it is in $\cell (K)$, and the map $\eta _V$
is in $\cell (K)$, so the composed vertical map $Z\cup ^XG(X)\rightarrt G(V)$ is in $\cell (K)$.
We conclude that the map $G(X)\rightarrt G(V)$ is in $\cell (K)$. This is the same as the composition
$\eta _V d \circ \eta _U b$ along the bottom of the previous diagram \eqref{rectangle}. 

This map is in $\cell (K)$ and goes between elements of $\Rr$, so it is a weak equivalence by Corollary \ref{maincor}.
Furthermore, the map $\eta _Ub$ on the bottom left of \eqref{rectangle} is in $\cell (K)$ and goes between elements
of $\Rr$, so it is a weak equivalence. We conclude by $3$ for $2$ that the map $G(U)\rightarrt^{\eta_V d}  G(V)$
is a weak equivalence. 

In our previous diagram \eqref{rectangle} the vertical maps are in $\cell (K)$, and in the $Y$ and $Z$ columns these
maps go between elements of $\Rr$, so the center and right vertical maps are weak equivalences. 
In the previous paragraph we have seen that the bottom of this rightward square is a weak equivalence, so by $3$ for $2$ we
conclude that $g:Y\rightarrt Z$ is a weak equivalence. This proves the lemma. 
\end{proof}

\begin{corollary}
\label{GcellK}
Suppose $f:X\rightarrt Y$ is a morphism in $\cell (K)$. Then $G(f): G(X)\rightarrt G(Y)$ is a weak equivalence. 
\end{corollary}
\begin{proof}
Apply the previous lemma to the triangle 
$$
\begin{diagram}
X & & \\
\downarr^{ \eta _X} & \rdTeXto  & \\
G(X) & \rightarr^{G(f)} & G(Y) .
\end{diagram}
$$
Naturality for the transformation $\eta$ says that the composition $G(f)\circ \eta _X$ is the same as $\eta _Y\circ f$.
We know that $\eta _Y\in \cell (K)$, and $f\in \cell (K)$ by hypothesis, so $G(f)\circ \eta _X\in \cell (K)$.
On the other hand, both $G(X)$ and $G(Y)$ are in $\Rr$, so Lemma \ref{triangle} applies to show that $G(f)$ is a
weak equivalence. 
\end{proof}

\begin{corollary}
\label{GtoGG}
For any object $X\in \mM$, both maps $G(\eta _X)$ and $\eta _{G(X)}$ from $G(X)$ to $G(G(X))$ are weak equivalences. 
\end{corollary}
\begin{proof}
The map $\eta _X$ is in $\cell (K)$ so the previous corollary shows that $G(\eta _X)$ is a weak equivalence. 
The map $\eta _{G(X)}$ is in $\cell (K)$ and goes between two objects in $\Rr$, so it is a weak equivalence by Corollary \ref{maincor}.
\end{proof}

\section{Weak monadic projection}

It is useful to look at the subcategory $\Rr$ and the functor $G$ in general terms.
One can axiomatize their properties, in a homotopy-theoretic analogue of the
discussion of Section \ref{sec-monadic}. 

Suppose $\mM$ is a model category, and $\Rr \subset \mM$ a full subcategory. 
We assume that $\Rr$ is stable under weak equivalences, which means that it comes from
a subset of isomorphism classes in $\Ho (\mM )$. The aim is to apply the general discussion of this section
to the subcategory $\Rr$ of the previous section; however, formally $\Rr$, and the functor $G$ which sometimes shows up below,
are not necessarily those of the previous section. 

A {\em weak monadic projection} from $\mM$ to $\Rr$ is a functor $F: \mM \rightarrt \mM$
together with a natural transformation $\eta _X: X\rightarrt F(X)$,
such that:
\newline
(WPr1)---$F(X)\in \Rr$ for all $X\in \mM$; 
\newline
(WPr2)---for any $X\in \Rr$, $\eta _X$ is a weak equivalence;
\newline
(WPr3)---for any $X\in \mM$, the map $F(\eta _X): F(X)\rightarrt F(F(X))$ is a weak equivalence;
\newline
(WPr4)---if $f:X\rightarrt Y$ is a weak equivalence between cofibrant objects then $F(f): F(X)\rightarrt F(Y)$ is a weak equivalence; and
\newline
(WPr5)---$F(X)$ is cofibrant for any cofibrant $X\in \mM$.

Note that the map in (WPr3) is different from the map $\eta _{F(X)}: F(X)\rightarrt F(F(X))$
which itself is a weak equivalence by (WPr1) and (WPr2). This differentiates the weak situation from Lemma \ref{monadicequal} for the
case of monadic projection considered in Chapter \ref{cattheor1}. Another notable detail is that 
in condition (WPr4) the objects are required
to be cofibrant; this is done in order to make the proof of Corollary \ref{Gweakmonadic} work below. 

Suppose $(F,\eta )$ is a weak monadic projection from $\mM$ to $\Rr$. Let $\Ho (\Rr )$ denote the 
image of $\Rr$ in $\Ho (\mM )$. Then we can construct a monadic projection  
$(\Ho (F), \Ho (\eta ))$ from $\Ho (\mM )$ to $\Ho (\Rr )$. 
Because of this restriction to cofibrant objects in (WPr4), we need to compose with a cofibrant replacement in order to define
$\Ho (F)$. This process, and hence the proof of Lemma \ref{howeakmon} below, can be simplified if $\mM$ is an injective
model category where all objects are cofibrant---and (WPr5) would be superfluous as well.

Let $P:\mM \rightarrt \mM$ be a functor with a natural transformation 
$\xi _X: P(X)\rightarrt X$ such that $P(X)$ is cofibrant and $\xi _X$ is a trivial fibration for all $X\in \mM$.
Then $FP$ is invariant under weak equivalences: if $f:X\rightarrt Y$ is a weak equivalence, we have a 
commutative diagram
$$
\begin{diagram}
P(X) & \rightarr & P(Y) \\
\downarr & & \downarr \\
X & \rightarr & Y
\end{diagram}
$$
such that three sides are weak equivalences. Therefore $P(X)\rightarrt P(Y)$ is a weak equivalence between cofibrant objects,
and by (WPr4), the map $FP(X) \rightarrt FP(Y)$ is again a weak equivalence. It follows that $FP$ descends to a functor 
which we denote by $\Ho (F): \Ho (\mM )\rightarrt \Ho (\mM )$. For any $X\in \mM$, consider the diagram 
$$
X \leftarr^{\xi _X} P(X) \rightarrt^{\eta _{P(X)}} FP(X).
$$
It projects in $\Ho (\mM )$ to a diagram where the first arrow is invertible; we can define
$$
\Ho (\eta )_X := \Ho (\eta _{P(X)})\circ \Ho (\xi _X)^{-1} : \Ho (X)\rightarrt \Ho (FP(X)) =: \Ho (F)(X).
$$
Here we denote also by $\Ho$ the functor from $\mM$ to $\Ho (\mM )$.

\begin{lemma}
\label{howeakmon}
Given a weak monadic projection $(F,\eta )$ and choosing a cofibrant replacement functor $(P,\xi )$,
the collection $\Ho (\eta )$ defined above is a natural transformation, and the pair
$(\Ho (F), \Ho (\eta ))$ is a monadic projection from $\Ho (\mM )$ to $\Ho (\Rr )$.  
\end{lemma}
\begin{proof}
For naturality of $\Ho (\eta )$, suppose $f: X\rightarrt Y$ is a morphism in $\mM$. Then we have a diagram
whose vertical arrows come from $f$,
$$
\begin{diagram}
X &\leftarr^{\xi _X} &P(X)& \rightarr^{\eta _{P(X)}} &FP(X) \\
\downarr & & \downarr & & \downarr \\
Y &\leftarr^{\xi _Y} &P(Y)& \rightarr^{\eta _{P(Y)}} &FP(Y) 
\end{diagram}
$$
which commutes by naturality of $\xi$ and $\eta$. The image of this is a commutative diagram in $\Ho (\mM )$ 
where the first horizontal arrows are isomorphisms; replacing them by their inverses and taking the outer
square we get a commutative diagram 
$$
\begin{diagram}
X &\rightarr^{\Ho (\eta )_X} &\Ho (F)(X) \\
\downarr & & \downarr  \\
Y &\rightarr^{\Ho (\eta )_Y} &\Ho (F)(Y)
\end{diagram}
$$
which shows naturality of $\Ho (\eta )$ with respect to morphisms coming from $\mM$. The same diagram implies naturality with respect to
the inverse of such a morphism, when $f$ is a weak equivalence. These two classes of morphisms generate the morphisms of $\Ho (\mM )$
so $\Ho (\eta )$ is a natural transformation. 

It remains to check the conditions (Pr1)--(Pr3) of a monadic projection. 
Condition (WPr1) implies that $\Ho (F)(X) = \Ho (FP(X))$ is in $\Ho (\Rr )$ for any $X$, giving (Pr1).

For property (Pr2), suppose $X\in \Rr$ (note that $\Rr$ and $\Ho (\Rr )$ have the same objects).
Then $P(X)\rightarrt X$ is a weak equivalence, so by 
the invariance of $\Rr$ with respect to weak equivalences $P(X)\in \Rr$, and condition (WPr2) says that
$\eta _{P(X)}$ is a weak equivalence. Therefore $\Ho (\eta )_X = \Ho (\eta _{P(X)})\circ \Ho (\xi _X)^{-1}$ is an isomorphism.

For property (Pr3), suppose $X\in \mM$. Recall that $\Ho (F)$ is obtained by descending the functor $FP$ to 
the homotopy category. In order to apply this to a composition such as $\Ho (\eta )_X = \Ho (\eta _{P(X)})\circ \Ho (\xi _X)^{-1}$,
we use the fact that $FP$ takes weak equivalences to weak equivalences to say that $FP(\xi _X)$ is a 
weak equivalence, and then 
$$
\Ho (F)(\Ho (\eta )_X) = \Ho (F)\left( \Ho (\eta _{P(X)})\circ \Ho (\xi _X)^{-1}\right)
$$
where the right side is defined to be $\Ho (FP (\eta _{P(X)})) \circ \Ho (FP(\xi _X))^{-1}$.
Consider the commutative diagram
obtained by applying $F$ to the naturality diagram for $\xi$ with respect to $\eta _{P(X)}$:
$$
\begin{diagram}
FP(P(X)) &\rightarr^{FP(\eta _{P(X)})}& FP(FP(X)) \\
\downarr^{ F(\xi _{P(X)})} & & \downarr  _{ F(\xi _{FP(X)})}\\
FP(X) &\rightarr^{F(\eta _{P(X)})} &F(FP(X))
\end{diagram} .
$$
By condition (WPr3), the bottom arrow $F(\eta _{P(X)})$ is a weak equivalence. 
On the other hand, $P(X)$ is cofibrant, so $\xi _{P(X)}$ is a weak equivalence between cofibrant objects; 
by (WPr4), $F(\xi _{P(X)})$ is a weak equivalence. Similarly, (WPr5) says that $FP(X)$ is
cofibrant, so $F(\xi _{FP(X)})$ is a weak equivalence. Three of the four sides of the square
are weak equivalences, so the fourth side $FP(\eta _{P(X)})$ is a weak equivalence. Hence
$\Ho (FP (\eta _{P(X)}))$ is an isomorphism, therefore $\Ho (F)(\Ho (\eta )_X)$ is an isomorphism. This proves (Pr3).
\end{proof}

\begin{proposition}
\label{weakmonadicunique}
Suppose $\Rr \subset \mM$ are as above, and $(F,\eta )$ and $(G,\varphi )$ are two weak monadic  projections
from $\mM$ to $\Rr$. Then, for any cofibrant object $X\in \mM$ the maps $F(\varphi _X): F(X)\rightarrt F(G(X))$
and $G(\eta _X): G(X)\rightarrt G(F(X))$ are weak equivalences, and the diagram of weak equivalences
$$
\begin{diagram}
F(X) & \rightarr^{F(\varphi _X)} & F(G(X)) \\
\downarr^{\varphi _{F(X)}} & & \uparr _{\eta _{G(X)}}\\
G(F(X)) & \leftarr^{G(\eta _X)} & G(X)
\end{diagram}
$$
becomes a commuting diagram of isomorphisms in the homotopy category $\Ho (\mM )$.
\end{proposition}
\begin{proof}
Use Proposition \ref{monadicunique} which is the same as the present
statement but for for monadic projections. In the case where all objects of $\mM$ are cofibrant,
we wouldn't need a cofibrant replacement in the construction for Lemma \ref{howeakmon} above
and the conclusion of the proposition follows directly. Otherwise, we need to unwind the occurrences of
the cofibrant replacement functor in the conclusion.

According to Lemma \ref{howeakmon}, $(\Ho (F), \Ho (\eta ))$ and $(\Ho (G), \Ho (\varphi ))$ are
monadic projections from $\Ho (\mM )$ to $\Ho (\Rr )$. Use the same cofibrant replacement
$(P,\xi )$ to define both of these.  Recall that
$$
\Ho (\varphi )_X := \Ho (\varphi _{P(X)})\circ \Ho (\xi _X)^{-1} 
$$
and 
$$
\Ho (F) (\Ho (\varphi )_X):= \Ho (FP (\varphi _{P(X)}))\circ \Ho (FP(\xi _X))^{-1} .
$$
Assume that $X$ is cofibrant. 
Apply Proposition \ref{monadicunique} to conclude that 
$\Ho (F) (\Ho (\varphi )_X)$ is an isomorphism in $\Ho (\mM )$. Hence the same for 
$\Ho (FP(\varphi _{P(X)})$ so  $FP(\varphi _{P(X)}$ is a weak equivalence. 
Look at $F$ applied to the diagram of naturality for $\xi$ with respect to $\varphi _{P(X)}$
$$
\begin{diagram}
FP(P(X)) & \rightarr^{FP(\varphi _{P(X)})} & FP(GP(X)) \\
\downarr^{ F(\xi _{P(X)})} & & \downarr_{ F(\xi _{GP(X)})}\\
FP(X) & \rightarr^{F(\varphi _{P(X)})} & F(GP(X))
\end{diagram} .
$$
The vertical arrows are obtained by applying $F$ to weak equivalences between cofibrant objects: 
use condition 
(WPr5) for $GP(X)$. The top is a weak equivalence as stated above, so we conclude that the
bottom arrow $F(\varphi _{P(X)})$ is a weak equivalence. Now look at the diagram 
$$
\begin{diagram}
FP(X) & \rightarr^{F(\varphi _{P(X)})} & F(GP(X)) \\
\downarr^{ F(\xi _X)} & & \downarr _{ FG(\xi _{X})}\\
F(X) & \rightarr^{F(\varphi _{X})} & F(G(X))
\end{diagram} .
$$
By the assumption that $X$ is cofibrant, and by (WPr5) for both $F$ and $G$, we get that
the vertical arrows are weak equivalences. The top is a weak equivalence by the previous square diagram,
so we conclude that the bottom arrow $F(\varphi _X)$ is a weak equivalence. 

The other statement that $G(\eta _X)$ is a weak equivalence, is obtained by symmetry.

The commutativity of the square diagram follows from the corresponding part of Proposition 
\ref{monadicunique},
using the same technique as above for removing occurences of $P$. 
\end{proof}

\section{New weak equivalences}

We now go back to the situation of the first section, with a 
combinatorial model category $\mM$ with subcategory $\Rr \subset \mM$ and 
set of morphisms $K$ satisfying axioms (A1)--(A6), and $K$-injective replacement functor $(G,\eta )$.

\begin{lemma}
\label{Gweakmonadic}
The pair $(G,\eta )$ is a weak monadic projection from $\mM$ to $\Rr$.
\end{lemma}
\begin{proof}
For (WPr1), if $X\in \mM$ then by definition of $G$ we have $G(X)\in \inj (K)$, but
$\inj (K)\subset \Rr$ by condition (A5). 

For (WPr2), suppose $X\in \Rr$. 
Then $\eta _X$ is a weak equivalence by Corollary \ref{determined}. 

For (WPr3), for any $X\in \mM$ we have $G(\eta _X)$ a weak equivalence by the previous corollary. 

For (WPr4), we first consider the case when $f$ is a trivial cofibration between cofibrant objects. 
Choose a factorization 
$$
X\rightarrt^{h} Z \rightarrt^{g} Y
$$
such that $h\in \cell (J)$ and $g\in \inj (J)$. Then, since $f\in \cof (J)$ there is a lifting
$s:Y\rightarrt Z$ such that $sf = h$ and $gs = 1_Y$. Thus, $f$ is a retract of $h$ in the category of objects under $X$. 
Applying the functor
$G$ we get that $G(f)$ is a retract of $G(h)$ in the category of objects under $G(X)$. On the other hand, $J \subset K$ by
hypothesis (A2), so $\cell (J)\subset \cell (K)$ and $h\in \cell (K)$. By  Corollary \ref{GcellK}, 
$G(h)$ is a weak equivalence. But weak equivalences are closed under retracts, so $G(f)$ is a weak equivalence. This
treats the case of a trivial cofibration.

Now for a general weak equivalence $f:X\rightarrt Y$ between cofibrant objects, 
consider the map $f\sqcup 1_Y : X\cup ^{\emptyset} Y \rightarrt Y$. Choose a factorization
$$
X\cup ^{\emptyset} Y \rightarrt^{h\sqcup s} Z \rightarrt^{g} Y
$$
such that $h\sqcup s$ is an $I$-cofibration and $g$ is a trivial fibration. Using the fact that $X$ and $Y$
are cofibrant objects, and also that $f$ is a weak equivalence,
we get that both maps $h: X\rightarrt Z$ and $s: Y\rightarrt Z$ are trivial cofibrations between cofibrant objects.
Hence $G(h)$ and $G(s)$ are weak equivalences. But $G(g)$, being a left inverse to $G(s)$, is therefore a weak equivalence,
so $G(f)=G(g)G(h)$ is a weak equivalence as desired to show (WPr4).

To show (WPr5), note simply that $\eta _X:X\rightarrt G(X)$ is in $\cell (K)$, so if $X$ is cofibrant then $G(X)$
is cofibrant too. 
\end{proof}

Fix a cofibrant replacement functor and natural transformation $(P,\xi )$. We say that a morphism $f:X\rightarrt Y$ is
a {\em new weak equivalence} if $GP(f):GP(X)\rightarrt GP(Y)$ is a weak equivalence. Using the theory of the previous
section, this notion depends only on $\Rr$ and not on $K$. In passing we also show that it doesn't depend on $P$. 

\begin{lemma} 
\label{nweiff}
Suppose $(H,\psi )$ is a monadic projection from $\Ho (\mM )$ to $\Ho (\Rr )$. Then $f$ is a new weak equivalence
if and only if $H(\Ho (f))$ is an isomorphism.
\end{lemma}
\begin{proof}
By Lemma \ref{howeakmon}, $(\Ho (G), \Ho (\eta ))$ defined using the original cofibrant replacement $(P,\xi )$ is also 
a monadic projection from $\Ho (\mM )$ to $\Ho (\Rr )$.
Applying Proposition \ref{monadicunique} to compare these, gives that $H(\Ho (f))$ is an isomorphism if and only if $\Ho (G)(f)$
is an isomorphism.
In turn, $\Ho (G)(f) = \Ho (GP(f))$ is an isomorphism if and only if $GP(f)$ is a weak equivalence.
\end{proof}

We can choose $(H,\psi )= (\Ho (G),\Ho (\eta ))$, so we can say that $f$ is a new weak equivalence if and only if
$\Ho (G)(\Ho (f))$ is an isomorphism.

If $(F,\varphi )$ is any other weak monadic projection and 
$(P',\xi ')$ is a cofibrant replacement functor not necessarily the same as $(P,\xi )$ then
$H=\Ho '(F)$ and $\psi = \Ho '(\varphi )$ defined as above using $(P', \xi ')$ form
a monadic projection 
from $\Ho (\mM )$ to $\Ho (\Rr )$. So we can also say that $f: X\rightarrt Y$ is a new weak equivalance 
if and only if $\Ho '(F)(\Ho (f))$ is an isomorphism, which is if and only if $FP'(f)$ is a weak equivalence.

\begin{corollary}
\label{nweindep}
The notion of new weak equivalence depends only on $\Rr\subset \mM$.
\end{corollary}
\begin{proof}
The condition of the lemma clearly depends only on $\Rr$. 
\end{proof}

\begin{corollary}
\label{oldwenew}
If $f$ is a weak equivalence in the original model structure for $\mM$, then it is a new weak equivalence.
\end{corollary}
\begin{proof}
In this case $\Ho (f)$ is an isomorphism, so after application of a functor $\Ho (G)(\Ho (f))$ remains a weak equivalence.
\end{proof}

\begin{corollary}
\label{newwecofib}
Suppose $f:X\rightarrt Y$ is a morphism between cofibrant objects. Then it is a new weak equivalence, if and only if $G(f): G(X)\rightarrt G(Y)$
is an old weak equivalence. If $(F,\varphi )$ is a different weak monadic projection then $f$ is a new weak equivalence if and only if $F(f)$
is a weak equivalence. 
\end{corollary}
\begin{proof}
If $X$ and $Y$ are cofibrant then $\Ho (F)(f)$ is isomorphic to the projection 
to $\Ho (\mM )$ of the  morphism $F(f)$. This applies in particular to $F=G$.
Conclude by using Lemma \ref{nweiff}. 
\end{proof}

The same can be said without specifying a full functor, but just looking at $X$ and $Y$. 

\begin{corollary}
\label{weloccrit}
Given a morphism $f:X\rightarrt Y$ between cofibrant objects and a diagram
$$
\begin{diagram}
X & \rightarr^{a} & A \\
\downarr^{ f} && \downarr _{ p} \\
Y & \rightarr^{b} & B
\end{diagram}
$$
with $a,b\in \cell (K)$ and $A,B\in \Rr$, the map $f$ is a new weak equivalence if and only if
$p$ is a weak equivalence for the original model structure. 
\end{corollary}
\begin{proof}
Apply $G$ to the
above diagram. Corollary \ref{GcellK} says that the resulting horizontal
arrows are old weak equivalences. Thus, by Corollary \ref{newwecofib},
the map $f$ is a new weak equivalence if and only if $G(A)\rightarrt^{G(p)} G(B)$ is an 
old weak equivalence. However, in the naturality square for $\eta$ with respect to $p$
$$
\begin{diagram}
A&\rightarr^{\eta _A}& G(A)\\
\downarr && \downarr \\
B&\rightarr^{\eta _B}& G(B)
\end{diagram}
$$
the horizontal arrows are old weak equivalences, since $A$ and $B$ are cofibrant
objects in $\Rr$ (see Corollary \ref{determined}). Putting these statements together
using 3 for 2 gives the proof that $f$ is a new weak equivalence if and only if 
$p$ is an old weak equivalence. 
\end{proof}

\begin{lemma}
\label{newretract32}
The notion of new weak equivalence is stable under retracts, and satisfies 3 for 2.
\end{lemma}
\begin{proof}
It is the pullback of the notion of isomorphism, by the monadic projection $\Ho (G)$.
\end{proof}

\section{Invariance properties}
\label{sec-invariance}

The transfinite consequence of the
left properness condition, Proposition \ref{prop-transfiniteleftproper}, allows us to prove invariance of the new trivial cofibrations with respect to
transfinite composition. Say that a morphism is
a {\em new trivial cofibration} if it is a cofibration and a new weak equivalence.

\begin{lemma}
\label{transfinite}
Suppose
$\{ X_n \} _{n\leq \beta}$ is a continuous
transfinite sequence indexed by an ordinal $\beta$, such that
for any limit ordinal $n$ we have $X_m = \colim _{n<m}X_n$, and for any $n$ with $n+1\leq \beta$ we have
that $X_n \rightarrt X_{n+1}$ is a new trivial cofibration. Then $X_0\rightarrt X_{\beta}$ is a new trivial cofibration.
\end{lemma}
\begin{proof}
Assume first that $X_0$ and hence all the $X_n$ are cofibrant. 
Choose a sequence $\{ Z _n \} _{n\leq \beta}$ and a morphism of sequences $X_n\rightarrt Z_n$ in the
following way: let $Z_0:= G(X_0)$; if $n+1\leq \beta$ and the morphisms up to $X_n\rightarrt Z_n$ are chosen,
let $Z_{n+1}:= G(Z_n \cup ^{X_n}X_{n+1})$; and 
if $m$ is a limit ordinal and the morphisms are chosen for all $n<m$, 
put $Z_m := \colim _{n<m}Z_n$. This defines the full sequence by transfinite induction. 
Furthermore, by induction we see that the morphisms $X_n\rightarrt Z_n$ are in $\cell (K)$.
This is clear at $n=0$ and at $n+1$ if we know it at $n$, by the definition. In the case where $m$ is a limit ordinal,
note that $Z_m = \colim _{n<m}(X_m \cup ^{X_n}Z_n)$ using the continuity 
property of $X_{\cdot}$. The transition maps in this colimit are 
of the form 
$$
\begin{diagram}
X_m \cup ^{X_n}Z_n  =  X_m \cup ^{X_{n+1}}(Z_n \cup ^{X_n}X_{n+1}) \\
\downarr \\
X_m \cup ^{X_{n+1}}Z_{n+1}= X_m \cup ^{X_{n+1}}G(Z_n \cup ^{X_n}X_{n+1}) ,
\end{diagram}
$$
which is in $\cell (K)$. By definition $\cell (K)$ is closed under transfinite composition, so the map
$X_m\rightarrt Z_m$ is in $\cell (K)$ as claimed. It follows 
(from Corollaries \ref{GcellK} and \ref{newwecofib} together using the fact that
$X_m$ are cofibrant) that the maps $X_n\rightarrt Z_n$ are new weak equivalences. 

We now argue by induction that the $Z_n$ are in $\Rr$ and the maps $Z_0\rightarrt Z_n$ are
old trivial cofibrations. If it is known for $Z_n$ then the map $Z_n\rightarrt Z_{n+1}$ is a new weak equivalence between fibrant objects in $\Rr$,
so by Corollary \ref{Rnweowe} it is an old weak equivalence, thus $Z_0\rightarrt Z_n$ is an old weak equivalence. If $m$ is a limit ordinal and we know the statement for
all smaller $n<m$,
then the transition maps in the system $\{ Z_n\} _{n<m}$ are old trivial cofibrations; so 
the map $Z_0\rightarrt Z_m = \colim _{n<m}Z_n$ is an old trivial cofibration (trivial cofibrations are closed under transfinite composition 
as can be seen from their characterization by the lifting property with respect to fibrations). 
Since $Z_0\in \Rr$ we get that $Z_n\in \Rr$ by the invariance property (A4). This completes the inductive proof for
the statement given at the start of the paragraph.

Applying it to $n=\beta$ we conclude that $Z_0\rightarrt Z_{\beta}$ is a new trivial cofibration, proving the lemma under the beginning assumption
that $X_0$ was cofibrant.

Recall that $\mM$ is left proper, which implies a left properness
statement for transfinite compositions also, Proposition \ref{prop-transfiniteleftproper}. 

In the arbitrary case, choose a sequence $\{ Y _n \} _{n\leq \beta}$ with a morphism of sequences 
consisting of old weak equivalences $p_n:Y_n\rightarrt X_n$,
such that the $Y_n$ are cofibrant and the  transition morphisms $Y_n\rightarrt Y_{n+1}$ are cofibrations. This is done by using the factorization
into $\inj (I)\circ \cell (I)$ at each stage. At a limit ordinal $m$ given the choice for all $n<m$, we set
$Y_m:= \colim _{n<m}Y_n$. By the transfinite left properness property $Y_m\rightarrt X_m$ is an old weak equivalence. 
This allows us to make the inductive choice of the sequence $Y_{\cdot}$. 
Then, the $Y_n\rightarrt Y_{n+1}$ are new trivial cofibrations  by Lemma \ref{newretract32}
and the objects $Y_n$ are cofibrant, so the first part of
the proof applies: $Y_0\rightarrt Y_{\beta}$ is a new weak equivalence. By \ref{newretract32} again, $X_0\rightarrt X_{\beta}$
is a new weak equivalence. 
\end{proof}

Using the left properness hypothesis on $\mM$ we can show 
that morphisms in $\cell (K)$ are
new weak equivalences, and also improve the earlier
criteria by removing the conditions that the morphism goes
between cofibrant objects.

\begin{proposition}
\label{cellKnwe}
A morphism $f:X\rightarrt Y$ in $\cell (K)$ is a new weak equivalence.
\end{proposition}
\begin{proof}
Corollaries \ref{GcellK} and \ref{newwecofib} together immediately imply this when $X$ is cofibrant. However this is not
sufficient in general, because new weak equivalences are defined using a cofibrant
replacement. Of course, in case $\mM$ satisfies the additional hypothesis that all objects 
are cofibrant, then this part of the argument (like many others) is considerably
simplified. 

To prove that $\cell (K)$ is contained in the new weak equivalences,
in view of the closure of new trivial cofibrations 
under transfinite composition (Lemma \ref{transfinite} above), 
it suffices to treat morphisms which are pushouts along
a single arrow in $K$. Suppose $X\leftarr U\rightarrt V$ is a 
diagram with $U\rightarrt V$ in $K$, then we need to show that
$X\rightarrt Y:=X\cup ^UV$ is a new
weak equivalence. 
The source $U$ 
is assumed to be cofibrant by Condition (A3). Choose a cofibrant replacement via a trivial fibration $P\rightarrt X$. 
The map from $U$ lifts to a map $U\rightarrt P$. Put $Z:= P\cup ^UV$,
then we have a cocartesian square
$$
\begin{diagram}
P& \rightarr & Z \\
\downarr & & \downarr \\
X & \rightarr & Y
\end{diagram}
$$
where the left vertical arrow is an old weak equivalence. By left properness of the
original model structure, the right vertical arrow is also an old weak equivalence.
The top  map is in $\cell (K)$ and goes between cofibrant objects, so as mentioned above 
it is a new weak equivalence. 
Now 3 for 2 given by Lemma \ref{newretract32} implies that the bottom map is a
new weak equivalence.  
\end{proof}

We can now obtain a criterion for a morphism to be a weak equivalence.
Notice that this condition coincides with the definition of weak equivalences (PG)
used in the construction of a model category by pseudo-generating sets in Section 
\ref{sec-pseudogen}. 

\begin{corollary}
\label{nwecriterion}
A morphism $f:X\rightarrt Y$ is a new weak equivalence if and only if there exists a diagram 
$$
\begin{diagram}
X & \rightarr^{a} & A \\
\downarr^{ f} && \downarr _{ g} \\
Y & \rightarr^{b} & B
\end{diagram}
$$
such that $a$ and $b$ are in $\cell (K)$ and $g \in \inj (I)$. 
\end{corollary}
\begin{proof}
If such a diagram exists, then $g$ is an old weak equivalence  hence a new one; also $a$ and $b$
are new weak equivalences by Corollary \ref{cellKnwe}, so by 3 for 2 for the new weak equivalences, given in the previous
corollary, we get that $f$ is a new weak equivalence.

Suppose $f$ is a new weak equivalence. 
Fix a cofibrant replacement functor and natural transformation $(P,\xi )$. Then 
$GP(f):GP(X)\rightarrt GP(Y)$ is an old weak equivalence. The map 
$P(X)\rightarrt GP(X)$ is in $\cell (K)$. Let $X'$ be the pushout in the cocartesian diagram
$$
\begin{diagram}
P(X) & \rightarr^{\eta _{P(X)}} & GP(X)\\
\downarr^{ \xi _X} && \downarr _{ r} \\
X & \rightarr^{u} & X'
\end{diagram}
$$
and similarly let
$Y'$ be the pushout in the cocartesian diagram
$$
\begin{diagram}
P(Y) & \rightarr^{\eta _{P(Y)}} & GP(Y)\\
\downarr^{ \xi _Y} && \downarr _{ s} \\
Y & \rightarr^{v} & Y' .
\end{diagram}
$$
The pushout maps $u$ and $v$ are in $\cell (K)$.
The maps $\xi _X$ and $\xi _Y$ are old weak equivalences, and $\eta _{P(X)}$ and 
$\eta _{P(Y)}$ are old cofibrations, so the left properness hypothesis on the old model
structure implies that the maps $r: GP(X)\rightarrt X'$ and 
$s:GP(Y)\rightarrt Y'$ are weak equivalences. The pushouts fit into a commutative cube,
and the condition that $GP(f):GP(X)\rightarrt GP(Y)$ is an old weak equivalence
implies, by 3 for 2 in the old model structure, that the map $X'\rightarrt Y'$ is an
old weak equivalence. We can factor this map as
$$
X'\rightarrt^{h} X'' \rightarrt^{g} Y'
$$
such that $h\in \cell (J)$ and $g$ is an old fibration; but it is an old weak
equivalence too so $g\in \inj (I)$. We have obtained the desired diagram
$$
\begin{diagram}
X & \rightarr^{hu} & X'' \\
\downarr^{ f} && \downarr _{ g} \\
Y & \rightarr^{v} & Y'
\end{diagram}
$$
with $hu\in \cell (K)$ because it is the composition of $u\in \cell (K)$ with 
$h\in \cell (J)\subset \cell (K)$, with $v\in \cell (K)$, and with $g\in \inj (I)$.
\end{proof}

We have the following criterion for new trivial cofibrations, which is the
same as in Proposition \ref{furtherinfo}. We repeat the proof here in order
to verify that it works in our current situation. 

\begin{corollary}
\label{ntcfcriterion}
A cofibration $X\rightarrt^f Y$ is a new trivial cofibration if and only if 
it fits into a diagram 
$$
\begin{diagram}
X & \rightarr^{a} & A \\
\downarr^{ f} && \uparr^s \downarr _{ g} \\
Y & \rightarr^{b} & B
\end{diagram}
$$
commutative using either arrow on the right, with $gs = 1_B$, such that
$a$ and $b$ are in $\cell (K)$.
\end{corollary}
\begin{proof}
The same as the first part of Proposition \ref{furtherinfo},
using some things that we already know: new weak equivalences
are closed under retracts and satisfy 3 for 2 by Lemma \ref{newretract32},
and they contain $\cell (K)$ by Proposition \ref{cellKnwe}.
\end{proof}

\begin{corollary}
\label{Rnweowe}
A morphism $f:X\rightarrt Y$ such that $X,Y\in \Rr$, is a new weak equivalence
if and only if it is a weak equivalence in the original model structure. 
\end{corollary}
\begin{proof}
The ``if'' direction is given by Corollary \ref{oldwenew}. 
Suppose given a new weak equivalence 
$f$ such that $X,Y\in \Rr$. Then $P(f):P(X)\rightarrt P(Y)$
is a new weak equivalence between cofibrant objects which are again in $\Rr$ (since
$\Rr$ is invariant under weak equivalences). In the diagram 
$$
\begin{diagram}
P(X) & \rightarr & P(Y) \\
\downarr & & \downarr \\
GP(X) & \rightarr & GP(Y) 
\end{diagram}
$$
the bottom arrow is an old weak equivalence by Corollary \ref{newwecofib}. 
The vertical arrows are old weak equivalences by Corollary \ref{determined}, (1)$\rightarrt$(2). By 3 for 2 in the original model structure, the top arrow
$P(f)$ is an old weak equivalence, hence $f$ is also.  
\end{proof}

\section{New fibrations}

Recall that a new trivial cofibration is a cofibration (which means the same thing in the
original and new structures) and also a new weak equivalence. 
Say that a morphism  is a {\em new fibration}
if it satisfies the right lifting property with respect to new trivial cofibrations.

\begin{lemma}
\label{newtrivcof}
The class of new trivial cofibrations is closed under composition and retracts, and contains
$\cof (K)$. In particular, trivial fibrations from the original model structure,
which are exactly $\cof (J)$, are also new trivial fibrations. 
\end{lemma}
\begin{proof}
By Lemma \ref{newretract32}, the new weak equivalences are closed under retracts and compositions; the
same holds for the cofibrations, so the intersection of these classes is closed under retracts and compositions.
Elements of $\cell (K)$ are cofibrations by (A3) and new weak equivalences by \ref{GcellK} so they are 
new trivial cofibrations, and $\cof (K)$ is the closure of $\cell (K)$ under right retracts. Finally,
since $J\subset K$ by assumption (A2) we get that $\cof (J)\subset \cof (K)$.
\end{proof}

The reader might hope that $\cof (K)$ is equal to the class of new trivial cofibrations; however, in general
it will be smaller and the next section below will be needed to remedy this problem.

\begin{corollary}
\label{newfibold}
A new fibration is also a fibration
in the old model structure. 
\end{corollary}
\begin{proof}
A new fibration satisfies right lifting with respect to 
the class of old trivial cofibrations, since that is contained in the class of new ones. 
\end{proof}

\begin{lemma}
A morphism which is a new fibration and a new weak equivalence, 
is a trivial fibration in the original model structure, in particular it satisfies the right lifting property
with respect to cofibrations. 
\end{lemma}
\begin{proof}
Suppose $f:X\rightarrt Y$ is a new fibration and a new weak equivalence. 
Choose a factorization $f=gh$ in the original model category structure $X\rightarrt^{h}Z\rightarrt^{g}Y$
with $h$ a cofibration and $g$ a trivial fibration. Then $g$ is an old weak equivalence, hence it is
a new weak equivalence by Corollary \ref{oldwenew}. By 3 for 2 for the new weak equivalences Lemma \ref{newretract32}
we  conclude that $h$ is a new weak equivalence, hence by definition it is a new trivial cofibration.
The condition that $f$ be a new fibration implies that it satisfies lifting with respect to $h$,
so there is a morphism $u:Z\rightarrt X$ such that $uh = 1_X$ and $fu = g$. In this way, $f$ becomes a retract of $g$
in the category of objects over $Y$. By closure of the original weak equivalences under retracts, $f$ is an original weak equivalence.
Since it is also a fibration in the original structure by Corollary \ref{newfibold}, we conclude that it is a trivial fibration
in the original model structure. 
\end{proof}

\begin{corollary}
\label{oldnewtfib}
The class of new trivial fibrations, defined as the intersection of the new fibrations and the new weak equivalences,
is equal to the original class of trivial fibrations.
\end{corollary}
\begin{proof}
One inclusion is given by the preceding lemma. 
In the other direction, suppose $f:X\rightarrt Y$ is an original trivial fibration. Then it satisfies lifting with
respect to cofibrations, in particular with respect to new trivial cofibrations. So it is a new fibration.
It is an old weak equivalence, hence a new one by Corollary \ref{oldwenew}. Thus it is a new trivial fibration. 
\end{proof}

We can sum up the preceding discussion as follows.

\begin{scholium}
The three classes of morphisms: new weak equivalences, new fibrations, and cofibrations which are the same as the old ones,
generate the notions of new trivial fibration and new trivial cofibration by intersection of new fibrations and cofibrations,
with new weak equivalences. 
All of these classes are closed under composition and retracts, and new weak equivalences satisfy 3 for 2.
If 
$$
\begin{diagram}
X & \rightarr^{a} & U \\
\downarr^{ i} && \downarr _{ p} \\
Y & \rightarr^{b} & V
\end{diagram}
$$
is a diagram, and if either $i$ is a new trivial cofibration and $p$ a new fibration; or $i$ a cofibration and
$p$ a new trivial fibration, then there exists a lifting $f:Y\rightarrt U$ with $fi=a$ and $pf =b$. 
\end{scholium}

We can add that the original model structure provides one of the two required factorizations:
if $f:X\rightarrt Y$ is any map, then it can be factored as 
$$
f=gh:X\rightarrt^{h}Z\rightarrt^{g}Y
$$
where $h$ is a cofibration and $g$ is an old or new trivial fibration (these being the same \ref{oldnewtfib}).

\section{Pushouts by new trivial cofibrations}

The class of new trivial cofibrations, defined to be the intersection of
the new weak equivalences with the cofibrations, is closed under pushout: 

\begin{lemma}
\label{pushout}
If 
$$
\begin{diagram}
X & \rightarr^{f} & Y \\
\downarr^{ a} && \downarr _{ b} \\
Z & \rightarr^{g} & W
\end{diagram}
$$
is a pushout square with $f$ a new trivial cofibration, then $g$ is a new trivial cofibration.
\end{lemma}
\begin{proof}
First reduce to the case when $X$  and hence $Y$ are cofibrant objects. 
Choose a cofibrant replacement $p: X'\rightarrt X$ with $p\in \inj (I)$, and a factorization of the map from $X'$ to 
$Y$ to give a commutative square 
$$
\begin{diagram}
X' & \rightarr^{f'} & Y' \\
\downarr^{ p} && \downarr _{ q} \\
X & \rightarr^{f} & Y
\end{diagram}
$$
such that $f'$ is a cofibration and $q\in \inj (I)$. Note that $f'$ is a  new trivial cofibration, by 3 for 2 for new weak equivalences
\ref{newretract32}. 

By the hypothesis that $\mM$ is left proper, the map $Y'\rightarrt X\cup ^{X'}Y'$ is an old weak equivalence, so by 3 for 2 the
map $X\cup ^{X'}Y'\rightarrt Y$ is an old weak equivalence. We get a map 
$$
Z\cup ^{X'}Y' = Z \cup ^X (X\cup ^{X'}Y')\rightarrt Z\cup ^{X}Y =W.
$$
The maps from $X$ to $X\cup ^{X'}Y'$ and $Y$ are cofibrations, so Lemma \ref{leftpropinvariance} applies:
the map $Z\cup ^{X'}Y' \rightarrt Z\cup ^{X}Y = W$ is an old weak equivalence. 
Assume we know the statement of the lemma for cofibrant objects; then $Z\rightarrt Z\cup ^{X'}Y'$ is a
new weak equivalence. The composition of these two maps is the same as the map $g:Z\rightarrt W$
so by  3 for 2 for new weak equivalences
\ref{newretract32}, we get that $g$ is a new weak equivalence. 

This reduces the statement of the lemma to the case where $X$ and $Y$ are cofibrant, which we now suppose.
Let $\eta _X: X\rightarrt G(X)$ be the map in $\cell (K)$ to $G(X)\in \inj (K)$. Let $V:= G(Y\cup ^X G(X))$.
The map $Y\rightarrt V$ is in $\cell (K)$ and $V\in \inj (K)$. Taking the pushout of the commutative square
$$
\begin{diagram}
X & \rightarr^{f} & Y \\
\downarr^{ \eta _X} && \downarr _{ s} \\
G(X) & \rightarr^{t} & V
\end{diagram}
$$
along the map $a:X\rightarrt Z$ gives a square 
$$
\begin{diagram}
Z & \rightarr^{g} & Z\cup ^X Y \\
\downarr^{ Z\cup ^X \eta _X} && \downarr _{ Z\cup ^X s} \\
Z\cup ^X G(X) & \rightarr^{Z\cup ^X t} & Z\cup ^X V
\end{diagram}
$$
and this fits into a cube with the previous diagram. Since $f$ is a new weak equivalence between cofibrant
objects, the vertical maps in the first square are in $\cell (K)$, and
the objects $G(X)$ and $V$ are in $\inj (K)\subset \Rr$, 
the map $t:G(X)\rightarrt V$ is an old weak equivalence by Corollary \ref{weloccrit}.
It is also a cofibration, so it is an old trivial cofibration. Its pushout 
$Z\cup ^X G(X) \rightarrt^{Z\cup ^X t}  Z\cup ^X V$ is therefore an old trivial cofibration.
On the other hand, the vertical maps in the second square, being pushouts of the vertical maps which were in $\cell (K)$
for the first square, are also in $\cell (K)$. In particular, these maps are in $\cof (K)$ hence new weak equivalences by 
Lemma \ref{newtrivcof}. 
The bottom map is an old,  hence a new weak equivalence. By 3 for 2, Lemma \ref{newretract32}
the map $g: Z\rightarrt Z\cup ^X Y = W$ is a new weak equivalence. This proves the lemma since $g$ is clearly a cofibration.
\end{proof}

\section{The model category structure}

We can now put together everything above to obtain the new model category structure on $\mM$,
as an application of Theorem \ref{recog} in the previous chapter,
which was our version of Smith's recognition theorem as exposed in \cite{Barwick}.

\begin{theorem}
\label{directLBL}
Suppose $\mM$ is a left proper tractable model category, and $(\Rr , K)$ is a directly localizing
system according to axioms (A1)--(A6) at the start of the chapter. 
The classes of original cofibrations, new weak equivalences, and new fibrations constructed above
provide $\mM$ with a
structure of closed model category, cofibrantly generated, and indeed 
combinatorial and even tractable. It is left proper. Furthermore,
this structure is the left Bousfield localization of $\mM$ by the original set of maps $K$.

The fibrant objects are the $K$-injective objects, and a morphism $W\rightarrt U$ 
to a fibrant object is a fibration if and only if it is in $\inj (K)$.  
\end{theorem}
\begin{proof}
We will be applying the discussion of Section \ref{sec-pseudogen} in the previous
chapter, about pseudo-generating sets. 
We are given subsets $I$ and $K$; the cofibrations are $\cof (I)$, and by Corollary \ref{nwecriterion}
the new weak equivalences are exactly the class defined by the condition (PG). 
The trivial cofibrations are as defined in Section \ref{sec-pseudogen}, 
hence the fibrations too. 
We verify the axioms
for pseudo-generating sets.  
\newline
(PGM1)---The hypothesis that $\mM$ is locally presentable, is contained in the tractability hypothesis of the present statement;
 $I$ is a small set as it is one of the generating sets for the old structure, and $K$ is a small set by (A1).
\newline
(PGM2)---the domains of arrows in $I$ are cofibrant, by the assumption that the old structure is tractable
(and $I$ is part of a tractable pair of generating sets); the domains of arrows in $K$ are cofibrant, and $K\subset \cof (I)$,  by (A3); 
\newline
(PGM3)---the class of weak equivalences is closed under retracts by Lemma \ref{newretract32};
\newline
(PGM4)---the class of weak equivalences satisfies 3 for 2 by Lemma \ref{newretract32};
\newline
(PGM5)---the class of trivial cofibrations is closed under pushouts by Lemma \ref{pushout};
\newline
(PGM6)---the class of trivial cofibrations is closed under transfinite composition by Lemma \ref{transfinite}.

By Theorem \ref{recog}, $\mM$ with these classes of morphisms is a tractable model category. 

To prove left properness, suppose given a cocartesian diagram 
$$
\begin{diagram}
X &\rightarr ^u & U \\
\downarr ^f & &\downarr _g \\
Y & \rightarr ^v & V
\end{diagram}
$$
where $f$ is a new weak equivalence and $u$ is a cofibration. Consider a fibrant
replacement $h:Y\rightarrt Y'$ and put $V':= Y'\cup ^YV = Y' \cup ^XU$.
Since $h$ is a new trivial cofibration, so is the map $V\rightarrt V'$.  
Factor the composed map $hf$ as $X\rightarr^i X'\rightarr^pY'$ where
$i$ is a new trivial cofibration and $p$ is a trivial fibration. 
In particular $p$ is a weak equivalence in the old model structure. 
The map $U\rightarrt X'\cup ^XU$ is a trivial cofibration by Lemma \ref{pushout},
and left properness of the original model structure implies that 
$$
X'\cup ^XU\rightarrt Y'\cup ^{X'}(X'\cup ^XU) = V'
$$
is an old weak equivalence (hence a new one). Composing these we get that $U\rightarrt V'$
is a new weak equivalence, and by 3 for 2 we conclude that $g$ is a new weak equivalence. 

The resulting model category is the left Bousfield localization of $\mM$ along the
subset $K$. The morphisms of $K$ go to weak equivalences in the new structure.
On the other hand, by the criterion
of Corollary \ref{nwecriterion}, given a model structure whose class of weak equivalences
contains the old ones plus $K$, that class must also contain the new weak equivalences.
So the class of new weak equivalences is the smallest one which can create a model
structure along with the old cofibrations. This says that the new structure is the
left Bousfield localization \cite{Hirschhorn} \cite{Barwick}.
Alternatively, the argument given at the end of the proof of Theorem \ref{transfer} below 
gives explicitly the left Bousfield property. 

The characterization of fibrant objects and of fibrations to fibrant objects
in the last paragraph, is given by Proposition \ref{furtherinfo}.   
\end{proof}

The left Bousfield localization depends only on the class $\Rr$ and not on the choice
of subset $K$. 

\begin{proposition}
\label{lblindep}
Suppose $(\Rr ,K)$ and $(\Rr , K')$ are two direct localizing systems for the same
class of objects $\Rr$ in a 
left proper tractable model category $\mM$. Then the two new model structures given by the previous
theorem are the same.
\end{proposition}
\begin{proof}
By Corollary \ref{nweindep} the classes of new weak equivalences are the same,
and by definition the classes of cofibrations are the same, so the classes of fibrations are
the same. 
\end{proof}

\section{Transfer along a left Quillen functor}
\label{sec-directtransfer}

Suppose $F:\mM \rBotharrow \mN : E$ is a Quillen pair of adjoint functors between model categories
$\mM$ and $\mN$, with $F$ the left adjoint and $E$ the right adjoint. Suppose that $\mM$ and $\mN$ are
both tractable and left proper. Let $I$ and $J$ be cofibrant generating sets for $\mM$,
and $I'$ and $J'$ cofibrant generating sets for $\mN$.  

Suppose that $\Rr \subset \mM$ is a full subcategory,
and $K\subset \Arr (\mM )$ is a small subset of arrows. Let $K':= F(K)\cup J'\subset \Arr (\mM )$ be the set consisting of $J'$ plus all the $F(f)$ for
$f\in K$. Let $\Rr ':= RE^{-1}(\Rr )\subset \mN$ be the full subcategory consisting of those objects $Y\in \mN$ such that,
for a fibrant replacement $Y\rightarrt Y_1$ we have $E(Y_1)\in \Rr$. By condition (A4) for $\Rr$ and the fact that $E$ preserves equivalences between
fibrant objets, membership of $Y$ in $\Rr$ is independent of the choice of fibrant replacement $Y_1$. 

\begin{theorem}
\label{transfer}
In the above situation, suppose that $(\Rr , K)$ is a direct localizing system in
$\mM$. Define $(\Rr ', K')$ as above. Then $(\Rr ', K')$ is a direct localizing system in $\mN$. The functor $F$ is a left Quillen functor from the left Bousfield localization
of $\mM$ along $K$, to the left Bousfield localization of $\mN$ along $K'$.
\end{theorem}
\begin{proof}
We verify the conditions (A1)--(A6). 
\newline
(A1)---$K'=F(K)\cup J$ is clearly a small subset.
\newline
(A2)---$J'\subset K'$ by definition.
\newline
(A3)---since $F$ is a left Quillen functor, it preserves cofibrations. But by hypothesis $K\subset \cof (I)$,
so $F(K)\subset \cof (I')$. Also $J'\subset \cof (I')$ so $K'\subset \cof (I')$. 
The domains of arrows in $F(K)$ are cofibrant again because $F$ preserves cofibrant objects. 
\newline
(A4)---if $X\in \Rr '$ and $X\cong Y$ in $\Ho (\mM )$, choose a cofibrant (resp. fibrant) replacement
by a weak equivalence  $X\leftarr X_1$ (resp. 
$Y\rightarrt Y_1$) in $\mN$. The isomorphism in the homotopy category lifts to a diagram
of the form $X_1\leftarr Z \rightarrt Y_1$ where both maps are weak equivalences, and the first is a fibration.
Hence $Z$ is fibrant, and $E(X_1)\leftarr E(Z)\rightarrt E(Y_1)$ is a diagram of weak equivalences in $\mM$.
The condition $X\in \Rr'$ means $E(X_1)\in \Rr$. Condition (A4) for $\Rr$ says that $E(Z), E(Y_1)\in \Rr$ hence $Y\in \Rr'$.
\newline
(A5)---Suppose $Y\in \inj (K')$. In particular $Y\in \inj (J')$ so $Y$ is fibrant. 
By the adjunction between $F$ and $E$, the lifting property for $Y$ with respect to $F(K)$ is equivalent
to the lifting property for $E(Y)$ with respect to $K$. Hence $E(Y)\in \inj (K)$ so $E(Y)\in \Rr$. 
Here $Y$ is its own fibrant replacement so this shows that $Y\in \Rr '$.

We  now get to the important part of the argument, which is the verification of 
(A6). Suppose  $X\in \Rr '$ and $X$ is fibrant (i.e. $J'$-injective), and that 
$X\rightarrt Y$ is a pushout by an element of $K'$.
If it is a pushout by an element of $J'$ then it is already a weak equivalence so we can take $Y=Z$
(the identity is by definition in $\cell (K')$). 

Suppose it is a pushout by an element 
$F(g):F(A)\rightarrt F(B)$ where $g\in K$. That is to say we have a map $i':F(A)\rightarrt X$ and
$Y=X\cup ^{F(A)}F(B)$. By adjunction this gives a map $i:A\rightarrt E(X)$. Since $X$ is its own fibrant replacement,
by definition of $\Rr '$ we have $E(X)\in \Rr$. Also $E$ is a right Quillen functor so it preserves
fibrant objects: $E(X)$ is fibrant.  Apply (A6) to the pushout $E(X)\rightarrt E(X)\cup ^A B$ in $\mM$.
That gives a map $h:E(X)\cup ^AB \rightarrt Z$ in $\cell (K)$ such that $E(X)\rightarrt Z$ is a weak equivalence.
Note that both $g$ and $h$ are cofibrations, so $F(g)$ and $F(h)$ are cofibrations, and 
$F(h): F(E(X)\cup ^AB)\rightarrt F(Z)$ is in $\cell (K)$. The pushout $E(X)\rightarrt E(X)\cup ^AB$ is a cofibration
so the composed  map $E(X)\rightarrt Z$ is a cofibration, hence a trivial cofibration. These are preserved by $F$ so
$F(E(X))\rightarrt F(Z)$ is a trivial cofibration. 

The adjunction map $F(E(X))\rightarrt X$ induces a map
$F(E(X))\cup ^AB)\rightarrt X\cup ^{F(A)}F(B)$ and the pushout of $F(h)$ along here gives a map 
$$
Y= X\cup ^{F(A)}F(B) \rightarrt ( X\cup ^{F(A)}F(B))\cup ^{F(E(X)\cup ^AB)} F(Z) =: Z'
$$
in $\cell (K)$. 
On the other hand,
$$
Z' = ( X\cup ^{F(A)}F(B))\cup ^{F(E(X)\cup ^AB)} F(Z) = X \cup ^{F(E(X))}F(Z),
$$
but that is a  pushout along the trivial cofibration $F(E(X))\rightarrt F(Z)$, so the map $X\rightarrt Z'$ is a trivial cofibration. 
This completes the verification of (A6).

To show that $F$ is a left Quillen functor, note that the cofibrations are the same in
both cases so we just have to show that $F$ preserves new trivial cofibrations.
Suppose $X\rightarrt ^f Y$ is a new trivial cofibration. It fits into a 
diagram as in Corollary \ref{ntcfcriterion}. Applying the functor $F$ 
takes $\cell (K)$ to $\cell (K')$ so $F(f)$, which is a cofibration by the original
left Quillen property of $F$, again fits into a diagram of the same form.
Thus by Corollary \ref{ntcfcriterion}, $F(f)$ is a new trivial cofibration. 
\end{proof}

Applying the last paragraph of Theorem \ref{directLBL}, which in turn comes from 
Proposition \ref{furtherinfo}, gives the following remark. It is the abstract
version appropriate to the present stage of our construction,
of Bergner's result characterizing fibrant Segal categories \cite{BergnerSegal}. 

\begin{remark}
\label{transferfurtherinfo}
In the situation of Theorem \ref{transfer}, an object $U\in \mN$ is fibrant for the new model structure corresponding to 
$(\Rr ', K')$, if and only if it is fibrant in the original structure of $\mN$
and $E(U)$ is $K$-injective. If $U$ is a new fibrant object of $\mN$ then
a morphism $W\rightarrt^p U$ is a new fibration if and only if it is an old fibration
and $E(p)$ is in $\inj (K)$. 
\end{remark}

Suppose now that we have a set $Q$ and a collection of tractable left proper model categories $\mM _q$ for $q\in Q$,
with a collection of Quillen functors $F_q:\mM _q\rBotharrow \mN : E_q$ to a fixed tractable left proper $\mN$.
Suppose $(\Rr _q, K_q)$ are direct localizing systems for the $\mM _q$. Let $(I',J')$ be generators for $\mN$,
and define $(\Rr ', K')$ as follows. First, $\Rr '\subset \mN$ is the full subcategory of objects $Y\in \mN$ such that
for a fibrant replacement $Y\rightarrt Y'$, we have $EE=_q(Y')\in \Rr _q$ for all $q\in Q$. Then let 
$$
K':= J' \cup \bigcup _{q\in Q}F_q(K_q).
$$

\begin{theorem}
\label{transferfamily}
In this situation, $(\Rr ', K')$ is a direct localizing system for $\mN$.
An object $U\in \mN$ (resp. a morphism $W\rightarrt ^pU$ to a fibrant object)
is fibrant for the resulting new model structure (resp. a new fibration), if and only
if it is fibrant in the old structure and if each $E_q(U)$ (resp. $E_q(p)$) is
fibrant in the direct left Bousfield localization corresponding to $(\Rr _q, K_q)$. 
\end{theorem}
\begin{proof}
The same as above. For (A6) note that we have to treat pushout by a single map in $K'$, which is either in $J'$ or
of the form $F_q(g)$ for $g\in K_q$. Use the argument of the previous proof for the functors $F_q$ and $E_q$. 

The characterization of fibrant objects and fibrations to fibrant objects, comes
from the corresponding part of Theorem \ref{directLBL} going back to
Proposition \ref{furtherinfo}. 
\end{proof}



\part{Generators and relations}


\chapter{Precategories}
\label{precat1}

This chapter introduces the main object of study, the notion of {\em $\mM$-precategory}. The terminology ``precategory'' has
been used in several different ways, notably by Janelidze \cite{Janelidze}. The idea of the word is to invoke a structure coming prior to the full structure of a
category. For us, a categorical structure means 
a category weakly enriched over $\mM$ following Segal's method. Then a ``precategory'' will be
a kind of simplicial object without imposing the Segal conditions. The passage from a precategory to a weakly enriched category consists
of enforcing the Segal conditions using the small object argument. The basic philosophy behind this construction is that the precategory
contains the necessary information for defining the category, by generators and relations. The small object argument then corresponds
to the calculus whereby the generators and relations determine a category, this operation being generically denoted $\Seg$. 
Splitting up things this way is motivated by the fact that
simplicial objects satisfying the Segal condition are not in any obvious way closed under colimits. When we take colimits we get to arbitrary
simplicial objects or precategories, which then have to generate a category by the $\Seg$ operation. 

The calculus of generators and relations is the main subject of several upcoming chapters. In the present chapter, intended as a reference, 
we introduce the definition
of precategory appropriate to our situation, and indicate the construction of some important examples which will be used later. 
For the purposes of the present chapter we don't need to be too specific about the hypotheses on $\mM$; it will generally be supposed to be a 
tractable left proper cartesian model category, but this will be discussed in detail in Chapter \ref{weakenr1}.

\section{Enriched precategories with a fixed set of objects}
\label{sec-fixedprecats}

Suppose $X$ is a set. Following Lurie \cite{LurieGC}, define the category $\Delta _X$ to have objects which are sequences of $x_i\in X$
denoted by $(x_0,\ldots , x_n)$ for any $n\in \Delta$,
and morphisms 
$$
(x_0,\ldots , x_n)\rightarrt ^{\phi} (y_0,\ldots , y_m)
$$ 
for any $\phi : n\rightarrt m$ in $\Delta$, whenever
$y_{\varphi (i)}=x_i$ for $i=0,\ldots , n$. 
Write $\Delta ^o_X$ for the opposite category.

If $f:X\rightarrt Y$ is a map of sets, we obtain a functor $\Delta _f:\Delta _X\rightarrt \Delta _Y$ defined by
$\Delta _f (x_0,\ldots , x_n):= (f(x_0),\ldots , f(x_n))$ and the corresponding functor on  opposite categories denoted
$\Delta ^o_f: \Delta ^o_X\rightarrt \Delta ^o_Y$.

The basic objects of study will be functors $\pF :\Delta ^o_X\rightarrt \mM$.
Such a functor specifies for each sequence of elements $x_0,\ldots , x_n\in X$, an object $\pF (x_0,\ldots , x_n)$
and for each arrow $\phi : n\rightarrt m$ and sequence $y_0,\ldots , y_m$ a morphism 
$$
\pF (y_0,\ldots , y_m)\rightarrt \pF (y_{\phi (0)}, \ldots , y_{\phi (n)})
$$
compatible with compositions and identities on the level of $\phi$. Recall that the category of such functors is denoted 
$\diag (\Delta ^o _X, \mM )$. 

For reasons which will become clear with the counterexample of Section \ref{sec-unitalnecessary} in the later chapter on products, we want to impose a {\em unitality condition},
saying that $\pF (y_0)=\ast$ for single-object sequences. 
The unitality condition corresponds, in a certain sense, to requiring strict units even though the composition of morphisms
is not yet specified. 

\begin{definition}
An {\em $\mM$-precategory over $X$} is a functor $\pF :\Delta ^o_X\rightarrt \mM$
such that $\pF (x_0)\cong  \ast$ is the coinitial object of $\mM$, for any $x_0\in X$.
Let $\precat (X; \mM )$ denote 
the category of $\mM$-precategories over $X$ with morphisms which are natural transformations
of diagrams. 
\end{definition}

We occasionally use various terminologies such as {\em weakly $\mM$-enriched precategory
with object set $X$},
or some subset of those words, for elements of $\precat (X; \mM )$.

In Section \ref{sec-unitality} of the next chapter we shall consider more generally a category of {\em diagrams with unitality condition}
denoted $\diag (\Phi /\Phi _0, \mM )$ whenever $\Phi$ is a small category and $\Phi _0$ a subset of objects or equivalently a full subcategory. 
In the present case $\Phi = \Delta _X^o$ and the
subset $\Phi _0$ consists of all sequences of length zero $(x_0)$. Note that $\Phi _0$ is isomorphic to the discrete category 
corresponding to the set $X$, so in the notation of Section \ref{sec-unitality}, 
$$
\precat (X; \mM ) := \diag (\Delta ^o_X /X, \mM ).
$$

It would of course be interesting to investigate further what would happen if we considered all objects in 
$\diag (\Delta ^o _X, \mM )$, with a weak unitality condition imposed by the Segal condition for $n=0$,
saying that $\pF (x_0)$ should be contractible. For the present purposes this doesn't lead directly to a cartesian model category,
indeed the degeneracies of the single points $\pF (x_0)=\ast$ play an important role in assuring the compatibility of weak
equivalences with products. It is likely that this problem could be worked around, and that the resulting theory would be closely
related to Kock's weak unit condition \cite{Kock} as well as to Paoli's special $Cat^n$-groups. 

For the present exposition it will be convenient to proceed as directly as possible without considering these other possibilities. 
When necessary, we refer to the objects of $\diag (\Delta ^o_X , \mM )$ as ``non-unital precategories''. 

\section{The Segal conditions}
\label{sec-segalcond}

Given an $\mM$-precategory $\pA $ over a set $X$, we can look at the {\em Segal maps}. If $(x_0,\ldots , x_n)$ is a sequence in $X$,
the principal edges of the $n$-simplex are maps in $\Delta _X$
$$
(x_i,x_{i+1})\rightarrt (x_0,\ldots , x_n) 
$$
and these give maps 
$$
\pA (x_0,\ldots , x_n)\rightarrt \pA (x_i,x_{i+1}).
$$
Put these together to get the {\em Segal map at $(x_0,\ldots , x_n)$}
$$
\pA (x_0,\ldots , x_n)\rightarrt \pA (x_0,x_{1})\times \cdots \times \pA (x_{n-1},x_{n}) .
$$

Note that for $n=0$ the Segal map at $(x_0)$ is 
$$
\pA (X_0) \rightarrt \ast .
$$
Thus the unitality condition says that the $n=0$ Segal maps are isomorphisms. 

Say that an $\mM$-precategory {\em satisfies the Segal condition} if the Segal maps are weak equivalences,
in other words they are contained in the subcategory of weak equivalences $\Ww \subset \mM$. 
A precategory satisfying this condition will be called an {\em weakly $\mM$-enriched category}
or just {\em weak $\mM$-category}. 
An $\mM$-precategory is said to be a {\em strict $\mM$-category} if the Segal maps are isomorphisms.

As was amply pointed out in Part I, the Segal condition at $n=2$ serves to define a weak composition
operation in the following sense. For any three objects $x_0,x_1,x_2\in X$
we  have a diagram 
$$
\begin{diagram}
\pA (x_0,x_1,x_2)& \rightarr & \pA (x_0,x_1)\times \pA (x_1,x_2)\\
\downarr & & \\
\pA (x_0,x_2)
\end{diagram}
$$
where the horizontal arrow is the Segal map; if it is a weak equivalence then
the vertical arrow projects to a map 
$$
\pA (x_0,x_1)\times \pA (x_1,x_2)\rightarrt \pA (x_0,x_2)
$$
in $\Ho (\mM )$. The Segal conditions for higher $n$ serve to fix the higher homotopy coherences necessary starting with associativity at $n=3$. It is necessary to include
all of the higher coherences in order to obtain a theory compatible with products,
see Chapter \ref{product1} below.

\section{Varying the set of objects}
\label{sec-varying1}

Up until now we have discussed the category of $\mM$-precategories on a fixed set of objects. This way of splitting off a first part of the argument
was originally suggested by Barwick \cite{BarwickThesis}, and that idea will be
continued in Section \ref{sec-pcxm} where we consider  model structures on
$\precat (X,\mM )$. 

The next step is to investigate what happens under maps between sets of objects,
and to define a category $\precat (\mM )$ of $\mM$-precategories without specified set of objects. 

If $f: X\rightarrt Y$ is a map of sets, we can {\em pull back} a structure of $\mM$-precategory on $Y$, to a structure of $\mM$-precategory 
on $X$, just by composing a diagram $\pA :\Delta ^o_Y\rightarrt \mM$ with the functor $\Delta ^o_f$: 
$$
(\Delta ^o_f)^{\ast}: \diag ( \Delta ^o_Y, \mM )\rightarrt \diag ( \Delta ^o_X, \mM ),
$$
$$
(\Delta ^o_f)^{\ast}(\pA )(x_0,\ldots , x_n)=\pA (f(x_0),\ldots , f(x_n)).
$$
This clearly preserves the unitality condition, so it restricts to a functor which, if no confusion
arises, will be denoted just by 
$$
f^{\ast}: \precat (Y; \mM )\rightarrt \precat (X; \mM ).
$$ 
We also get a left adjoint to $f^{\ast}$ denoted $f_! : \precat (X; \mM )\rightarrt \precat (Y; \mM )$. If 
$f: X\hookrightarrow Y$ is an inclusion, then we can write $Y= f(X)\sqcup Z$ where $Z$ is the complement of the image.
In this case, $f_! (\pA )(y_0,\ldots , y_n) = \pA (x_0,\ldots , x_n)$ if all $y_i\in f(X)$ and $x_i$ is the preimage of $y_i$.
If $y_0= \ldots = y_n = z\in Z$ then $f_! (\pA )(y_0,\ldots , y_n)=\ast$, and $f_! (\pA )(y_0,\ldots , y_n)=\emptyset$ in all
other cases. Thus, $f_!(\pA )$ is the precategory $\pA $, extended by adding on the discrete set $Z$ considered as a discrete 
$\mM$-enriched category.

\begin{lemma}
\label{pfpb}
Suppose $f:X\rightarrt Y$ is a map of sets. If $\pA $ is an $\mM$-precategory on object set $Y$, which satisfies
the Segal condition, then $f^{\ast}(\pA )$ satisfies the Segal condition as an $\mM$-precategory on $X$.

Assume at least condition (DCL) which is part of the cartesian condition \ref{def-cartesian} on $\mM$. 
If $f$ is injective and if $\pB $ is an $\mM$-precategory on $X$ satisfying the Segal condition, then
$f_!(\pB )$ satisfies the Segal condition as an $\mM$-precategory on $Y$.
\end{lemma}
\begin{proof}
The Segal maps for $f^{\ast}(\pA )$ are some among the Segal maps for $\pA $, so the Segal conditions for $\pA $ imply the same for $f^{\ast}(\pA )$. Suppose $f$ is injective and $\pB \in \precat (X,\mM )$ satisfies the Segal conditions. Use the notations 
$Y= f(X)\sqcup Z$ of the paragraph before the lemma. Given a sequence $(y_0,\ldots , y_n)$
of objects $y_i\in Y$, if $y_i=f(x_i)$ for all $i$ then 
the Segal map for $f_!(\pB )$ at $(y_0,\ldots , y_n)$ is the same as that for $\pA $
at $(x_0,\ldots , x_n)$. If $y_0= \ldots = y_n = z\in Z$ then the Segal map for $f_!(\pB )$ at
$(y_0,\ldots , y_n)$ is the identity of $\ast$. If $(y_0,\ldots , y_n)$ is a sequence
which is not constant and which contains at least one element of $Z$, then one of the
adjacent pairs $(y_{i-1},y_{i})$ has to be nonconstant and contain an element of $Z$.
Using condition (DCL) via Lemma \ref{emptyempty}, the direct product of anything
with $\emptyset$ is again $\emptyset$, so in these cases the Segal maps are the identity
of $\emptyset$. In all cases, the Segal maps remain weak equivalences. 
\end{proof}

\begin{lemma}
\label{pbwe}
Suppose $f:X\rightarrt Y$ is a map of sets. 
The pullback $f^{\ast}$ preserves levelwise weak equivalences, that is it takes levelwise weak equivalences in $\precat (Y,\mM )$
to levelwise weak equivalences in $\precat (X,\mM )$. Similarly, $f^{\ast}$ preserves
levelwise cofibrations and trivial cofibrations. 
\end{lemma}
\begin{proof}
The pullback $f^{\ast}$ of diagrams will preserve any levelwise properties. 
\end{proof}

\section{The category of precategories}
\label{sec-varying2}

We would  like to define a notion of $\mM$-precategory on an unspecified or variable set of objects, just as for 
usual categories. 
An {\em $\mM$-precategory} $\pA $ consists of a set denoted $\Ob (\pA )$, and an $\mM$-precategory over $\Ob (\pA )$ denoted 
$$
(x_0,\ldots , x_n)\mapsto \pA (x_0,\ldots , x_n)\in \mM .
$$

A {\em morphism} $f:\pA \rightarrt \pB $ between two $\mM$-precategories consists of a map of sets $\Ob (f):\Ob (\pA )\rightarrt \Ob (\pB )$,
and for any $x_0,\ldots , x_n\in \Ob (\pA )$ a morphism 
$$
f(x_0,\ldots , x_n): \pA (x_0,\ldots , x_n)\rightarrt \pB (\Ob (f)(x_0),\ldots , \Ob (f)(x_n));
$$
these are required to satisfy naturality with respect to maps $\phi : n\rightarrt m$ in $\Delta$. 
Henceforth, if no confusion is possible, we denote by $f(x_i):= \Ob(f)(x_i)$ the action on objects, and use $f$ also to denote $f(x_0,\ldots , x_n)$. 

In terms of the notation given in Section \ref{sec-varying1},
we can think of $f$ as a morphism of $\mM$-precategories over object set $\Ob (\pA )$, denoted 
$$
\Mor _f : \pA \rightarrt \Ob (f) ^{\ast}(\pB ). 
$$

Composition of morphisms is defined in the obvious way. 

Let $\precat (\mM )$ denote the category of $\mM$-precategories thusly defined. Taking the ``underlying set of objects'' is a functor 
$$
\Ob : \precat (\mM )\rightarrt \Sets .
$$
This is a fibered category as we now explain. Given a map of sets $g:X\rightarrt Y$ this induces a functor $\Delta _g: \Delta _X\rightarrt \Delta _Y$
hence an adjoint pair
$$
\Delta _{g,!}: \diag (\Delta _X^o; \mM )\rBotharrow \diag (\Delta _Y^o; \mM ): \Delta _g^{\ast}.
$$
We have $\Delta _g^{\ast}(\pF )(x_0,\ldots , x_n) = \pF (g(x_0),\ldots , g(x_n))$. The adjoint pair on the categories of $\mM$-precategories is
$$
\precat _{g,!}: \precat (X; \mM )\rBotharrow \precat ( Y; \mM ): \Delta _g^{\ast}.
$$
with $\precat _{g,!} = U_{Y,!}\circ \Delta _{g,!,u}$. 

Consider the fibered category $\Ff \rightarrt \Sets$ whose fiber over $X$ is $\precat( X; \mM )$ and whose pullback maps are the
$\Delta _g^{\ast}$. An object of $\Ff$ is by definition a pair $(X,\pA )$ where $X\in \Sets$ and $\pA \in \precat(X; \mM )$.
A morphism from $(X,\pA )$ to $(Y,\pB )$ is a pair $(g,h)$ where $g:X\rightarrt Y$ is a morphism in $\Sets$ and $h: \pA \rightarrt g^{\ast}(\pB )= \Delta _g^{\ast}$
is a morphism in $\precat(Y; \mM )$. By inspection, this is the same structure as defined above, that is $\Ff = \precat (\mM )$,
with functor $\Ob$ being the projection to $\Sets$. A morphism $f= (\Ob (f), f)$ is cartesian if and only if $\Mor _f$ induces
isomorphisms
$$
f(x_0,\ldots , x_n): \pA (x_0,\ldots , x_n)\rightarrt ^{\cong} \pB (\Ob (f)(x_0),\ldots , \Ob (f)(x_n)).
$$

\section{Basic examples}
\label{sec-precatexamples}

In this section we indicate many of the basic examples of precategories which will be important at various places in
the subsequent chapters. Rather than give a full definition and discussion here, we refer to the appropriate places when necessary.

Suppose $X$ is a set. The {\em discrete precategory} $\disc (X)$ is defined to have $\Ob (\disc (X)):=X$ and
$$
\disc (X) (x_0,\ldots , x_n) := \left\{  
\begin{array}{ll}
\ast & \mbox{if } x_0= \cdots = x_n \\
\emptyset & \mbox{otherwise} 
\end{array}
\right.
$$
Here $\emptyset$ and $\ast$ denote the initial and coinitial objects of $\mM$ respectively. The functor $X\mapsto \disc (X)$
provides a functor $\Sets \rightarrt \precat (\mM )$, fully faithful and left adjoint to $\Ob : \precat (\mM )\rightarrt \Sets$.

The {\em codiscrete precategory} $\codisc (X)$ is defined to have $\Ob (\codisc (X)):= X$ and 
$$
\codisc (X) (x_0,\ldots , x_n) := \ast 
$$
for all sequences $x_{\cdot}$. The functor $\codisc :\Sets \rightarrt \precat (\mM )$ is fully faithful and right adjoint to $\Ob$. 

Suppose $k\in \nn$. Denote by 
$$
[k]:= \{ \upsilon _0,\ldots , \upsilon _k\},\;\;\; \upsilon _0 < \upsilon _1 < \cdots < \upsilon _k
$$
the standard linearly ordered set with $k+1$ elements. The $[k]$ are the objects of $\Delta$. It will be important to consider various different
precategories whose underlying set of objects is $[k]$. These examples will all be {\em ordered precategories} in the sense that 
if $(y_0,\ldots , y_p)$ is any sequence of elements of the set $[k]$ with 
$y_j = \upsilon _{i_j}$, and if the sequence $y_{\cdot}$
is strictly decreasing at any place (i.e. if it is not an increasing sequence),
then
$$
\pA (y_0,\ldots , y_p) = \emptyset.
$$ 

For $B\in \mM$ we have the precategory $h([k],B)$ which has the following concrete description
(see Proposition \ref{hdescription}):
suppose $(y_0,\ldots , y_p)$ is any sequence of elements of the set $[k]$ with 
$y_j = \upsilon _{i_j}$. Then:
\newline
---if $(y_0,\ldots , y_p)$ is increasing but not constant i.e. $i_{j-1}\leq i_j$ but 
$i_0 < i_p$ then
$$
h([k]; B)(y_0,\ldots , y_p) = B; 
$$
---if $(y_0,\ldots , y_p)$ is  constant i.e. $i_0=i_1=\ldots = i_p$ then
$$
h([k]; B)(y_0,\ldots , y_p) = \ast; 
$$
and otherwise, that is if there exists $1\leq j \leq p$ such that $i_{j-1}>i_j$ then
$$
h([k]; B)(y_0,\ldots , y_p) = \emptyset.
$$ 
For any precategory $\pA \in \precat (\mM )$, a map $h([k];B)\rightarrt \pA $ is the same thing as a collection of elements
$x_0,\ldots , x_k\in \Ob (\pA )$ together with a map $B\rightarrt \pA (x_0,\ldots , x_k)$ in $\mM$, so we can think of 
$h([k],B)$ as being a ``representable'' object in a certain sense. 

It has as ``boundary'' the precategory
$h(\partial [k],B)$, defined using the skeleton construction in Chapter \ref{cofib1}, 
with the following description (see Lemma \ref{hdelbardescription}):
\newline
---if $(y_0,\ldots , y_p)$ is increasing but not constant i.e. 
$i_{j-1}\leq i_j$ but $i_0 < i_p$, and
if there is any $0\leq m\leq k$ such that $i_j\neq m$ for all $0\leq j\leq k$, then
$$
h(\partial [k]; B)(y_0,\ldots , y_p) = B; 
$$
---if $(y_0,\ldots , y_p)$ is  constant i.e. $i_0=i_1=\ldots = i_p$ then
$$
h(\partial [k]; B)(y_0,\ldots , y_p) = \ast; 
$$
and otherwise, that is if either there exists $1\leq j \leq p$ such that $i_{j-1}>i_j$
or else if the map $j\mapsto y_j$ is a surjection from $\{ 0,\ldots , p\}$ to $[k]$,  then
$$
h(\partial [k]; B)(y_0,\ldots , y_p) = \emptyset.
$$

More generally, the pushouts
$$
h([k],\partial [k]; A \rightarrt ^f B):= 
h([k]; A)\cup ^{h(\partial [k];A)} \partial h(\partial [k];B)
$$
are also useful, being the generators of the Reedy cofibrations in $\precat (\mM )$.

In Section \ref{sec-upsilon} and later in Chapter \ref{freecat1} we consider precategories
$\Upsilon (B_1,\ldots , B_k)$ with the same set of objects $[k]$, depending on $B_1,\ldots , B_k\in \mM$. 
The basic idea is to put $B_i$ in as space of morphisms from $\upsilon _{i-1}$ to $\upsilon _i$. Thus,
the main part of the structure of precategory is given by 
$$
\Upsilon (B_1,\ldots , B_k)(\upsilon _{i-1}, \upsilon _i ) := B_i.
$$
This is extended whenever there is a constant string of points on either side:
$$
\Upsilon (B_1,\ldots , B_k)(\upsilon _{i-1},\ldots , \upsilon _{i-1}, \upsilon _i ,\ldots , \upsilon _i) := B_i.
$$
The unitality condition on the diagram $\Delta _{\{ \upsilon _0,\ldots , \upsilon _k\}} \rightarrt \mM$ means that for $0\leq i \leq k$ we have
$$
\Upsilon (B_1,\ldots , B_k)(\upsilon _{i}, \ldots ,  \upsilon _i ) := \ast .
$$
In all other cases, 
$$
\Upsilon (B_1,\ldots , B_k)(x_0, \ldots ,  x_n ) := \emptyset .
$$
If $\pA \in \precat (\mM )$, a map $\Upsilon (B_1,\ldots , B_k)\rightarrt \pA $ is the same thing as a collection of objects $x_0,\ldots , x_k\in \Ob (\pA )$
together with maps $B_i\rightarrt \pA (x_{i-1},x_i)$ in $\mM$. 
 
The categorification of $\Upsilon (B_1,\ldots , B_k)$
can be described explicitly; it is denoted by
$\Upstild (B_1,\ldots , B_k)$ in Section \ref{basicseqfree}.
For any sequence $\upsilon _{i_0},\ldots , \upsilon _{i_n}$ with
$i_0\leq \ldots \leq i_n$, we put 
\begin{equation}
\label{upstildform}
\Upstild _k(B_1,\ldots , B_k)(\upsilon _{i_0},\ldots , \upsilon _{i_n}):= 
B_{i_0+1}\times B_{i_0+2 }\times \cdots \times B_{i_n -1}\times B_{i_n}.
\end{equation}
For any other sequence, that is to say any sequence which is not increasing, the value is $\emptyset$. The value on a constant sequence is $\ast$. 

As a rough approximation, the calculus of generators and relations can be understood as being the small object argument applied using the
inclusions $\Upsilon (B_1,\ldots , B_k)\hookrightarrow \Upstild (B_1,\ldots , B_k)$. Starting with a precategory $\pA $, for every
map $\Upsilon (B_1,\ldots , B_k)\rightarrt \pA $, take the pushout with $\Upstild (B_1,\ldots , B_k)$. Keep doing this and eventually
we get to a precategory which satisfies the Segal conditions.

\section{Limits, colimits and local presentability}
\label{sec-limits}

It will be useful to have explicit descriptions of limits and colimits in $\precat (\mM )$. We can then show local 
presentability. 

Suppose $\{ \pA _i\} _{i\in \alpha}$ is a diagram of objects $\pA _i\in \precat (\mM )$, that is a functor $\alpha \rightarrt \precat (\mM )$. 
Let $X_i:= \Ob (\pA _i)$ denote the object sets so we can consider $\pA _i\in \precat (X_i; \mM )$. 
For any $f:i\rightarrt j$ in $\alpha$
denote by $\phi _{f}: X_i\rightarrt X_j$ the transition map on object sets, then 
$\rho _f: \pA _i\rightarrt \phi _f^{\ast} \pA _j$ the transition maps on the level of precategories.

Start by constructing the limit in $\precat (\mM )$. The object set of the limit will be
$$
X := \mylim _{i\in \alpha} X_i.
$$
We have maps $p_i:X\rightarrt X_i$, hence $p_i^{\ast}(\pA _i)\in \precat (X; \mM )$. These are provided with transition
maps, indeed for $f:i\rightarrt j$ in $\alpha$, $\phi _fp_i=p_j$ so $p_i^{\ast}(\phi _f^{\ast}\pA _j)=p_j^{\ast}\pA _j$
and $p_i^{\ast}(\rho _f): p_i^{\ast}(\pA _i)\rightarrt p_j^{\ast}\pA _j$ provide the transition maps for
the system of $p_i^{\ast}(\pA _i)$ considered as a diagram $\alpha \rightarrt \precat (X; \mM )$.  Put
$$
\pA := \mylim ^{\precat (X; \mM )}_{i\in \alpha} p_i^{\ast}(\pA _i) \;\; \in \; \precat (X; \mM ).
$$
We claim that this is the limit of the diagram $\{ \pA _i\} _{i\in \alpha}$ in $\precat (\mM )$. 

Suppose $\pB \in \precat (\mM )$ with $\Ob (\pB )=Y$ and suppose given a system of maps $\pB \rightarrt \pA _i$.
These correspond to maps $r_i:Y\rightarrt X_i$ and $\varphi _i: \pB \rightarrt r_i^{\ast}(\pA _i)$.
The collection of $r_i$ gives a uniquely determined map $r:Y\rightarrt X$, and $r^{\ast}p_i^{\ast}(\pA _i)= r_i^{\ast}(\pA _i)$
so the $\varphi  _i$ correspond to a collection of maps $\pB \rightarrt r^{\ast}(p_i^{\ast}(\pA _i))$. This gives a uniquely determined map
$\varphi : \pB \rightarrt r^{\ast}(\pA )$ whose composition with the projections of the limit expression for $\pA $, are the
$\varphi _i$. The pair $(r,\varphi )$ is a map $\pB \rightarrt \pA $ in $\precat (\mM )$ uniquely solving the universal problem 
to show that $(X,\pA )$ is the limit of the $\pA _i$. 

The limit $\pA $ can be described explicitly as an element of $\precat (X;\mM )$, by the discussion preceding Lemma \ref{colimexpr}:
for any $x_0,\ldots , x_n\in X$,
$$
\pA (x_0,\ldots , x_n) = \mylim _{i\in \alpha} \pA _i (p_ix_0,\ldots , p_ix_n).
$$

Turn now to the construction of the colimit. The object set will be
$$
Z:= \colim _{i\in \alpha} X_i
$$
with maps $q_i:X_i\rightarrt Z$. These give $q_{i,\ast}(\pA _i)\in \precat (X; \mM )$, with transition maps defined as follows.
If $f:i\rightarrt j$ is an arrow in $\alpha$ then $q_j\phi _f=q_i$ so $q_{j,!}(\phi _{f,!}(\pA _i))= q_{i,!}(\pA _i)$. The adjunction
between $\phi _{f,!}$ and $\phi _f^{\ast}$ means that the transition map $\rho _f$ for the system of $\pA _i$ may be viewed as a map denoted
$\tilde{\rho}_f:\phi _{f,!}(\pA _i)\rightarrt \pA _j$, in particular we get 
$$
q_{j,!}(\tilde{\rho}_f) : q_{i,!}(\pA _i) = q_{j,!}(\phi _{f,!}(\pA _i))\rightarrt q_{j,!}(\pA _j).
$$
These provide transition maps for the diagram $\{ q_{i,!}(\pA _i)\}_{i\in \alpha}$ with values in $\precat (Z; \mM )$ and set
$$
\pC := \colim ^{\precat (Z; \mM )} _{i\in \alpha}q_{i,!}(\pA _i) \;\; \in \; \precat (Z; \mM ).
$$
The object $(Z,\pC )$ is the colimit of the diagram of $\pA _i$, for the same formal reason as in the previous discussion of the limit. 

The colimit to define $\pC $ is taken in the category $\precat (Z; \mM )$ which presupposes in general applying the operation $U_!$
if $\alpha$ is disconnected. However, when the definition is unwound explicitly this phenomenon disappears, being absorbed 
in the calculation of the value of $\pC $ on a sequence of points $(z_0,\ldots , z_n)$ via the introduction of a new category
$\alpha /(z_0,\ldots , z_n)$. 

Suppose $(z_0,\ldots , z_n)$ is a sequence of elements of $Z$. Let $\alpha / (z_0,\ldots , z_n)$ denote the 
set of pairs $(i,(x_0,\ldots , x_n))$ where $i\in \alpha$ and $(x_0,\ldots , x_n)$ is a sequence of points in $X_i$ such that
$q_i(x_k)= z_k$ for $k=0,\ldots , n$. The association $(i,(x_0,\ldots , x_n))\mapsto \pA _i(x_0,\ldots , x_n)$ is a diagram from
$\alpha / (z_0,\ldots , z_n)$ to $\mM$.

\begin{lemma}
\label{colimcalc}
In the above situation, the value of the colimiting object $\pC $ on the sequence $(z_0,\ldots , z_n)$ is calculated as
a colimit of the diagram $\alpha / (z_0,\ldots , z_n)\rightarrt \mM$:
$$
\pC (z_0,\ldots , z_n) = \colim _{(i,(x_0,\ldots , x_n))\in \alpha / (z_0,\ldots , z_n)}\pA _i(x_0,\ldots , x_n) .
$$
\end{lemma}
\begin{proof}
Let $q_{i,!}^{wu}: \diag (\Delta ^o_{X_i}, \mM )\rightarrt \diag (\Delta ^o_{Z};\mM )$ denote the pushforward in the world of non-unital
precategories. It is the pushforward for diagrams valued in $\mM$ along the functor 
$$
\Delta ^o_{X_i}\rightarrt \Delta ^o_Z.
$$
If $(z_0,\ldots , z_n)$ is a sequence of points in $Z$, then $q_{i,!}^{wu}(\pA _i)(z_0,\ldots , z_n)$ is the colimit of the
$\pA _i(x_0,\ldots , x_k)$ over the category of pairs
$(u,(x_0,\ldots , x_k))$ where $u: (q_ix_0,\ldots , q_ix_k)\rightarrt (z_0,\ldots , z_n)$ is a map in $\Delta ^o_Z$
or equivalently $(z_0,\ldots , z_n)\rightarrt (q_ix_0,\ldots , q_ix_k)$ is a map in $\Delta _Z$. Such a map factors
through a unique map $(z_0,\ldots , z_n)\rightarrt (q_ix_{i_0},\ldots , q_ix_{i_n})$ where $(i_0,\ldots , i_n)$
is a multiindex representing a map $[n]\rightarrt \nocom [k]$ in $\Delta$. Hence the category of pairs in question is a disjoint union of
categories having coinitial objects. This yields the expression of $q_{i,!}^{wu}(\pA _i)(z_0,\ldots , z_n)$ as the disjoint sum
of $\pA _i(x_0,\ldots , x_n)$ over all sequences $(x_0,\ldots , x_n)$ of objects in $X_i$ such that $q_i(x_j)=z_j$. 
Putting these together over all $i\in \alpha$ gives 
$$
\colim _{i\in \alpha}(q_{i,!}^{wu}(\pA _i)(z_0,\ldots , z_n)) = \colim _{(i,(x_0,\ldots , x_n))\in \alpha / (z_0,\ldots , z_n)}\pA _i(x_0,\ldots , x_n).
$$
Now $U_!$ is a left adjoint so it preserves colimits, in particular
$$
\pC = \colim ^{\precat (Z; \mM )}_{i\in \alpha}q_{i,!}(\pA _i) = \colim _{i\in \alpha}U_!(q_{i,!}^{wu}(\pA _i)) 
$$
$$
= 
U_! \left( \colim _{i\in \alpha}q_{i,!}^{wu}(\pA _i)\right)  = U_!(\pC ')
$$
where  $\pC '\in \diag (\Delta ^o_{Z}; \mM )$ is the object defined by
$$
\pC '(z_0,\ldots , z_n):= \colim _{(i,(x_0,\ldots , x_n))\in \alpha / (z_0,\ldots , z_n)}\pA _i(x_0,\ldots , x_n).
$$
To finish the proof note that $\pC '$
is already in $\precat (Z; \mM )$. Indeed for a single element $z_0\in Z$ the category 
$\alpha /(z_0)$ of pairs $(i,x_0)$ with $q_i (x_0)= z_0$ is connected. This is factoid about colimits of sets. 
Since $\pA _i$ is unital 
we have $\pA _i(x_0)=\ast$, thus
$$
\pC '(z_0) = \colim _{(i,x_0)\in \alpha / (z_0)}\ast = \ast ,
$$
in other words $\pC '$ is unital. Therefore $\pC = U_!(\pC ')=\pC '$ which is the statement of the lemma. 
An alternative proof would be to note that $\pC '$ satisfies the required universal property for defining a colimit. 
\end{proof}

\begin{corollary}
\label{constobjcolim}
If $\alpha$ is a connected category and $\pA _{\cdot}:\alpha \rightarrt \precat (X; \mM )$ is a diagram of 
$\mM$-precategories with a common object set $X$, then the natural map 
$$
\colim ^{\precat ( \mM )}_{i\in \alpha} \pA _i \rightarrt \colim ^{\precat ( X;\mM )}_{i\in \alpha} \pA _i 
$$
is an isomorphism. 
\end{corollary}
\begin{proof}
The explicit descriptions of both sides are the same.
\end{proof}

It is worthwhile to discuss explicitly some special cases. For example $\Ob (\pA \times \pB ) = \Ob (\pA )\times \Ob (\pB )$
and for any sequence $((x_0,y_0),\ldots , (x_n,y_n))$ of elements of $\Ob (\pA )\times \Ob (\pB )$ we have
$$
\pA \times \pB  ((x_0,y_0),\ldots , (x_n,y_n)) = \pA (x_0,\ldots , x_n)\times \pB (y_0,\ldots , y_n).
$$
This extends to fiber products: if $\pA \rightarrt \pC $ and $\pB \rightarrt \pC $ are maps, then
$$
\Ob (\pA \times _{\pC }\pB ) = \Ob (\pA )\times _{\Ob ({\pC })}\Ob (\pB )
$$
and again for any sequence $((x_0,y_0),\ldots , (x_n,y_n))$ of elements in the fiber product of object sets we have
$$
\pA \times _{\pC } \pB  ((x_0,y_0),\ldots , (x_n,y_n)) = \pA (x_0,\ldots , x_n)\times _{{\pC }(z_0,\ldots , z_n)}\pB (y_0,\ldots , y_n)
$$
where the $z_i$ are the common images of $x_i$ and $y_i$ in $\Ob ({\pC })$. 

For colimits, the disjoint sum $\pA \sqcup  \pB $ has object set $\Ob (\pA )\sqcup \Ob (\pB )$, and
for a sequence of elements $(z_0,\ldots , z_n)$ in here we have
$$
(\pA \sqcup  \pB )(z_0,\ldots , z_n)= \left\{ \begin{array}{ll} 
\pA (z_0,\ldots , z_n) &\mbox{if all }z_i\in \Ob (\pA ) \\
\pB (z_0,\ldots , z_n) &\mbox{if all }z_i\in \Ob (\pB ) \\
\emptyset &\mbox{otherwise }
\end{array}
\right.  .
$$
Look at the coproduct or pushout of
two maps $u:{\pC }\rightarrt \pA $ and $v:{\pC }\rightarrt \pB $, supposing that one of the maps say $v$ is injective on
the set of objects. This will usually be the case in our applications because we usually look at pushouts along cofibrations.
The coproduct $\pA \cup ^{\pC } \pB $ has object set
$$
\Ob (\pA )\cup ^{\Ob ({\pC })} \Ob (\pB ) = \Ob (\pA )\sqcup (\Ob (\pB )-\Ob ({\pC })).
$$
The category $\alpha$ indexing the colimit has three objects denoted $a,b,c$ with arrows $a\leftarrow c \rightarrt b$.
Given a sequence $(z_0,\ldots , z_n)$ of elements in here, the category $\alpha / (z_0,\ldots , z_n)$ 
has objects of three kinds denoted $a$, $b$ and $c$, and objects of a given kind correspond to sequences in $\Ob (\pA )$,
$\Ob (\pB )$ or $\Ob ({\pC })$ respectively mapping to the given $z_{\cdot}$. Our general formula of Lemma \ref{colimcalc}
expresses $(\pA  \cup ^{\pC }\pB )(z_0,\ldots , z_n)$ as the pushout of disjoint sums
$$
\coprod _{x_{\cdot}\mapsto z_{\cdot}} \pA (x_0,\ldots , x_n) 
\cup ^{
\coprod _{w_{\cdot}\mapsto z_{\cdot}} {\pC }(w_0,\ldots , w_n)
}
\coprod _{y_{\cdot}\mapsto z_{\cdot}} \pB (y_0,\ldots , y_n).
$$

Another important case is 
that of filtered colimits. Suppose $\alpha$ is a filtered (resp. $\kappa$-filtered) category
and $\{ \pA _i\} _{i\in \alpha}$ is a diagram in $\precat ( \mM )$. Put $X_i:= \Ob (\pA _i)$ and
$Z:= \colim _{i\in \alpha}X_i$. Then for any sequence $(z_0,\ldots , z_n)$ of elements of $Z$,
the category $\alpha / (z_0,\ldots , z_n)$ is again filtered (resp. $\kappa$-filtered),
so 
$$
\pA (z_0,\ldots , z_n) = \colim _{(i,(x_0\ldots , x_n)\in \alpha / (z_0,\ldots , z_n)} \pA (x_0\ldots , x_n)
$$
is a filtered (resp. $\kappa$-filtered) colimit in $\mM$.

\begin{lemma}
\label{kappapres}
Suppose $\mM$ is locally presentable.
An object $\pA \in \precat (\mM )$ is $\kappa$-presentable if and only if $X:=\Ob (\pA )$ is a set of cardinality $<\kappa$,
and each $\pA (x_0,\ldots , x_p)$ is a $\kappa$-presentable object of $\mM$.
\end{lemma}
\begin{proof}
Suppose $X:=\Ob (\pA )$ is a set of cardinality $<\kappa$,
and each $\pA (x_0,\ldots , x_p)$ is a $\kappa$-presentable object of $\mM$. Suppose 
$\{ \pB _i\} _{i\in \beta}$ is a diagram in $\precat ( \mM )$ indexed by a $\kappa$-filtered category $\beta$,
and suppose given a map 
$$
\pA  \rightarrt \pB := \colim _{i\in \beta} \pB _i .
$$
In particular we get a map $\Ob (\pA )\rightarrt \Ob (\pB )=\colim _{i\in \beta} \Ob (\pB _i)$
and the condition $|\Ob (\pA )| < \kappa$ implies that this map factors through a map $\Ob (\pA )\rightarrt \Ob (\pB _j)$ for
some $j\in \beta$. Given a sequence of objects $(x_0,\ldots , x_n)\in \Ob (\pA )$, let $(y_0,\ldots , y_n)$
denote the corresponding sequence of objects in $\Ob (\pB _j)$ and $(z_0,\ldots , z_n)$ the sequence in $\Ob (\pB )$.
Let $j\backslash \alpha$ denote the category of objects under $j$ in $\alpha$. Given $j\rightarrt i$
the image of $(y_0,\ldots , y_n)$ is a sequence denoted $(y^i_0,\ldots , y^i_n)$ in $\Ob (\pB _i)$ mapping to $(z_0,\ldots , z_n)$. 
This gives a functor 
$$
j\backslash \alpha \rightarrt \alpha / (z_0,\ldots , z_n),
$$
and the $\kappa$-filtered property of $\alpha$ implies that this functor is cofinal. 
Hence by Lemma \ref{colimcalc} and the invariance of colimits under cofinal functors, 
$$
\pB (z_0,\ldots , z_n) = \colim _{(j\rightarrow i)\in j\backslash \alpha} \pB _i(y^i_0,\ldots , y^i_n).
$$
The category $j\backslash \alpha $ is $\kappa$-filtered, so the map 
$$
\pA (x_0,\ldots , x_n)\rightarrt \pB (z_0,\ldots , z_n)
$$
factors through one of the $\pB _i(y^i_0,\ldots , y^i_n)$ for $j\rightarrt i$ in $\alpha$.
The cardinality of the set of possible sequences $(x_0,\ldots , x_n)$  is $<\kappa$,
so the $\kappa$-filtered property says that we can choose a single $i$. Then the standard kind of argument
shows that by going further along, the maps $\pA (x_0,\ldots , x_n)\rightarrt \pB _i(y^i_0,\ldots , y^i_n)$
can be assumed to all fit together into a natural transformation in terms of $(x_0,\ldots , x_n)\in \Delta ^o_X$.
Thus we get a factorization of our map $\pA \rightarrt \pB $ through some $\pA \rightarrt \pB _i$. This shows that 
$\pA $ is $\kappa$-presentable.

Suppose on the other hand that $\pA $ is $\kappa$-presentable. We note first of all that $|X|<\kappa$. Indeed,
if $|X|\geq \kappa$ then we could consider a $\kappa$-filtered system of subsets $Z_i\subset X$ with 
$|Z_i|<\kappa$, but $\colim Z_i =X$. Let $\codisc (Z_i)$ and $\codisc (X)$  denote the codiscrete  precategories
on these object sets, that is the precatgories whose value is $\ast$ on any sequence of objects. Then 
$\colim _i \codisc (Z_i) =\codisc (X)$ in $\precat (\mM )$ (see Lemma \ref{connstar}), but the identity map on underlying object sets
gives a map $\pA \rightarrt \codisc (X)$ not factoring through any $\codisc (Z_i)$, contradicting the assumed
$\kappa$-presentability of $\pA $. Hence we may assume that $|X|<\kappa$. 

We claim that $\pA $ is $\kappa$-presentable when considered as an object in $\precat (X; \mM )$. Indeed,
if $\{ \pB _i\} _{i\in \beta}$ is a $\kappa$-filtered
diagram in $\precat (X; \mM )$ then since $\beta$ is connected, any map 
$$
\pA \rightarrt \colim ^{\precat ( X; \mM )}_{i\in \beta} \pB _i
$$
is also a map $\pA \rightarrt \colim ^{\precat ( \mM )}_{i\in \beta} \pB _i$ by 
Corollary \ref{constobjcolim}. By the assumed $\kappa$-presentability of $\pA $ this
would have to factor through one of the $\pB _i$, necessarily as a map inducing the identity on underlying object sets $X$.
This shows that $\pA $ is $\kappa$-presentable when considered as an object in $\precat (X; \mM )$.
Now Corollary \ref{constobjpres} tells us that the $\pA (x_0,\ldots , x_p)$ are $\kappa$-presentable objects of $\mM$.
\end{proof}

\begin{proposition}
\label{precatlocpres}
Suppose $\mM$ is locally presentable. Then the category of $\mM$-precategories $\precat(\mM )$ is also locally presentable.
\end{proposition}
\begin{proof}
Let $\kappa$ be a regular cardinal such that $\mM$ is locally $\kappa$-presentable. 
Note that, for any set $X$ the category $\precat (X,\mM )$ is locally $\kappa$-presentable
by its adjunction with the diagram category $\diag (\Delta _X^o,\mM )$ which in turn
is locally $\kappa$-presentable by Lemma \ref{diagpres}. 

Existence of arbitrary colimits in $\precat (\mM )$
was shown at the start of the section. It is clear
from the description of $\kappa$-presentable objects in Lemma \ref{kappapres}
that the isomorphism classes of $\kappa$-presentable objects of $\precat (\mM )$
form a set. Furthermore,
any object $(Z,\pA )$ is $\kappa$-accessible. Indeed, consider the category of triples $(X,\pB ,u)$ where
$X\subset Z$ is a subset of cardinality $|X| < \kappa$, $\pB \in \precat (X,\mM )$ is a $\kappa$-presentable object, and
$u:\pB \rightarrt \pA |_X$ is a morphism to the pullback of $\pA $ along the inclusion $X\subset Z$.
Using the local $\kappa$-presentability of each $\precat (X,\mM )$
and the expression of $Z$ as a $\kappa$-filtered union of the subsets $X$,
the category of triples is $\kappa$-filtered and the colimit of the tautological 
functor is $(Z,\pA )$. 
\end{proof}

\section{Interpretations as presheaf categories}
\label{sec-interpretations}

With some additional hypotheses on $\mM$, the trick of introducing the categories $\Delta _X$ becomes unnecessary. This corresponds more
closely with some previous references \cite{Tamsamani} \cite{svk} \cite{limits}, and will be useful in establishing notation for iterated
$n$-precatgories. Starting from this first discussion we show later on that if $\mM$ itself is a presheaf category then $\precat (\mM )$ is
a presheaf category. 

This consideration will not usually enter into our argument although it does provide a convenient
change of notation for Chapter \ref{secat1}. Nevertheless, many arguments become considerably simpler in the case of a presheaf category---as we have already seen for cell complexes
in Chapter \ref{cattheor1}. As we shall now see, the passage from $\mM$ to $\precat (\mM )$ preserves the condition of being a presheaf category, and many important initial 
cases such as the model category of simplicial sets $\mK$ satisfy this condition. 
So it should be helpful throughout the book to be able to think 
of the case of presheaf categories.

The first part of our discussion follows Pelissier \cite{Pelissier}. 

Suppose $\mM$ is a locally presentable category. Define a functor $\disc : \Sets \rightarrt \mM$ by 
$$
\disc (U):= \coprod _{u\in U}\ast .
$$
If necessary we shall denote by $\ast _u$ the term corresponding to $u\in U$ in the
coproduct. The discrete object functor 
has a right adjoint $\disc ^{\ast}:\mM \rightarrt \Sets$ defined by
$\disc ^{\ast}(X):= \Hom _{\mM} (\ast , X)$, with
$$
\Hom _{\mM}(\disc (U),X)= \Hom _{\Sets}(U,\disc ^{\ast}(X)).
$$
In particular, $\disc$ commutes with colimits. 

Given an expression $X=\coprod _{u\in U}X_u$ the universal property of the coproduct
applied to the maps $X_u\rightarrt \ast _u$ yields a map $X\rightarrt \disc (U)$.
On the other hand, given a map $f:X\rightarrt \disc (U)$ we can put
$$
f^{-1}(u):= X\times _{\disc (U)}\ast _u.
$$
We have a natural map $\coprod _{u\in U}f^{-1}(u)\rightarrt X$ compatible with 
the maps to $\disc (U)$. 

Consider the following hypothesis on $\mM$, saying that disjoint unions (i.e. colimits over discrete categories) behave well.
This kind of hypothesis was introduced for the same purpose by Pelissier \cite[Definition 1.1.4]{Pelissier}.

\begin{condition}[DISJ]
\label{disj}
$ \; $
\newline
(a)---If $f:X\rightarrt \disc (U)$ is a map, then the natural map is an
isomorphism $\coprod _{u\in U}f^{-1}(u)\cong X$. If $X=\coprod _{u\in U}X_u$ is
a coproduct expression then $X_u=f^{-1}(u)$ for the corresponding map $f:X\rightarrt
\disc (U)$. 
\newline
(b)---The coinitial object $\ast$ is indecomposable, that is to say it cannot be written
as a coproduct of two nonempty objects.
\newline
(c)---The category $\mM$ has more than just a single object up to isomorphism. 
\end{condition}

To see why condition (b) is  necessary, 
note for example that if $\mM = \presh (\Phi )$ is the category of presheaves on a category $\Phi$ which
has more than one connected component (for example the discrete category
with two objects) then $\ast$ is decomposable. 

Condition (c) is required to rule out the trivial category $\mM =\ast$ which satisfies all of
our other hypotheses. 

\begin{lemma}
\label{disjointproperties}
Assume that $\mM$ is locally presentable, 
and satisfies Condition (DCL) of the cartesian condition \ref{def-cartesian}.
Assume that $\mM$ satisfies Condition (DISJ) \ref{disj} above. Then we have the
following further properties:
\newline
(1)---For any object $X\in \mM$, giving an expression $X\cong \coprod _{u\in U}X_u$ is equivalent to giving a map $X\rightarrt \disc (U)$. 
\newline
(2)---If $Y\rightarrt \coprod _{u\in U}X_u$ is a map from a single object to a coproduct, then setting $Y_u:= Y\times _{\coprod _{u\in U}X_u} X_u$
the natural map $\coprod _{u\in U}Y_u\rightarrt Y$ is an isomorphism. 
\newline
(3)---The map $\emptyset \rightarrt \ast$ is not an isomorphism. 
\newline
(4)---For any set $U$ the adjunction map $U\rightarrt \disc ^{\ast}\disc (U)$ is
an isomorphism. 
\newline
(5)---The functor $\disc$ is fully faithful, and $\disc ^{\ast}$ gives an  inverse
on its essential image.
\newline
(6)---If $X\cong \coprod _{u\in U}X_u$ and $Y\cong \coprod _{u\in U}Y_u$ are decompositions corresponding to maps $X,Y\rightarrt \disc (U)$ then 
$$
X\times _{\disc (U)}Y = \coprod _{u\in U}X_u\times Y_u .
$$
(7)---Coproducts are disjoint: if $\{ X_u\}_{u\in U}$ is a collection of objects 
and $u\neq v$ then 
$$
X_u\times _{\coprod _{u\in U}X_u}X_v = \emptyset .
$$
(8)---The functor $\disc$ preserves finite limits. 
\end{lemma}
\begin{proof}
Condition (1) is just a restatement of the first part of (DISJ). 

For (2),
suppose $Y\rightarrt X=\coprod _{u\ni U}X_u$ is a map. Compose with
the map $X\rightarrt \disc (U)$ given by (1), to get $Y\rightarrt \disc (U)$.
In turn this corresponds to a decomposition $Y=\coprod _{u\in U}Y_u$ with
$$
Y_u = Y\times _{\disc (U)}\ast _u,
$$
but $X_u=X\times _{\disc (U)}\ast _u$ so $Y_u = Y\times _XX_u$ which is the desired statement.

For (3), if $\emptyset \rightarrt \ast$ were an isomorphism, then
we would have $\emptyset \times X = X$ for all objects $X$, but in view of
Lemma \ref{emptyempty} following from (DCL) this would imply that all objects are
$\emptyset$ contradicting the nontriviality hypothesis (DISJ) (c) on $\mM$. 

For (4) suppose first we are given a map $f:\ast \rightarrt \disc (U)$. 
By (1) this corresponds to a decomposition $\ast = \coprod _{u\in U}\ast \times _{\disc (U)}\ast _u$. By Condion (DISJ)(b), all but one of the summands must be $\emptyset$.
From (3) it follows that 
one of the summands is different, so there is unique one of the summands
which is $\ast$. We get that our map factors through a unique $\ast _u$. This shows that
the adjunction map $U\rightarrt \disc ^{\ast}\disc (U)$ is an isomorphism. 

For (5), suppose given a map $f:\disc (U)\rightarrt \disc (V)$. By property (4)
there is a unique map $U\rightarrt^g V$ compatible with $\disc ^{\ast}(f)$ by the
adjunction isomorphisms. Uniqueness shows that $\disc $ is faithful. 
Applying $\disc$ to this comptibility diagram and composing with
the naturality square for the other adjunction map gives a square
$$
\begin{diagram}
\disc (U) &\rightarr ^{\disc (g)}& \disc (V) \\
\downarr && \downarr \\
\disc (U) &\rightarr ^{f}& \disc (V) 
\end{diagram}
$$
where the vertical maps are the adjunction compositions 
$$
\disc (U)\rightarrt \disc \disc ^{\ast}\disc (U)\rightarrt \disc (U)
$$
and the same for $V$. These are the identities so $f=\disc (g)$. This shows that $\disc$ is
fully faithful, and the $\disc ^{\ast}$ gives an essential inverse by (4). 

In the situation of (6), given $Z\rightarrt^h \disc (U)$ corresponding to
$Z=\coprod _{u\in U}Z_u$, a map $Z\rightarrt X$ compatible with $h$ is the
same thing as a collection of maps $Z_u\rightarrt X_u$ by (2). Similarly 
for a map to $Y$. It follows that  $\coprod _{u\in U}X_u\times Y_u$ satisfies the
universal property for the fiber product, giving (6). 

For (7), it suffices to show that $\ast _u\times _{\disc (U)}\ast _v = \emptyset$
for $v\neq u$,
in view of the consequence  Lemma \ref{emptyempty} of Condition (DCL) saying that
no nonempty object can map to $\emptyset$. But  
$\ast _u = \coprod _{v\in U}(\ast _u)\times _{\disc (U)}\ast _v$
and, as was seen in the proof of (4),
condition (DISJ) (b) implies that all but one of these terms must be $\emptyset$.
Since there is a diagonal map from $\ast$ to the term $v=u$, the terms for $v\neq u$ must
be $\emptyset$. 

For (8), the functor $\disc $ preserves finite direct products: given maps 
$Z\rightarrt \disc (U)$ and $Z\rightarrt \disc (V)$ they correspond to decompositions
as in (1). The maps $Z_{u}\rightarrt \disc (V)$ correspond to decompositions
$$
Z_u = \coprod _{v\in V}Z_{u,v},\;\;\; Z_{u,v} = Z_u\times _{\disc (V)}\ast _v.
$$
Putting these together over all $u\in U$ we get a decomposition of $Z$ corresponding
to a unique map $Z\rightarrt \disc (U\times V)$. This shows that $\disc (U\times V)$
satisfies the universal property to be the product $\disc (U)\times \disc (V)$.
A similar argument shows that $\disc$ preserves equalizers. 
\end{proof}

Suppose now that $\mM$ satisfies condition (DISJ) in addition to being tractable left proper
and cartesian, 
so the properties of the preceding lemma apply. Given $\pA \in \precat (\mM )$ define
a functor $\Delta ^o \rightarrt \mM$ denoted $[n]\mapsto \pA _{n/}$ by 
$$
\pA _{n/}:= \coprod _{(x_0,\ldots , x_n)\in \Ob (\pA )^{n+1}} \pA (x_0,\ldots , x_n).
$$
The functoriality maps are defined using those of $\pA $. This has the property that $\pA _{0/}$ is a discrete object, indeed the unitality conditions say that $\pA (x_0)=\ast$
so there is a natural 
isomorphism $\pA _{0/}\cong \disc (\Ob (\pA ))$, and we identify $\pA _{0/}$ with the set $\Ob (\pA )$, sometimes using the notation $\pA _0$. 

The property that $\pA _{0/}$ is a discrete object is called the {\em constancy condition}
\cite{Tamsamani}, closely related to the globular nature of the theory
of $n$-categories. Let $\diag ([0]\subset \Delta ^o, \Sets \subset \mM  )$ denote the full subcategory of $\diag (\Delta ^o, \mM )$ consisting of functors
which satisfy this constancy condition. 

Suppose on the other hand that $n\mapsto \pA _{n/}$ is a functor $\Delta ^o \rightarrt \mM$ which satisfies the constancy condition.
Set $X:= \Ob (\pA ):= \disc ^{\ast}( \pA _{0/})$, then $\pA _{0/}= \disc (\Ob (\pA ))$. 
For any sequence $x_0,\ldots , x_n\in X=\Ob (\pA )$,
we get a map 
$$
(x_0,\ldots , x_n):\ast \rightarrt \pA _{0/}\times \cdots \times \pA _{0/}
$$
and define 
$$
\pA (x_0,\ldots , x_n):= \pA _{n/}\times _{\pA _{0/}\times \cdots \times \pA _{0/}}\ast
$$
where the map $\pA _{n/}\rightarrt \pA _{0/}\times \cdots \times \pA _{0/}$ is obtained using the $n+1$ vertices of the simplex $[n]$. Notice that 
$$
\pA _{0/}\times \cdots \times \pA _{0/} = \disc (\Ob (\pA )\times \cdots \times \Ob (\pA ))
$$
since $\disc$ preserves finite products (Lemma \ref{disjointproperties} (8)).
The simplicial maps for $\pA _{\cdot /}$ provide transition maps to make $\pA (\cdots )$ into a functor $\Delta ^p _X\rightarrt \mM$.

\begin{theorem}
\label{interp1}
Suppose $\mM$ satisfies Condition (DISJ) \ref{disj}.
The above constructions provide essentially inverse functors 
$$
\precat (\mM )\rBotharr \diag ([0]\subset \Delta ^o, \Sets \subset  \mM  ).
$$
The full subcategory $\diag ([0]\subset \Delta ^o, \Sets \subset  \mM  )\subset \diag (\Delta ^o, \mM  )$ is closed under limits and colimits,
and the above essentially inverse functors preserve limits and colimits. 
\end{theorem}
\begin{proof}
This follows from Lemma \ref{disjointproperties}: suppose fixed the set of objects $X$,
then giving $\pA _{n/}$ with map
to the direct product
$$
\pA _{n/}\rightarrow \disc (X\times \cdots X) = \disc (X)\times \cdots \times \disc (X)
$$
is the same as giving the pieces of the decomposition 
$$
\pA _{n/}=\coprod _{(x_0,\ldots , x_n)}\pA (x_0,\ldots , x_n).
$$
\end{proof}

For iteration of our basic construction, the following lemmas show that the above notational reinterpretation is very reasonable.

\begin{lemma}
\label{precatdisj}
Suppose $\mM$ satisfies the hypothesis (DCL) of Definition 
\ref{def-cartesian}. Then the category $\precat (\mM )$ satisfies condition (DISJ).
\end{lemma}
\begin{proof}
The discrete precategories are constructed as follows. 
The coinitial object $\ast\in \precat (\mM )$ has a single element $\Ob (\ast ) = \ast = \{ x\}$ and $\ast (x,\ldots , x)=\ast \in \mM $.
Thus, if $U$ is any set then $\disc (U)$ calculated in $\precat (\mM )$ has object set $\disc (U)\cong U$,
and $\disc (U)(x_0,\ldots , x_n)= \ast$ if $x_0= \cdots = x_n$, otherwise $\disc (U)(x_0,\ldots , x_n)= \emptyset$ if the sequence
is not constant.

For part (a), suppose $\pA \in \precat (\mM )$ and $\pA \rightarrt \disc (U)$ is a map. Put $X:= \Ob (\pA )$, with $X\rightarrt U$
which corresponds to a decomposition $X=\coprod _{u\in U}X_u$. Let $\pA _u\subset \pA $ be the pullback of $\pA $ along $X_u\rightarrt X$,
that is it is the ``full sub-precategory'' with object set $X_u\subset X$. This is indeed $\pA \times _{\disc (U)}\ast _u$. The map 
$$
\coprod _{u\in U}\pA _u\rightarrt \pA 
$$
is an isomorphism. This follows from the fact that $\disc (U)(x_0,\ldots , x_n)=\emptyset$ for any nonconstant sequence,
plus the consequence Lemma \ref{emptyempty} of (DCL). On the other hand given a decomposition in coproduct
$\pA =\coprod _{u\in U}\pA _u$ we get a map $\pA \rightarrt \disc (U)$ for which the components
$\pA _u$ are the fibers. 

Conditions (b) and (c) are easy. 
\end{proof}

\begin{lemma}
\label{startdisj}
Suppose $\Psi$ is a connected category (i.e. its nerve is a connected space, or equivalently any two objects
are joined by a zig-zag of arrows). Then the category $\presh (\Psi )= \diag (\Psi ^o; \Sets )$ satisfies condition (DISJ), with the
discrete objects being those equivalent to constant functors. 
In particular, 
the model category $\Sets$ where all morphisms are fibrations and cofibrations, and the weak equivalences are isomorphisms,
satisfies (DISJ). And the model category $\mK$ of simplicial sets satisfies (DISJ).
\end{lemma}
\begin{proof}
If $U$ is a set, the discrete presheaf $\disc (U)$ is the constant presheaf
$x\mapsto U$ on all $x\in \Psi$. For (a), a map $\pA \rightarrt \disc (U)$ is the same
thing as a decomposition $\pA =\coprod _{u\in U}\pA _u$, as can be seen levelwise. Condition (b)
follows from the connectedness of $\Psi$ (indeed it is equivalent), and Condition (c) is easy. 
\end{proof}

We now turn to the further situation where  $\mM$ is a  presheaf category, say $\mM = \presh (\Psi )= \diag (\Psi ^o; \Sets )$.
In this case
$$
\diag (\Delta ^o_X ; \mM ) = \diag (\Delta ^o_X \times \Psi ^o; \Sets ) = \presh (\Delta _X \times \Psi )
$$
is again a presheaf category. 

We construct a new category denoted $\scone (\Psi )$, which looks somewhat like a
``cone'' on $\Phi$. It is defined by 
contracting $\{ [0]\} \times \Psi \subset  \Delta  \times \Psi$ to a single object denoted $0$. Thus, the objects of $\scone (\Psi )$
are of the form $([n],\psi )$ for $n\geq 1$ and $\psi \in \Psi$, or the object $0$. The morphisms are as follows:
\newline
---there is a single identity morphism between $0$ and itself;
\newline
---for $n\geq 1$ and $\psi \in \Psi$ there is a unique morphisms from any $([n],\psi )$ to $0$;
\newline
---the morphisms from $0$ to $([n],\psi )$ are the same as the morphisms $[0]\rightarrt \nocom [n]$ in $\Delta$ (i.e. there are $n+1$ of them); and
\newline
---the morphisms from $([n],\psi )$ to $([n'],\psi ')$ are of two kinds: either $(a,f)$ where $a:[n]\rightarrt \nocom [n']$ is
a morphism in $\Delta$ such that $a$ doesn't factor through $[0]$ and $f:\psi \rightarrt \psi '$ is a morphism in $\Psi$;
or else $(a)$ where $a:[n]\rightarrt \nocom [n']$ is
a morphism in $\Delta$ which factors as $[n]\rightarrt \nocom [0]\rightarrt \nocom [n']$.
\newline
Composition of morphisms is defined in an obvious way. Notice that the composition of anything with a morphism which factors through $[0]$ will again factor through $[0]$,
which allows us to define compositions of the form $(a)\circ (a',f')$ or $(a,f)\circ (a')$ in the last case. The division of the morphisms
from $([n],\psi )$ to $([n'],\psi ')$ into two cases is necessary in order to define composition.

\begin{proposition}
\label{upresheafcat}
Suppose $\mM = \presh (\Psi )$ is a presheaf category. 
There is a natural isomorphism between 
$\presh (\scone (\Psi ))$ and the category of unital $\mM$-precategories $\precat(\mM )$. Thus, $\precat(\mM )$ 
is again a presheaf category. 
\end{proposition}
\begin{proof}
Suppose $\pF  :\scone (\Psi )^o\rightarrt \Sets $ is a presheaf. Let $X:= \pF  (0 )$. For $n\geq 1$ and any $(x_0,\ldots , x_n)\in \Delta _X$,
let $\pA (x_0,\ldots , x_n)$ be the presheaf on $\Psi$ which assigns to $\psi \in\Psi $ the subset of elements of $\pF  ([n],\psi )$
which project to $(x_0,\ldots , x_n)$ under the $n+1$ projection maps $\pF ([n],\psi )\rightarrt \pF (\iota )=X$ 
corresponding to the $n+1$ maps $0 \rightarrt ([n],\psi )$. At $n=0$, set $\pA (x_0):= \{x_0 \}$.
The pair $(X,\pA )$ is an element of $\precat (\mM )$. Conversely,
given any $(X,\pA )\in \precat (\mM )$, define a presheaf $\pF :\scone (\Psi )^o\rightarrt \Sets $
by setting $\pF (0 ):= X$ and
$$
\pF ([n],\psi ):= \coprod _{(x_0,\ldots , x_n)\in X^{n+1}}\pA (x_0,\ldots , x_n)(\psi ).
$$
These constructions are inverses.
\end{proof}

Let $\sconemap _{\Psi}: \Delta \times \Psi \rightarrt \scone (\Psi )$ denote the projection. 
If $\Psi$ is connected, then 
the category of presheaves $\presh (\scone (\Psi ))$ may be identified, via $\sconemap _{\Psi}^{\ast}$, with the full subcategory of $\presh (\Delta \times \Psi )$ consisting
of presheaves $\pA $ which satisfy {\em Tamsamani's constancy condition} that
$\pA (0,\psi )$ is a constant set independent of $\psi \in \Psi$.

The construction $\Psi \mapsto \scone (\Psi )$ was the iterative step in the construction of the sequence of categories denoted $\Theta ^n$ in \cite{svk}.
In that notation, $\Theta ^0 = \ast$ and $\Theta ^{n+1} = \scone (\Theta ^n)$. However, the notation $\Theta ^n$ was subsequently used by Joyal \cite{JoyalTheta}
to denote a related but different sequence of categories; in order to avoid confusion we will use the cone notation $\scone$. 

For the theory of non-unital $\mM$-precategories, 
one can also construct a category denoted $\sconeplus (\Psi )$ with the property that 
$\presh (\sconeplus  (\Psi ))$
is the category of non-unital $\mM$-precategories. This is in fact even more straightforward
than the construction of $\scone$ for the unital theory. The following discussion is
optional, but serves to put the previous discussion of $\scone$ into a better perspective. 

Recall that $\precat (\mM )$ was defined as a fibered category over $\Sets$ whose fiber over a set $X$ was $\precat (X,\mM )$.
In the same way, the {\em category of non-unital $\mM$-precategories} is the fibered category over $\Sets$ whose
fiber over $X$ is $\diag (\Delta ^o_X, \mM )$. 

The construction of $\sconeplus   (\Psi )$ 
is to formally add an object denoted $\iota$ to $\Delta \times \Psi$. 
The objects of $\sconeplus   (\Psi )$ are of the form either $\iota$ or $([n],\psi )$ for $n\in \Delta$ and $\psi \in \Psi$.
The morphisms of $\sconeplus   (\Psi )$ are defined as follows: 
\newline
---there is a single identity morphism between $\iota$ and itself;
\newline
---there are no morphisms from $([n],\psi )$ to $\iota$;
\newline
---the morphisms from $([n],\psi )$ to $([n'],\psi ')$ are the same as the morphisms in $\Delta \times \Psi$; and
\newline
---the morphisms from $\iota$ to $([n],\psi )$ are the same as the morphisms 
$[0]\rightarrt \nocom [n]$ in $\Delta$ (i.e. there are $n+1$ of them).
\newline
Composition of morphisms is defined so that the single automorphism of $\iota$ is the identity, so that
composition of morphisms within $\Delta \times \Psi$ is the same as from that category, and a composition of the form
$$
\iota \rightarrt ([n],\psi ) \rightarrt ([n'],\psi ')
$$
is given by the corresponding composition $[0]\rightarrt \nocom [n]\rightarrt \nocom [n']$.

\begin{proposition}
\label{presheafcat}
Suppose $\mM = \presh (\Psi )$ is a presheaf category. 
There is a natural isomorphism between $\presh (\sconeplus   (\Psi ))$ and the category of non-unital $\mM$-precategories.
In particular, the latter is a presheaf category. 
\end{proposition}
\begin{proof}
Suppose $\pF :\sconeplus   (\Psi )^o\rightarrt \Sets $ is a presheaf. Let $X:= \pF (\iota )$. For any $(x_0,\ldots , x_n)\in \Delta _X$,
let $\pA (x_0,\ldots , x_n)$ be the presheaf on $\Psi$ which assigns to $\psi \in\Psi $ the subset of elements of $\pF ([n],\psi )$
which project to $(x_0,\ldots , x_n)$ under the $n+1$ projection maps $\pF ([n],\psi )\rightarrt \pF (\iota )=X$ 
corresponding to the $n+1$ maps $\iota \rightarrt ([n],\psi )$. Then pair $\pA $ is an element of $\diag (\Delta ^o_X, \mM )$. Conversely,
given any $\pA \in \diag (\Delta ^o_X, \mM )$, define a presheaf $\pF :\sconeplus   (\Psi )^o\rightarrt \Sets $
by setting $\pF (\iota ):= X$ and
$$
\pF ([n],\psi ):= \coprod _{(x_0,\ldots , x_n)\in X^{n+1}}\pA (x_0,\ldots , x_n)(\psi ).
$$
These constructions are inverses.
\end{proof}

The inclusion functor $U^{\ast}$ from $\precat(\mM )$ to the category of non-unital precategories, can be viewed as a pullback. Indeed, there is a projection functor
$$
\pi : \sconeplus   (\Psi )\rightarrt \scone (\Psi )
$$
defined by $\pi ([n],\psi ) = ([n],\psi )$ if $n\geq 1$, $\pi ([0],\psi ) = 0$ and $\pi (\iota ) = 0$. The pullback $U^{\ast}$ is just the pullback functor 
$$
\pi ^{\ast} : \precat (\mM )\cong \presh (\scone (\Psi ))\rightarrt \presh (\sconeplus   (\Psi )).
$$
The left adjoint of the pull back or inclusion, denoted pushforward $U_{!}$ above, is also the pushforward $\pi _{!}$ for the functor $\pi$. 

Using the operations $\sconeplus   (\Psi )$ and $\scone (\Psi )$ we may stay essentially entirely within the realm of presheaf categories.
The only place where we go outside of there is when we speak of the category  of unital $\mM$-precategories over a fixed set of objects $X$;
even if $\mM$ is a presheaf category, $\precat(X; \mM )$ will not generally be a presheaf category. On the other hand, the non-unital version 
$\diag (\Delta ^o_X,  \mM )= \presh (\Delta _X \times \Psi )$ remains a presheaf category.


\chapter{Algebraic theories in model categories}
\label{algtheor1}

In this chapter we consider algebraic diagram theories 
consisting of a collection of finite product conditions imposed on
diagrams $\Phi \rightarrow  \mM$. This is motivated by the situation considered
in the previous Chapter \ref{precat1}. There we defined the notion of precategory
on a fixed set of objects, which is a diagram $\pA :\Delta ^o_X\rightarrow \mM$. The
Segal conditions require that certain maps be weak equivalences. 
Imposing these conditions amounts to a homotopical analogue of the ``finite product theories''
often considered in category theory \cite{Lawvere} \cite{LawvereThesis},
see further historical remarks in \cite[pp 171-172]{AdamekRosicky}.
The homotopical analogue, whose origins go back to the various theories of
$H$-spaces and loop spaces,
was considered by Badzioch in \cite{Badzioch}, and his treatment was
used by Bergner for Segal categories in \cite{BergnerThreeModels} \cite{BergnerSMSC},
and also generalized to models in simplicial categories in \cite{BergnerRigidification}. 
Rosicky carried these ideas further in \cite{RosickyHomotopyVarieties} and he points
out more references. 

Let $\epsilon (n)$ denote the category 
$$
{\setlength{\unitlength}{.5mm}
\begin{picture}(60,60)
\put(0,30){\ensuremath{\xi _0}}
\put(10,30){\circle*{1.5}}
\put(55,55){\circle*{1.5}}
\put(60,55){\ensuremath{\xi _1}}
\put(55,45){\circle*{1.5}}
\put(60,45){\ensuremath{\xi _2}}
\put(55,35){\circle*{1.5}}
\put(60,35){\ensuremath{\xi _3}}

\put(40,25){\ensuremath{\vdots}}

\put(55,15){\circle*{1.5}}
\put(60,15){\ensuremath{\xi _{n-1}}}
\put(55,5){\circle*{1.5}}
\put(60,5){\ensuremath{\xi _n}}

\qbezier(10,30)(30,42)(50,53)
\qbezier(10,30)(30,37)(50,44)
\qbezier(10,30)(30,32)(50,35)
\qbezier(10,30)(30,23)(50,16)
\qbezier(10,30)(30,18)(50,6)

\put(50,53){\vector(3,2){0}}
\put(50,44){\vector(3,1){0}}
\put(50,35){\vector(4,1){0}}
\put(50,16){\vector(3,-1){0}}
\put(50,6){\vector(3,-2){0}}
\end{picture}
}
$$
The Segal maps are obtained by pulling back $\pA $ along functors 
$\epsilon (n)\rightarrow \Delta ^0_X$. This sets up a localization problem
which can be phrased in more general terms. We treat the general situation in the
present chapter. Even ignoring any possible other applications, that simplifies
notations for the general aspects of the problem of enforcing the given collection
of finite product conditions. 

By an {\em algebraic theory} we mean a category $\Phi$ provided with a collection of 
``direct product diagrams'', that is diagrams with the shape of a direct product, which are
functors $\epsilon (n)\rightarrt^P \Phi$.  
A realization of the theory in a classical $1$-category $\Cc$
is a functor $\Phi \rightarrt \Cc$ which sends
these diagrams to direct products in $\Cc$. Many of the easiest kinds of structures can be written
this way, although it is well understood that to get more complicated structures, one needs to go
to the notion of {\em sketch} which is a category provided with more generally shaped limit diagrams. 
For our purposes, it will suffice to consider direct product diagrams.

Suppose $\mM$ is an appropriate kind of model category. Then a homotopy realization of the theory in $\mM$
is a functor $\Phi \rightarrt \mM$ which sends the direct product diagrams, to homotopy direct products
in $\mM$. 

These  notions lead to a ``calculus of generators and relations'' where we start with an arbitrary functor
$\Phi \rightarrt \Cc$ (resp. $\Phi \rightarrt \mM$) and try to enforce the direct product (resp. homotopy 
direct product) condition. The main work of this chapter will be to do this for the case of homotopy realizations
in a model category $\mM$.

\section{Diagrams over the categories $\epsilon (n)$}

The first task is to take a close look at the categories indexing direct products. 
Let $\epsilon (n)$ denote the category with objects $\xi _0, \xi _1,\ldots , \xi _n$;
and whose only morphisms apart from the identities, are single morphisms $\rho _i: \xi _0 \rightarrt \xi_i$.
We also let $\xi _i$  denote the functors $\xi _i:\{ \ast \}\rightarrt \epsilon (n)$ sending the point $\ast$ 
to the object $\xi _i(\ast ) = \xi _i$. 

This includes the case $n=0$ where $\epsilon (0)$ is the discrete category with one object $\xi_0$. 

Suppose $\Cc$ is some other category. A functor $\pA :\epsilon (n)\rightarrt \Cc$ corresponds
to a collection of objects $a _0,a _1,\ldots ,a _n\in ob (\Cc )$ together with maps 
$p_i:a _0\rightarrt a_i$
for $1\leq i\leq n$. We sometimes write 
$$
\pA = (a_0,\ldots , a_n; p_1,\ldots , p_n).
$$
Say that $\pA $ is a {\em direct product diagram} if the collection of maps $p_i$ expresses $a_0$ as
a direct product of the $a_1,\ldots , a_n$ in $\Cc$. 

Suppose $\mM$ is a model category, in particular direct products exist. We say that a diagram $\pA = (a_0,\ldots , a_n; p_1,\ldots , p_n)$
from $\epsilon (n)$ to $\mM$ is a {\em homotopy direct product diagram} if the morphism 
$$
(p_1,\ldots , p_n): a_0\rightarrt a_1\times \cdots \times a_n
$$
is a weak equivalence.

Consider now the diagram category $\diag (\epsilon (n), \mM )$. We describe explicitly its injective and projective model structures,
which are well-known to exist (see Hirschhorn \cite{Hirschhorn}, Barwick \cite{Barwick}). 

A morphism of diagrams 
$$
\pA = (a_0,\ldots , a_n; p_1,\ldots , p_n)\rightarrt  
\pB = (b_0,\ldots , b_n; q_1,\ldots , q_n)
$$ 
is a collection of maps $g=(g_0,\ldots , g_n)$ with $g_i:a_i\rightarrt b_i$ such that
$q_ig_0 = g_ip_i$. We can also think of a morphism $g$ as consisting of the maps $(g_1,\ldots , g_n)$
plus a map 
$$
g_P : a_0 \rightarrt b_0 \times _{b_1\times \cdots \times b_n} (a_1\times \cdots \times a_n)
$$
such that the second projection is the structural map for $\pA $. 
We then write $g= \langle g_1,\ldots , g_n; g_P\rangle $.

A morphism $g$ is fibrant in the projective structure, if and only if each $g_0,\ldots , g_n$ is fibrant in $\mM $.
Similarly, $g$ is cofibrant in the injective structure, if and only if each $g_0,\ldots , g_n$ is cofibrant in $\mM$. 

To study fibrant maps in the injective
structure, suppose 
$$
\pC = (c_0,\ldots , c_n; r_1,\ldots , r_n), \;\;\;\; 
\pD = (d_0,\ldots , d_n; s_1,\ldots , s_n)
$$
are two diagrams, with morphisms forming a square
$$
\begin{diagram}
\pA  & \rightarr^{u} & \pC  \\
\downarr^g & & \downarr _h\\
\pB  & \rightarr^{v} & \pD  &.
\end{diagram} 
$$
These give diagrams
$$
\begin{diagram}
a_i & \rightarr^{u_i} & c_i \\
\downarr^{g_i} & & \downarr _{h_i}\\
b_i & \rightarr^{v_i} & d_i & .
\end{diagram} 
$$
We look for a lifting $f: \pB \rightarrt \pC $ such that $fg = u$ and $hf = v$.
This amounts to asking for liftings $f_i:b_i\rightarrt c_i$ such that $f_ig_i = u_i$ and $h_if_i = v_i$,
and also such that $s_if_0 = f_ir_i$. 

Suppose given already the liftings $f_1,\ldots , f_n$. Then to specify a full lifting $f$ as above we need to find
$f_0: b_0\rightarrt c_0$ such that the diagram
$$
\begin{diagram}
a_0 & \rightarr & c_0  \\
\downarr ^{g_0}& \ruTeXto^{f_0}& \downarr \\
b_0 & \rightarr^{v_P} & d_0 \times _{d_1\times \cdots \times d_n} (c_1\times \cdots \times c_n)
\end{diagram} 
$$
commutes. 

\begin{lemma}
\label{injfib}
A morphism $h: \pC \rightarrt \pD $ with notations as above, is fibrant in the injective model structure on $\diag (\epsilon (n), \mM )$
if and only if each $h_1,\ldots , h_n$ is fibrant in $\mM$, and the map 
$$
h_P : c_0 \rightarrt d_0 \times _{d_1\times \cdots \times d_n} (c_1\times \cdots \times c_n)
$$
is fibrant in $\mM$. 
\end{lemma}
\begin{proof}
If $h$ satisfies the conditions stated in the lemma, and $f$ is a levelwise
cofibration i.e.\ each $f_i$ is a cofibration,
then we can first
choose the liftings $f_i$ for $i=1,\ldots , n$, then choose $f_0$ lifting $h_P$.
Therefore $h$ is fibrant. Suppose on the other hand that $h$ is fibrant. 
For any cofibration $a\rightarrt b$ and any $i=1,\ldots , n$ we have
an injective cofibration between objects with $a$ (resp. $b$) placed at $i$ and
the remaining places filled in with $\emptyset$. Since $h$ satisfies lifting
along any such cofibration, it follows that $h_i$ is fibrant in $\mM$.
Similarly, given a cofibration $a\rightarrt b$ with $a\rightarrt c_0$
and $b\rightarrt d_0 \times _{d_1\times \cdots \times d_n} (c_1\times \cdots \times c_n)$,
we get a square diagram as above with $\pA =(a,b,b\ldots , b)$ and $\pB =(b,b,\ldots , b)$.
The map $\pA \rightarrow \pB $ is a levelwise cofibration, so the fibrant condition
for $h$ implies existence of a lifting, which gives a lifting
for the map $h_P$.
\end{proof}

Consider now the same question in the other direction: choose first the lifting $f_0$. 
For any $i\geq 1$ we get a map $b_0 \sqcup ^{a_0} a_i  \rightarrt  c_i$, and
the choices of $f_i$ for $i=1,\ldots , n$ correspond to choices of lifting in the diagrams
$$
\begin{diagram}
b_0 \sqcup ^{a_0} a_i & \rightarr & c_i  \\
\downarr & \ruTeXto^{f_i} & \downarr_{h_i} \\
b_i & \rightarr^{v_i} & d_i
\end{diagram} .
$$

\begin{lemma}
\label{projcofib}
A morphism $g: \pA \rightarrt \pB $ with notations as above, is cofibrant in the projective model structure on $\diag (\epsilon (n), \mM )$
if and only if $g_0$ is cofibrant in $\mM$, and for each $i=1,\ldots , n$ the map 
$b_0 \sqcup ^{a_0} a_i  \rightarrt  b_i$ is cofibrant in $\mM$.
In particular an object $\pB $ is cofibrant if and only if $b_0$ is cofibrant and
each $b_0\rightarrt b_i$ is a cofibration. 
\end{lemma}
\begin{proof}
Similar to the previous proof. 
\end{proof}

One can describe explicit generating sets for the cofibrations and trivial cofibrations in
both structures $\diag _{\rm inj}(\epsilon (n), \mM )$ and $\diag _{\rm proj}(\epsilon (n), \mM )$.

Recall the standard adjunctions for the functors $\xi _i: \ast \rightarrt \mM$. If
$X\in \mM$ is an object, we obtain a diagram $\xi _{i,!}(X): \epsilon (n)\rightarrt \mM$.
Explicitly, if $i=0$ then $\xi _{0,!}(X)$ is the constant diagram $(X,X,\ldots , X; 1_X, \ldots , 1_X)$
with values $X$. If $i\geq 1$ then $\xi _{0,!}(X)$ is the diagram $(\emptyset , \ldots , \emptyset , X ,\emptyset \ldots ; \iota , \ldots , \iota )$
where $\iota$ denote the  unique maps from $\emptyset$ to anything else, and here $X$ is at the $î$-th place. 
The adjunction says that for any diagram $\pA = (a_0,\ldots , a_n; p_1,\ldots , p_n)$, a
morphism $\xi _{i,!}(X)\rightarrt \pA $ is the same thing as a morphism $X\rightarrt \xi _i^{\ast}(\pA )=a_i$.

Suppose given generating sets $I$ for the cofibrations of $\mM$ and $J$ for the trivial cofibrations. 

For $f:X\rightarrt Y$ in $I$, we obtain a cofibration $\xi _{0,!}(f): \xi _{0,!}(X)\rightarrt \xi _{0,!}(Y)$
in the projective model structure. To see this, use Lemma \ref{projcofib} on $g= \xi _{0,!}(f)$
and note that $g_0$ is just $f$ so it is cofibrant; and $b_0 \sqcup ^{a_0} a_i = Y\sqcup ^X X = Y$ maps to $b_i=Y$
by a cofibration. If $f\in J$ then $\xi _{0,!}(f)$ is a trivial cofibration: over each object of $\epsilon (n)$ we just
get back the map $f$ so it is an levelwise weak equivalence. 

For $f\in I$ as above, at any $i\geq 1$, $\xi _{i,!}(f)$ is a cofibration, again using Lemma \ref{projcofib} on $g= \xi _{0,!}(f)$.
The map $g_0$ is the identity of $\emptyset$, and the maps 
$b_0 \sqcup ^{a_0} a_j  \rightarrt  b_j$ are either the identity of $\emptyset$ for $j\neq i$, or $f$ when $j=i$,
so these are cofibrations. Again, if $f\in J$ then $\xi _{i,!}(f)$ is a trivial cofibration: over each object 
of $\epsilon (n)$ it gives either the identity of $\emptyset$ which is automatically a weak equivalence, or else $f$.

Define $I_{\epsilon (n)}$ (resp. $J_{\epsilon (n)}$) 
to be the set consisting of diagrams of the form $\xi _{i,!}(f)$ for $0\leq i\leq n$ and $f\in I$ (resp. $f\in J$).

\begin{proposition}
\label{projepsilongen}
The sets $I_{\epsilon (n)}$  and $J_{\epsilon (n)}$ are generators for the projective model category structure
$\diag _{\rm proj}(\epsilon (n), \mM )$.
\end{proposition}
\begin{proof}
By the adjunction, a morphism $g$ in $\diag (\epsilon (n), \mM )$ satisfies the right lifting property with
respect to $I_{\epsilon (n)}$ (resp. $J_{\epsilon (n)}$) if and only if $\xi _i^{\ast}(g)$ satisfies
the right lifting property with respect to $I$ (resp. $J$) for all $0\leq i\leq n$. Since $I$ (resp. $J$)
is a set of generators for the cofibrations (resp. trivial cofibrations) of $\mM$, this lifting
property is equivalent to saying that each $\xi _i^{\ast}(g)$ is a trivial fibration (resp. a fibration).
By the definition of the projective model structure, this is equivalent to saying that $g$ is a 
trivial fibration (resp. a fibration). Hence $\inj (I_{\epsilon (n)})$ is the class of trivial fibrations
and $\inj (J_{\epsilon (n)})$ is the class of fibrations, so $\cof (I_{\epsilon (n)})$ is the class of
cofibrations and $\cof (J_{\epsilon (n)})$ is the class of
trivial cofibrations. 
\end{proof}

To get generators for the injective model structure, we need to add a new kind of injective cofibration. 
If $f:X\rightarrt Y$ is a cofibration, consider the diagrams 
\begin{equation}
\label{varpidef}
\varpi (f):= \xi _{0,!}(X) \cup ^{\coprod _i \xi _{i,!}(X)}\coprod _i \xi _{i,!}(Y) = (X, Y, \ldots , Y; f,\ldots , f)
\end{equation}
and $\xi _{0,!}(Y)$. We have a map $\varpi (f)\rightarrt \xi _{0,!}(Y)$ which is $f$ at the object $\xi _0(\ast )$ and $1_Y$ at
the other objects. Denote this map by $\rho (f)=(f,1,\ldots , 1)$.

\begin{proposition}
\label{injepsilongen}
Let $I^+_{\epsilon (n)}$ denote the union of $I_{\epsilon (n)}$ with the set of 
maps of the form $\rho (f)=(f,1,\ldots , 1)$ for $f\in I$. 
Let $J^+_{\epsilon (n)}$ denote the union of $J_{\epsilon (n)}$ with the set of 
maps of the form $\rho (f)$ for $f\in J$.  Then
$I^+_{\epsilon (n)}$ and $J^+_{\epsilon (n)}$ are generating sets for 
the injective model category structure
$\diag _{\rm inj}(\epsilon (n), \mM )$.
\end{proposition}
\begin{proof}
If $f: X\rightarrt Y$ is a cofibration, then 
$(f,1,\ldots , 1): \pU \rightarrt \pV $ is a cofibration with the notations as above.  Suppose that $g:\pA \rightarrt \pB $ satisfies 
right lifting with respect to $f$, where $g= (g_0,\ldots , g_n)$ goes from
$\pA = (a_0,\ldots , a_n; p_1,\ldots , p_n)$ to
$\pB = (b_0,\ldots , b_n; q_1,\ldots , q_n)$. The lifting property says that for any map
$u_0:X\rightarrt a_0$ and maps $u_i:Y\rightarrt a_i$ for $i\geq 1$, and maps $v_i: Y\rightarrt b_i$ for $i\geq 0$
such that $p_iu_0 = u_i f$, $v_i = g_iu_i$ for $i\geq 1$, and $v_0f = g_0 u_0$, then there should exist
a map $u'_0: Y\rightarrt a_0$ such that $u'_0 f = u_0$, $g_0 u'_0 = v_0$, and $p_iu'_0 = u_i$. 
This is the same as the right lifting property for the square
$$
\begin{diagram}
X & \rightarr & a_0 \\
\downarr & & \downarr \\
Y & \rightarr & b_0 \times _{b_1\times \cdots \times b_n} (a_1\times \cdots \times a_n)
\end{diagram}
$$
so the condition that $g$ satisfies right lifting with respect to any $(f,1,\ldots , 1)$ for $f\in J$ is equivalent to the
condition that $a_0\rightarrt b_0 \times _{b_1\times \cdots \times b_n} (a_1\times \cdots \times a_n)$
is a fibration. Thus, $\inj (J^+_{\epsilon (n)})$ consists of maps which are levelwise fibrations
(because of lifting with respect to $J_{\epsilon (n)}$) and such that the map
of Lemma \ref{injfib} is a fibration. 
By Lemma \ref{injfib}, $\inj (J^+_{\epsilon (n)})$ is the class of fibrations. 

Similarly, the fact that $g$ satisfies right lifting with respect to any $(f, 1,\ldots , 1)$ for $f\in I$ is
equivalent to the conditoin that 
$$
a_0\rightarrt b_0 \times _{b_1\times \cdots \times b_n} (a_1\times \cdots \times a_n)
$$
be a trivial fibration.  Thus, $\inj (I^+_{\epsilon (n)})$ consists of maps which are levelwise trivial fibrations 
(because of lifting with respect to $I_{\epsilon (n)}$) and such that the map
of Lemma \ref{injfib} is a trivial fibration. 

We claim that this is equal to the class of fibrations. 
If $g\in \inj (I^+_{\epsilon (n)})$, then it is a fibration for the injective structure
since $J^+_{\epsilon (n)}\subset I^+_{\epsilon (n)}$, and also levelwise a weak equivalence, so
it is a trivial fibration. If $g$ is a trivial fibration, then it is levelwise a trivial fibration, and
the map $a_0\rightarrt b_0 \times _{b_1\times \cdots \times b_n} (a_1\times \cdots \times a_n)$ is a fibration.
However, the maps $a_i\rightarrt b_i$ are trivial fibrations, so 
$a_1\times \cdots \times a_n\rightarrt b_1\times \cdots \times b_n$ is a trivial fibration. Thus, the map 
$b_0 \times _{b_1\times \cdots \times b_n} (a_1\times \cdots \times a_n)\rightarrt b_0$
is a trivial fibration, and by 3 for 2 we conclude that 
$a_0\rightarrt b_0 \times _{b_1\times \cdots \times b_n} (a_1\times \cdots \times a_n)$ is a weak equivalence.
Hence it is a trivial fibration as required for the claim.

This identifies $\inj (I^+_{\epsilon (n)})$ and $\inj (J^+_{\epsilon (n)})$ with the classes of trivial fibrations and fibrations
respectively,
so $\cof (I^+_{\epsilon (n)})$ and $\cof (J^+_{\epsilon (n)})$ are the classes of cofibrations and trivial cofibrations respectively. 
\end{proof}

\begin{scholium}
If $(\mM , I, J)$ is a cofibrantly generated (resp. combinatorial, tractable, left proper)
model category, then 
$\diag _{\rm inj}(\epsilon (n), \mM )$ and $\diag _{\rm proj}(\epsilon (n), \mM )$ are cofibrantly generated 
(resp. combinatorial, tractable, left proper) model categories. 
\end{scholium}

\noindent
See Theorem \ref{injectiveprojective}.

\section{Imposing the product condition}

Assume that $\mM$ is a tractable left proper cartesian model category.
Recall that the cartesian condition \ref{def-cartesian}
implies that
for any $X\in \mM$ and weak equivalence $f:Y\rightarrt Z$ the induced map $X\times Y\rightarrt X\times Z$ is a
weak equivalence. 

We are going to apply the direct left Bousfield localization theory of
Chapter \ref{direct1}. 
Say that an object $\pA \in \diag (\epsilon (n), \mM )$ is {\em product-compatible} if the map 
$\pA (\xi _0) \rightarrt \pA (\xi _1)\times \ldots \times \pA (\xi _n)$ is a weak equivalence. Let $\Rr _{\epsilon (n)}\subset \diag (\epsilon (n), \mM )$
denote the full subcategory of product-compatible objects. 

\begin{lemma}
\label{Rinvariant}
Suppose $\mM$ is cartesian. 
The subcategory $\Rr_{\epsilon (n)}$ is invariant under weak equivalence: if $\pA $ is product-compatible and 
its image  in $\Ho  \diag (\epsilon (n), \mM )$
is isomorphic to the image of another object $\pB $, then $\pB $ is also product-compatible.
\end{lemma}
\begin{proof}
If $\pA \rightarrow \pB $ is a levelwise weak equivalence, then the horizontal arrows in the
diagram
$$
\begin{diagram}
\pA (\xi _0) & \rightarr & \pB (\xi _0) \\
\downarr & & \downarr \\
\pA (\xi _1)\times \ldots \times \pA (\xi _n) &\rightarr & \pB (\xi _1)\times \ldots \times \pB (\xi _n)
\end{diagram}
$$
are weak equivalences, using the cartesian condition for the bottom arrow. Hence, one
of the vertical arrows is a weak equivalence if and only if the other one is. 
\end{proof}

\subsection{Direct localization of the projective structure}
\label{sec-dirlocproj}

We now define the direct localizing system which goes with the full subcategory of product-compatible objects in the projective model category structure. 
For each generating cofibration $f:X\rightarrt Y$ in $I$,
recall the morphism $\rho (f) = (f,1,\ldots , 1)$ defined above (see \eqref{varpidef})
$$
\varpi (f)= (X, Y, \ldots , Y; f,\ldots , f) \rightarr^{\rho (f)} 
\xi _{0,!}(Y) = (Y,\ldots , Y; 1,\ldots , 1).
$$
It is obviously an injective cofibration, and the domain $\varpi (f)$ is
projectively cofibrant (apply Lemma \ref{projcofib}). However, $\rho (f)$ will
not in general be a projective cofibration. Choose a factorization 
$$
Y\cup ^X Y\rightarrt^{a(f)} Z\rightarrt ^{b(f)}Y
$$
such that $a(f)$ is a cofibration and $b(f)$ is a trivial fibration. 
Denoting by $i_2$ the first inclusion $Y\rightarrt Y\cup ^XY$ we get a
trivial cofibration $a(f)i_2:Y\rightarrt Z$. 
Let
$$
\varpi (f)= (X, Y, \ldots , Y; f,\ldots , f)\rightarrt^{\zeta (f)} 
\psi (f):= (Y,Z,\ldots , Z; a(f)i_2,\ldots , a(f)i_2)
$$
be the map given by $f$ over $\xi _0$ and $a(f)i_1$ over $\xi _1,\ldots , \xi _n$.
Then the map occuring in Lemma \ref{projcofib} is exactly $a(f)$, so 
$\zeta (f)$ is a projective cofibration.

Recall the explicit generating set $J_{\epsilon (n)}$ for the trivial cofibrations in the projective model structure. Put 
$$
K^{\rm proj}_{\epsilon (n)}:= J_{\epsilon (n)} \cup \{ \zeta (f)\} _{f\in I} .
$$

\begin{theorem}
\label{diagpi}
The pair $(\Rr _{\epsilon (n)}, K^{\rm proj}_{\epsilon (n)})$ is a direct localizing system for the
model category $\diag _{\rm proj}(\epsilon (n), \mM )$ of $\epsilon (n)$-diagrams in $\mM$ 
with the projective model structure. Let $\diag _{{\rm proj}, \Pi }(\epsilon (n), \mM )$ denote the
left Bousfield localized model structure constructed in Chapter \ref{direct1}. 
A morphism $\pA \rightarrt \pB $ is a weak equivalence if and only if it induces 
weak equivalences levelwise over the objects $\xi _1,\ldots , \xi _n$. 
An object $\pA $ is fibrant in the 
localized structure if
and only if it is product-compatible and each $\pA (\xi _i)$ is fibrant. 
\end{theorem}
\begin{proof}
We verify the properties (A1)--(A6). Properties (A1) and (A2) are immediate. 
For (A3) we have chosen $\zeta (f)$ so as to be a projectively cofibration, 
and its domain is projectively cofibrant. Condition (A4) is given by Lemma \ref{Rinvariant}.

Suppose a diagram $\pA $ is in $\inj (K^{\rm proj}_{\epsilon (n)})$. In particular it is
$J_{\epsilon (n)}$-injective, which is to say fibrant in the projective model structure.
This means that each $\pA (\xi _i)$ is a fibrant object in $\mM$. 
The lifting property along the $\zeta (f)$ implies the following homotopy lifting
property of the map
$$
p:\pA (\xi _0)\rightarrt \pA (\xi _1)\times \cdots \times \pA (\xi _n).
$$
along a generating cofibration $f:X\rightarrow Y$ in $I$. 
Recall that $Y\cup ^XY\rightarrt^a Z \rightarrt^b Y$ was chosen above.  If we are given a
diagram 
$$
\begin{diagram}
X & \rightarr^{u} & \pA (\xi _0) \\
\downarr^f & & \downarr_p \\
Y & \rightarr^{v} & \pA (\xi _1)\times \cdots \times \pA (\xi _n)
\end{diagram}
$$
then there is a map $Y\rightarrt^w \pA (\xi _0)$ such that $wf=u$, and a map $Z\rightarrt^h
\pA (\xi _1)\times \cdots \times \pA (\xi _n)$ such that $hai_2 = v$ but $hai_1 = pw$. 
This homotopy lifting property implies that $p$ is a weak equivalence,
see Lemma \ref{homotopylifting}.
Thus, $\pA \in \Rr _{\epsilon (n)}$ which gives (A5).

For condition (A6), suppose $\pA $ is in $\Rr _{\epsilon (n)}$ and
$\pA \rightarrt \pB $ is a pushout along an element of $K^{\rm proj}_{\epsilon (n)}$.
Applying the small object argument we can find a map $\pB \rightarrt \pC $ in 
$\cell (K^{\rm proj}_{\epsilon (n)})$ such that
$\pC \in \inj (K^{\rm proj}_{\epsilon (n)})$. Then $\pA \rightarrt \pC $ is also in
$\cell (K^{\rm proj}_{\epsilon (n)})$. Notice, however, that the elements of
$K^{\rm proj}_{\epsilon (n)}$ are levelwise trivial cofibrations over the
objects $\xi _1,\ldots , \xi _n\in \epsilon (n)$ (but not over $\xi  _0$).
Therefore the map $\pA \rightarrow \pC $ induces weak equivalences over each $\xi _1,
\ldots , \xi _n$. Using the cartesian condition on $\mM$, this implies that
the right vertical map in the square
$$
\begin{diagram}
\pA (\xi _0) & \rightarr & \pA (\xi _1)\times \cdots \times \pA (\xi _n) \\
\downarr & &\downarr \\
\pC (\xi _0) & \rightarr & \pC (\xi _1)\times \cdots \times \pC (\xi _n)
\end{diagram}
$$
is a weak equivalence. The hypothesis that $\pA \in \Rr _{\epsilon (n)}$
says that the top map is a weak equivalence, and the fact that 
$\pC \in \inj (K^{\rm proj}_{\epsilon (n)})$ and part (A5) proved above
say that $\pC \in \Rr _{\epsilon (n)}$, so the bottom map is a weak equivalence.
By 3 for 2 the left vertical arrow is a weak equivalence, showing that $\pA \rightarrt \pC $ is
a levelwise weak equivalence of diagrams. This shows (A6). 

Our direct localizing system leads to a left Bousfield localization by Theorem 
\ref{directLBL}. 

We now look at the characterizations of new weak equivalences. 
As seen above, the elements of $\cell (K^{\rm proj}_{\epsilon (n)})$ are
weak equivalences levelwise over the $\xi _1,\ldots , \xi _n$. Using the 
characterizatin of Corollary \ref{nwecriterion} we see that all new weak equivalences are
levelwise weak equivalences over the $\xi _1,\ldots , \xi _n$. Suppose
$\pA \rightarrt^f \pB $ is a morphism inducing a weak equivalence over each
$\xi _1,\ldots , \xi _n$. Choose $\pB \rightarrt^b \pB '$ in 
$\cell (K^{\rm proj}_{\epsilon (n)})$ such that $\pB '$ is in 
$\inj (K^{\rm proj}_{\epsilon (n)})$, in particular it is product-compatible. 
Factor the composed map as
$$
\pA \rightarrt ^{a}\pA '\rightarrt ^g \pB '
$$
where $a\in \cell (K^{\rm proj}_{\epsilon (n)})$ and 
$g\in \inj (K^{\rm proj}_{\epsilon (n)})$. All of the above maps are levelwise
weak equivalences over the objects $\xi _1,\ldots , \xi _n$, however $\pA '$ and $\pB '$
are product-compatible. It follows that $g$ is a levelwise weak equivalence.
The criterion of Corollary \ref{nwecriterion} implies that $f$ is a new
weak equivalence. 

The characterization of fibrant objects is a first version in the simplified
situation of a single product diagram, of Bergner's characterization of fibrant
Segal categories \cite{BergnerSegal}. 
The new fibrant objects are in $\inj (K^{\rm proj}_{\epsilon (n)})$ so
they are in $\Rr _{\epsilon (n)}$ i.e. product-compatible, and levelwise fibrant. 
Suppose $\pA $ is product-compatible and levelwise fibrant. Suppose $\pU \rightarrt^f \pV $
is a new trivial cofibration, and we are given a map $\pU \rightarrt^u \pA $. 
By the previous paragraph it induces a levelwise trivial 
cofibration over the objects $\xi _1,\ldots , \xi _n$. Hence the components $u_1,\ldots , u_n$
extend to maps $\pV (\xi _i)\rightarrt^{v'_i} \pA (\xi _i)$. Putting these together, the composition
$$
\pV (\xi _0)\rightarrt \pV (\xi _1)\times \cdots \times \pV (\xi _n)\rightarrt
\pA (\xi _1)\times \cdots \times \pA (\xi _n)
$$
gives the bottom arrow of the diagram
$$
\begin{diagram}
\pU (\xi _0) & \rightarr^{u_0} & \pA (\xi _0) \\
\downarr^{f_0} & & \downarr \\
\pV (\xi _0) & \rightarr & \pA (\xi _1)\times \cdots \times \pA (\xi _n) .
\end{diagram}
$$
The right vertical arrow is a weak equivalence between fibrant objects, so by Lemma 
\ref{homotopylifting} there is a homotopy lifting relative $\pU (\xi _0)$,
in other words a map $\pV (\xi _0)\rightarrt^{v_0} \pA (\xi _0)$ such that $rf_0=u_0$,
and the other triangle commutes up to a homotopy relative $\pU (\xi _0)$. 
Lemma \ref{htyext} 
says we can change the maps $v'_i$ to maps $v_i$, still restricting to $u_i$ on
$\pU (\xi _i)$, but compatible with $v_0$. We have now constructed the required extension
$\pV \rightarrt \pA $, showing that $\pA $ is a new fibant object. 
\end{proof}

\subsection{Direct localization of the injective structure}
\label{sec-dirlocinj}

When possible, it is more convenient to use the injective model structure on $\diag (\epsilon (n), \mM )$. 
Consider the explicit generating set $J^+_{\epsilon (n)}$ for the trivial cofibrations in the injective model structure, given by Proposition \ref{injepsilongen}. 
We can define two different sets of cofibrations, the first extending 
$K^{\rm proj}_{\epsilon (n)}$:
$$
K^{{\rm inj}+}_{\epsilon (n)}:=K^{\rm proj}_{\epsilon (n)}\cup J^{+}_{\epsilon (n)}
=J^{+}_{\epsilon (n)} \cup \{ \zeta (f)\} _{f\in I};
$$
and the second defined using the simpler maps $\rho (f)$ which were already
injective cofibrations: 
$$
K^{\rm inj}_{\epsilon (n)}:= J^{+}_{\epsilon (n)} \cup \{ \rho (f)\} _{f\in I} .
$$

\begin{theorem}
\label{injdiagpi}
The pairs $(\Rr _{\epsilon (n)}, K^{{\rm inj}+}_{\epsilon (n)})$ 
and $(\Rr _{\epsilon (n)}, K^{\rm inj}_{\epsilon (n)})$
are both direct localizing systems for the
model category $\diag _{\rm inj}(\epsilon (n), \mM )$ of $\epsilon (n)$-diagrams in $\mM$ 
with the injective model structure. Let $\diag _{{\rm inj}, \Pi }(\epsilon (n), \mM )$ denote the
left Bousfield localized model structure, which is the same in both cases. 
The weak equivalences are the same
as for the projective structure. 
An object $\pA $ is fibrant in the 
localized structure if
and only if it is product-compatible and satisfies the fibrancy criterion of Lemma
\ref{injfib} for the injective model structure. 
The identity functor is a left Quillen functor
$$
\diag _{{\rm proj}, \Pi }(\epsilon (n), \mM )\rightarrt
\diag _{{\rm inj}, \Pi }(\epsilon (n), \mM )
$$
from the new projective to the new injective model structure. 
\end{theorem}
\begin{proof}
The functor 
$\diag _{{\rm proj}}(\epsilon (n), \mM )\rightarrt
\diag _{{\rm inj} }(\epsilon (n), \mM )$
is a left Quillen functor, whose corresponding right Quillen functor (both being the identity on underlying categories) preserves the class $\Rr _{\epsilon (n)}$
of product-compatible diagrams. In other words the transfered class is the same.
The subset 
$K^{{\rm inj}+}_{\epsilon (n)}$ is the transfered subset given in Theorem \ref{transfer},
so by that theorem $(\Rr _{\epsilon (n)}, K^{{\rm inj}+}_{\epsilon (n)})$ is
a direct localizing system and the identity functor is a left Quillen functor
from the previous new projective model structure to the
resulting left Bousfield localization of the injective structure 
$$
\diag _{{\rm proj}, \Pi }(\epsilon (n), \mM )\rightarrt
\diag _{{\rm inj}, \Pi }(\epsilon (n), \mM ) .
$$

The proof that $(\Rr _{\epsilon (n)}, K^{{\rm inj}}_{\epsilon (n)})$ is
a direct localizing system is the same as in the proof of the previous
Theorem \ref{diagpi}, but in fact easier since an object which satisfies lifting
with respect to the $\rho (f)$ has the stronger property that 
$$
\pA (\xi _0)\rightarrt \pA (\xi _1)\times \cdots \times \pA (\xi _n)
$$
is in $\inj (I)$, that is it is a trivial fibration. So in this case we don't need
to rely on the notion of homotopy lifting property as was done in the previous proof.
We get conditions (A1)--(A6) and also the same description of weak equivalences,
and the corresponding description of fibrant objects. 

The two model structures given by Theorem \ref{directLBL} applied to 
$(\Rr _{\epsilon (n)}, K^{{\rm inj}}_{\epsilon (n)})$ 
and $(\Rr _{\epsilon (n)}, K^{{\rm inj}+}_{\epsilon (n)})$ are the same,
by Proposition \ref{lblindep}.
\end{proof}

Suppose $\pA \in \diag (\epsilon (n), \mM )$. Suppose given a factorization 
$$
\pA (\xi _0)\rightarr^{e_0} E_0 \rightarr^{p} \pA (\xi _1)\times \cdots \times \pA (\xi _n)
$$
in $\mM$. Let $E_i:= \pA (\xi _i)$ for $i=1,\ldots , n$. The structural map $p$ gives a structure of $\epsilon (n)$-diagram
to the collection $(E_0,\ldots , E_n)$, call it $\pE$. The map $e$ gives a map  $e:\pA \rightarrt \pE$. 
If $e_0$ is a cofibration in $\mM$ then $e$ is a cofibration in $\diag _{\rm inj}(\epsilon (n), \mM )$.

\begin{lemma}
\label{standardinjtcf}
In the above situation, if $e_0$ is a cofibration and $p$ is a weak equivalence in $\mM$ then $e:\pA \rightarrt \pE$ is a trivial 
cofibration in $\diag _{{\rm inj}, \Pi }(\epsilon (n), \mM )$.
\end{lemma}
\begin{proof}
The map $e$ is levelwise cofibrant by construction. It is a weak equivalence 
since it induces a weak equivalence levelwise over the objects $\xi _1,\ldots , \xi _n$. 
\end{proof}

\subsection{Transfering these structures}

Putting together the above analysis of diagrams over $\epsilon (n)$ with the transfer along a Quillen functor
gives the following general picture. Suppose we are given a set $Q$, integers $n(q)\geq 0$ for $q\in Q$,
a family of tractable left proper cartesian model categories
$\mM _q$ for $q\in Q$, a tractable left proper model category $\mN$,
and a family of Quillen functors $F_q : \diag _{\rm proj}(\epsilon (n(q)), \mM _q )\rBotharrow \mN : G_q$.

Let $\Rr '\subset \mN$ be the full subcategory of objects $Y$ such that, for a fibrant replacement $Y\rightarrt Y'$,
the diagrams $G_q(Y'): \epsilon (n(q))\rightarrt \mM _q$ are product-compatible. 
Let $(I^q,J^q)$ be generating sets for $\mM _q$, and $(I',J')$ generators for $\mN$. 

\begin{corollary}
\label{dirlocgen}
Let $K'$ be the union of $J'$, of the set of morphisms of the form $F_q (g)$ for $g\in J^q_{\epsilon (n(q))}$,
and of the set of morphisms of the form $F_q (\zeta (f))$ for $f\in I^q$. Then $(\Rr ',K')$ is a direct localizing system
for $\mN$. 

If furthermore $F_q$ are left Quillen functors from 
$\diag _{\rm proj}(\epsilon (n(q)), \mM _q )$ to $\mN$, then we can consider
$K^{\rm inj}$, the union of $J'$ with the 
set of morphisms of the form $F_q (\rho (f))$ for $f\in I^q$, and 
$(\Rr ',K^{\rm inj})$ is a direct localizing system for $\mN$ giving the same model
structure as  $(\Rr ',K')$. 
\end{corollary}
\begin{proof}
Let $K''$ be the union of $K'$ with the set of morphisms of the form $F_q (g)$ for $g\in J^q_{\epsilon (n(q))}$. Then 
$(\Rr ',K'')$ is a direct localizing system
for $\mN$, by Theorem \ref{transferfamily} applied to the direct localizing systems of Theorem \ref{diagpi}. However, the $F_q (g)$ for $g\in J^q_{\epsilon (n(q))}$
are trivial cofibrations between cofibrant objects in the 
original model structure of $\mN$, so they could be included in a bigger
generating set $J''$ for the original trivial cofibrations of $\mN$. 
But one can note that in the construction of a direct localizing system by 
adding on some new morphisms to the original generating set, the properties
are independent of the choice of original generating set. So $K'$ works as well as $K''$. 

Suppose now that $F_q$ remain left Quillen functors when we use the 
injective model structures on their sources. Then, with a similar discussion for 
leaving out the images of the morphisms in $J^{q,+}_{\epsilon (n(q))}$, Theorem
\ref{transferfamily} applies to the direct localizing systems of
Theorem \ref{injdiagpi} to conclude that $(\Rr ',K^{\rm inj})$ is a direct localizing system.
As pointed out in Proposition \ref{lblindep}, the resulting model structure 
is the same as for $(\Rr ',K')$. 
\end{proof}

\section{Algebraic diagram theories}

Classically, an ``algebraic theory'' is given by a small category $\Phi$ and a collection of product diagrams $P_q: \epsilon (n(q))\rightarrt \Phi$.
The objects of the theory are the functors $\pA :\Phi \rightarrt Set$ with the property that $p_q^{\ast}(\pA )$ is a direct product, that is
{\em product-compatible} in the above terminology. Of course this theory has since been much generalized, to include the notion of ``finite limit sketches''
among other things. However, for our purposes it will be sufficient to consider just the basic version of the theory, and to give it a weak-enriched
counterpart using the notion of direct left Bousfield localization we have developped so far. 

So, suppose $\Phi$ is a small category, $Q$ is a small set, we have integers $n(q)\geq 0$ for $q\in Q$, and suppose given 
functors $P_q : \epsilon (n(q))\rightarrt \Phi$ for $q\in Q$.

For the coefficients, fix a 
tractable left proper model category $\mM$ satisfying condition (PROD). Let $(I,J)$ be a set of generators for $\mM$.
Playing the role of the model category $\mN$
will be the category of $\Phi$-diagrams in $\mM$ with its projective or injective model structure, $\mN =\diag _{\rm proj}(\Phi , \mM )$
(resp. $\mN =\diag _{\rm proj}(\Phi , \mM )$).
Let $(I_{\Phi ,  {\rm proj}}, J_{\Phi,  {\rm proj}})$ and 
$(I_{\Phi ,  {\rm inj}}, J_{\Phi,  {\rm inj}})$ be the sets of generators for the projective and 
injective model structures respectively, see the
discussion of references for Theorem \ref{injectiveprojective}. 
Recall that the identity functor is a Quillen adjunction between the projective and injective diagram categories
$$
{\bf 1}: \diag _{\rm proj}(\Phi , \mM )\rBotharrow 
\diag _{\rm inj}(\Phi , \mM ): {\bf 1} .
$$

\begin{lemma}
With the above notations, for each $q\in Q$ we get a Quillen adjunction 
$$
P_{q,!} : \diag _{\rm proj}(\epsilon (n(q)), \mM )\rBotharrow \diag _{\rm proj}(\Phi , \mM ): P_q^{\ast} .
$$
This composes with the identity functor to give a Quillen adjunction 
$$
{\bf 1}P_{q,!} : \diag _{\rm proj}(\epsilon (n(q)), \mM )\rBotharrow \diag _{\rm inj}(\Phi , \mM ): P_q^{\ast} {\bf 1}.
$$
\end{lemma}
\begin{proof}
This is
the standard Quillen adjunction between Bousfield projective model structures coming from
the functor $P_q:\epsilon (n(q))\rightarrt \Phi$.
\end{proof}

In the above situation, let $\Rr (\Phi , P_{\cdot}, \mM )$ denote the full subcategory of $\diag (\Phi , \mM )$
consisting of diagrams $\pA $ such that for a fibrant replacement $\pA \rightarrt \pA '$ in
the injective model structure (which is also a fibrant replacement in the projective model structure),
for all $q\in Q$, $P_q^{\ast}(\pA ')\in \diag (\epsilon (n(q)),\mM )$ is product-compatible.
Let $K_{{\rm proj}/{\rm inj}}(\Phi , P_{\cdot}, \mM , I,J)$ be the sets given by Corollary \ref{dirlocgen} for the  
projective/injective structure, consisting of
the elements of $J_{\Phi , {\rm proj}/{\rm inj}}$, of the $P_{q,!}(J_{\epsilon (n(q))})$, and of the $P_{q,!}(\rho _{\epsilon (n(q))}(f))$ for $f\in I$. 
Here a choice of subscript ${\rm proj}$ or ${\rm inj}$ is proposed when necessary. 

\begin{theorem}
\label{diagprod}
The pair $(\Rr (\Phi , P_{\cdot}, \mM ), K_{{\rm proj}/{\rm inj}}(\Phi , P_{\cdot}, \mM , I,J))$ is a direct localizing system. Define the {\em model category
of weak $(\Phi , P_{\cdot})$-algebras in $\mM$} denoted by $\Alg _{{\rm proj}/{\rm inj}}(\Phi , P_ {\cdot};  \mM )$, to be
the direct left Bousfield localization of the projective or injective diagram model category $\diag _{{\rm proj}/{\rm inj}}(\Phi , \mM )$ with
respect to the full subcategory $\Rr (\Phi , P_{\cdot}, \mM )$ and sets $K_{{\rm proj}/{\rm inj}}(\Phi , P_{\cdot}, \mM , I,J)$. 
The cofibrations are levelwise cofibrations in the injective structure, and projective diagram cofibrations in the projective structure. 
The fibrant objects of $\Alg _{{\rm proj}/{\rm inj}}(\Phi , P_ {\cdot};  \mM )$ are the diagrams $\pA :\Phi \rightarrt \mM $ such that $\pA $
is levelwise fibrant (for the projective structure) or fibrant in the injective diagram structure, 
and for each $q\in Q$ the pullback $P_q^{\ast}(\pA ):\epsilon (n(q))\rightarrt \mM $ is product-compatible.
Given a map $f:\pA \rightarrt \pB $ of $\Phi$-diagrams in $\mM$, the following conditions are equivalent:
\newline
---$f$  is a weak equivalence in $\Alg _{{\rm proj}/{\rm inj}}(\Phi , P_ {\cdot};  \mM )$
\newline 
---for any square diagram 
$$
\begin{diagram}
\pA  & \leftarr & \pA ' & \rightarr & \pA '' \\
\downarr & & \downarr & & \downarr \\
\pB  & \leftarr & \pB ' & \rightarr & \pB ''
\end{diagram}
$$
such that the left horizontal arrows are projective/levelwise cofibrant replacements,  and $\pA '',\pB '' \in \Rr (\Phi , P_{\cdot}, \mM )$, 
then the morphism $\pA ''\rightarrt \pB ''$ is an levelwise weak equivalence; 
\newline
---there exists a square diagram as above with $\pA ''\rightarrt \pB ''$ an levelwise weak equivalence; 
\newline
---there exists a square diagram as above with $\pA ''\rightarrt \pB ''$ an levelwise weak equivalence,
but without the requirement $\pA '',\pB '' \in \Rr (\Phi , P_{\cdot}, \mM )$.
\newline
These projective and injective model categories of weak algebras are tractable and left proper. 
\end{theorem}
\begin{proof}
Apply the construction of the direct left Bousfield localization given in the previous chapter,
starting with either $\diag _{{\rm proj}}(\Phi , \mM )$ or
$\diag _{{\rm inj}}(\Phi , \mM )$ of Theorem \ref{injectiveprojective},
and transfering the direct localizing systems as in Corollary \ref{dirlocgen}. 

For the characterization of fibrant objects, see the successive statements
of Proposition \ref{furtherinfo}, Theorem \ref{directLBL},
Remark \ref{transferfurtherinfo} and Theorem \ref{transferfamily}.
Apply these together with the characterizations of fibrant objects
in the product-compatible $\epsilon (n)$-diagram model structures.  
\end{proof}

\section{Unitality}

\label{sec-unitality}

Suppose given a full subcategory $\Phi _0\subset \Phi$. Typically, these will be the $P_q(\xi _0)$ for $q\in Q$ such that $n(q)=0$.
We would like to consider diagrams $\pA :\Phi \rightarrt \mM$ such that $\pA (x)=\ast $ for $x\in \Phi _0$. Call such a diagram {\em unital along $\Phi _0$}.

Let $\diag (\Phi /\Phi _0, \mM )\subset \diag (\Phi , \mM )$ denote the full subcategory of diagrams which are unital along $\Phi _0$. 
Denote by 
$$
U_{\Phi _0}^{\ast}: \diag (\Phi /\Phi _0, \mM )\rightarrt \diag (\Phi , \mM ) 
$$
the identity inclusion functor. 

The idea for this notation is that $\Phi /\Phi _0$ represents the contraction of $\Phi _0$ to a point, the result being a pointed category
i.e. a category with distinguished object; and $\diag (\Phi /\Phi _0, \mM )$ is the category of pointed functors from here to the category
$\mM$ pointed by distinguishing the coinitial object $\ast$. 

The present discussion plays an important role in the theory of weakly enriched precategories: the unitality condition corresponds to Tamsamani's constancy
condition in the case of $n$-precategories, corresponding to the idea of having a globular theory in which the objects form a discrete set. The motivating
example, first introduced in Section \ref{sec-fixedprecats} above, is when $\Phi =  \Delta _X$ and $\Phi _0=X$ is the subcategory of sequences of length $0$. 

For any $x\in \Phi$ denote by $\Phi _0/x$ the category of arrows $z\rightarrt x$ with $z\in \Phi_0$.

\begin{theorem}
\label{unitaltheorylocpres}
If $\mM$ is a locally presentable category, then the category $\diag (\Phi /\Phi _0, \mM )$ is locally presentable
and $U_{\Phi _0}^{\ast}$ has a left adjoint $U_{\Phi _0, !}$. The left
adjoint is given as follows: if $\pA \in \diag (\Phi , \mM )$ then
$U_{\Phi _0, !}\pA $ is the diagram which sends an object $x\in \Phi$ to the coproduct of $\pA (x)$ and $\colim _{\Phi _0 / x}\ast$
over $\colim _{z\in \Phi _0 / x}\pA (z)$. The adjunction map $U_{\Phi _0, !}U_{\Phi _0}^{\ast}\pB \rightarrt \pB $ is the 
identity for any $\pB \in \diag (\Phi /\Phi _0, \mM )$, so $U_{\Phi _0, !}$ is a monadic projection in the terminology of Section \ref{sec-monadic}. 

The full subcategory $\diag (\Phi /\Phi _0, \mM )\subset \diag (\Phi , \mM )$ is closed under 
small limits and over colimits with small nonempty
connected index sets, in particular it is closed under coproducts, filtered colimits, and transfinite composition. 

For any regular cardinal $\kappa$ with $\mM$ being locally $\kappa$-presentable and
$|\Phi |<\kappa$, an object  $\pA \in \diag (\Phi /\Phi _0, \mM )$ is $\kappa$-presentable if and only if each of the $\pA (x)$ are
$\kappa$-presentable in $\mM$. 
\end{theorem}
\begin{proof}
Put
$$
U_{\Phi _0, !}(\pA )(x):= \pA (x)\cup ^{\colim _{z\in \Phi _0 / x}\pA (z)} \colim _{\Phi _0 / x}\ast .
$$
Given a morphism $x\rightarrt y$, we obtain morphisms 
$$
\colim _{\Phi _0 /x}\ast \rightarrt \colim _{\Phi _0 / y}\ast
$$
and 
$$
\colim _{z\in \Phi _0 / x}\pA (z) \rightarrt \colim _{z\in \Phi _0 / y}\pA (z).
$$
These are compatible with the maps in the above coproduct so they give $U_{\Phi _0, !}(\pA )$ a structure of diagram (i.e. functor).
If $x\in \Phi _0$ then $x$ is the coinitial object of $\Phi _0/x$,
and we get $U_{\Phi _0, !}(\pA )(x)=\pA (x)\cup ^{\pA (x)}\ast = \ast$. Thus $U_{\Phi _0, !}(\pA )\in \diag (\Phi /\Phi _0, \mM )$. 

To show adjunction, suppose $\pB \in \diag (\Phi /\Phi _0, \mM )$. Given a map $\pA \rightarrt U_{\Phi _0}^{\ast}\pB $
then for any $z\in \Phi _0$ the map $\pA (z)\rightarrt \pB (z)$ factors through $\ast$ (since indeed $\pB (z)=\ast$).
For any $x$ this gives a factorization
$$
\begin{diagram}
\colim _{z\in \Phi _0 / x}\pA (z) &\rightarr &\colim _{z\in \Phi _0 / x}\ast \\
\downarr  & & \downarr \\
\pA (x) & \rightarr &\pB (x)
\end{diagram}
$$
so our map of diagrams factors through a unique map $U_{\Phi _0, !}(\pA )\rightarrt \pB $. 

If $\pA = U_{\Phi _0}^{\ast}\pB $ is a diagram with $\pA (z)=\ast$ for $z\in \Phi _0$ already, then the second map in 
the coproduct defining $U_{\Phi _0, !}(\pA )$ is the identity, so $U_{\Phi _0, !}(\pA )=\pA $, i.e. the adjunction is a monadic projection.

Closure under arbitrary small limits is automatic since $U_{\Phi_0}^{\ast}$
is a right adjoint. 
For closure under connected colimits, suppose $\alpha$ is an index category with connected nerve. Then $\colim _{\alpha} \ast = \ast$
where $\ast$ is the coinitial object of $\mM$ and the colimit is taken over the constant functor $\alpha \rightarrt \mM$
(Lemma \ref{connstar}).
As colimits in $\diag (\Phi , \mM )$ are calculated levelwise, it follows that colimits over $\alpha$ preserve the condition 
for inclusion in $\diag (\Phi /\Phi _0, \mM )$ which is that the diagram take values $\ast$ levelwise over $\Phi  _0$. 
Note that this property says that connected colimits in $\diag (\Phi /\Phi _0, \mM )$ are calculated levelwise. That wouldn't be true,
however, for disconnected colimits such as disjoint sums. 

We now identify the $\kappa$-presentable objects of $\diag (\Phi /\Phi _0, \mM )$. 
If each  
$\pA (x)$ is a $\kappa$-presentable object of $\mM$, then
by Lemma \ref{diagpres}
$\pA $ is $\kappa$-presentable in 
$\diag (\Phi , \mM )$. If we are given 
a $\kappa$-filtered
system $\{ \pB _i\} _{i\in \beta}$ in $\diag (\Phi /\Phi _0,  \mM )$, any map
$$
\pA \rightarrt \colim ^{\diag (\Phi /\Phi _0 \mM )}_{i\in \beta} \pB _i
$$
is also a map to $\colim ^{\diag (\Phi , \mM )}_{i\in \beta} \pB _i$ by 
the closure under connected colimits; hence it factors through one of the $\pB _i$.
This factorization is a morphism in the full subcategory $\diag (\Phi /\Phi _0, \mM )$,
which shows that $\pA $ is $\kappa$-presentable in 
$\diag (\Phi /\Phi _0, \mM )$. 

Suppose on the other hand that $\pA $ is $\kappa$-presentable 
in $\diag (\Phi /\Phi _0, \mM )$. 
Given our assumption that $\mM$ is locally $\kappa$-presentable and
$|\Phi |<\kappa$, the category $\diag (\Phi , \mM )$ is locally $\kappa$-presentable,
and its $\kappa$-presentable objects are exactly the diagrams $\pB $ such that $\pB (x)$ is
$\kappa$-presentable in $\mM$. This was stated as Lemma \ref{diagpres}
with reference to \cite{AdamekRosicky}. In particular, we can express $\pA $ as a colimit
in the category $\diag (\Phi , \mM )$
$$
\pA = \colim _{i\in \beta} \pB _i
$$
with $\pB _i(x)$ being $\kappa$-presentable, and indexed by a $\kappa$-filtered category $\beta$.
The unitalization functor $U_{\Phi _0, !}$ being a left adjoint, we get
\begin{equation}
\label{Acolimexpr}
\pA =\colim _{i\in \beta} U_{\Phi _0, !}(\pB _i)\mbox{  in  }\diag (\Phi /\Phi _0 \mM ).
\end{equation}
The hypothesis that $\pA $ is $\kappa$-presentable 
in $\diag (\Phi /\Phi _0, \mM )$, applied to the identity map of $\pA $, says that
the identity factors through a map $\pA \rightarrt U_{\Phi _0, !}(\pB _i)$.
On the other hand, the explicit description of 
$U_{\Phi _0, !}$ shows that each $U_{\Phi _0, !}(\pB _i)(x)$ is $\kappa$-presentable.
We have a retraction 
$$
\pA (x)\rightarrt  U_{\Phi _0, !}(\pB _i)(x) \rightarrt \pA (x)
$$
the composition being the identity of $\pA (x)$. It easily follows that $\pA (x)$ is $\kappa$-presentable in $\mM$.
This completes the proof of the identification of $\kappa$-presentable objects of $\diag (\Phi /\Phi _0 \mM )$.

It is clear from this description that the $\kappa$-presentable objects form a small set. 
The above argument shows that any  $\pA \in \diag (\Phi /\Phi _0 \mM )$ is a $\kappa$-filtered 
colimit of $\kappa$-presentable objects, indeed we obtained the expression \eqref{Acolimexpr}.
\end{proof}

We can construct the injective model category structure.

\begin{proposition}
\label{unitalinjective}
Suppose $\mM$ is tractable. 
There exists a tractable injective model category structure $\diag _{\rm inj}(\Phi /\Phi _0, \mM )$ where the weak equivalences are levelwise
weak equivalences, and cofibrations are levelwise cofibrations.
If $\mM$ is left proper then so is the model category
$\diag _{\rm inj}(\Phi /\Phi _0, \mM )$.
\end{proposition}
\begin{proof}
Weak equivalences are clearly closed under retracts and satisfy 3 for 2. 
The class of trivial cofibrations, that is the intersection of the classes of cofibrations and weak equivalences,
is closed under pushout and transfinite composition since these colimits are calculated levelwise. 

The sets of injective cofibrations and injective trivial cofibrations have generating sets, as was shown in 
Theorem \ref{lurie-unital} using Lurie's technique of Theorems \ref{lurie-thm1} and \ref{lurie-thm2}. This 
gives the necessary accessibility argument which allows to apply Smith's recognition lemma
to obtain the model structure, such as described in \cite{Barwick}. If $\mM$ is left proper, the colimits involved in this condition are connected 
so they are computed  levelwise, hence the same condition holds for $\diag _{\rm inj}(\Phi /\Phi _0, \mM )$. 
\end{proof}

And the projective structure. 

\begin{proposition}
\label{unitalprojective}
If $\Phi _0\subset \Phi$ and 
$\mM$ is a tractable left proper model category, 
we get a projective model structure $\diag _{\rm proj}(\Phi /\Phi _0, \mM )$ which is a tractable left proper model category.
Furthermore the identity on the underlying category constitutes a left Quillen functor
$$
\diag _{\rm proj}(\Phi /\Phi _0, \mM )\rightarrt 
\diag _{\rm inj}(\Phi /\Phi _0, \mM ),
$$
and the unitalization construction is a left Quillen functor 
$$
\diag _{\rm proj}(\Phi , \mM )\rightarr^{U_{\Phi _0, !}}
\diag _{\rm proj}(\Phi /\Phi _0, \mM ).
$$
\end{proposition}
\begin{proof}
The weak equivalences in $\diag _{\rm proj}(\Phi /\Phi _0, \mM )$ are defined to be the levelwise weak equivalences.
These satisfy 3 for 2 and are closed under retracts. 
The fibrations are defined to be the levelwise fibrations. The trivial fibrations are the intersection of these classes.
The cofibrations are determined by the left lifting property with respect to trivial fibrations.  
We construct explicitly a generating set, by a small variant of Bousfield's original construction. 

Choose a generating set $I$ for the cofibrations of $\mM$.
It leads to the set $I_{\Phi}$ of generators for cofibrations in the projective model structure
$\diag _{\rm proj}(\Phi , \mM )$ discussed in Theorem \ref{injectiveprojective}. Recall 
that $I_{\Phi}$ consists of all morphisms of the 
form $i_{x,!}(f)$ where $f:A\rightarrt B$ is in $I$, where $i_x:\{ x\}\rightarrt \Phi$ is the inclusion of a discrete single object
and where 
$$
i_{x,!}: \mM = \diag (\{ x\} , \mM )\rightarrt \diag (\Phi , \mM )
$$ 
is the corresponding left adjoint functor. 
This is just Bousfield's classic generating set for projective diagram cofibrations. 
Set 
$$
I_{\Phi /\Phi _0}:= U_{\Phi _0, !}(I_{\Phi}). 
$$
Notice that 
$$
U_{\Phi _0, !}i_{x,!}: \mM = \diag (\{ x\} , \mM )\rightarrt \diag (\Phi /\Phi _0, \mM )
$$
is the left adjoint functor for inducing unital diagrams from objects of $\mM$ placed over $x\in \Ob (\Phi )$.
The set $I_{\Phi /\Phi _0}$ consists of all $U_{\Phi _0, !}i_{x,!}(f)$ where $f$ runs through the set $I$ of
generating cofibrations for $\mM$ and $x$ runs through $\Ob (\Phi )$. 

For $x$ and $f:A\rightarrt B$ fixed,
$$
U_{\Phi _0, !}i_{x,!}(f):U_{\Phi _0, !}i_{x,!}(A)\rightarrt U_{\Phi _0, !}i_{x,!}(B)
$$ 
has the following explicit description. For any object $y\in \Phi$, let $\Phi _{\rm nf}(x,y)$ denote the set of arrows
from $x$ to $y$ which don't factor through an objects of $\Phi _0$, and let 
$\Phi _{\rm f}(x,y)$ denote the set of arrows which factor through an element of $\Phi_0$. 
Thus $\Phi (x,y)= \Phi _{\rm nf}(x,y)\sqcup \Phi _{\rm f}(x,y)$. 
Then, 
$$
U_{\Phi _0, !}i_{x,!}(A)(y) = \coprod _{\Phi _{\rm nf}(x,y)}A, \;\; \mbox{  if  }\Phi _{\rm f}(x,y)=  \emptyset ; 
$$
$$
U_{\Phi _0, !}i_{x,!}(A)(y) = \ast \; \sqcup \coprod _{\Phi _{\rm nf}(x,y)}A,  \;\; \mbox{  if  }\Phi _{\rm f}(x,y) \neq  \emptyset .
$$
For an arrow $y\rightarrt z$, composition induces $\Phi _{\rm f}(x,y)\rightarrt \Phi _{\rm }(x,z)$ but
only $\Phi _{\rm nf}(x,y)\rightarrt \Phi _{\rm nf}(x,y)\sqcup \Phi _{\rm f}(x,y)$. The morphisms of functoriality for the diagram 
$U_{\Phi _0, !}i_{x,!}(A)$ are either the identity on $A$ or on $\ast$, or else the projection $A\rightarrt \ast$ in the case of an arrow
in $\Phi _{\rm nf}(x,y)$ which composes with $y\rightarrt z$ to give an arrow in $\Phi _{\rm f}(x,z)$.

Note that if $u\in \Phi _0$ then $\Phi _{\rm nf}(x,u)=\emptyset$ and $U_{\Phi _0, !}i_{x,!}(A)(u)=\ast$
so the above formula defines a unital diagram. One can check by hand that the explicit construction 
described above is adjoint to the
functor $\diag (\Phi /\Phi _0, \mM )\rightarrt \mM$ of evaluation at $x$, which serves to show that
the explicit construction is indeed $U_{\Phi _0, !}i_{x,!}$.

The same description holds for $U_{\Phi _0, !}i_{x,!}(B)$ and the map $U_{\Phi _0, !}i_{x,!}(f)$
is obtained by applying either $f$ or $1_{\ast}$ on the various factors.

A few things are immediate from this description:
\newline
(1) if $f$ is any  cofibration in $\mM$ then $U_{\Phi _0, !}i_{x,!}(f)$ is an injective i.e. levelwise cofibration in 
$\diag (\Phi /\Phi _0; \mM )$, indeed the $U_{\Phi _0, !}i_{x,!}(f)(y)$ are disjoint unions of copies of
$f$ and of the isomorphism $1_{\ast}$; and
\newline
(2) if $f$ is a trivial cofibration in $\mM$ then $U_{\Phi _0, !}i_{x,!}(f)$ is an injective i.e. levelwise trivial cofibration in 
$\diag (\Phi /\Phi _0; \mM )$, for the same reason. 

On the other hand, the adjunction formula says that $U_{\Phi _0, !}i_{x,!}$ is left adjoint to the restriction
$$
i_x^{\ast} : \diag (\Phi /\Phi _0, \mM )\rightarrt \mM 
$$
i.e. the evaluation at $x\in \Phi$. Hence, a morphism $g$ of diagrams in $\diag (\Phi /\Phi _0, \mM )$
satisfies right lifting with respect to $U_{\Phi _0, !}i_{x,!}(f)$, if and only if $i_x^{\ast}(g)=g(x)$ satisfies
the right lifting property with respect to $f$. 

The previous paragraph implies that $\inj (I_{\Phi /\Phi _0})$ is equal to the class of levelwise trivial fibrations, hence
$\cof (I_{\Phi /\Phi _0})$ is the class of cofibrations. Thus $I_{\Phi /\Phi _0}$ is a set of generators for the cofibrations.

If we had started with a set $J$ generating the trivial cofibrations of $\mM$, then defining $J_{\Phi /\Phi _0}$
to be the set of all $U_{\Phi _0, !}i_{x,!}(f)$ for $f\in J$, gives by the same argument 
$\inj (J_{\Phi /\Phi _0})$ equal to the class of levelwise fibrations. We claim that $\cof (J_{\Phi /\Phi _0})$ is then
equal to the class of trivial cofibrations. By property (2) above, $J_{\Phi /\Phi _0}$ and hence $\cof (J_{\Phi /\Phi _0})$
consist of levelwise weak equivalences, so they are contained in the class of trivial cofibrations.
Suppose $g: R\rightarrt S$ is a trivial cofibration. By the small object argument it can be factored as $g=ph$ where
$p\in \inj (J_{\Phi /\Phi _0})$ and $h\in \cell (J_{\Phi /\Phi _0})$. In particular $p$ is an levelwise fibration,
but it is also an levelwise weak equivalence by 3 for 2, so it is a trivial fibration hence satisfies lifting with
respect to cofibrations. As $g$ is assumed to be a cofibration, there is a lifting which shows $g$ to be a retract of $h$.
Thus $g\in \cof (J_{\Phi /\Phi _0})$. We have now shown conditions (CG1)--(CG3b) for
$I$ and $J$ with respect to the given three classes of morphisms, and we know that 
weak equivalences are closed under retracts and satisfy 3 for 2. These give a cofibrantly
generated model structure (see Proposition \ref{cofibgenrec}). It is tractable since the elements of the generating sets
have cofibrant domains. Left properness is checked levelwise. 

For the statements about left Quillen functors, note that both functors in question are left adjoints.
Furthermore, they preserve cofibrations and trivial cofibrations, indeed by (1) and (2) above the
generating cofibrations of $\diag _{\rm proj}(\Phi /\Phi _0, \mM )$ are also injective cofibrations; 
and by our construction the generating sets $I_{\Phi}$ and $J_{\Phi}$ for cofibrations and trivial cofibrations
in $\diag _{\rm proj}(\Phi , \mM )$ are mapped by $U_{\Phi _0, !}$ to the generating sets for 
cofibrations and trivial cofibrations in $\diag _{\rm proj}(\Phi /\Phi _0, \mM )$.
\end{proof}

\begin{remark}
\label{unfortunately}
Unfortunately $U_!$ is not necessarily a left Quillen functor between the injective
structures. 
\end{remark}

In the next chapter when our discussion is applied to the special case $\Delta ^o_X/X$
we can impose an additional condition \ref{cond-inj} 
on $\mM$ so that it works, or alternatively
use the Reedy structure on the category of diagrams. 

\section{Unital algebraic diagram theories}
\label{unitalgtheory}

Combine the previous discussions: 
suppose $\Phi$ is a small category, $\Phi _0\subset \Phi$ is a full subcategory,
$Q$ is a small set, we have integers $n(q)\geq 0$ for $q\in Q$, and suppose given 
functors $P_q : \epsilon (n(q))\rightarrt \Phi$ for $q\in Q$. Suppose that $\mM$ is a
tractable left proper cartesian model category with generating sets $I$ and $J$. We obtain
left Quillen functors (the leftmost varying in a family
indexed by $q\in Q$):
$$
\begin{diagram}
\diag _{\rm proj}(\epsilon (n(q)), \mM ) & 
\rightarr^{P_{q,!}} &  
\diag _{\rm proj}(\Phi , \mM )& \\
&& \downarr_{U_{\Phi _0, !}}& \\
& & \diag _{\rm proj}(\Phi /\Phi _0, \mM )&\rightarr^{1}&
\diag _{\rm inj}(\Phi /\Phi _0, \mM ).
\end{diagram}
$$
By the property of transfer of families of direct localizing systems 
along Quillen functors (Theorem \ref{transferfamily}), 
we obtain left Bousfield localizations of 
$\diag _{\rm proj}(\Phi /\Phi _0, \mM )$ and 
$\diag _{\rm inj}(\Phi /\Phi _0, \mM )$ along the 
images of the $\zeta _{\epsilon (n(q))}(f)$ for $f\in I$. 
Denote these respectively by 
$\Alg _{{\rm proj}}(\Phi /\Phi _0, P_ {\cdot};  \mM )$ 
and $\Alg _{{\rm inj}}(\Phi /\Phi _0 , P_ {\cdot};  \mM )$. 
They are tractable left proper model categories whose underlying categories are
the unital diagram category $\diag (\Phi /\Phi _0, \mM )$. 
In the projective structure, the cofibrations are generated by $U_{\Phi _0,!}(I_{\Phi})$
whereas in the injective structure the cofibrations are the levelwise cofibrations. 
The fibrant objects in the projective structure are the levelwise fibrant objects whose pullback to
each $\epsilon (n(q))$ satisfies the product condition.

In the projective case this compares with the non-unital algebraic diagram theories by a
Quillen adjunction
$$
U_{\Phi _0,! } : \Alg _{{\rm proj}}(\Phi , P_ {\cdot};  \mM )\rBotharrow \Alg _{{\rm proj}}(\Phi /\Phi _0, P_ {\cdot};  \mM ): U_{\Phi _0}^{\ast} .
$$
Indeed, the direct left Bousfield localization of the unital theory can be
seen as coming 
from the localization of the non-unital
theory given in Theorem \ref{diagprod}, by transfer along the left Quillen functor 
$U_{\Phi _0,! }$, and in the situation of Theorem \ref{transfer} we still get a 
left Quillen functor. 

\begin{remark}
\label{iffurthermore}
If furthermore we know that $U_{\Phi _0,!}$ gives a left Quillen functor on the injective
diagram structures, then this completes to a Quillen adjunction 
$$
U_{\Phi _0,! } : \Alg _{{\rm inj}}(\Phi , P_ {\cdot};  \mM )\rBotharrow \Alg _{{\rm inj}}(\Phi /\Phi _0, P_ {\cdot};  \mM ): U_{\Phi _0}^{\ast} .
$$ 
\end{remark}

\begin{lemma}
\label{productreplacementlevelwise}
Suppose $r:\pA \rightarrt \pA '$ is a trivial cofibration towards a fibrant object,
in either of the projective or injective model structures 
on $\Alg (\Phi /\Phi _0, P_ {\cdot};  \mM )$. If $\pA $ satisfies the product condition,
then $r$ is levelwise a weak equivalence, that is $r(x):\pA (x)\rightarrt \pA '(x)$ is
a weak equivalence in $\mM$ for any $x\in \Phi$. 
\end{lemma}
\begin{proof}
This follows from Corollary \ref{Rnweowe}, noting that the model 
structures on  $\Alg (\Phi /\Phi _0, P_ {\cdot};  \mM )$ are obtained
from direct left Bousfield localizing systems with $\Rr$ being the class
of objects satisfying the product condition. 
\end{proof}

\section{The generation operation}
\label{sec-gendef}

Suppose $\Phi$ is a small category, $\Phi _0\subset \Phi$ is a full subcategory,
$Q$ is a small set, we have integers $n(q)\geq 0$ for $q\in Q$, and suppose given 
functors $P_q : \epsilon (n(q))\rightarrt \Phi$ for $q\in Q$. Suppose that $\mM$ is a
tractable left proper model category with generating sets $I$ and $J$.

make the following assumption:
\newline
(INJ)---the functors $U_{\Phi _0, !}P_{q,!}$ send cofibrations (resp. trivial cofibrations)
in $\diag _{\rm inj}(\epsilon (n_q),\mM )$ to levelwise cofibrations (resp. levelwise
trivial cofibrations). 

In other words we have left Quillen functors 
$$
\diag _{\rm inj}(\epsilon (n_q),\mM )\rightarrt^{U_{\Phi _0, !}P_{q,!}}diag _{{\rm inj}}(\Phi /\Phi _0,  \mM ).
$$
In this case we can use the generating set for the new model structure on 
$\diag _{\rm inj}(\epsilon (n_q),\mM )$ made from the simpler 
cofibrations $\rho (f)$.

Suppose $\pA \in \diag (\Phi /\Phi _0, \mM )$ and $q\in Q$. Let $n:= n_q$. 
Define a trivial cofibration $\pA \rightarrt \Gen (\pA ; q)$ as follows: 
choose a factorization 
$$
P_q^{\ast}(\pA )(\xi _0) \rightarr^{e_0} E_0 \rightarr^{p} P_q^{\ast}(\pA )(\xi _1)\times \cdots \times P_q^{\ast}(\pA )(\xi _n)
$$
with $e_0$ a cofibration and $p$ a weak equivalence, in $\mM$. 
This gives a trivial cofibration $P_q^{\ast}(\pA ) \rightarr^{e}\pE$ in $\diag _{{\rm inj},\Pi}(\epsilon (n),\mM )$. Using condition (INJ) we obtain a cofibration 
\begin{equation}
\label{cofibtogen}
\pA \rightarrt \Gen (\pA ; q):= \pA  \cup ^{U_{\Phi _0, !}P_{q,!}(\pA )}U_{\Phi _0, !}P_{q,!}(E)
\end{equation}
in the injective model structure $\diag _{{\rm inj}}(\Phi /\Phi _0,  \mM )$. 

Note that $\Gen (\pA ,q)$ doesn't depend, up to equivalence, on the choice of factorization $E$.
If necessary we can include the factorization in the notation $\Gen (\pA ,q; e_0,p)$. 

The weak monadic projection from $\diag (\Phi /\Phi _0, \mM )$
to the class of objects satisfying the product condition, may be thought of as
a transfinite iteration of the operation $\pA \mapsto \Gen (\pA ,q)$ over all $q\in Q$. 



\section{Reedy structures}

In the main situation where the theory of this chapter will be applied, the underlying category $\Phi$ is a {\em Reedy category},
and we can give the category of diagrams $\Phi \rightarrt \mM$ the Reedy model category structure denoted by  $\diag _{\rm Reedy}(\Phi , \mM )$. 

Assume that $\Phi _0$ consists of objects of the bottom degree in the Reedy
structure. Then there is a corresponding model structure denoted $\diag _{\rm Reedy}(\Phi /\Phi _0, \mM )$ on the category of unital diagrams,
such that $(U_!, U^{\ast})$ remains a Quillen adjunction. 
The Reedy structures $\diag _{\rm Reedy}(\Phi , \mM )$ and 
$\diag _{\rm Reedy}(\Phi /\Phi _0, \mM )$
can again by localized by direct left Bousfield localization, to give model categories denoted 
$\Alg _{\rm Reedy}(\Phi , P_ {\cdot};  \mM )$ and $\Alg _{\rm Reedy}(\Phi /\Phi _0, P_ {\cdot};  \mM )$. These fit in between the projective and the injective
structures above. 



\chapter{Weak equivalences}
\label{weakenr1}

This chapter continues the study of weakly enriched categories using Segal's method. 
We use the model category for algebraic theories, developped in the previous chapter, to 
get model structures for Segal precategories on a fixed set of objects. This structure will be studied in detail later, to
deal with the passage from a Segal precategory to the Segal category it generates. 

Then we consider the full category of Segal precategories, with movable sets of objects, giving various definitions and notations. 
Constructing a model structure in this case is the main subject of the subsequent chapters. 

The reader will note that this division of the global argument into two pieces,
was present already in Dwyer-Kan's treatment of the model category for 
simplicial categories. They discussed the model category for simplicial categories
on a fixed set of objects in a series of papers \cite{DK1} \cite{DK2} \cite{DK3}; 
but it wasn't until some time later
with their unpublished manuscript with Hirschhorn \cite{DwyerHirschhornKan},
which subsequently became \cite{DwyerHirschhornKanSmith}, and then Bergner's paper \cite{BergnerModel} that the global case was treated. 

For the theory of weak enrichment following Segal's method, the corresponding division
and introduction of the notion of left Bousfield localization for the first part,
was suggested in Barwick's thesis \cite{BarwickThesis}. 

Assume throughout that $\mM$ is a tractable left proper cartesian model category. 
See Chapter \ref{cartmod1} for an explanation and first consequences of the cartesian
condition.

\section{The model structures on $\precat (X,\mM )$}
\label{sec-pcxm}

The Segal conditions (Section \ref{sec-segalcond}) for $\mM$-precategories can be expressed in terms of the algebraic diagram theory of the previous chapter,
which was the  motivation for introducing that notion. 

Let $\Phi := \Delta _X^o$, and let $\Phi _0 =\disc (X)$ be the discrete subcategory
on object set $X$, considered as a subcategory by letting $x\in X$ correspond to the
sequence $(x)$. An $\mM$-precategory $\pA \in \precat (X,\mM )$ is 
by definition the same thing
as a functor $\Phi \rightarrow \mM$ sending the objects of $\Phi _0$ to $\ast$, which
is to say
$$
\precat (X,\mM )= \diag (\Delta ^o_X/X, \mM ).
$$

The Segal conditions are a collection of finite product conditions
as was considered in the previous chapter. 
The set of product conditions
$Q$ consists of the full set of objects of $\Delta ^o_X$.
For $q= (x_0,\ldots , x_n)$ the integer $n(q)$ is equal to $n$, and we define a functor 
$$
P_{(x_0,\ldots , x_n)}: \epsilon (n) \rightarrt \Delta ^o_X
$$
by 
$$
P_{(x_0,\ldots , x_n)}(\xi _0) := (x_0,\ldots , x_n), 
$$
$$
P_{(x_0,\ldots , x_n)}(\xi _i) := (x_{i-1},x_i) \mbox{   for } 1\leq i \leq n.
$$
The images of the projection maps in $\epsilon (n)$ are the opposites of the inclusion maps $(x_{i-1},x_i)\hookrightarrow (x_0,\ldots , x_n)$ in $\Delta _X$.
An $\mM$-precategory is an $\mM$-enriched Segal category, if and only if it satisfies the product condition with respect to the collection of 
functors $P_{\cdot}$, indeed the two conditions are identically the same. Unitality gives the product
condition whenever $n=0$, whereas the product condition is automatically true whenever $n=1$ because the Segal maps are the identity in this case.
Thus, this condition needs only to be imposed for $n\geq 2$. 

We obtain adjoint functors
$$
P_{(x_0,\ldots , x_n),!}: \diag (\epsilon (n),\mM )\rBotharrow \diag (\Delta ^o_X,\mM ):  P_{(x_0,\ldots , x_n)}^{\ast}
$$
and the unital versions
$$
U_!P_{(x_0,\ldots , x_n),!}: \diag (\epsilon (n),\mM )\rBotharrow \precat (X,\mM ):  U_!P_{(x_0,\ldots , x_n)}^{\ast}
$$
where the right adjoint is just the pullback of a diagram $\pA :\Delta ^o_X\rightarrow \mM$
to the category $\epsilon (n)$. An $\mM$-precategory $\pA $ satisfies the Segal conditions,
if and only if $U_!P_{(x_0,\ldots , x_n)}^{\ast}(\pA )$ is a product-compatible diagram $\epsilon (n)\rightarrt \mM$, for each sequence $(x_0,\ldots , x_n)$.

A direct application of the construction of Chapter \ref{algtheor1} gives two
model structures (projective and injective) 
on $\precat (X,\mM )$ such that the fibrant objects satisfy the
Segal condition. We add a third {\em Reedy} structure since $\Delta ^o_X$ is
a Reedy category. 

The consideration of these model structures is intermediate with respect to 
our main goal of constructing global model
structures on $\precat (\mM )$: the maps in 
$\precat (X,\mM )$ are ones which induce the identity on the set of objects $X$.
Nonetheless, these easier model structures on $\precat (X,\mM )$ will be
very useful in numerous arguments later. 

The idea of introducing the intermediate  model category 
$\precat (X,\mM )$, and of expressing it as a left Bousfield localization,
is due to Barwick \cite{BarwickThesis}. 

Start by discussing the case of non-unital diagrams. 
If we fix generating sets $(I,J)$ for $\mM$, then we obtain generating sets for the projective model structure 
$\diag _{\rm proj}(\Delta ^o_X, \mM )$ and injective model structure
$\diag _{\rm inj}(\Delta ^o_X, \mM )$, recalled in Theorem \ref{injectiveprojective}. 
Notice that if $\mM$ is tractable and $I$ and $J$ consist of arrows with cofibrant domains, the explicit generators for the projective model structure 
also have cofibrant domains. Thus $\diag _{\rm proj}(\Delta ^o_X, \mM )$ will again be tractable.
For the injective model structure the construction of generating sets of Barwick and Lurie
\cite{Barwick} \cite{LurieTopos} (as discussed in Theorem \ref{lurie-thm2} above) 
was complicated, and it doesn't seem clear
whether we can choose generators with cofibrant domains. This problem can be bypassed later with the Reedy model structures where again
the generators become explicit. 

Within the projective or injective model categories of non-unital diagrams, we obtain a direct left Bousfield localizing system $(\Rr ^{\rm nu}, K^{\rm nu})$ where $\Rr$ is the 
class of non-unital diagrams satisfying the Segal conditions, and $K^{\rm nu}$ is given by the generating trivial cofibrations plus the maps of the form 
$P_{(x_0,\ldots , x_n),!}(\zeta _n(f))$ for $f$ in the generating set $I$ of cofibrations of $\mM$.
Here $\zeta _n(f)$ is the diagram $\epsilon (n)\rightarrt \mM$ considered in 
Section \ref{sec-dirlocproj}. 
This yields the direct left Bousfield localized model structures which were designated by the notation $\Alg (\ldots )$
in the previous chapter. Denote these model categories of weakly unital precategories now by
$$
\Alg _{\rm proj}(\Delta ^o_X, P_{\cdot}; \mM ),
\;\;\; \Alg _{\rm inj}(\Delta ^o_X, P_{\cdot}; \mM ).
$$
The underlying categories are both the same $\diag (\Delta ^o _X, \mM )$. 
The fibrant objects are diagrams which are fibrant in the projective or injective model structures for diagrams, and which are Segal categories.

The same will work for for $\mM$-precategories where the unitality condition is imposed. 
The model structure on $\mM$-precategories for a fixed set of objects is given by the following theorem. Note the introduction of
the notation $\pA \rightarrt \Seg (\pA )$ for a choice of fibrant replacement in either of the model categories. This notation will be used later,
but can mean that a choice is made at each usage, rather than fixing a global choice once and for all. Most constructions will be independent of the choice,
up to equivalence.

\begin{theorem}
\label{modstrucs}
Supose $\mM$ is a tractable left proper cartesian  
model category, then
there are left proper combinatorial model category structures on the unital $\mM$-precategories
$$
\precat_{\rm proj}(X; \mM ):= \Alg _{\rm proj}(\Delta ^o_X /X, P_{\cdot}; \mM ),
$$
$$
\precat_{\rm inj}(X; \mM ):= \Alg _{\rm inj}(\Delta ^o_X /X, P_{\cdot}; \mM ).
$$
The fibrant objects are unital fibrant diagrams which satisfy the Segal condition. 
The cofibrations are the projective or injective cofibrations in the unital diagram category
$\diag (\Delta ^o_X /X; \mM )$. The weak equivalences are the same in both structures.

Let $\pA \rightarrt \Seg (\pA )$ denote a trivial cofibration towards a 
fibrant replacement of $\pA $ in the projective model structure $\precat_{\rm proj}(X; \mM )$; 
this can be chosen functorially. A map $\pA \rightarrt \pB $ is a weak equivalence
if and only if $\Seg (\pA )\rightarrt \Seg (\pB )$ is a levelwise weak equivalence when 
considered as a map of diagrams
$\Delta ^o_X\rightarrt \mM$. 
\end{theorem}
\begin{proof}
Propositions \ref{unitalprojective} and \ref{unitalinjective}
give projective and injective diagram
model structures $\diag _{\rm proj}(\Delta ^o_X/X,\mM )$ and
$\diag _{\rm inj}(\Delta ^o_X/X,\mM )$ on the category of unital diagrams,
which is the same underlying category as $\precat (X,\mM )$. 

We get direct left Bousfield localizing systems $(\Rr , K_{\rm proj})$
and  $(\Rr , K_{\rm inj})$ for these model structures,
by transfering the direct localizing systems for $\epsilon (n)$-diagrams
of Theorem \ref{diagpi}, as was discussed in
Section \ref{unitalgtheory} using Theorem \ref{transferfamily}. 

In both cases $\Rr$ is the class of $\mM$-precategories which satisfy the Segal conditions;
then $K_{\rm proj}$ (resp. $K_{\rm inj}$)
is the union of the set of generators for trivial cofibrations in the projective
(resp. injective) diagram structure, plus the
morphisms of the form  $U_!P_{(x_0,\ldots , x_n),!}(\zeta _n(f))$ for $f$ in the
generating set $I$ of cofibrations of $\mM$. 
Note that the images $PU_{(x_0,\ldots , x_n),!}(g)$ of generating trivial cofibrations for the
diagram categories on $\epsilon (n(q))$ are already trivial cofibrations in
$\diag (\Delta ^o_X/X,\mM )$ so we don't need to include them again.

Now
Theorem \ref{directLBL} applies to give the required model structures. 
The characterization of weak equivalences comes from Lemma \ref{Rnweowe}.
\end{proof}

\section{Unitalization adjunctions}
\label{sec-unadj}

The projective model structure
of Theorem \ref{modstrucs} is related to the non-unital version by a
Quillen adjunction of unitalization
$$
U_{!}: \Alg _{\rm proj}(\Delta ^o_X , P_{\cdot}; \mM )\rBotharrow \precat_{\rm proj}(X; \mM ) : U^{\ast} .
$$
This follows from the application of Theorem \ref{directLBL}. 

It is useful to describe the
unitalization operation for $\Delta ^o_X/X$.

\begin{lemma}
\label{unitaldescrip}
Suppose $\pA :\Delta ^o_X\rightarrow \mM $ is a functor. Then $U_!\pA $ has the
following explicit description: 
\newline
---if $x_0=\ldots = x_n=x$ is a constant sequence then the full degeneracy gives a map $\pA (x)\rightarrt \pA (x,\ldots ,x)$ and 
$$
(U_!\pA )(x,\ldots , x) = \pA (x,\ldots , x)\cup ^{\pA (x)}\ast ;
$$
---otherwise, if $x_0,\ldots , x_n$ is not a constant sequence then
$$
(U_!\pA )(x_0,\ldots , x_n) = \pA (x_0,\ldots , x_n).
$$
\end{lemma}
\begin{proof}
The object explicitly constructed in this way is again a functor $\Delta ^o_X\rightarrt \mM$
because if $(x_0,\ldots , x_n)$ is a constant sequence and 
$(y_0,\ldots , y_k)\rightarrt (x_0,\ldots , x_n)$ is any map in $\Delta _X$
then $y_{\cdot}$ must also be a constant sequence, so the map 
$\pA (x_0,\ldots , x_n)\rightarrt \pA (y_0,\ldots , y_k)$ passes to the quotients.
The resulting functor satisfies the required left adjunction property  
with respect to $U^{\ast}$ on the right, so it must be $U_!\pA $. 
\end{proof}

In the injective case, $U_!$ will not in general be a left Quillen functor; we need
to impose an additional hypothesis. 

\begin{condition}[INJ]
\label{cond-inj}
Suppose 
\begin{diagram}
A & \rightarr & X& \rightarr & A \\
\downarr & & \downarr & & \downarr \\
B & \rightarr & Y& \rightarr & B
\end{diagram}
is a diagram in $\mM$ such that the vertical arrows are cofibrations
(resp. trivial cofibrations)
and the horizontal compositions are the identity. Then 
$X\cup ^A \ast \rightarrow Y\cup ^B\ast$ is a cofibration (resp. trivial cofibration). 
\end{condition}

The following observations allow us to use the injective model categories in
many cases. In fact, in these cases the injective structure also coincides with the
Reedy structure, however it seems comforting to be able to use the injective structure
which is conceptually simpler, instead. 

\begin{lemma}
\label{condinjpresh}
Suppose  $\mM$ is a presheaf category, is left proper, and the class of cofibrations is the
class of monomorphisms of presheaves. Then Condition \ref{cond-inj} holds.
\end{lemma}
\begin{proof}
Note first that given a diagram of sets as in Condition \ref{cond-inj} where the
vertical maps are injections, then the map $X\cup ^A \ast \rightarrow Y\cup ^B\ast$
is an injection of sets. Indeed if $x\in X$ maps to an element of $B$ then the image
of $x$ by the projection to $A$, maps to the same element of $B$. By injectivity of
the map $X\rightarrt Y$it follows that $x\in A$, which shows injectivity of the
map on quotients. 

Now if $\mM$ is a presheaf category and the cofibrations are the monomorphisms,
applying the previous paragraph levelwise we obtain the desired result for cofibrations.
Suppose given a diagram whose vertical arrows are trivial cofibrations. 
The split injections $A\rightarrt X$ and 
$B\rightarrt Y$ are monomorphisms,  hence cofibrations, so Corollary \ref{pushoutequivcor}
applies to conclude that the pushout map is a weak equivalence, hence it is a trivial
cofibration.  
\end{proof}

\begin{lemma}
\label{injectiveunitalization}
If $\mM$ satisfies Condition \ref{cond-inj} then the unitalization functors
give a Quillen adjunction between injective diagram structures
$$
U_{!}: \Alg _{\rm inj}(\Delta ^o_X , P_{\cdot}; \mM )\rBotharrow \precat_{\rm inj}(X; \mM ) : U^{\ast}
$$
where $U^{\ast}$ is just the identity inclusion of unital precategories in all precategories. 
\end{lemma}
\begin{proof}
We first show that unitalization is a Quillen adjunction between levelwise
injective diagram categories
$$
U_!:\diag _{\rm inj}(\Delta ^o_X,\mM )
\rBotharrow \diag _{\rm inj}(\Delta ^o_X/X,\mM ):U^{\ast} .
$$
Suppose $\pA \rightarrt \pB $ is a levelwise cofibration (resp. trivial cofibration)
of diagrams, the claim is that
$U_!\pA  \rightarrt U_!\pB $ is a levelwise cofibration (resp. trivial cofibration). 
In view of the description of Lemma \ref{unitaldescrip} it suffices to
look at the values over a constant sequence $(x,\ldots , x)$. We have a diagram
$$
\begin{diagram}
\pA (x) & \rightarr & \pA (x,\ldots , x)& \rightarr & \pA (x) \\
\downarr & & \downarr & & \downarr \\
\pB (x) & \rightarr & \pB (x,\ldots , x)& \rightarr & \pB (x)
\end{diagram}
$$
where the vertical arrows are cofibrations (resp. trivial cofibrations), 
and the second horizontal arrows are,
say, the projections corresponding to the first object of the sequence. 
Condition \ref{cond-inj} now says exactly that 
$$
\begin{diagram}
U_!\pA (x,\ldots , x) = \pA (x,\ldots , x)\cup ^{\pA (x)} \ast \\
\downarr \\
U_!\pB (x,\ldots , x) = \pB (x,\ldots , x)\cup ^{\pB (x)} \ast
\end{diagram}  
$$
is a cofibration (resp. trivial cofibration). This shows the Quillen adjunction for
the diagram categories. 

Now, the categories in question for the lemma are obtained by direct left Bousfield
localization using the sets of generators plus the morphisms of the form $P_{q,!}(\zeta _{n(q)}(f))$ for $q\in Q$ and $f$ in a generating set for cofibrations of $\mM$.
The left Quillen functor passes to a left Quillen functor between localizations
by Theorem \ref{transfer}.  
\end{proof}

\section{The Reedy structure}
\label{sec-reedyXstruc}

The category $\Delta ^o_X$ is a Reedy category, using the 
subcategories of injective and surjective maps of finite ordered sets, as 
direct and inverse subcategories, and the
length function $(x_0,\ldots , x_n)\mapsto n$. 
This leads to a {\em Reedy model structure} on the category of diagrams,
using levelwise weak equivalences,
denoted $\diag _{\rm Reedy}(\Delta ^o_X,\mM )$ (Proposition \ref{reedystructure}).

The notion of Reedy structure on diagram categories is a slightly more technical area
of the theory of model categories, but these structures are very natural and turn
out to be the best ones for our theory of precategories. In many useful examples 
the Reedy structure coincides with the injective structure. This is the case
for example if $\mM=\presh (\Psi )$ is a presheaf category and the cofibrations of $\mM$
are the monomorphisms of presheaves, see Proposition \ref{reedyinjective}.

There is a unital version of the Reedy model structure. 

\begin{proposition}
\label{unitalreedy}
For any tractable left proper model category $\mM$, 
the unital diagram category has a tractable left proper model structure denoted 
$\diag _{\rm Reedy}(\Delta ^o_X/X,\mM )$, related to the non-unital Reedy structure
by a Quillen adjunction
$$
U_!:\diag _{\rm Reedy}(\Delta ^o_X,\mM )
\rBotharrow \diag _{\rm Reedy}(\Delta ^o_X/X,\mM ):U^{\ast}
$$
using the levelwise weak equivalences. The fibrations (resp. cofibrations) of $\diag _{\rm Reedy}(\Delta ^o_X/X,\mM )$ are exactly the maps $f$ such that $U^{\ast}(f)$ is a Reedy fibration (resp. cofibration) in
$\diag _{\rm Reedy}(\Delta ^o_X,\mM )$, in particular they have the same description
in terms of latching and matching objects. 
The generating sets for the unital Reedy structure are obtained by applying $U_!$ to the
generating sets for the regular Reedy diagram structure. 

The Reedy structure lies in between the projective
and injective structures with a diagram of left Quillen functors
$$
\begin{diagram}
\diag _{\rm proj}(\Delta ^o_X,\mM )&\rightarr & \diag _{\rm Reedy}(\Delta ^o_X,\mM )
& \rightarr & \diag _{\rm inj}(\Delta ^o_X,\mM ) \\
\downarr ^{U_!} & & \downarr_{U_!} \\
\diag _{\rm proj}(\Delta ^o_X/X,\mM )&\rightarr & \diag _{\rm Reedy}(\Delta ^o_X/X,\mM )
& \rightarr & \diag _{\rm inj}(\Delta ^o_X/X,\mM )
\end{diagram}
$$
where the horizontal rows are identity functors.  If $\mM$ satisfies Condition 
\ref{cond-inj} then this can be completed by putting in the rightmost vertical arrow. 
\end{proposition}
\begin{proof}
The first thing to show is that if $\pA \rightarrt ^f \pB $ is a Reedy cofibration
of diagrams on $\Delta ^o_X$, then $U^{\ast}U_!\pA  \rightarrt ^{U^{\ast}U_!f}
U^{\ast}U_!\pB $ is again a Reedy cofibration. The latching objects over a sequence
$(x_0,\ldots , x_n)\in \Delta ^o_X$ involve maps which correspond to
surjections $(x_0,\ldots , x_n)\rightarrt^{\sigma} (y_0,\ldots , y_k)$ of sequences, with
$\sigma : [n]\rSurjarrow  \sqb{k}$ with
$y_{\sigma (i)} = x_i$. In particular, any $(y_0,\ldots , y_k)$ involved in the latching map at a non-constant sequence $(x_0,\ldots , x_n)$, is also not constant. Over these sequences
the relative latching map for $U^{\ast}U_!f$ is the same as for $f$ so the 
Reedy condition is preserved. We may therefore concentrate on the case of a constant
sequence. Let $x^n:= (x_0,\ldots , x_n)$ with $x_i=x$. Let $\pA ':= U^{\ast}U_!\pA $,
$\pB ':= U^{\ast}U_!\pB $, and $f':= U^{\ast}U_!f$. Thus
$$
\pA '(x^n)=\pA (x^n)\cup ^{\pA (x)}\ast , \;\;\; \pB '(x^n)=\pB (x^n)\cup ^{\pB (x)}\ast .
$$
The latching object for $\pA $ is expressed as the pushout
$$
\begin{diagram}
\coprod _{0< i < j \leq n}\pA (x^{n-2}) & \rightarr & \coprod _{0< i \leq n}\pA (x^{n-1})\\
\downarr & & \downarr \\
\coprod _{0< i \leq n}\pA (x^{n-1}) &\rightarr & \latch (\pA ,x^n) 
\end{diagram}
$$
and similarly for $\latch (\pB ,x^n)$.  The coproducts may be considered as being taken over
contractions of adjacent pairs in the sequence $(x_0,\ldots , x_n)$ with all $x_i=x$.

We get a map from $\pA (x)$ into all elements of the
above cocartesian diagram, and then the same diagram for $\pA '$ gives 
$$
\latch (\pA ',x^n) = \latch (\pA ,x^n) \cup ^{\pA (x)}\ast .
$$
Similarly 
$$
\latch (\pB ',x^n) = \latch (\pB ,x^n) \cup ^{\pB (x)}\ast .
$$
Recall that 
$$
{\latch}(f,x^n):= {\latch}(\pB ,x^n)\cup ^{{\latch}(\pA ,x^n)}\pA (x^n)\rightarrt \pB (x^n).
$$
Thus ${\latch}(f',x^n)$ is obtained by contracting out $\pB (x)$ in the domain and
range, in particular
$$
{\latch}(\pB ',x^n)\cup ^{{\latch}(\pA ',x^n)}\pA '(x^n) = \ast \cup ^{\pB (x)}{\latch}(\pB ,x^n)\cup ^{{\latch}(\pA ,x^n)}\pA (x^n)
$$
as can be seen by using the universal property of coproducts. In the diagram
$$
\begin{diagram}
\pB (x)&\rightarr & {\latch}(\pB ,x^n)\cup ^{{\latch}(\pA ,x^n)}\pA (x^n)& \rightarr & \pB (x^n) \\
\downarr & & \downarr & & \downarr \\
\ast & \rightarr & {\latch}(\pB ',x^n)\cup ^{{\latch}(\pA ',x^n)}\pA '(x^n)& \rightarr & \pB '(x^n)
\end{diagram}
$$
the left square is cocartesian, and the outer square is cocartesian,
therefore the right square is cocartesian. 
The upper right map is assumed to be a cofibration, so the bottom right
map is also a cofibration. 

This shows that $U^{\ast}U_!f$ is a Reedy cofibration whenever $f$ is. 
The same argument shows that if $f$ is a Reedy trivial
cofibration then  $U^{\ast}U_!f$ is a Reedy trivial
cofibration.

To construct the model structure on the unital diagram
category, define the weak equivalences (resp. cofibrations, fibrations)
to be those arrows $f$ in 
$\diag _{\rm Reedy}(\Delta ^o_X/X,\mM )$ 
such that $U^{\ast}f$ is a weak equivalence (resp. cofibration,
fibration). In view of this definition, the same holds for the intersection classes of
trivial fibrations and trivial cofibrations. Since $U^{\ast}$ takes compositions to 
compositions and retracts to retracts, the classes are all closed under retracts, and
weak equivalences (which are levelwise weak equivalences) satisfy 3 for 2. 

Suppose $I_R$ and $J_R$ are generating sets for the Reedy structure on non-unital diagrams
$\diag _{\rm Reedy}(\Delta ^o_X,\mM )$, and let $I'_R:= U_!I_R$ and $J'_R:=U_!J_R$.
Apply Proposition \ref{cofibgenrec}. Note that elements of $I'_R$ are cofibrations (CG2a)
and elements of $J'_R$ are trivial cofibrations (CG3a),
by the above arguments. All objects are small (CG1). 

To complete, we need to show that
a map $g:U\rightarrt V$ is a fibration (resp. trivial fibration) if and only if
it satisfies lifting with respect to $I'_R$ (resp. $J'_R$). But this is an immediate
consequence of the definitions and the adjunction property. For example,
$g$ is a fibration $\Leftrightarrow$ $U^{\ast}g$ is a fibration $\Leftrightarrow$
$U^{\ast}g$ satisfies right lifting with respect to $J_R$ $\Leftrightarrow$
$g$ satisfies right lifting with respect to $U_!J_R=J'_R$. The proof for trivial 
fibrations is the same. 

Applying Proposition \ref{cofibgenrec} we get a cofibrantly generated model category,
from which it also follows that 
$\cof (I'_R)$ is the class of cofibrations, and $\cof (J'_R)$
is the class of trivial cofibrations. 

The Reedy diagram structure is tractable, so we may assume that the domains
of elements of $I_R$ and $J_R$ are cofibrant, from which it follows that the
same is true of $I'_R$ and $J'_R$. Reedy cofibrations are injective cofibrations.
But Reedy-unital cofibrations (i.e. cofibrations in 
$\diag _{\rm Reedy}(\Delta ^o_X/X,\mM )$) are maps whose $U^{\ast}$ are Reedy diagram
cofibrations, hence they are levelwise cofibrations. Also weak equivalences are 
defined levelwise; so left properness may be verified levelwise. 

For the second paragraph of the theorem,
the upper sequence of left Quillen functors 
$$
\diag _{\rm proj}(\Delta ^o_X,\mM )\rightarr  \diag _{\rm Reedy}(\Delta ^o_X,\mM )
\rightarr  \diag _{\rm inj}(\Delta ^o_X,\mM ) 
$$
is classical. From this the lower sequence follows too. To go from the projective
to the Reedy structure, the generators of 
the projective unital diagram structure are obtained by applying $U_!$ to generators
of the projective non-unital structure, but these go to Reedy cofibrations (resp.
trivial cofibrations) in the non-unital Reedy structure,
which in turn get sent to the same in the unital Reedy structure.
To go from the Reedy to the injective structure, note that the unital Reedy
cofibrations are also non-unital Reedy cofibrations by definition, and these are
levelwise cofibrations by the upper sequence. The same is true for trivial
cofibrations because the weak equivalences are the same in all cases. 

If $\mM$ satisfies Condition \ref{cond-inj} then $U_!$ is a left Quillen functor
by Lemma \ref{injectiveunitalization} and it clearly makes the diagram commute. 
This completes the proof. 
\end{proof}

\begin{theorem}
\label{reedyprecat}
Suppose $\mM$ is a tractable left proper cartesian model category. 
The model structure of $\diag _{\rm Reedy}(\Delta ^o_X/X,\mM )$
admits a direct left Bousfield localization denoted $\precat _{\rm Reedy}(X,\mM )$
by a system $(\Rr , K_{\rm Reedy})$ where $\Rr$ is the
class $\mM$-precategories satisfying the Segal conditions, 
and $K_{\rm Reedy}$ is the generating set for trivial
Reedy cofibrations, plus the set of cofibrations of the
form $U_!P_{(x_0,\ldots , x_n),!}(\rho _n(f))$ for $f$ a generating cofibration of $\mM$. 
The cofibrations of $\precat _{\rm Reedy}(X,\mM )$ are the
maps which are Reedy cofibrations in the diagram category 
$\diag _{\rm Reedy}(\Delta ^o_X,\mM )$. The fibrant 
objects of $\precat _{\rm Reedy}(X,\mM )$ are the Reedy-fibrant diagrams
which satisfy the Segal conditions. The Reedy structure lies between the 
projective and injective structures with the identity functors being left Quillen
$$
\precat _{\rm proj}(X,\mM )\rightarrt\precat _{\rm Reedy}(X,\mM )
\rightarrt\precat _{\rm inj}(X,\mM ).
$$
The weak equivalences are the same for the three structures. 
\end{theorem}
\begin{proof}
This is analogous to the proof of Theorem \ref{modstrucs}
using 
Theorem \ref{directLBL} and Theorem \ref{transferfamily}.
The only difference worth pointing out here is that we can use
the morphisms $\rho _n(f)$ defined just before Proposition \ref{injepsilongen},
in place of the $\zeta _n(f)$ from Section \ref{sec-dirlocproj}.

We claim that 
$$
U_!P_{(x_0,\ldots , x_n),!}(\varpi _n(f))\rightarrt^{U_!P_{(x_0,\ldots , x_n),!}(\rho _n(f))}
U_!P_{(x_0,\ldots , x_n),!}(\xi _{0,!}(V))
$$
is a Reedy cofibration
if $U\rightarrt ^fV$ is a cofibration in $\mM$. 
Recall from Equation \eqref{varpidef}
that 
$$
\varpi _n(f):= \xi _{0,!}(U) \cup ^{\coprod _i \xi _{i,!}(U)}\coprod _i \xi _{i,!}(V) = (U, V, \ldots , V; f,\ldots , f).
$$
It will be more convenient to consider an arbitrary $\pA \in \precat (X,\mM )$
with a map $U_!P_{(x_0,\ldots , x_n),!}(\varpi _n(f))\rightarrt^a \pA $
and let 
$$
\pB := \pA \cup ^{U_!P_{(x_0,\ldots , x_n),!}(\varpi _n(f))}
U_!P_{(x_0,\ldots , x_n),!}(\xi _{0,!}(V))
$$
be the pushout along 
$U_!P_{(x_0,\ldots , x_n),!}(\rho _n(f))$. We would like to show that $\pA \rightarrt \pB $
is a cofibration. 

The map $a$ corresponds to a collection of maps $a_i:V\rightarrt \pA (x_{i-1},x_i)$
and $a_0:U\rightarrt \pA (x_0,\ldots , x_n)$ such that the diagrams
$$
\begin{diagram}
U & \rightarr & \pA (x_0,\ldots , x_n)\\
\downarr && \downarr \\
V& \rightarr & \pA (x_{i-1},x_i)
\end{diagram}
$$
commute. 
We can now describe $\pB $ explicitly. For any sequence $(y_0,\ldots , y_k)$,
$$
\pB (y_0,\ldots , y_k)= \pA (y_0,\ldots , y_k)\cup ^{ \coprod _{y_{\cdot}\rightarrow x_{\cdot}} U}
\left( \coprod _{y_{\cdot}\rightarrow x_{\cdot}} V \right) ,
$$
where the coproduct is taken over all arrows $y_{\cdot}\rightarrow x_{\cdot}$
which don't factor through one of the adjacent pairs 
$(x_{i-1},x_i)\rightarrow (x_0,\ldots , x_n)$.
The left map in the coproduct uses $a_0$, whereas
the functoriality maps to define the functor $\pB :\Delta ^o_X\rightarrt \mM$
are made using the $a_i$ for $i=1,\ldots , n$. 
Suppose $(z_0,\ldots , z_l)\rightarrt (y_0,\ldots , y_k)$ is a map. For any $h:y_{\cdot}\rightarrow x_{\cdot}$
which doesn't factor through an adjacent pair, if the composition 
$z_{\cdot}\rightarrt x_{\cdot}$ does factor through an adjacent pair $(x_{i-1},x_i)$
then 
we use the map $a_i: V\rightarrt \pA (x_{i-1},x_i)\rightarrow \pA (z_0,\ldots , x_l)$
on the component $V$ of $\pB (y_0,\ldots , y_k)$ corresponding to $h$. 
The other maps of functoriality are straightforward. 

From this description we already see easily that $\pA \rightarrow \pB $ is an injective i.e.
levelwise cofibration. 
To show that it is a Reedy cofibration, consider $(y_0,\ldots , y_k)$.
The latching objects are colimits over the category of surjections $(y_0,\ldots , y_k)\rSurjarrow^{\sigma} (z_0,\ldots , z_l)$ with $l<k$.
For any $\sigma$, the arrows $z_{\cdot}\rightarrow x_{\cdot}$ which don't factor
through an adjacent pair, form a subset of the similar arrows  $y_{\cdot}\rightarrow x_{\cdot}$. 
We can write
$$
\latch (\pB ,y_{\cdot}) = \latch (\pA ,y_{\cdot})\cup ^{ \coprod _{y_{\cdot}\rightarrow x_{\cdot},{\rm fact}} U}
\left( \coprod _{y_{\cdot}\rightarrow x_{\cdot},{\rm fact}} V \right)
$$
where the coproduct is now taken over $y_{\cdot}\rightarrow x_{\cdot}$
which factor through some surjection $y_{\cdot}\rightarrow z_{\cdot}$ in the latching
category. 
The relative latching object for the map $\pA \rightarrow \pB $ is then 
$$
\pA (y_{\cdot})\cup ^{ \coprod _{y_{\cdot}\rightarrow x_{\cdot},{\rm fact}} U}
\left( \coprod _{y_{\cdot}\rightarrow x_{\cdot},{\rm fact}} V \right)
\rightarrow 
\pA (y_{\cdot})\cup ^{ \coprod _{y_{\cdot}\rightarrow x_{\cdot}} U}
\left( \coprod _{y_{\cdot}\rightarrow x_{\cdot}} V \right) .
$$
Thus, the relative latching map consists just of adding on some additional pushouts
along $U\rightarrt V$. It is cofibrant, which shows that $\pA \rightarrow \pB $ is a Reedy
cofibration. 

Now the same proof as for Theorem \ref{modstrucs} shows that 
we can get a direct localizing system for the class $\Rr$ of
$\mM$-enriched Segal categories by starting with the generating set for
trivial cofibrations in $\diag _{\rm Reedy}(\Delta ^o_X/X,\mM )$
and adding the $U_!P_{(x_0,\ldots , x_n),!}(\rho _n(f))$ to get $K_{\rm Reedy}$.

We could also have added the
$U_!P_{(x_0,\ldots , x_n),!}(\zeta _n(f))$, or both. These all give the same
left Bousfield localization (Proposition  \ref{lblindep}). 
Using the $U_!P_{(x_0,\ldots , x_n),!}(\zeta _n(f))$, Theorem \ref{transfer} gives the
left Quillen functors from the projective to the Reedy and then to the injective
structure. 

We can remark that the above proof also shows that we could have used the 
$U_!P_{(x_0,\ldots , x_n),!}(\rho _n(f))$ to form a direct localizing system 
for the injective structure. 
\end{proof}

\section{Some remarks}
\label{sec-someremarks}

The model structures on $\precat (X; \mM )$ are the main intermediate structure to be considered \cite{BarwickThesis} before getting the model structures for
moving object sets.    
In the remainder of the book, when we speak of weak equivalences or the injective
or projective model structures on $\precat (X; \mM )$, we mean the 
structures of Theorem \ref{modstrucs}. In case it is necessary, we use the
terminology {\em levelwise weak equivalence} to speak of
weak equivalences in the unital diagram category $\diag (\Delta ^o_X/X, \mM )$.

\begin{lemma}
\label{localwesegal}
If $\pA ,\pB \in \precat (X; \mM )$ satisfy the Segal conditions, then 
a map $f:\pA \rightarrt \pB $ in $\precat (X; \mM )$ is a weak equivalence 
if and only if it is a levelwise weak equivalence. 

If $\pB $ satisfies the
Segal conditions and $\pA \rightarrt \pB $ is a trivial fibration in
either the projective or injective model structures on 
$\precat (X; \mM )$, then $\pA $ also satisfies the Segal conditions and $f$ is
a levelwise weak equivalence.
\end{lemma}
\begin{proof}
Consider the square
$$
\begin{diagram}
\pA  & \rightarr & \Seg (\pA )\\
\downarr ^f& & \downarr _{\Seg (f)}\\
\pB  & \rightarr & \Seg (\pB ).
\end{diagram}
$$
If $\pA $ and $\pB $ both satisfy the Segal conditions, then the horizontal arrows are
levelwise weak equivalences. By the definition of the model structure
of Theorem \ref{modstrucs} the map $f$ is 
a weak equivalence if and only if $\Seg (f)$ is a levelwise weak equivalence,
and by 3 for 2 this is equivalent to requiring that $f$ be a levelwise weak equivalence. 

Suppose $\pB $ satisfies the Segal conditions, and $f:\pA \rightarrt \pB $ is a trivial fibration
in the projective model structure on $\precat (X,\mM )$. This means that $f$ satisfies
right lifting with respect to cofibrations, but the cofibrations are the
same as those of the unital diagram category. Thus, $f$ is a trivial fibration
in $\diag (\Delta ^o_X/X,\mM )$, in particular it is a levelwise weak equivalence. 
Under our hypotheses on $\mM$,  Lemma \ref{prodconsequence} says that direct product 
is invariant under weak equivalences in $\mM$, from which it follows that $\pA $ satisfies
the Segal conditions too. 
\end{proof}

In the diagram category $\diag (\Delta ^o _X, \mM )$ limits and colimits are computed
levelwise: if
$\{ \pA _i\} _{i\in \alpha}$ is a diagram of objects $\pA _i\in \diag (\Delta ^o _X, \mM )$ then 
$$
(\colim _{i\in \alpha}\pA _i)(x_0,\ldots , x_n) = \colim _{i\in \alpha}(\pA _i(x_0,\ldots , x_n)),
$$
and
$$
(\lim _{i\in \alpha}\pA _i)(x_0,\ldots , x_n) = \lim _{i\in \alpha}(\pA _i(x_0,\ldots , x_n)).
$$
The right adjoint functor $U^{\ast}$ being the identity means that the same holds true for all limits in $\precat (X; \mM )$:
if $\{ \pA _i\} _{i\in \alpha}$ is a diagram of objects $\pA _i\in \precat (X;\mM )$ then 
$$
(\lim _{i\in \alpha}\pA _i)(x_0,\ldots , x_n) = \lim _{i\in \alpha}(\pA _i(x_0,\ldots , x_n)).
$$
On the other hand, for colimits in $\precat (X;\mM )$ we have to reapply the functor $U_!$.
Write this with superscripts $\precat$ or $\diag$ to indicate in which category the colimit is taken:
$$
\colim ^{\precat}_{i\in \alpha}\pA _i = U_!( \colim ^{\diag}_{i\in \alpha}\pA _i ).
$$
In the special case where $\colim ^{\diag}_{i\in \alpha}\pA _i$ is already in $\precat (X,\mM )$ then the functor $U_!$ acts as
the identity. The condition that $\colim ^{\diag}_{i\in \alpha}\pA _i$ already be in $\precat (X,\mM )$, says that the value on
sequences of length zero $(x_0)$ is $\ast$. As the $\colim ^{\diag}$ is calculated levelwise, it is sufficient to require
that $\colim _{\alpha}\ast = \ast$ which holds when $\alpha$ is a connected category by Lemma \ref{connstar}. This condition
is necessary and sufficient unless $X$ is the empty set. 

\begin{lemma}
\label{colimexpr}
Suppose $\alpha$ is a connected category and $\{ \pA _i\} _{i\in \alpha}$ is a diagram of objects $\pA _i\in \precat (X;\mM )$.
Then the colimit in $\precat (X;\mM )$ may be calculated levelwise:
$$
(\colim ^{\precat}_{i\in \alpha}\pA _i)(x_0,\ldots , x_n) = \colim _{i\in \alpha}(\pA _i(x_0,\ldots , x_n)).
$$
\end{lemma}
\begin{proof}
Above (also, this is the same as the corresponding statement in Theorem \ref{unitaltheorylocpres} of
the previous chapter). 
\end{proof}

\begin{corollary}
\label{constobjpres}
Suppose $\mM$ is locally $\kappa$-presentable, and $X$ is a set with $<\kappa$ elements. 
An object $\pA \in \precat (X; \mM )$ is $\kappa$-presentable if and only if
each $\pA (x_0,\ldots , x_p)$ is a $\kappa$-presentable object of $\mM$.
The category $\precat (X; \mM )$ is locally $\kappa$-presentable. 
\end{corollary}
\begin{proof}
Apply Theorem \ref{unitaltheorylocpres}.
\end{proof}

For colimits over disconnected categories, in particular for disjoint sums, one has to apply the functor $U_!$. 
This will not be used very often, because when we work in the category $\precat (\mM )$ with variable set of objects
which will be defined next, disjoint sums lead to disjoint unions of sets of objects, returning a reasonable description (Section \ref{sec-limits}).

\section{Global weak equivalences}
\label{sec-globalwe}

We turn now from consideration of the various model categories for $\mM$-precategories
with a fixed set of objects, to the global category $\precat (\mM )$
of $\mM$-precategories with arbitrary variable set of objects. We'll get
back to a closer analysis of the generation operation $\Seg$ on $\precat (X,\mM )$
in Chapter \ref{genrel1} on the calculus of generators and relations, but
first we consider weak equivalences in the global case, and in the next Chapter
\ref{cofib1}, various types of cofibrations in the global case. 

In this section we formalize Tamsamani's induction step for defining equivalences of $n$-categories, as applied in the $\mM$-enriched case.
For classical Segal categories which is the $\mK$-enriched case, this notion of equivalence goes back at least to Dwyer and Kan and is
known as ``Dwyer-Kan equivalence'' \cite{Rezk} \cite{BergnerSurvey}. 
The general case requires a functor $\tau _{\leq 0}:\mM \rightarrt \Sets$ in order to be able to talk about the essential surjectivity
condition; this issue was explored by Pelissier, but we have modified his situation a little bit by introducing sets $X$ of objects external to $\mM$. Our discussion here
specializes the more general discussion in Section \ref{sec-compequivtrunc}. 

Define the {\em $0$-th truncation functor} $\tau _{\leq 0}:\mM \rightarrt \Sets$ by
$$
\tau _{\leq 0}(\pA ):= \Hom _{\Ho (\mM )}(\ast , \pA ).
$$
Then define a {\em truncation functor} $\tau _{\leq 1}: \precat (\mM ) \rightarrt \Cat$ as follows: 
for $\pA \in \mM$ with $\Ob (\pA )= X$, choose a weak equivalence in $\precat (X,\mM )$
from $\pA $ to an $\mM$-enriched Segal category $\Seg (\pA )$ (i.e. a precategory which satisfies the Segal conditions). Consider the functor 
from $\Delta _X$ to $\Sets$ defined by
$$
\tau _{\leq 1}(\pA )(x_0,\ldots , x_n):= \tau _{\leq 0}(\Seg (\pA )(x_0,\ldots , x_n)) 
$$
$$
= \Hom _{\Ho (\mM )}(\ast , \Seg (\pA )(x_0,\ldots , x_n)).
$$
The Segal maps for this functor are isomorphisms (proven in the next lemma), 
so it defines a $1$-category whose object set is $X$ and such that the set of morphisms from $x_0$ to $x_1$
is $\tau _{\leq 1}(\pA )(x_0,x_1) = \tau _{\leq 0}(\Seg (\pA )(x_0,x_1)$. 

One good choice for $\Seg (\pA )$ would be to take a fibrant replacement, however other smaller choices can be useful too (see Chapter \ref{genrel1} below). 

\begin{lemma}
\label{truncationproperties}
The truncation $\tau _{\leq 0}: \mM \rightarrt \Sets$ sends weak equivalences to isomorphisms and products to products. 
The $\Sets$-precategory $\tau _{\leq 1}(\pA )$ defined above satisfies the Segal conditions so it corresponds to a $1$-category,
and this category is independent of the choice of $\Seg (\pA )$.
It gives a functor $\tau _{\leq 1}: \precat (\mM ) \rightarrt \Cat$, and for a fixed object set $X$ the functor 
takes weak equivalences in $\precat (X, \mM )$ to isomorphisms. 
\end{lemma}
\begin{proof}
The $\tau _{\leq 0}$ factors through the projection to the homotopy category
so it clearly preserves weak equivalences. It preserves products since $\mM$ is
assumed to be cartesian. In particular the Segal condition for $\Seg (\pA )$
yields the Segal condition for the $\Delta ^o_X$-set $\tau _{\leq 1}(\pA )$ 
to be the nerve of a category, independent of the choice of $\Seg (\pA )$ because
$\tau _{\leq 0}$ sends weak equivalences to isomorphisms. We get a functor 
$\tau _{\leq 1}$. If $\pA \rightarrt \pB $ is a weak equivalence 
in $\precat (X,\mM )$ then $\Seg (\pA )\rightarrt \Seg (\pB )$
is a levelwise weak equivalence (see Theorem \ref{modstrucs})
so the resulting $\tau _{\leq 1}(\pA )\rightarrow \tau _{\leq 1}(\pB )$ is fully faithful,
but in the present case it is also the identity on the set $X$ of objects so it 
is an isomorphism. 
\end{proof}

\begin{lemma}
\label{truncationproducts}
If $\pA$ and $\pB$ are $\mM$-precategories satisfying the Segal condition 
then there is a natural isomorphism
$$
\tau _{\leq 1}(\pA \times \pB ) \cong \tau _{\leq 1}(\pA )
\times \tau _{\leq 1}(\pB ).
$$
\end{lemma}
\begin{proof}
As $\tau _{\leq 0}$ is compatible with products this follows from the definition of
$\tau _{\leq 1}$. 
\end{proof}

The truncation operation is crucial to Tamsamani's construction of $n$-categories \cite{Tamsamani}, since it allows us to define equivalences between
enriched category objects. In the case of $n$-categories, $\mM$ corresponds to the model category for $n-1$-categories and $\Ho (\mM )$
can be identified with the homotopy category of Segal $n-1$-categories. For Segal $n$-categories, the truncation as we have defined it here coincides with
Tamsamani's truncation operation which was constructed by induction on $n$. The general iteration procedure
requires a general definition of truncation starting from a model category $\mM$. This more general case was considered by 
Pelissier \cite[Definition 1.4.1]{Pelissier}.

We now proceed with the definition of global weak equivalence, combining two conditions: ``fully faithful'' is defined by asking for
weak equivalences in $\mM$ between morphism objects, and ``essentially surjective'' is defined using the above truncation operation.

A map $f:(X,\pA ) \rightarrt (Y,\pB )$ consists of a map $\Ob (f):X\rightarrt Y$ which we sometimes denote by $f$ for short, togegher with maps 
$f:\pA (x_0,\ldots , x_p)\rightarrt \pB (f(x_0),\ldots , f(x_p))$ in $\mM$
for each $(x_0,\ldots , x_p)\in \Delta _{X}$, compatible with restrictions along maps in $\Delta _X$.

Suppose that $(X,\pA )$ and $(Y,\pB )$ are $\mM$-precategories satisfying the Segal condition. A map $f: (X,\pA )\rightarrt (Y,\pB )$
is {\em fully faithful} if, for every sequence $(x_0,\ldots , x_n)$ of objects in $X$ the map
$$
f(x_0,\ldots , x_n): \pA (x_0,\ldots , x_n) \rightarrt \pB (f(x_0),\ldots , f(x_n))
$$
is a weak equivalence in $\mM$.

If $(X,\pA )$ and $(Y,\pB )$ are only $\mM$-precategories,
let $\pA \rightarrt \Seg (\pA )$ and $\pB \rightarrt \Seg (\pB )$ denote injective
trivial cofibrations towards objects satisfying the Segal conditions---for example, fibrant replacements in either the injective model structures  $\precat _{\rm inj}(X,\mM )$ and $\precat _{\rm inj}(Y,\mM )$. 
Recall that $\Ob (f)^{\ast}(\pB )$ is an $\mM$-precategory on $X= \Ob (\pA )$, as is $\Ob (f)^{\ast}(\Seg (\pB ))$. 
Furthermore, the map $\Ob (f)^{\ast}(\pB )\rightarrt \Ob (f)^{\ast}(\Seg (\pB ))$
is again a weak equivalence in the model structure $\precat _{\rm inj}(X,\mM )$ (see Lemma \ref{pbwe}). 
It is also again a levelwise cofibration, so it is again a trivial cofibration. 

By making an appropriate choice of  
construction of the operation $\Seg (\cdot )$ we can obtain that $f$ induces a map $\Seg (\pA )\rightarrt \Ob (f)^{\ast}(\Seg (\pB ))$; however, 
in any case we could modify the choice of $\Seg (\pA )$ in order to obtain such a factorization. So we shall assume that this factorization is given. 
This results in a map $\Seg (\pA )\rightarrt \Seg (\pB )$. 

We say  that $f:\pA \rightarrt \pB $ is {\em fully faithful} if
the map $\Seg (\pA )\rightarrt \Seg (\pB )$ is fully faithful in the previous sense, i.e. it induces a levelwise equivalence 
$$
\Seg (\pA )(x_0,\ldots , x_p)\stackrel{\sim}{\rightarrt} \Seg (\pB )(f(x_0),\ldots , f(x_p))
$$
for each sequence of objects $(x_0,\ldots , x_p)\in \Delta _X{\Ob (\pA )}$. Since both $\Seg (\pA )$ and $\Seg (\pB )$ satisfy the 
Segal condition, it is easy to see that it is sufficient to check this on sequences of length one, that is 
it suffices to require that 
$$
\Seg (\pA )(x,y)\stackrel{\sim}{\rightarrt} \Seg (\pB )(f(x), f(y))
$$
for every pair of objects $x,y\in \Ob (\pA )$. 

We say that a map $f: (X,\pA )\rightarrt (Y,\pB )$ is {\em essentially surjective} if the functor of categories obtained by the truncation operation 
$\tau _{\leq 1}(f): \tau _{\leq 1}(X,\pA )\rightarrt \tau _{\leq 1}(Y,\pB )$ is an essentially surjective map of categories, in other words if the set of
isomorphism classes $\Iso \tau _{\leq 1}(\pA )$ of $\tau _{\leq 1}(X,\pA )$ surjects to the set of isomorphism classes 
$\Iso \tau _{\leq 1}(\pB )$ of $\tau _{\leq 1}(Y,\pB )$.

We say that $f: (X,\pA )\rightarrt (Y,\pB )$ is a {\em global weak equivalence} if it is fully faithful 
and essentially surjective. 

\begin{lemma}
\label{levelwiseglobalwe}
A morphism $f:\pA \rightarrt \pB$ which induces an isomorphism on sets of objects
$\Ob (f):\Ob (\pA )\cong \Ob (\pB )$ is a global weak equivalence if and only
if the corresponding map $f_{\sharp}: \Ob (f)_!\pA \rightarrt \pB$ in
$\precat (\Ob (\pB ),\mM )$ is a weak equivalence in the model structures of 
Theorems \ref{modstrucs} and \ref{reedyprecat}. 

In particular, if $\Ob (f)$ is an isomorphism and $f_{\sharp}$ is a levelwise weak equivalence
of $\Delta ^o_{\Ob (\pB )}$-diagrams in $\mM$, then $f$ is a weak equivalence. 
\end{lemma}
\begin{proof}
This is immediate from the definitions, since the essential surjectivity is automatically
guaranteed by the condition that $\Ob (f)$ be an isomorphism. For the last paragraph
note that Corollary \ref{oldwenew} applies in the direct left Bousfield localizations
of Theorems \ref{modstrucs} and \ref{reedyprecat}.
\end{proof}

\section{Categories enriched over $\Ho (\mM )$}

As in the proof of the preceding propsition,
using the small object argument for the pseudo-generating sets for trivial cofibrations in 
$\precat _{\rm proj}(X,\mM )$, we can obtain a functorial construction which associates to any $\pA \in \precat (\mM )$
another $\Seg (\pA )\in \precat (\Ob (\pA ),\mM )$ with a natural
transformation $\pA \rightarrt \Seg (\pA )$ which is a weak equivalence in $\precat (\Ob (\pA ),\mM )$ such that 
$\Seg (\pA )$ satisfies the Segal condition. For the present section, fix one such construction, although the discussion is left invariant
if we make a different choice.

The category $\Ho (\mM )$ admits finite direct products, and these are calculated by direct products in $\mM$ using conditions (DCL) and (PROD) (Lemma \ref{prodconsequence}).
Therefore it makes sense to look at the category $\Cat (\Ho (\mM ))$ of $\Ho (\mM )$-enriched categories, see Section \ref{sec-enriched}. 
Recall that a $\Ho (\mM )$-enriched category $C$ consists of a set  $\Ob (C)$ and
for any $x,y\in \Ob (C)$, a morphism object $C(x,y)\in \Ho  (\mM )$, plus composition maps satisfying strict associativity. 

\begin{lemma}
Suppose $\pA \in \precat (\mM )$. Then we obtain a $\Ho (\mM )$-enriched category denoted $\hec (\pA )$ with the same set of objects as $\pA $,
by putting $\hec (\pA )(x,y)$ equal to the class of $\Seg (\pA )(x,y)$ in $\Ho (\mM )$. The remainder of the simplicial diagram structure of $\Seg (\pA )$
(up to $\Seg (\pA )(x_0,x_1,x_2,x_3)$) provides this with composition maps which are unital and strictly associative. This gives a functor from 
$\precat (\mM )$ to $\Cat (\Ho (\mM ))$. 
\end{lemma}
\begin{proof}
Composing $\Delta ^o_X\rightarrt^{\Seg (\pA )} \mM $ with the projection $\mM \rightarrt 
\Ho (\mM )$ we get a functor $\hec (\pA ): \Delta ^o _X\rightarrt \Ho (\mM )$. 
The homotopy Segal conditions for $\Seg (\pA )$ imply that the Segal maps for $\hec (\pA )$
are isomorphisms. By the interpretation of Theorem \ref{enriched-interp} (the reader is
encouraged by now to have looked at the proof of this theorem, which was left as an exercise),
this says that $\hec (\pA )$ corresponds to a $\Ho (\mM )$-enriched category. 
\end{proof}

Recall that the truncation $\tau _{\leq 0}$  on
$\mM$ factors through a functor denoted the same way $\tau _{\leq 0}: \Ho (\mM )\rightarrt \Sets$ compatible with direct products,
so and applying that to the enrichment category gives a functor 
$$
\tau _{\leq 1}^h: \Cat (\Ho (\mM ))\rightarrt \Cat .
$$
The image of $\hec (\pA )$ under this functor is exactly $\tau _{\leq 1}(\pA )$ according to its construction, that is to say 
\begin{equation}
\label{tausame}
\tau _{\leq 1}^h(\hec (\pA ))=\tau _{\leq 1}(\pA ).
\end{equation}
Recall from Section 
\ref{sec-enriched} that a functor $g:C\rightarrt C'$ in $\Cat (\Ho (\mM ))$ is an {\em equivalence of categories} if $g$ essentially surjective i.e. 
$$
\Iso \tau ^h_{\leq 1}(g):\Iso \tau ^h_{\leq 1}(C) \rightarrt \Iso \tau ^h_{\leq 1}(C')
$$
is surjective, and if $g$ is fully faithful i.e. for any $x,y\in \Ob (C)$ the map $C(x,y)\rightarrt C'(g(x),g(y))$ is an isomorphism 
in $\Ho (\mM )$. As in Lemma \ref{transferequiv}, if $g$ is an equivalence of 
categories then $\Iso \tau ^h_{\leq 1}(g)$ is an isomorphism. 

\begin{proposition}
\label{conservative}
A morphism of $\mM$-precategories 
$f:\pA \rightarrt \pB $ is a global weak equivalence in $\precat (\mM )$ if and only if the corresponding morphism 
$\hec (f): \hec (\pA )\rightarrt \hec (\pB )$ is an equivalence of 
$\Ho (\mM )$-enriched categories.
\end{proposition}
\begin{proof}
By the identification \eqref{tausame} between truncations, the
essential surjectivity condition is the same. The fully faithful conditions
are the same too, because a morphism in $\mM$ is an equivalence if and only 
if its image in $\Ho (\mM )$ is an isomorphism.  
\end{proof}

This point of view allows us to apply the arguments of Section \ref{sec-enriched}, which concerned a simpler $1$-categorical notion of enrichment,
to obtain some important first properties of the class of global weak equivalences. 

\begin{lemma}
\label{levelwiseff}
Suppose $f$ is a global weak equivalence. Then $\tau _{\leq 1}(f): \tau _{\leq 1}(\pA )\rightarrt \tau _{\leq 1}(\pB )$ is an equivalence of categories, in particular
it induces an isomorphism of sets $\Iso \tau _{\leq 1}(\pA )\cong \Iso \tau _{\leq 1}(\pB )$. 
\end{lemma}
\begin{proof}
Use Lemma \ref{transferequiv}.
\end{proof}

\begin{proposition}
\label{globalretract32}
The class of global weak equivalences is closed under retracts and satisfies 3 for 2.
Furthermore if $f:\pA \rightarrt \pB$ and $g:\pB \rightarrt \pA$ are morphisms 
such that $fg$ and $fg$ are global weak equivalences, then $f$ and $g$ are global
weak equivalences. 
\end{proposition}
\begin{proof}
Using Proposition \ref{conservative}, this becomes an immediate consequence of Theorem \ref{enrichedretract32}.
\end{proof}

\section{Change of enrichment category}
\label{changeenrichment}

The following discussion is the key to resolving the issue that was raised by Pelissier in \cite{Pelissier}: he found a mistake in \cite{svk},
in what will correspond to our construction of interval objects in Chapter \ref{interval1} below.
Pelissier solved the problem for $\mM  =\mK$, but one can go pretty easily from there to any $\mM$ by investigating what happens under change of enrichment category.
This strategy will be used to construct interval objects in Chapter \ref{interval1}.

Suppose $F:\mM \rightarrt \mM '$ is a left Quillen functor between tractable left
proper cartesian model categories; denote its right adjoint by $F^{\ast}$. Then $F$ 
it induces left Quillen functors between the levelwise diagram model categories
$$
F^{\Delta ^o_X}: \diag _{\rm inj}(\Delta ^o_X, \mM )\rightarrt \diag _{\rm inj}(\Delta ^o_X, \mM ')
$$
and
$$
F^{\Delta ^o_X}: \diag _{\rm proj}(\Delta ^o_X, \mM )\rightarrt \diag _{\rm proj}(\Delta ^o_X, \mM ')
$$
whose right adjoints are obtained by applying the right adjoint $F^{\ast}$ levelwise. 
The Quillen property is verified levelwise on the left for the injective structure
and levelwise on the right for the projective structure. 

Define the functor
$$
\precat (X,F) : \precat (X,\mM )\rightarrt \precat (X,\mM ')
$$
by putting 
$$
\precat (X,F)(\pA ):= U_!(F^{\Delta ^o_X}(U^{\ast}\pA )).
$$
The reader may note that if we assume furthermore $F(\ast )=\ast$ then 
$\precat (X,F)$ is just the restriction of $F^{\Delta ^o_X}$ to the
unital diagrams, and the following discussion could be simplified. 
We don't make that assumption in general. 

\begin{lemma}
The functor $\precat (X,F)$ is left adjoint to the functor 
$\precat (X,F^{\ast})$ obtained by restricting the levelwise adjoint to the
unital diagram categories. It is compatible with $F^{\Delta ^o_X}$ via the diagram
$$
\begin{diagram}
\diag (\Delta ^o_X, \mM )&\rightarr^{F^{\Delta ^o_X}}& \diag (\Delta ^o_X, \mM ') \\
\downarr^{ U_!} & & \downarr _{ U_!} \\ 
\precat (X,\mM )&\rightarr^{\precat (X,F)}& \precat (X,\mM ').
\end{diagram}
$$
\end{lemma}
\begin{proof}
Note that $F^{\ast}$ commutes with limits so $F^{\ast}(\ast )=\ast$,
in particular applying $F^{\ast}$ levelwise preserves the unitality condition,
defining a functor $\precat (X,\mM ')\rightarrt^{\precat (X,F^{\ast})}\precat (X,\mM )$.  
We verify that $\precat (X,F)$ as defined above is its left adjoint. 
Suppose $\pA \in \precat (X,\mM )$ and $\pB \in \precat (X,\mM ')$. 
A map $\precat (X,F)(\pA )\rightarrt \pB $ in $\precat (X,\mM ')$ is the same
as a map $F^{\Delta ^o_X}(U^{\ast}\pA )\rightarrt U^{\ast}\pB $ which in turn is
the same as a map $U^{\ast}\pA \rightarrt F^{\ast , \Delta ^o_X}U^{\ast}\pB $,
but this is the same as a map $\pA \rightarrt \precat (X,F^{\ast})\pB $. 

For the compatibility diagram note that 
$$
\precat (X,F)(U_!\pA ) = U_!(F^{\Delta ^o_X}(U^{\ast}U_!\pA ))=U_!(F^{\Delta ^o_X}(\pA )).
$$
This can be seen by the explicit description of the unitalization operation
$U^{\ast}U_!$
as a coproduct at the constant sequences 
$$
U_!(F^{\Delta ^o_X}(U^{\ast}U_!\pA )) (x^n) = 
F(U^{\ast}U_!\pA (x^n))\cup ^{F(U^{\ast}U_!\pA (x))}\ast 
$$
$$
=
F(\pA (x^n)\cup ^{\pA (x)}\ast)\cup ^{F(\ast )}\ast 
=
F(\pA (x^n))\cup ^{F(\pA (x))}\ast = U_! (F^{\Delta ^o_X}(\pA ))(x^n).
$$
Recall that the unitalization is trivial at nonconstant sequences. 
\end{proof}

Assume that $\mM$ and $\mM '$ are tractable left proper cartesian model categories.
If $F:\mM \rightarrt \mM '$ is a left Quillen functor, then for any finite
collection of objects $A_1,\ldots , A_m$ (including the empty collection with $m=0$) we have
$$
F(A_1\times \cdots \times A_m)\rightarrow F(A_1)\times \cdots \times F(A_m).
$$
We say that $F$ is {\em product-compatible} if these maps are weak equivalences
for any sequence of objects. In particular
for the case $m=0$ this says that $F(\ast )$ is contractible.

Strangely enough, this natural condition does not seem to be 
needed for the following proposition. 

\begin{proposition}
\label{precatF}
Suppose that  $\mM \rightarrt ^F \mM'$ 
is a left Quillen functor between two 
tractable left proper cartesian model categories.
Then the functor 
$\precat (X; F)$ is a left Quillen functor between the projective (resp. Reedy) model category structures on
$\precat (X,\mM )$ and $\precat (X,\mM ')$. If $\mM '$ furthermore
satisfies Condition \ref{cond-inj} then 
$\precat (X; F)$ is also a left Quillen functor between the injective structures. 
\end{proposition}
\begin{proof}
We first note that $\precat (X,F)$ is a left Quillen functor between unital
diagram categories 
$$
\diag (\Delta ^o_X/X,\mM )\rightarrt \diag (\Delta ^o_X/X,\mM ')
$$
in the projective and Reedy structures (and the injective structures if Condition \ref{cond-inj} holds). 
 
For the projective structure, the fibrations and trivial fibrations are defined
levelwise, but the right adjoint $\precat (X,F^{\ast})$ is also defined
levelwise so it preserves levelwise (trivial) fibrations. For the Reedy structure,
recall that the diagrammatic Reedy structure is defined using the inclusion 
$U^{\ast}$ from unital to all diagrams: a morphism $f$ of unital diagrams is a 
Reedy (trivial) cofibration if and only if $U^{\ast}f$ is.
So $F^{\Delta ^o_X}(U^{\ast}$ sends Reedy (trivial) cofibrations 
in $\diag _{\rm Reedy}(\Delta ^o_X/X,\mM )$ to 
Reedy (trivial) cofibrations in $\diag _{\rm Reedy}(\Delta ^o_X,\mM )$.
But now $U_!$ is left Quillen for the Reedy structure by Proposition \ref{unitalreedy}, so 
$\precat (X,F)=U_!F^{\Delta ^o_X}(U^{\ast}$ preserves Reedy (trivial) cofibrations. 

If we assume furthermore Condition \ref{cond-inj} then by 
Lemma \ref{injectiveunitalization}, $U_!$ is left Quillen for the injective structures
so the same argument works. This completes the proof of what is claimed in the first
paragraph. 

Now $\precat (X;F)$ passes to a left Quillen functor
between the direct left Bousfield localizations.
We may assume that $F$ sends the generating sets for $\mM$ into ones for $\mM '$. Then
the functor $\precat (X;F)$ sends the pseudo-generators 
$P_{(x_0,\ldots , x_n),!}(\zeta _n (f))$ for the
direct left Bousfield localized structure on $\precat (X,\mM )$, to 
pseudo-generators for $\precat (X,\mM ')$; in particular
to trivial cofibrations. By the left Bousfield localization property
$\precat (X;F)$ are left Quillen functors between the direct left Bousfield localizations. 
\end{proof}


\chapter{Cofibrations}
\label{cofib1}

In this chapter, continuing the construction of model structures for the category 
$\precat (\mM )$ of $\mM$-enriched precategories 
with variable set of objects, we define and discuss various classes of cofibrations. 
The corresponding classes of trivial cofibrations
are then defined as the intersection of the cofibrations with the weak equivalences 
defined in the preceding chapter.  
The fibrations are defined by the lifting property with respect
to trivial cofibrations. It will have to be proven later that the class of trivial 
fibrations, defined as the intersection of the fibrations
and the weak equivalences, is also defined by the lifting property with respect to cofibrations. 

In reality, we consider  three model structures, called the {\em injective}, 
the {\em projective}, and the {\em Reedy} structures.
These will be denoted by subscripts $I,P,R$ respectively when necessary. They share the 
same class of weak equivalences, but the cofibrations are different. 
It turns out that the Reedy structure is the best one for the purposes of iterating the 
construction. That is a somewhat subtle point 
because the Reedy structure coincides with the injective structure when $\mM$ lies in a 
wide range of model categories where the monomorphisms are 
cofibrations, see Proposition \ref{reedyinjective}. When they are different it is better to take the Reedy route.

\section{Skeleta and coskeleta}

The definition of Reedy cofibrations, as well as the study of the projective ones, is 
based on consideration of 
the skeleta and coskeleta of objects in $\precat (\mM )$. 
The Reedy structure will be defined in an explicit way without making use of the general 
definition of Reedy category.
The category $\precat (\mM )$ consists of precategories with various different 
object sets $X$, so it isn't enough to just invoke the Reedy structure of each $\Delta _X$.  
Considering instead the relevant structures explicitly
has the added advantage of making the discussion accessible for readers who wish to 
avoid plunging into the full theory of Reedy categories. Those already
familiar with those kinds of things can try to view the discussion in a more general 
language 
(see Barwick \cite{BarwickReedy}, Berger-Moerdijk \cite{BergerMoerdijk} and others). 

For a fixed set $X$, consider the subcategory $\Delta _{X,m}\subset \Delta _X$ 
consisting of all sequences $(x_0,\ldots , x_p)$ 
of length $p\leq m$. For example, $\Delta _{X,0}$ consists of sequences of length $0$ 
i.e. $(x_0)$ only. 
Let $\sigma (X,m):\Delta _{X,m}\hookrightarrow \Delta _X$ denote the inclusion functor. 
Then we have a functor
$$
\sk _m :=  \sigma (X,m)_! \circ \sigma (X,m)^{\ast} : \diag (\Delta _X, \mM )
\rightarrt \diag (\Delta _X, \mM ).
$$
We call $\sk _m (\pA )=  \sigma (X,m)_! ( \sigma (X,m)^{\ast}(\pA ))$ the {\em $m$-skeleton} 
of $\pA $. There is a natural
map $\sk _m(\pA )\rightarrt \pA $.

Given $m$ and $\pA $ we would like to calculate $\sk _m(\pA )$. 
Suppose $x_{\cdot}=(x_0,\ldots , x_k)$ is a sequence of elements of $X$ of length $k$.
Consider the category of all surjections $x_{\cdot}\twoheadrightarrow y_{\cdot}$ where
$y_{\cdot}= (y_0,\ldots , y_p)$ is a sequence of length $p\leq m$. Morphisms are
morphisms of objects $y_{\cdot}$ in $\Delta _X$ commuting with the maps 
from $x_{\cdot}$. Note that
$y_{\cdot}\in \Delta _{X,m}$ so the pullback 
of $\pA $ to $\Delta _{X,m}$ restricts, by the projection functor to the variable $y_{\cdot}$,
to a contravariant functor from our category of surjections to $\mM$. We claim that 
\begin{equation}
\label{skformula}
\sk _m(\pA )(x_{\cdot}) = \colim _{x_{\cdot}\twoheadrightarrow y_{\cdot}}
\pA (y_{\cdot}).
\end{equation}
Indeed, this is almost exactly the same as the
expression of Lemma \ref{fshriek} for the pushforward $\sigma (X,m)_!$.
In the general expression, the colimit is taken over all morphisms from 
$x_{\cdot}$ to some $y_{\cdot}\in \Delta _{X,m}$. However, any such morphism
canonically factors as 
$x_{\cdot}\twoheadrightarrow y'_{\cdot}\hookrightarrow y_{\cdot}$
with $y'_{\cdot}\in \Delta _{X,m}$ too, so the category of
surjections is cofinal and suffices for calculating the colimit. 

From this expression we conclude that if $x_{\cdot}=(x_0,\ldots , x_k)$
is a sequence of length $k\leq m$ then the natural map
$$
\sk _m(\pA )(x_{\cdot})\rightarrt \pA (x_{\cdot})
$$
is an isomorphism, in other words the restrictions of $\pA $ and $\sk _m(\pA )$ to
$\Delta _{X,m}$ are the same. It follows that $\sk _m$ is a monad, in other
words $\sk _m(\sk _m(\pA ))=\sk _m(\pA )$ (or rather they are isomorphic in a natural
way). We say that $\pA $ is {\em $m$-skeletic} if the map $\sk _m(\pA )\rightarrt \pA $
is an isomorphism; and $\sk _m$ is a monadic projection to the full subcategory
of $m$-skeletic objects. 

If $n\geq m$ and $\pA $ is $m$-skeletic then it is
also $n$-skeletic, as $\sk _n(\sk _m(\pA ))=\sk _m(\pA )$. This formula comes from
the fact that the restriction to $\Delta _{X,n}$ of $\sk _m(\pA )$ is the same as
the pushforward from $\Delta _{X,m}$ to $\Delta _{X,n}$ of the restriction
$\sigma (X,m)^{\ast}(\pA )$, as can be seen by the expressions of the pushforwards
as colimits over surjective maps; then the composition of pushforwards
from $\Delta _{X,m}$ to $\Delta _{X,n}$ and then to $\Delta _X$,
is the same as the pushforward $\sigma (X,m)_!$ to give the claimed formula.

\begin{lemma}
The functor $\sk _m$ preserves unitality, that is if $\pA \in \precat (X,\mM )$ then 
$\sk _m(\pA )\in \precat (X,\mM )$.
It has a right adjoint $\csk _m : \precat (X,\mM )\rightarrt \precat (X,\mM )$ 
called the {\em $m$-coskeleton}.
Both of these are compatible with changing $X$ so they give functors 
$\precat (\mM )\rightarrt \precat (\mM )$. 
The natural morphism induces an isomorphism
$\colim _m\sk _m(\pA )\cong \pA $.
The morphisms $\sk _{m-1}(\pA )\rightarrt \sk _m(\pA )
\rightarrt \pA $ are monomorphisms if $\mM$ is a presheaf category. 
\end{lemma}
\begin{proof}
Suppose $(x_0)$ is a sequence of length zero. Then $(x_0)\in \Delta _{X,m}$ for any
$m\geq 0$ so the category of surjections $(x_0)\twoheadrightarrow y_{\cdot}$
considered above has as initial object $(x_0)$ itself. The formula \eqref{skformula}
then says 
$$
\sk _m(\pA )(x_0) = \pA (x_0).
$$
It follows that if $\pA $ is unital, so is $\sk_m(\pA )$.

The functors $\sigma (X,m)_!$ and $\sigma (X,m)^{\ast}$ have
right adjoints $\sigma (X,m)^{\ast}$ and $\sigma (X,m)_{\ast}$ respectively.
So, considered as a functor between diagram categories, $\sk_m$ 
has as right adjoint:
$$
\csk_m (\pB ) = \sigma (X,m)_{\ast}\sigma (X,m)^{\ast}(\pB ).
$$
As before we can note that this functor has the expression
$$
\csk _m(\pB ) (x_{\cdot}) = \lim _{y_{\cdot}\hookrightarrow x_{\cdot}} \pB (y_{\cdot})
$$
where the limit is taken over inclusions $y_{\cdot}\hookrightarrow x_{\cdot}$
in $\Delta _X$ such that $y_{\cdot}$ is in $\Delta _{X,m}$ i.e. has length $\leq m$.
Then again note that if $x_{\cdot}=(x_0)$ is a sequence of length zero, the
index category for the limit has an initial object $(x_0)$ itself, so 
$$
\csk _m(\pB )(x_0)=\pB (x_0).
$$
Thus $\csk _m$ preserves the unitality condition and gives an endofunctor of
$\precat (X,\mM )$, right adjoint to $\sk _m$. 

Suppose $\psi : X\rightarrt Y$ is a morphism of sets. 
For any $\pA \in \precat (X,\mM )$ and $\pB \in \precat (Y,\mM )$
together with a morphism $f:\pA \rightarrt \psi ^{\ast}\pB $ in $\precat (X,\mM )$
we would like to get a map $\sk _m(\pA )\rightarrt \psi ^{\ast}\sk _m(\pB )$.
In view of the definition of $\sk _m$ this is equivalent to giving a map
of diagrams over $\Delta _{X,m}$
\begin{equation}
\label{skfunc}
\sigma (X,m)^{\ast}(\pA )\rightarrt \sigma (X,m)^{\ast}\psi ^{\ast}\sk _m(\pB ).
\end{equation}
Let $\psi _m:\Delta _{X,m}\rightarrt \Delta _{Y,m}$ denote the
induced functor, then $\psi \circ \sigma (X,m) = \sigma (Y,m)\circ \psi_m$
so 
$$
\sigma (X,m)^{\ast}\psi ^{\ast}\sk _m(\pB ) = \psi _m^{\ast}\sigma (Y,m)^{\ast}\sk _m(\pB ) 
$$
$$
= \psi _m ^{\ast}\sigma (Y,m)^{\ast}(\pB ) = \sigma (X,m)^{\ast}\psi ^{\ast}(\pB )
$$
and functoriality for $\sigma (X,m)^{\ast}$ applied to $f$ gives a map
$$
\sigma (X,m)^{\ast}(f) : \sigma (X,m)^{\ast}(\pA )\rightarrt 
\sigma (X,m)^{\ast}\psi ^{\ast}(\pB ).
$$
Which yields the desired map \eqref{skfunc}. We also have a commutative square
$$
\begin{diagram}
\sk _m(\pA ) & \rightarr & \psi ^{\ast}\sk _m(\pB ) \\
\downarr & & \downarr \\
\pA & \rightarr & \psi ^{\ast}\pB .
\end{diagram}
$$
So, given a morphism $\pA \rightarrt \pB $
in $\precat (\mM )$ this construction has given $\sk _m(\pA )\rightarrt \sk_m(\pB )$.
It is compatible with identities and compositions so $\sk _m$ is an endofunctor
of $\precat(\mM )$, with
natural transformation $\sk _m(\pA )\rightarrt \pA $ again giving a structure of 
monadic projection to the full subcategory of $m$-skeletic precategories.

Similarly for the coskeleton we look for a map
\begin{equation}
\label{cskfunc}
\psi _!\csk _m(\pA )\rightarrt \sigma (Y,m)_{\ast}\sigma (Y,m)^{\ast}(\pB ).
\end{equation}
This is equivalent to looking for 
$$
\sigma (Y,m)^{\ast}\psi _!\csk _m(\pA ) \rightarrt \sigma (Y,m)^{\ast}(\pB ).
$$
But in general 
$$
\psi _!(\pC )(y_{\cdot}) = \colim _{y_{\cdot}\rightarrt \psi x_{\cdot}}\pC (x_{\cdot})
$$
and the colimit can be taken over surjective maps 
$y_{\cdot}\twoheadrightarrow \psi x_{\cdot}$. Thus if $y_{\cdot}$ has length $m$
then it only depends on the restriction of $\pC $ to $\Delta _{X,m}$, in particular
$$
\sigma (Y,m)^{\ast}\psi _!\csk _m(\pA ) =\sigma (Y,m)^{\ast}\psi _!(\pA ) 
$$
and we can look for a map $\psi _!\pA \rightarrt \csk _m(\pB )$ or equivalently
by adjunction,
$\pA \rightarrt \psi ^{\ast}\csk _m(\pB )$. But the natural map $\pB \rightarrt \csk _m(\pB )$
gives by pullback under $\psi$ and composition 
$$
\pA \rightarrt \psi ^{\ast}(\pB )\rightarrt \psi ^{\ast}\csk _m(\pB )
$$
as required. We get the functoriality maps for defining 
$$
\csk _m : \precat (\mM )\rightarrt \precat (\mM ), 
$$
again with a natural
transformation $\pA \rightarrt \csk _m(\pA )$. 

Suppose $\pA \in \precat (X,\mM )$. 
For any $x_{\cdot}$ a sequence of length $p$, then for $m\geq p$ we have
$\sk _m(\pA )(x_{\cdot})= \pA (x_{\cdot})$. Hence the morphism 
$$
\colim _m\sk _m(\pA )\rightarrt \pA 
$$
is an isomorphism on each object $x_{\cdot}$ of $\Delta _X$. As colimits in diagram categories
are computed objectwise and the unitalification operation $U_!$ preserves colimits,
we have $\pA \cong \colim _m\sk _m(\pA )$ in $\precat (X,\mM )$. Recall from Corollary 
\ref{constobjcolim} that
connected colimits in $\precat (X,\mM )$ also are colimits in $\precat (\mM )$ so
we get the same formula in $\precat (\mM )$. 

One can remark that dually and for the same reason, 
$$
\pA \rightarrt^{\cong} \lim _m\csk _m(\pA ).
$$

Finally we note that if $\mM$ is a presheaf category then we have a family of functors
$h_u:\mM \rightarrt \Sets$, the evaluations on objects of the underlying category for
the presheaves, such that $h_u$ preserves colimits. It follows that $h_u\circ (\sk _m(\pA ))=
\sk _m(h_u\circ \pA )$, however $h_u\circ \pA $ is a set-valued diagram over $\Delta _X$.
As is well-known, the inclusion of the $m$-skeleton is an injection of simplicial sets,
and this works equally well for $\Delta _X^o$-sets. Thus, the map
$$
h_u\circ (\sk _m(\pA ))=\sk _m(h_u\circ \pA )\rightarrt h_u\circ \pA 
$$
is a monomorphism. Therefore the map $\sk _m(\pA )\rightarrt \pA $ induces a monomorphism
upon application of all of the functors $h_u$; but as $\mM $ is a presheaf category
this implies that our map is a monomorphism as claimed in the last statement of the lemma.
\end{proof}

Notice that $\sk _0(\pA )=\disc (\Ob (\pA ))$ is the set $\Ob (\pA )$ considered as a discrete precategory,
that is one whose $p$-fold morphism objects are all $\ast$ for constant sequences or 
$\emptyset$ otherwise (see Section \ref{sec-precatexamples}). Indeed if $x_{\cdot}$ is a sequence then 
there is a surjection to a unique $(y_0)$ of length $0$, if and only if $x_{\cdot}$ 
is constant, so $\sk_0(\pA )(x_{\cdot})=\pA (y_0)=\ast$ if $x_{\cdot}$ is constant of value
$y_0$, and $\sk_ 0(\pA )(x_{\cdot})= \emptyset$ if $x_{\cdot}$ is nonconstant.

If $\pA \in \precat (\mM )$ and $(x_0,\ldots , x_m)$ is a sequence of objects with $m\geq 1$, 
define the {\em degenerate subobject}
$$
\dgt (\pA ; x_0,\ldots , x_m):= \sk _{m-1}(\pA )(x_0,\ldots , x_m).
$$
It has a map denoted 
$$
\delta (\pA ; x_0,\ldots , x_m): \dgt (\pA ; x_0,\ldots , x_m)\rightarrt \pA ( x_0,\ldots , x_m),
$$
and we can express $D$ as a colimit:

\begin{lemma}
\label{degeneratesub}
With the above notations, 
$$
\dgt (\pA ; x_0,\ldots , x_m) = \colim _{x_{\cdot}\rightarrow y_{\cdot}} \pA (y_0,\ldots , y_p)
$$
where the colimit is taken over surjective maps $x_{\cdot}\rightarrow y_{\cdot}$ in 
$\Delta _{\Ob (\pA )}$
such that $y_{\cdot}$ are sequences of length $p<m$. The map $\delta (\pA ; x_0,\ldots , x_m)$ 
is obtained from the restriction maps of $\pA $
via the universal property of the colimit. The map $\delta (\pA ; x_0,\ldots , x_m)$
is a monomorphism if $\mM$ is a presheaf category. 
\end{lemma}
\begin{proof}
This comes from the definition of $D$ and the colimit expression
for $\sk_{m-1}$. The monomorphism property comes from the previous lemma. 
\end{proof}

\section{Some natural precategories}
\label{sec-naturalprecats}

Consider the ordered sets $[k]\in \Ob (\Delta )$ with the notation
$$
[k] = \{ \upsilon _0,\ldots , \upsilon _k \} , \;\;\; \upsilon _0 < \cdots < \upsilon _k .
$$

For any $B\in \mM$ and any $[k]\in \Ob (\Delta )$ define the precategory 
$h([k]; B)\in \precat ([k],\mM )$
as follows.
Let ${\bf t}_k\in \Ob (\Delta _{[k]})$ denote the tautological object 
${\bf t}_{k}:= ( \upsilon _0,\ldots , \upsilon _k)$ of length $k$, 
and let 
$$
\bfi \{ {\bf t}_k\} : \{ {\bf t}_k\} \hookrightarrow \Delta _{[k]}
$$
denote the inclusion
of the discrete category on a single point ${\bf t}_k\in \Ob (\Delta _{[k]})$. 
Let  $B_{{\bf t}_k}$ denote the constant diagram with value $B$ on the one-point category 
$\{ {\bf t}_k \}$, and put
$$
h([k]; B):= U_!\bfi \{ {\bf t}_k\}_! (B_{{\bf t}_k}).
$$
Recall that $U_!$ is the unitalization operation, necessary here because the $\Delta_X$-diagram 
$\bfi \{ {\bf t}_k\}_! (B_{{\bf t}_[k]})$ will not in general be unital. The following 
lemma gives a concrete
description of $h([k]; B)$ and could be taken as its definition. 

\begin{lemma}
\label{hdescription}
Suppose $(y_0,\ldots , y_p)$ is any sequence of elements of the set $[k]$ with 
$y_j = \upsilon _{i_j}$. Then:
\newline
---if $(y_0,\ldots , y_p)$ is increasing but not constant i.e. $i_{j-1}\leq i_j$ but 
$i_0 < i_p$ then
$$
h([k]; B)(y_0,\ldots , y_p) = B; 
$$
---if $(y_0,\ldots , y_p)$ is  constant i.e. $i_0=i_1=\ldots = i_p$ then
$$
h([k]; B)(y_0,\ldots , y_p) = \ast; 
$$
and otherwise, that is if there exists $1\leq j \leq p$ such that $i_{j-1}>i_j$ then
$$
h([k]; B)(y_0,\ldots , y_p) = \emptyset.
$$
The pullback maps giving $h([k]; B)$ a structure of diagram are all either the unique 
maps of the form $\emptyset \rightarrt B$,
$\emptyset \rightarrt \ast$, or $B\rightarrt \ast$, or the identity $B\rightarrt B$.
\end{lemma}
\begin{proof}
The functor $B\mapsto h([k]; B)$ is by construction left adjoint to the functor 
$\bfi \{ {\bf t}_k\}^{\ast}$
from $\precat ([k],\mM )$ to $\mM$ which sends $\pA $ to $\pA (\upsilon _0,\ldots , \upsilon _k)$.
On the other hand we can check by hand (as in the proof of the next lemma below) 
that the functor sending $B$ to
the precategory defined explicitly
in the statement of the lemma, is also adjoint to the same functor. 
\end{proof}

\begin{lemma}
\label{universalh}
If $B\in \mM$ and $\pC \in \precat (\mM )$ then a morphism $f:h([k], B)\rightarrt \pC $ is the 
same thing as a
sequence of objects $x_0,\ldots , x_k\in \Ob (\pC )$ together with a map 
$\varphi: B\rightarrt \pC (x_0,\ldots , x_k)$ in $\mM$.
\end{lemma}
\begin{proof}
Use the explicit description of the previous lemma. Given $f$ we get $x_i:= f(\upsilon _i)$
and the map $\varphi$ is given by $f_{\upsilon _0,\ldots , \upsilon _k}$.
On the other hand, given $x_{\cdot}$ and $\varphi$,  the restrictions
$\pC (x_0,\ldots , x_k)\rightarrt \pC (x_{i_0},\ldots , x_{i_p})$ lead to maps
$$
B\rightarrt \pC (x_{i_0},\ldots , x_{i_p})
$$
for any sequence as in the first part of the previous lemma; if $i_0=\cdots = i_p$
then this map factors through $\ast$ by the unitality condition for $\pC $, treating
the second part of the previous lemma; and nothing is needed for defining a map in the 
third case of the previous lemma. This defines the required map $f$. 
\end{proof}

Note that the above construction is functorial in $B$, that is a map 
$A \rightarrt B$ induces $h([k];A )\rightarrt h([k]; B)$.
Define the ``boundary'' of  $h([k]; B)$ by the skeleton operation:
$$
h(\partial [k]; B) := \sk _{k-1}h([k]; B),
$$
with the natural inclusion
$$
h(\partial [k]; B)\rightarrt h([k]; B).
$$
It is also functorial in $B$. 

\begin{lemma}
\label{hdelbardescription}
This boundary object has the following concrete description:
\newline
---if $(y_0,\ldots , y_p)$ is increasing but not constant i.e. 
$i_{j-1}\leq i_j$ but $i_0 < i_p$, and
if there is any $0\leq m\leq k$ such that $i_j\neq m$ for all $0\leq j\leq k$, then
$$
h(\partial [k]; B)(y_0,\ldots , y_p) = B; 
$$
---if $(y_0,\ldots , y_p)$ is  constant i.e. $i_0=i_1=\ldots = i_p$ then
$$
h(\partial [k]; B)(y_0,\ldots , y_p) = \ast; 
$$
and otherwise, that is if either there exists $1\leq j \leq p$ such that $i_{j-1}>i_j$
or else if the map $j\mapsto y_j$ is a surjection from $\{ 0,\ldots , p\}$ to $[k]$,  then
$$
h(\partial [k]; B)(y_0,\ldots , y_p) = \emptyset.
$$
\end{lemma}
\begin{proof}
Use the description of Lemma \ref{hdescription} and the formula
$$
h(\partial [k]; B)(y_{\cdot}) = \colim _{y_{\cdot}\twoheadrightarrow z_{\cdot}}
h(B)(z_{\cdot}).
$$
\end{proof}

If $f:A \rightarrt B$ is a cofibration in $\mM$, put
$$
h([k],\partial [k]; A\stackrel{f}{\rightarrow} B):= 
h([k]; A)\cup ^{h(\partial [k];A)} h(\partial [k];B).
$$
We therefore obtain two natural maps coming from $f$, the first is 
$$
P([k];f): h([k]; A)\rightarrt h([k]; B)
$$
and the second is \label{Rdefpage}
$$
R([k];f): h([k],\partial [k]; A\stackrel{f}{\rightarrow} B)\rightarrt h([k]; B) .
$$
The maps $P([k]; f)$ will form the generators for the projective cofibrations, while the 
$R([k]; f)$ form the generators for the Reedy cofibrations.
 Unfortunately we don't have an easy way to describe generators for the injective cofibrations. 

\section{Projective cofibrations}
The different model structures are characterized and 
differentiated by their notions of cofibrations. In the projective structure, the generating 
set is easiest to describe but on the other hand there is no 
easy criterion for being a cofibration. 

A map $f:\pA \rightarrt \pB $ is a {\em projective cofibration} 
if on the set of objects it is an injective map of sets
$\Ob (f):\Ob (\pA )\rightarrt \Ob (\pB )$, and if the map
$\Ob (f)_!(\pA )\rightarrt \pB $ is a cofibration in the projective 
model category structure $\precat _{\rm proj}(\Ob (\pB ),\mM )$.
Recall that $\Ob (f)_!(\pA )$ is the precategory $\pA $ transported 
to the subset $f(\Ob (\pA ))\subset \Ob (\pB )$, then extended by adding on the discrete
or initial precategory over the complementary subset $\Ob (\pB )-\Ob (f)(\Ob (\pA ))$. 

Say that $f$ is a {\em projective trivial cofibration} if it is a 
projective cofibration, and a global weak equivalence.

For now, we say that a morphism $u:U\rightarrt V$ in $\precat (\mM )$
is a {\em projective fibration} if it satisfies the right lifting property
with respect to all projective trivial cofibrations; we say that $u$ is a
{\em projective trivial fibration} if it is a projective fibration and a
global weak equivalence. 

For the time being we also need a separate notation: say that $u:U\rightarrt V$
in $\precat (\mM )$ is an {\em apparent projective trivial fibration} if
it satisfies the right lifting property with respect to all projective cofibrations.
One of our tasks in future chapters will be to identify the class of
projective trivial fibrations with the apparent ones. For now we can characterize
and use the apparent ones. 

\begin{lemma}
\label{apparentdual}
A morphism 
$u:U\rightarrt V$ in $\precat (\mM )$ is an apparent projective trivial fibration 
if and only if
$\Ob (u)$ is surjective and $U\rightarrt \Ob (u)^{\ast}(V)$
is an objectwise trivial fibration. A morphism $f:\pA \rightarrt \pB $ is a projective
cofibration if and only if it satisfies the left lifting property
with respect to the apparent projective trivial fibrations. 
\end{lemma}
\begin{proof}
Consider the class $\mA$ of morphisms $u$ such that 
$\Ob (u)$ is surjective and $U\rightarrt \Ob (u)^{\ast}(V)$
is an objectwise trivial fibration. If $f$ is a projective cofibration
then it breaks down as a composition $f=f'\circ d$ where $f':\Ob (f)_!(\pA )\rightarrt \pB $
is a morphism in $\precat (\Ob (\pB ),\mM )$ and $d$ is the extension by adding on
the discrete precategory $\Ob (\pB ) - \Ob (f)(\Ob (\pA ))$. Similarly, if $u\in \mA$ then
$u=p\circ u'$ where $u':U\rightarrt \Ob (u)^{\ast}(V)$ is an objectwise trivial
fibration in $\precat (\Ob (U),\mM )$ and $p$ is the tautological map
$\Ob (u)^{\ast}(V)\rightarrt V$. Notice that $p$ satisfies the right lifting property
with respect to any morphism which induces an isomorphism on sets of objects,
and also with respect to extensions by adding discrete sets; and dually $f'$ satisfies
the left lifting property with respect to tautological maps such as $p$, while
$d$ satisfies the left lifting property with respect to any map surjective on objects.
All told, the lifting property with $f$ on the left and $u$ on the right, is
equivalent to the lifting property with $f'$ on the left and $u'$ on the right.
This holds whenever
$u'$ is an objectwise trivial fibration and $f'$ is
a projective cofibration in $\precat (\Ob (\pB ),\mM )$. 

These considerations show that $\mA$ is contained in the class of apparent projective trivial
fibrations. Furthermore, if $u$ satisfies right lifting for any projective cofibration
$f$ then in particular $\Ob (u)$ must be surjective (using for $f$ any extension by a
nonempty discrete set); and $u'$ must satisfy right lifting with respect to 
projective trivial cofibrations inducing isomorphisms on sets of objects, so 
$u'$ is an objectwise trivial fibration (by the projective
model structure on the $\precat (X,\mM )$). This shows that the apparent projective
trivial cofibrations are contained in $\mA$ so the two classes coincide.

Then, similarly, if $f$ satisfies left lifting with respect to $\mA$ then first of all
it must be injective on sets of objects, as seen by considering $u$ which are surjective
maps of
codiscrete objects (i.e. objects $U$ with $U(y_{\cdot})=\ast$ for all sequences $y_{\cdot}$).
And then decomposing $f=f'\circ d$ where $d$ is the extension by the complementary
subset, we see that $f'$ should be a projective cofibration in $\precat (\Ob (\pB ),\mM )$.
Thus $f$ is a projective cofibration by definition. 
But all projective cofibrations satisfy lifting with respect
to the apparent projective trivial fibrations, by the definition of this latter class,
and equality with the class $\mA$ shows that the class of projective cofibrations
is exactly that which satisfies left lifting with respect to $\mA$.
\end{proof}

We have stability under pushouts and retracts.

\begin{corollary}
\label{Pcofibclosure}
Assume that $\mM$ is a tractable left proper cartesian model category. 
Suppose $f:\pA \rightarrt \pB $ and $g:\pA \rightarrt \pC $ are morphisms in 
$\precat (\mM )$ such that $f$ is a projective cofibration.
Then the map 
$$
u:\pC \rightarrt \pB \cup ^{\pA }\pC 
$$
is a projective cofibration. Furthermore the projective cofibrations 
are stable under retracts and transfinite composition. 
The class of  projective trivial cofibrations is 
closed under transfinite composition and retracts. 
\end{corollary}
\begin{proof}
Stability under pushouts, retracts and transfinite composition come from 
the characterization of the projective cofibrations as dual to the class of 
apparent projective trivial fibrations. 
The last sentence now follows immediately from Proposition \ref{globalretract32}.
\end{proof}

The advantage of the notion of projective cofibration 
is that the generating set is very easy to describe. 

\begin{proposition}
Fix a generating set $I$ for the cofibrations of $\mM$. 
Then for any $k$ and any $f:A \rightarrt B$ in $I$,
consider the map
$$
P([k];f): h([k]; A)\rightarrt h([k]; B).
$$
The collection of these for integers $k$ and all $f\in I$, 
together with the map $\emptyset \rightarrt \ast$,
forms a generating set for the class of projective cofibrations in 
$\precat (\mM )$.
\end{proposition}
\begin{proof}
Any projective cofibration $f:\pA \rightarrt \pB $ factors as $f=f'd$
where
$$
\pA \rightarrt^{d} \pA '\sqcup \disc (Z) \rightarrt^{f'} \pB 
$$
where $Z=\Ob (\pB )=\Ob (f)(\Ob (\pA ))$, where $\pA '=\Ob (f)_{!}(\pA )$
is the precategory on the set $\Ob (f)(\Ob (\pA ))$ obtained
by transport of structure, $\disc (Z)$ is the discrete precategory
on the set $Z$, and $d$ is the extension morphism. In this situation
furthermore, $f'$ is a
projective cofibration in $\precat (X,\mM )$ where $X=\Ob (\pB )=\Ob (\pA '\sqcup \disc (Z))$.

Now $d$ is obtained by successive pushout along $\emptyset \rightarrt \ast$,
while the $w_!P([k];f)$ for various maps $w:[k]\rightarrt X$ form
the set of generators for the projective cofibrations in $\precat (X,\mM )$. Thus,
$f'$ is a retract of a transfinite pushout along the $w_!P([k];f)$. Putting these
together gives the expression of $f$ as a retract of a transfinite pushout
of morphisms in our generating set. 
\end{proof}

\section{Injective cofibrations}
It is easier to see whether a given map is an injective cofibration, 
since this condition is defined objectwise,
but the only way to get a generating set in general is by an accessibility argument.

A morphism $f:\pA \rightarrt \pB $ is an {\em injective cofibration} 
if $\Ob (f)$ is an injective map of sets, and if the map
$\Ob (f)_!(\pA )\rightarrt \pB $ is a cofibration in the injective 
model category structure $\precat _{\rm inj}(\Ob (\pB ),\mM )$.
Since injective cofibrations are defined objectwise, we can be 
more explicit about this condition: it means that for any 
sequence $(x_0,\ldots , x_p)$ in $\Ob (\pA )$, the map 
$\pA (x_0,\ldots , x_p)\rightarrt \pB (f(x_0),\ldots , f(x_p))$ is a cofibration in $\mM$,
for any constant sequence $(y,\ldots , y)$ at a point 
$y\in Y-f(X)$ the map $\ast \rightarrt \pB (y,\ldots , y)$ is
a cofibration in $\mM$, and for any non-constant sequence 
$(y_0,\ldots , y_p)$ such that at least one of the $y_i$ is not in $f(X)$,
the map $\emptyset \rightarrt \pB (y_0,\ldots , y_p)$ is a 
cofibration i.e. $\pB (y_0,\ldots , y_p)$ is a cofibrant object. 

Say that $f$ is an {\em injective trivial cofibration} if 
it is an injective cofibration, and a global weak equivalence.

We again have the stability proposition: 

\begin{proposition}
\label{Icofibclosure}
Assume that $\mM$ is a tractable left proper cartesian model category. 
Suppose $f:\pA \rightarrt \pB $ and $g:\pA \rightarrt \pC $ are morphisms in 
$\precat (\mM )$ such that $f$ is an injective cofibration.
Then the map 
$$
\pC \rightarrt \pB \cup ^{\pA }\pC 
$$
is an injective cofibration. Furthermore the injective cofibrations 
are stable under retracts and transfinite composition. 
The class of  injective trivial cofibrations is closed under transfinite 
composition and retracts. 
\end{proposition}
\begin{proof}
The explicit conditions given in the paragraph defining the injective cofibrations,
are defined objectwise over $\Delta _Z$ where $Z$ is the set of objects of the target
precategory. The conditions are preserved by pushouts, retracts and transfinite
composition. 
Again, the last sentence now follows immediately from Proposition \ref{globalretract32}.
\end{proof}

\begin{lemma}
\label{injcofibpushout}
Suppose $\mM$ is a tractable model category. Then the class of injective 
cofibrations in $\precat (\mM )$ admits
a small set of generators.
\end{lemma}
\begin{proof}
Follow the argument given in Barwick \cite{Barwick} for the proof using a general 
accessibility argument. 
\end{proof}

Unfortunately, we don't get
very much information on the set of generators. In Section \ref{sec-reedyrelation} 
below, if $\mM$ satisfies some further hypotheses
which hold for presheaf categories, then the injective cofibrations are the 
same as the Reedy cofibrations and
the generators will be described explicitly. 

A similar problem with the injective structure is that if $\mM$ is tractable, we don't know
whether the generating cofibrations for $\precat (\mM )$ have cofibrant domains.
This question for diagram
categories isn't treated by Lurie in \cite{LurieTopos} and seems to remain an open question.

For our purposes this question is one further reason for introducing the Reedy model structure
which has an explicit set of generators; in many cases (such as when $\mM$ is a
presheaf category, Proposition \ref{reedyinjective}) the Reedy and injective structures will coincide.

\section{A pushout expression for the skeleta}
\label{exprskeleta}

The main observation crucial for understanding the Reedy cofibrations, 
is an expression for the successive skeleta as pushouts along maps of the form
$R([k]; f)$. If $x_0,\ldots , x_m$ is a sequence of objects, 
recall the notation 
$$
\dgt (\pA ; x_0,\ldots , x_m)\rightarrt ^{\delta (\pA ,x_{\cdot})}\pA (x_0,\ldots , x_m)
$$
from Section \ref{sec-naturalprecats}. The identity map
$$
\dgt (\pA ; x_{\cdot})\rEqualarr \sk _{m-1}(\pA )(x_0,\ldots , x_m)
$$
corresponds by the universal property of $h([m]; \cdot )$ to a map
$$
h([m]; \dgt (\pA ; x_{\cdot}))\rightarrt \sk _{m-1}(\pA )
$$
in $\precat (\mM )$.
On the other hand, by the definition of $h(\partial [m];\cdot )$ and 
functoriality of the skeleton operation the map 
$$
h([m];\pA (x_{\cdot}))\rightarrt \pA 
$$
yields a map
$$
h(\partial [m];\pA (x_{\cdot}))\rightarrt \sk _{m-1}(\pA ).
$$
These two maps agree on $h(\partial [m];\dgt (\pA ; x_{\cdot}))$ so they give a map 
defined on the coproduct 
$$
h([m],\partial [m]; \dgt (\pA ; x_{\cdot})\stackrel{\delta (\pA ; x_{\cdot})}{\longrightarrow} 
\pA (x_{\cdot}) ) 
$$
$$
= 
h([m]; \dgt (\pA ; x_{\cdot}))\cup ^{h(\partial [m];\dgt (\pA ; x_{\cdot}))}  h(\partial [m];\pA (x_{\cdot})),
$$
giving the top map in the commutative square
$$
\begin{diagram}
h([m],\partial [m];\dgt (\pA ; x_{\cdot})\rightarrow
\pA (x_{\cdot})) & \rightarr & \sk _{m-1}(\pA )\\
\downarr^{ R([m];\delta (\pA ; x_{\cdot}))} && \downarr \\
h([m];\pA (x_{\cdot}))& \rightarr & \sk _m(\pA ).
\end{diagram}
$$
The map on the bottom is given by adjunction from the equality
$$
\pA (x_0,\ldots , x_m)= \sk _m(\pA )(x_0,\ldots , x_m). 
$$
Putting these together over all sequences $x_0,\ldots , x_m$ of length $m$ 
we get an expression for $\sk _m(\pA )$.

\begin{proposition}
\label{expression}
For any $m$ we have an expression of $\sk _m(\pA )$ as a pushout of 
$\sk _{m-1}(\pA )$ by copies of the standard maps $R([m];\cdot )$
indexed by sequences $x_{\cdot}=(x_0,\ldots , x_m)$:
$$
\sk _m(\pA ) = \sk _{m-1}(\pA )\cup ^{\coprod _{x_{\cdot}}
h([k],\partial [k]; \delta (\pA ; x_{\cdot}))}
\coprod _{x_{\cdot}}h([m],\pA (x_{\cdot})).
$$
This is a pushout, and uses coproducts, in the category $\precat (\mM )$. 
\end{proposition}
\begin{proof}
This is a classical fact about simplicial objects.
\end{proof}

\section{Reedy cofibrations}
\label{sec-Reedycofib}

Consider a map $f:\pA \rightarrt \pB $ in $\precat (\mM )$, giving for
each $k$ a map on skeleta $\sk _k(\pA )\rightarrt \sk _k(\pB )$. Define the {\em relative
skeleton of $f$}
$$
\pA \cup ^{\sk _k(\pA )}\sk _k(\pB )\rightarrt^{\sk _k^{\rm rel}(f)} \pB .
$$
We say that $f$ is a {\em Reedy cofibration} if the relative skeleton maps are 
injective cofibrations for every $k$. 

\begin{lemma}
The class of Reedy cofibrations is closed under pushout, transfinite composition, and retracts.
\end{lemma}
\begin{proof}
The skeleton operation is given by a pushforward which is a kind of colimit, 
so it commutes with colimits over connected categories which are computed levelwise. 
The relative skeleton map of a retract is again a retract so closure under
retracts comes from the same property in $\mM$ levelwise. 
\end{proof} 

\begin{theorem}
\label{reedycriterion}
Suppose $f:\pA \rightarrt \pB $ is map such that $\Ob (f)$ is 
injective and view $\Ob (\pA )$ as a subset of $\Ob (\pB )$. 
The following conditions are equivalent:
\newline
(a)---$f$ is a Reedy cofibration;
\newline
(b)---for any $m\geq 1$ the map 
\begin{equation}
\label{reedymapb}
\sk _m (\pA )\cup ^{\sk_{m-1}(\pA )}\sk _{m-1}(\pB )\rightarrt \sk _m (\pB )
\end{equation}
is an injective cofibration; 
\newline
(c)---for 
any sequence $(x_0,\ldots , x_p)$ of objects in $\Ob (\pA )$, the map
\begin{equation}
\label{reedymapc1}
\pA (x_0,\ldots , x_p)\cup ^{\dgt (\pA ;x_0,\ldots , x_p)}\dgt (\pB ;x_0,\ldots , x_p) 
\rightarrt \pB (x_0,\ldots , x_p)
\end{equation}
is a cofibration in $\mM$, and for any sequence $(y_0,\ldots , y_p)$ 
of objects in $\Ob (\pB )$ not all in $\Ob (\pA )$,
the map $\dgt (\pB ;y_0,\ldots , y_p) \rightarrt \pB (y_0,\ldots , y_p)$ 
is a cofibration in $\mM$;
\newline
(d)---letting $X:=\Ob (\pB )$,
the map $\Ob (f)_!(\pA )\rightarrt \pB $ is Reedy cofibrant in 
the model structure $\precat _{\rm Reedy}(X,\mM )$ of Theorem \ref{reedyprecat}. 
\end{theorem}
\begin{proof}
First note that (a) implies (c), indeed (c) is the statement of (a) 
for $k=p-1$ at the sequence of objects $(x_0,\ldots , x_p)$ 
or $(y_0,\ldots , y_p)$. 

Next we show that (b) implies the following more general statement: 
for any $0\leq n \leq m$ the map 
\begin{equation}
\label{reedymapb1}
\sk _m (\pA )\cup ^{\sk_{n}(\pA )}\sk _{n}(\pB )\rightarrt \sk _m (\pB )
\end{equation}
is an injective cofibration. This is tautological for $m=n$. 
Let $n$ be fixed, and suppose we know this
statement for some $m\geq n$. Then 
$$
\sk _{m+1} (\pA )\cup ^{\sk_{n}(\pA )}\sk _{n}(\pB ) = 
\sk _{m+1}(\pA )\cup ^{\sk _m(\pA )}(\sk _m (\pA )\cup ^{\sk_{n}(\pA )}\sk _{n}(\pB )).
$$
Injective cofibrations are stable under pushout (Lemma \ref{injcofibpushout}), 
and our inductive hypothesis says that
$$
\sk _m (\pA )\cup ^{\sk_{n}(\pA )}\sk _{n}(\pB )\rightarrt \sk _m(\pB )
$$
is an injective cofibration. 
Take the pushout of this by $\sk _{m+1}(\pA )$ over $\sk _m(\pA )$ and use the previous 
identification, to get that 
$$
\sk _{m+1} (\pA )\cup ^{\sk_{n}(\pA )}\sk _{n}(\pB )\rightarrt 
\sk _{m+1} (\pA )\cup ^{\sk_{m}(\pA )}\sk _{m}(\pB )
$$
is an injective cofibration. On the other hand condition (b) says that 
$$
\sk _{m+1} (\pA )\cup ^{\sk_{m}(\pA )}\sk _{m}(\pB ) \rightarrt \sk _{m+1} (\pB )
$$
is an injective cofibration, so composing these gives that 
$$
\sk _{m+1} (\pA )\cup ^{\sk_{n}(\pA )}\sk _{n}(\pB )\rightarrt \sk _{m+1}(\pB )
$$
is an injective cofibration. This proves by induction that the maps \eqref{reedymapb1}
are injective cofibrations. 

This statement now implies condition (a), indeed if $(x_0,\ldots , x_p)$ 
is any sequence of objects in $\Ob (\pB )$ 
then for any $m\geq p, \; m\geq k$ we have 
$$
\pB (x_0,\ldots , x_p)= \sk _{m}(\pB )(x_0,\ldots , x_p).
$$ 
The same is true for $\pA $ if the $x_i$ are all in $\Ob (\pA )$. 
Hence 
$$
\pA \cup ^{\sk _k(\pA )}\sk _k(\pB ) (x_0,\ldots , x_p)= \sk _m (\pA )
\cup ^{\sk_{k}(\pA )}\sk _{k}(\pB )(x_0,\ldots , x_p),
$$
so the map 
$$
\pA \cup ^{\sk _k(\pA )}\sk _k(\pB ) (x_0,\ldots , x_p)\rightarrt \pB  (x_0,\ldots , x_p)
$$
is the same as the map 
$$
\sk _m (\pA )\cup ^{\sk_{k}(\pA )}\sk _{k}(\pB )(x_0,\ldots , x_p)\rightarrt \pB  (x_0,\ldots , x_p).
$$
This latter is just the value of \eqref{reedymapb1} on the sequence $(x_0,\ldots , x_p)$
so it is a cofibration in $\mM$, as needed to show the Reedy condition (a).
This shows that (b) implies (a). 

To complete the proof we need to show that (c) implies (b).
For this, use the expression of Proposition \ref{expression} for 
$\sk _m(\pB )$ as a pushout of $\sk _{m-1}(\pB )$ and 
the standard inclusions $R([m], \delta (\pB ;y_0,\ldots , y_p))$ and similarly for $\pA $.
We have a diagram
$$
\begin{diagram}
\sk _{m-1}(\pA ) & \rightarr & \sk _m(\pA ) \\
\downarr && \downarr \\
\sk _{m-1}(\pB ) & \rightarr & \sk _m(\pB )
\end{diagram}
$$
where the top arrow is pushout along the $R([m], \delta (\pA ;x_0,\ldots , x_p))$ for sequences
of objects $(x_0,\ldots ,x_p)$ of $\pA $, and the bottom arrow is pushout along
the $R([m], \delta (\pB ;x_0,\ldots , x_p))$ for sequences
of objects $(x_0,\ldots ,x_p)$ of $\pB $. 
Taking the pushout of the upper left corner of the diagram gives the expression
$$
\sk _m (\pA )\cup ^{\sk_{m-1}(\pA )}\sk _{m-1}(\pB ) = 
$$
$$
\sk _{m-1}(\pB ) 
\cup ^{\coprod _{x_{\cdot}} h([m],\partial [m];\delta (\pA ; x_{\cdot}))}
\coprod _{x_{\cdot}}h([m],\pA (x_{\cdot})).
$$
The coproducts are over sequences 
$x_{\cdot} = (x_0,\ldots , x_m) $ of length $m$ of objects of $\pA $, 
however it can be extended to a coproduct over sequences of
objects of $\pB $ by setting $\pA (x_0,\ldots , x_m):= \emptyset$ as well as 
$\dgt (\pA ; x_0,\ldots , x_m):=\emptyset$ 
if any of the $x_i$ are not in $\Ob (\pA )$.
On the other hand, 
$$
\sk _m (\pB ) =\sk _{m-1}(\pB ) 
\cup ^{\coprod _{x_{\cdot}}\partial h([m],\partial [m];\delta (\pB ;x_{\cdot}))}
\coprod _{x_{\cdot}}h([m],\pB (x_{\cdot})).
$$
Putting these two together, we conclude that 
$$
\sk _m (\pB ) = \left( \sk _m (\pA )\cup ^{\sk_{m-1}(\pA )}\sk _{m-1}(\pB ) \right)
\cup ^{\coprod _{x_{\cdot}} \pC (x_{\cdot})} \coprod _{x_{\cdot}}h([m],\pB (x_{\cdot}))
$$
where
$$
\pC (x_{\cdot}):= h([m],\pA (x_{\cdot}))
\cup ^{ h([m],\partial [m];\delta (\pA ; x_{\cdot}))}  
h([m],\partial [m];\delta (\pB ; x_{\cdot})).
$$
Hence, to prove that the map $\sk _m (\pA )\cup ^{\sk_{m-1}(\pA )}\sk _{m-1}(\pB )\rightarrt \sk _m (\pB )$
is an injective cofibration, using stability of injective cofibrations under pushouts,
it suffices to show that each of the maps 
\begin{equation}
\label{cmap}
\pC (x_0,\ldots , x_m) \rightarrt 
h([m],\pB ( x_0,\ldots , x_m))
\end{equation}
is an injective cofibration. For both sides, the set of objects is now our 
standard set $[m]=\{ \upsilon _0,\ldots , \upsilon _m\}$.
Consider the value on a sequence of objects $(y_0,\ldots , y_p)$ in $[m]$. 
There are several possible cases:
\newline
(i)---If it is a constant sequence, then
$$
h([m],\pA ( x_{\cdot}))(y_{\cdot}) =  h([m],\partial [m];\delta (\pA ; x_{\cdot}))(y_{\cdot})
$$
$$
 = h([m],\partial [m];\delta (\pB ;x_{\cdot}))(y_{\cdot}) = h([m],\pB (x_{\cdot}))(y_{\cdot}) =\ast
$$
and the map \eqref{cmap} is the identity. 
\newline
(ii)---If the sequence is somewhere decreasing i.e. there is some $j$ with $y_j<y_{j-1}$ then 
$$
h([m],\pA ( x_{\cdot}))(y_{\cdot}) =  h([m],\partial [m];\delta (\pA ; x_{\cdot}))(y_{\cdot})
$$
$$
 = h([m],\partial [m];\delta (\pB ;x_{\cdot}))(y_{\cdot}) = h([m],\pB (x_{\cdot}))(y_{\cdot})=\emptyset
$$
and again the map \eqref{cmap} is the identity.
\newline
(iii)---If the sequence is nondecreasing, but
there is some $\upsilon _i$ not contained in
$y_{\cdot}$ then 
$$
h([m],\pA ( x_{\cdot}))(y_{\cdot}) =  
h([m],\partial [m];\delta (\pA ; x_{\cdot}))(y_{\cdot}) = \pA (x_{\cdot}), 
$$
and 
$$
h([m],\pB ( x_{\cdot}))(y_{\cdot}) =  
h([m],\partial [m];\delta (\pB ; x_{\cdot}))(y_{\cdot}) = \pB (x_{\cdot}), 
$$
so in this case 
$$
\pC (x_{\cdot})(y_{\cdot}) = \pA (x_{\cdot})\cup ^{\pA (x_{\cdot})}\pB (x_{\cdot}) = \pB (x_{\cdot})
$$
and once again \eqref{cmap} is the identity.
\newline
(iv)---If the sequence is nondecreasing and surjects onto the full set of objects, then 
$$
h([m],\pA ( x_{\cdot}))(y_{\cdot}) = \pA (x_{\cdot}), \;\;\; 
h([m],\partial [m];\delta (\pA ; x_{\cdot}))(y_{\cdot}) = \dgt (\pA ; x_{\cdot}),
$$
and 
$$
h([m],\pB ( x_{\cdot}))(y_{\cdot}) = \pB (x_{\cdot}), \;\;\; 
h([m],\partial [m];\delta (\pB ; x_{\cdot}))(y_{\cdot}) = \dgt (\pB ; x_{\cdot}),
$$
so 
$$
\pC (x_{\cdot})(y_{\cdot}) = \pA (x_{\cdot})\cup ^{\dgt (\pA ;x_{\cdot})}\dgt (\pB ; x_{\cdot})
$$
so the map
$$
\pC (x_{\cdot})(y_{\cdot})\rightarrt h([m],\pB ( x_{\cdot}))(y_{\cdot}) 
$$
is exactly the map \eqref{reedymapc1} 
$$
\pA (x_{\cdot})\cup ^{\dgt (\pA ;x_{\cdot})}\dgt (\pB ; x_{\cdot})\rightarrt \pB (x_{\cdot})
$$
which is known to be a cofibration in $\mM$ because we are assuming condition (c) of the theorem.
This completes the proof that the map \eqref{reedymapb} 
is an injective cofibration, showing (c)$\Rightarrow$(b).
This completes the proof of the equivalence of (a), (b) and (c). 

For the equivalence with (d), note that $\sk _m(\Ob (f)_!\pA ) = \Ob (f)_!\sk _m(\pA )$
by commutation of pushforwards. Now
$$
\pA \cup ^{\sk _m(\pA )}\pB  = \Ob (f)_!(\pA )\cup ^{\Ob (f)_!\sk _m(\pA )}\pB 
$$
from the definition of colimits in $\precat (\mM )$ (Section \ref{sec-limits}),
so $\pA \rightarrt \pB $ is Reedy cofibrant if and only if $\Ob (f)_!\pA \rightarrt \pB $
is. However, this latter induces an isomorphism on sets of objects, and
for such maps the criterion (c) above is the same as the Reedy condition that
the relative latching maps be cofibrant in $\mM$. This shows that (d) is
equivalent to (a),(b),(c).
\end{proof}

The map \eqref{cmap} occuring above is of the form $R([m],g)$ as shown in the following lemma.

\begin{lemma}
\label{twomaps}
Suppose 
$$
\begin{diagram}
E & \rightarr^{a}& F \\
\downarr^u && \downarr_v \\
U & \rightarr^{b}& V
\end{diagram}
$$
is a diagram in $\mM$. Consider the induced map $g: U\cup ^E F \rightarrt V$. 
Then the two maps 
$$
h([m],U)\cup ^{ h([m],\partial [m];E\stackrel{u}{\rightarrow} U)}  
h([m],\partial [m];F\stackrel{u}{\rightarrow} V)
\rightarrt h([m],V)
$$
and
$$
R([m],g): h([m],\partial [m];g)\rightarrt h([m],V)
$$
are the same. 
\end{lemma}
\begin{proof}
Look at the values on any sequence $y_{\cdot} = (y_0,\ldots , y_p)$ 
of objects in $[m]= \{ \upsilon _0,\ldots , \upsilon _m\}$. There are several possibilities:
\newline
---if the sequence is constant then both maps are $\ast \rightarrt \ast$;
\newline
---if the sequence is anywhere decreasing then both maps are $\emptyset \rightarrt \emptyset$; 
\newline 
---if the sequence is nondecreasing but misses some element $\upsilon _j$, 
then we are in the boundary $\partial [m]$ and the first map 
is
$$
U \cup ^U V \rightarrt V
$$
and the second map is $U\rightarrt V$, these are the same; 
\newline
---if the sequence is nondecreasing and surjects onto $[m]$ then both maps are 
$$
U\cup ^E F \rightarrt V.
$$
We need to point out that these identifications of the maps, and in 
particular of their sources, are functorial in the restriction 
maps for diagrams over $\Delta ^o_{[m]}$, so they give an identification 
$$
h([m],U)\cup ^{ h([m],\partial [m];E\stackrel{u}{\rightarrow} U)}  
h([m],\partial [m];F\stackrel{u}{\rightarrow} V) 
\cong h([m],\partial [m];g)
$$
and both maps from here to $h([m],V)$ are the same. 
\end{proof}

\begin{corollary}
\label{Rreedy}
For an object $\pA \in \precat (\mM )$, the following are equivalent:
\newline
(a)---$\pA $ is Reedy cofibrant;
\newline
(b)---for any $m\geq 1$ the map $\sk _{m-1}(\pA )\rightarrt \sk _m(\pA )$ is an injective cofibration;
\newline
(c)---for any sequence of objects $(x_0,\ldots , x_p)$ the map 
$$
\sk _{p-1}(\pA )(x_0,\ldots , x_p)\rightarrt \pA (x_0,\ldots , x_p)
$$
is a cofibration in $\mM$;
(d)---$\pA $ is a Reedy cofibrant object in $\precat _{\rm Reedy}(X,\mM )$
where $X=\Ob (\pA )$. 
\end{corollary}
\begin{proof}
Apply the previous proposition to the map $\emptyset \rightarrt \pA $. 
Note that $\Ob (\emptyset )=\emptyset$
and $\sk _m(\emptyset )= \emptyset$. Condition (c) here is the part of 
condition (c) of the proposition,
concerning sequences of objects not all coming from the source. 
\end{proof}

\begin{corollary}
\label{RisReedy}
If $f:A\rightarrt B$ is a cofibration in $\mM$, then $R([k]; f)$ (page \pageref{Rdefpage})
is a Reedy cofibration. 
\end{corollary}
\begin{proof}
Recall that 
$$
h([k],\partial [k],f) = h([k],A)\cup ^{\sk_{k-1}h([k],A)}
\sk_{k-1}h([k],B)
$$
and $R([k];f)$ is the map from here to $h([k],B)$. 
If $m\leq k-1$ then $\sk _mh([k],\partial [k],f) =\sk _m([k],B)$
so the map occuring in condition (b) of Theorem \ref{reedycriterion} 
is the identity in this case.
For $m\geq k$, 
$$
\sk _mh([k],\partial [k],f) \cup ^{\sk _{m-1}h([k],\partial [k],f)}
\sk _{m-1}h([k],B) =
$$
$$ 
\sk _mh([k],A)\cup ^{\sk_{m-1}h([k],A)}
\sk_{m-1}h([k],B)
$$
so the map occuring in condition (b) is
$$
\sk _mh([k],A)\cup ^{\sk_{m-1}h([k],A)}
\sk_{m-1}h([k],B)\rightarrt \sk _m(h[k],B).
$$ 
As for the equivalence with condition (c), it suffices to note that
this induces a cofibration over sequences $x_0,\ldots , x_m$ of length $m$
in $\{ \upsilon _0,\ldots , \upsilon _k\}$ (similar to the proof of
Proposition \ref{projreedyinj} below).  
\end{proof}

\begin{corollary}
\label{RfaceReedy}
Many maps between the $h([k],B)$ are Reedy cofibrations due to the previous corollary.
For example if $f:A\rightarrt B$ is a cofibration in $\mM$ then the  map
$$
h([k-1],A)\rightarrt h([k],B)
$$
induced by applying $f$ at one of the faces of the $k$-simplex $[k-1]\subset [k]$,
is a Reedy cofibration.
\end{corollary}
\begin{proof}
Either calculate directly the skeleta, or use the previous corollary inductively.
\end{proof}

The following statement is somewhat similar to giving a set of generators for the 
Reedy cofibrations.
On the one hand it refers to the full class of morphisms of the form $R([m],g)$ 
for cofibrations $g$, while on the other hand
giving a stronger expression without refering to retracts. 

\begin{proposition}
\label{reedypushoutR}
A morphism $f:\pA \rightarrt \pB $ is a Reedy cofibration if and only if it is a 
transfinite composition of a disjoint union with a discrete set,
and then pushouts 
along morphisms of the form $R([m],g)$ for cofibrations $g$ in $\mM$.
\end{proposition}
\begin{proof}
A Reedy cofibration $f:\pA \rightarrt \pB $ can be expressed
as the countable composition of the morphisms $\pA \cup ^{\sk _{m-1}(\pA )}\sk _{m-1}(\pB )  
\rightarrt \pA \cup ^{\sk _m(\pA )}\sk _m(\pB ) $
which are themselves Reedy cofibrations. 
At the start, $\pA \cup ^{\sk _0(\pA )}\sk _0(\pB )$ is the disjoint union of $\pA $ with the 
discrete set $\Ob (\pB )-\Ob (\pA )$.
Then, as we have seen in the proof of Theorem \ref{reedycriterion},
at each stage the morphism in question is obtained by simultaneous pushout along 
morphisms \eqref{cmap} of the form 
\begin{equation}
\label{ourmap}
h([m],U)\cup ^{ h([m],\partial [m];E\stackrel{u}{\rightarrow}U)}  
h([m],\partial [m];F\stackrel{v}{\rightarrow}V) \rightarrt 
h([m],V)
\end{equation}
where from the notations of \eqref{cmap} we put $U:= \pA (x_0,\ldots , x_m)$, 
$V:= \pB (x_0,\ldots , x_m)$,
$E:= \dgt (\pA ; x_0,\ldots , x_m)$ and $F:=\dgt (\pB ; x_0,\ldots , x_m)$, with
$$
u:=\delta (\pA ; x_0,\ldots , x_m) : E\rightarrt U
$$
and 
$$
v:=\delta (\pB ; x_0,\ldots , x_m) : F\rightarrt V.
$$
The maps $u$ and $v$, as well as the maps $U\rightarrt V$ and $E\rightarrt F$ 
all fitting into a commutative square. 
Condition (c) of the theorem says that the map $g:U\cup ^E F\rightarrt V$ is a cofibration. 
By Lemma \ref{twomaps}, the map \eqref{ourmap} is the same as $R([m]; g)$. 
\end{proof}

The the notion of Reedy cofibration is similar to that of projective cofibration, 
in that we can give explicitly
the set of generators.

\begin{lemma}
\label{reedycell}
Suppose $I$ is a set of maps in $\mM$. Let $R(I)\subset \Arr (\precat (\mM ))$ denote the set of all arrows of the form
$R([k]; g)$ for $k\in \nn$ and $g\in I$. 
If $f\in \cell (I)$ and $m\in \nn$ then $R([m]; f)\in \cell (R(I))$.
\end{lemma}
\begin{proof}
Look at the behavior of
$R([k]; f)$ under pushouts and transfinite composition. 
Suppose 
$$
\begin{diagram}
A & \rightarr^{f} & B \\
\downarr^g && \downarr_{ u} \\
C & \rightarr^{v} & P 
\end{diagram}
$$
is a pushout diagram in $\mM$, that is $P=B \cup ^A C$. This induces a diagram 
$$
\begin{diagram}
h([k],\partial [k];g) & \rightarr & h([k],\partial [k];u) \\
\downarr_{R([k];g)} && \downarr_{R([k];u)}\\
h([k]; C) & \rightarrt & h([k];P).
\end{diagram}
$$
We claim that this second diagram is then also a pushout in $\precat (\mM )$. 
In fact it is a diagram in $\precat ([k],\mM )$,
and connected colimits of diagrams in $\precat ([k],\mM )$ are the same as the 
corresponding colimits in $\precat (\mM )$,
also in turn they are the same as the corresponding colimits in 
$\diag (\Delta _{[k]}^o, \mM )$(Section \ref{sec-limits}), so the
pushout of the second diagram can be computed objectwise. 
Then it is easy to see that it is a pushout, using the
explicit description of the values of $h(([k],\partial [k]);\cdot )$ 
and $h([k];\cdot )$.
The conclusion from this discussion, is that any pushout along $R([k]; u)$ 
will also be a pushout along $R([k]; g)$. 

Consider now a transfinite composition: suppose we have a series
$$
\ldots \rightarrt A_i \rightarrt^{f_{i,i+1}} A_{i+1} \rightarrt \ldots 
$$
in $\mM$
indexed by $i\in \beta$ for some ordinal $\beta$. To treat  limit ordinals  
we need also to consider the 
transition maps $f_{i,j}:A_i\rightarrt A_j$ for any $i<j$. Assume that if $j$ 
is a limit ordinal then
$A_j=\colim _{i<j}A_i$, and let $A_{\beta}:= \colim _{i\in \beta} A_i$. Consider the map 
$$
f:A_0\rightarrt A_{\beta}.
$$
We would like to express 
$R([k]; f)$ as a transfinite composition of pushouts along  $R([k]; \cdot )$.
Consider the series
$$
G_i:= h([k],\partial [k];A_i\rightarrow A_{\beta}),
$$
$$
\ldots \rightarrt G_i \rightarrt^{g_{i,i+1}}
G_{i+1} \rightarrt \ldots
$$
with the more general transition maps $g_{i,j}: G_i\rightarrt G_j $.
This is still a transfinite series: if $j$ is a limit ordinal then $G_j= \colim _{i<j}G_i$,
and $G_{\beta}:= \colim _{i\in \beta}G_i$ is equal to 
$h([k],\partial [k];1_{ A_{\beta}})=h([k]; A_{\beta})$.
These can be seen by calculating the colimits objectwise over $\Delta _{[k]}$. 
The map
$$
G_0\rightarrt G_{\beta}
$$
is equal to $R([k]; f)$. Furthermore, $G_{i+1}$ is the pushout of $G_i$ along the map 
$R([k]; f_{i,i+1})$. Thus, $R([k]; A_0\rightarrt A_{\beta})$ is a 
transfinite composition of pushouts along
the $R([k]; f_{i,i+1})$.  
In turn, if $f_{i,i+1}$ is a pushout along an element of $I$ then 
by the discussion of pushouts above,
$R([k]; f_{i,i+1})$ is a pushout along an element of $R(I)$. 
Thus, if our series gives an expression for 
$f$ as an element of $\cell (I)$, then $R([k]; f)$ is seen to be in $\cell (R(I))$. 
\end{proof}

A similar statement is needed for retracts.
Suppose 
$$
\begin{diagram}
& & B \\
& \ruTeXto(2,2)^g & \downarr^p \uparr _s \\
A & \rightarr^{f} & C 
\end{diagram}
$$
is a retract diagram in $\mM$, that is $f=pg$, $g = sf$, and $ps=1$. This gives a diagram
$$
\begin{diagram}
h([k],\partial [k]; g)& \rightarr^{R([k]; g)} &  h([k]; B) \\
\downarr \uparr & & \downarr \uparr \\
h([k],\partial [k]; f) & \rightarr^{R([k]; f)} & h([k]; C) 
\end{diagram}
$$
where the vertical arrows are induced by $p$ and $s$. 

\begin{lemma}
\label{reedyretract}
In the above situation, suppose we are given a pushout diagram 
$$
\begin{diagram}
h([k],\partial [k]; f) & \rightarr & \pU \\
\downarr^{ R([k]; f)} && \downarr  \\
h([k]; C) & \rightarr & \pV .
\end{diagram}
$$
Then $\pU\rightarrt \pV$ is a retract of the pushout of $\pU$ along $R([k]; g)$,
in the category of objects under $\pU$.
\end{lemma}
\begin{proof}
Using the left downward map in the previous diagram, we get a map 
$h([k],\partial [k]; g)\rightarrt \pU$ so we can form the pushout $\pV'$ of
$\pU$ along the map $R([k]; g)$. The vertical maps on the right of the previous diagram
induce maps $\pV\rightarrt \pV'$ and $\pV'\rightarrt \pV$, compatible with the maps from $\pU$,
and the composition $\pV\rightarrt \pV'\rightarrt \pV$ is the identity as desired. 
To see all of these
things, note that in 
$$
\begin{diagram}
& & h([k],\partial [k]; g) \\
& \ldTeXto(2,2) & \downarr \uparr \\
\pU & \leftarr & h([k],\partial [k];f)
\end{diagram}
$$
both triangles, obtained by using the upward or the downward arrows, commute. 
The required statements are then obtained
by functoriality of pushout diagrams using the identity on $\pU$. 
\end{proof}

\begin{corollary}
\label{reedyretractcor}
If $f\in \cof (I)$ then $R([k]; f)\in \cof (R(I))$.
\end{corollary}
\begin{proof}
Write $f$ as a retract of $g\in \cell (I)$, and apply the previous Lemma \ref{reedyretract}
to $R([k]; f)$ seen as a pushout of itself and the identity; thus $R([k]; f)$ is a 
retract of a pushout along
$R([k]; g)$. On the other hand, Lemma \ref{reedycell} shows that $R([k]; g)\in \cell (R(I))$.
Thus, $R([k]; f)$ is  a retract of a map in $\cell (R(I))$, so it is in $\cof (R(I))$.
\end{proof}

\begin{proposition}
\label{reedygenerators}
Fix a generating set $I$ for the cofibrations of $\mM$. Then for any $k$ and 
any $f:A\rightarrt B$ in $I$,
consider the map
$$
R([k];f): h(([k],\partial [k]); A\stackrel{f}{\rightarrow} B)\rightarrt h([k]; B) .
$$
The collection of these for integers $k$ and all $f\in I$, forms a 
generating set $R(I)$ for the class of Reedy cofibrations in 
$\precat (\mM )$. If the elements of $I$ have cofibrant domains, then the
elements of $R(I)$ have Reedy-cofibrant domains. 
\end{proposition}
\begin{proof}
We have seen in Proposition \ref{reedypushoutR} 
that any Reedy cofibration can be written as a 
successive pushout by maps of the form $R([k]; f)$ for
various $k$ and various cofibrations $f$ in $\mM$. As $\cof (R(I))$ is 
closed under pushout and transfinite composition,
is suffices to show that for any cofibration $f$, the map $R([k]; f)$ is in $\cof (R(I))$. 
This is exactly the statement of the previous Corollary \ref{reedyretractcor}.
If furthermore $f$ has cofibrant domain then $R([k]; f)$
will have a Reedy cofibrant domain.
\end{proof}

If $\mM$ is cartesian, then the Reedy cofibrations also satisfy the cartesian property. 
This is one of the
main reasons for introducing the Reedy objects. This is closely related to the
corresponding result for diagrams over a Reedy category, see for
example Barwick \cite{BarwickReedy} and  Berger-Moerdijk \cite{BergerMoerdijk}.
I would like to thank several people including Clemens Berger, 
Clark Barwick, Ieke Moerdijk, and Mark Johnson, 
for replying to a query about this on the topology mailing list. 

\begin{proposition}
\label{reedyproduct}
Suppose $\mM$ is a cartesian model category. 
Consider morphisms $A\rightarrt^{f} B$ and
$U\rightarrt^{g} V$ in $\mM$. Assume that they are cofibrations. 
Then for any $k,m$ the morphism $\xi$ from
$$
\pU:=h([k],\partial [k];f)\times h([m]; V)\cup ^{h([k],\partial [k];f)\times h([m],\partial [m];g)}
h([k]; B)\times h([m],\partial [m];g)
$$
to $\pF:= h([k]; B)\times h([m]; V)$ is a Reedy cofibration.
\end{proposition} 
\begin{proof}
Use the criterion (c) of Theorem \ref{reedycriterion}.
Consider any sequence of objects $z_{\cdot}=((x_0,y_0),\ldots , (x_p,y_p))$
where $x_i\in [k]$ and $y_j\in [m]$.  We first need to calculate 
$\sk _{p-1}(\pU )(z_{\cdot})$ and $\sk _{p-1}(\pF)(z_{\cdot})$.

If either one of the sequences $x_{\cdot}$ or $y_{\cdot}$ is decreasing at any index,
then (using that the product of anything with $\emptyset$ is again $\emptyset$)
$$
\sk _{p-1}(\pU )(z_{\cdot})= \sk _{p-1}(\pF)(z_{\cdot}) = \pU (z_{\cdot})=\pF(z_{\cdot})=\emptyset ,
$$
and the map \eqref{reedymapc1} 
for $\xi$ is an isomorphism hence a cofibration. 
Similarly, if $z_{\cdot}$ is constant then 
$$
\sk _{p-1}(\pU )(z_{\cdot})= \sk _{p-1}(\pF)(z_{\cdot}) = 
\pU (z_{\cdot})=\pF(z_{\cdot}) = \ast ,
$$
so again the map \eqref{reedymapc1} is an isomorphism hence a cofibration.

Thus, we may assume that both sequences $x_{\cdot}$ and $y_{\cdot}$ are nondecreasing and
at least one of them is nonconstant. 

However, if one or the other of the sequences is constant, then the morphism in question
becomes the same as the map \eqref{reedymapc1} for the other side, and we know that $R([k];f)$ or
$R([m]; g)$ are Reedy cofibrations by \ref{Rreedy}.

So, we may now assume that both sequences are nonconstant. 

If $z_{\cdot}$ is strictly increasing, then it has no quotient 
$z_{\cdot}\rightarrt w_{\cdot}$ of length $\leq p-1$,
so 
$$
\sk _{p-1}(\pF )(z_{\cdot}) = \sk_{p-1}(\pU )(z_{\cdot})=\emptyset .
$$
So in this case the map \eqref{reedymapc1} is just the map 
$$
\pU  (z_{\cdot})\rightarrt \pF(z_{\cdot}).
$$
First note that $\pF(z_{\cdot})=B\times V$. The calculation of $\pU (z_{\cdot})$ breaks into several
cases. If $x_{\cdot}$ lies in $\partial [k]$ (i.e. it misses at least one object of $[k]$)
and $y_{\cdot}$ lies in $\partial [m]$ then  
$$
\pU  (z_{\cdot}) = B\times V \cup ^{B\times V} B\times V = B\times V,
$$
so the map from here to $\pF(z_{\cdot})$ is an isomorphism hence a cofibration. If 
$x_{\cdot}$ lies in $\partial [k]$ but $y_{\cdot}$ surjects onto $[m]$ then
$$
\pU  (z_{\cdot}) = B\times V \cup ^{B\times U} B\times U = B\times V,
$$
so again the map to $\pF(z_{\cdot})$ is an isomorphism. Similarly, if 
$x_{\cdot}$ surjects onto $[k]$
but $y_{\cdot}$ lies in $\partial [m]$, then 
$$
\pU  (z_{\cdot}) = A\times V \cup ^{A\times V} B\times V = B\times V,
$$
and the map to $\pF (z_{\cdot})$ is an isomorphism. Finally, suppose that 
$x_{\cdot}$ surjects onto $[k]$
and $y_{\cdot}$ surjects onto $[m]$. Then 
$$
\pU  (z_{\cdot}) = A\times V \cup ^{A\times U} B\times U,
$$
so the map from here to $B\times V$ is a cofibration by the 
cartesian axiom for $\mM$ and the assumption that $f$ and $g$ were
cofibrations of $\mM$. 
This completes the proof that the map \eqref{reedymapc1} for $\xi$ is a 
cofibration in the case where 
the sequence $z_{\cdot}$ is strictly increasing.

Assume therefore that $z_{\cdot}$ is not strictly increasing, i.e. it has at 
least one adjacent pair of objects which are equal. 
We have 
$$
\sk _{p-1}(\pF )(z_{\cdot}) = \dgt (\pF ; z_{\cdot}) = 
\colim _{z_{\cdot}\rightarrt w_{\cdot}}\pF (w_{\cdot})
$$
$$
\sk _{p-1}(\pU  )(z_{\cdot}) = \dgt (\pU  ;z_{\cdot}) = 
\colim _{z_{\cdot}\rightarrt w_{\cdot}}\pU  (w_{\cdot})
$$
where the colimits are taken over surjective maps 
$z_{\cdot}\rightarrt w_{\cdot}$ such that $w_{\cdot}$ has length less than or
equal to $p-1$, see Lemma \ref{degeneratesub}. 
The category of quotients $z_{\cdot}\rightarrt w_{\cdot}$ 
of length $q\leq p-1$ is nonempty, because we are assuming that 
$z_{\cdot}$ is not strictly increasing. The opposite category of
this category of quotients has an initial object, corresponding to the 
quotient of minimal length obtained by identifying
all adjacent equal objects. 

The diagram which to 
$z_{\cdot}\rightarrt w_{\cdot}$ associates $\pF (w_{\cdot})$ is constant, 
taking values $B\times V$. This uses the
definitions of $h([k]; B)$ and $h([m]; V)$ and the fact that both sequences 
$x_{\cdot}$ and $y_{\cdot}$ are nondecreasing and nonconstant. Hence, 
$$
\sk _{p-1}(\pF )(z_{\cdot})=\colim _{z_{\cdot}\rightarrt 
w_{\cdot}}B\times V = B\times V = \pF (z_{\cdot})
$$
since the colimit of a constant diagram over a category with 
initial object, is equal to the constant value
of the diagram. 

Next, look at 
$$
\sk _{p-1}(\pU  )(z_{\cdot}) = \dgt (\pU  ; z_{\cdot}) = 
\colim _{z_{\cdot}\rightarrt w_{\cdot}}\pU  (w_{\cdot}).
$$
Suppose given a quotient $z_{\cdot}\rightarrt w_{\cdot}$ of 
length $q\leq p-1$. Write $w_{\cdot} = ((r_0,s_0),\ldots , (r_q,s_q))$.
Note that $r_{\cdot}$ is a quotient sequence of $x_{\cdot}$ and $s_{\cdot}$ 
is a quotient sequence of $y_{\cdot}$. The question of whether
the sequence of first elements $r_{\cdot}$ lies in $\partial [k]$ or $[k]$, or 
whether the sequence of second elements $s_{\cdot}$
lies in $\partial [m]$ or
$[m]$, is independent of the choice of quotients and depends only on 
$x_{\cdot}$ or $y_{\cdot}$. Hence, 
the values  $h([k],\partial [k]; f)(r_{\cdot})$ and $h([m],\partial [m]; g)(s_{\cdot})$ 
are independent of the choice of
$w_{\cdot}$, and the colimit defining $\sk _{p-1}(\pU  )(z_{\cdot})$ is equal to its 
constant value on any of the objects $w_{\cdot}$.
This breaks into exactly the same cases as considered previously, and by the same 
reasoning we see that 
$$
\sk _{p-1}(\pU  )(z_{\cdot}) = \pU  (z_{\cdot}).
$$
Now the map \eqref{reedymapc1} for $\xi$ is written as
$$
\pU  (z_{\cdot})\cup ^{\sk _{p-1}(\pU  )(z_{\cdot})} 
\sk _{p-1}(\pF )(z_{\cdot}) \rightarrt \pF (z_{\cdot}),
$$
but in view of the identifications given above this map is just the 
identity of $\pF (z_{\cdot})=B\times V$
so it is a cofibration. This completes the proof of the proposition. 
\end{proof}

\begin{corollary}
\label{reedycofcart}
Suppose $f:\pA\rightarrt \pB $ and $g:\pU\rightarrt \pV$ are Reedy 
cofibrations in $\precat (\mM )$. Then
$$
\pA \times \pV \cup ^{\pA \times \pU}\pB \times \pU \rightarrt \pB \times \pV
$$
is a Reedy cofibration in $\precat (\mM )$.
\end{corollary}
\begin{proof}
Both $f$ and $g$ may be expressed as transfinite compositions of 
pushouts along elementary Reedy cofibrations of the form $R([k], h)$.
The previous proposition gives the cartesian property for these. Since 
Reedy cofibrations are closed under pushout and transfinite
composition, we get the cartesian property for any $f$ and $g$. 
\end{proof}

One of the main steps in our proof will be to give the same property for 
trivial Reedy cofibrations, in Chapter \ref{product1} below.

\section{Relationship between the classes of cofibrations}
\label{sec-reedyrelation}

\begin{proposition}
\label{projreedyinj}
A projective cofibration is a Reedy cofibration, and a Reedy cofibration is 
an injective cofibration.
\end{proposition}
\begin{proof}
Starting from generators for the projective cofibrations,
if $f:A\rightarrt B$ is a cofibration in $\mM$ then $h([k],f)$ is a
Reedy cofibration. Indeed, the set of objects is $\{ \upsilon _0,\ldots , \upsilon _k\}$.
Suppose given a sequence of the form $x_{\cdot}=(\upsilon _{i_0},\ldots , \upsilon _{i_m})$.
If the sequence is increasing and nonconstant, then the same is true of
any surjective image, and we get 
$$
h([k],A)(x_{\cdot}) = A, \;\; \sk _m h([k],A)(x_{\cdot}) = A \mbox{ or } \emptyset
$$
with $\emptyset$ occuring if there are no surjections to sequences of length $\leq m$. 
Similarly for $B$. The relative skeleton map for $h([k],f)$ at $x_{\cdot}$ 
is either $f$
or the identity of
$B$ in this case. If the sequence is constant, then the same is true of any surjective
image and the relative skeleton map is the identity of $\ast$. If the
sequence is anywhere strictly decreasing, again the same is true of any surjective image
and the relative skeleton map is the identity of $\emptyset$. Thus, 
$h([k],f)$ is a Reedy cofibration. 

Using $m=0$ in the definition of Reedy cofibrations, we get that they are
levelwise cofibrations. 
\end{proof}

\begin{proposition}
\label{reedyinjective}
Suppose $\mM$ is a presheaf category and the cofibrations are the monomorphisms of
$\mM$. Then the Reedy and injective cofibrations of $\precat (\mM )$ coincide, and
if $I$ is a generating set of cofibrations for $\mM$ then the set $R(I)$ 
consisting of the
$R([k]; f)$ for $f\in I$,
is a generating set of cofibrations for the injective cofibrations.
\end{proposition}
\begin{proof}
If $\mM  =\presh (\Phi )$ we can verify the cofibrant property of the
relative skeleton maps levelwise over $\Phi$. It then reduces to the classical
statement that the skeleton of a simplicial set is a simplicial subset. 
\end{proof}

\begin{theorem}
\label{injglobalwe}
Suppose a map $f:\pA \rightarrt \pB $ in $\precat (\mM )$ satisfies the 
right lifting property with respect to the class
of projective (resp. injective, Reedy) cofibrations. Then $f$ is a global weak equivalence.
\end{theorem}
\begin{proof}
It suffices to treat the case of projective cofibrations, since the 
other ones contain this class.

Recall that any set $X\in \Sets$ 
corresponds to a discrete precategory 
$\disc (X)\in \precat (\mM )$ whose object set is $X$ itself,
and whose morphism objects are defined by
$$
\disc (X)(x_0,\ldots , x_n) = \left\{ \begin{array}{ll} 
\ast & \mbox{if } x_0= \cdots = x_n \\
\emptyset & \mbox{otherwise.}
\end{array}
\right.
$$

Included among the projective cofibrations 
is $\emptyset \rightarrt \disc (\{ x\} )$. 
The morphism $f:\pA \rightarrt \pB $ in $\precat (\mM )$ satisfies the right 
lifting property with respect to
$\emptyset \rightarrt \{ x\}$, if and only if $\Ob (\pA )\rightarrt \Ob (\pB )$ is surjective. 
In particular $f$ is essentially surjective. 

Next, suppose $g$ is a generating cofibration of $\mM$.
If $f$ satisfies the right lifting property with respect to a given $h([k], g)$, 
Lemma \ref{universalh}
implies that for any $[k]$-sequence of objects $x_0,\ldots , x_k\in \Ob (\pA )$, the map
$$
\pA (x_0,\ldots , x_k)\rightarrt \pB (f(x_0),\ldots , f(x_k))
$$
satisfies the right lifting property with respect to $g$. If $f$ satisfies the 
right lifting property with
respect to the generators of projective cofibrations of $\precat (\mM )$ which is to say
all the $h([k], g)$ as $g$ runs over a generating set of cofibrations of $\mM$ 
and $k$ is any positive integer,
it follows that 
for any $x_0,\ldots , x_k\in \Ob (\pA )$, the map
$$
\pA (x_0,\ldots , x_k)\rightarrt \pB (f(x_0),\ldots , f(x_k))
$$
is a trivial fibration in $\mM$---in particular it is a weak equivalence. 
This shows that $f$ is fully faithful. 
\end{proof}


\chapter{Calculus of generators and relations}
\label{genrel1}

In this chapter we look more closely at the specific calculus of generators and relations 
corresponding to the direct left Bousfield localization of $\precat (X; \mM )$ discussed in Chapter \ref{weakenr1}.
Throughout, the model category $\mM$ is assumed to be tractable left proper and cartesian.

\section{The $\Upsilon$ precategories}
\label{sec-upsilon}
 
Recall that $\emptyset$ denotes the initial object and $\ast$ the coinitial object of $\mM$. 
We introduce some $\mM$-enriched precategories $\Upsilon _k(B_1,\ldots , B_k)$ via an adjunction. 

If $X$ is a set, it may be considered as a discrete precategory $\disc (X)$
with object set itself, and morphism
objects  $\disc (X)(x_0,\ldots , x_n):= \emptyset$ when some $x_i\neq x_j$ but 
$\disc (X)(x_0,\ldots , x_n):= \ast$ when $x_0= \cdots = x_n$. We can consider
the category $\disc (X)/\precat (\mM )$ of arrows $\disc (X)\rightarrt \pA$ in $\precat (\mM )$. Such an arrow is equivalent
to giving an $\mM$-enriched  precategory $\pA$ together with a map of sets $X\rightarrt \Ob (\pA )$. 

Recall that $[k]$ denotes the ordered set $\{ \upsilon _0,\ldots , \upsilon_k\}$, so we can
form the category
$\disc ([k])/\precat (\mM )$. 

Looking at the morphism objects
between adjacent elements of $[k]$ gives a functor
$$
(\Ee _1,\ldots , \Ee _k): \disc ([k])/\precat (\mM ) \rightarrt \mM \times \cdots \times \mM 
$$
defined by
$$
\left( f:[k]\rightarrt \Ob (\pA )\right) 
\mapsto \left( \pA (f(0), f(1)),\pA (f(1),f(2)), \ldots , \pA (f(k-1),f(k))\right) .
$$
It has a left adjoint, consisting of a functor
$$
\Upsilon _k: \mM ^k \rightarrt \precat (\mM )
$$
together with a natural transformation $\upsilon : \disc ([k])\rightarrt \Upsilon _k(B_1,\ldots , B_k)$
so that the resulting functor 
$$
\mM ^k \rightarrt \disc ([k])/\precat (\mM ), \;\;\; (B_1,\ldots , B_k)\mapsto
\left( \disc ([k])\rightarrt^{\upsilon}\Upsilon _k (B_1,\ldots , B_k) \right)
$$
is left adjoint to $(\Ee _1,\ldots , \Ee _k)$. 

After this abstract introduction, we can describe explicitly $\Upsilon (B_1,\ldots , B_k)$, and the reader could well
skip the above discussion at first reading and consider just the following construction.
The main part of the structure of precategory is given by 
$$
\Upsilon (B_1,\ldots , B_k)(\upsilon _{i-1}, \upsilon _i ) = B_i.
$$
This is extended whenever there is a constant string of points on either side:
$$
\Upsilon (B_1,\ldots , B_k)(\upsilon _{i-1},\ldots , \upsilon _{i-1}, \upsilon _i ,\ldots , \upsilon _i) = B_i.
$$
The unitality condition on the diagram $\Delta ^o_{\{ \upsilon _0,\ldots , \upsilon _k\}} \rightarrt \mM$ implies, by minimality of $\Upsilon _k$, that for $0\leq i \leq k$ we have
$$
\Upsilon (B_1,\ldots , B_k)(\upsilon _{i}, \ldots ,  \upsilon _i ) = \ast .
$$
In all other cases, 
and
$$
\Upsilon (B_1,\ldots , B_k)(x_0, \ldots ,  x_n ) = \emptyset .
$$
The reader is invited to check that the obvious maps turn $\Upsilon (B_1,\ldots , B_k)$ as defined above, into a functor from 
$\Delta ^o_{\{ \upsilon _0,\ldots , \upsilon _k\}}$ to $\mM$, which is unital, in other words it is an element of $\precat (\mM )$. 
The tautological map 
$$
\upsilon : [k]\rightarrt \Ob (\Upsilon (B_1,\ldots , B_k)) = \{\upsilon _0,\ldots , \upsilon _k\} , \;\;\; i\mapsto \upsilon _i
$$
provides $\Upsilon (B_1,\ldots , B_k)$ with a structure of element of $\disc ([k])/\precat (\mM )$. 

The construction is clearly functorial in $(B_1,\ldots , B_k)$, and $\upsilon$ is a natural transformation. 
We get a functor $\Upsilon : \mM ^k \rightarrt \disc ([k])/\precat (\mM )$. 

\begin{lemma}
\label{upsilonadjunction}
The explicitly constructed $\Upsilon$ is left adjoint to the functor $(\Ee _1,\ldots , \Ee _k)$. Furthermore it satisfies
a more global adjunction property:
if $\pR \in \precat (\mM )$ then a morphism 
$\Upsilon (B_1,\ldots , B_k)\rightarrt \pR$ is the same thing as a string of objects $x_0,\ldots , x_k\in \Ob (\pR )$,
and maps $B_i \rightarrt \pR (x_{i-1},x_i)$. 
\end{lemma}
\begin{proof} 
Given a string of objects $x_0,\ldots , x_k\in \Ob (\pR )$,
and maps $B_i \rightarrt \pR (x_{i-1},x_i)$, the functoriality maps for the diagram 
$\pR : \Delta _{\Ob (\pR )}^o \rightarrt \mM$ provide the required maps to define $\Upsilon (B_1,\ldots , B_k)\rightarrt \pR $
and this is inverse to the obvious construction in the other direction.
\end{proof}

\begin{corollary}
\label{upsilonfromh}
There is also an expression as a pushout of the standard representables:
$$
\Upsilon (B_1,\ldots , B_k) = 
h([1],B_1)\cup ^{\upsilon _1} h([1],B_2)\cup ^{\upsilon _2}\cdots \cup ^{\upsilon _{k-1}}
h([1],B_k)
$$
where the maps in the coproducts go alternately from the single point precategories
to the second or first objects of $h([1],-)$.
\end{corollary}
\begin{proof}
The pushout expression satisfies the same universal property as 
given for $\Upsilon$ in the previous lemma.
\end{proof}

The construction $\Upsilon$ sends sequences of cofibrations in $\mM$ to
Reedy cofibrations. 

\begin{lemma}
\label{upsilonreedy}
If $f_i:A_i\rightarrt B_i$ are cofibrations in $\mM$ then the induced map
$$
\Upsilon (A_1,\ldots , A_k)\rightarrt \Upsilon (B_1,\ldots , B_k)
$$
is a Reedy cofibration.
\end{lemma}
\begin{proof}
Put together in a string the cofibrations $R([1],f_i)$ of Corollary \ref{RisReedy}.
\end{proof}

One could form a more complicated Reedy cofibration out of the $f_i$,
for example when $k=2$
$$
\Upsilon (B_1,A_2)\cup ^{\Upsilon (A_1,A_2)}\Upsilon (A_1,B_2)\rightarrt \Upsilon (B_1,B_2)
$$
is a Reedy cofibration. We leave it to the reader to elucidate notation for the
most general possibility.

As a particular case of the $\Upsilon$ construction which will enter into our 
discussion in the next section,
note that for $(B_1,\ldots , B_k)=(B,\ldots , B)$ 
there is a tautological map $\Upsilon (B,\ldots , B)\rightarrt h([k],B)$. In terms of universal properties, if $\pR \in \precat (\mM )$
and we are given a sequence of objects $x_0,\ldots , x_n\in \Ob (\pR )$ and $B\rightarrt \pR (x_0,\ldots , x_n)$, which corresponds to
$h([k],B)\rightarrt \pR $, then
composing with the principal edge maps (i.e. those which make up the Segal map) we get $B\rightarrt \pR (x_{i-1},x_i)$
and this collection corresponds to the map $\Upsilon (B,\ldots , B)\rightarrt \pR $.

\section{Some trivial cofibrations}
\label{sec-trivgen}

There is an important link between the pseudo-generating set used to define global weak equivalences in Chapters \ref{algtheor1} and  \ref{weakenr1},
and the objects $\Upsilon $ defined above.  

A morphism of sets $g:X\rightarrt Y$ induces a functor $\Delta _g^o:\Delta _X^o \rightarrt \Delta _Y^o$.
Consider a sequence of objects $x_0,\ldots , x_n\in X$ and the resulting Segal functor
$P_{(x_0,\ldots , x_n)}:\epsilon (n)\rightarrt \Delta _X^o$. Then the composition with $\Delta _g^o$
is equal to the functor corresponding to the sequence $g(x_0),\ldots , g(x_n)\in Y$:
$$
\Delta _g^o\circ P_{(x_0,\ldots , x_n)}= P_{(g(x_0),\ldots , g(x_n))}:\epsilon (n)\rightarrt \Delta _Y^o.
$$
We obtain a diagram 
$$
\begin{diagram}
\diag (\epsilon (n),\mM ) & \rightarr & \diag (\Delta _X^o,\mM )& \rightarr^{g_!} & \diag (\Delta _Y^o,\mM ) \\
& \rdTeXto(2,2) & \downarr & & \downarr \\
& & \precat (X,\mM )& \rightarr^{g_!} & \precat (Y,\mM )
\end{diagram}
$$
where the vertical arrows are the unitalization operators $U_!$. 

If $f:A \rightarrt B$ is a generating cofibration for $\mM$, 
look at what happens to the map
$$
\varpi _n(f)= (A, B, \ldots , B; f,\ldots , f) \rightarrt^{\rho _n(f)} 
\xi _{0,!}(B) = (B,\ldots , B; 1,\ldots , 1)
$$
in $\diag (\epsilon (n),\mM )$ which was considered in Chapter \ref{algtheor1}. 
The image in $\precat (X,\mM )$ is  the  pseudo-generator 
$$
U_!  P_{(x_0,\ldots , x_n),!}(\rho _n(f)) : U_!  P_{(x_0,\ldots , x_n),!}(\varpi _n (f))\rightarrt U_!  P_{(x_0,\ldots , x_n),!}(\xi _{0,!}(B))
$$
for the Reedy model structure considered in Theorem \ref{reedyprecat}.  When we project to $\precat (Y,\mM )$ using $g_!$ this gives the corresponding pseudo-generator 
in $\precat (Y,\mM )$ for the sequence $g(x_0),\ldots , g(x_n)$:
$$
g_!\left( U_!  P_{(x_0,\ldots , x_n),!}(\rho _n(f)) \right) = U_!  P_{(g(x_0),\ldots , g(x_n)),!}(\rho _n(f)) .
$$

This remark is particularly useful when $X= [n]$ is the universal set having a sequence of objects
$\upsilon _0,\ldots , \upsilon _n\in [n]$. For any sequence $y_0,\ldots , y_n\in Y$ there is a unique map 
$g_{(y_0,\ldots , y_n)}: [n]\rightarrt Y$ sending $\upsilon _i$ to $y_i$, and the standard generator for the cofibration $f$ at the
sequence $(y_0,\ldots , y_n)$ is then expressed as
\begin{equation}
\label{gshriekstd}
U_!  P_{(y_0,\ldots , y_n),!}(\rho _n(f)) = g_{(y_0,\ldots , y_n),!}\left( U_!  P_{(\upsilon _0,\ldots , \upsilon _n),!}(\rho _n(f)) \right) .
\end{equation}

So, in order to understand these generators it suffices to look at the corresponding ones in $\precat ([n],\mM )$.

Using $\zeta _n(f)$ instead of $\rho _n(f)$ gives the pseudo-generators for 
the projective model structure (Theorem \ref{modstrucs}). The pushouts by  
$U_!  P_{(\upsilon _0,\ldots , \upsilon _n),!}(\zeta _n (f))$ will not however have
such a canonical description in view of the additional choice necessary to define $\zeta _n(f)$. The reader is invited to
modify the following discussion accordingly in that case.

\begin{lemma}
\label{pseudogens}
Fix $n$ and suppose $f:A\rightarrt B$ is a cofibration in $\mM $. Define the arrow
$$
\Psi ([n], f):= \left( h([n],A )\cup ^{\Upsilon (A,\ldots , A)}\Upsilon (B,\ldots , B) \rightarrt h([n],B) \right) 
$$
in $\precat ([n],\mM )$.
Then
$$
U_!  P_{(\upsilon _0,\ldots , \upsilon _n),!}(\rho (f)) = \Psi ([n], f) .
$$
\end{lemma}
\begin{proof}
The source of $\rho (f)$ is the object $\varpi (f)= (A, B, \ldots , B; f,\ldots , f)$ in 
$\diag (\epsilon (n),\mM )$. Its image in $\diag (\Delta _{[n]}^o, \mM )$ is the universal object for diagrams $\pR $
with maps $A\rightarrt \pR (\upsilon _0,\ldots , \upsilon _n)$ and $B\rightarrt \pR (\upsilon _{i-1},\upsilon _i)$
for $i=1,\ldots , n$, making commutative diagrams 
$$
\begin{diagram}
A & \rightarr & B \\
\downarr & & \downarr \\
\pR (\upsilon _0,\ldots , \upsilon _n) & \rightarr & \pR (\upsilon _{i-1},\upsilon _i) .
\end{diagram}
$$
After applying $U_!$ the image in $\precat ([n],\mM )$ is again the universal object for unital diagrams $\pR $ as above. 

The collection of maps $B\rightarrt \pR (\upsilon _{i-1},\upsilon _i)$ corresponds by adjunction to a map $\Upsilon (B,\ldots , B)\rightarrt \pR $,
the map $A\rightarrt \pR (\upsilon _0,\ldots , \upsilon _n)$  corresponds to a map $h([n],A)\rightarrt \pR $, and 
the commutative diagrams amount to requiring that they lead to the same map $\Upsilon (A,\ldots , A)\rightarrt \pR $.
Hence, our universal diagram is the pushout
$$
U_!  P_{(\upsilon _0,\ldots , \upsilon _n),!}(\varphi (f))= h([n],A )\cup ^{\Upsilon (A,\ldots , A)}\Upsilon (B,\ldots , B).
$$
Similarly, the target of $\rho (f)$ is the universal object for unital diagrams $S$ with  maps $B\rightarrt S(\upsilon _0,\ldots , \upsilon _n)$,
thus
$$
U_!  P_{(\upsilon _0,\ldots , \upsilon _n),!}(\xi _{0,!}(B)) = h([n],B).
$$
The map between them is the natural one induced by $h([n],A)\rightarrt h([n],B)$ and 
$\Upsilon (B,\ldots , B)\rightarrt h([n],B)$. 
\end{proof}

\begin{corollary}
If $y_0,\ldots , y_n$ is any sequence of objects in $Y$ and $f$ is a cofibration in $\mM$ then
$U_!  P_{(y_0,\ldots , y_n),!}(\rho (f))$ is obtained by applying  
$g_{(y_0,\ldots , y_n),!} $ to the map of the lemma. 
\end{corollary}
\begin{proof}
Apply \eqref{gshriekstd} to the previous lemma. 
\end{proof}

\begin{lemma}
Suppose $f:A\rightarrt B$ is a  cofibration in $\mM$ and $n\geq 0$. 
Then $\Psi ([n],f)$ is a Reedy trivial cofibration in $\precat ([n],\mM )$ and a globally trivial Reedy cofibration in $\precat (\mM )$. 
\end{lemma}
\begin{proof}
This follows from Theorem \ref{reedyprecat}. 
\end{proof}

For brevity, the source of $\Psi ([n],f)$ will be denoted
$$
\src \Psi ([n],f) := h([n],A )\cup ^{\Upsilon (A,\ldots , A)}\Upsilon (B,\ldots , B).
$$
The target is $\targ \Psi ([n],f)=h([n],B)$. 

\begin{proposition}
Suppose $f:A\rightarrt B$ is a  cofibration in $\mM$, and suppose $\pR \in \precat (\mM )$. 
To give a map from the source of $\Psi ([n],f)$ to $\pR $ is the same as to give a sequence of objects $x_0,\ldots , x_n\in \Ob (\pR )$
together with a commutative diagram 
$$
\begin{diagram}
A & \rightarr & B \\
\downarr & & \downarr \\
\pR (x_0,\ldots , x_n)&\rightarr & \pR (x_0,x_1)\times \cdots \times \pR (x_{n-1},x_n)
\end{diagram}
$$
with $f$ on the top and the Segal map on the bottom.
The pushout $\pR \cup ^{\src \Psi ([n],f)}h([n],B)$ computed in $\precat (\mM )$ has a set of objects isomorphic to $\Ob (\pR )$ and
transporting by this identification (which will usually be made tacitly) this pushout is the same as the pushout
of $\pR $ along $U_!  P_{(y_0,\ldots , y_n),!}(\rho (f))$ in the category $\precat (\Ob (\pR ),\mM )$.
\end{proposition}
\begin{proof}
The description of maps $\src \Psi ([n],f)\rightarrt \pR $ was done in the proof of Lemma \ref{pseudogens}. Since $\Psi ([n],f)$ is an
isomorphism on sets of objects, pushout along it preserves the set of objects up to canonical isomorphism. 

In general if $g:X\rightarrt Y$ is a map of sets,  if $\pT \rightarrt \pT '$ is a morphism in $\precat (X,\mM )$ and $\pR \in \precat (Y,\mM )$
with a map $g_!(\pT )\rightarrt \pR $, then the pushout $\pR \cup ^\pT  \pT '$ in $\precat (\mM )$ corresponds (under the canonical isomorphism between
$Y$ and the pushout of $Y$ along $1_X$) to the pushout $\pR \cup ^{g_!(\pT )} g_!(\pT ')$ in $\precat (Y,\mM )$. Indeed they both satisfy the
same universal property in $\precat (Y, \mM )$ because for any $\pR '\in \precat (Y,\mM )$ a map $g_!(\pT )\rightarrt \pR '$ is the same thing as a map
$\pT \rightarrt \pR '$ in $\precat (\mM )$ inducing $g$ on sets of objects; and the same for $\pT '$. 

Apply this general fact to pushouts along the map $\Psi ([n],f)$ to get the last statement of the proposition. 
\end{proof}

\begin{corollary}
\label{psitoseg}
If $\pA  \in \precat (\mM )$ then there is a global trivial Reedy cofibration $\pA  \rightarrt \Seg (\pA  )$ obtained as a transfinite 
composition of pushouts along morphisms either of the form $\Psi ([n],f)$, or levelwise trivial Reedy cofibrations, 
such that $\Seg (\pA  )$ satisfies the Segal condition. 
\end{corollary}
\begin{proof}
If $X=\Ob (\pA  )$ then there is a transfinite composition of pushouts along elements of the standard generating set for the direct left Bousfield
localization considered in Chapters \ref{algtheor1} \ref{weakenr1},
see Theorem \ref{reedyprecat}. The elements of the standard generating set are either 
generating trivial cofibrations for the unital diagram theory without the product condition, i.e. generating levelwise trivial projective cofibrations;
or else maps of the form $U_!  P_{(x_0,\ldots , x_n),!}(\rho (f))$. We have seen above that pushout 
in $\precat (X,\mM )$ along $U_!  P_{(x_0,\ldots , x_n),!}(\rho (f))$ is the same as pushout along $\Psi ([n],f)$ in the global category $\precat (\mM )$.
\end{proof}

Notice from our discussion that any transfinite composition of pushouts along maps as used in the corollary, is a global trivial cofibration,
indeed it is a trivial cofibration in the direct Bousfield localized projective model structure constructed in Chapters \ref{algtheor1} and \ref{weakenr1}.
In particular, if $\pA  \rightarrt \pA  ''$ is some such transfinite composition, then applying the corollary to $\pA  ''$ we obtain another 
such transfinite composition $\pA  '' \rightarrt \Seg (\pA  '')$ such that $\Seg (\pA  '')$ satisfies the Segal condition. In this case all the maps
$$
\pA  \rightarrt \pA  '' \rightarrt \Seg (\pA  '')
$$
are global trivial cofibrations. 

\section{Pushout by isotrivial cofibrations}

One of the main problems is to prove that pushout by a global weak equivalence is again a global weak equivalence. An important first case
is pushout by a global weak equivalence which is an isomorphism on objects. 
In the following discussion we use the generic term ``cofibration'' for either
a projective, Reedy or injective cofibration: the statements come in 
three versions one for each of the model structures of 
Theorem \ref{modstrucs} or Theorem \ref{reedyprecat}.

An {\em isotrivial cofibration} is a cofibration $\pA\rightarrt^f \pB$
(in whichever of the projective, Reedy, or injective structures we are using), 
such that $\Ob (f)$ induces
an isomorphism $\Ob (\pA )\cong \Ob (\pB )$, and $f$ is a global weak equivalence.

\begin{lemma}
\label{isotrivialdef}
A morphism $\pA\rightarrt^f \pB$ is an isotrivial cofibration, if and only if
$\Ob (f)$ is an isomorphism and $\Ob (f)_!\pA \rightarrt ^{f_{\sharp}}\pB$
is a trivial cofibration in the appropriate model structure
of Theorem \ref{modstrucs} or Theorem \ref{reedyprecat} on 
$\precat (X,\mM )$, where $X=\Ob (\pB )$.
\end{lemma}
\begin{proof}
The condition that $\Ob (f)$ is an isomorphism is tautologically necessary so we assume it.
By transport of structure using $\Ob (f)_!$ we may assume 
that $f$ is a morphism in $\precat (X,\mM )$. Considered as a morphism of $\mM$-precategories
$f$ is then automatically essentially surjective. We get the diagram 
$$
\begin{diagram}
\pA & \rightarr & \pB \\
\downarr && \downarr \\
\Seg (\pA )& \rightarr & \Seg (\pB )
\end{diagram}
$$
in $\precat (X,\mM )$, where the vertical arrows are weak equivalences
in the model structures of \ref{modstrucs} or \ref{reedyprecat}.
The fully faithful condition for $f$ is by definition the condition that the
bottom map be a levelwise weak equivalence; but this is equivalent to being
a weak equivalence since the objects satisfy the Segal conditions.
By 3 for 2 this condition is equivalent to the top map being a weak equivalence
in $\precat (X,\mM )$ for the model structures
of \ref{modstrucs} or \ref{reedyprecat}.
\end{proof}

\begin{lemma}
\label{objectwisetrivial}
Suppose $\pA \in \precat (\mM )$ and $f:\pB  \rightarrt \pC $ is an isotrivial cofibration.
Suppose furthermore that $f$ is a levelwise weak equivalence of diagrams, which means
that for any sequence of objects $(x_0,\ldots , x_p)\in \Ob (\pB  )$ the map 
$$
\pB  (x_0,\ldots , x_p )\rightarrt \pC  (f(x_0),\ldots , f(x_p))
$$
is a weak equivalence in $\mM$.  
Suppose given a map $g:\pB  \rightarrt \pA $. Then 
$$
\pA  \rightarrt \pA  \cup ^{\pB } \pC  
$$
also induces an isomorphism on sets of objects, and is a levelwise weak equivalence of diagrams. In particular, it is
an isotrivial cofibration.
\end{lemma}
\begin{proof}
By transport of structure we  may assume that
$\Ob (f)$ is the identity of  $X=\Ob (\pB  ) = \Ob (\pC )$ and think of $\pB  , \pC  \in \precat (X,\mM )$. Let $Y=\Ob (\pA  )$
and denote also by $g:X\rightarrt Y$ the map induced by $g$ on sets of objects. Then 
$$
\pA  \cup ^{\pB } \pC   = \pA  \cup ^{g_!(\pB  )} g_!(\pC  ) \mbox{  in  } \precat (Y, \mM )
$$
but $g_!(\pB  )\rightarrt g_!(\pC  )$ is a levelwise trivial cofibration,
so the pushout is a levelwise trivial cofibration, hence in particular a trivial cofibration
in the model structure
of Theorem \ref{modstrucs} or Theorem \ref{reedyprecat}.
\end{proof}

\begin{theorem}
\label{obisopushout}
Suppose $\pA \in \precat (\mM )$ and $f:\pB  \rightarrt \pC $ is an isotrivial cofibration (in the projective, Reedy or injective structures). 
Suppose given a map $g:\pB  \rightarrt \pA $. Then 
$$
\pA  \rightarrt \pA  \cup ^{\pB } \pC  
$$
is an isotrivial cofibration (in the projective, Reedy or injective structures respectively).
\end{theorem}
\begin{proof}
Using Corollary \ref{psitoseg} let 
$$
\pA  \rightarrt \pA  ', \;\; \pB  \rightarrt \pB  ', \pC  \rightarrt \pC  '
$$
be global trivial cofibrations towards objects which satisfy the Segal property,
obtained by series of pushouts along maps which are either levelwise trivial cofibrations, or 
of the form $\Psi ([n],u)$ where $u$ are generating cofibrations of $\mM$. 
Make the choice for $\pB $ first, and let $\pA  '' = \pA  \cup ^{\pB } \pB  '$ and $\pC  '' = \pC  \cup ^{\pB } \pB  '$. 
Then as in the remark after Corollary \ref{psitoseg}, we can continue with $\pA  '' \rightarrt \pA  '$ and
$\pC '' \rightarrt \pC  '$. That way, there are maps 
$$
\pA  ' \leftarr \pB  ' \rightarrt \pC  '
$$
the second one still being a global trivial cofibration (in any of the projective, Reedy or injective structures), and  
$$
\pA  \cup ^{\pB } \pC  \rightarrt \pA  '\cup ^{\pB  '} \pC  '
$$
is itself obtained by pushout along
maps of the same form. In particular this latter map is also a global trivial cofibration, so by the 3 for 2
property for global weak equivalences, it suffices to show that the map
$$
\pA  ' \rightarrt \pA  '\cup ^{\pB  '} \pC  '
$$
is a global weak equivalence. But now, the map $\pB  '\rightarrt \pC  '$ satisfies the hypothesis of Lemma \ref{objectwisetrivial},
exactly by the fully faithful condition. 
So, by that lemma, the map $\pA  ' \rightarrt \pA  '\cup ^{\pB  '} \pC  '$ is a global trivial cofibration as required.
\end{proof}

Cofibrant pushouts are invariant under 
global weak equivalences inducing isomorphisms on sets of objects.

\begin{lemma}
\label{obisopushoutinvariance}
Suppose given a diagram 
$$
\begin{diagram}
\pA  & \leftarr & \pB  & \rightarr & \pC  \\
\downarr & & \downarr && \downarr \\
\pA  '& \leftarr & \pB  '& \rightarr & \pC  '
\end{diagram}
$$
such that the left horizontal arrows are injective cofibrations, and the
vertical arrows are global weak equivalences inducing isomorphisms on sets of objects.
Then the induced map
$$
\pA  \cup ^{\pB  }\pC  \rightarrt \pA  '\cup ^{\pB  '}\pC  '
$$
is a global weak equivalence inducing an isomorphism on sets of objects.
\end{lemma}
\begin{proof}
Choose a diagram
$$
\begin{diagram}
\pA  & \leftarr & \pB  & \rightarr & \pC  \\
\downarr & & \downarr && \downarr \\
\pA  ''& \leftarr & \pB  ''& \rightarr & \pC  ''
\end{diagram}
$$
such that the vertical maps are injective isotrivial cofibrations, and the precategories along the bottom row satisfy the Segal condition. 
We can do this by first choosing $\pB  '':= \Seg (\pB  )$ which is a transfinite composition
of pushouts along standard morphisms in the direct localizing system $K_{\rm inj}$ 
of Theorem 
\ref{modstrucs}, which in the global context of $\precat (\mM )$ may be interpreted as
the standard maps of Corollary \ref{psitoseg}. Let 
$$
\pA  '' := \Seg (\pA  \cup ^{\pB } \pB  '')
$$
and similarly for $\pC  ''$. That way the map $\pB  ''\rightarrt \pA  ''$ is again an injective cofibration.
Now put
$$
\pB  ^3 := \Seg (\pB  ' \cup ^{\pB } \pB  '' ),
$$
then
$$
\pA  ^3 := \Seg \left( \pB^3 \cup ^{\pB ' \cup ^{\pB}\pB''}\pA '\cup ^{\pA} \pA '' \right) 
$$
and
$$
\pC  ^3 := \Seg \left( \pB^3 \cup ^{\pB ' \cup ^{\pB}\pB''}\pC '\cup ^{\pC} \pC '' \right) .
$$
We have the diagram 
$$
\begin{diagram}
\pB ' \cup ^{\pB} \pB '' & \rightarr & \pA '\cup ^{\pA } (\pA \cup ^{\pB }\pB '') &\rightarr & \pA '\cup ^{\pA} \pA '' \\
\downarr & & \downarr && \downarr \\
\pB ^3 & \rightarr & \pB^3 \cup ^{\pB ' \cup ^{\pB}\pB''} \pA ' \cup ^{\pB}\pB'' & 
\rightarr & \pB^3 \cup ^{\pB ' \cup ^{\pB}\pB''}\pA '\cup ^{\pA} \pA '' \\
\dEqualarr && \downarr && \downarr \\
\pB ^3 &\rightarr & \Seg (\ldots ) & \rightarr & \pA ^3 
\end{diagram}
$$ 
in which the two upper squares are cartesian, the vertical maps
are obtained by pushouts along standard maps of Corollary \ref{psitoseg},
and the top vertical arrows are injective cofibrations.
It follows that $\pB ^3 \rightarrt \pA^3$ is again an injective cofibration. 
There is a similar diagram for $\pC^3$ (but the horizontal
arrows are not necessarily cofibrant). 

Furthermore,
$\pA^3$, $\pB^3$ and $\pC^3$ are pushouts of $\pA '$, $\pB'$ and $\pC'$ respectively
along the standard morphisms of Corollary \ref{psitoseg}. 
They fit into a diagram
$$
\begin{diagram}
\pA  ''& \leftarr & \pB  ''& \rightarr & \pC  ''\\
\downarr & & \downarr && \downarr \\
\pA  ^3& \leftarr & \pB  ^3& \rightarr & \pC  ^3
\end{diagram}
$$
consisting of objects satisfying the Segal condition. The 
vertical arrows are global weak equivalences
inducing isomorphisms on objects. Indeed on the left for example,
pushout along the isotrivial cofibration $\pA  \rightarrt \pA  ''$ 
is again an isotrivial cofibration (Lemma \ref{obisopushout}), from which it follows
that the map
$\pA  '\rightarrt \pA  ^3$ is a global weak equivalence. By 3 for 2, the map $\pA  '' \rightarrt \pA  ^3$
is a global weak equivalence.   

We get that the vertical arrows in the previous diagram
are levelwise weak equivalences of diagrams. 
These are preserved by pushout when one of the maps in the pushout is 
a levelwise cofibration, as is verified levelwise (recall that $\mM$ is assumed to be
left proper).
Therefore the map
$$
\pA  ''\cup ^{\pB  ''}\pC  ''\rightarrt \pA  ^3\cup ^{\pB  ^3}\pC  ^3
$$
is a levelwise weak equivalence of diagrams, so it is a global weak equivalence. 

Now, the map 
$$
\pA  \cup ^{\pB  }\pC  \rightarrt \pA  ''\cup ^{\pB  ''}\pC  ''
$$
is obtained by a transfinite composition of pushouts along the standard maps of 
Corollary \ref{psitoseg}, so it is a global weak equivalence. 
On the other hand, the map 
$$
\pA  '\cup ^{\pB  '}\pC  ' \rightarrt \pA  ^3\cup ^{\pB  ^3}\pC  ^3
$$
is also obtained by a transfinite composition of pushouts along the standard maps, so it is a global weak equivalence. 
By 3 for 2 we conclude that the map 
$$
\pA  \cup ^{\pB  }\pC  \rightarrt \pA  '\cup ^{\pB  '}\pC  '
$$
is a global weak equivalence. 
\end{proof}

A similar argument gives closure under transfinite composition. A cofibrancy condition can be avoided in $\precat (\mM )$ 
under the condition that it could be avoided in $\mM $ already.

\begin{lemma}
\label{globaltranfinite}
The notion of global weak equivalence in $\precat (\mM )$ is closed under transfinite compositions such that the transition
maps are any kind of cofibrations (injective ones being the weakest).
If weak equivalences of $\mM$ are closed under transfinite composition, 
then the notion of global weak equivalence in $\precat (\mM )$
is also closed under transfinite composition. 
\end{lemma}
\begin{proof}
Suppose we are given a transfinite sequence $\{ \pA _i\} _{i\in \alpha}$ indexed by an ordinal $\alpha$, with
continuity at limit ordinals $<\alpha$. Put $\pA _{\alpha}:= \colim _{i<\alpha}\pA _i$. Suppose that the transition
maps $\pA _i\rightarrt \pA _{i+1}$ are global weak equivalences; and suppose either
\newline
(i)---that the transition maps are injective cofibrations,
or else 
\newline
(ii)---that weak equivalences in $\mM$ are closed under transfinite composition. 

We want to prove that $\pA _0\rightarrt \pA _{\alpha}$ is a
global weak equivalence. By induction on $\alpha$ we may assume that this statement is known for all sequences indexed by
strictly smaller ordinals, in particular we get that the transition maps $\pA _i\rightarrt \pA _j$ are global weak equivalences for
all $i<j<\alpha$. 

Set $X_i:= \Ob (\pA _i)$ and $X:= \Ob (\pA _{\alpha})$. Then $X$ is the filtered colimit of the $X_i$. 
Operating by induction on $i$, we will choose a sequence of morphisms $\pA _i\rightarrt \pA '_i$ where
$\pA '_i\in \precat (X_i, \mM )$ is a weak equivalent replacement of $\pA _i$ in the injective model structure
$\precat _{\rm inj}(X_i, \mM )$ for each $i$, such that $\pA '_i$ satisfies the Segal conditions and such that these are compatible in the sense that for 
$i<j<\alpha $ we have a commutative diagram 
$$
\begin{diagram}
\pA _i & \rightarr & \pA _j\\
\downarr & & \downarr \\
\pA '_i & \rightarr & \pA '_j .
\end{diagram}
$$
To make this choice, suppose it is done for all $i<j$. 

If $j$ is a limit ordinal, set $\pA '_j:= \colim _{i<j}\pA '_i$.
By an argument which we will use below for the colimit at $\alpha$ (which we don't repeat here because the notations will be
more comfortable later), the map $\pA _j\rightarrt \pA '_j$ is a weak equivalence in $\precat (X_j, \mM )$ and $\pA '_j$ satisfies the Segal conditions.
This treats the case of a limit ordinal. 

Suppose $j=i+1$ is a successor ordinal. Consider the diagram 
$$
\begin{diagram}
\pA _i & \rightarr & \pA _j\\
\downarr & & \downarr  \\
\pA '_i & \rightarr & P_j
\end{diagram}
$$
where $ \pP_j:= \pA '_i\cup ^{\pA _i}\pA _j$ is defined as the pushout. 
Now, the left vertical map is a weak equivalence, by the inductive hypothesis on the statement we are trying to prove,
and it is also an injective cofibration in $\precat (X_{i}, \mM )$. By Lemma \ref{obisopushout}, 
pushout along a trivial cofibration which induces an isomorphism on sets of objects, is again a trivial cofibration. Thus
$\pA _j\rightarrt \pP_j$ is a trivial cofibration, and we can let $\pP_j\rightarrt \pA '_j$ be a $K$-injective replacement of $\pP_j$,
which will also be a $K$-injective replacement of $\pA _j$ accepting compatible maps from all the $\pA '_i$ for $i<j$. Here $K$ denotes the
pseudo-generating set for $\precat (X_j,\mM )$ considered in Chapter \ref{weakenr1}. 

This completes the choice of the sequence $\pA '_i$, which by construction is continuous at limit ordinals. Put $\pA '_{\alpha}:= \colim _{i<\alpha}\pA '_i$.
As promised above, we show that $\pA '_{\alpha}$ satisfies the Segal conditions and $\pA _{\alpha}\rightarrt \pA '_{\alpha}$ is a weak equivalence 
in $\precat (X,\mM )$.

Suppose given a sequence of objects $x_0,\ldots , x_n\in X$.
It comes from a sequence $x^k_0,\ldots , x^k_n\in X_k$ for some $k<\alpha$ (which depends on the sequence), and denote by 
$x^i_0,\ldots , x^i_n\in X_i$ the image sequence for any $i\geq k$. The colimit morphism objects for $\pA '_{\alpha}$ are
$$
\pA '_{\alpha}(x_0,\ldots , x_n) = \colim _{k\leq i < \alpha} \pA '_i(x^i_0,\ldots , x^i_n).
$$
Commutation of direct products with colimits, condition (DCL) on $\mM$, now tells us that 
$$
\pA '_{\alpha}(x_0,x_1)\times \cdots \times \pA '_{\alpha}(x_{n-1},x_n)
$$
$$
= \colim _{k\leq i < \alpha} \pA '_i(x^i_0,x^i_1)\times \cdots \times \pA '_i(x^i_{n-1}, x^i_n).
$$
Transfinite composition of a sequence of weak equivalences is again a weak equivalence, so the Segal maps for $\pA ' _{\alpha}$ are weak equivalences. 
Thus, $\pA '_{\alpha}$ satisfies the Segal conditions. 

By the construction above, the map $\pA _{\alpha}\rightarrt \pA '_{\alpha}$ is obtained by a transfinite composition of pushouts along 
pseudo-generating trivial cofibrations for the various $\precat (X_i,\mM )$, thus it is a trivial cofibration.

We can now show that $\pA _0\rightarrt \pA _{\alpha}$ is a global weak equivalence. Use $\pA _{\alpha}\rightarrt 
\pA '_{\alpha}$ as replacement satisfying the Segal conditions, so the problem is to show that $\pA '_0\rightarrt \pA '_{\alpha}$
is a weak equivalence. 

It is essentially surjective, indeed given any $x\in X=\Ob (\pA '_{\alpha})$ there is some $i<\alpha$ such that $x$ comes from $x^i\in X_i$.
Hence the isomorphism class of $x$ in the truncation, is in the image of 
$$
\Iso \tau _{\leq 1}(\pA '_i)\rightarrt \Iso \tau _{\leq 1}(\pA '_{\alpha}).
$$
However, since $\pA '_0\rightarrt \pA '_i$ is a global weak equivalence, it is essentially surjective in other words
$$
\Iso \tau _{\leq 1}(\pA '_0)\rightarrt \Iso \tau _{\leq 1}(\pA '_{i}).
$$
is surjective. It follows that the isomorphism class of $x$ is in the image of 
$$
\Iso \tau _{\leq 1}(\pA '_0)\rightarrt \Iso \tau _{\leq 1}(\pA '_{\alpha}).
$$
This shows essential surjectivity. 

To show that the map is fully faithful, let $x_0,\ldots , x_n$ be a sequence of objects in $X_0= \Ob (\pA '_0)$. 
Let $x^i_j$ denote their images in $X_i$ including $i=\alpha$. Now
$$
\pA '_{\alpha}(x^{\alpha}_0,\ldots , x^{\alpha}_n) = \colim _{i<\alpha} \pA '_i(x^i_0,\ldots , x^i_n).
$$
To show full faithfulness we need to show that 
\begin{equation}
\label{zeroalpha}
\pA '_0(x_0,\ldots , x_n) \rightarrt \colim _{i<\alpha} \pA '_i(x^i_0,\ldots , x^i_n)
\end{equation}
is a weak equivalence.

Recall that we are assuming either (i) or (ii) above. In case (i), the transition maps $\pA _i\rightarrt \pA _j$
are injective cofibrations, and by construction of $\pA '_i$ the same is true of the transition maps
$\pA '_i\rightarrt \pA '_j$. Hence the colimit expression above is a sequential colimit whose transition maps
are trivial cofibrations, thus the map \eqref{zeroalpha}
is a trivial cofibration, which shows that $\pA '_0\rightarrt \pA '_{\alpha}$ is fully faithful. 

In case (ii), we just know that the transition maps in the colimit are weak equivalences, but the hypothesis (ii)
says again that the map \eqref{zeroalpha} is a weak equivalence.

In either case, $\pA '_0\rightarrt \pA '_{\alpha}$ is fully faithful.
\end{proof}

\section{An elementary generation step $\Gen$}
\label{sec-gen}

Corollary  \ref{psitoseg} looks at the replacement $\pA \rightarrt \Seg (\pA )$  from the
point of view of direct left Bousfield localizing systems 
of Chapters \ref{algtheor1} and \ref{weakenr1}. It will also be useful to have an approach
which is more homotopically canonical, in other words some kind
of process which looks canonical when viewed in the homotopy category 
$\Ho (\diag (\Delta ^o_X/X,\mM ))$ of unital diagrams up to levelwise weak equivalences.
The full process will be broken up into elementary steps denoted $\pA\rightarrt \Gen (\pA , q)$. These should be thought of as ``calculating generators and relations 
at the Segal map corresponding to $q$''. 

Fix a set of objects $X$. Recall that the underlying category used  is $\Phi = \Delta _X^o$;
the subcategory on which the unitality condition will be imposed is $\Phi _0 = \{ (x_0)\} _{x_0\in X}$, and the set $Q$ used to determine the algebraic theory of Segal categories over $X$,
is just the object set of $\Delta _X$ that is the set of sequences $(x_0,\ldots , x_n)$ of $x_i\in X$.
For $q=(x_0,\ldots , x_n)$, we have $n_q:=n$; the functor $P_q:\epsilon (n)\rightarrt \Delta _X^o$ 
sends $\xi _0$ to $q= (x_0,\ldots , x_n)$ and
for $1\leq j\leq n$ it sends $\xi _j$ to the adjacent pair $(x_{j-1},x_j)$. The structural
maps for $P_q$ come from inclusions of each adjacent pair in the full sequence; they are the maps
which together make up the Segal maps. 

An $\pA \in \precat (X;\mM )$ is a functor $\pA : \Delta _X^o\rightarrt \mM$ such that $\pA (x)= \ast$ for any
single element sequence $(x)$. 

For any $q= (x_0,\ldots , x_n)\in Q=\Ob (\Delta _X^o)$, the functor 
$$
U_!P_{q,!}: \diag (\epsilon (n), \mM )\rightarrt \precat (X,\mM )
$$
sends injective cofibrations to Reedy cofibrations.
Recall that it sends the standard $\rho _n(f)$ to the cofibrations $\Psi ([n],f)$
see Corollary \ref{psitoseg}. 
Abbreviate the right adjoint $P_q^{\ast}U^{\ast}$ by just $P_q^{\ast}$.

For $\pA \in \precat (X,\mM )$ and
$q= (x_0,\ldots , x_n)$ we have considered in Section \ref{sec-gendef} the cofibration 
$\pA \rightarrt \Gen (\pA ,q)$ 
which depends on a choice of factorization
$$
P_q^{\ast}(\pA )_0 \rightarrt E \rightarrt P_q^{\ast}(\pA )_1\times \cdots \times P_q^{\ast}(\pA )_n.
$$
Note that $P_q^{\ast}(\pA )_0 = \pA (x_0,\ldots , x_n)$ while $P_q^{\ast}(\pA )_j= \pA (x_{j-1},x_j)$ for
$1\leq j\leq n$. Thus, the factorization we need to choose may be written as
$$
\pA (x_0,\ldots , x_n) \rightarrt^{e} E \rightarrt^{(p_1,\ldots , p_n)}
\pA (x_0,x_1)\times \cdots \times \pA (x_{n-1},x_n).
$$
Given such a factorization, which will generally be chosen so that $e$ is a
cofibration and $(p_1,\ldots , p_n)$ is a weak equivalence, we get a map
$$
\gen (\pA ,q): \pA \rightarrt \Gen (\pA ,q)[E,e, p_1,\ldots , p_n] .
$$
If no confusion arises, we abbreviate the right side to $\Gen (\pA ,q)$. 
In the notation of the previous section, 
\begin{equation}
\label{gendef}
\Gen (\pA ,q)[E,e, p_1,\ldots , p_n] = \pA \cup ^{\src \Psi ([n],e)}h([n],E).
\end{equation}

\begin{lemma}
Let $\pA \in \precat (X,\mM )$ and
$q= (x_0,\ldots , x_n)$, and suppose given a choice of factorization $E,e, p_1,\ldots , p_n$ as above. 
Suppose $e$ is a cofibration and $(p_1,\ldots , p_n): E\rightarrt \pA (x_0,x_1)\times \cdots \times \pA (x_{n-1},x_n)$
is a weak equivalence in $\mM$. Then $\gen (\pA ,q)$ is a trivial cofibration from $\pA $ to 
$\Gen (\pA ,q)=\Gen (\pA ,q)[E,e, p_1,\ldots , p_n]$ in the model structure $\precat _{\rm Reedy}(X; \mM )$ constructed in Theorem \ref{reedyprecat}.
\end{lemma}
\begin{proof}
Indeed $\gen (\pA ,q)$ is a pushout along $\Psi ([n],e)$ by \eqref{gendef}. 
\end{proof}

We can describe explicitly the structure of $\Gen (\pA ;x_0,\ldots , x_n)$ depending on the choice of $E,e, p_1,\ldots , p_n$. 
For any sequence $(z_0,\ldots , z_m)$, let 
$$
\Delta _X^{NS}((z_0,\ldots , z_m), (x_0,\ldots , x_n))\subset \Delta _X((z_0,\ldots , z_m), (x_0,\ldots , x_n))
$$
be the subset of maps $\phi : (z_0, \ldots , z_m)\rightarrt (x_0,\ldots ,x_m)$ in $\Delta _X$
which don't factor through any of the adjacent pairs $(x_{j-1},x_j)$. 
Abbreviate this by $\Delta _X^{NS}(z_{\cdot}, x_{\cdot})$ when convenient. 
Then
$\Gen (\pA ;x_0,\ldots , x_n)(z_0,\ldots ,z_m)$ is a pushout in the diagram
$$
\begin{diagram}
\coprod _{\phi\in  \Delta _X^{NS}(z_{\cdot}, x_{\cdot})} \pA (x_0,\ldots , x_n) & \rightarr &
\coprod _{\phi\in  \Delta _X^{NS}(z_{\cdot}, x_{\cdot})} E \\
\downarr & & \downarr \\
\pA (z_0,\ldots , z_m) & \rightarr & \Gen (\pA ;x_0,\ldots , x_n)(z_0,\ldots ,z_m).
\end{diagram}
$$
In order to define the functoriality it is convenient to rewrite this as follows: for any object $B\in \mM$ 
and any $\phi \in \Delta _X((z_0,\ldots , z_m), (x_0,\ldots , x_n))$ denote
by ${\bf alt}(\phi, B)$ either: $B$ if $\phi \in \Delta _X^{NS}(z_{\cdot}, x_{\cdot})$;
or $\ast$ if $\phi$ factors through a
singleton
$$
(z_0,\ldots , z_m)\rightarrt (x_j)\rightarrt (x_0,\ldots ,x_n);
$$
or $\pA (x_{j-1},x_j)$ if $\phi$ factors through a map 
$$
(z_0,\ldots , z_m)\rightarrt (x_{j-1},x_j)\rightarrt (x_0,\ldots ,x_n)
$$
but not through a singleton (in which case the choice of $j$ is unique). Now, for either $B=\pA (x_0,\ldots , x_n)$ or $B=E$,
we have structural maps $B\rightarrt \pA (x_{j-1},x_j)$. Using these, 
given a map $\psi : (y_0,\ldots , y_p)\rightarrt (z_0,\ldots , z_m)$
then for any $\phi \in \Delta _X((z_0,\ldots , z_m), (x_0,\ldots , x_n))$ we get a map
$$
{\bf alt}(\phi, B) \rightarrt \nocom {\bf alt} (\phi \psi , B)
$$
and these are compatible with composition in $\psi$. Indeed, if $\phi$ factors through an adjacent pair or a singleton,
then so does $\phi \psi$. If $\phi \psi$ factors through a singleton then the map from 
${\bf alt}(\phi, B) $ is the unique map to ${\bf alt} (\phi \psi , B)=\ast$. If $\phi \psi$ factors through an
adjacent pair but $\phi$ didn't factor through an adjacent pair, then the map 
$$
{\bf alt}(\phi, B) = \pA (x_0,\ldots , x_n)\mbox{ or }E \rightarrt \nocom {\bf alt} (\phi \psi , B)=\pA (x_{j-1},x_j)
$$
is the given structural map. 

Now with these notations, $\Gen (\pA ;x_0,\ldots , x_n)(z_0,\ldots ,z_m)$ can also be expressed as a pushout of the form
$$
\begin{diagram}
\coprod _{\phi\in  \Delta _X(z_{\cdot}, x_{\cdot})} {\bf alt}(\phi,\pA (x_0,\ldots , x_n)) & \rightarr &
\coprod _{\phi\in  \Delta _X(z_{\cdot}, x_{\cdot})} {\bf alt}(\phi,E) \\
\downarr & & \downarr \\
\pA (z_0,\ldots , z_m) & \rightarr & \Gen (\pA ;x_0,\ldots , x_n)(z_0,\ldots ,z_m).
\end{diagram}
$$
Note that for a map $\phi : z_{\cdot}\rightarrt x_{\cdot}$ factoring through $(x_{j-1},x_j)$ (resp. $(x_j)$) we get a map 
$(z_0,\ldots , z_m)\rightarrt (x_{j-1},x_j)$ (resp. to $(x_j)$) and hence a map 
$\pA (x_{j-1},x_j) \rightarrt \pA (z_0,\ldots , z_m)$ (resp. a map $\ast = \pA (x_j)\rightarrt \pA (z_0,\ldots , z_m)$).
These combine to give the left vertical map. 

Now, in case of a map $\psi : (y_0,\ldots , y_p)\rightarrt (z_0,\ldots , z_m)$ we get maps 
$$
\coprod _{\phi\in  \Delta _X(z_{\cdot}, x_{\cdot})} {\bf alt}(\phi,\pA (x_0,\ldots , x_n)) \rightarrt
\coprod _{\phi\in  \Delta _X(z_{\cdot}, x_{\cdot})} {\bf alt}(\phi\psi ,\pA (x_0,\ldots , x_n))
$$
$$ 
\rightarrt 
\coprod _{\zeta \in  \Delta _X(y_{\cdot}, x_{\cdot})} {\bf alt}(\zeta ,\pA (x_0,\ldots , x_n))
$$
and 
$$
\coprod _{\phi\in  \Delta _X(z_{\cdot}, x_{\cdot})} {\bf alt}(\phi,E) \rightarrt
\coprod _{\phi\in  \Delta _X(z_{\cdot}, x_{\cdot})} {\bf alt}(\phi\psi ,E) 
\rightarrt 
\coprod _{\zeta \in  \Delta _X(y_{\cdot}, x_{\cdot})} {\bf alt}(\zeta ,E)
$$
which induce the map of functoriality
$$
\Gen (\pA ;x_0,\ldots , x_n)(z_0,\ldots ,z_m)\rightarrt 
\Gen (\pA ;x_0,\ldots , x_n)(y_0,\ldots , y_p).
$$

This whole discussion can be simplified considerably in the case when the $x_0,\ldots , x_n$ are all distinct
in other words $\forall 0\leq i\neq j\leq n,\; x_i\neq x_j$. Then, for any sequence of objects $(z_0,\ldots , z_m)\in \Delta _X$
there is at most one  map $\phi : (z_0,\ldots , z_m)\rightarrt (x_0,\ldots , x_n)$. 
It factors through an adjacent pair if
and only if $z_{\cdot}$ is of the form $(x_{j-1},\ldots , x_{j-1},x_j,\ldots , x_j)$ and it factors through a singleton if and only
if $z_{\cdot}= (x_j,\ldots , x_j)$. There is a map, but not factoring through an adjacent pair or a singleton, if and only if
$$
z_{\cdot}= (x_{i_0},\ldots , x_{i_m})
$$
where $0\leq i_0\leq i_1\leq \cdots \leq i_m\leq n$ and $i_0+1 < i_m$. We state the result in this case as a lemma.

\begin{lemma}
\label{gendescription-lem}
Suppose $\pA \in \precat (X,\mM )$ and
$(x_0,\ldots , x_n)$ is a sequence of pairwise disjoint objects. Suppose given a choice of factorization $E,e, p_1,\ldots , p_n$ as above. 
Then for any increasing sequence of indices $0\leq i_0\leq i_1\leq \cdots \leq i_m\leq n$ with $i_0+1 < i_m$,
$$
\Gen (\pA ; x_0,\ldots , x_n) (x_{i_0},\ldots , x_{i_m}) = \pA (x_{i_0},\ldots , x_{i_m})\cup ^{\pA (x_0,\ldots , x_n)}E.
$$
For any other sequence of objects $(z_0,\ldots , z_m)$ we have 
$$
\Gen (\pA ; x_0,\ldots , x_n)(z_0,\ldots , z_m) = \pA (z_0,\ldots , z_m).
$$
These expressions are compatible with the natural maps from $\pA (z_0,\ldots , z_m)$, and the maps of functoriality are given by the
structural maps for $E$ in the case of a map $\psi : (z_0,\ldots , z_m)\rightarrt (x_{i_0},\ldots , x_{i_m})$ which factors through a
sequence projecting to an 
adjacent pair in $(x_0,\ldots , x_n)$. The maps of functoriality are given by the degeneracy maps of $\pA $ in case $\psi$ factors through a repeated singleton. 
\end{lemma}
\begin{proof}
Left to the reader. 
\end{proof}

\section{Fixing the fibrant condition locally}

Applying the operation $\Gen$ may destroy the fibrancy condition in the plain diagram
Reedy structure. Let $\fix :\pA \rightarrt \Fix (\pA )$ be a replacement by a fibrant
object in the Reedy model structure on the plain diagram category $\diag _{\rm Reedy}(\Delta ^o_X/X,\mM )$
where $X=\Ob (\pA )$. The map $\fix (\pA )$
is a pushout along elements of $K_{\rm Reedy}$ namely the generators for the
trivial cofibrations of $\diag _{\rm Reedy}(\Delta ^o_X/X,\mM )$. In particular
it is an isotrivial cofibration and levelwise a trivial cofibration. This doesn't
change the homotopy type in $\Ho (\diag (\Delta ^o_X/X,\mM ))$.

\section{Combining generation steps}

We can put together the elementary generation steps defined above, to obtain a model for the replacement to a Segal object. 

\begin{theorem}
\label{transgen}
Suppose $\pA \in \precat (X, \mM )$. There is a transfinite composition sequence 
$\{ \pA _i\} _{i\in [\alpha ]}$ indexed by an ordinal $[\alpha ] = \alpha +1$, with $\pA _0=\pA $, such that $\pA _{\cdot}$ is continuous at limit ordinals,
such that $\pA _{\alpha}$ satisfies the Segal conditions, and such that each $\pA _i\rightarrt \pA _{i+1}$ has the form 
$\fix \gen (\pA _i,q_i)$ for some $q_i= (x_{i,0},\ldots , x_{i,n_i})$. In other words, 
$\pA _{i+1}= 
\Fix \Gen (\pA _i,q_i)[E_i,e_i, p_{i,1},\ldots , p_{i,n_i}]$ as considered above. 
Furthermore the transition maps are all trivial cofibrations in $\precat _{\rm Reedy}(X; \mM )$, in particular $\pA \rightarrt \pA _{\alpha}$ is a trivial cofibration.
Denote the end result by $\Seg ^{FG}(\pA )$; it is a fibrant object in
the Reedy model structure $\precat (X,\mM )$ and equivalent to $\Seg (\pA )$ 
in the Reedy model structure on the plain 
diagram category $\diag _{\rm Reedy}(\Delta ^o_X/X,\mM )$,
indeed it is a possible choice for the fibrant replacement $\Seg (\pA )$.
\end{theorem}
\begin{proof}
The $\gen (\pA _i, q_i)$ are trivial cofibrations which
may be chosen to contain a pushout along any particular cofibration
of the form $\Psi ([n],f)$. Similarly the maps $\fix$ may be chosen to
contain any of the pushouts along generating trivial cofibrations
for the level structure. Thus, the maps envisioned above 
may contain all pushouts along elements of the pseudo-generating
set $K_{\rm Reedy}$. Choose the sequence using the small object argument
so that the end result $\pA _{\alpha}$ is in $\inj (K_{\rm Reedy})$.
In particular it satisfies the Segal conditions.  
Note that $K_{\rm Reedy}$-injective objects are also fibrant, see Theorem \ref{reedyprecat}
where that statement really comes from Theorem \ref{directLBL}. Thus,
$\pA \rightarrt \Seg ^{FG}(\pA )$ is a trivial cofibration towards a fibrant object
in $\precat _{\rm Reedy}(X,\mM )$. 
\end{proof}

An important corollary of this procedure is the following statement, tautological for maps between Segal categories.  

\begin{corollary}
\label{loctoglobwe}
Suppose a map $f:\pA \rightarrt \pB $ in $\precat (\mM )$ induces a weak equivalence
$$
\pA (x_0,\ldots , x_p)\rightarrt^{\sim} \pB (f(x_0),\ldots , f(x_p))
$$
for any sequence $(x_0,\ldots , x_p)\in \Delta _{\Ob (\pA )}$. Then $f$ is fully faithful.
\end{corollary}
\begin{proof}
In the process of Theorem \ref{transgen} this condition is preserved at each step. 
It follows that $\Seg ^{FG}(\pA )\rightarrt \Seg ^{FG}(\pB )$ is fully faithful,
which is by definition the full faithfulness criterion for $f$. 
\end{proof}

\section{Functoriality of the generation process}
\label{funcgen}

Suppose $\mM$ and $\mN$ are tractable left proper cartesian model categories. Suppose $F:\mM \rightarrt \mN$ is a left Quillen functor.
We say that $F$ is {\em compatible with products} if, for any $A,B \in \mM$ the natural map $F(A \times B)\rightarrt F(A )\times F(B)$ is
a weak equivalence in $\mN$. 

The functor $F$ induces functors $\precat (X,F): \precat (X,\mM )\rightarrt \precat (X,\mN )$ and hence, as $X$ varies,
it induces a functor $\precat (F):\precat (\mM )\rightarrt \precat (\mN )$. For $\pA \in \precat (X,\mM )$ and a sequence of objects $x_0,\ldots , x_n$ in 
$X$,
$$
F(\pA )(x_0,\ldots , x_n):= F(\pA (x_0,\ldots , x_n)). 
$$

\begin{lemma}
In the above situation, suppose $\pA \in \precat (X,\mM )$, choose
$q= (x_0,\ldots , x_n)$, and suppose given a choice of factorization $E,e, p_1,\ldots , p_n$ as above. 
Suppose $e$ is a cofibration and $(p_1,\ldots , p_n): E\rightarrt \pA (x_0,x_1)\times \cdots \times \pA (x_{n-1},x_n)$
is a weak equivalence in $\mM$. Then $F(e)$ is a cofibration and 
$$
(Fp_1,\ldots , Fp_n): F(E)\rightarrt F(\pA )(x_0,x_1)\times \cdots \times F(\pA )(x_{n-1},x_n)
$$
is a weak equivalence in $\mN$, so $(F(E), F(e), F(p_1),\ldots , F(p_n))$ constitute data for defining 
$$
F(\pA )\rightarrt^{\gen (F(\pA ),q)} \Gen (F(\pA ),q)=\Gen (F(\pA ),q)[F(E), F(e), F(p_1),\ldots , F(p_n)]
$$
in $\precat (X,\mN )$. In these terms we have $\gen (F(\pA ),q) = F(\gen (\pA ,q))$ and
$$
\Gen (F(\pA ),q)[F(E), F(e), F(p_1),\ldots , F(p_n)] = F(\Gen (\pA ,q)[E,e, p_1,\ldots , p_n] )
$$
written more succinctly as $\Gen (F(\pA ),q)= F(\Gen (\pA ,q))$.
\end{lemma}
\begin{proof}
By inspection. 
\end{proof}

\begin{corollary}
In our above situation of a left Quillen functor $F:\mM \rightarrt \mN$ compatible with products, between two tractable left proper cartesian model categories, 
suppose $\pA \in\precat (X,\mM )$ and let $\{ \pA _i\} _{i\in [\alpha ]}$ be a transfinite composition of elementary generation steps in 
$\precat (X,\mM )$ with $\pA =\pA _0\rightarrt \pA _{\alpha}$ a trivial cofibration in $\precat _{\rm inj}(X,\mM )$ and 
$\pA _{\alpha}$ satisfying the Segal conditions. Then $\{ F(\pA _i)\} _{i\in [\alpha ]}$ is a transfinite composition of elementary generation
steps in $\precat (X,\mN )$, the map $F(\pA ) \rightarrt F(\pA _{\alpha})$ is a trivial cofibration in $\precat _{\rm inj}(X,\mN )$,
and $F(\pA _{\alpha})$ satisfies the Segal conditions as an $\mN$-enriched precategory. 
\end{corollary}
\begin{proof}
Just apply $F$ to the sequence of Theorem \ref{transgen}. 
\end{proof}

We state the next corollary as a theorem, since it is the main statement we need from this section.
Recall from Proposition \ref{precatF} that given a left Quillen functor
between tractable left proper cartesian model categories, it induces
a left Quillen functor on the Reedy model categories of precategories. 
One didn't even need to suppose that $F$ preserves products. If that is the case,
the statement can be slightly improved to say that $\precat (X,F)$ preserves
weak equivalences (rather than just trivial cofibrations). 

\begin{proposition}
\label{precatFpreserveswe}
Suppose
$F:\mM \rightarrt \mN$ is a left Quillen functor, preserving weak equivalences and
compatible with products,
between two tractable left proper cartesian model categories. Then $\precat (X; F) :\precat _{\rm Reedy}(X,\mM )\rightarrt \precat _{\rm Reedy}(X,\mN )$ 
is a left Quillen functor which takes weak equivalences to weak equivalences. 
\end{proposition}
\begin{proof}
Suppose $\pA \rightarrt \pB$ is a weak equivalence in 
$\precat _{\rm Reedy}(X,\mM )$. This means that
$\Seg ^{FG}(\pA )\rightarrt \Seg ^{FG}(\pB )$ is a levelwise weak equivalence
of diagrams, so applying $\precat (X,F)$ gives a new levelwise weak equivalence of
diagrams towards $\mN$. By the previous corollary, when we apply $\precat (X,F)$
we get $\Seg ^{FG}(\precat (X,F)\pA )\rightarrt \Seg ^{FG}(\precat (X,F)\pB )$,
it follows that $\precat (X,F)\pA )\rightarrt \precat (X,F)\pB $
is a weak equivalence in $\precat _{\rm Reedy}(X,\mN )$. 
\end{proof}

\section{Example: generators and relations for $1$-categories}
\label{genrel1cat}

It is interesting and instructive to consider the case where $\mM = \Sets$ is the model category of sets. Here,
the weak equivalences are isomorphisms, and the fibrations and cofibrations are arbitrary maps. It is easily seen to be tractable, left proper and cartesian. 
The category of $\Sets$-enriched precategories $\precat (\Sets )$
may be identified with the category of simplicial sets. For $\pA \in \precat (\Sets  )$, if we make this identification
then the $0$-simplices correspond to the objects; the $1$-simplices give generators for the morphisms, and the $2$-simplices 
give relations among the morphisms. The Segal conditions for
$\pA \in \precat (\Sets )$ are exactly the classical conditions which are equivalent to
stating that $\pA$ is the nerve of a category. The calculus of generators and relations
constructed above, reduces in this case to the classical calculus of generators and relations for $1$-categories.  

Fix a set $X$. For this section, the main part of the data of $\pA \in \precat (X,\Sets )$
consists of sets $\pA (x,y)$ for each pair $x,y\in X$, and $\pA (x,y,z)$ for 
each triple $x,y,z\in Z$. Recall that $\pA (x)$ is a singleton, so the degeneracy
maps yield elements denoted $1_x\in \pA (x,x)$. Similarly, for any $f\in \pA (x,y)$
the degeneracies yield elements denoted $[1_yf]\in \pA (x,y,y)$ and $[f1_x]\in \pA (x,x,y)$.
In the case $f=1_x$ there is no confusion as both elements denoted $[1_x1_x]$ correspond
to the same element of $\pA (x,x,x)$ obtained from the singleton $\pA (x)$ by 
degeneracy. We have projections 
$$
i_{01}^{\ast}:\pA (x,y,z)\rightarrt \pA (x,y), \;\;\; 
i_{12}^{\ast}:\pA (x,y,z)\rightarrt \pA (y,z), 
$$
$$
i_{02}^{\ast}:\pA (x,y,z)\rightarrt \pA (x,z).
$$
These are compatible with the degeneracies in an obvious way, for example 
$i_{02}^{\ast}[1_yf] = i_{02}^{\ast}[f1_x] = f$.

The elements of $\pA (x,y)$ are the ``generating arrows'' between $x$ and $y$.
An element $r\in \pA (x,y,z)$ corresponds to an ``elementary relation'' of the form
$f=gh$, where 
$$
f=i_{02}^{\ast}(r)\in \pA (x,z), \;\;\; g=i_{12}^{\ast}(r)\in \pA (y,z),
\;\;\; h=i_{01}^{\ast}(r)\in \pA (x,y).
$$
This can be seen by looking more closely 
at what happens under the generation step $\Gen (\pA , q)$ for a triple 
$q=(x,y,z)$. Such a step involves a choice of factorization
$$
\pA (x,y,z)\rightarrt ^e E \rightarrt ^{(p_1,p_2)} \pA (x,y)\times \pA (y,z).
$$
All maps are cofibrations in $\Sets$ so there is no restriction on $e$;
on the other hand the weak equivalences are the isomorphisms, so $(p_1,p_2)$
has to be an isomorphism and we may in effect suppose $E=\pA (x,y)\times \pA (y,z)$.

We can describe explicitly the resulting precategory $\pA ':= \Gen (\pA ,q)$.
Suppose $(u,v)$ is a pair of elements of $X$. We  need to consider the maps 
$\phi : (u,v)\rightarrt (x,y,z)$ in $\Delta _X$. There are $6$ possibilities:
$$
(u,v)=(x,y)\rightarrt ^{i_{01}}(x,y,z),
$$
$$
(u,v)=(y,z)\rightarrt ^{i_{12}}(x,y,z),
$$
$$
(u,v)=(x,z)\rightarrt ^{i_{02}}(x,y,z),
$$
$$
(u,v)=(x,x)\rightarrt ^{i_{00}}(x,y,z),
$$
$$
(u,v)=(y,y)\rightarrt ^{i_{11}}(x,y,z),
$$
$$
(u,v)=(z,z)\rightarrt ^{i_{22}}(x,y,z).
$$
In case of coincidences among the $x,y,z$ these can overlap in the sense that
the same $(u,v)$ could have several maps to $(x,y,z)$. However, in the notation
of Section \ref{sec-gen}, the only one of these maps which 
is in $\Delta ^{NS}_X((u,v),(x,y,z))$ is
$i_{02}$. Thus, using $E=\pA (x,y)\times \pA (y,z)$ we have a pushout diagram
$$
\begin{diagram}
\pA (x,y,z) & \rightarr & \pA (x,y)\times \pA (y,z) \\
\downarr & & \downarr \\
\pA (x,z) & \rightarr & \pA ' (x,z).
\end{diagram}
$$
In other words, $\pA '(x,z)$ is obtained by taking the old $\pA (x,z)$
and adding the symbols $gh$ for $g\in \pA (y,z)$ and $h\in \pA (x,y)$,
subject to the relations that $f=gh$ whenever there is an 
element of $\pA (x,y,z)$ mapping to 
$f\in \pA (x,z)$, $h\in \pA (x,y)$ and $g\in \pA (y,z)$.

Similarly, if $(u,v,w)$ is a triple of elements, the only maps 
$$
(u,v,w)\rightarrt (x,y,z)
$$
which are in $\Delta ^{NS}_X((u,v),(x,y,z))$ are the identity $(x,y,z)\rightarrt (x,y,z)$,
and the degeneracies of the previous map $i_{02}$ 
$$
(x,z,z)\rightarrt (x,y,z),\;\;\; (x,x,z)  \rightarrt (x,y,z).
$$
The pushout expression along the identity says that
for any such symbol $gh$ which is added to $\pA '(x,z)$ following the previous
discussion, there will also be a corresponding
element of $\pA '(x,y,z)$ saying that it is the composition of $g$ and $h$,
in other words the formal composition $gh$ will be recorded as a composition of $g$ and $h$. 
The degeneracies of $i_{02}$ take care of keeping track of left and right identities
for the new morphisms $gh$ which are added. 

We leave as an exercise to verify that the generation steps $\Gen (\pA , q)$ at
quadruples $q=(w,x,y,z)$ correspond to enforcing the associativity axiom.

Putting this all together, we see that after repeating the generation steps infinitely
many times for each uplet of objects, the resulting precategory $\Seg ^{FG}(\pA )$
is the nerve of a category, and it is the category generated by the original 
arrows $\pA (x,y)$ subject to the original relations $\pA (x,y,z)$.


\chapter{Generators and relations for Segal categories}
\label{secat1}

In this chapter, we consider one of the main examples of the theory: when $\mM =\mK$ is the Kan-Quillen
model category of simplicial sets. 

\begin{theorem}
\label{kanquillen}
The Kan-Quillen model category $\mK$ of simplicial sets, is a 
tractable left proper cartesian model category.
\end{theorem}
\begin{proof}
Cofibrations are just monomorphisms of simplicial sets, so the cartesian and tractability conditions are easy to check. 
Condition (PROD) of Definition \ref{def-cartesian} for trivial cofibrations is Eilenberg-Zilber. 
\end{proof}

\begin{corollary}
\label{precatKstrucs}
For a fixed set of objects $X$, we get model categories  of $\mK$-enriched precategories over $X$ with either
the projective structure $\precat _{\rm proj}(X; \mK )$, or
the Reedy structure $\precat _{\rm Reedy}(X; \mK )$ which is the same as the
injective structure $\precat _{\rm inj}(X; \mK )$.
\end{corollary}
\begin{proof}
Apply Theorems \ref{modstrucs} and \ref{reedyprecat}. Note that $\mK$ is
a presheaf category and the cofibrations are monomorphisms, so by Proposition 
\ref{reedyinjective} the Reedy and injective cofibrations coincide. The weak
equivalences are the same in all structures so the Reedy and injective structures are the
same.   
\end{proof}

A $\mK$-precategory is called a {\em Segal precategory}; those satisfying the Segal conditions are called {\em Segal categories}. 
We have defined the classes of global weak equivalences, and two
flavors of cofibrations (hence fibrations), on the way to constructing the
projective and Reedy$=$injective model 
structures 
on the category $\precat (\mK )$ of Segal precategories 
in the next part. 
The notations of Section \ref{sec-interpretations} apply, so a Segal precategory may be considered as a bisimplicial
set ($m\mapsto \pA _{m/}\in \mK$ or $m,n\mapsto \pA _{m,n}\in \Sets$, see below)
whose simplicial set $\pA _{0/}$ in degree $0$ is the discrete set $\Ob (\pA )$. 

The theory of Segal categories goes back a long time. 
They were first considered by Dwyer, Kan, Smith \cite{DKS} and
Schw\"{a}nzl and Vogt \cite{Vogt}, but of course the definition goes essentially back to Segal \cite{Segal} who just restricted to the case when there is only one object.
Other references are  \cite{Adams} \cite{Thomason} \cite{CordierPorter}.
Bergner \cite{BergnerThreeModels}
and Pelissier \cite{Pelissier} have already given complete
constructions of the global model category structure. 
The model categories $\precat _{\rm proj}(\mK )$ and $\precat _{\rm Reedy}(\mK )$ which
we are in the process of constructing, are the same as those of Bergner \cite{BergnerThreeModels}.

\section{Segal categories}
\label{secatreal}

Recall that if $\pA $ is a Segal category, then the {\em truncation} $\tau _{\leq 1}(\pA )$
was defined in Section \ref{sec-globalwe}. It is the usual $1$-category
whose nerve is the simplicial set $m\mapsto \pi _0(\pA _{m/})$. 
A {\em Segal groupoid}
is a Segal category such that $\tau _{\leq 1}(\pA )$ is a groupoid. 

The category $\mK$ satisfies the additional Condition \ref{disj} (DISJ) on disjoint unions. Therefore $\precat (\mK )$ may be viewed as a presheaf
category. More precisely, as $\mK$ is the category of presheaves on $\Delta$,
by Theorem \ref{interp1} we have a natural identification 
\begin{equation}
\label{seprecatident}
\precat (\mK ) \cong \Sets ^{\scone (\Delta )} \subset \Sets ^{\Delta \times \Delta}.
\end{equation}
A Segal precategory for us is a pair $(X,\pA )$
where $X\in \Sets$ and $\pA $ is a collection of simplicial sets $\pA (x_0,\ldots , x_n)$
for sequences of $x_i\in X$. 
Recall that $\scone (\Delta )$ is a quotient of $\Delta \times \Delta$,
thus the subset relation in \eqref{seprecatident}.
Under the identification \eqref{seprecatident}
this corresponds to a bisimplicial set given by the formula 
$$
\widetilde{\pA }:(n,m)\in \scone (\Delta ) \mapsto \coprod _{(x_0,\ldots , x_n)\in X^{n+1}}
\pA (x_0,\ldots , x_n)_m.
$$
If $n=0$ then this is constant in $m\in \Delta $, which is {\em Tamsamani's constancy condition} in this case. The constancy condition characterizes the bisimplicial
sets which are presheaves on the quotient $\scone (\Delta )$. 

In what follows, we use the following notation 
of Section \ref{sec-interpretations} for a Segal precategory $\pA$:
\newline
---$\pA _0  = \Ob (\pA )=\widetilde{\pA}(0,m)$ for any $m$;
\newline
---$\pA _{n/}$ is the simplicial set $m\mapsto \widetilde{\pA }(n,m)$. 

The assignment $n\mapsto \pA _{n/}$ is a functor $\Delta ^o\rightarrt \mK$,
which is the same as the bisimplicial set $\widetilde{\pA }$.
We denote the first component by $\pA _0$ rather than $\pA _{0/}$ to
emphasize the constancy condition, that it is a set considered as a
constant simplicial set. 

This notation is convenient
because
colimits (resp. limits) in $\precat (\mK )$ correspond to levelwise colimits (resp. limits)
of the
functors $\Delta ^o\rightarrt \mK$. 

For a more pictorial point of view, it is often more intuitive to
replace the model category $\mK$ by the category $Top$ of topological spaces, 
and to think of Segal (pre)-categories as functors $\Delta _X^o\rightarrt \Top$.
As usual, for technical details it is more convenient to use simplicial sets, and
we don't treat the thorny questions surrounding model category structures on $Top$. 
Below we will sometimes replace $\mK$ by $\Top$ and leave to the reader to insert the appropriate realization and 
singular complex functors between these two.

A Segal precategory may therefore
be thought of as a functor $\pA :\Delta ^o \rightarrt \Top$ denoted $n\mapsto \pA _{n/}$,
which is to say a simplicial space, satisfying the ``contancy'' or ``globular'' condition that
$\pA _0=\Ob (\pA )$ is a discrete set. 

\section{The Poincar\'e-Segal groupoid}
\label{PoincareSegal}

Given a Kan simplicial set $X$, that is a fibrant object in $\mK$, we can define its  
{\em Poincar\'e-Segal groupoid} $\Pi _S(X)$ which is a Segal groupoid, i.e. a Segal category
whose $\tau _{\leq 1}$ is a groupoid. 
It is constructed as 
the right adjoint of the diagonal realization functor, which we consider first.

The {\em diagonal} ${\bf d}:\Delta \rightarrt \Delta \times \Delta$
provides a pullback functor 
$$
{\bf d}^{\ast} : \Sets ^{\Delta \times \Delta}\rightarrt \Sets ^{\Delta} =\mK ;
$$
so composing with \eqref{seprecatident} we obtain the {\em realization functor}
$$
| \;\; | : \precat (\mK )\rightarrt \mK , 
$$
$$ 
|(X,\pA )|:= {\bf d}^{\ast}\widetilde{\pA } = \left( m\mapsto \widetilde{\pA }_{m,m}\right) .
$$
Recall that $\sconemap _{\Delta}: \Delta \times \Delta \rightarrt \scone (\Delta )$
denotes the projection map. We obtain a functor
$$
\sconemap _{\Delta}\circ {\bf d}: \Delta \rightarrt \scone (\Delta ),
$$
and the realization functor is just the pullback 
\begin{equation}
\label{realpullbackform}
| \pA  | = (\sconemap _{\Delta}\circ {\bf d})^{\ast}(\pA ).
\end{equation}
from $\precat (\mK )\cong \Sets ^{\scone (\Delta )}$ to $\mK = \Sets ^{\Delta}$. 
It follows that $| \;\; |$ preserves limits and colimits, indeed it has both left and
right adjoints.  

Taking the right adjoint of the expression \eqref{realpullbackform}, define
$$
\Pi _S := (\sconemap _{\Delta}\circ {\bf d})_{\ast} : \mK \rightarrt \precat (\mK ).
$$
Note that $\Pi _S$ commutes with limits.
For a fibrant simplicial set $X$, call $\Pi _S(X)$ the {\em Poincar\'e-Segal groupoid}
of $X$. 

We now state the general theorem relating the homotopy theory in $\precat (\mK )$
of Segal groupoids, with the classical homotopy theory of simplicial sets. 
It is essentially due to Segal, although 
Segal's arguments were mostly stated for the
situation of categories with a single  object. These were reviewed in the more general
categorical setting by Tamsamani \cite{Tamsamani}. Tamsamani then furthermore iterated
the result to obtain an equivalence between the theory of $n$-truncated homotopy types,
and their Poincar\'e $n$-groupoids. We don't delve into the details of this, refering the
reader to \cite{Tamsamani} instead. 

The theorem essentially says that we have a Quillen adjunction, although at the
present point in our argument we haven't yet shown that the given classes of morphisms
provide a model structure for $\precat (\mK )$. Of course that statement is also
contained in the references \cite{BergnerThreeModels} \cite{Pelissier}. 

\begin{theorem}
\label{poincaresegal}
The realization functor $| \;\; |$ sends injective
cofibrations (resp. injective global trivial cofibrations) to 
cofibrations (resp. trivial cofibrations) of simplicial sets. 
The Poincar\'e-Segal groupoid
functor sends fibrations in $\mK$ to new fibrations in $\precat (\mK )$. 
If $X$ is a fibrant simplicial set then $\Pi _S(X)$ is fibrant, in particular it
is a Segal category. It is, in fact, a Segal groupoid, and furthermore the
adjunction  map 
$$
| \Pi _S(X)|\rightarrt X
$$
is a weak equivalence. Conversely, if $\pA \in \precat (\mK )$ is a Segal groupoid 
and $|\pA |\rightarrt Y$ is a fibrant replacement in $\mK$, 
then the map obtained by adjunction
$$
\pA  \rightarrt \Pi _S(Y)
$$ 
is a global equivalence of Segal categories. 
\end{theorem}
\begin{proof}
We refer to Tamsamani \cite{Tamsamani} for the interpretation of Segal's original
results \cite{Segal} in the context of many-object groupoids. 
See also  
Bergner \cite{BergnerThreeModels}, Dwyer, Kan and Smith \cite{DKS},
Berger \cite{BergerCellularNerve}, Duskin \cite{Duskin}, and others. 
\end{proof}

Tamsamani defines the homotopy groups of an $n$-nerve in \cite{Tamsamani}.
Taken at the first level, we can similarly define the homotopy groups
of a Segal groupoid $\pA$ at any $x\in \Ob (\pA )$ by putting
$$
\pi _0(\pA ):= {\rm Iso}\tau _{\leq 1}\pA 
$$
and for $i\geq 1$, 
$$
\pi _i(\pA , x):= \pi _{i-1}(\pA (x,x), 1_x)
$$
where $1_x:\ast \rightarrt \pA (x,x)$ is the degeneracy map using $\pA (x)=\ast$.

From the Segal groupoid condition it follows that we can
define the relation of homotopy on $\Ob (\pA )$ by saying that $x\sim y$ if
and only if $\pA (x,y)\neq \emptyset$, and then  $\pi _0(\pA ) = \Ob (\pA )/\sim$. 

\begin{lemma}
A morphism $f:\pA \rightarrow \pB$ between Segal groupoids is a global weak
equivalence, if and only if $\pi _0(\pA )\rightarrt \pi _0(\pB )$ is surjective
(resp. an isomorphism),
and
for each $i\geq 1$ and each
object $x\in \Ob (\pA )$, the induced  maps $\pi _i(\pA ,x)\rightarrt \pi _i(\pB , f(x))$
are isomorphisms. 
\end{lemma}
\begin{proof}
This follows directly from the definition using the criterion that a map between simplicial
sets is a weak equivalence if and only if it induces an isomorphism on homotopy groups.
\end{proof}

\begin{proposition}
\label{segalpi}
If $\pA$ is  a Segal groupoid then $\pi _0(\pA )=\pi _0(|\pA |)$ and 
for any $i\geq 1$ and $x\in \Ob (\pA )$, $\pi _i(\pA , x)=\pi _i(|\pA |,|x|)$.
Similarly if $Z$ is a simplicial set satisfying Kan's fibrancy condition then
$\pi _0(\Pi _S(Z))=\pi _0(Z)$ and for any vertex $z\in Z_0$, $\pi _i(\Pi _S(Z),z)= \pi _i(Z,z)$. 

A morphism $f:\pA \rightarrt \pB$ between Segal groupoids is a global weak equivalence
if and only if the induced map on realizations is a weak equivalence $|\pA |\sim |\pB |$.
\end{proposition}
\begin{proof}
The first part states some of the essential facts about Segal's construction 
\cite{Segal}, entering into the proof of Theorem \ref{poincaresegal},
see \cite{Tamsamani}. The second paragraph follows immediately from 
the first part together with the previous lemma.
\end{proof}

\section{The calculus}

After this first part, the remainder of the chapter is devoted to following out the
calculus of generators and relations in the case of Segal precategories. 
In the spirit of Segal's ``delooping machine'', this process gives
a good explanation of how the technical machinery introduced in the previous chapters
works. First is a detailed description of how to go from a Segal 
precategory to a Segal category, keeping track of the connectivity properties
of the intervening spaces. 
In Section \ref{computingLoop}, we show how the
process leads, in principle, to a calculation of the loop space of a space.  
In Section \ref{s2example} below we will be able to follow along what happens as
one of the first nontrivial homotopy groups $\pi _3(S^2)$ appears out of this process. 
This chapter
constitutes a reworked version of the preprint \cite{effective}. 

In view of the topological motivation, we concentrate
on the case of Segal groupoids and more specifically the case when there is
only one object and $\pA _{1/}$ is connected. In this case, a stronger finiteness
property will hold: to arrange things up to a certain level of connectivity,
it will suffice to do a finite number of elementary generation operations
which we denote here by $Arr$ (or later $Arr2$). Of course, if $\pA _{1/}$
is not connected then we may be in the presence of a fundamental group and 
it requires, in principle, an infinite and even undecideable number of operations
to compute the group. Some remarks on this aspect are given near the end of the chapter.

\subsection{Arranging in degree $m$}

Suppose $\pA $ is a Segal precategory with $\pA _0 = \ast$. We say that $\pA $ is {\em
$(m,k)$-arranged} if the Segal map
$$
\pA _{m/} \rightarrt \pA _{1/} \times  \ldots \times \pA _{1/}
$$
induces isomorphisms on $\pi _i$ for $i< k$ and a surjection on $\pi _k$.

Note that for $l\geq k$, adding $l$-cells to $\pA _{m/}$ or $l+1$-cells to
$\pA _{1/}$ doesn't affect this property.

We now define an operation where we try to ``arrange'' $\pA $ in degree $m$.
This operation is inspired by the operation $\Gen$ considered in Section \ref{sec-gen}. 
In the present version it will be called
$$
\pA  \mapsto Arr (\pA , m).
$$
Fix $m$ in what follows.
Let $C$ be the mapping cone of the Segal map
$$
\pA _{m/} \rightarrt \pA _{1/}\times \ldots \times \pA _{1/}.
$$
To be precise, as a bisimplicial set
$$
C = (I \times \pA _{m/}) \cup ^{\{ 1\} \times \pA _{m/}} (\pA _{1/}\times \ldots \times
\pA _{1/}),
$$
where $I$ is the standard simplicial interval, and the notation is coproduct
of bisimplicial sets (note also that the globular condition is preserved, so it
is a coproduct of Segal precategories). Note that $\{ 1 \} \times \pA _{m/}$ denotes the
second endpoint of the interval crossed with $\pA _{m/}$. We have morphisms
$$
\pA _{m/} \rHookarrow ^a C
\rightarrt ^b \pA _{1/}\times \ldots \times
\pA _{1/}, $$
the morphism $a$ being the inclusion of $\{ 0\} \times \pA _{m/}$ into
$I \times \pA _{m/}$ (thus it is a cofibration i.e. injection of simplicial
sets) and the second morphism $b$ coming from the projection  $I\times \pA _{m/}
\rightarrt \pA _{m/}$. The second morphism $b$ is a weak equivalence.

We now define $Arr (\pA , m)$: for any $p$, let
$$
Arr (\pA ,m)_{p/} := \pA _{p/} \cup ^{(\bigcup \pA _{m/})} \left( \bigcup _{p\rightarrt
m}
 C \right)
$$
be the combined coproduct of $\pA _{p/}$ with several copies of the
morphism $a:\pA _{m/} \rightarrt C$ , one copy for each map $p\rightarrt m$ not
factoring through a principal edge (see below for further discussion of this
condition), these maps inducing $\pA _{m/} \rightarrt \pA _{p/}$.

We need to define $Arr (\pA ,m)$ as a Segal precategory, i.e. as a bisimplicial set.
For this we need morphisms of functoriality
$$
Arr (\pA , m)_{p/} \rightarrt Arr (\pA , m)_{q/}
$$
for any $q\rightarrt p$. These are defined as follows.
We consider a component of $Arr (\pA , m)_{p/}$ which is a copy of $C$
attached along a map $\pA _{m/}\rightarrt \pA _{p/}$ corresponding to $p\rightarrt
m$ which doesn't factor through a principal edge.  If the composed map
$q\rightarrt p\rightarrt m$ doesn't factor through a principal edge
then the component $C$ maps to the corresponding component of
$Arr (\pA , m)_{1/}$. If the map does factor through a principal edge
$q\rightarrt 1\rightarrt m$ then we obtain a map
$C \rightarrt \pA _{1/}$ (the component of the map $b$ corresponding to this
principal edge). Compose with the map $\pA _{1/}\rightarrt \pA _{q/}$ to obtain a
map $C\rightarrt \pA _{q/}$. Note that if the map further factors
$$
q\rightarrt 0 \rightarrt 1 \rightarrt m
$$
then the map $\pA _{1/}\rightarrt \pA _{q/}$ factors through the basepoint
$$
\pA _{1/}\rightarrt \pA _0 \rightarrt \pA _{q/},
$$
and our map on $C$ factors through the basepoint. This factorization doesn't
depend on choice of principal edge containing the map $0\rightarrt m$.

One can verify that this prescription defines a functor $p\mapsto Arr
(\pA ,m)_{p/}$ from $\Delta$ to simplicial sets. This verification will be a
consequence of the more conceptual description which follows.

Let $h(m)$ denote the simplicial set representing the standard $m$-simplex;
it is the contravariant functor on $\Delta$ represented by the object $m$.
Let $\Sigma (m)\subset h(m)$ be the subcomplex which is the union of the
principal edges.

\begin{definition}
\label{boxtimesdef}
If $X$ is a simplicial set and $B$ is another simplicial set
denote by $X\boxtimes B$ the bisimplicial set exterior product, defined by
$$
(X\boxtimes B)_{p,q} := X_p \times B_q.
$$
\end{definition}

If $B$ is any simplicial set then putting $h(m)$ or $\Sigma (m)$ in the first
variable, we obtain an inclusion of bisimplicial sets which we denote
$$
\Sigma (m)\boxtimes B \hookrightarrow h(m)\boxtimes B .
$$
Note that these bisimplicial sets are not Segal precategories because they don't
satisfy the globular condition (they are not constant over $0$ in the first
variable). However, that the morphism of simplicial sets
$$
(\Sigma (m)\boxtimes B)_{0/}\hookrightarrow  (h(m)\boxtimes B)_{0/}
$$
is an isomorphism because $\Sigma$ contains all of the vertices.

If $\pA $ is a Segal precategory then a morphism $h(m)\boxtimes B
\rightarrt \pA $ is the same thing as a morphism $B\rightarrt \pA _{m/}$. Similarly,
a morphism $$
\Sigma (m)\boxtimes B \rightarrt \pA 
$$
is the same thing as a morphism
$$
B\rightarrt \pA _{1/} \times _{\pA _0} \ldots \times _{\pA _0}\pA _{1/}.
$$

The morphism of realizations
$$
|\Sigma (m)\boxtimes B | \rightarrt | h(m) \boxtimes B|
$$
is a weak equivalence. To see this note that it is the product of $|B|$ and
$$
|\Sigma (m)| \rightarrt | h(m) |,
$$
and this last morphism is a weak equivalence (it is the inclusion from the
``spine'' of the $m$-simplex to the $m$-simplex; both are contractible).

Suppose $B'\subset B$ is an injection of simplicial sets. Put
$$
\pU := \left( \Sigma (m)\boxtimes B\right) \cup ^{\Sigma (m) \boxtimes B'}
\left( h(m)\boxtimes
B'\right) ,
$$
and
$$
\pV := h(m)\boxtimes B.
$$
We have an injection $\pU \hookrightarrow \pV $.
If $\pA $ is a Segal precategory then a map $\pU \rightarrt \pA $ consists of a
commutative diagram
$$
\begin{diagram}
B' & \rightarr & B \\
\downarr && \downarr \\
\pA _{m/} & \rightarr & \pA _{1/} \times _{\pA _0}\ldots \times _{\pA _0} \pA _{1/} .
\end{diagram}
$$
The inclusion
$$
\Sigma (m) \boxtimes B \rightarrt \pU 
$$
induces a weak equivalence of realizations, because of the fact that the
inclusion
$\Sigma (m) \boxtimes B'\rightarrt h(m)\boxtimes B'$ does. Therefore the morphism
$|\pU |\rightarrt |\pV |$ is a weak equivalence.

We can now interpret our operation $Arr (\pA , m)$ in these terms. Applying the
previous paragraph to the inclusion $\pA _{m/} \hookrightarrow C$, we obtain
an inclusion of bisimplicial sets
$\pU \hookrightarrow \pV $. We get a map $\pU \rightarrt \pA $ corresponding to the diagram
$$
\begin{diagram}
\pA _{m/} & \rightarr & C \\
\downarr && \downarr \\
\pA _{m/} & \rightarr & \pA _{1/} \times _{\pA _0}\ldots \times _{\pA _0} \pA _{1/} .
\end{diagram}
$$
The left vertical arrow is the identity map, the top arrow is $a$ and the
right vertical arrow is $b$. The bottom arrow is the Segal map.

It is easy to see that
$$
Arr (\pA ,m) = \pA  \cup ^{\pU } \pV .
$$
In passing, this proves associativity of the previous formulas for
functoriality of $Arr (\pA ,m)$.

We get
$$
|Arr (\pA ,m) |= |\pA | \cup ^{|\pU |} |\pV |.
$$
Since $|\pU |\rightarrt |\pV |$ is a weak equivalence, this implies the

\begin{lemma}
\mylabel{triviality}
The morphism induced by the above inclusion on realizations,
$$
|\pA | \hookrightarrow |Arr (\pA ,m)|
$$
is a weak equivalence of spaces.
\end{lemma}
\eop

The key observation is the following proposition.

\begin{proposition}
\mylabel{arrangement}
Suppose $\pA $ is a Segal precategory
with $\pA _0 = \ast$ and $\pA _{1/}$ connected. Suppose that $\pA$ is
$(m, k-1)$-arranged and $(p,k)$-arranged for some
$p\neq m$. Then  $Arr (\pA , m)$ is $(p,k)$-arranged and $(m,k)$-arranged.
\end{proposition}
\begin{proof}
Keep the hypotheses of the proposition.  Let $C$ be the cone occuring in the
construction $Arr (\pA , m)$. Denote $\pB := Arr (\pA , m)$. Then the map
$$
a:\pA _{m/} \hookrightarrow C
$$
is weakly equivalent to a map obtained by adding cells of dimension $\geq k$ to
$\pA _{m/}$.  This is by the condition that $\pA $ is $(m, k-1)$-arranged. Let
$h_1,\ldots , h_u$ be the $k$-cells that are attached to $\pA _{m/}$ to give $C$.

We first show that $\pB $ is $(p,k)$-arranged.
Note that $\pB _{p/}$ is obtained from $\pA _{p/}$ by attaching a certain number of
$k$-cells, $h^{p\rightarrt m}_i$ for $i=1,\ldots , u$ indexed by the maps
$p\rightarrt m$ not factoring through the principal edges of $m$; plus some
cells of dimension $\geq k+1$. The higher-dimensional  cells don't have any
effect on the question of whether $\pA $ is $(p,k)$-arranged.

On the other hand, $\pB _{1/}$ is obtained from $\pA _{1/}$ by attaching cells
$h^{1\rightarrt m}_i$ for $i=1,\ldots , u$ and indexed by the maps
$1\rightarrt m$ not factoring through the principal edges. Note here that
these maps cannot be degenerate, thus they are the non-principal edges.

Now $\pB _{1/} \times \ldots \times \pB _{1/}$ (product of $p$-copies) is obtained
from  $\pA _{1/} \times \ldots \times \pA _{1/}$ by adding $k$-cells indexed as
$\nu _{1\rightarrt p}(h^{1\rightarrt m}_i)$ where the indexing
$1\rightarrt p$ are principal edges and $1\rightarrt m$ are non-principal
edges. Then by adding cells of dimension $\geq k+1$ which have no effect on the
question. The notation $\nu _{1\rightarrt p}$ refers to the map
$$
\pA _{1/} \rightarrt \pA _{1/}\times \ldots \times \pA _{1/}
$$
putting the base point (i.e. the degeneracy of the unique point in $\pA _0$) in all
of the factors except the one corresponding to the map $1\rightarrt p$.

For every principal edge $1\rightarrt p$ there is a unique degeneracy
$p\rightarrt 1$ inducing an isomorphism $1\rightarrt 1$ and this establishes
a bijection between principal edges and degeneracies. Thus we may rewrite our
indexing of the $k$-cells attached to the product above as $\nu _{p\rightarrt
1}(h^{1\rightarrt m}_i)$.

Now for every pair $(p\rightarrt 1, 1\rightarrt m)$
the composition $p\rightarrt m$ is a degenerate morphism, not factoring
through a principal edge; and these degenerate morphisms are all different for
different pairs $(p\rightarrt 1, 1\rightarrt m)$. Thus $\pB _{p/}$ contains a
$k$-cell $h^{p\rightarrt m}_i$ for each $i=1,\ldots , u$ and each of these
maps $p\rightarrt m$ (plus possibly other cells for other maps $p\rightarrt
m$ but we don't use these).  Take such a cell $h^{p\rightarrt m}_i$, and look
at its image in $\pB _{1/} \times \ldots \times \pB _{1/}$ by the Segal map.
The projection to any factor $1\rightarrt p$ other than the one which splits
the degeneracy $p\rightarrt 1$, is totally degenerate coming from a
factorization $1\rightarrt 0\rightarrt 1$, hence goes to the unique basepoint.
The projection to the unique factor which splits the degeneracy is just the
cell $h^{1\rightarrt m}_i$. Thus the projection of our cell
$h^{p\rightarrt m}_i$ to the product is exactly the cell
$\nu _{p\rightarrt
1}(h^{1\rightarrt m}_i)$.   This shows that all of the new $k$-cells which
have been added to the product, are lifted as new $k$-cells in $\pB _{p/}$.
Together with the fact that $\pA _{p/} \rightarrt \pA _{1/} \times \ldots \times
\pA _{1/}$ was an isomorphism on $\pi _i$ for $i< k$ and a surjection for $i=k$,
we obtain the same property for $\pB _{p/} \rightarrt \pB _{1/}\times \ldots \times
\pB _{1/}$.  Note that the further $k$-cells which are attached to $\pB _{p/}$ by
morphisms $p\rightarrt m$ other than those we have considered above, don't
affect this property. (In general, attaching $k$-cells to the domain of a map
doesn't affect this property, but attaching cells to the range can affect it,
which was why we had to look carefully at the cells attached to $\pB _{1/}$).
This completes the proof that $\pB $ remains $(p,k)$-arranged.

We now prove that $\pB $ becomes $(m,k)$-arranged.
Note that $\pB _{m/}$ is obtained by first adding on $C$ to $\pA _{m/}$ via the
identity map $m\rightarrt m$; then adding some other stuff which we treat in a
minute. The Segal map for $\pB $ maps this copy of $C$ directly into
$\pA _{1/} \times \ldots \pA _{1/}$. The fact that $C$ is a
mapping cone for the Segal map means that the map
$$
C\rightarrt \pA _{1/} \times \ldots \pA _{1/}
$$
is a homotopy equivalence. In particular, it is bijective on $\pi _i$ for $i<k$
and surjective for $i=k$.

Now $\pB _{m/}$ is obtained from $C$ by adding various cells to $C$ along
degenerate
maps $m\rightarrt m$. The new $k$-cells which are added to $\pB _{1/}
\times \ldots \times \pB _{1/}$ ($m$-factors this time) are lifted to cells in
$\pB _{m/}$ added to $C$ via the degeneracies $m\rightarrt m$ which factor
through a principal edge. The argument is the same as above and we don't repeat
it. We obtain that $\pB $ is $(m,k)$-arranged.

This completes the proof of the proposition.
\end{proof}

\subsection{Connectivity properties of $(m,k)$-arranged precategories}

Suppose
$\pA $ is a Segal precategory with $\pA _0= \ast$ and $\pA _{1/}$ connected, such that
$\pA $ is $(m,k)$-arranged for all $m+k \leq n$. We would like to obtain a map
$\pA \rightarrt \pA '$ to a Segal category, which induces an equivalence in the domain
covered by the $(m,k)$-arrangement conditions which are already known.

Let $\Lambda (m)\subset h(m)$ denote some ``horn'', i.e. a union of all faces but one.
Use the notation of Definition \ref{boxtimesdef}. For any inclusion of simplicial sets
$B '\subset B $ we can set $$ \pU :=(\Lambda (m) \boxtimes \pB ) \cup ^{\Lambda (m) \boxtimes
B '} (h(m) \boxtimes B ')
$$ and
$$
\pV := h(m)\boxtimes B,
$$
we obtain an inclusion of bisimplicial sets $\pU \subset \pV $ such that
$|\pU |\rightarrt |\pV |$ is a weak equivalence. If $\pA $ is a Segal precategory and
$\pU \rightarrt \pA $ a morphism then set $\pA ':= \pA  \cup ^{\pU }\pV $. The morphism
$|\pA |\rightarrt |\pA ' |$ is a weak equivalence. Note that $\pA '$ is again a Segal
precategory because the morphism $\Lambda (m)\rightarrt h(m)$ includes all of the
vertices of $h(m)$.  Finally, note that for $p\leq m-2$
the morphism
$$
\pA _{p/}\rightarrt \pA '_{p/}
$$
is an isomorphism (because the same is true of $\Lambda (m)\rightarrt h(m)$).
This last property allows us to conserve the homotopy type of the smaller
$\pA _{p/}$.

We would like to use operations of the above form, to $(m,k)$-arrange $\pA $.
In order to do this we analyze what a morphism from $\pU $ to $\pA $ means.

For any simplicial set $X$, we can form the simplicial set $\leftbrack X, \pA ]$ with the
property that a map $B\rightarrt \leftbrack X,\pA ]$ is the same thing as a map $X\boxtimes
B\rightarrt \pA $.  If $X\hookrightarrow Y$ is an inclusion obtained by adding
on an $m$-simplex $h(m)$ over a map $Z\rightarrt X$ where $Z\subset h(m)$ is
some subset of the boundary, then
$$
\leftbrack Y,\pA ] = \leftbrack X , \pA ] \times _{\leftbrack Z,\pA ]} \pA _{m/}.
$$
In this way, we can reduce $\leftbrack X,\pA ]$ to a gigantic iterated fiber product of the
various components $\pA _{p/}$.

\begin{lemma}
\label{claimlemma}
In the above situation, 
there exists a cofibration $\pA \rightarrt \pA '$ such that $\pA _{p/}
\rightarrt \pA '_{p/}$ is a weak equivalence for all $p$, and such that 
for any cofibration of simplicial sets $X\subset
Y$ including all of the vertices of $Y$, the morphism
$$
\leftbrack Y,\pA '] \rightarrt \leftbrack X,\pA ']
$$
is a Kan fibration of simplicial sets. 
\end{lemma}
\begin{proof}
In order to construct $\pA '$, we just
``throw in'' everything that is necessary. More precisely, suppose
$B'\subset B$ is a trivial fibration of simplicial sets. A diagram
$$
\begin{diagram}
B' & \rightarr & B \\
\downarr && \downarr \\
\leftbrack Y,\pA '] &\rightarr & \leftbrack X,\pA ']
\end{diagram}
$$
corresponds to a diagram
$$
\begin{diagram}
\pU & \rightarr & X \boxtimes B \\
\downarr && \downarr \\
\pA ' &\rightarr & \pA '
\end{diagram}
$$
with $\pU = (Y\boxtimes B')\cup ^{X\boxtimes B'} (X\boxtimes B)$.
The morphism of bisimplicial sets $\pU \rightarrt X \boxtimes B$ is a weak
equivalence on each vertical column $(\, )_{p/}$.  Therefore we can throw in
to $\pA '$ the pushout along this morphism, without changing the weak equivalence
type of the $\pA _{p/}$.  Note that the new thing is again a Segal precategory
because of the assumption that $X\subset Y$ contains all of the vertices. Keep
doing this addition over all possible diagrams, an infinite number of times,
until we get the required Kan fibration condition to prove the lemma. 
\end{proof}

\begin{theorem}
\mylabel{bound1}
If $\pA $ is a Segal precategory with $\pA _0= \ast$ and $\pA _{1/}$ connected, such that
$\pA $ is $(m,k)$-arranged for all $m+k \leq n$ then there exists a morphism
$\pA \rightarrt \pA '$ such that:
\newline
(1)\,\, the morphism $|\pA |\rightarrt |\pA '|$ is a weak equivalence;
\newline
(2)\,\, $\pA '$ is a Segal
groupoid; and
\newline
(3)\,\, the map of simplicial sets $\pA _{m/} \rightarrt \pA '_{m/}$ induces an
isomorphism on $\pi _i$ for $i+m < n$.
\end{theorem}
\begin{proof}
The answer $\pA '$ is the result of a procedure which we now describe.
At each step of
the procedure,  the construction of the previous
lemma will be applied without necessarily saying so everywhere. 
Thus  we may always assume that
our Segal precategory satisfies the condition of Lemma \ref{claimlemma}.

We have already described above an arranging operation $\pA \mapsto Arr (\pA ,m)$.
We now describe a second arranging operation, under the hypothesis that
$\pA $ satisfies the conclusion of Lemma \ref{claimlemma}.
Fix $m$ and fix a horn $\Lambda (m)\subset h(m)$ (complement of all but one of
the faces, and the face that is left out should be neither the first nor the
last face).  Let $C$ be the cone on the map
$$
\pA _{m/} = \leftbrack h(m), \pA ] \rightarrt \leftbrack \Lambda (m), \pA ].
$$
Thus we have a diagram
$$
\pA _{m/} \rightarrt C \rightarrt \leftbrack \Lambda (m), \pA ].
$$
Note that $\Lambda (m)$ is a gigantic iterated fiber product of various $\pA _{p/}$
for $p<m$. This diagram corresponds to a map $\pU \rightarrt \pA $
where
$$
\pU := (\Lambda (m) \boxtimes C )\cup ^{\Lambda (m)\boxtimes \pA _{m/}} h(m)\boxtimes \pA _{m/}.
$$
Letting $\pV := h(m)\boxtimes C$ we set
$$
Arr2(\pA , m):= \pA \cup ^{\pU } \pV .
$$
Notice first of all that by the previous discussion, the map
$$
|\pA | \rightarrt | Arr2(\pA , m)|
$$
is a weak equivalence of spaces.

We have to try to figure out what effect
$Arr2(\pA ,m)$ has.
We do this under the following hypothesis on the utilisation of this operation:
that $\pA $ is $(m, k-1)$-arranged, and $(p,k)$-arranged for all $p<m$.

The first step is to notice that the fiber product in the expression of
$\leftbrack \Lambda (m), \pA ]$ is a homotopy fiber product, because of the 
condition of Lemma \ref{claimlemma} which is 
imposed on $\pA $. Furthermore the elements in this fiber
product all satisfy the Segal condition up to $k$ (bijectivity on $\pi _i$
for $i<k$ and surjectivity for $\pi _k$). Thus the morphism
$$
\leftbrack \Lambda (m) , \pA ]
\rightarrt \pA _{1/} \times _{\pA _0} \ldots \times _{\pA _0} \pA _{1/}
$$
(the Segal fiber product for $m$, on the right)
is an isomorphism on $\pi _i$ for $i<k$ and a surjection for $i=k$.
Thus when we add to $\pA _{m/}$ the cone $C$, we obtain the condition of being
$(m,k)$-arranged.

By hypothesis, $\pA $ is $(m,k-1)$-arranged, in particular the map
$\pA _{m/} \rightarrt C$ is an isomorphism on $\pi _i$ for $i<k-1$ and
surjective for $\pi _{k-1}$.  Thus $C$ may be viewed as obtained from $\pA _{m/}$
by adding on cells of dimension $\geq k$.  Therefore, for all $p$ the morphisms
$\pU _{p/}\rightarrt \pV _{p/}$ are homotopically obtained by addition of cells of
dimension $\geq k$.

From the previous paragraph, some extra cells of dimension $\geq k$ are added to
various $\pA _{p/}$ in the process. This doesn't spoil the condition of being
$(p,k)$-arranged wherever it exists.  However, the major advantage
of this second operation is that the $\pA _{p/}$ are left unchanged for $p\leq m-2$.  This
is because all $p$-faces of the $m$-simplex are then contained in the horn
$\Lambda (m)$.

We review the above results. First, the hypotheses on $\pA $ were:
\newline
(a)\,\, that $\pA $ satisfies the lifting condition of Lemma \ref{claimlemma};
\newline
(b)\,\, that $\pA $ is $(p,k)$-arranged for $p<m$, and $(m,k-1)$-arranged.
We then obtain a construction $Arr2(\pA ,m)$ with the following properties:
\newline
(1)\,\, the map
$\pA _{p/}\rightarrt Arr2(\pA ,m)_{p/}$ is an isomorphism for $p\leq m-2$;
\newline
(2)\,\, for any $p$ the map $\pA _{p/}\rightarrt Arr2(\pA ,m)_{p/}$ induces an
isomorphism on $\pi _i$ for $i<k-1$;
\newline
(3)\,\, if $\pA $ is $(p,k)$-arranged for any $p$ then $Arr2(\pA ,m)$ is also
$(p,k)$-arranged;
\newline
(4)\,\, and $Arr2(\pA ,m)$ is $(m,k)$-arranged.

{\em Remark:} at $m=2$ the operations $Arr(\pA ,2)$ and $Arr2(\pA ,2)$ coincide.

With an infinite series of applications of the construction $Arr2(\pA ,m)$ and the
replacement operation of Lemma \ref{claimlemma}
we can prove Theorem \ref{bound1}. The reader may
do this as an exercise or else read the explanation below.

Take an array of dots, one for each $(p,k)$. Color the dots green if $\pA $ is
$(p,k)$-arranged, red otherwise (note that one red dot in a column implies
red dots everywhere above).  We do a sequence of operations of the form
of Lemma \ref{claimlemma} (which doesn't change the homotopy type levelwise) 
and then $Arr2(\pA ,m)$.
When we do this, change the colors of the dots appropriately.

Also mark an
$\times$ at any dot $(p,k)$ such that the $\pi _i(\pA _{p/})$ change for
any $i \leq k-1$.  (Keep any $\times$ which are marked, from one step to
another).  If a dot $(p,k)$ is never marked with a $\times$ it means that
the $\pi _i(\pA _{p/})$ remain unchanged for $i<k$.

We don't color the dots $(1,k)$ but we still might mark an
$\times$.

Suppose the dot $(m,k)$ is red,
the dots $(p,k)$ are green for $p<m$ and the dot $(m,k-1)$ is green.
Then apply the replacement operation of Lemma \ref{claimlemma} 
and the operation $Arr2(\pA 	,m)$.
This has the following effects.  Any green dot $(p,j)$ for $p\leq m-2$
(and arbitrary $j$) remains green. The dot $(m-1,k)$ remains green.
However, the dot $(m-1, k+1)$ becomes red.
The dot $(m,k)$ becomes green.
The dots $(m-1, k)$ and $(m,k)$, as well as all $(p,k)$ for $p>m$, are marked
with an $\times$. The dots above these are also marked with an $\times$ but
no other dots are (newly)  marked with an $\times$.

In the situation of Theorem \ref{bound1}, we start with green dots at $(p,k)$
for $p+k \leq n$. We may as well assume that the rest of the dots are colored
red. Start with $(m,k)=(n+1, 0)$ and apply the procedure of the previous
paragraph. The dot $(n+1,0)$ becomes green, the dot $(n,0)$ stays green,
and the dots $(n, 0), (n+1, 0), \ldots $ are marked with an $\times$. Continue
now at $(n,1)$ and so on. At the end we have made all of the dots $(p,k)$ with
$p+k=n+1$ green, and we will have marked with an $\times$ all of the dots
$(p,k)$ with $p+k=n$ (including the dot $(1,n-1)$; and also all of the dots
above this line).

We can now iterate the procedure. We successively get green dots on each of the
lines $p+k= n+j$ for $j=1,2,3,\ldots$. Furthermore, no new dots will be marked
with a $\times$.  After taking the union over all of these iterations, we obtain
an $\pA '$ which is $(p,k)$-arranged for all $(p,k)$.
Thus $\pA '$ is a Segal category.

Note that the morphism $|\pA |\rightarrt |\pA '|$ is  a weak equivalence of spaces.

By looking at which dots are
marked with an $\times$, we find that the morphisms
$$
\pA _{p/} \rightarrt \pA '_{p/}
$$
induce isomorphisms on $\pi _i$ whenever $i < n-p$.   This completes the proof
of the theorem.
\end{proof}

\subsection{Iteration}

The following corollary to Theorem \ref{bound1} says
that in order to calculate the $n$-type of
$\Omega |\pA |$ we just have to change $\pA $ by pushouts preserving the weak
equivalence type of $|\pA |$ in such a way that $\pA $ is $(m,k)$-arranged for all
$m+k \leq n+2$.

\begin{corollary}
\mylabel{bound2}
Suppose $\pA $ is a Segal precategory with $\pA _0= \ast$ and $\pA _{1/}$ connected, such that
$\pA $ is $(m,k)$-arranged for all $m+k \leq n$. Then the natural morphism
$$
|\pA _{1/}| \rightarrt \Omega | \pA |
$$
induces an isomorphism on $\pi _i$ for $i < n-1$.
\end{corollary}
{\em Proof:}
Use Theorem \ref{bound1} to obtain a morphism $\pA \rightarrt \pA '$ with the
properties stated there (which we refer to as (1)--(3)).  We have a diagram
$$
\begin{diagram}
|\pA _{1/}| & \rightarr & \Omega | \pA | \\
\downarr & & \downarr \\
|\pA '_{1/}| & \rightarr & \Omega | \pA '|
\end{diagram} .
$$
By property (1) the vertical morphism on the right is a weak equivalence.
By property (2) and Theorem \ref{segal} the morphism on the bottom is a weak
equivalence. By property (3) the vertical morphism on the right induces
isomorphisms on $\pi _i$ for $i<n-1$. This gives the required statement.
\eop

\begin{corollary}
\mylabel{thinkIcan}
Fix $n$, and suppose $\pA $ is a Segal precategory with $\pA _0 = \ast$ and $\pA _{1/}$
connected. By applying the operations $\pA  \mapsto Arr (\pA , m)$ for various $m$, a
finite number of times (less than $(n+2)^2$) in a predetermined way,
we can effectively get to a morphism of Segal precategories $\pA \rightarrt \pB $ such
that
$$
|\pA | \rightarrt |\pB |
$$
is a weak equivalence of spaces, and such that $\pB $ is $(m,k)$-arranged for all
$m+k \leq n+2$. Furthermore $\pB _0 = \ast$ and $\pB _{1/}$ is
connected. 
\end{corollary}
{\em Proof:}
By Corollary \ref{triviality} any successive application of the operations
$\pA \mapsto Arr (\pA , m)$ yields a morphism $|\pA |\rightarrt |\pB |$ which is a weak
equivalence of spaces. By Proposition \ref{arrangement} it suffices, for
example,
to successively apply $Arr (\pA ,i)$ for $i= 2, 3,\ldots , n+2$ and to repeat this
$n+2$ times. These operations preserve connectedness of the pieces in degree $1$,
so $\pB _{1/}$ is connected.  
\eop

\begin{corollary}
\mylabel{OmegaDone}
Fix $n$, and suppose $\pA $ is a Segal precategory with $\pA _0 = \ast$ and $\pA _{1/}$
connected. Let $\pB $ be the result of the operations of Corollary
\ref{thinkIcan}. Then the $n$-type of the simplicial set $\pB _{1/}$ is equivalent
to the $n$-type of $\Omega |\pA  |$.
\end{corollary}
{\em Proof:}
Apply Corollaries \ref{bound2} and \ref{thinkIcan}.
\eop

\section{Computing the loop space}
\label{computingLoop}

Suppose $X$ is a simplicial set with $X_0 = X_1 = \ast$, and with finitely
many nondegenerate simplices.  Fix $n$. We will obtain, by iterating an
operation
closely related to the operation $\gen $ of Chapter \ref{genrel1}, a finite complex
representing the $n$-type of $\Omega X$. 

Let $\pA $ be $X$
considered as a Segal precategory constant in the second variable, in other words
$$
\pA _{p,k} := X_p.
$$
Apply the arrangement process, iterated as in Corollary \ref{thinkIcan}. 

\begin{corollary}
\mylabel{OmegaDone2}
Fix $n$, and suppose $X$ is a simplicial set with finitely many nondegenerate
simplices, with $X_0=X_1 = \ast$. Let $\pA $ be $X$
considered as a Segal precategory. Let $\pB $ be the result of the operations of
Corollary \ref{thinkIcan}. Then the $n$-type of the simplicial set $\pB _{1/}$ is
equivalent to the $n$-type of $\Omega X$.
 \end{corollary}
{\em Proof:}
An immediate restatement of \ref{OmegaDone}.
\eop

{\em Remark:} Any finite region of the Segal precategory $\pB $ is effectively
computable. In fact it is just an iteration of operations pushout and mapping
cone, arranged in a way which depends on combinatorics of simplicial sets.
Thus the $n+1$-skeleton of the simplicial set $\pB _{1/}$ is effectively
calculable (in fact, one could bound the number of simplices in $\pB _{1/}$).

\begin{corollary}
Fix $n$, and suppose $X$ is a simplicial set with finitely many
nondegenerate simplices, with $X_0=X_1 = \ast$. Then we
can effectively calculate $H_i (\Omega |X|, \zz )$ for $i\leq n$.
\end{corollary}
{\em Proof:}
Immediate from above.
\eop

In some sense this corollary is the ``most effective'' part of the present
argument, since we can get at the calculation after a bounded number of
easy steps of the form $\pA \mapsto Arr (\pA ,m)$.

We describe how to use the above description of $\Omega X$ inductively to
obtain the $\pi _i(X)$.  This seems to be a new algorithm, different from those
of E. Brown \cite{EBrown} and Kan--Curtis \cite{Kan1} \cite{Kan2}
\cite{Curtis}.

There is an unboundedness to the resulting algorithm, coming essentially from a
problem with $\pi _1$ at each stage. Even though we know in advance that the
$\pi _1$ is abelian, we would need to know ``why'' it is abelian in a
precise way
in order to specify a strategy for making $\pA _{1/}$ connected at the appropriate
place in the loop.  In the absence of a particular description of the proof we
are forced to say ``search over all proofs'' at this stage. See Subsection \ref{sec-godement}
for further discussion.

\subsection{Getting $\pA _{1/}$ to be connected}
\label{gettingconnected}

In the general situation, we have to tackle the problem of computation of a
fundamental group using generators and relations, known to be undecideable in 
general. Some sub-cases can still be treated effectively. 

The first question is how to arrange $\pA $ on the level of $\tau
_{\leq 1}(\pA )$. 

We define operations $Arr ^{0\, {\rm only}}(\pA ,m)$ and
$Arr ^{1\, {\rm only}}(\pA ,m)$. These consist of doing the operation $Arr (\pA ,m)$
but instead of using the entire  mapping cone $C$, only adding on $0$-cells
to $\pA _{m/}$ to get a surjection on $\pi _0$; or only adding on $1$-cells to get
an injection on $\pi _0$.  Note in the second case that we {\em don't} add
extra $0$-cells. This is an important point, because if we added further
$0$-cells every time we added some $1$-cells, the process would never stop.

To define  $Arr ^{0\, {\rm only}}(\pA ,m)$, use the same construction as for $Arr
(\pA ,M)$ but instead of setting $C$ to be the mapping cone, we put
$$
C':= \pA _{m/} \cup sk_0(\pA _{1/} \times \ldots \times \pA _{1/}.
$$
Here $sk_0$ denotes the $0$-skeleton of the simplicial set, and $im$ means the
image under the Segal map. Let $C\subset C'$ be a subset where we choose
only one point for each  connected component of the product.
With this $C$ the
same construction as previously gives $Arr ^{0\, {\rm only}}(\pA ,m)$.

With the subset $C\subset C'$ chosen as above (note that this choice can
effectively be made) the resulting simplicial set
$$
p\mapsto \pi _0\left( Arr ^{0\, {\rm only}}(\pA ,m)_{p/} \right)
$$
may be described only in terms of the simplicial set
$$
p\mapsto \pi _0(\pA _{p/}).
$$
That is to say, this operation $Arr ^{0\, {\rm only}}(\pA ,m)$ commutes with the
operation of componentwise applying $\pi _0$. We formalize this as
$$
\tau _{\leq 1}Arr ^{0\, {\rm only}}(\pA ,m)
=\tau _{\leq 1}Arr ^{0\, {\rm only}}(\tau _{\leq 1}\pA ,m).
$$

To define $Arr ^{1\, {\rm only}}(\pA ,m)$, let $C$ be the cone of the map from
$\pA _{m/}$ to
$$
im(\pA _{m/}) \cup sk_1(\pA _{1/} \times \ldots \times \pA _{1/})^o
$$
where  $sk_1(\pA _{1/} \times \ldots \times \pA _{1/})^o$ denotes the union of
connected components of the $1$-skeleton of the product, which touch
$im(\pA _{m/})$. In this case, note that the inclusion
$$
\pA _{m/} \hookrightarrow C
$$
is $0$-connected (all connected components of $C$ contain elements of $\pA _{m/}$).
Using this $C$ we obtain the operation $Arr ^{1\, {\rm only}}(\pA ,m)$.
It doesn't introduce any new connected components in the
new simplicial sets $\pA '_{p/}$, but may connect together some components which
were disjoint in $\pA _{p/}$.

Again, the operation $Arr ^{1\, {\rm only}}(\pA ,m)$ commutes with truncation:
we have
$$
\tau _{\leq 1}Arr ^{1\, {\rm only}}(\pA ,m)
=\tau _{\leq 1}Arr ^{1\, {\rm only}}(\tau _{\leq 1}\pA ,m).
$$

Our goal in this section is to find a sequence of operations which makes $\tau
_{\leq 1}(\pA )_1$ become trivial (equal to $\ast$).  In view of this, and the
above commutations, we may henceforth work with simplicial sets (which we
denote $U = \tau _{\leq 1}\pA $ for example) and use the above operations
followed by
the truncation $\tau _{\leq 1}$ as modifications of the simplicial set $U $.
We try to obtain $U _1=\ast$.  This corresponds to making $\pA _{1/}$  connected.

Our operations have the following interpretation. The operation
$$
U \mapsto \tau _{\leq 1}Arr ^{0\, {\rm only}}(U ,2)
$$
has the effect of formally adding to $\pU _1$ all binary products of pairs of
elements in $U _1$.  (We say that a binary product of $u,v\in U _1$ is defined if
there is an element $c\in U_2$ with principal edges $u$ and $v$ in $U_1$; the
product is then the image $w$ of the third edge of $c$).

The operation
$$
U\mapsto \tau _{\leq 1}Arr ^{1\, {\rm only}}(U,2)
$$
has the effect of identifying $w$ and $w'$ any time both $w$ and $w'$ are
binary products of the same elements $u,v$.

The operation
$$
U\mapsto \tau _{\leq 1}Arr ^{0\, {\rm only}}(U,3)
$$
has the effect of introducing, for each triple $(u,v,w)$, the
various binary products one can make (keeping the same order) and giving  a
formula
$$
(uv)w = u(vw)
$$
for certain of the binary products thus introduced.

It is somewhat unclear whether blindly applying the composed operation
$$
U\mapsto \tau _{\leq 1}
Arr ^{1\, {\rm only}}(\tau _{\leq 1}Arr ^{0\, {\rm only}}(U,3),2)
$$
many times must automatically lead to $U_1=\ast$ in case the actual fundamental
group is trivial.   This is because in the process of adding
the associativity, we also add in some new binary products; to which
associativity might then have to be applied in order to get something trivial,
and so on.

If the above doesn't work, then we may need a slightly revised version of the
operation $Arr ^{0\, {\rm only}}(U,3)$ where we add in only certain triples
$u,v,w$.  This can be accomplished by choosing a subset of the original $C$ at
each time.  Similarly for the
$Arr ^{0\, {\rm only}}(U,2)$ for binary products.
We now obtain a situation where
we have operations which effect the appropriate changes on $U$ corresponding to
all of the various possible steps in an elementary proof that the associative
unitary monoid generated by generators $U_1$ with relations $U_2$, is trivial.
Thus if we have an elementary proof that the associative unitary monoid
generated by $U_1$ with relations $U_2$ is trivial, then we can read off
from the
steps in the proof, the necessary sequence of operations to apply to get
$U_1=\ast$.  On the level of $\pA $ these same steps will result in a new $\pA $ with
$\pA _{1/}$ connected.

In our case we are interested in the {\em group completion} of the monoid: we
want to obtain the condition of being a Segal groupoid not just a Segal
category. It is possible that the simplicial set $X$ we start with would yield
a monoid which is not a group, when the above operations are applied. To fix
this, we take note of another operation which can be applied to $\pA $ which
doesn't affect the weak type of the realization, and which guarantees that, when
the monoid $U$ is generated, it becomes a group.

Let $I$ be the category with two objects and one morphism $0\rightarrt 1$, and
let $\overline{I}$ be the category with two objects and an isomorphism between
them. Consider these as Segal categories (taking their nerve as
bisimplicial sets constant in the second variable).  Note that $|I|$ and $|
\overline{I}|$ are both contractible, so the obvious inclusion
$I\hookrightarrow \overline{I}$ induces an equivalence of realizations.

The bisimplicial set $\overline{I}$ is just that which is represented by
$(1,0)\in \Delta \times \Delta$. Thus for a Segal precategory $\pA $,
if $f\in \pA _{1,0}$ is an object of $\pA _{1/}$ (a ``morphism'' in $\pA $) then it
corresponds to a morphism $I\rightarrt \pA $.
Set
$$
\pA ^f:= \pA  \cup ^{I} \overline{I}.
$$
Now the morphism $f$ is strictly invertible in the precategory $\tau _{\leq
1}(\pA ^f)$ and in particular, when we apply the operations described above, the
image of $f$ becomes invertible in the resulting category.  If $\pA _0=\ast$
(whence $\pA ^f_0=\ast$ too) then
the image of $f$ becomes invertible in the resulting monoid.
Note finally that
$$
|\pA  | \rightarrt | \pA ^f| = | \pA | \cup ^{|I|} |\overline{I}|
$$
is a weak equivalence.
In fact we want to invert all of the $1$-morphisms.
Let
$$
\pA ' := \pA  \cup ^{\bigcup _fI}\left( \bigcup _f \overline{I} \right)
$$
where the union is taken over all $f\in \pA _{1,0}$. Again
$|\pA |\rightarrt | \pA '|$ is a weak equivalence. Now, when we apply the previous
procedure to $\tau _{\leq 1} (\pA ')$ giving a category $U$ (a monoid if $\pA $ had
only one object), all morphisms coming from $\pA _{1,0}$ become invertible. Note
that the morphisms in $\pA '$, i.e. objects of $\pA '_{1,0}$, are either morphisms in
$\pA $ or their newly-added inverses. Thus all of the morphisms coming from
$\pA '_{1,0}$ become invertible in the category $U$. But it is clear from the
operations described above that $U$ is generated by the morphisms in
$\pA '_{1,0}$. Therefore $U$ is a groupoid. In the case of only one object,
$U$ becomes a group.

By Segal's theorem we then have $U= \pi _1(| \pA |)$. If we know for some reason
that $|\pA |$ is simply connected, then $U$ is the trivial group. More
precisely, search for a proof that $\pi _1= 1$, and when such a proof is found,
apply the corresponding series of operations to
$\tau _{\leq 1} (\pA ')$ to obtain $U=\ast$. Applying the
operations to $\pA '$ upstairs, we obtain a new $\pA ''$ with $|\pA ''|\cong |\pA '|\cong
|\pA |$ and  $\pA ''_{1/}$ connected.

Another way of looking at this is to say that every time one needs to take the
inverse of an element in the proof that the group is trivial, add on a copy of
$\overline{I}$ over the corresponding copy of $I$.

\subsection{The case of finite homotopy groups}

We first present our algorithm for the case of finite homotopy groups.
Suppose we want to calculate $\pi _n(X)$. We assume known that the
$\pi _i(X)$ are finite for $i\leq n$.

\noindent
{\em Start:}
Fix $n$ and
start with a simplicial set $X$ containing a finite number of nondegenerate
simplices. Suppose we know that $\pi
_1(X,x)$ is a given finite group; record this group, and set $Y$ equal to
the corresponding covering space of $X$. Thus $Y$ is simply connected. Now
contract out a maximal tree to obtain $Z$ with $Z_0=\ast$.

\noindent
{\em Step 1.} Let $\pA _{p,k}:= Z_p$ be the
corresponding Segal precategory. It has only one object.

\noindent
{\em Step 2.} Let $\pA '$ be the coproduct of $\pA $ with one copy of the nerve of
the category $\overline{I}$ (containing two isomorphic objects), for each
morphism $I\rightarrt \pA $ (i.e. each point in $\pA _{1,0}$).

\noindent
{\em Step 3.} Apply the procedure of Subsection \ref{gettingconnected} 
to obtain a morphism $\pA '\rightarrt
\pA ''$ with $\pA ''_{1/}$ connected, and inducing a weak equivalence on
realizations.  (This step can only be bounded if we have a specific proof that
$\pi _1(Y,y)=1$).

\noindent
{\em Step 4.} Apply the procedure of Corollary \ref{OmegaDone2}
and Theorem \ref{bound1} to obtain a morphism $\pA ''
\rightarrt  \pB $ inducing a weak equivalence on realizations, such that $\pB $ is a
Segal groupoid.
Note that the
$n-1$-type of $\pB _{1/}$ is effectively calculable (the non-effective parts of the
proof of Theorem \ref{bound1} served only to prove the properties in question).  
By Segal's theorem,
$$
|\pB _{1/}| \sim \Omega | \pB | \sim \Omega | Y|,
$$
which in turn is the
connected component of $\Omega | X|$. Thus
$$
\pi _n(|X|)= \pi _{n-1}(|\pB _{1/}|).
$$

\noindent
{\em Step 5.} Go back to the {\em Start} with the new $n$ equal to the old
$n-1$, and the new $X$ equal to the simplicial set $\pB _{1/}$ above. The new
fundamental group is known to be abelian (since it is $\pi _2$ of the previous
$X$). Thus we can calculate the new fundamental group as $H_1(X)$ and, under our
hypothesis, it will be finite.

Keep repeating the procedure until we get down to $n=1$ and have recorded the
answer.

\subsection{How to get rid of free abelian groups in $\pi _2$}

In the case where the higher homotopy groups are infinite (i.e. they contain
factors of the form $\zz ^a$) we need to do something to get past these
infinite groups. If we go down to the case where $\pi _1$ is infinite, then
taking the universal covering no longer results in a finite complex. We prefer
to avoid this by tackling the problem at the level of $\pi _2$, with a
geometrical argument. Namely, if $H^2(X, \zz)$ is nonzero then we can take a
class there as giving a line bundle, and take the total space of the
corresponding $S^1$-bundle. This amounts to taking the fiber of a map
$X\rightarrt K(\zz , 2)$.  This can be done explicitly and effectively,
resulting again in a calculable finite complex. In the new complex we will have
reduced the rank of $H_2(X, \zz )= \pi _2(X)$ (we are assuming that $X$ is
simply connected).

The original method of E. Brown \cite{EBrown} for effectively calculating the
$\pi _i$ was basically to do this at all $i$. The technical problems in
\cite{EBrown} are caused by the fact that one doesn't have a finite complex
representing $K(\zz , n)$.  In the case $n=2$ we don't have these technical
problems because we can look at circle fibrations and the circle is a finite
complex.  For this section, then, we are in some sense reverting to an easy
case of \cite{EBrown} and not using the Seifert-Van Kampen technique.

Suppose $X$ is a simplicial set with finitely many nondegenerate simplices,
and suppose $X_0=X_1 = \ast$.  We can calculate $H^2(X, \zz )$ as
the kernel of the differential
$$
d: \zz ^{X'_2} \rightarrt \zz ^{X'_3}.
$$
Here $X'_i$ is the set of nondegenerate $i$-simplices.
(Note that a basis of this kernel can effectively be computed using Gaussian
elimination). Pick an element $\beta$ of this basis, which is a collection of
integers $b_t$ for each $2$-simplex (i.e. triangle) $t$.  For each triangle $t$
define an $S^1$-bundle $L_t$ over $t$ together with trivialization
$$
L_t |_{\partial t} \cong \partial t \times S^1.
$$
To do this, take $L_t = t \times S^1$ but change the trivialization along the
boundary by a bundle automorphism
$$
\partial t \times S^1 \rightarrt \partial t \times S^1
$$
obtained from a map $\partial t \rightarrt S^1$ with winding number $b_t$.
 Let $L^{(2)}$ be the $S^1$-bundle over the $2$-skeleton of $X$ obtained by
glueing together the $L_t$ along the trivializations over their boundaries. We
can do this effectively and obtain $L^{(2)}$ as a simplicial set with a finite
number of nondegenerate simplices.

The fact that $d(\beta )=0$ means that for a $3$-simplex $e$, the restriction
of $L^{(2)}$ to $\partial e$ (which is topologically an $S^2$) is a trivial
$S^1$-bundle. Thus $L^{(2)}$ extends to an $S^1$-bundle $L^{(3)}$ on
the $3$-skeleton of $X$. Furthermore, it can be extended across any simplices
of dimension $\geq 4$ because all $S^1$-bundles on $S^k$ for $k\geq 3$, are
trivial ($H^2(S^k, \zz )= 0$). We obtain an $S^1$-bundle $L$ on $X$.
By subdividing things appropriately (including possibly subdividing $X$)
we can assume that $L$ is a simplicial set with a finite number of
nondegenerate simplices.  It depends on the choice of basis element $\beta$, so
call it $L(\beta )$.

Let
$$
T=L(\beta _1)\times _X \ldots \times _X L(\beta _r)
$$
where $\beta _1,\ldots , \beta _r$ are our basis elements found above.
It is a torus bundle with fiber $(S^1)^r$. The long exact homotopy sequence for
the map $T\rightarrt X$ gives
$$
\pi _i (T)= \pi _i (X), \;\;\; i \geq 3;
$$
and
$$
\pi _2(T) = \ker (\pi _2(X) \rightarrt \zz ^r).
$$
Note that $\zz ^r$ is the dual of $H^2(X , \zz )$ so the kernel $\pi _2(T)$ is
finite. Finally, $\pi _1(T)=0$ since the map  $\pi _2(X) \rightarrt \zz ^r$
is surjective.

Note that we have a proof that $\pi _1(T)=0$.

\subsection{The general algorithm}

Here is the general situation. Fix $n$. Suppose $X$ is a simplicial set with
finitely many nondegenerate simplices, with $X_0=\ast$ and with a proof that $\pi
_1(X)=\{1\}$.  We will calculate $\pi _i(X)$ for $i\leq n$.

\noindent
{\em Step 1.}\,
Calculate (by Gaussian elimination) and record $\pi _2(X) = H_2(X, \zz )$.

\noindent
{\em Step 2.}\,
Apply the operation described in the previous subsection above, to obtain a new
$T$ with $\pi _i(T)=\pi _i(X)$ for $i\geq 3$, with $T_0 = \ast$, with $\pi
_1(T)=1$, and with $\pi _2(T)$ is finite.

Let $\pA $ be the  Segal precategory corresponding to $T$.

\noindent
{\em Step 3.}\,
Use the discussion of Section \ref{gettingconnected} 
to obtain a morphism $\pA \rightarrt \pA '$ inducing a
weak equivalence of realizations, such that $\pA '_{1/}$ is connected. For this
step we need a proof that $\pi _1(T)=1$.  In the absence of a specific (finite)
proof, search over all proofs.

\noindent
{\em Step 4.}\,
Use Corollary \ref{thinkIcan} to replace $\pA '$ by a Segal  precategory $\pB $ with
$|\pA '|\rightarrt |\pB |$ a weak equivalence, such that
the $n-1$-type of $\pB _{1/}$ is equivalent to $\Omega |\pB |$ which in turn is
equivalent to $\Omega |X|$. Let $Y= \pB _{1/}$ as a simplicial set.

Note that $Y$ is connected and $\pi _1(Y)$ is finite, being equal to
$\pi _2(T)$. We have $\pi _i(X) = \pi _{i-1}(Y)$ for
$3\leq i \leq n$.

\noindent
{\em Step 5.}\,
Choose a universal cover of $Y$, and mod out by a maximal
tree in the $1$-skeleton to obtain a simplicial set $Z$, with finitely many
nondegenerate simplices, with $Z_0 = \ast$, and with a proof that $\pi _1(Z)=
1$.  We have $\pi _i(X) = \pi _{i-1}(Z)$ for
$3\leq i \leq n$.

Go back to the beginning of the algorithm and plug in $(n-1)$
and $Z$. Keep doing this until, at the step where we calculate $\pi _2$ of the
new object, we end up having calculated $\pi _i(X)$ as desired.

\subsection{Proofs of Godement}
\label{sec-godement}

We pose the following question: how could one obtain, in the process of
applying the above algorithm, an explicit proof that at each stage the
fundamental group (of the universal cover $Z$ in step 5) is trivial? This could
then be plugged into the machinery above to obtain an explicit strategy,
thus we would avoid having to try all possible strategies. To do this we would
need an explicit proof that $\pi _1(Y)$ is finite in step 4, and this in turn
would be based on a proof that $\pi _1(Y)=\pi _2(T)$ as well as a proof of the
Godement property that $\pi _2(T)$ is abelian.

\section{Example: $\pi _3(S^2)$}
\label{s2example}

The story behind the preprint \cite{effective} 
was that Ronnie Brown came by Toulouse for Jean Pradines'
retirement party, and we were discussing Seifert-Van Kampen. He pointed out that
the result of \cite{svk} didn't seem to lead to any actual calculations.
After that, I tried to use that technique (in its simplified
Segal-categoric version) to calculate $\pi _3(S^2)$. It was apparent from this
calculation that the process was effective in general.

We describe here what happens for calculating $\pi _3(S^2)$.
We take as simplicial model a simplicial set
with the basepoint as unique $0$-cell $\ast$ and with one nondegenerate simplex
$e$ in degree $2$.  Note that this leads to many degenerate simplices in degrees
$\geq 2$ (however there is only one degenerate simplex which we denote $\ast$ in
degree $1$).

We follow out what happens in a language of cell-addition. Thus we don't feel
required to take the whole cone $C$ at each step of an operation $Arr (\pA ,m)$;
we take any addition of cells to $\pA _{m/}$ lifting cells in $\pA _{1/}\times \ldots
\times \pA _{1/}$.

We keep the notation $\pA $ for the result of each operation (since our discussion
is linear, this shouldn't cause too  much confusion).

The first step is to $(2,0)$-arrange $\pA $. We do this by adding a $1$-cell
joining the two $0$-cells in $\pA _{2/}$, in an operation of type $Arr (\pA ,
2)$. Note that both $0$-cells map to the same point $\pA _{1/}\times \pA _{1/} =
\ast$.
The first result of this is to add on $1$-cells in the $\pA _{m/}$ connecting all
of the various degeneracies of $e$, to the basepoint. Thus the $\pA _{m/}$ become
connected.
Additionally we get a new $1$-cell added onto to $\pA _{1/}$
corresponding to the third face $(02)$.  Furthermore, we obtain all images of
this cell by degeneracies $m\rightarrt 1$. Thus we get $m$ circles
attached to  the pieces which became connected in the first part of this
operation. Now each $\pA _{m/}$ is a wedge of $m$ circles.

In particular note that $\pA $ is now $(m,1)$-arranged for all $m$.

The next step is to $(2,2)$-arrange $\pA $. To do this, note that the Segal map
is
$$
S^1 \vee S^1 = \pA _{2/} \rightarrt \pA _{1/} \times \pA _{1/} = S^1 \times S^1.
$$
To arrange this map we have to add a $2$-cell to $S^1 \vee S^1$
with attaching map the commutator relation. Again, this has the result of
adding on $2$-cells to all of the $\pA _{m/}$ over the pairwise commutators of the
loops. Furthermore, we obtain an extra $2$-cell added onto $\pA _{1/}$ via the
edge $(02)$. The attaching map here is the commutator of the generator with
itself, so it is homotopically trivial and we have added on a $2$-sphere.
(Note in passing that this $2$-sphere is what gives rise to the class of the
Hopf map).
 Again, we obtain the  images of this $S^2$ by all of the degeneracy
maps $m\rightarrt 1$. Now
$$
\pA _{1/} = S^1 \vee S^2,
$$
$$
\pA _{2/} = (S^1 \times S^1)\vee S^2\vee S^2,
$$
and in general $\pA $ is $(m, 2)$-arranged for all $m$.  Looking forward to the next section, we see that adding
$3$-cells to $\pA _{m/}$ for $m\geq 3$  in the appropriate way as described in the
proof of \ref{bound1}, will end up resulting in the addition of $4$-cells (or
higher) to $\pA _{1/}$ so this no longer affects the $2$-type of $\pA _{1/}$. Thus
(for the purposes of getting $\pi _3(S^2)$) we may now ignore the $\pA _{m/}$ for
$m\geq 3$.

The remaining operation is to $(2,3)$-arrange $\pA $. For this, look at the
Segal map
$$
\pA _{2/} = (S^1 \times S^1)\vee S^2\vee S^2 \rightarrt
$$
$$
\pA _{1/}\times \pA _{1/} = (S^1 \vee S^2)\times (S^1 \vee S^2).
$$
Let $C$ be the mapping cone on this map. Then we end up attaching one copy
of $C$ to $\pA _{1/}$ along the third edge map $\pA _{2/} \rightarrt \pA _{1/}$.
This gives the answer for the $2$-type of $\Omega S^2$:
$$
\tau _{\leq 2}(\Omega S^2) =
\tau _{\leq 2}\left(
(S^1 \vee S^2)\cup ^{(S^1 \times S^1)\vee S^2\vee S^2}C \right) .
$$
To calculate $\pi _2(\Omega S^2)$ we revert to a homological formulation
(because it isn't easy to ``see'' the cone $C$).  In homology of degree $\leq
2$,  the above Segal map
$$
(S^1 \times S^1)\vee S^2\vee S^2 \rightarrt
(S^1 \vee S^2)\times (S^1 \vee S^2)
$$
is an isomorphism. Thus the map $\pA _{2/} \rightarrt C$ is an isomorphism on
homology in degrees $\leq 2$, and adding in a copy of $C$ along $\pA _{2/}$
doesn't change the homology. Thus
$$
H_2(\Omega S^2) = H_2 (S^1 \vee S^2) = \zz .
$$
Noting that (as we know from general principles) $\pi _1(\Omega S^2) = \zz
$ acts
trivially on $\pi _2(\Omega S^2)$ and $\pi _1$ itself has no homology in degree
$2$, we get that $\pi _2(\Omega S^2)= H_2(\Omega S^2) = \zz$.

{\em Exercise:} Calculate $\pi _4(S^3)$ using the above method.

{\em Remark:} our above recourse to homology calculations suggests that it
might be interesting to do pushouts and the operation $Cat$ in the context of
simplicial chain complexes.

\subsection{Seeing Kan's simplicial free groups}
Using the above procedure, we can actually see how Kan's simplicial free groups
arise in the calculation for an arbitrary simplicial set $X$. They
arise just from a first stage where we add on $1$-cells. Namely, if in doing the
procedure $Arr (\pA ,m)$ we replace $C$ by a choice of $1$-cell joining any two
components of $\pA _{m/}$ which go to the same component under the Segal map, then
applying this operation for various $m$, we obtain a simplicial space whose
components are connected and homotopic to wedges of circles.  (We have to start
with an $X$ having $X_1=\ast$). The resulting simplicial space has the same
realization as $X$. If $X$ has only finitely many nondegenerate simplices then
one can stop after a finite number of applications of this operation.  Taking
the fundamental groups of the component spaces (based at the degeneracy of the
unique basepoint) gives a simplicial free group. Taking the classifying
simplicial sets of these groups in each component we obtain a bisimplicial set
whose realization is equivalent to $X$. This bisimplicial set actually satisfies
$\pA _{p,0}= \pA _{0,k}=\ast$, in other words it satisfies the globular condition in
both directions!  We can therefore view it as  a Segal precategory in two ways. The
second way, interchanging the two variables, yields a Segal precategory where the
Segal maps are {\em isomorphisms} (because at each stage it was the classifying
simplicial set for a group).  Thus, viewed in this way, it is a Segal groupoid
and Segal's theorem implies that the simplicial set $p\mapsto \pA _{p,1}$, which
is the underlying set of a simplicial free group, has the homotopy type of
$\Omega X$.



\part{The model structure}


\chapter{Sequentially free precategories}
\label{freecat1}

In this chapter, we continue the study of weakly enriched categories by looking at some basic objects:
these are the categories with an ordered set of objects $x_0,\ldots , x_n$ and
morphisms other than the  identity from $x_i$ to $x_j$ only when $i<j$. More precisely we consider the {\em free}
categories of this type obtained by specifying an object $B_i\in \mM$ of morphisms from $x_{i-1}$ to $x_i$ 
for $1\leq i\leq n$; then the object of morphisms from $x_i$ to $x_j$ should be the product of the $B_k$
for $k=i+1,\ldots , j$. One of the main tasks is to look at a notion of precategory which corresponds to this
notion of category. At the end will be our main calculation, which is what happens when one takes the product
of two such categories. 

Throughout, the notion of weak equivalence on $\precat (X,\mM )$ is the one
given by the model structure of Theorem \ref{modstrucs}.

\section{Imposing the Segal condition on $\Upsilon$}
\label{basicseqfree}

Recall $\mM$-precategories $\Upsilon (B_1,\ldots , B_k)$ defined in Sections \ref{sec-precatexamples} and \ref{sec-upsilon}. These can be strictly categorified, which is to say that we can construct corresponding strict $\mM$-categories which will be weakly
equivalent (Theorem \ref{upstildequiv} below).

Define a precategory $\Upstild _k(B_1,\ldots , B_k)$ as follows: the set of objects is
the same as for $\Upsilon _k(B_1,\ldots , B_k)$, that is  $[k]=\{ \upsilon _0,\ldots , \upsilon _k\}$. For any sequence $\upsilon _{i_0},\ldots , \upsilon _{i_n}$ with
$i_0\leq \ldots \leq i_n$, we put 
\begin{equation}
\label{upstildformRecall}
\Upstild _k(B_1,\ldots , B_k)(\upsilon _{i_0},\ldots , \upsilon _{i_n}):= 
B_{i_0+1}\times B_{i_0+2}\times \cdots \times B_{i_n-1} \times B_{i_n}.
\end{equation}
This includes the unitality condition 
$\Upstild _k(B_1,\ldots , B_k)(\upsilon _i,\ldots , \upsilon _i)=\ast$. 
For any other sequence, that is to say any sequence which is not increasing, the value is $\emptyset$. 
Note in particular that 
\begin{equation}
\label{adjident}
\Upstild _k(B_1,\ldots , B_k)(\upsilon _{i-1},
\upsilon _{i}) = B_i = \Upsilon _k(B_1,\ldots , B_k)(\upsilon _{i-1},\upsilon _{i}).
\end{equation}
By the adjunction property of $\Upsilon _k(B_1,\ldots , B_k)$ there is a unique map 
\begin{equation}
\label{upsinclusion}
\Upsilon _k(B_1,\ldots , B_k)\rightarrt \Upstild _k(B_1,\ldots , B_k)
\end{equation}
inducing the identity \eqref{adjident} on adjacent pairs of objects. In fact, we can consider $\Upsilon _k(B_1,\ldots , B_k)$
as a subobject of $\Upstild _k(B_1,\ldots , B_k)$ with this map as the inclusion. 

\begin{lemma}
The precategory $\Upstild _k(B_1,\ldots , B_k)$ is a stict $\mM$-category, in particular it is a Segal $\mM$-category.
\end{lemma}
\begin{proof}
If $\upsilon _{i_0},\ldots , \upsilon _{i_n}$ is an increasing sequence i.e. $i_0\leq \ldots \leq i_n$, and if $1\leq j \leq n-1$ then
by the formula \eqref{upstildformRecall} the natural maps obtained by splitting the sequence of objects at $\upsilon _{i_j}$ give an isomorphism
\begin{diagram}
\Upstild _k(B_1,\ldots , B_k)(\upsilon _{i_0},\ldots , \upsilon _{i_n}) \\  
\downarr_{\cong}  \\
\Upstild _k(B_1,\ldots , B_k)(\upsilon _{i_0},\ldots , \upsilon _{i_j})\times \Upstild _k(B_1,\ldots , B_k)(\upsilon _{i_j},\ldots , \upsilon _{i_n}).
\end{diagram}
By induction it follows that the Segal maps are isomorphisms. 
\end{proof}

One should think of $\Upstild _k(B_1,\ldots , B_k)$ as being the {\em free $\mM$-category generated by morphism objects $B_i$ going from $\upsilon _{i-1}$
to $\upsilon _i$}. Here the objects are linearly ordered, and the generating objects of $\mM$ are placed between adjacent objects in the ordering. 

We would like to make precise this intuition by showing that it is the Segal $\mM$-category generated by the precategory $\Upsilon _k (B_1,\ldots , B_k)$
as stated in the following theorem. 

\begin{theorem}
\label{upstildequiv}
For any $k$ and any sequence of objects $B_1,\ldots , B_k$, the inclusion \eqref{upsinclusion} is a weak equivalence in $\precat ( \{ \upsilon _0,\ldots , \upsilon _k\} , \mM )$. If each $B_i$ is cofibrant, then
it is a trivial cofibration in the injective model structure. 
\end{theorem}

\section{Sequentially free precategories in general}
\label{genseqfree} 

Before getting to the proof of the previous theorem, which will be done at the
end of the chapter on page \pageref{upstildequivproof}, 
it is useful to generalize the above situation by giving a criterion for when a precategory will lead to a free
Segal category with linearly ordered object set and generators $B_i$ between adjacent objects. 

\begin{definition} 
\label{seqfreedef}
A {\em sequentially free $\mM$-precategory} consists of a 
finite linearly ordered set $X=\{ x_0,\ldots , x_k\}$ together with a structure of 
$\mM$-precategory $A\in \precat (X,\mM )$, satisfying the following properties:
\newline
(SF1)---if $x_{i_0},\ldots , x_{i_n}$ is a sequence of objects which are not increasing, i.e. there is some $0<j\leq n$ with $i_{j-1}>i_j$,
then $A(x_{i_0},\ldots , x_{i_n})=\emptyset$; and 
\newline
(SF2)---if $x_{i_0},\ldots , x_{i_n}$ is a sequence of objects in increasing order i.e. $i_{j-1}\leq i_j$ for all $0<j\leq n$, then
the outer map for the $n$-simplex provides a weak equivalence 
\begin{equation}
\label{sf2sim}
A(x_{i_0},\ldots , x_{i_n})\rightarrt^{\sim} A(x_{i_0},x_{i_n}).
\end{equation}
\end{definition}

\begin{remark}
\label{seqfreeconstcontr}
Note that condition (SF2) for a sequence $(x_i)$ of length $n=0$ says that the map 
$\ast = A(x_i)\rightarrt A(x_i,x_i)$ is a weak equivalence. 
\end{remark}

We often say that $A$ is sequentially free {\em with respect to a given order on $X$} if $(X,A)$ is sequentially free for the
ordering in question. 

\begin{lemma}
Both $\Upsilon _k(B_1,\ldots , B_k)$ and $\Upstild _k(B_1,\ldots , B_k)$ are sequentially free $\mM$-precategories with respect to the ordering on $[k]$. 
\end{lemma}
\begin{proof}
Both clearly satisfy (SF1). For $\Upstild_k(B_1,\ldots , B_k)$ the condition (SF2)
holds by construction. For $\Upsilon _k(B_1,\ldots , B_k)$, suppose we have
an increasing
sequence of objects $\upsilon _{i_0}\leq \cdots \leq \upsilon _{i_n}$.
If $i_n=i_0$ then the values on both sides of \eqref{sf2sim}
are $\ast$. If $i_n=i_0+1$ then the values on both sides are $B_{i_n}$,
whereas if $i_n>i_0+1$ the values on both sides are $\emptyset$.
In all three cases the map \eqref{sf2sim} is an isomorphism, {\em a fortiori} a weak equivalence. 
\end{proof}

\begin{lemma}
Suppose $X$ is a fixed linearly ordered finite set, $\beta$ is an ordinal, and
$\{ A(b )\} _{b\in \beta}$ is a transfinite sequence of $\mM$-precategories,
that is to say a functor $\beta \rightarrt \precat (X,\mM )$. 
Suppose that each $(X,A(b))$ is a sequentially free  $\mM$-precategory with respect
to the given order on $X$, and that for any $b\leq b'$  the transition map $A(b)\rightarrow A(b')$  
is an injective cofibration in $\precat (X,\mM )$.
Then the colimit $(X,\colim _{b\in \beta }A(b))$ is again a sequentially free  $\mM$-precategory
with respect to the same ordering. 
\end{lemma}
\begin{proof}
The colimit is calculated levelwise as a $\Delta _X^o$-diagram in $\mM$, since
filtered colimits preserve the unitality condition. A colimit of objects $\emptyset$
is again $\emptyset$, so the colimit satisfies (SF1). The left properness 
hypotheses on $\mM$ implies
transfinite left properness (Proposition \ref{prop-transfiniteleftproper}),
so the weak equivalence \eqref{sf2sim} is preserved in the colimit whose transition
maps are cofibrations, giving (SF2). 
\end{proof}

We now come to one of the main steps where we gain some control over the process of generators and relations. Recall from Section \ref{sec-gen}
the operation $\Gen$ consisting of applying one step of the calculus of generators and relations. 

\begin{lemma}
\label{sfgenstep}
Suppose $(X,A)$ is a sequentially free  $\mM$-precategory with linearly ordered object set $X=\{ x_0,\ldots , x_m\}$, and $x_{a},x_{a+1}, \ldots , x_{b}$ is a strictly increasing
sequence of adjacent objects with $0\leq a < b\leq m$. 
Then  the new $\mM$-precategory $\Gen ((X,A); x_{a},\ldots , x_{b})$ is also sequentially free
with the same ordered set of objects $X$, and furthermore for any $a<j\leq b$ the map 
$$
A(x_{j-1},x_j)\rightarrt \Gen ((X,A); x_{a},\ldots , x_{b})(x_{j-1},x_j)
$$
is an isomorphism (hence a weak equivalence) in $\mM$. 
\end{lemma}
\begin{proof}
Fix a factorization $E,e,p_1,\ldots , p_{b-a}$ as used for the construction of $\Gen ((X,A); x_{a},\ldots , x_{b})$.
Note that 
$$
p:E\rightarrt^{\sim} A(x_a,x_{a+1})\times \cdots \times A(x_{b-1},x_b)
$$
is a weak equivalence. 
The sequence of objects $x_a,\ldots , x_b$ is disjoint, so we can use the description of 
$\Gen ((X,A); x_{a},\ldots , x_{b})$ given in Lemma \ref{gendescription-lem}. 
That says that for any sequence
of the form $x_{i_0},\ldots , x_{i_p}$ if $a\leq i_0\leq \cdots \leq i_p \leq b$
with $i_0 + 2\leq i_p$
then
\begin{equation}
\label{casewithin}
\Gen (A; x_a,\ldots , x_b) (x_{i_0},\ldots , x_{i_p}) = A(x_{i_0},\ldots , x_{i_p})\cup ^{A(x_a,\ldots , x_b)}E.
\end{equation}
but for any other sequence, 
\begin{equation}
\label{caseother}
\Gen (A; x_a,\ldots , x_b) (x_{i_0},\ldots , x_{i_p}) = A(x_{i_0},\ldots , x_{i_p}).
\end{equation}
We can now check the conditions (SF1) and (SF2). If the sequence of objects $(x_{i_0},\ldots , x_{i_p})$ is
not increasing, then it falls into the second case \eqref{caseother}, and $A(x_{i_0},\ldots , x_{i_p})=\emptyset$ by (SF1) for $A$,
which gives (SF1) for $\Gen (A; x_a,\ldots , x_b)$. 
Suppose that $(x_{i_0},\ldots , x_{i_p})$ is increasing. We need to check (SF2). If either ${i_0}< a$ or
${i_p}>b$ then we again fall into case \eqref{caseother} for both the full sequence $(x_{i_0},\ldots , x_{i_p})$ and also
the pair of endpoints $(x_{i_0}, x_{i_p})$. Thus, by (SF2) for $A$ we have a weak equivalence 
\begin{diagram}
\Gen (A; x_a,\ldots , x_b) (x_{i_0},\ldots , x_{i_p}) = A(x_{i_0},\ldots , x_{i_p})
\\
\downarr_{\sim} 
\\
A(x_{i_0}, x_{i_p})= 
\Gen (A; x_a,\ldots , x_b) (x_{i_0},x_{i_p})
\end{diagram}
giving this case of (SF2) for $\Gen (A; x_a,\ldots , x_b)$. 
Suppose on the other hand that $a\leq {i_0}\leq {i_p}\leq b$. If $i_p\leq i_0+1$ then we again are in case \eqref{caseother}
for both the full sequence and the pair of endpoints, so we get the condition (SF2) as above. Suppose therefore that
$i_0 + 2\leq i_p$. In the diagram
$$
\begin{diagram}
A(x_{i_0},\ldots , x_{i_p})& \leftarr & A(x_a,\ldots , x_b)& \rightarr & E \\
\downarr && \dEqualarr && \dEqualarr \\
A(x_{i_0}, x_{i_p})& \leftarr & A(x_a,\ldots , x_b)& \rightarr & E
\end{diagram}
$$
the left vertical arrow is an equivalence by (SF2) for $A$. By left properness of $\mM$ 
via Corollary \ref{pushoutequivcor}, the vertical maps therefore induce
a weak equivalence from the pushout of the top row to the pushout of the bottom row. By equations \eqref{casewithin}
which apply both to the full sequence $(x_{i_0},\ldots , x_{i_p})$ and the pair of endpoints $(x_{i_0},x_{i_p})$,
these pushouts are respectively $\Gen (A; x_a,\ldots , x_b) (x_{i_0},\ldots , x_{i_p})$ and
$\Gen (A; x_a,\ldots , x_b) (x_{i_0},x_{i_p})$. Thus, the map from the one to the other is an equivalence, which gives condition (SF2) in
this last case.  This proves that $\Gen ((X,A); x_{a},\ldots , x_{b})$ is again a sequentially free $\mM$-precategory.

For the last statement, note that the sequence $(x_{j-1},x_j)$ falls into case \eqref{caseother} because the space between the endpoints is
only $1$. Hence the formula \eqref{caseother} says that the map 
$$
A(x_{j-1},x_j)\rightarr  \Gen ((X,A); x_{a},\ldots , x_{b})(x_{j-1},x_j)
$$
is an isomorphism.
\end{proof}

Recall that the map $A\rightarr \Gen ((X,A); x_{a},\ldots , x_{b})$ is a weak equivalence in  $\precat (X; \mM )$. 

\begin{lemma}
\label{sfgeniterate}
Suppose $(X,A)$ is a sequentially free  $\mM$-precategory with ordered object set $X=\{ x_0,\ldots , x_m\}$.
By iterating a series of operations of the form
$A\mapsto \Gen ((X,A); x_{a},\ldots , x_{b})$ we can obtain a weak equivalent sequentially free $\mM$-precategory
$(X,A)\rightarrt (X,A')$ such that $A'$ satisfies the Segal conditions.
\end{lemma}
\begin{proof}
We show how to obtain the Segal condition for strictly increasing sequences of adjacent objects. At the end of the proof we go from
here to the Segal condition for general sequences. 

Fix an integer $n_0$ and
suppose $(X,A)$ satisfies the Segal condition for all adjacent sequences $(x_c,x_{c+1},\ldots , x_d)$ with $d-c>n_0$, 
and a certain number of adjacent sequences $(x_{c_v},x_{c_v+1},\ldots , x_{d_v})$ for $(c_v,d_v)$ indexed by $v\in V$ for some set $v$,
with $d_v-c_v = n_0$. 
Suppose $0\leq a < b \leq m$ with $b-a= n_0$. Then $\Gen ((X,A); x_{a},\ldots , x_{b})$ satisfies the 
Segal condition for all adjacent sequences $(x_c,x_{c+1},\ldots , x_d)$ with $d-c>n_0$, 
for the given adjacent sequences $(x_{c_v},x_{c_v+1},\ldots , x_{d_v})$ with $v\in V$, and also for the sequence
$(x_a,\ldots , x_b)$. Indeed, if $(x_c,x_{c+1},\ldots , x_d)$ with $d-c>n_0$
then the terms entering into the Segal map for this sequence are all covered by the situation \eqref{caseother}
in the explicit description of $\Gen ((X,A); x_{a},\ldots , x_{b})$ used in the previous proof (see Lemma \ref{gendescription-lem}). 
Note that the terms in the product on the right hand side of the Segal map are of the form
$\Gen ((X,A); x_{a},\ldots , x_{b})(x_{j-1},x_j)= A(x_{j-1},x_j)$. Here we use the condition that we are only looking at adjacent sequences.
By the recurrence hypothesis on $A$, the Segal map for this sequence is an equivalence. 

Similarly, if $(x_{c_v},x_{c_v+1},\ldots , x_{d_v})$  is one of our given sequences with $d_v-c_v = n_0$, and if it is different
from the sequence $(x_a,\ldots , x_b)$, then either $c_v<a$ or $d_v>b$ so again everything entering into the Segal map 
for this sequence is the same as for $A$, by \eqref{caseother}. Thus, again the inductive hypothesis on $A$ implies the
Segal condition for $\Gen ((X,A); x_{a},\ldots , x_{b})(x_{c_v},x_{c_v+1},\ldots , x_{d_v})$. 

At the sequence $(x_a,\ldots , x_b)$ equation \eqref{casewithin} gives
$$
\Gen (A; x_a,\ldots , x_b) (x_{a},\ldots , x_{b}) = A(x_{a},\ldots , x_{b})\cup ^{A(x_a,\ldots , x_b)}E = E
$$
and the map 
$$
E\rightarrt A(x_{a},x_{a+1})\times \cdots \times A(x_{b-1},x_b)
$$
is an equivalence by hypothesis on the choice of $E$. Therefore, the Segal condition holds at the sequence 
$(x_a,\ldots , x_b)$ too and we can add this to our collection $V$ of good sequences of lenght $n_0$.

So, the inductive procedure is to start with the maximal sequence $x_0,\ldots , x_m$ of length $m$, impose
the Segal condition here by changing $A$ to $\Gen (A; x_0,\ldots , x_m)$; then successively impose the Segal condition
on all sequences of length $m-1$, then $m-2$ and so on, using the above inductive observation. At each step $A$
is replaced by $\Gen (A; x_a,\ldots , x_b)$ and there is a map from the old to the new $A$ which is a weak equivalence
in $\precat (X,\mM )$. By Lemma \ref{sfgenstep}, the new $A$ is always again a sequentially free $\mM$-precategory.
Combining these steps (there are a finite number) down to $n_0=2$ we arrive
at a map $(X,A)\rightarrt (X,A')$ which is a weak equivalence in $\precat (X,\mM )$, 
and such that $A'$
is a sequentially free $\mM$-precategory which
satisfies the Segal condition for all sequences of the form $(x_c,\ldots , x_d)$. 

We claim that this implies that $A'$ satisfies the Segal condition for any  sequence 
of the form $x_{i_0},\ldots ,x_{i_p}$. Note that if $p=0$ then 
the Segal condition is automatic since $A'(X)=\ast$. Suppose $p\geq 1$ and 
consider our sequence $(x_{i_0},\ldots ,x_{i_{p-1}}, x_{i_p})$.
If any $i_{j-1}>i_j$ then $A'(x_{i_0},\ldots ,x_{i_{p-1}}, x_{i_p}) = \emptyset$ and
$A'(x_{i_{j-1}}, x_{i_j})=\emptyset$. The second statement, plus Lemma \ref{emptyempty} on the direct product with $\emptyset$,
imply that the right hand side of the Segal map is $\emptyset$, which is the same as the left side by the first statement. 
So in this case, the Segal condition is automatic. Hence we may assume that $i_{j-1}\leq i_j$ for all $1\leq j\leq p$.

Let $a:= i_0$ and $b:= i_p$, and denote by $(x_a,\ldots , x_b)$ the full sequence of adjacent elements (each counted once)
going from $x_a$ to $x_b$. There is a unique map $\sigma : (x_{i_0},\ldots ,x_{i_{p-1}}, x_{i_p})\rightarrt (x_a,\ldots , x_b)$
in $\Delta _X$, sending each $x_{i_j}$ to the same object at the unique place it occurs in $(x_a,\ldots , x_b)$.
Note that we are using here the reduction of the previous paragraph that $i_{j-1}\leq i_j$ for all $1\leq j\leq p$,
guaranteeing that each $x_{i_j}$ occurs in the list $(x_a,\ldots , x_b)$. 
This $\sigma$ induces a map 
\begin{equation}
\label{topmap}
\sigma ^{\ast}: A'(x_a,\ldots , x_b)\rightarrt^{\sim} A'(x_{i_0},\ldots ,x_{i_{p-1}}, x_{i_p}).
\end{equation}
That this is a weak equivalence, can be seen by considering the diagram with maps to the spanning pair
$(x_a,x_b)=(x_{i_0},x_{i_p})$: 
$$
\begin{diagram}
A'(x_a,\ldots , x_b) & \rightarr & A'(x_{i_0},\ldots ,x_{i_{p-1}}, x_{i_p})\\
\downarr && \downarr \\
A'(x_{a},x_{b}) & \rEqualarr & 
A'(x_{i_0},x_{i_p}) .
\end{diagram}
$$
The vertical maps are weak equivalences by the sequentially free condition, so the top map is a weak equivalence by 3 for 2. 

On the other hand, for each pair in the original sequence, 
consider the Segal map for it:
$$
\sigma _{i_0,i_1}: A'(x_{i_0},\ldots , x_{i_1})\rightarrt  A'(x_{i_0},x_{i_0+1})\times \cdots \times A'(x_{i_1-1},x_{i_1}) ,
$$
$$
\vdots 
$$
$$
\sigma _{i_{p-1},i_p}: A'(x_{i_{p-1}},\ldots , x_{i_p})\rightarrt  A'(x_{i_{p-1}},x_{i_{p-1}+1})\times \cdots \times A'(x_{i_p-1},x_{i_p}) .
$$
The mini-sequences appearing on the left hand sides, are the sequences of length $i_j-i_{j-1}$ going from $i_{j-1}$ to $i_j$ by intervals of step $1$.
Whenever $i_{j-1}= i_j$, we have a sequence of length zero and both sides of the map are equal to $\ast$.  
The case of the Segal condition which we already know, says that the $\sigma _{i_{j-1},i_j}$ are weak equivalences. 

Putting these all together, we get a weak equivalence whose target is the full Segal product for the sequence $(x_a,\ldots , x_b)$ considered above: 
$$
A'(x_{i_0},\ldots , x_{i_1})\times \cdots \times A'(x_{i_{p-1}},\ldots , x_{i_p}) \rightarrt^{\sim} A'(x_a,x_{a+1})\times \cdots \times A'(x_{b-1},x_b).
$$
On the other hand, the map to the spanning interval gives a map for each mini-sequence
\begin{equation}
\label{minispan}
A'(x_{i_{j-1}},\ldots , x_{i_j})\rightarrt^{\sim}A'(x_{i_{j-1}},x_{i_j})
\end{equation}
which is a weak equivalence, by the condition that $A'$ is sequentially free. The mini-sequences map into the full adjacent sequence $(x_a,\ldots , x_b)$,
giving maps 
$$
A'(x_a,\ldots , x_b)\rightarrt A'(x_{i_{j-1}},\ldots , x_{i_j}).
$$
The Segal map for $(x_a,\ldots , x_b)$ thus factors as
\begin{equation}
\begin{diagram}
\label{leftvert}
A'(x_a,\ldots , x_b)&\rightarr^{\sim} &A'(x_{i_0},\ldots , x_{i_1})\times \cdots \times A'(x_{i_{p-1}},\ldots , x_{i_p}) \\
& & \downarr^{\sim}\\
& &  A'(x_a,x_{a+1})\times \cdots \times A'(x_{b-1},x_b).
\end{diagram}
\end{equation}
The second arrow is a weak equivalence as pointed out previously, and the composition is a weak equivalence by the Segal condition for 
$(x_a,\ldots , x_b)$ which we already know, so the first map is a weak equivalence by 3 for 2. 

Now we have a commutative diagram 
$$
\begin{diagram}
A'(x_a,\ldots , x_b) & \rightarr & A'(x_{i_0},\ldots ,x_{i_{p-1}}, x_{i_p})\\
\downarr && \downarr \\
A'(x_{i_0},\ldots , x_{i_1})\times \cdots \times A'(x_{i_{p-1}},\ldots , x_{i_p}) & 
\rightarr & 
A'(x_{i_0},x_{i_1})\times \cdots \times A'(x_{i_{p-1}},x_{i_p})
\end{diagram}
$$
where the top arrow is a weak equivalence as seen above \eqref{topmap}, the left vertical map is a weak equivalences by \eqref{leftvert},
and the bottom map is a weak equivalence by combining together the equivalences \eqref{minispan}. By 3 for 2, it follows that the right vertical
map is a weak equivalence, which is the Segal condition for the sequence $(x_{i_0},\ldots ,x_{i_{p-1}}, x_{i_p})$. 
This completes the proof of the lemma. 
\end{proof}

\begin{corollary}
\label{seqfreeinvariance}
Suppose $(X,A)$ is a sequentially free ordered $\mM$-precategory, and let $r:(X,A)\rightarrt (X,A')$ be a fibrant replacement in
either the projective or  injective model structure on $\precat (X;\mM )$ constructed in 
Theorem \ref{modstrucs}. Then $(X,A')$ is a
sequentially free ordered $\mM$-precategory for the same order on $X$, and 
for any $0<j\leq n$ the map $A(x_{j-1},x_j)\rightarrt A'(x_{j-1},x_j)$
is a weak equivalence. 

The same conclusions holds if, instead of a fibrant replacement,
$r$ is a weak equivalence to an object $A'$ which satisfies the Segal conditions. 
\end{corollary}
\begin{proof}
Use either the projective or the injective structure in what follows. 
Note that the sequentially free condition is preserved by 
levelwise weak equivalences of
unital diagrams on $\Delta ^o_X$. 

By Lemma \ref{sfgeniterate}, there is a map $A\rightarrow A''$, weak equivalence
in $\precat (X,\mM )$, such that $A''$ satisfies the Segal conditions and is still
sequentially free. Consider furthermore a map $s:A''\rightarrow A^3$ which is
a trivial cofibration to a fibrant object in $\precat (X,\mM )$. Then $A^3$ also
satisfies the Segal conditions, so Lemma \ref{localwesegal} says that 
$s$ is a levelwise weak equivalence. It follows that $A^3$ is again sequentially free.

Our different fibrant replacement $r:A\rightarrow A'$ is a trivial cofibration,
so there is a map $g:A'\rightarrow A^3$ compatible with the maps from $A$.
By 3 for 2 this map $g$ is a weak equivalence in $\precat (X,\mM )$, so again by 
Lemma \ref{localwesegal} it is a levelwise weak equivalence, hence $A'$ is sequentially
free. Furthermore, $A(x_{j-1},x_j)\rightarrt A'' (x_{j-1},x_j)$ is a weak
equivalence in $\mM$; it follows from the above levelwise weak equivalences compatible
with the maps from $A$, that the same is true for $A^3$, then $A'$. 

For the last paragraph, suppose $r:A\rightarrt A'$ is a weak equivalence
in $\precat (X,\mM )$ and $A'$ satisfies the Segal conditions. Then 
choosing first a fibrant replacement $A^3$ of $A'$, then a factorization,
we can get a square 
$$
\begin{diagram}
A & \rightarr & A' \\
\downarr & & \downarr \\
A'' & \rightarr & A^3
\end{diagram}
$$
such that the vertical arrows are fibrant replacements and all arrows are weak
equivalences. By Lemma \ref{localwesegal}, the right vertical and bottom maps are
levelwise weak equivalences since $A'$ and $A''$ satisfy the Segal conditions.
By the first part of the present lemma, $A''$ is sequentially free, 
it follows that $A^3$ and then $A'$ are sequentially free. On adjacent pairs of
objects the left vertical map induces an equivalence, so by 3 for 2 the top map 
does too. This shows that the required statements hold for $A'$. 
\end{proof}

In the situation of the previous corollary,
suppose $(x_{i_0},\ldots , x_{i_n})$ is an increasing  sequence of objects. Then we have weak equivalences 
$$
A'(x_{i_0}, x_{i_n})
\rightarrt^{\sim}
A'(x_{i_0},\ldots , x_{i_n}) \rightarrt^{\sim}
A'(x_{i_0}, x_{i_1})\times \cdots \times A'(x_{i_{n-1}}, x_{i_n}).
$$
On the other hand if the sequence is not increasing, the corresponding morphism space is $\emptyset$.

Suppose that our sequence of objects is a sequence of adjacent objects, that is to say look at a sequence
of the form $(x_a, x_{a+1},\ldots , x_b)$ for $a\leq b$. The condition of the corollary says furthermore that we have weak equivalences
$A(x_{j-1},x_j)\rightarrt^{\sim} A'(x_{j-1},x_j)$ for $a<j\leq b$. These equivalences go together to give
an equivalence on the level of the direct product, by the product conditions (PROD) and 
(DCL) for $\mM$ as in Lemma \ref{prodconsequence}. Thus, the pair of weak equivalences in
the previous paragraph extends to a chain of equivalences
\begin{equation}
\label{chainofeq}
\begin{diagram}
A'(x_{a},\ldots , x_{b}) & \rightarr^{\sim}& 
A'(x_{a}, x_{a+1})\times \cdots \times A'(x_{b-1}, x_{b})\\
\downarr^{ \sim} & & \uparr _{ \sim} \\
A'(x_{a}, x_{b}) && A(x_{a}, x_{a+1})\times \cdots \times A(x_{b-1}, x_{b}).
\end{diagram}
\end{equation}
All in all, this says that $(X,A')$ looks up to  homotopy very much like the $\Upstild _k(B_1,\ldots , B_k)$ defined at the
start, with $B_i = A(x_{i-1},x_i)$. 

\begin{corollary}
\label{seqfreemorphism}
Suppose $f:(X,A)\rightarrt (X,B)$ is a morphism in $\precat (X; \mM )$ between $\mM$-precategories which are both sequentially free for
the same ordering of the underlying set of objects $X$. Suppose that for any adjacent objects $x_{j-1},x_j$ in the ordering,
$f$ induces a weak equivalence $A(x_{j-1},x_j)\rightarrt^{\sim} B(x_{j-1},x_j)$. Then $f$ is a weak equivalence in $\precat (X; \mM )$. 
\end{corollary}
\begin{proof}
View the object set as numbered $X=\{ x_0,\ldots , x_k\}$. Let $(X,A')$ and $(X,B')$ be fibrant replacements for $(X,A)$ and
$(X,B)$ respectively, and we may assume that $f$ extends to a map $f':(X,A')\rightarrt (X,B')$. 
For any $0\leq a \leq b \leq k$, $f'$ and $f$ induce a morphism between the chains of equivalences considered in \eqref{chainofeq},
for $A,A'$ and $B,B'$. The hypothesis of the present corollary says that the induced map on the bottom right corner is a weak
equivalence; thus the induced map $f':A'(x_a,x_b)\rightarrt B'(x_a,x_b)$ is a weak equivalence. Note on the other hand that
if $a>b$ then $A'(x_a,x_b)= \emptyset$ and $B'(x_a,x_b)=\emptyset$. Thus $f'$ induces a weak equivalence on the morphism space for
any pair of objects. By the Segal condition, it follows that $f'$ is an objectwise weak equivalence in the category of $\Delta _X^o$-diagrams in $\mM$,
so it is a weak equivalence in the model structures. Since $f'$ is a fibrant replacement for $f$, this implies that $f$ was a weak equivalence. 
\end{proof}

{\em Proof of Theorem \ref{upstildequiv}:} 
\label{upstildequivproof}
the map 
$$
f:\Upsilon _k(B_1,\ldots , B_k)\rightarrt \Upstild _k(B_1,\ldots , B_k)
$$
satisfies the hyptheses of the previous corollary, so we conclude that it is a weak equivalence. If the $B_i$ are cofibrant, then the map in question
is a cofibration in the injective
model structure, because is is an objectwise cofibration when viewed in the injective model structure: the induced maps are
either identities, or maps from $\emptyset$ into a product of copies of $B_i$. 
A product of copies of $B_i$ is cofibrant, because of the
cartesian
conditions on $\mM$ (see Lemma \ref{prodproperties}).
\eop


\chapter{Products}
\label{product1}

In this chapter, we consider the direct product of two $\mM$-enriched precategories.
On the one hand, we would like to maintain the cartesian or product condition (PROD) for the new model category we are constructing.
On the other hand, compatibility with direct product provides the main technical tool we need in order to study
pushouts by weak equivalences in $\precat (\mM )$. Recall that we already have
the injective and projective model structures of Theorem \ref{modstrucs} on $\precat (X,\mM )$
at  our disposal. 

Here is the place where we really use the full structure of the category $\Delta$, as well as the unitality condition.
At the end of the chapter, we'll discuss some counterexamples showing why these aspects are necessary. 

\section{Products of sequentially free precategories}
Suppose $X = \{ x_0,\ldots , x_m\} $ and $Y= \{ y_0,\ldots , y_n\} $ are finite linearly ordered sets, and
$(X,\pA )$ and $(Y,\pB )$ are sequentially free $\mM$-precategories, with respect to the orderings. The product of these objects
considered in the category $\precat (\mM )$, has the form $(X\times Y,\pA \boxtimes \pB )$ where
$$
(\pA \boxtimes \pB )((x_{i_0}, y_{j_0}),\ldots , (x_{i_p}, y_{j_p}) := \pA (x_{i_0},\ldots , x_{i_p})\times \pB (y_{j_0},\ldots , y_{j_p}).
$$

Let $\pA \rightarrt \pA '$ and
$\pB \rightarrt \pB '$ be weak equivalences towards objects satisfying the Segal condition and which are therefore also sequentially free (Corollary \ref{seqfreeinvariance}). 
We obtain a map of $\mM$-precategories
$\pA \boxtimes \pB \rightarrt \pA '\boxtimes \pB '$ on the object set $X\times Y$. We would like to show that this is a weak equivalence in $\precat (X\times Y; \mM )$. 

For any subset $S\subset X\times Y$, let $\bfi (S):S\hookrightarrow X\times Y$ denote the inclusion map, and
$$
\bfi (S)^{\ast}: \precat (X\times Y, \mM )\rightarrt \precat (S,\mM )
$$
the pullback map on $\mM$-precategory structures. 
This can be pushed back to a $\mM$-precategory structure with object set $X\times Y$, and the resulting functor
will be denoted as $\bfi (S) ^{\ast}_!$ for brevity: 
$$
\bfi (S) ^{\ast}_!:= \bfi (S)_! \bfi (S)^{\ast}: \precat (X\times Y, \mM )\rightarrt \precat (X\times Y, \mM ).
$$

This may be described explicitly as follows. Suppose $\pC \in \precat (X\times Y, \mM )$. 
Recall from the first paragraph of Section \ref{sec-varying1}, that 
$\bfi (S)_! \bfi (S)^{\ast}(\pC )$  is the same as the precategory $\bfi (S)^{\ast}(\pC )$ over set of objects $S$,
extended by adding on the discrete $\mM$-enriched category on the complementary set $X\times Y -S$. Thus,
for
any sequence of objects $((x_{i_0}, y_{i_0}),\ldots , (x_{i_p}, y_{i_p}))$ we have:
$$
\bfi (S) ^{\ast}_!(\pC )((x_{i_0}, y_{i_0}),\ldots , (x_{i_p}, y_{i_p})) = \pC  ((x_{i_0}, y_{i_0}),\ldots , (x_{i_p}, y_{i_p}))
$$
if all the $(x_{i_j},y_{i_j})$ are in $S$; if any of the pairs is not in $S$ and the sequence is not constant then
$$
\bfi (S) ^{\ast}_!(\pC )((x_{i_0}, y_{i_0}),\ldots , (x_{i_p}, y_{i_p})) =\emptyset ,
$$
while the unitality condition requires that 
$$
\bfi (S) ^{\ast}_!(\pC )((x_{i_0}, y_{i_0}),\ldots , (x_{i_0}, y_{i_0})) =\ast ,
$$
for a constant sequence at any point $(x_{i_0}, y_{i_0})\in  X\times Y$. 

\begin{lemma}
If $\pC $ satisfies the Segal condition, then so does $\bfi (S) ^{\ast}(\pC )$ (on object set $S$) and
$\bfi (S) ^{\ast}_!(\pC )$ (on object set $X\times Y$).
\end{lemma}
\begin{proof}
This follows immediately from Lemma \ref{pfpb}.
\end{proof}

The above discussion is usually applied to an object of the form $\pC =\pA \boxtimes \pB \in \precat (X\times Y; \mM )$. 

Here are the special types of subsets $S$ we will consider. 
A {\em box} will be a subset of the form 
$$
{\bf B}_{a,b}:= \{x_0,\ldots , x_a\} \times \{ y_0,\ldots , y_b\} \subset X\times Y,
$$
which can be pictured for example with $(a,b)=(5,3)$ as
$$
{\setlength{\unitlength}{.5mm}
\begin{picture}(100,100)
\put(0,10){\line(0,1){60}}
\put(0,10){\line(1,0){100}}
\put(100,10){\line(0,1){60}}
\put(0,70){\line(1,0){100}}
\multiput(0,10)(20,0){6}{\circle*{2}}
\multiput(0,30)(20,0){6}{\circle*{2}}
\multiput(0,50)(20,0){6}{\circle*{2}}
\multiput(0,70)(20,0){6}{\circle*{2}}
\put(-15,2){\ensuremath{(x_0,y_0)}}
\put(-15,75){\ensuremath{(x_0,y_b)}}
\put(95,2){\ensuremath{(x_a,y_0)}}
\put(95,75){\ensuremath{(x_a,y_b)}}
\end{picture}
} 
$$
A {\em notched box} is a subset 
$$
{\bf B}_{a,b}\supset {\bf B}_{a,b}^{\nu} := {\bf B}_{a,b-1}\cup \{ (x_{\nu}, y_b),\ldots , (x_a,y_b)\} ,
$$
defined when $b\geq 1$. 
Thus, a  notched box ${\bf B}_{a,b}^{\nu}$ is a box ${\bf B}_{a,b}$ minus a part of the top row $\{ (x_{0}, y_b),\ldots , (x_{\nu -1},y_b)\}$. 
For example a notched box with $(a,b)=(5,3)$ and $\nu = 3$ looks like
$$
{\setlength{\unitlength}{.5mm}
\begin{picture}(100,100)
\put(0,10){\line(0,1){40}}
\put(0,10){\line(1,0){100}}
\put(100,10){\line(0,1){60}}
\put(0,50){\line(1,0){60}}
\put(60,50){\line(0,1){20}}
\put(60,70){\line(1,0){40}}

\multiput(0,10)(20,0){6}{\circle*{2}}
\multiput(0,30)(20,0){6}{\circle*{2}}
\multiput(0,50)(20,0){6}{\circle*{2}}
\multiput(60,70)(20,0){3}{\circle*{2}}
\put(-15,2){\ensuremath{(x_0,y_0)}}
\put(-15,55){\ensuremath{(x_0,y_{b-1})}}
\put(95,2){\ensuremath{(x_a,y_0)}}
\put(45,75){\ensuremath{(x_{\nu},y_b)}}
\put(45,42){\ensuremath{(x_{\nu},y_{b-1})}}
\put(95,75){\ensuremath{(x_a,y_b)}}
\end{picture}
} 
$$
The boxes and notched boxes are defined for $0\leq a \leq m$ and $0\leq b\leq n$, and $0\leq \nu \leq a$. 

If $a=m,b=n$ then ${\bf B}_{a,b}=X\times Y$, and if $\nu = 0$ then ${\bf B}_{a,b}^{\nu} = {\bf B}_{a,b}$. For any $\nu$, the upper right corner 
of ${\bf B}_{a,b}^{\nu}$ is equal to the point $(x_a,y_b)$. The lower left corner is $(x_0,y_0)$.
The {\em corner of the notch} of ${\bf B}_{a,b}^{\nu}$ is
the point $(x_{\nu}, y_b)$. Note that 
$$
{\bf B}_{a,b}^{\nu} - \{ (x_{\nu}, y_b)\} = {\bf B}_{a,b}^{\nu +1}
$$
whenever $0\leq \nu < a$. For $\nu = a$ we have ${\bf B}_{a,b}^{\nu} - \{ (x_{\nu}, y_b)\} = {\bf B}_{a,b-1}$. 

A {\em tail} is a subset of the form 
$$
T^{i_0,\ldots , i_p}_{j_0,\ldots , j_p} := \{ (x_{i_0},y_{j_0}), \ldots , (x_{i_p},y_{j_p}) \}
$$
which will be considered whenever 
$$
i_0\leq i_1\leq \cdots \leq i_p, \;\;\; i_{u}\leq i_{u-1}+1,
$$
$$
j_0\leq j_1\leq \cdots \leq j_p, \;\;\; j_{u}\leq j_{u-1}+1,
$$
and $(i_u,j_u)\neq (i_{u-1}, j_{u-1})$.
Put another way, at each place in the pair of sequences, there are three possibilities for $(i_u,j_u)$ in terms of $(i_{u-1}, j_{u-1})$:
$$
(i_u,j_u)= (i_{u-1} + 1, j_{u-1}) \mbox{ or } (i_{u-1} , j_{u-1}+1) \mbox{ or } (i_{u-1} + 1, j_{u-1}+1). 
$$
Thus, the sequence  moves in single horizontal, vertical or diagonal steps. 
We say that the tail {\em goes from $(x_{i_0},y_{j_0})$ to $(x_{i_p},y_{j_p})$}. A {\em full tail} is one which goes
from $(x_0,y_0)$ to $(x_m,y_n)$. 

A {\em dipper} is a subset $S$ which is the union of a notched box ${\bf B}_{a,b}^{\nu}$, with a tail 
$T^{i_0,\ldots , i_p}_{j_0,\ldots , j_p}$ going from $(x_{i_0},y_{j_0})= (x_a,y_b)$ to $(x_{i_p},y_{j_p})= (x_m,y_n)$.
These overlap at the point $(x_a,y_b)$. For example, a dipper starting
with the previous notched box, and going out to $(m,n)=(9,7)$ could look like this:
$$
{\setlength{\unitlength}{.5mm}
\begin{picture}(100,85)
\put(0,5){\line(0,1){20}}
\put(0,5){\line(1,0){50}}
\put(50,5){\line(0,1){30}}
\put(0,25){\line(1,0){30}}
\put(30,25){\line(0,1){10}}
\put(30,35){\line(1,0){20}}

\multiput(0,5)(10,0){6}{\circle*{1.5}}
\multiput(0,15)(10,0){6}{\circle*{1.5}}
\multiput(0,25)(10,0){6}{\circle*{1.5}}
\multiput(30,35)(10,0){3}{\circle*{1.5}}
\put(32,58){\ensuremath{(x_a,y_b)}}
\qbezier(43,55)(42,45)(48,37)
\put(48,37){\vector(2,-3){0}}

\put(10,46){\ensuremath{(x_{\nu},y_b)}}
\qbezier(22,43)(22,41)(28,37)
\put(28,37){\vector(4,-3){0}}

\put(50,35){\line(1,0){10}}
\put(60,35){\circle*{1.5}}

\put(60,35){\line(0,1){10}}
\put(60,45){\circle*{1.5}}

\put(60,45){\line(1,1){10}}
\put(70,55){\circle*{1.5}}

\put(70,55){\line(1,0){10}}
\put(80,55){\circle*{1.5}}

\put(80,55){\line(0,1){10}}
\put(80,65){\circle*{1.5}}

\put(80,65){\line(1,1){10}}
\put(90,75){\circle*{1.5}}
\put(85,78){\ensuremath{(x_m,y_n)}}
\end{picture}
} 
$$

For $b\geq 1$, the union of ${\bf B}_{a,b}^{b}$ with a tail $T$ going from 
$(x_a,y_b)$ to $(x_m,y_n)$, is equal to the union of ${\bf B}_{a,b-1}^0$ with a tail $T'= \{ (x_a,x_{b-1})\}$
going from $(x_a,x_{b-1})$ to $(x_m,y_n)$, with the union overlapping at the point $(x_a,x_{b-1})$. 
For example, in the previous picture if we set $\nu = b$ it becomes
$$
{\setlength{\unitlength}{.5mm}
\begin{picture}(100,85)
\put(0,5){\line(0,1){20}}
\put(0,5){\line(1,0){50}}
\put(50,5){\line(0,1){30}}
\put(0,25){\line(1,0){50}}

\multiput(0,5)(10,0){6}{\circle*{1.5}}
\multiput(0,15)(10,0){6}{\circle*{1.5}}
\multiput(0,25)(10,0){6}{\circle*{1.5}}
\put(50,35){\circle*{1.5}}
\put(32,58){\ensuremath{(x_a,y_b)}}
\qbezier(43,55)(42,45)(48,37)
\put(48,37){\vector(2,-3){0}}

\put(70,23){\ensuremath{(x_a,y_{b-1})}}
\put(67,25){\vector(-1,0){14}}

\put(50,35){\line(1,0){10}}
\put(60,35){\circle*{1.5}}

\put(60,35){\line(0,1){10}}
\put(60,45){\circle*{1.5}}

\put(60,45){\line(1,1){10}}
\put(70,55){\circle*{1.5}}

\put(70,55){\line(1,0){10}}
\put(80,55){\circle*{1.5}}

\put(80,55){\line(0,1){10}}
\put(80,65){\circle*{1.5}}

\put(80,65){\line(1,1){10}}
\put(90,75){\circle*{1.5}}
\put(85,78){\ensuremath{(x_m,y_n)}}
\end{picture}
} 
$$

In this situation if furthermore $b-1=0$, then the subset $S$ becomes a full tail going from $(x_0,y_0)$ to $(x_m,y_n)$.
Note also that a box with $a=0$ plus a tail, is again equal to a full tail. 

In view of this, we consider the variables $a,b,\nu$ for the notched box in a dipper $S$, only in the
range $0\leq \nu < a$ and $0<b$. The corner of the notch $(x_{\nu },y_b)$ is well-defined, and
$S-  \{ (x_{\nu },y_b)\}$ is again either a dipper, or a full tail. 

The product of the linear orders on $X$ and $Y$ is a partial order on $X\times Y$ where
$(x_a,y_b)< (x_{a'},y_{b'})$ whenever $a\leq a'$, $b\leq b'$ and $(a,b)\neq (a',b')$. Given a dipper
$S$, let $(x_{\nu },y_b)$ be the corner of the notch and define new subsets 
$$
S_{>(x_{\nu },y_b)<} := \{ (x',y')\in S, \mbox{ either } (x',y')< (x_{\nu },y_b) \mbox{ or } (x_{\nu },y_b) < (x',y') \}
$$
and
$$
S_{\geq (x_{\nu },y_b)\leq } := \{ (x',y')\in S, \mbox{ either } (x',y')\leq  (x_{\nu },y_b) \mbox{ or } (x_{\nu },y_b) \geq (x',y') \}.
$$
Note that
$$
S_{\geq (x_{\nu },y_b)\leq } = S_{>(x_{\nu },y_b)<} \cup \{ (x_{\nu },y_b) \} .
$$

\newpage

For example, for the dipper pictured previously with $a=5,b=3,\nu = 3,m=9,n=7$, we have the picture
for $S_{>(x_{\nu },y_b)<}$
$$
{\setlength{\unitlength}{.5mm}
\begin{picture}(100,85)
\put(0,5){\line(0,1){20}}
\put(0,5){\line(1,0){30}}
\put(30,5){\line(0,1){20}}
\put(0,25){\line(1,0){30}}

\multiput(0,5)(10,0){4}{\circle*{1.5}}
\multiput(0,15)(10,0){4}{\circle*{1.5}}
\multiput(0,25)(10,0){4}{\circle*{1.5}}

\put(50,35){\circle*{1.5}}
\put(40,35){\circle*{1.5}}

\put(32,58){\ensuremath{(x_a,y_b)}}
\qbezier(43,55)(42,45)(48,37)
\put(48,37){\vector(2,-3){0}}

\put(10,46){\ensuremath{(x_{\nu +1},y_b)}}
\qbezier(22,43)(22,39)(38,35)
\put(38,35){\vector(3,-1){0}}

\put(50,23){\ensuremath{(x_{\nu},y_{b-1})}}
\put(47,25){\vector(-1,0){14}}

\put(30,25){\line(1,1){10}}
\put(40,35){\line(1,0){10}}

\put(50,35){\line(1,0){10}}
\put(60,35){\circle*{1.5}}

\put(60,35){\line(0,1){10}}
\put(60,45){\circle*{1.5}}

\put(60,45){\line(1,1){10}}
\put(70,55){\circle*{1.5}}

\put(70,55){\line(1,0){10}}
\put(80,55){\circle*{1.5}}

\put(80,55){\line(0,1){10}}
\put(80,65){\circle*{1.5}}

\put(80,65){\line(1,1){10}}
\put(90,75){\circle*{1.5}}
\put(85,78){\ensuremath{(x_m,y_n)}}
\end{picture}
} 
$$
and the picture for $S_{\geq(x_{\nu },y_b)\leq}$
$$
{\setlength{\unitlength}{.5mm}
\begin{picture}(100,85)
\put(0,5){\line(0,1){20}}
\put(0,5){\line(1,0){30}}
\put(30,5){\line(0,1){20}}
\put(0,25){\line(1,0){30}}

\multiput(0,5)(10,0){4}{\circle*{1.5}}
\multiput(0,15)(10,0){4}{\circle*{1.5}}
\multiput(0,25)(10,0){4}{\circle*{1.5}}

\put(50,35){\circle*{1.5}}
\put(40,35){\circle*{1.5}}
\put(30,35){\circle*{1.5}}

\put(32,58){\ensuremath{(x_a,y_b)}}
\qbezier(43,55)(42,45)(48,37)
\put(48,37){\vector(2,-3){0}}

\put(10,46){\ensuremath{(x_{\nu},y_b)}}
\qbezier(22,43)(22,41)(28,37)
\put(28,37){\vector(4,-3){0}}

\put(30,25){\line(0,1){10}}
\put(30,35){\line(1,0){20}}

\put(50,35){\line(1,0){10}}
\put(60,35){\circle*{1.5}}

\put(60,35){\line(0,1){10}}
\put(60,45){\circle*{1.5}}

\put(60,45){\line(1,1){10}}
\put(70,55){\circle*{1.5}}

\put(70,55){\line(1,0){10}}
\put(80,55){\circle*{1.5}}

\put(80,55){\line(0,1){10}}
\put(80,65){\circle*{1.5}}

\put(80,65){\line(1,1){10}}
\put(90,75){\circle*{1.5}}
\put(85,78){\ensuremath{(x_m,y_n)}}
\end{picture}
} 
$$

\begin{lemma}
If $S$ is a dipper with notched box ${\bf B}^{\nu}_{a,b}$ for $0\leq \nu <a$ and $0<b$, then 
either $\nu > 0$ and $b>1$ in which case $S_{>(x_{\nu },y_b)<}$ and 
$S_{\geq (x_{\nu },y_b)\leq }$ are dippers with
notched box ${\bf B}^0_{\nu , b-1}$ and tails going from $(x_{\nu}, y_{b-1})$ to $(x_m,y_n)$;
or else either $\nu =0$ or $b=1$ in which case $S_{>(x_{\nu },y_b)<}$ and $S_{\geq (x_{\nu },y_b)\leq }$ are full tails.
\end{lemma}
\begin{proof}
Look at the above pictures. 
\end{proof}

The idea of the proof of the product property 
is to consider subsets $S$ which are dippers or full tails,
and prove that $\bfi (S)^{\ast}(\pA \boxtimes \pB )\rightarrt \bfi (S)^{\ast}(\pA '\boxtimes \pB ')$ 
is a weak equivalence in $\precat (S;\mM )$. 

We first discuss the case of a full tail. Then the case of dippers will be treated by induction 
on $a,b,\nu$, eventually getting to the case $a=m,b=n$ which gives the desired theorem. 
The case of full tails, treated in the following proposition, encloses the case of all increasing 
single-step paths going from $(0,0)$ to $(m,n)$. This formalizes the intuition that to understand the
product we should understand what happens on each path. This is similar to what is going on in the 
decomposition of a product of simplices: the product has a decomposition into simplices indexed by the
same collection of paths---although the simplices of maximal dimension correspond to paths without
diagonal steps.

\begin{proposition}
\label{tailcase}
Suppose $(X,\pA )$ and $(Y,\pB )$ are sequentially free $\mM$-precategories.
Suppose $T= T^{i_0,\ldots , i_p}_{j_0,\ldots , j_p}$ is a full tail.
The induced ordering on $T$
is a linear order, and
$\bfi (T)^{\ast}(\pA \boxtimes \pB )$ is a sequentially free $\mM$-precategory
on the linearly ordered object set $T$. Furthermore, suppose $\pA \rightarrt \pA '$ and
$\pB \rightarrt \pB '$ are weak equivalences towards sequentially free precategories
satisfying the Segal condition, in the model categories of Theorem \ref{modstrucs}. 
Then 
$$
\bfi (T)^{\ast}(\pA \boxtimes \pB )\rightarrt \bfi (T)^{\ast}(\pA '\boxtimes \pB ')
$$
is a weak equivalence whose target satisfies the Segal condition. 
The same is also true of $\bfi (T)^{\ast}_!(\pA \boxtimes \pB )$ whose object set is $X\times Y$.
\end{proposition}
\begin{proof}
An increasing sequence of objects of $T$ has the form 
$(z_0,\ldots , z_r)$ where $z_k= (x_{u(k)},y_{v(k)})$
with $u(k)=i_{a(k)}$ and $v(k)=j_{a(k)}$ for 
$0\leq a(0)\leq \ldots \leq  a(r)\leq p$ an increasing sequence in the set
of indices for the objects of $T$. 
The tail condition means that for any such increasing sequence, the
sequences $x_{u(0)},\ldots , x_{u(r)}$ and $y_{u(0)},\ldots , y_{u(r)}$
are increasing sequences in $X$ and $Y$ respectively. 
Now
$$
\bfi (T)^{\ast}(\pA \boxtimes \pB )(z_0,\ldots , z_r) = \pA (x_{u(0)},\ldots , x_{u(r)})\times  \pB (y_{u(0)},\ldots , y_{u(r)})
$$
whereas 
$$
\bfi (T)^{\ast}(\pA \boxtimes \pB )(z_0, z_r) = \pA (x_{u(0)}, x_{u(r)})\times  
\pB (y_{u(0)}, y_{u(r)})
$$
so the fact that $\pA $ and $\pB $ are sequentially free implies, via the cartesian condition for $\mM$,  that 
the map 
$$
\bfi (T)^{\ast}(\pA \boxtimes \pB )(z_0,\ldots , z_r)\rightarrt \bfi (T)^{\ast}(\pA \boxtimes \pB )(z_0, z_r)
$$
is a weak equivalence. 

Suppose $\pA \rightarrt \pA '$ and $\pB \rightarrt \pB '$ are weak equivalences
towards sequentially free precategories satisfying the Segal condition. 
It follows from Corollary \ref{seqfreeinvariance} that for any
adjacent pair of objects $x_{i-1},x_i\in X$, the map
$$
\pA (x_{i-1},x_i)\rightarrt \pA '(x_{i-1},x_i)
$$
is a weak equivalence. The sequentially free condition implies that $\pA (x_i,x_i)$
and $\pA '_(x_i,x_i)$ are contractible (see Remark \ref{seqfreeconstcontr}). In particular
for any object $x_i$ the map 
$$
\pA (x_{i},x_i)\rightarrt \pA '(x_{i},x_i)
$$
is a weak equivalence. The same two statements hold for $\pB $. But now the tail
property of $T$ says that an adjacent pair
of objects in $T$ is of the form $(x_i,y_j),(x_k, y_l)$
where $x_i,x_k$ are either adjacent objects or the same object in $X$,
and $y_j,y_l$ are either adjacent objects or the same object in $Y$.
It follows that the map 
$$
\pA (x_i,x_i)\times \pB (y_j,y_l)\rightarrt  \pA '(x_i,x_i)\times \pB '(y_j,y_l)
$$
is a weak equivalence. By Corollary \ref{seqfreemorphism}, the map
$$
\bfi (T)^{\ast}(\pA \boxtimes \pB )\rightarrt \bfi (T)^{\ast}(\pA '\boxtimes \pB ')
$$
is a weak equivalence. The set of objects of $T$ injects into $X\times Y$,
so upon pushing forward to precategories on object set $X\times Y$, again
$$
\bfi (T)^{\ast}_!(\pA \boxtimes \pB )\rightarrt \bfi (T)^{\ast}_!(\pA '\boxtimes \pB ')
$$
is a weak equivalence. 
\end{proof}

The next step is to note that when we remove the corner of the notch from a dipper $S$
we have a pushout diagram.
In order to give the statement with some generality, 
say that $\pC \in \precat (X\times Y,\mM )$ is {\em ordered} if 
$\pC ((x_{i_0},y_{j_0}),\ldots , (x_{i_p},y_{j_p}))=\emptyset$ whenever we have a
nonincreasing sequence $(x_{i_0},y_{j_0}),\ldots , (x_{i_p},y_{j_p})$ in the product order,
that is to say whenever there is some $k$ such that either $i_{k-1}>i_k$ or $j_{k-1}>j_k$. 
Let $\precat ^{\rm ord}(X\times Y; \mM )$ denote the category of ordered $\mM$-enriched precategories over $X\times Y$. 

\begin{proposition}
\label{spushout}
Suppose $S = {\bf B}^{\nu}_{a,b}\cup T$ is a dipper and $(x_{\nu} ,y_b)$ the corner of the notch.
Assume $0\leq \nu < a$ and $0<b$. Then the pushout expression for subsets of $X\times Y$
$$
S = (S-\{ (x_{\nu} ,y_b) \}) \cup ^{S_{>(x_{\nu },y_b)<}}S_{\geq (x_{\nu },y_b)\leq }
$$
extends to a pushout expression for any ordered $\pC \in \precat ^{\rm ord}(X\times Y; \mM )$: the square
$$
\begin{diagram}
\bfi (S_{>(x_{\nu },y_b)<}) ^{\ast}_!(\pC ) & \rightarr & \bfi (S_{\geq (x_{\nu },y_b)\leq }) ^{\ast}_!(\pC ) \\
\downarr & & \downarr \\
\bfi (S-\{ (x_{\nu} ,y_b) \}) ^{\ast}_!(\pC ) & \rightarr & \bfi (S) ^{\ast}_!(\pC )
\end{diagram}
$$
is a pushout square in the category $\precat (X\times Y; \mM)$. 
\end{proposition}
\begin{proof}
It suffices to verify this levelwise on $\Delta _{X\times Y}$, in other words
we have to verify it for every sequence of objects $(x_{i_0},y_{i_0}),\ldots , (x_{i_p},y_{i_p})$.
If either of $x_{\cdot}$ or $y_{\cdot}$ is not increasing then it is trivially
true, so we may assume that both are 
increasing. If none of the elements in the sequence are $(x_{\nu},y_b)$ then
it is again trivially true. So we may assume that there is a $j$ with 
$(x_{i_j},y_{i_j})= (x_{\nu},y_b)$. Then the full sequence is contained
in the region $S_{\geq (x_{\nu },y_b)\leq }$. It follows that both vertical
maps in the above diagram, at the level of the sequence $(x_{\cdot}, y_{\cdot})$,
are isomorphisms. This implies that the diagram is a pushout. 
\end{proof}

Looking at this before putting everything back onto the same object set $X\times Y$,
one can also say that the square
$$
\begin{diagram}
\bfi (S_{>(x_{\nu },y_b)<}) ^{\ast}(\pC ) & \rightarr & \bfi (S_{\geq (x_{\nu },y_b)\leq }) ^{\ast}(\pC ) \\
\downarr & & \downarr \\
\bfi (S-\{ (x_{\nu} ,y_b) \}) ^{\ast}(\pC ) & \rightarr & \bfi (S) ^{\ast}(\pC )
\end{diagram}
$$
is a pushout square in $\precat (\mM)$. This version of the statement perhaps explains better what is
going on, but is less useful to us since we are currently working in the model category structure
on $\precat (X\times Y; \mM)$ for a fixed set of objects $X\times Y$. 

\begin{corollary}
\label{prodcor}
Suppose $g:\pC \rightarrt \pC '$ is a map of ordered $\mM$-enriched precategories over object set $X\times Y$,
such that both $\pC $ and $\pC '$ are levelwise cofibrant, that is cofibrant in $\precat _{\rm inj}(X\times Y, \mM )$. 
Suppose that for any full tail $T$ going from $(0,0)$ to $(m,n)$ the map 
$$
\bfi (T)^{\ast}_! (\pC )\rightarrt \bfi (T)^{\ast}_! (\pC ')
$$
is a weak equivalence. Then $g$ is a weak equivalence.
\end{corollary}
\begin{proof}
We show by induction that for any dipper of the form $S = {\bf B}^{\nu}_{a,b}\cup T$, the map 
$$
\bfi (S)^{\ast}_! (\pC )\rightarrt \bfi (S)^{\ast}_! (\pC ')
$$
is a weak equivalence. The induction is by $(b,a-\nu )$ in lexicographic order. In the initial case $b=1$ and $\nu = a$,
$S$ is a tail so the inductive statement is one of the hypotheses. Suppose the statement is known for all dippers $S'$
corresponding to $(a',b',\nu ')$ with $b'<b$ or $b'=b$ and $a'-\nu ' < a-\nu$. Use the expression of Proposition \ref{spushout}:
$\bfi (S) ^{\ast}_!(\pC )$ and  $\bfi (S) ^{\ast}_!(\pC ')$ are respectively the pushouts of the top and bottom rows in the diagram
$$
\begin{diagram}
\bfi (S-\{ (x_{\nu} ,y_b) \}) ^{\ast}_!(\pC )&\leftarr &\bfi (S_{>(x_{\nu },y_b)<}) ^{\ast}_!(\pC )  &\rightarr & \bfi (S_{\geq (x_{\nu },y_b)\leq }) ^{\ast}_!(\pC )\\
\downarr && \downarr && \downarr \\
\bfi (S-\{ (x_{\nu} ,y_b) \}) ^{\ast}_!(\pC ')&\leftarr &\bfi (S_{>(x_{\nu },y_b)<}) ^{\ast}_!(\pC ')  &\rightarr & \bfi (S_{\geq (x_{\nu },y_b)\leq }) ^{\ast}_!(\pC ').
\end{diagram}
$$
The vertical maps in the diagram induce the given map $\bfi (S) ^{\ast}_!(\pC )\rightarrt \bfi (S) ^{\ast}_!(\pC ')$.

The horizontal
maps are cofibrations in the injective model structure $\precat _{\rm inj}(X\times Y, \mM )$. Indeed, on any sequence 
$(x_{i_0},y_{j_0}),\ldots , (x_{i_p},y_{j_p})$ of points in $X\times Y$, each of the horizontal maps is either the identity,
or the inclusion from $\emptyset$ to $\pC (x_{i_0},y_{j_0}),\ldots , (x_{i_p},y_{j_p})$ or $\pC '(x_{i_0},y_{j_0}),\ldots , (x_{i_p},y_{j_p})$.
By hypothesis, $\pC $ and $\pC '$ are levelwise cofibrant, so the inclusions from $\emptyset$ are cofibrations. This shows that the horizontal
maps are cofibrations. 

The inductive hypothesis applies to each of the vertical maps. This is seen by noting that the invariants $(b',a'-\nu ')$ for the
dippers $S-\{ (x_{\nu} ,y_b) \}$, $S_{>(x_{\nu },y_b)<}$ and $S_{\geq (x_{\nu },y_b)\leq }$ are strictly smaller than $(b,a-\nu )$
in lexicographic order:
\newline
---for $S-\{ (x_{\nu} ,y_b) \}$ we have $a'= a$, $\nu ' = \nu + 1$ and $b'=b$, unless $\nu = a-1$ in which case
$b' = b-1$; 
\newline
---for $S_{>(x_{\nu },y_b)<}$ and $S_{\geq (x_{\nu },y_b)\leq }$ we have $b'=b-1$.

Hence, by the inductive hypothesis, each of the vertical maps is a weak equivalence. Now 
$\precat _{\rm inj}(X\times Y,\mM )$ is left proper because it is
a left Bousfield localization cf Theorem \ref{modstrucs}, which
implies by Corollary \ref{pushoutequivcor} that cofibrant pushouts are preserved by weak equivalences (see alternatively Lemma \ref{obisopushoutinvariance}). 
This shows that the map on pushouts
$\bfi (S) ^{\ast}_!(\pC )\rightarrt \bfi (S) ^{\ast}_!(\pC ')$ is a weak equivalence. This completes the inductive proof.

At the last case $S=X \times Y$ we obtain the conclusion of the corollary, that $g:\pC \rightarrt \pC '$ is a weak equivalence. 
\end{proof}

This gives the first main result of this chapter. 

\begin{theorem}
\label{sfprodthm}
Suppose $(X,\pA )$ and $(Y,\pB )$ are  sequentially free $\mM$-precategories, cofibrant in $\precat _{\rm inj}(X\times Y, \mM )$.
Suppose $\pA \rightarrt \pA '$ and $\pB \rightarrt \pB '$ are trivial cofibrations in $\precat _{\rm inj}(X, \mM )$
and $\precat _{\rm inj}(Y, \mM )$ respectively,
towards sequentially free $\mM$-precategories. Then the map
$$
(X\times Y, \pA \boxtimes \pB ) \rightarrt (X\times Y, \pA '\boxtimes \pB ')
$$
is a trivial cofibration in  $\precat _{\rm inj}(X\times Y, \mM )$.
\end{theorem}
\begin{proof}
Suppose first that $\pA '$ and $\pB '$ satisfy the Segal condition. In this case, 
Corollary \ref{prodcor} applies with $\pC =\pA \boxtimes \pB $ and $\pC '=\pA '\boxtimes \pB '$. The hypothesis on full tails used in Corollary \ref{prodcor}
is provided by Proposition \ref{tailcase}.

If $\pA '$ and $\pB '$ do not themselves satisfy the Segal condition, we can choose further morphisms $\pA '\rightarrt \pA ''$ and 
$\pB '\rightarrt \pB ''$ which are trivial cofibrations in $\precat _{\rm inj}(X, \mM )$
and $\precat _{\rm inj}(Y, \mM )$ respectively, such that $\pA ''$ and $\pB ''$ satisfy the Segal condition. The first case of this proof then
applies to the maps from $\pA $ and $\pB $, and also to the maps from $\pA '$ and $\pB '$. These show that 
$$
(X\times Y, \pA \boxtimes \pB ) \rightarrt (X\times Y, \pA ''\boxtimes \pB '')
$$
and
$$
(X\times Y, \pA '\boxtimes \pB ') \rightarrt (X\times Y, \pA ''\boxtimes \pB '')
$$
are trivial cofibrations in $\precat _{\rm inj}(X\times Y, \mM )$. By 3 for 2 it follows that 
$$
(X\times Y, \pA \boxtimes \pB ) \rightarrt (X\times Y, \pA '\boxtimes \pB ')
$$ 
is a weak equivalence in $\precat _{\rm inj}(X\times Y, \mM )$, and it is a trivial cofibration by the cartesian property of $\Mm$ (conditions (PROD) and (DCL), applied
as in Lemma \ref{prodproperties}).
\end{proof}

\section{Products of general precategories}

The next step is to extend the result of Theorem \ref{sfprodthm} from the sequentially free case, to the product of arbitrary $\mM$-enriched precategories.

Recall that we  have defined morphisms 
$$
\Upsilon _k(B_1,\ldots , B_k)\rightarrt \Upstild _k(B_1,\ldots , B_k)
$$
of sequentially free $\mM$-precategories, the target satisfies the Segal conditions, and the morphism is a weak equivalence in $\precat ([k],\mM )$ by Theorem \ref{upstildequiv}.
Use the notation
$$
\Sigma ([k]; B):= \Upsilon _k(B,\ldots , B)
$$
where the same object $B$ occurs $k$ times. 

\begin{lemma}
In the case where $B_1=\ldots =B_k=B\in \mM$, the map $\Upsilon \rightarrt \Upstild$ factors as
$$
\Sigma ([k]; B)\rightarrt h([k]; B)\rightarrt \Upstild _k(B,\ldots , B),
$$
and both maps are global weak equivalences between sequentially free $\mM$-enriched precategories.
\end{lemma}
\begin{proof}
All three are sequentially free, and Corollary \ref{seqfreemorphism} applies. 
For $\Sigma ([k]; B)$
and $\Upstild _k(B,\ldots , B)$ this is exactly what we said
in the proof of Theorem \ref{upstildequiv}; for $h([k]; B)$
see the explicit description in Section \ref{sec-precatexamples}. 
\end{proof}

\begin{proposition}
\label{sigmahtoseg}
Suppose $\pA\in \precat (\mM )$. Then we can obtain a global weak equivalence $\pA\rightarrt \pA'$ such that $\pA'$ satisfies the Segal conditions
i.e. $\pA'\in \Rr$, by taking a transfinite composition of pushouts along morphisms of the form 
\begin{equation}
\label{sigmahgens}
\Sigma ([k]; V)\cup ^{\Sigma ([k]; U)}h([k];U)
\rightarrt h([k]; V)
\end{equation}
for generating cofibrations $U\rightarrt^f V$ in $\mM$.
\end{proposition}
\begin{proof}
The maps \eqref{sigmahgens} are the same as the $\Psi ([k],f)$ considered
in Corollary \ref{psitoseg} and which make up the new pieces in $K_{\rm Reedy}$.
Note that this collection is missing the piece of $K_{\rm Reedy}$ consisting
of the generators for levelwise trivial Reedy cofibrations.
However, for the present statement that piece is not needed:
if we apply the small object argument to the present collection of morphisms
we can obtain a map $\pA\rightarrt \pA'$ which is a transfinite composition of
pushouts along morphisms of the form \eqref{sigmahgens},
such that $\pA'$ satisfies the left lifting property with respect to this collection.
The pushouts in question preserve weak equivalences, indeed the maps are a part of 
$K_{\rm Reedy}$ so that follows from the construction of
the model structure of Theorem \ref{reedyprecat} by direct left Bousfield
localization; or else one could apply
Theorem \ref{obisopushout} and Lemma \ref{globaltranfinite} which is
really saying pretty much the same thing.
Now, an object which satisfies the left lifting property with respect to
the  $\Psi ([k],f)$, satisfies the Segal conditions because the product maps
satisfy lifting along any generating cofibration $U\rightarrt V$ for $\mM$, thus they
are trivial fibrations in $\mM$, which shows that $\pA '$ satisfies the Segal
condition. 
\end{proof}

Notice that in the  construction of the proposition, the resulting
$\pA'$ will not in general be even levelwise fibrant, one would have to include
pushouts along morphisms of the form $h([k],f)$ for $f$ a generating
trivial cofibration of $\mM$. 
 
\begin{theorem}
\label{AtimesSigma}
Suppose $\pA  \in \precat (\mM )$ is Reedy cofibrant, $k\in \nn$ and $B\in \mM$ is a cofibrant object. Then the map
$$
\pA \times \Sigma ([k]; B)\rightarrt \pA \times h([k]; B)
$$
is a global weak equivalence.  
\end{theorem}
\begin{proof}
The proof goes in several steps. 
\newline
(i)\,\, Suppose $\pA$ is a sequentially free $\Mm$-enriched precategory. Note that 
$\Sigma ([k]; B)\rightarrt  h([k]; B)$ is a trivial cofibration of sequentially free $\Mm$-enriched precategories, inducing an order-preserving
isomorphism on objects. Apply Theorem \ref{sfprodthm} to this map and the identity of $\pA$, to conclude that
$$
\pA \times \Sigma ([k]; B)\rightarrt \pA \times h([k]; B)
$$
is a global weak equivalence.  This completes the proof when $\pA$ is sequentially free. This applies in particular to
the $h([k]; B)$ which are sequentially free.

\noindent
(ii)\,\, 
Recall from Theorem \ref{obisopushout}, that we already know that any pushout along a global trivial cofibration inducing an isomorphism on sets of objects,
is again a global trivial cofibration. 

\noindent
(iii)\,\, 
Suppose we know the statement of the theorem for $\pA$, $\pA '$ and $\pA ''$ and suppose given a diagram in which one of the arrows is at least an injective (i.e. levelwise)
cofibration
$$
\pA' \leftarrow \pA \rightarrt \pA '',
$$
then we claim that the statement of  the theorem is true for $\pQ := \pA '\cup ^{\pA} \pA ''$. Indeed,
$$
\pQ \times \Sigma ([k]; B) =
(\pA '\times \Sigma ([k]; B))\cup ^{\pA\times \Sigma ([k]; B)} (\pA ''\times \Sigma ([k]; B))
$$
by commutation of colimits and direct products in $\precat (\Mm )$ (which is part (DCL)
of the cartesian condition \ref{def-cartesian}).
The same is true for the product with $h([k]; B)$. The map
\begin{equation}
\label{Qmap}
\pQ \times \Sigma ([k]; B)\rightarrt \pQ \times h ([k]; B)
\end{equation}
is therefore obtained by functoriality of the pushout of the columns in 
$$
\begin{diagram}
\pA ' \times \Sigma ([k]; B) &\rightarr & \pA '\times h([k]; B)\\
\uparr & & \uparr \\
\pA \times \Sigma ([k]; B) &\rightarr & \pA \times h([k]; B)\\
\downarr & &\downarr \\
\pA ''\times \Sigma ([k]; B) &\rightarr & \pA ''\times h([k]; B) .
\end{diagram}
$$
We are supposing that we know that each of the horizontal maps is a global weak equivalence,
also they induce isomorphisms on objects. By the cartesian property for $\mM$
applied levelwise, the same one of the vertical maps is a levelwise cofibration.

By Lemma \ref{obisopushoutinvariance}, the induced map on pushouts \eqref{Qmap}
is a global weak equivalence. 
This proves the claim for step (iii).

\noindent
(iv) \,\, Suppose given a sequence $\pA _i$ indexed by an ordinal $\beta$, with injective 
cofibrant transition maps. 
Suppose the statement of the theorem is true for each $\pA _i$,
then it is true for $\pQ := \colim _{i\in \beta } \pA _i$.  Indeed, just as in the previous part $\pQ \times h([k]; B)$ can be expressed
as a transfinite composition of pushouts of $\pQ \times \Sigma ([k]; B)$ along maps which are by hypothesis global trivial cofibrations
which induce isomorphisms on objects. By Theorem \ref{obisopushout} and 
Lemma \ref{globaltranfinite}, the composition is
a global weak equivalence. 

\noindent
(v)\,\, We show by induction on $m\in \nn$ that if $\sk _m (\pA ) \cong \pA$ then the statement of the theorem holds
for $\pA$. It is easy to see in case $m=0$ because then $\pA$ is just a discrete set.
Suppose this is known for any $m\leq n$, and suppose $\pA = \sk _n(\pA )$. By 
Proposition \ref{expression} we can express $\pA $ as a transfinite composition of
pushouts of
$\sk _{n-1}(\pA )$ along maps of the form $h([n],\partial [n]; U\rightarrt V)\rightarrt h([n]; V)$. On the other hand,
$$
h([n],\partial [n]; U\rightarrt V) = h([n]; U)\cup ^{h(\partial [n];U)} \partial h(\partial [n];V),
$$
and $h(\partial [n];U) = \sk _{n-1}h([n]; U)$.  By the inductive hypothesis 
the statement of the theorem is known for $h(\partial [n];U)$ and similarly for $h(\partial [n];V)$.
It is known for $h([n]; U)$ by (i). So by (iii) the statement of the theorem is known for 
$h([n],\partial [n]; U\rightarrt V)$. Furthermore it is known for $\sk _{n-1}(\pA )$ by the inductive hypothesis. 
Again by (iii) and (iv) we conclude the statement for $\pA$. 

\noindent
(vi) \,\, Any Reedy cofibrant $\pA$ can be expressed as a transfinite composition of the
maps $\sk _m(\pA )\rightarrt \sk_{m-1}(\pA )$, so by (iv) and (v)
we get the statement of the theorem for any Reedy cofibrant $\pA$. This completes the proof.
\end{proof}

Recall from Corollary \ref{psitoseg} and the remark at the
beginning of the proof of Proposition \ref{sigmahtoseg} above, 
for any cofibration $f:U\rightarrt V$ we have the notation
$$
\src \Psi ([k],f) = \Sigma ([k]; V)\cup ^{\Sigma ([k]; U)}h([k];U)
$$
and the map $\Psi ([k],f)$ goes from here to $h([k]; V)$.

\begin{corollary}
\label{AtimesPsi}
Suppose $\pA  \in \precat (\mM )$ is Reedy cofibrant, $k\in \nn$, and 
$f:U\rightarrt V$ is a cofibration in $\mM$. Then the map
$$
\pA \times \src \Psi ([k],f)\rightarrt \pA \times h([k]; v)
$$
is a global weak equivalence. 
\end{corollary}
\begin{proof}
In the cocartesian diagram
$$
\begin{diagram}
\pA \times \Sigma ([k]; U) & \rightarr & \pA \times h([k];U) \\
\downarr && \downarr \\
\pA \times \Sigma ([k]; V) & \rightarr & \pA \times \src \Psi ([k],f)
\end{diagram}
$$
the upper arrow is a global weak equivalence inducing an isomorphism on the set of
objects, by the previous theorem. Furthermore it is a Reedy cofibration by 
Proposition \ref{reedyproduct}, so 
it is a Reedy isotrivial cofibration. By Theorem \ref{obisopushout}, 
the bottom map is a global weak equivalence. The statement of the corollary
now follows by again using the previous Theorem \ref{AtimesSigma} as well as 3 for 2.
\end{proof}

\begin{theorem}
\label{AtimesB}
Assume $\mM$ is a tractable left proper cartesian model category.
For any $\pA , \pB \in \precat (\mM )$, the map
$$
\pA \times \pB \rightarrt \Seg (\pA )\times \Seg (\pB )
$$
is a global weak equivalence.
\end{theorem}
\begin{proof}
We suppose first that $\pA$ and $\pB$ are Reedy cofibrant. 
There is a map $\pA \rightarrt \pA '$ to an object satisfying the Segal conditions,
which is a transfinite composition of pushouts along 
morphisms of the form $\src \Psi ([k],f)\rightarrt^{\Psi ([k],f)} h([k]; V)$. Each of these pushouts is a global trivial Reedy cofibration
inducing an isomorphism on the set of objects, by Theorem \ref{obisopushout}. The map 
$$
\pA \times \pB \rightarrt \pA '\times \pB 
$$
is the corresponding transfinite composition of pushouts along morphisms of the form 
$$
\src \Psi ([k],f)\times \pB \rightarrt h([k]; U)\times \pB.
$$
This is because part of the cartesian hypothesis for $\mM$ is (DCL) commutation of direct products and colimits. 
By Corollary \ref{AtimesPsi}, the morphisms 
$\src \Psi ([k],f)\times \pB \rightarrt h([k]; U)\times \pB$ are global weak equivalences; they are also
Reedy cofibrations by Proposition \ref{reedyproduct} because we assumed that $\pB$ is Reedy cofibrant.
These maps are again isomorphisms on objects, so we can 
apply Theorem \ref{obisopushout} which says that global trivial cofibrations which induce isomorphisms
on the set of objects are preserved under pushout (and see Lemma \ref{globaltranfinite} for the transfinite composition). 
Therefore $\pA \times \pB \rightarrt \pA '\times \pB $ is a global weak equivalence. 

Arguing in the same way for the product of a map $\pB \rightarrt \pB '$ with $\pA'$,
then composing the two equivalences we conclude that the map
$$
\pA \times \pB \rightarrt \pA '\times \pB '
$$
is a global weak equivalence. On the other hand, $\pA ' \rightarrt \Seg (\pA ')$ and
$\pB ' \rightarrt \Seg (\pB ')$ are levelwise weak equivalences,
so 
$$
\pA '\times \pB '\rightarrt \Seg (\pA ')\times \Seg (\pB ')
$$
is a levelwise  weak equivalence. Similarly, the fact that $\pA \rightarrt \pA'$ is
a global weak equivalence inducing an isomorphism on objects, implies that 
$$
\Seg (\pA )\rightarrt \Seg (\pA ')
$$
is a levelwise weak equivalence, and by the same remark for $\pB$ then taking the product, we get that
$$
\Seg (\pA )\times \Seg (\pB )\rightarrt \Seg (\pA ')\times \Seg (\pB ')
$$
is a levelwise weak equivalence. 
Thus we obtain a diagram
$$
\begin{diagram}
\pA \times \pB & \rightarr & \pA '\times \pB ' \\
\downarr & & \downarr \\
\Seg (\pA )\times \Seg (\pB )& \rightarr & \Seg (\pA ')\times \Seg (\pB ')
\end{diagram}
$$
where the top map is a global weak equivalence, and the right vertical and bottom maps
are levelwise hence global weak equivalences (Lemma \ref{levelwiseglobalwe}). 
By 3 for 2 it follows that the left vertical
map is a global weak equivalence as required for the theorem. This completes the proof for
the case of Reedy cofibrant objects. 

Now suppose $\pA$ and $\pB$ are general objects of $\precat (\mM )$. Consider Reedy cofibrant replacements 
$\pA '\rightarrt \pA$ and $\pB '\rightarrt \pB$; these may be chosen as levelwise equivalences of diagrams,
which are then global weak equivalences by Lemma \ref{levelwiseglobalwe}.
In particular,  $\Seg (\pA ')\rightarrt \Seg (\pA )$ is a levelwise weak equivalence and the same for $\pB '$ so 
$$
\Seg (\pA ')\times \Seg (\pB ') \rightarrt  \Seg (\pA )\times \Seg (\pB )
$$
is a levelwise weak equivalence. 
Note also that 
$$
\pA ' \times \pB ' \rightarrt \pA \times \pB 
$$
is a levelwise weak equivalence of diagrams over $\Delta ^o _{\Ob (\pA )\times \Ob (\pB )}$, so it is a global
weak equivalence by Lemma \ref{levelwiseglobalwe}.
The first part of the proof treating the Reedy cofibrant case shows that
$$
\pA ' \times \pB ' \rightarrt  \Seg (\pA ')\times \Seg (\pB ')
$$
is a global weak equivalence. In the square diagram
$$
\begin{diagram}
\pA ' \times \pB ' & \rightarr & \pA \times \pB \\
\downarr & & \downarr \\
\Seg (\pA ')\times \Seg (\pB ') & \rightarr & \Seg (\pA )\times \Seg (\pB )
\end{diagram}
$$
the top, bottom and left vertical arrows are global weak equivalences, so by 3 for 2 the right vertical arrow is
a global weak equivalence. This completes the proof.
\end{proof}

\begin{corollary}
\label{prodwe}
Suppose $\pA \rightarrt \pB$ and $\pC \rightarrt \pD$ are global weak equivalences. Then the map
$$
\pA \times \pC \rightarrt \pB \times \pD
$$
is a global weak equivalence. 
\end{corollary}
\begin{proof}
Suppose first of all that $\pA$, $\pB$, $\pC$ and $\pD$ are  objects satisfying the Segal conditions. 
Then the  products also satisfy the Segal conditions. Truncation of these is compatible with direct products, by Lemma \ref{truncationproducts},
so the map in question is essentially surjective. By looking at the morphism objects we see that it is fully faithful,
so it is a global weak equivalence by following the definition. 

Next suppose that $\pA$, $\pB$, $\pC$ and $\pD$ are
any objects, and
look at the diagram 
$$
\begin{diagram}
\pA \times \pC & \rightarr & \pB \times \pD \\
\downarr & & \downarr \\
\Seg (\pA )\times \Seg (\pC )& \rightarr & \Seg (\pB)\times \Seg (\pD ).
\end{diagram}
$$
The vertical maps are global weak equivalences by Theorem \ref{AtimesB}, while the bottom map is a global weak equivalence
by the first paragraph of the proof. By 3 for 2, the top map is a global weak equivalence.
\end{proof}

\section{The role of unitality, degeneracies and higher coherences}
\label{sec-role}

In this section, we point out why we need to impose the unitality condition $A(x_0)=\ast$,
to include the degeneracy maps in $\Delta$ (which also correspond to some sort of
unit condition), and why we can't truncate $\Delta$ by, say, dropping the objects
$[n]$ for $n\geq 4$. These all have to do with the arguments of this chapter about
products. In some sense it goes back to the Eilenberg-Zilber theorem; our product
condition can be viewed as a generalization to
the present context where the information of direction of arrows is retained.

\subsection{The unitality condition}
\label{sec-unitalnecessary}

Suppose we tried to use non-unital precategories. These would be pairs $(X,\pA )$
where $\pA :\Delta ^o_X\rightarrt \mM$ is an arbitrary functor. The Segal condition
would include, for sequences of length $n=0$, the fact that $\pA (x_0)\rightarrt \ast$ 
should be a weak equivalence, in other words $\pA (x_0)$ is weakly contractible.
So, this would constitute a weak version of the unitality condition. We would proceed much
as above, imposing the Segal condition by the small object argument in an operation
denoted $\pA \mapsto \Seg ^n (\pA )$. It seems likely
that this would lead to a model category, conjecturally Quillen equivalent to the
model categories on unital precategories which we are constructing here. 

However, even if the model structure existed, it could not be cartesian. The 
reason for this occurs at some very degenerate objects: consider the non-unital
precategory with object set $\Ob (\pB )= \{ y\}$ a singleton, but with functor 
the constant functor with values the initial object:
$$
\pB (y,\ldots , y):= \emptyset .
$$
This includes the case of sequences of length $0$: $\pB (y)=\emptyset$,
so $\pB$ doesn't satisfy the unitality condition. Now $\Seg ^n( \pB )$
would be some kind of $\mM$-enriched category with a single object;
it seems clear that it would be the coinitial $\ast$ but in any case has
to contain $\ast$ as a retract.

Suppose $\pA $ is another non-unital precategory (which might in fact be
unital). Consider $\pA \times \pB$. The object set is $\Ob (\pA )\times \{ y\}\cong \Ob (\pA )$. But for any sequence $((x_0,y),\ldots , (x_n,y))$ we have
$$
(\pA \times  \pB )((x_0,y),\ldots , (x_n,y)) =
\pA (x_0,\ldots ,x_n)\times \pB (y,\ldots , y) 
$$
$$
=\pA (x_0,\ldots ,x_n)\times \emptyset =\emptyset .
$$
In particular, the structure of $\pA \times \pB$ depends only on $\Ob (\pA )$ and
not on $\pA$ itself. This would be incompatible with the cartesian condition (for any reasonable choice of $\mM$), because
$$
\Seg ^n(\pA )\times \ast \rightarrt \Seg ^n(\pA )\times \Seg ^n(\pB )
$$
is contained as a retract, but $\Seg ^n(\pA \times \pB )$ is essentially trivial.

To make the last step of the above argument precise we would need to investigate
$\Seg ^n$ explicitly. In the case $\mM = \Sets$, the same discussion as 
in Section \ref{genrel1cat} applies, expressing $\Seg ^n( \pA )$ as the
category generated by $\pA$ considered as a system of generators and relations; the
first step would be to impose the Segal condition at $n=0$ which, for $\mM =\Sets$,
is exactly the unitality condition; from there the rest is the same. In this case
we see that $\Seg ^n(\pB )$ is really just $\ast$, $\Seg ^n(\pA \times \pB )$ is
the discrete category on object set $\Ob (\pA )\times \{ y\}$, and it cannot
contain $\Seg ^n(\pA )$ as a retract in general.

\subsection{Degeneracies}

Let $\Phi \subset \Delta$ denote the category consisting only of face maps, in other words the objects of $\Phi$ are the 
nonempty finite linearly ordered sets $[k]$ for $k\in \nn$, but the maps are the injective order-preserving maps. 
We could try to create a theory of weak categories based on $\Phi$ rather than $\Delta$. It would be appropriate to call these 
``weak semicategories'', because the degeneracy maps in $\Delta$ correspond to inserting identity morphisms into a composable sequence. 
This theory is undoubtedly interesting and important, and has not been fully worked out as far as I know. 

This would undoubtedly be
related to the work of J. Kock on weakly unital higher categories \cite{Kock} \cite{Kock2}, as one could start by considering
weak semicategories, then impose a weak unitality condition.  

Unfortunately the theory of products again doesn't work if $\Delta$ is replaced by $\Phi$.
Indeed, a slight modification of the example of the previous subsection again provides
a counterexample. Let $\pB$ be the precategory with a single object $y$, 
with $\pB (y)=\ast$
so it is unital, but with $\pB (y,\ldots , y)=\emptyset$ for any sequence of length $n\geq 1$ (that is to say, with $n+1$ elements). This will still be a valid functor from $\Phi _{\{ y\}}$ to $\mM$, however, taking the product $\pA \times \pB$ will destroy the structure
of $\pA$. 

As in the previous subsection, this can be made precise in the case $\mM = \Sets$.
The non-unital precategories may then be considered as systems of generators and relations
for a category, but the system doesn't contain the degeneracies. 

It is interesting to look more closely at how this works in the case of systems of
generators for a monoid, that is to say for a category with a single object. 
The $1$-cells of a precategory $\pA$ 
correspond to generators of the monoid, and the $2$-cells correspond to
relations of the form $f=gh$ among the generators.  In this case, the system of generators
and relations is just given by two sets $\pA (x,x)$ and $\pA (x,x,x)$ with three maps
$$
\begin{diagram}
\pA (x,x,x) &\pile{\rightarr \\ \rightarr \\ \rightarr }& \pA (x,x).
\end{diagram}
$$
The process of generators and relations would have to include the addition of identities.

In this case, for example, if $\pA$ has a single generator and no relations, then
its product with itself $\pA \times \pA$ will again be a system with a single generator
and no relations; but $\pA$ generates the monoid $\nn$ and $\nn \times \nn$ is different
from $\nn$, so the product of systems of generators and relations doesn't generate the
product of the corresponding monoids. 

One can see in this simple example how the degeneracies come to the rescue. 
A system of generators and relations with unitality and 
degeneracies corresponds to a diagram of
the form 
$$
\begin{diagram}
\pA (x,x,x) &\pile{\rightarr \\ \leftarr \\ \rightarr \\ \leftarr \\ \rightarr }& \pA (x,x) & \pile{\rightarr \\
\leftarr \\ \rightarr} & \pA (x) = \ast .
\end{diagram}
$$
This means that there is an explicit element $1$ among the generators, with the
relations $1\cdot f = f$ and $f\cdot 1 = f$ for any other generator $f$.

Now let's look again at $\nn$ generated by a system $\pA$ consisting of
a single generator. We have $\pA (x,x)=\{ f, 1\}$ with relations
corresponding to the left and right identities for both $f$ and $1$ itself. 
The product now has
generators 
$$
\pA \times \pA (x,x) = \{ f\times f, f\times 1, 1\times f, 1\times 1\}
$$
where we have noted the single object of $\pA\times \pA$ as $x$ again rather than $(x,x)$. 
The unit generator is $1\times 1$. 
The relations include the left and right identities with the unit  generator $1\times 1$,
plus two new relations of the form
$$
(f\times 1) \cdot (1\times f)=f\times f
$$
and
$$
(1\times f) \cdot (f\times 1)=f\times f.
$$
The first of these two relations serves to eliminate the generator $f\times f$ so we get to a
monoid with two generators $f\times 1$ and $1\times f$, then the second relation
gives the commutativity
$(f\times 1)\cdot (1\times f) = (1\times f)\cdot (f\times 1)$. So, the monoid generated
by the system $\pA \times \pA$ is indeed $\nn \times \nn$.

Working out this example demonstrates how the degeneracies of $\Delta$ enter into
the cartesian condition in an important
way.

\section{Why we can't truncate $\Delta$}

The above examples could all be done in the case of $\mM = \Sets$, where the
passage from precategories to categories is the process of generating a category
by generators and relations. For that, we didn't need to consider the part of
$\Delta$ involving $[n]$ for $n\geq 4$ (the case $n=3$ being needed for the
associativity condition). 

On the other hand, 
for weak enrichment in a general model category $\mM$,
we can't replace $\Delta$ by any finite truncation, that is by a subcategory $\Delta ^{\leq m}$
of finite ordered sets of size $\leq m$. 

This can be seen by the requirement that there should be a higher Poincar\'e groupoid
construction; in the case when $\mM = \mK$ is the model category of simplicial sets,
the $\mK$-enriched precategories should be realizable into arbitrary homotopy types;
and in particular the Segal groupoids should be eqivalent to homotopy types by a
pair of functors including Poincar\'e-Segal category and realization. These should be
compatible with homotopy groups in a way similar to that described in Chapters \ref{ncats1}
and \ref{nonstrict1}.  

If we impose these conditions, it becomes easy to see that 
the Segal groupoids defined using only
$\Delta ^{\leq m}$ (say for $m\geq 3$) don't model all homotopy types.
This means that, for the program we are pursuing here, the category $\Delta ^{\leq m}$
cannot be sufficient.

It should also be possible to show that $\Delta ^{\leq m}$ can't be used to model
homotopy $n$-types for $n>m$, in any way at all; however, it doesn't seem completely
clear how to formulate a good statement of this kind. 

If for $1$-categories it suffices to look at 
$\Delta ^{\leq 3}$, we could expect more generally that in order to 
consider $n$-categories it would suffice to look at $\Delta ^{\leq n+2}$, indeed this
showed up in the explicit example of Chapter \ref{secat1} and will show up again in
our discussion of stabilization in a future version (see \cite{BBDSH}).


\chapter{Intervals}
\label{interval1}

Given our tractable, left proper and cartesian model category $\mM$,
the main remaining problem in order to construct the global model structure on 
$\precat (\mM )$ is to consider the notion of {\em interval} which should be
an $\mM$-precategory (to be called $\Xi (N|N')$ in our notations below),
weak equivalent 
to the usual category $\overline{\pE}$ with two isomorphic objects $\upsilon _0,\upsilon _1\in \overline{\pE}$,
and with a single morphism between any pair of objects. 

If $\pA \in \precat (\mM )$ is a weakly $\mM$-enriched category, an {\em internal equivalence}
between $x_0,x_1\in \Ob (\pA )$ is a ``morphism from $x_0$ to $x_1$'' (see \eqref{mordef} below), which projects to an isomorphism in the truncated category $\tau _{\leq 1}(\pA )$.
This terminology was introduced by Tamsamani in \cite{Tamsamani}. It plays a 
vital role in the study of global weak equivalences.
Essential surjectivity of a morphism $f:\pA \rightarrt \pB$
means (assuming that $\pB$ is levelwise fibrant)
that for any object $y\in\Ob (\pB )$, there is an object $x\in \Ob (\pA )$
and an internal equivalence between $f(x)$ and $y$.

Unfortunately, an internal equivalence between $x_0$ and $x_1$ in $\pA$
doesn't necessarily
translate into the existence of a morphism $\overline{\pE}\rightarrt \pA$.
This will work after we have established the model structure on $\precat (\mM )$
if we assume that $\pA$ is a fibrant object. However, in order to finish
the construction of the model structure, we should start with the weaker 
hypothesis\footnote{As observed by Bergner \cite{BergnerSegal} this hypothesis 
will be equivalent to fibrancy in the global projective model structure, once we know
that it exists.}
that $\pA$ satisfies the Segal conditions and is levelwise fibrant. 
The ``interval object'' $\Xi (N|N')$ should be contractible, and
have the versality property that
whenever  $x_0$ and $x_1$ are internally equivalent, there is a morphism 
$\Xi (N|N')\rightarrt \pA$ relating them. 

The construction of such a versal interval was the subject of an error in \cite{svk},
found and corrected by Pelissier in \cite{Pelissier}. 
This was somewhat
similar to a mistake in Dwyer-Hirschhorn-Kan's original construction of the model category structure for simplicial categories (\cite{DwyerHirschhornKan}, now \cite{DwyerHirschhornKanSmith}),
pointed out by Toen and subsequently fixed for simplicial categories by Bergner
\cite{BergnerModel}. Pelissier fixed this problem for the model category of Segal categories 
by constructing an explicit interval object and verifying its topological properties using the comparison between Segal $1$-groupoids and spaces.
Drinfeld has constructed an interval object for differential graded categories \cite{DrinfeldIntervalDG}. 

Pelissier's correction as written
covered only the case of $\mK$-enriched weak categories, and one of our purposes here is
to point out that 
his argument serves to construct the required intervals in general, by functoriality with respect to a left Quillen functor 
$\mK \rightarrt \mM$. 
For the main result which is contractibility of 
the $\Xi (N|N')$ we proceed therefore in two steps: first considering the problem for the case of Segal categories
i.e. $\Kk$-enriched weak categories as was done in \cite{Pelissier}; then going to the case of $\mM$-enriched weak categories 
by transfer along $\Kk \rightarrt \mM$. 
Sections \ref{sec-directtransfer}, \ref{changeenrichment},
and \ref{funcgen}
about transfer along a left Quillen functor 
were motivated by this movement. 
The possibility of doing that is 
one of the advantages of the fully iterative point of view originally suggested by
Andr\'e Hirschowitz in Pelissier's thesis topic, in which $\mM$ is a general input
into the construction. 
It should also be possible to adapt Pelissier's correction directly
to the original $n$-nerves considered in \cite{Tamsamani} \cite{svk},
by using Tamsamani's theorems on the topological realization of weak $n$-groupoids
which in turn applied Segal's original results in a partially iterative way.
That would be more geometrically motivated, 
but for the present treatment the fully 
iterative approach is both more general and more direct. 

I would like to thank Regis Pelissier for finding and correcting this problem.

\section{Contractible objects and intervals in $\mM$}

An object $A\in \mM$ is {\em contractible} if the unique morphism $A\rightarrt \ast$ is a weak equivalence.
An {\em interval object} is a triple $(B,i_0,i_1)$ where $B\in \mM$ and $i_0,i_1:\ast \rightarrt B$
such that $B$ is contractible and $i_0\sqcup i_1:\ast \cup ^{\emptyset} \ast \rightarrt B$ is a cofibration. 

Assumption (AST) in the cartesian condition \ref{def-cartesian} 
says that $\ast$ is a cofibrant object, so an interval object is itself cofibrant. 

A {\em morphism between intervals} from $(B,i_0,i_1)$ to $(B',i'_0,i'_1)$ means a morphism $f:B\rightarrt B'$
such that $f\circ i_0=i'_0$ and $f\circ i_1=i'_1$. Since $B$ and $B'$ are contractible, a morphism $f$ is automatically
a weak equivalence. 

\begin{lemma}
\label{intervalbeyond}
Suppose  $(B,i_0,i_1)$ and $(B',i'_0,i'_1)$ are two interval objects. Then there is a third one
$(B'',i''_0,i''_1)$ and morphisms of intervals $f:B\rightarrt B''$ and $f':B'\rightarrt B''$. 
These may be assumed to be trivial cofibrations. 
\end{lemma}
\begin{proof}
Put $A:= B\cup ^{\ast \cup ^{\emptyset}\ast} B'$ and choose a factorization
$$
A\rightarrt ^f B'' \rightarrt \ast
$$
where the first morphism is a cofibration and the second morphism is a weak equivalence. 
Now $i_0$ and $i'_0$ are the same when considered as maps $\ast \rightarrt A$ 
because of the coproduct in the definition of $A$. Thus $f\circ i_0=f\circ i'_0$ gives a map
$i''_0:\ast \rightarrt B''$. Similarly $i''_1:= f\circ i_1=f\circ i'_1$. 
The map $\ast \cup ^{\emptyset}\ast\rightarrt A$ is cofibrant, and since $f$ is cofibrant
the composition into $B''$ is cofibrant. Note that the maps $B\rightarrt A$ and $B'\rightarrt A$ are cofibrations,
so the same is true of the maps to $B''$, and since the intervals are weakly equivalent to $\ast$ these maps
are trivial cofibrations. 
\end{proof}

Recall that we defined in Chapter \ref{weakenr1} 
a functor $\tau _{\leq 0}: \mM \rightarrt \Sets$ by $\tau _{\leq 0}(A):= \Hom _{\Ho (\mM )}(\ast , A')$ where $A\rightarrt A'$ is a
fibrant replacement. 

\begin{lemma}
\label{tauinterval}
Suppose $A$ is a fibrant object and $a,b:\ast \rightarrt A$. The following conditions are equivalent:
\newline
(a)---The classes of $a$ and $b$ in $\tau _{\leq 0}(A)$ coincide;
\newline
(b)---for any interval object $(B,i_0,i_1)$ there exists a map $B\rightarrt A$ sending $i_0$ to $a$ and $i_1$ to $b$;
\newline
(c)---there exists an interval object $(B,i_0,i_1)$ and a map $B\rightarrt A$ sending $i_0$ to $a$ and $i_1$ to $b$.
\end{lemma}
\begin{proof}
This is an exercise in Quillen's theory of the homotopy category of a model category, which we do for the reader's convenience. 

Note that $\mbox{(b)}\Rightarrow \mbox{(c)}\Rightarrow \mbox{(a)}$ easily. Assume that $A$ is also cofibrant. 
To prove that $\mbox{(a)}\Rightarrow\mbox{(c)}$,
suppose that the classes of $a$ and $b$ coincide in $\tau _{\leq 0}(A)$. This is equivalent to saying that the two maps
$a,b:\ast \rightarrt A$ project to the same map in $\Ho (\mM )$. Recall from Quillen \cite{Quillen} that $\Ho (\mM )$ is
also the category of fibrant and cofibrant objects of $\mM$, with homotopy classes of maps. As $\ast$ is automatically fibrant, and
cofibrant by hypothesis; and we are assuming that $A$ is cofibrant and fibrant, then condition (a) says that the two maps $a$ and $b$
are homotopic in the sense of Quillen \cite{Quillen}, which says exactly condition (c). For the implication 
$\mbox{(a)}\Rightarrow\mbox{(c)}$, but with $A$ assumed only to be fibrant, choose a trivial fibration from a cofibrant object $A'\rightarrt A$.
Lift to maps $a',b':\ast \rightarrt A'$. Since $A'\rightarrt A$ projects to an isomorphism in $\Ho (\mM )$, the maps $a'$ and $b'$ are
equivalent in $\tau _{\leq 0}(A')$ so by condition (c) proven for $A'$ previously, there exists an interval object $(B,i_0,i_1)$ and an 
extension of $a'\sqcup b'$ to $B\rightarrt A'$. Composing gives the required map $B\rightarrt A$.

To finish the proof it suffices to show $\mbox{(c)}\Rightarrow\mbox{(b)}$. Suppose $(B,i_0,i_1)$ and $(B',i'_0,i'_1)$ are two
interval objects. Applying Lemma \ref{intervalbeyond} there is an interval object $(B'', i''_0,i''_1)$ with trivial cofibrations from both $B$ and $B'$.
If $a\sqcup b: \ast \cup ^{\ast}\ast \rightarrt A$ extends to a map $B\rightarrt A$, and if $A$ is fibrant, then
the lifting property for $A$ gives the extension to $B''\rightarrt A$ which then restricts to a map $B'\rightarrt A$ as required to show
$\mbox{(c)}\Rightarrow\mbox{(b)}$.
\end{proof}

Using the assumption that $\mM$ is cartesian,
we can make a similar statement explaining the relation of homotopy between
morphisms using an interval, if the target is a fibrant object of $\mM$. 

\begin{lemma}
\label{homotopicequiv}
Suppose $\mM$ is a cartesian model category. Suppose $A$ is a cofibrant object and
$C$ is a fibrant object. Then, for two morphisms $f,g:A\rightarrt B$ the
following statements are equivalent:
\newline
(a)---$f$ and $g$ are homotopic in Quillen's sense, meaning that 
the classes of $f$ and $g$ in $\Hom _{\Ho (\mM )}(A,C)$ coincide;
\newline
(b)---for any interval object $(B,i_0,i_1)$ there exists a map $h:A\times B\rightarrt C$
such that $h\circ (1_A\times i_0) = f$ and $h\circ (1_A\times i_1)=g$;
\newline
(c)---there exists an interval object $(B,i_0,i_1)$ and a map $h:A\times B\rightarrt C$
such that $h\circ (1_A\times i_0) = f$ and $h\circ (1_A\times i_1)=g$.
\end{lemma}
\begin{proof}
The cartesian property of $\mM$ implies that for any interval object $B$, the diagram 
$$
A\times (\ast \cup ^{\emptyset} \ast )=
A\cup ^{\emptyset} A\rightarrt A\times B\rightarrt A
$$
is an $A\times I$-object in Quillen's sense, and so can be used to measure
homotopy between our two maps.
\end{proof}

\section{Intervals for $\mM$-enriched precategories}

Let $\pE:= \Upsilon (\ast )$ denote the category with two objects $\upsilon _0,\upsilon _1$ and a single
morphism between them. Thus, $\pE(\upsilon _0,\ldots , \upsilon _0)=\ast$,
$\pE(\upsilon _1,\ldots , \upsilon _1)=\ast$, $\pE(\upsilon _0,\ldots ,\upsilon _0,\upsilon _1,\ldots , \upsilon _1)=\ast$
and the remaining values are $\emptyset$. This is the image of the usual category $[0\rightarrt 1]$
under the map $\precat (\Sets )\rightarrt \precat (\mM )$. 

An alternative description of $\pE$ in terms of the representable object notation of Section \ref{sec-precatexamples} is $\pE=h([1],\ast )$. 
If $\pA$ is any $\mM$-enriched precategory,
a map $\pE\rightarrt \pA$ is the same thing as
a triple $(x_0,x_1,a)$ where $x_0,x_1\in \Ob(\pA )$ and $a:\ast \rightarrt \pA(x_0,x_1)$
is an element of the ``set of morphisms from $x_0$ to $x_1$''. This ``set of morphisms''
may be denoted by
\begin{equation}
\label{mordef}
\Mor ^1_{\pA}(x_0,x_1):= \Hom _{\mM}(\ast , \pA (x_0,x_1)) = \Hom _{\precat (\mM )}^{x_0,x_1}(\pE , \pA )
\end{equation}
where the superscript on the right designates the subset of maps $\pE \rightarrt \pA$ sending
$\upsilon _0$ to $x_0$ and $\upsilon _1$ to $x_1$. 

Let $\overline{\pE}$ denote the image of the category with two isomorphic objects
under the map $\precat (\Sets )\rightarrt \precat (\mM )$. We think of $\overline{\pE}$ 
as containing
$\pE$ as a subcategory. Thus $\overline{\pE}$ again has objects $\upsilon _0,\upsilon _1$, but $\pE(x_0,\ldots , x_p)=\ast $
for any sequence of objects. One can also view it as 
the codiscrete precategory with two objects,
$\overline{\pE}=\codisc ([1])=\codisc (\{ \upsilon _0, \upsilon _1\} )$
in the notation of Section \ref{sec-precatexamples}. 

If $\pA\in \precat (\mM )$, a map $\overline{\pE}\rightarrt \pA$ is sure to correspond
to an internal
equivalence between the images of the two endpoints $\upsilon _0,\upsilon _1$.
Say that a map $\pE\rightarrt \pA$ which extends to $\overline{\pE}\rightarrt
\pA$, is {\em strongly invertible}. An important little observation is that the
identity morphisms (i.e. those given by the images of the degeneracies
$\ast = \pA (x_0)\rightarrt \pA(x_0,x_0)$) are strongly invertible. 

Unfortunately, given a general morphism from $x_0$ to $x_1$ in $\pA$,
the corresponding map $\pE\rightarrt \pA$ will not in general extend to 
$\overline{\pE}\rightarrt \pA$. That is to say, not all internal equivalences
will be strongly invertible. This is why we need to do some further work to
construct versal interval objects.

Assuming that $\pA$ satisfies the Segal condition and is levelwise fibrant,
suppose $x_0,x_1\in\Ob (\pA )$ and suppose $a:\ast \rightarrt \pA (x_0,x_1)$
is a morphism from $x_0$ to $x_1$. The condition of $a$ 
being an internal equivalence means that there
should be morphisms $b$ and $c$ from $x_1$ to $x_0$, such that
$ba$ is homotopic to the identity of $x_0$ and $ac$ is homotopic to the identity of $x_1$.
In turn, these homotopies can be represented by maps from interval objects in $\mM$ which we
shall denote by $N$ and $N'$ respectively.  
We will build up a big coproduct representing this collection of data. 

It turns out to be convenient to relax slightly the conditions that the homotopies
go between $ab$ and the identity (resp. $ca$ and the identity). Instead, we say that the
homotopies go between $ab$ or $ca$ and strongly invertible morphisms. In particular,
the source of $c$ could be an object $x'_1$ different from $x_1$ and the target of $b$
could be an object $x'_0$ different from $x_0$. 

This situation can be represented diagramatically by
$$
{\setlength{\unitlength}{.5mm}
\begin{picture}(60,100)
\put(51,95){\ensuremath{x'_1}}
\put(51,27){\ensuremath{x_1}}
\put(3,72){\ensuremath{x_0}}
\put(3,3){\ensuremath{x'_0}}

\put(50,64){\ensuremath{\Leftrightarrow}}
\put(50,56){\ensuremath{N'}}
\put(5,40){\ensuremath{\Leftrightarrow}}
\put(5,32){\ensuremath{N}}
\put(30,88){\ensuremath{c}}
\put(28,55){\ensuremath{a}}
\put(30,10){\ensuremath{b}}

\put(50,95){\vector(-2,-1){40}}
\put(10,70){\vector(1,-1){40}}
\put(50,27){\vector(-2,-1){40}}

\put(50,33){\vector(1,-2){0}}
\put(55,33){\vector(-1,-2){0}}
\put(55,93){\vector(-1,2){0}}

\put(10,10){\vector(-1,-2){0}}
\put(5,10){\vector(1,-2){0}}
\put(5,70){\vector(1,2){0}}

\qbezier(50,93)(35,63)(50,33)
\qbezier(8,69)(25,40)(10,10)

\linethickness{.4mm}
\qbezier(55,93)(70,63)(55,33)
\qbezier(5,70)(-10,40)(5,10)

\end{picture}
} .
$$
The goal in this section is to construct a precategory $\Xi (N|N')$ such that a map $\Xi (N|N')\rightarrt \pA$
is the same thing as such a diagram. Notice that the diagram may be divided into two
triangles which are independent except for the fact that they share the same edge
labeled $a$. Our $\Xi (N|N')$ will be a coproduct of two precategories
$\Xi (N)$ and $\Xi (N')$ along $\pE$, where each of the pieces represents
a triangular diagram.

So to start, suppose 
given an interval object $(N,i_0,i_1)$ with $N\in \mM$ and $i_0,i_1:\ast \rightarrt N$.
We will construct a precategory $\Xi (N)$ such that
a map from $\Xi (N)$ to $\pA$ is the same thing as a diagram of the form
$$
{\setlength{\unitlength}{.5mm}
\begin{picture}(60,80)
\put(54,38){\ensuremath{y_1}}
\put(3,73){\ensuremath{y_0}}
\put(3,3){\ensuremath{y_2}}

\put(5,40){\ensuremath{\Leftrightarrow}}
\put(5,32){\ensuremath{N}}

\put(10,70){\vector(3,-2){42}}
\put(52,35){\vector(-3,-2){40}}

\put(10,10){\vector(-1,-2){0}}
\put(5,10){\vector(1,-2){0}}
\put(5,70){\vector(1,2){0}}

\qbezier(8,69)(25,40)(10,10)

\linethickness{.4mm}
\qbezier(5,70)(-10,40)(5,10)

\end{picture}
} 
$$
in $\pA$. There are three pieces. The part involving $N$ is
a map to $\pA$ from 
an  $\mM$-enriched precategory of the form $\Upsilon (N)$ (see Section \ref{sec-upsilon} and Chapter \ref{freecat1}),
which comes with two maps $\Upsilon (i_0)$ and 
$\Upsilon (i_1)$ from $\pE$ to $\Upsilon (N)$. 
The commutative triangle corresponds to a map from a representable precategory of the form $h([2],\ast )$ to $\pA$. The strongly invertible morphism on the left corresponds to an extension
of one of the $\pE\rightarrt \Upsilon (N)\rightarrt \pA$ to 
a map 
$\overline{\pE}\rightarrt \pA$. 

Motivated by this picture,
define the $\mM$-enriched precategory $\Xi (N)$ to be the coproduct of three terms
corresponding to these three pieces:
$$
\Xi  (N):= \overline{\pE}\cup ^{\Upsilon (i_0)(\pE )} \Upsilon (N) 
\cup ^{\Upsilon (i_1)(\pE )} h([2],\ast ).
$$
The map at the end of the coproduct notation is $\pE =h([1],\ast )\rightarrt h([2],\ast )$
corresponding to the edge $[1]\rightarrt \nocom [2]$ sending $0$ to $0$ and $1$ to $2$. 
The objects of $\Xi (N)$ will be denoted $\xi _0,\xi _1,\xi _2$. These correspond to the three objects of $h([2],\ast )$.  In case of a map $\Xi (N)\rightarrt \pA$ corresponding
to a diagram as above, the images of $\xi _i$ are the objects labeled $y_i$ above. 
Thus the two objects $\upsilon _0,\upsilon _1$ of both copies of $\pE$ as well as $\Upsilon (N)$ correspond to $\xi _0$ and $\xi _2$ respectively.

\begin{lemma}
\label{mapfromXiN}
Suppose $\pA\in \precat (\mM )$. Then a map $\Xi (N)\rightarrt \pA$ corresponds to giving three objects 
$x_0,x_1,x_2\in \Ob (\pA )$, to giving an element $t:\ast \rightarrt \pA(x_0,x_1,x_2)$, to giving a map 
$b:N\rightarrt A(x_0,x_2)$ and to giving a map $g:\overline{\pE}\rightarrt \pA$ such that 
$b\circ \Upsilon (i_1)= \partial _{02}(t)$, and $b\circ \Upsilon (i_0)=g(e_{01})$ where $e_{01}:\ast \rightarrt \overline{\pE}(\upsilon _0,\upsilon _1)$
is the unique map. 
\end{lemma} 
\begin{proof}
This comes from the coproduct description for $\Xi (N)$. 
\end{proof}

We think of $t:\ast \rightarrt \pA(x_0,x_1,x_2)$ as corresponding to a commutative triangle 
with maps $\partial _{01}(t)$ and $\partial _{12}(t)$ whose ``composition'' is 
$$
\partial _{12}(t)\circ \partial _{01}(t)=\partial _{02}(t).
$$
Then $N$ can be a homotopy from $\partial _{02}(t)$ to the map $g(e_{01})$ (in our application $N$ will be contractible).
Then the extension of this map to $g$ defined on $\overline{\pE}$ says that $g(e_{01})$ is strictly invertible. 
So, roughly speaking when we look at a map $\Xi (N)\rightarrt \pA$ we are looking at two morphisms 
whose composition $\partial _{12}(t)\circ \partial _{01}(t)$ is equivalent to an invertible map. 

The two different maps in question correspond to maps $\zeta _{01}, 
\zeta _{12}:\pE \rightarrt \Xi (N)$ with $\zeta _{01}$ corresponding to 
$\partial _{01}(t)$ and $\zeta _{12}$ to $\partial _{12}(t)$. 

The construction $\Xi $ also works for the other half of $\Xi (N|N')$. We distinguish
the two interval objects which are used here, for clarity of notation. Obviously one
could choose the same on both sides. 

Given two interval objects $N$ and $N'$, we can form 
$$
\Xi (N|N'):= \Xi (N)\cup ^{\pE}\Xi (N')
$$
where the map $\pE \rightarrt \Xi (N)$ is $\zeta _{01}$ and the map $\pE\rightarrt \Xi (N')$ is $\zeta _{12}$.
These become the same map denoted $\eta : \pE\rightarrt \Xi (N|N')$. 
Denote the four objects of $\Xi (N|N')$ by $\xi _{|0}$, $\xi _{0|1}$, $\xi _{1|2}$ and $\xi _{2|}$, these corresponding with
the objects of $\Xi (N)$ or $\Xi (N')$ by saying that $\xi _{i|j}$ corresponds
to $\xi _i$ in the left piece $\Xi (N)$ and to $\xi _j$ in the right 
piece $\Xi (N')$ to give the following picture of $\Xi (N|N')$:
$$
{\setlength{\unitlength}{.5mm}
\begin{picture}(60,100)
\put(51,95){\ensuremath{\xi _{|0}}}
\put(51,27){\ensuremath{\xi _{1|2}}}
\put(3,72){\ensuremath{\xi _{0|1}}}
\put(3,3){\ensuremath{\xi _{2|}}}

\put(43,61){\ensuremath{\Upsilon (N')}}
\put(-2,39){\ensuremath{\Upsilon (N)}}

\put(64,48){\ensuremath{\overline{\pE}}}
\put(-8,50){\ensuremath{\overline{\pE}}}

\put(17,71){\ensuremath{h([2],\ast )}}
\put(19,29){\ensuremath{h([2],\ast )}}
\put(29,55){\ensuremath{\eta}}

\put(50,95){\vector(-2,-1){40}}
\put(10,70){\vector(1,-1){40}}
\put(50,27){\vector(-2,-1){40}}

\put(50,33){\vector(1,-2){0}}
\put(55,33){\vector(-1,-2){0}}
\put(55,93){\vector(-1,2){0}}

\put(10,10){\vector(-1,-2){0}}
\put(5,10){\vector(1,-2){0}}
\put(5,70){\vector(1,2){0}}

\qbezier(50,93)(35,63)(50,33)
\qbezier(8,69)(25,40)(10,10)

\linethickness{.4mm}
\qbezier(55,93)(70,63)(55,33)
\qbezier(5,70)(-10,40)(5,10)

\end{picture}
} 
$$
displaying $\eta$ as a morphism from 
$\xi _{0|1}$ to $\xi _{1|2}$ i.e. an element of the set
$\Mor ^1_{\Xi (N|N')}( \xi _{0|1},\xi _{1|2})$
defined in \eqref{mordef}. 

\begin{lemma}
\label{mapfromXiNN}
Suppose $(N,i_0,i_1)$ and $(N',i'_0,i'_1)$ are interval objects of $\mM$. If $\pA$ is
an $\mM$-enriched precategory, then a map $\Xi (N|N')\rightarrt \pA$
corresponds to the data of a morphism $(x_0,x_1,a)$ in $\pA$,
of two other objects $x'_0$ and $x'_1$,
together with commutative triangles 
$$
s:\ast \rightarrt \pA (x_0,x_1,x'_0),
\;\;\;
t:\ast \rightarrt \pA (x'_1,x_0,x_1),
$$
with maps $h:\Upsilon (N)\rightarrt \pA$ and $h':\Upsilon (N')\rightarrt \pA$
and two maps $u,v:\overline{\pE}\rightarrt \pA$ such that various maps $\pE\rightarrt \pA$
induced by these data coincide (see the diagram in the proof below).  
\end{lemma}
\begin{proof}
This comes from the coproduct description of $\Xi (N|N')$ and the corresponding properties for 
$\Xi (N)$ and $\Xi (N')$. 
The objects $x_0$ and $x_1$ are the images of 
$\xi _{0|1}$ and $\xi _{1|2}$ while $x'_0$ is the image of
$\xi _{2|}$ and $x'_1$ is the image of $\xi _{|0}$.

The maps which are supposed to coincide may be read off
from the diagram
$$
{\setlength{\unitlength}{.5mm}
\begin{picture}(60,100)
\put(51,95){\ensuremath{x'_1}}
\put(51,27){\ensuremath{x_1}}
\put(3,72){\ensuremath{x_0}}
\put(3,3){\ensuremath{x'_0}}

\put(50,64){\ensuremath{\Leftrightarrow}}
\put(50,56){\ensuremath{h'}}
\put(5,40){\ensuremath{\Leftrightarrow}}
\put(5,32){\ensuremath{h}}
\put(25,71){\ensuremath{t}}
\put(29,71){\ensuremath{\circlearrowleft}}
\put(28,55){\ensuremath{a}}
\put(24,31){\ensuremath{s}}
\put(29,31){\ensuremath{\circlearrowright}}
\put(-10,41){\ensuremath{u}}
\put(68,65){\ensuremath{v}}

\put(50,95){\vector(-2,-1){40}}
\put(10,70){\vector(1,-1){40}}
\put(50,27){\vector(-2,-1){40}}

\put(50,33){\vector(1,-2){0}}
\put(55,33){\vector(-1,-2){0}}
\put(55,93){\vector(-1,2){0}}

\put(10,10){\vector(-1,-2){0}}
\put(5,10){\vector(1,-2){0}}
\put(5,70){\vector(1,2){0}}

\qbezier(50,93)(35,63)(50,33)
\qbezier(8,69)(25,40)(10,10)

\linethickness{.4mm}
\qbezier(55,93)(70,63)(55,33)
\qbezier(5,70)(-10,40)(5,10)

\end{picture}
} 
$$
\label{fullpicture}
in which the $2$-cells represent maps from $\Upsilon (N)$ or $\Upsilon (N')$,
the thick lines represent maps from $\overline{\pE}$, and the triangles
represent maps from $h([2],\ast )$. 
For example the boundary $\partial _{02}\circ t:\ast \rightarrt \pA (x'_1,x_1)$
corresponds to a map $\pE\rightarrt \pA$ sending $\upsilon _0$ to $x'_1$ and $\upsilon _1$ to
$x_1$. This should be the same as the map 
$$
\pE = \Upsilon (\ast )\rightarrt^{\Upsilon (i'_0)}\Upsilon (N')\rightarrt^{h'}\pA .
$$
The other identifications are similar. 
\end{proof}

\begin{lemma}
\label{xireedycof}
For any two interval objects $N$ and $N'$, the $\mM$-precategory $\Xi (N|N')$
is  Reedy cofibrant, and indeed the inclusion 
$$
\disc \{\xi _{0|1}, \xi _{1|2} \} \rightarrt \Xi (N|N')
$$
is Reedy cofibrant hence injectively cofibrant. 
\end{lemma}
\begin{proof}
Use Corollaries \ref{RisReedy} and \ref{RfaceReedy}, and Lemma \ref{upsilonreedy}.
\end{proof}

However, $\Xi (N|N')$ is not projectively cofibrant, because 
the inclusions of edges $\pE \rightarrow h([2],N)$ are Reedy but not projective cofibrations.
This issue will be addressed further in the comments after Remark \ref{XiNNforproj} below. 

Record here what happens when we change the intervals used in the construction. 

\begin{lemma}
\label{changeofinterval}
Suppose $f:N\rightarrt P$  is a  morphism between interval objects $(N,i_0,i_1)$
and $(P,j_0,j_1)$, that is $f\circ i_0=j_0$, $f\circ i_1=j_1$. Suppose similarly
$f':N'\rightarrt P'$ is a  morphism between interval objects $(N',i'_0,i'_1)$
and $(P',j'_0,j'_1)$. Then these induce a global weak equivalence 
$$
\Xi (N| N')\rightarrt \Xi (P | P').
$$
\end{lemma}
\begin{proof}
It is a levelwise weak equivalnce,
being a pushout of maps which are levelwise a weak equivalences.
\end{proof}

\section{The versality property}

From the universal property of Lemma \ref{mapfromXiNN}, we obtain the versality property
of $\Xi (N|N')$. 

\begin{theorem}
\label{xiinterval}
Suppose $\pA\in \precat (X,\mM )$ satisfies the Segal condition 
and is fibrant in the Reedy diagram model structure $\diag _{\rm Reedy}(\Delta^o_X/X,\mM )$. 
Suppose that $x,y\in X=\Ob (\pA )$ and $a:\ast \rightarrt \pA (x,y)$ is an element of
$\Mor ^1_{\pA }(x,y)$. Suppose that $a$ is an inner equivalence, in other
words the image of $a$ in 
the truncated category $\tau _{\leq 1}(\pA )\in \Cat$, is invertible. Then for any interval objects $N$ and $N'$ in $\mM$ there exists
a morphism $\Xi (N|N')\rightarrt \pA$ sending $\xi _{0|1}$ and $\xi _{2|}$ to $x$ and $\xi _{1|2}$ and $\xi _{|0}$ to $y$,
sending the tautological morphism $\eta$ to $a$, and sending the
two copies of $\overline{\pE}$ to the identities of $x$ and $y$ respectively. 
\end{theorem}
\begin{proof}
Since $\pA$ satisfies the Segal condition, the truncation $\tau _{\leq 1}(\pA )$
may be defined using $\pA $ itself, that is to say that the truncation is the $1$-category
with $x+\Ob (\pA )$ as set of objects, and whose nerve relative to this set is the functor
$$
\Delta ^o_{x}\rightarrt \Sets , \;\;\; (x_0,\ldots , x_n)\mapsto \tau _{\leq 0}\pA (x_0,\ldots , x_n).
$$
The Reedy fibrant condition for $\pA$ implies that it is levelwise fibrant,
which means that for any sequence
$(x_0,\ldots , x_n)\in \Delta ^o_{X}$ the image $\pA(x_0,\ldots , x_n)$ is a fibrant object of $\mM$, so 
$$
\tau _{\leq 0}\pA (x_0,\ldots , x_n) = \Hom _{\mM }(\ast , \pA (x_0,\ldots , x_n))/\sim 
$$
where $\sim$ is the relation of homotopy occuring in Lemma \ref{tauinterval}.

The fact that $a$ maps to an isomorphism in $\tau _{\leq 1}\pA$ therefore
means that there is an inverse
$b\in \tau _{\leq 0}\pA (y,x)$; and by the levelwise fibrant condition it can be
represented by $b:\ast \rightarrt \pA (y,x)$. 

By the Segal condition the morphism 
$$
\pA (x,y,x)\rightarrt \pA (x,y)\times \pA (y,x)
$$
is a weak equivalence. On the other hand, the Reedy fibrant condition in the diagram
category means that the matching map at $(x,y,z)$ is a fibration in $\mM$, which
in turn implies that the Segal map above is a fibration. Hence it is a trivial fibration,
in particular the element $(a,b): \ast \rightarrt  \pA (x,y)\times \pA (y,x)$
lifts to a map $\ast \rightarrt \pA (x,y,x)$. This gives a diagram 
$$
s:h([2],\ast )\rightarrt \pA
$$
representing ``the composition $b\circ a$'', fitting into the lower left triangle in
the picture on page \pageref{fullpicture}. The image of $s$ 
in the nerve of $\tau _{\leq 1}\pA$
is the commutative triangle for the composition of the images of $b$ and $a$.
We chose $b$ as representing an inverse to $a$ in the truncated category, so the $02$
edge of $s$ is homotopic to the identity of $x$; by Lemma \ref{tauinterval}
there exists a map $N\rightarrow \pA (x,x)$ representing this homotopy, or by adjunction
$$
h:\Upsilon (N)\rightarrt \pA .
$$
This gives the lower or $\Xi (N)$ part of the required diagram. A similar discussion
using the fact that $a\circ b$ is homotopic to the identity of $y$, gives the upper
or $\Xi (N')$ part, and these glue together to give the required map 
$\Xi (N|N')\rightarrt \pA$.
\end{proof}

\begin{remark}
\label{XiNNforproj}
Let $\Xi ^{\rm proj}(N|N')$ denote a cofibrant replacement for $\Xi (N|N')$
in the projective diagram category
$\diag _{\rm Reedy}(\Delta^o_X/X,\mM )$. Then it has the same versality property
with respect to any $\pA$ which is levelwise fibrant and satisfies the Segal 
condition. 
\end{remark}

The projectively cofibrant version $\Xi ^{\rm proj}(N|N')$
could be  constructed explicitly 
by inserting objects of the form $\Upsilon (L)$  and $\Upsilon (L')$
in between $\zeta$ and the
two triangles, for intervals $L$ and $L'$, according to the picture
$$
{\setlength{\unitlength}{.5mm}
\begin{picture}(60,100)
\put(51,95){\ensuremath{\xi _{|0}}}
\put(51,27){\ensuremath{\xi _{1|2}}}
\put(3,72){\ensuremath{\xi _{0|1}}}
\put(3,3){\ensuremath{\xi _{2|}}}

\put(49,64){\ensuremath{N'}}
\put(4,40){\ensuremath{N}}

\put(25,44){\ensuremath{L}}
\put(30,50){\ensuremath{L'}}

\put(50,95){\vector(-2,-1){40}}
\put(10,70){\vector(1,-1){40}}
\put(50,27){\vector(-2,-1){40}}

\put(53,33){\vector(1,-2){0}}
\put(55,33){\vector(-1,-2){0}}
\put(55,93){\vector(-1,2){0}}

\put(7,10){\vector(-1,-2){0}}
\put(5,10){\vector(1,-2){0}}
\put(5,70){\vector(1,2){0}}

\qbezier(10,70)(42,62)(50,30)
\qbezier(10,70)(18,38)(50,30)

\qbezier(53,93)(42,63)(53,33)
\qbezier(7,69)(18,40)(7,10)

\linethickness{.4mm}
\qbezier(55,93)(63,63)(55,33)
\qbezier(5,70)(-3,40)(5,10)

\end{picture} .
} 
$$
This corresponds to the step in the proof of the theorem where we used
the Reedy fibrant property to lift $(a,b)$ to an element $s:\ast \rightarrt \pA (x,y,x)$.
If $\pA$ is assumed only fibrant in the projective structure (i.e. levelwise fibrant)
then $(a,b)$ only lifts up to a homotopy given by $L\rightarrt \pA (x,y)\times \pA (y,x)$.
The second term may be neglected since we don't care about the bottom arrow
of the big diagram, and the piece $L\rightarrt \pA (x,y)$ corresponds to a map $\Upsilon (L)\rightarrt \pA$. We presented the Reedy version in our main discussion above 
because the diagrams are easier to picture.

Let 
$\widetilde{\Xi}(N|N')\subset \Seg (\Xi (N|N'))$ be
the full subcategory containing only the two objects
$\xi _{|0}$ and $\xi _{2|}$. In the situation of the theorem,
we get by restriction plus functoriality of $\Seg$ a map 
$$
\widetilde{\Xi}(N|N')\rightarrt \Seg (\pA ).
$$
Similarly if $\widetilde{\Xi}^{\rm proj}(N|N')\subset \Seg (\Xi ^{\rm proj}(N|N'))$ is
the full subcategory containing only 
$\xi _{|0}$ and $\xi _{2|}$, then in the situation of the remark
we get a map as stated in the following corollary.

\begin{corollary}
Suppose $\pA$ is levelwise fibrant and satisfies the Segal conditions,
and suppose $x,y\in \Ob (\pA )$ are two internally equivalent objects. Then
there is a map 
$$
\widetilde{\Xi}^{\rm proj}(N|N')\rightarrt \Seg (\pA )
$$
sending the two objects of $\widetilde{\Xi}^{\rm proj}(N|N')$ to $x$ and $y$
respectively.
\end{corollary}
\begin{proof}
As in the above remark we get a map 
${\Xi}^{\rm proj}(N|N')\rightarrt \pA $, hence by functoriality of $\Seg$ 
and composition with the inclusion, 
$$
\widetilde{\Xi}^{\rm proj}(N|N')\subset 
\Seg (\Xi ^{\rm proj}(N|N'))\rightarrt \Seg (\pA )
$$
as required.
\end{proof}

It remains to be seen that $\Xi (N|N'))$ and thus $\widetilde{\Xi}(N|N')$
are contractible. 

\section{Contractibility of intervals for $\mK$-precategories}

Given the above construction, the main problem is to prove that
$\Xi (N|N')$  is contractible. In this section we do that 
for enrichment over the Kan-Quillen model category $\mK$ of simplicial sets. 

\begin{theorem}
\label{xicontractibleK}
Suppose $N,i_0,i_1$ and $N',i'_0,i'_1$ are two interval objects in the Kan-Quillen model category of
simplicial sets $\mK$. Then 
$\Xi (N|N')$ is contractible in $\precat (\mK )$, that is $\Xi (N|N')\rightarrt \ast$ is a global weak equivalence. 
We have a map $\Xi (N|N')\rightarrt \overline{\pE }\times \overline{\pE }$ which is a global weak equivalence and an isomorphism
on the sets of objects. In particular the map 
$$
\Seg (\Xi (N|N'))\rightarrt \overline{\pE }\times \overline{\pE }
$$
induces an objectwise weak equivalence, which is to say that 
$$
\Seg (\Xi (N|N'))(x_0,\ldots , x_p) \mbox{ is contractible in }\mM
$$
for any sequence of objects $x_0,\ldots , x_p\in \{ \xi _{|0}, \xi _{0|1}, \xi _{1|2}, \xi _{2|}\}$. 
\end{theorem}
\begin{proof}
This was treated in the last chapter of 
\cite{Pelissier} and our 
present version is only slightly different in that we have expanded somewhat $\Xi$ as something with $4$ objects. Our present picture is perhaps closer to Drinfeld's intervals
for DG-categories \cite{DrinfeldIntervalDG}.

Elements of $\precat (\mK )$ may be considered as certain kinds of bisimplicial sets (see 
Section \ref{sec-interpretations} and Chapter \ref{secat1}),
and this commutes with coproducts. Similarly the diagonal realization from bisimplicial sets to simplicial sets
commutes with coproducts and takes Reedy or injective cofibrations\footnote{The
Reedy and injective cofibrations are the same in  $\precat (\mK )$
by Proposition \ref{reedyinjective} as was pointed out in Corollary \ref{precatKstrucs}.}
in $\precat (\mK )$
to cofibrations in $\mK$ (which are just the monomorphisms).
Call the composition of these two operations $|\cdot | : \precat (\mK )\rightarrt \mK$.
Note that $|\pE |$, $|\overline{\pE } |$, and $|h([2],\ast )|$ are contractible simplicial sets,
and if $N$ is an interval object in $\mK$ then $|\Upsilon (N)|$ is contractible. Thus, 
$|\Xi (N)|$ is a successive cofibrant coproduct of contractible objects over contractible objects, so it is contractible.
Similarly the coproduct of two of these over the contractible $|\pE |$ (mapping into both sides by cofibrations)
is contractible, so $|\Xi (N|N')|$ is
contractible. In general a map $A\rightarrt \Seg (A)$ induces a weak equivalence of simplicial sets $|A|\rightarrt^{\sim}
|\Seg (A)|$. Thus in our case, $|\Seg (\Xi (N|N')|$ is contractible. On the other hand, all of the $1$-morphisms
in $\Xi (N|N')$ go to invertible morphisms in $\Seg (\Xi (N|N')$, in effect
the main middle morphism $\eta$ has by construction a left and a right inverse up to equivalence; so it goes
to an equivalence, and its inverses go to equivalences too. Thus $\Seg (\Xi (N|N')$ is a Segal groupoid. 
Now, a Segal groupoid whose realization is contractible, is contractible (see Proposition \ref{segalpi}). Thus $\Seg (\Xi (N|N')$ is contractible, which proves the theorem in the case of $\mK$.
\end{proof}

\section{Construction of a left Quillen functor $\mK\rightarrow \mM$}

In order to transfer the above contractibility result for $\Xi (N|N')$ in
the $\mK$-enriched case, to the general case, we
explain in this section the essentially well-known construction of
a left Quillen functor $\mK\rightarrt \mM$. 
In Hovey \cite{Hovey} was 
explained the intuition that every monoidal model category is a module over $\mK$,
and even without the monoidal structure there is a left Quillen functor from $\mK$ into
$\mM$.
The construction is based on a choice
of contractible cosimplicial object in $\mM$, or more precisely
a {\em cosimplicial resolution} in the sense of Hirschhorn \cite{Hirschhorn}.
That means  a functor $R: \Delta \rightarrt \mM$ which 
is cofibrant in the Reedy model structure  $\diag _{\rm Reedy}(\Delta , \mM )$.

Recall that an object $A\in \mM$ is {\em contractible} if the unique morphism $A\rightarrt \ast$ is
a weak equivalence. We say that a cosimplicial object $\resol :\Delta \rightarrt \mM$ is
{\em levelwise contractible} if $\resol ([n])$ is contractible for each object $[n]\in \ob (\Delta )$. 

\begin{lemma}
There exists a choice of Reedy-cofibrant levelwise contractible cosimplicial object $\resol :\Delta \rightarrt \mM$.
\end{lemma}
\begin{proof}
See \cite{Hirschhorn}.
\end{proof}

We fix one such choice, from now on. The objects $\resol ([n])$ may be thought of as the ``standard $n$-simplices'' in $\mM$. 
If $A\in \mM$, define $\resol ^{\ast}(A): \Delta ^o\rightarrt \Sets$ to be the functor 
$$
\resol ^{\ast}(A): [n]\mapsto \Hom _{\mM}(\resol ([n]), A).
$$

\begin{theorem}
\label{resolthm}
If $\resol $ is a Reedy-cofibrant levelwise contractible cosimplicial object, then $\resol ^{\ast}$
is a right Quillen functor from $\mM$ to $\mK = \diag (\Delta ^o, \Sets )$. Its left adjoint 
$$
\resol _! : \mK \rightarrt \mM
$$
is a left Quillen functor given by the usual formula for the topological realization of a simplicial set,
but with the standard $n$-simplex replaced by  $\resol ([n])\in \mM$. 
\end{theorem}
\begin{proof}
See \cite{Hirschhorn}.
\end{proof}

\begin{corollary}
The realization functor induces, for every object set $X$, a functor 
$$
\precat (X; \resol _!): \precat _c(X, \mK )\rightarrt \precat _c(X,\mM ).
$$
This is a left Quillen functor for $c$ denoting either the projective, injective
or Reedy model structures on $\mK$ and $\mM$-enriched
precategories over $X$. It is compatible with change of set $X$, and gives a functor 
$$
\precat (\resol _!): \precat (\mK )\rightarrt \precat (\mM )
$$
from the  Segal precategories to the  $\mM$-enriched precategories. 
\end{corollary}
\begin{proof}
Combine Theorem \ref{resolthm}, with the discussion of Proposition \ref{precatF}. 
\end{proof}

The corresponding right Quillen functor $\precat (\resol ^{\ast}): \precat (\mM )\rightarrt \precat (\mK )$
should be applied to $\pA \in \precat (\mM )$ only after taking a fibrant replacement $\pA \rightarrt \pA '$.
Define $\Int _{\resol} (\pA ):= \precat (\resol ^{\ast})(\pA ')$. We call this the {\em $\resol $-interior} of $\pA $, since it
is obtained by looking at maps from the standard simplices $\resol ([n])$ into $\pA (x_0,\ldots , x_n)$
so it measures $\pA $ ``from the inside''. This construction is compatible with truncation:

\begin{lemma}
For $\pA \in \precat (\mM )$ we have an isomorphism of categories $\tau _{\leq 1}(\pA ) \cong \tau _{\leq 1}(\Int _{\resol }(\pA ))$. 
\end{lemma}
\begin{proof}
This follows from the definition of $\tau _{\leq 1}$ in Section \ref{sec-globalwe}. 
\end{proof}

By its construction as a colimit, $\precat (\resol _!)$ preserves coproducts, preserves constructions $\Upsilon$ and $h$, preserves the
various notions of cofibrancy. Since $\mM$ is cartesian, Proposition \ref{precatFpreserveswe}
says that
$\precat (\resol _!)$  takes weak equivalences in $\precat (X,\mK )$ to
weak equivalenes in $\precat (X, \mM )$. Furthermore since it preserves truncations, 
$\precat (\resol _!)$
preserves global weak equivalences, and preserves the
notion of interval objects.

\section{Contractibility in general}

We can now use the functor $\precat (X; \resol _!)$ to transfer the 
the contractibility result for $\mK$-enriched precategories, to 
$\precat (\mM )$. 
This yields the first main theorem of the present chapter saying that $\Xi (N|N')$ is a good interval object in $\precat (\mM )$. 
This was the 
last step missing in Pelissier's \cite{Pelissier} correction of \cite{svk}, 
but which is actually straightforward from a Quillen-functorial point of view. 

The contractibility statement is made before we have completely finished the construction of the model structure, although it is the penultimate step.
Some care is still therefore necessary in using only the parts of the model structure
which are already known. 

\begin{theorem}
\label{xicontractible}
Suppose $N,i_0,i_1$ and $N',i'_0,i'_1$ are two interval objects. Then 
$\Xi (N|N')$ is contractible in $\precat (\mM )$, that is $\Xi (N|N')\rightarrt \ast$ is a global weak equivalence. 
We have a map $\Xi (N|N')\rightarrt \overline{\pE}\times \overline{\pE}$ which is a global weak equivalence and an isomorphism
on the sets of objects. In particular the map 
$$
\Seg (\Xi (N|N'))\rightarrt \overline{\pE }\times \overline{\pE }
$$
induces an objectwise weak equivalence, which is to say that 
$$
\Seg (\Xi (N|N'))(x_0,\ldots , x_p) \mbox{ is contractible in }\mM
$$
for any sequence of objects $x_0,\ldots , x_p\in \{ \xi _{|0}, \xi _{0|1}, \xi _{1|2}, \xi _{2|}\}$. 
\end{theorem}
\begin{proof}
First notice that the statement of the theorem is independent of the choice of interval object: if $N\rightarrt P$ and
$N'\rightarrt P'$ are maps of interval objects then the statement of the theorem for $(N,N')$ is equivalent to the
statement for $(P,P')$. See Lemma \ref{changeofinterval}. 
In particular it suffices to prove the statement for one pair of  intervals. 

Theorem \ref{xicontractibleK} gives the same statement for precategories enriched over the Kan-Quillen model category $\mK$ of simplicial sets. 
Then choose a left Quillen functor $\resol _!:\mK \rightarrt \mM$. This gives a functor $\precat (\resol _!):\precat (\mK )\rightarrt 
\precat (\mM )$ which preserves coproducts. 
Suppose $(B,i_0,i_1)$ is an interval object in $\mK$. Then $\resol _!(B)$ is an interval object in $\mM$ and 
$$
\precat (\resol _!)(\Xi (B|B) = \Xi (\resol _!(B)|\resol _!(B)).
$$
Now since $\precat (\resol _!)$ preserves global weak equivalences, we obtain the statement of the theorem 
for the pair of interval objects $\resol _!(B)|\resol _!(B)$. By the invariance discussed in the first paragraph of the proof,
this implies the statement of the theorem for all $N,N'$. 
\end{proof}

Recall that $\widetilde{\Xi}(N|N')\subset \Seg (\Xi (N|N'))$ is 
the full subcategory containing only the two objects
$\xi _{|0}$ and $\xi _{2|}$. Since all objects of $\Seg (\Xi (N|N'))$ are equivalent, the inclusion
$$
\widetilde{\Xi}(N|N')\hookrightarrow \Seg (\Xi (N|N'))
$$
is a global weak equivalence (it is by definition fully faithful and both sides satisfy the Segal conditions).
It follows from Theorem \ref{xicontractible} and the 3 for 2 property of global weak equivalences, that the functor 
$$
p_{N,N'}:\widetilde{\Xi}(N|N')\rightarrt \overline{\pE}
$$
is a global weak equivalence; furthermore this induces isomorphisms on the sets of objects (there are exactly two
objects on each side), and both sides satisfy the Segal conditions, so $p_{N,N'}$ is an objectwise weak equivalence of 
diagrams.

\section{Pushout of trivial cofibrations}

These interval objects allow us to analyze pushouts along trivial cofibrations which are not isomorphisms on objects.
In this discussion, we use injective cofibrations since this encompasses the Reedy and projective cofibrations too. 

We start by considering the pushout along the standard interval $\overline{\pE }$.

\begin{lemma}
\label{pushoutE}
Suppose $\pA \in \precat (\mM )$, and $y\in \Ob (\pA )$. Then the pushout
morphism 
$$
a:\pA \rightarrt \pA \cup ^{\{ y\} } \overline{\pE }
$$
obtained by identifying $\upsilon _0$ and $y$, is a global weak equivalence.
\end{lemma}
\begin{proof}
By Corollary \ref{prodwe} applied to the identity of $\pA$ and the map $p:\overline{\pE}\rightarrt \ast$,
the map 
$$
1_{\pA}\times p: \pA \times \overline{\pE}\rightarrt \pA
$$
is a global weak equivalence. Let $i_0,i_1:\ast \rightarrt \overline{\pE}$ be the two inclusions of objects
$\upsilon _0$ and $\upsilon _1$. The two maps 
$$
1_{\pA}\times i_0,1_{\pA}\times i_1: \pA \rightarrt \pA \times \overline{\pE}
$$
are global weak equivalences, as can be seen by composing with $1_{\pA}\times p$ and using 3 for 2.

Now, consider the morphism $g:\overline{\pE}\times \overline{\pE}\rightarrt \overline{\pE}\times \overline{\pE}$
equal to the identity on $\overline{\pE}\times \{ \upsilon _0\}$ and sending $\overline{\pE}\times \{ \upsilon _1\}$
to the single object $(\upsilon _0,\upsilon _1)$. Set 
$$
\pB := \pA \cup ^{\{ y\} } \overline{\pE }.
$$
Let $q:\pB \rightarrt \pA$ denote the projection obtained by sending all of $\overline{\pE }$ to the single object $y\in \Ob (\pA ^)$.
Then 
$$
\pB \times \overline{\pE } = \left( \pA \times \overline{\pE } \right) \cup ^{\{ \upsilon _0\} \times  \overline{\pE }}
\left( \overline{\pE} \times \overline{\pE } \right) .
$$
Apply the map $g$ to the second factor of this pushout, to obtain a map 
$$
f:\pB \times \overline{\pE }\rightarrt \pB \times \overline{\pE }
$$
such that $f$ restricts to the identity on $\pB\times \{ \upsilon _0\}$, while $f|_{\pB \times \{ \upsilon _1\}}$ is the projection
$q:\pB \rightarrt \pA$. By the first paragraph of the proof, the maps induced by $f$,
$$
\pB  \times \{ \upsilon _0 \} \rightarrt \pB  \times \overline{\pE }
$$
and
$$
\pB  \times \{ \upsilon _1 \} \rightarrt \pB  \times \overline{\pE }
$$
are global weak equivalences. The composition of $1_{\pB}\times p : \pB \times \overline{\pE }\rightarrt \pB$
with the morphism $f$ considered above, is a morphism 
$$
(1_{\pB}\times p)\circ f:  \pB \times \overline{\pE }\rightarrt \pB
$$
such that the composition $(1_{\pB}\times p)\circ f \circ (1_{\pB}\times i_0)$ is the identity of $\pB$, and the
composition $(1_{\pB}\times p)\circ f \circ (1_{\pB}\times i_0)$ is the composition
$$
\pB \rightarrt^{q}\pA \rightarrt^{a} \pB .
$$
The facts that $(1_{\pB}\times p)\circ f \circ (1_{\pB}\times i_0)$ is the identity of $\pB$, and that 
$(1_{\pB}\times i_0)$ is a global weak equivalence, imply by 3 for 2 that $(1_{\pB}\times p)\circ f$
is a global weak equivalence. But then, composing with the global weak equivalence $(1_{\pB}\times i_1)$
we see that $(1_{\pB}\times p)\circ f \circ (1_{\pB}\times i_0)$ is a global weak equivalence, in other words
that the composition $aq$
is a global weak equivalence. In the other direction, the composition
$$
\pA \rightarrt^{a} \pB \rightarrt^{q}\pA
$$
is the identity of $\pA$. Thus, we conclude from 
the last sentence of Theorem \ref{globalretract32}
that both $q$ and the inclusion $\pA \rightarrt \pB$ are global weak equivalences.
This last statement is what we are supposed to prove. 
\end{proof}

\begin{corollary}
\label{pushoutgeninterval}
Suppose $\pB$ is an $\mM$-enriched precategory with two objects $b_0,b_1$. Suppose $\pB$  satisfies the Segal conditions, and is contractible, that is
the map $\pB\rightarrt \ast$ is a global weak equivalence. Then for any $\pA \in \precat (\mM )$ and any object $y\in \Ob (\pA )$
the map
$$
\pA \rightarrt \pA \cup ^{\{ b_0\}}\pB 
$$
obtained by identifying $b_0$ to $y$, is a global weak equivalence.
\end{corollary}
\begin{proof}
There is a unique map $f:\pB \rightarrt \overline{\pE}$ sending $b_0$ to $\upsilon _0$ and $b_1$ to $\upsilon _1$.
This map is a global weak equivalence, as seen by applying 3 for 2 to the composition
$$
\pB \rightarrt \overline{\pE}\rightarrt \ast .
$$
But since it induces an isomorphism on objects, and both sides satisfy the Segal conditions, it is an objectwise
weak equivalence of diagrams. 
Applying $f$ to the second piece of the given pushout, we get a map 
$$
g:\pA \cup ^{\{ b_0\}}\pB \rightarrt \pA \cup ^{\{ y\} } \overline{\pE }
$$
to the pushout considered in the previous corollary. However, $g$ is an objectwise weak equivalence of diagrams,
so it is a global weak equivalence. Note that $g$ commutes with the maps from $\pA$, so by the previous corollary and
3 for 2 we conclude that the map of the present statement is a global weak equivalence. 
\end{proof}

Suppose $f:\pA \rightarrt \pB$ is an injective trivial cofibration, and suppose $\pB$ is levelwise fibrant and satisfies the
Segal conditions. Let $Z:= \Ob (\pB ) - f(\Ob (\pA ))$. For each $z\in Z$  choose $e(z)\in \Ob (A)$ and 
$$
a(z)\in \pB (f(e(z)), z)
$$
such that the image of $a(z)$ is invertible in the truncated category. This is possible by the definition 
of essential surjectivity of $\pA \rightarrt \pB$. 

Applying Theorem \ref{xiinterval}
There exist collections of interval objects 
$N_z,N'_z$ indexed by $z\in Z$, and functors $t_i:\Xi (N_z| N'_x)\rightarrt \pB$ sending $\xi _{|0}$ to $e(z)$,
$\xi _{2|}$ to $z$, and sending the tautological morphism $\eta$ to $a(z)$. By functoriality of the construction $\Seg$
we get
$$
\Seg (\Xi (N_z| N'_x))\rightarrt \Seg (\pB ),
$$
and restricting this gives $\tilde{t}_i:\widetilde{\Xi}(N_z|N_z')\rightarrt \Seg (\pB )$. Now, $\tilde{t_i}$ sends
the first object to $e(z)\in \pA $ and the second object to $z$. Putting these all together we get a morphism
in $\precat (\mM )$,
$$
\pA \cup ^{\coprod _{z\in Z} \{ \tilde{t}_i\xi (|0)\} } 
\coprod _{z\in Z}\widetilde{\Xi}(N_z|N_z') \rightarrt^{T}
\Seg (\pB ),
$$
and now $T$ induces an isomorphism on sets of objects. It is no longer necessarily a cofibration. 
We would like to show that $T$ is a global weak equivalence. We start by considering pushouts along the interval $1$-category
$\overline{\pE}$.  For this proof we make essential use of the cartesian property of $\mM$ and the discussion of products in Chapter \ref{product1}.

\begin{corollary}
\label{pushoutinterval}
Suppose $\pA \in \precat (\mM )$, and $y\in \Ob (\pA )$. Suppose $N,N'$ are interval objects in $\mM$.
Then the pushout
morphism 
$$
\pA \rightarrt \pA \cup ^{\{ y\} } \widetilde{\Xi}(N_z| N'_x)
$$
obtained by identifying $\xi (|0)$ and $y$, is a global weak equivalence.
\end{corollary}
\begin{proof}
Apply Corollary \ref{pushoutgeninterval}
with $\pB = \Seg (\Xi (N_z| N'_x))$.
\end{proof}

\begin{corollary}
\label{Tisoeq}
In the situation described above the preceding corollary, the morphism 
$$
\pA \rightarrt \pA \cup ^{\coprod _{z\in Z} \{ \tilde{t}_i\xi (|0)\} } 
\coprod _{z\in Z}\widetilde{\Xi}(N_z|N_z')
$$
is a global weak equivalence. Given that $\pA \rightarrt \pB$ was a global weak equivalence, the functor
$$
T: \pA \cup ^{\coprod _{z\in Z} \{ \tilde{t}_i\xi (|0)\} } 
\coprod _{z\in Z}\widetilde{\Xi}(N_z|N_z')
\rightarrt \Seg (\pB )
$$
is a global weak equivalence inducing an isomorphism on sets of objects.
\end{corollary}
\begin{proof}
Choose a well-ordering of $Z$, giving an exhaustion of $Z$ by subsets $Z_i$ indexed by an ordinal $i\in \beta$. 
Let $\pB _i\subset \pB$ be the full subobject whose object set is $f(\Ob (A))\cup Z_i$.
By transfinite induction we obtain that the functors
$$
T_i: \pA \cup ^{\coprod _{z\in Z_i} \{ \tilde{t}_i\xi (|0)\} } 
\coprod _{z\in Z_i}\widetilde{\Xi}(N_z|N_z')
\rightarrt \Seg (\pB _i)
$$ 
are global weak equivalences, using the previous corollary at each step. At the end of the induction we obtain the
required statement. 
\end{proof}

We are now ready to show the preservation of global trivial cofibrations under pushouts. 

\begin{theorem}
\label{globalpushout}
Suppose $\pA \rightarrt \pB$ is an injective trivial cofibration. Suppose $\pA \rightarrt \pC$ is any morphism in $\precat (\mM )$.
Then the pushout morphism
$$
\pC \rightarrt \pC \cup ^{\pA } \pB 
$$
is a global weak equivalence.
\end{theorem}
\begin{proof}
We first show this statement assuming that all three objects $\pA$, $\pB$ and $\pC$ satisfy
the Segal conditions. Noting that $\pB \rightarrt \Seg (\pB )$ is an isomorphism on sets of objects and applying Lemma \ref{obisopushoutinvariance},
it suffices to show that the map 
$$
\pC \rightarrt \pC \cup ^{\pA } \Seg (\pB ) 
$$
is a global weak equivalence. Define
$$
\pF := \pA \cup ^{\coprod _{z\in Z} \{ \tilde{t}_i\xi (|0)\} } 
\coprod _{z\in Z}\widetilde{\Xi}(N_z|N_z') ,
$$
and consider the map $T:\pF \rightarrt \pB$ defined above. By Corollary \ref{Tisoeq}, $T$ is a global
weak equivalence inducing an isomorphism on sets of objects. By Lemma \ref{obisopushoutinvariance} it follows that the map
$$
\pC \cup ^{\pA } \pF \rightarrt \pC \cup ^{\pA } \Seg (\pB ) 
$$
is a global weak equivalence, so by 3 for 2 it suffices to show that 
$$
\pC \rightarrt \pC \cup ^{\pA } \pF
$$
is a global weak equivalence. But the map $\pA\rightarrt \pF$ is obtained as a transfinite composition of pushouts along
things of the form $\{ \xi (|0)\} \rightarrt \widetilde{\Xi}(N_z|N_z')$, and by
Corollary \ref{pushoutinterval} these pushouts are global weak equivalences. Thus, the map 
$\pC \rightarrt \pC \cup ^{\pA } \pF$ is a global weak equivalence, which finishes this part of the proof. 

Starting with $\pC \leftarr \pA \rightarrt \pB$ in general, let $\pA ':= \Seg (\pA )$, then
$$
\pB ' := \Seg (\pA ' \cup ^{\pA }\pB ), \;\;\; 
\pC ' := \Seg (\pA ' \cup ^{\pA }\pC ).
$$
We get a diagram 
$$
\begin{diagram}
\pC & \leftarr & \pA & \rightarr & \pB \\
\downarr & & \downarr && \downarr \\
\pC '& \leftarr & \pA '& \rightarr & \pB '
\end{diagram}
$$
such that the bottom row satisfies the hypothesis for the first part of the proof (all objects satisfy the Segal condition and the second
map is a global trivial cofibration), and such that the vertical arrows are global weak equivalences inducing isomorphisms on sets of objecs.
By Lemma \ref{obisopushoutinvariance}, the bottom map in the diagram 
$$
\begin{diagram}
\pC & \rightarr & \pC ' \\
\downarr & & \downarr \\
\pC \cup ^{\pA } \pB  & \rightarr & \pC ' \cup ^{\pA '} \pB '
\end{diagram}
$$
is a global weak equivalence. The top vertical map is a global weak equivalence by construction of $\pC '$ and the
right vertical map is one too, by the first part of the proof above. By 3 for 2 we conclude that the left vertical map
is a global weak equivalence, as required. 
\end{proof}

\section{A versality property}

The versality properties for the intervals constructed above, yield a similar versality property for any cofibrant
replacement of $\overline{\pE }$ if the target $\pA$  is fibrant in the 
diagram structure $\diag _c(\Delta ^i_{\Ob (\pA )}/\Ob (\pA ),\mM )$.  

\begin{theorem}
Suppose $\pA\in \precat (\mM )$, and suppose $\pA$ is fibrant as an object of $\precat _{c}(\Ob (\pA ), \mM )$ 
where $c$ indicates either the projective, the Reedy or the injective structures. Let $\pB \rightarrt \overline{\pE}$ be
a cofibrant replacement in $\precat _c([1], \mM )$, so $\Ob (\pB )$ is still $[1]=\{ \upsilon _0, \upsilon _1\}$.  
Then if $x,y\in \Ob (\pA )$ are two
objects, they project to equivalent objects in $\tau _{\leq 1}(\pA )$ if and only if 
there exists a morphism $\pB \rightarrt \pA$
sending $\upsilon _0$ to $x$ and $\upsilon _1$ to $y$.
\end{theorem}
\begin{proof}
Since $\tau _{\leq 1}\pB = \overline{\pE}$ is the category with two isomorphic objects,
existence of a map $\pB\rightarrt \pA$ sending $\upsilon _0$ to $x$ and $\upsilon _1$ to $y$
implies that $x$ and $y$ are internally equivalent in $\pA$. 

Suppose $x$ and $y$ are internally equivalent. 
If $\pA$ is a fibrant object for the Reedy or injective model 
structures relative to $\Ob (\pA )$, there is a morphism
$\Xi (N|N)\rightarrt \pA$ given by Theorem \ref{xiinterval}. For the
projective structure use $\Xi ^{\rm proj}(N|N')$ given by Remark \ref{XiNNforproj}
instead.
Denote either of these maps by $\pC\rightarrt \pA$. 
Let $\widetilde{\pC}\subset \Seg (\pC )$ be the
full subcategory consisting of only the two main objects, but
identify $\Ob (\widetilde{\pC})$
with the two element set $[1]=\{ \upsilon _0,\upsilon _1\} =\Ob (\pE )$. 
The map $\pC \rightarrt \Seg (\pC )$ is a isotrivial cofibration
so 
$$
\pA \rightarrt \pA \cup ^{\pC}\Seg (\pC )
$$
is an isotrivial cofibration by Theorem \ref{obisopushout}, it follows that
our map extends to $\Seg (\pC )\rightarrt \pA$. This now restricts to a  map
$$
\widetilde{\pC} \rightarrt \pA
$$
sending $\upsilon _0$ to $x$ and $\upsilon _1$ to $y$. 
By contractibility, Theorem \ref{xicontractible}, 
$$
\widetilde{\pC}\rightarrt \overline{\pE}
$$
is a weak equivalence inducing an isomorphism on sets of objects. 
Choose a factorization
$$
\widetilde{\pC}\rightarrt^i \widetilde{\pC}'\rightarrt^p \overline{\pE}
$$
where $i$ is a trivial cofibration 
and 
$p$ is a trivial fibration in $\precat ([1],\mM )$. Again our map
extends to $\widetilde{\pC}'\rightarrt \pA$, but now 
since $\pB$ is cofibrant and $p$ is a trivial fibration there is
a lifting $\pB \rightarrt \widetilde{\pC}'$ inducing the identity
on the set of objects. We get the required map $\pB \rightarrt \pA$. 
\end{proof}

The importance of this versality property is that it allows us to replace a global weak equivalence by one which is
surjective on sets of objects. 

\begin{corollary}
\label{versality}
Let $\pB \rightarrt \overline{\pE}$ be
a cofibrant replacement in $\precat _c([1], \mM )$
Suppose $f:\pA\rightarrt \pC$ is a global weak equivalence, and suppose $\pC$ is a fibrant object in 
$\precat _{c}(\Ob (\pC ), \mM )$. Then there exists a pushout $\pA \rightarrt \pA '$ by a collection of copies of 
$\{ \upsilon _0 \} \hookrightarrow \pB$, and a map $\pA '\rightarrt \pC$ which is a global weak equivalence and a
surjection on sets of objects. 
\end{corollary}
\begin{proof}
For each object $y\in \Ob (\pC )$ choose $x\in \Ob (\pA )$ such that $f(x)$ is
internally equivalent to $y$. For each such pair we get a map $\pB \rightarrt \pC$
sending $\upsilon _0$ to $f(x)$ and $\upsilon _1$ to $y$; attaching a copy of $\pB$
to $\pA$ by sending $\upsilon _0$ to $x$ and then doing this for all objects $y$
we obtain the required pushout $\pA '$ and extension of the map. 
\end{proof}


\chapter{The model category of $\mM$-enriched precategories}
\label{mproof1}

In this chapter, we finish the proof that the category $\precat (\mM )$ of $\mM$-enriched precategories, 
with variable set of objects, has natural injective and projective model structures. Given the 
product theorem of Chapter \ref{product1} and the discussion of intervals in Chapter \ref{interval1}, the proof presents no further obstacles. 
We also show that the Reedy structure
$\precat _{\rm Reedy}(\mM )$ is
again tractable, left proper and cartesian,
allowing us to iterate the operation. 

\section{A standard factorization}

Follow up on Corollary \ref{versality} of the previous chapter, by analyzing further the case of maps which are
surjective on the set of objects. 

\begin{lemma}
\label{standardfactor}
Suppose $f:\pA \rightarrt \pB$ is a morphism in $\precat (\mM )$ such that $\Ob (f)$ is surjective. 
Let $\Ob (f)^{\ast}(\pB )\in \precat (\Ob (\pA ); \mM )$
be the precategory obtained by pulling back along $\Ob (f):\Ob (\pA )\rightarrt \Ob (\pB )$. Then $f$ factors as
$$
\pA \rightarrt \Ob (f)^{\ast}(\pB ) \rightarrt \pB ,
$$
where the first map is an isomorphism on sets of objects, and  
the second map $\Ob (f)^{\ast}(\pB ) \rightarrt \pB $ satisfies the right lifting property with respect to any
morphism $g:\pU \rightarrt \pV$ such that $\Ob (g)$ is injective. 
\end{lemma}
\begin{proof}
Given a diagram 
$$
\begin{diagram}
\pU & \rightarr & \Ob (f)^{\ast}(\pB ) \\
\downarr & & \downarr \\
\pV & \rightarr & \pB
\end{diagram}
$$ 
in order to get a lifting it suffices to choose a lifting on the level of objects
$$
\Ob (\pV )\rightarrt
\Ob (\Ob (f)^{\ast}(\pB ))=\Ob (\pA ). 
$$
This is possible since $\Ob (q)$ is
injective and $\Ob (f)$ surjective. 
\end{proof}

\begin{corollary}
\label{surjlift}
In the situation of the lemma, if $I$ is a set of morphisms in $\precat (\mM )$
which are all injective on the level of objects, and if the first map $\pA \rightarrt \Ob (f)^{\ast}(\pB )$ is
in $\inj (I)$ then $f\in \inj (I)$. 
\end{corollary}
\begin{proof}
Combine the lifting property for $\inj (I)$ with the one of the previous lemma. 
\end{proof}

\section{The model structures}

We will be applying Theorem \ref{recog} of Chapter \ref{modcat1} to construct the model structure on $\precat (\mM )$.

We fix a class of cofibrations denoted generically by $c$, with $c={\rm proj}$ or 
$c= {\rm Reedy}$.
This choice determines the corresponding notions of cofibrations in $\precat _c(\mM )$ or $\precat _c(X; \mM )$. 
Let $I$ be a set of generators for the $c$-cofibrations in $\precat _c(\mM )$, as discussed in Chapter \ref{cofib1}. 
We can choose $I$ to consist of maps with $c$-cofibrant domains, for the
Reedy and projective structures, see Chapter \ref{cofib1}. 

Let $K_{\rm loc}$ denote a set of morphisms which are pseudo-generators for the local weak equivalences in $\precat _c([k],\mM )$, as from Theorem \ref{modstrucs} or Theorem \ref{reedyprecat}.  We may assume that they are $c$-cofibrations,
with $c$-cofibrant domains if $c$ is Reedy or projective. 

Recall that $\pE = \Upsilon (\ast )$ is the category with a single non-identity morphism $\upsilon _0\rightarrt \upsilon _1$,
and $\overline{\pE }$ is the category obtained by inverting this map, that is with a single isomorphism $\upsilon _0\cong \upsilon _1$. 
Consider the morphism $ \{ \upsilon _0 \} \hookrightarrow \overline{\pE }$. Choose a $c$-cofibrant replacement 
$$
\{ \upsilon _0,\upsilon _1 \}\rightarrt \pP \rightarrt \overline{\pE }
$$
in the model category $\precat _c([1], \mM )$, and let 
$$
\{ \upsilon _0 \} \rightarrt^{i_0} \pP
$$
denote the inclusion morphism of a single object. This is still a $c$-cofibration in $\precat _c(\mM )$ (because of Condition (AST) in Definition \ref{def-cartesian}). 

Let $K_{\rm glob} := K_{\rm loc} \cup \{ i_0 \}$. Note that the domain of $i_0$ is the single object precategory
$\{ \upsilon _0 \}$ which is $c$-cofibrant for any of the $c$.

\begin{theorem}
\label{globalmodel}
The class of global weak equivalences is pseudo-generated by $K_{\rm glob}$ in the sense of construction (PG) of Section \ref{sec-pseudogen}.
Furthermore, $I$ and $J$ satisfy axioms (PGM1)--(PGM6), so they define tractable left proper model structures by Theorem \ref{recog}. 
For these model structures, the weak equivalences are the global weak equivalences; the cofibrations are the 
projective  (resp. injective resp. Reedy) cofibrations, and the fibrations are the 
projective  (resp. injective resp. Reedy) global fibrations. 
\end{theorem}
\begin{proof}
We have to show that $K_{\rm glob}$ leads to the class of global weak equivalences via prescription (PG). This amounts to showing
that a map $f:\pX\rightarrt \pY$ is a global weak equivalence if and only if there exists a diagram 
$$
\begin{diagram}
\pX & \rightarr & \pA \\
\downarr & & \downarr \\
\pY & \rightarrt & \pB
\end{diagram}
$$ 
with the horizontal maps in $\cell (K_{\rm glob})$ and the right vertical map in $\inj (I)$.

The maps in $K_{\rm glob}$ are trivial cofibrations in the projective structure, and the global trivial cofibrations are preserved
by pushout (Theorem \ref{globalpushout}) and transfinite composition 
(Lemma \ref{globaltranfinite}), so the maps in $\cell (K_{\rm glob})$ are global trivial cofibrations.
By 3 for 2 for global trivial cofibrations (Proposition \ref{globalretract32}) 
it follows that if there exists a square diagram as above then $f$ is a global weak
equivalence.

Suppose $f$ is a global weak equivalence, we would like to construct a square as above. Let $r:\pY\rightarrt \pB$ be the map given
by applying the small object argument to $\pY$ with respect to $K_{\rm loc}$. Thus $\pB$ is $K_{\rm loc}$-injective. 
It follows that it satisfies the Segal conditions, and is a fibrant object in $\precat _{c}(\Ob (\pB ), \mM )$ (see Theorems \ref{modstrucs} and \ref{reedyprecat}).

By Corollary \ref{versality} of the preceding Chapter \ref{interval1}, there exists a pushout $\pX \rightarrt \pX '$ by 
a collection of copies of the map $i_0:\{\upsilon _0\} \rightarrt \pP$ and an extension of $rf$ to a global weak equivalence 
$g:\pX '\rightarrt \pB$ which is surjective on the set of objects. Note that $\pX\rightarrt \pX '$ is in $\cell (K_{\rm glob})$. Consider the factorization 
$$
\pX ' \rightarrt \Ob (g) ^{\ast} \pB \rightarrt \pB
$$
of Lemma \ref{standardfactor} above. 
The first map is an isomorphism on the set of objects, so it can be considered as a map 
in $\precat (\Ob (\pX '),\mM )$. 
Apply the small object argument for the set $K_{\rm loc}$, to the first map to yield a factorization
$$
\pX ' \rightarrt \pA \rightarrt \Ob (g) ^{\ast} \pB
$$
such that $\pX ' \rightarrt \pA $ is in $\cell (K_{\rm loc})$ and 
$\pA \rightarrt \Ob (g) ^{\ast} \pB$ is a fibration in $\precat (\Ob (\pX '),\mM )$.
However, the composed map $\pX ' \rightarrt \Ob (g) ^{\ast} \pB $ is a local weak equivalence, so by 3 for 2 in the local 
model structure, the map $\pA \rightarrt \Ob (g) ^{\ast} \pB$ is a trivial fibration; hence it is in $\inj (I)$. 
Apply now Corollary \ref{surjlift}. Note that the factorization of Lemma \ref{standardfactor} for the map $\pA \rightarrt \pB$ is
just 
$$
\pA \rightarrt \Ob (g) ^{\ast} \pB \rightarrt \pB
$$
where the first map is the same as previously; thus the first map is in $\inj (I)$ and by Corollary \ref{surjlift} the full map
$\pA \rightarrt \pB$ is in $\inj (I)$. The composition 
$$
\pX \rightarrt \pX ' \rightarrt \pA
$$
of two maps in $\cell (K_{\rm glob})$ is again in $\cell (K_{\rm glob})$. This completes the verification that our global weak equivalence
$f$ satisfies the condition (PG). 

We now verify axioms (PGM1)--(PGM6) needed to apply Theorem \ref{recog} of Chapter \ref{modcat1}. 
\newline
(PGM1)---by hypothesis $\mM$ is locally presentable, and $I$ and $K_{\rm glob}$ are chosen to be small sets of morphisms;
\newline
(PGM2)---we  have chosen the domains of arrows in $I$ and $K_{\rm glob}$ to be cofibrant, and $K_{\rm glob}$ consists
of $c$-cofibrations in other words it is contained in $\cof (I)$; 
\newline
(PGM3)---the class of global weak equivalences is closed under retracts by Proposition \ref{globalretract32} of Chapter \ref{weakenr1};
\newline
(PGM4)---the class of global weak equivalences satisfies 3 for 2 again by Proposition \ref{globalretract32};
\newline
(PGM5)---the class of global trivial $c$-cofibrations is closed under pushouts by  
Theorem \ref{globalpushout};
\newline
(PGM6)---the class of global trivial $c$-cofibrations is closed under transfinite composition, indeed the cofibrations are
closed unter transfinite composition since they have generating sets,
see Chapter \ref{cofib1}, and the weak equivalences are too by Lemma \ref{globaltranfinite}.

Theorem \ref{recog} now applies to show that $\precat (\mM )$ with the given classes of $c$-cofibrations,
global weak equivalences, hence global trivial $c$-cofibrations whence global $c$-fibrations, is a tractable left proper closed model category. 
\end{proof}

\section{The cartesian property}

\begin{lemma}
\label{cartesianproperty}
Suppose $\pA \rightarrt \pB$ and $\pC \rightarrt \pD$ are global Reedy cofibrations, with the first one being a weak equivalence. Then the map
$$
\pA \times \pD \cup ^{\pA\times \pC} \pB \times \pC \rightarrt \pB \times \pD
$$
is a global trivial Reedy cofibration.
\end{lemma}
\begin{proof}
It is a Reedy cofibration by Corollary \ref{reedycofcart}. 
We just have to show that it is a global weak equivalence. 
By Corollary \ref{prodwe}, the vertical maps in the diagram 
$$
\begin{diagram}
\pA \times \pC & \rightarr & \pA \times \pD \\
\downarr & & \downarr \\
\pB \times \pC & \rightarr & \pB \times \pD
\end{diagram}
$$ 
are global weak equivalences. Applying
Theorem \ref{globalpushout} to pushout along the left vertical 
global trivial cofibration, then using 3 for 2, 
it follows that the map
$$
\pA\times \pD\rightarrt \pA \times \pD \cup ^{\pA\times \pC} \pB \times \pC
$$
is a global weak equivalence. Then by 3 for 2 using the right vertical global weak equivalence,
the map in the statement of the lemma is a global weak equivalence.
\end{proof}

\begin{theorem}
\label{maintheorem}
Suppose $\mM$ is a tractable left proper cartesian  model category. Then 
the model category $\precat _{\rm Reedy}(\mM )$ of $\mM$-enriched precategories with Reedy cofibrations is again 
a tractable left proper cartesian  model category. 
\end{theorem}
\begin{proof}
Observe first of all that direct product commutes with colimits in $\precat (\mM )$, as can be seen from the explicit description of products and 
colimits and using the corresponding condition for $\mM$.

Next, note that the map $\emptyset \rightarrt \ast$ is a Reedy cofibration, from the definition. 

Proposition \ref{reedyproduct} gives the cofibrant property of the pushout-product map; and the previous
Lemma \ref{cartesianproperty} gives the trivial cofibration property. This shows that $\precat _{\rm Reedy}(\mM )$ is cartesian.

To finish the proof we need to note that it is tractable. This can be seen by inspection of the generating cofibrations for the Reedy structure,
given in Proposition \ref{reedygenerators}.
\end{proof}

Of course, the projective model structure is definitely not cartesian. 
On the other hand, one can hope to treat the injective model structure. 
There is already a problem with tractability: 
Lurie and Barwick don't mention if their constructions of injective model categories preserve the tractability condition, at least until a most recent version of a paper in
which Barwick states this property.
That would clearly be an important result, giving tractability for 
$\precat _{\rm inj}(\mM )$ in general. In the case of presheaf categories with monomorphisms as cofibrations, of course
this condition becomes automatic. Similarly, it doesn't seem immediately clear whether $\precat _{\rm inj}(\mM )$ will satisfy
condition (PROD) in general, although again this is relatively easy to see for presheaf categories with monomorphisms as cofibrations. 

\section{Properties of fibrant objects}

Julie Bergner made the very interesting observation \cite{BergnerSegal}
that one could give an
explicit characterization for the fibrant objects in the case of Segal categories 
$\mM = \mK$. We get the same kind of property in general.  

\begin{proposition}
\label{fibprop}
Let $c={\rm proj}$ or $c={\rm Reedy}$. 
In the model category $\precat _{c}(\mM )$ 
constructed above, an object $\pA$ with $\Ob (\pA )=X$ 
is fibrant if and only if it is fibrant
when considered as an object of the model category $\precat _{c}(X,\mM )$ of Theorem \ref{modstrucs} or \ref{reedyprecat}.
In turn this condition is equivalent to saying that $\pA$ satisfies the
Segal conditions, and is fibrant as an object of the unital diagram model category
$\diag _{c}(\Delta ^o_X/X,\mM )$. 
\end{proposition}
\begin{proof}
Left to the reader in the current version.
\end{proof}

For $c={\rm proj}$ then, an $\mM$-precategory $\pA$ is fibrant if and only if
it satisfies the Segal conditions, and is levelwise fibrant. For $c={\rm Reedy}$
the fibrancy condition is also pretty easy to check: it just means that the
standart matching maps are fibrations in $\mM$.  


\section{The model category of strict $\mM$-enriched categories}
\label{strictMcats}

Dwyer and Kan proposed, in a series of papers, a model category structure on the category of strict simplicial categories.
Their program was finished by Bergner \cite{BergnerModel}. Lurie then generalized this to construct a model category of strict $\mM$-enriched
categories in \cite[Appendix]{LurieTopos}, then used that to construct the model category of weakly $\mM$-enriched precategories
as we have done above. 

\begin{theorem}
Suppose $\mM$ is a tractable left proper cartesian model category. Let $\Cat (\mM )$ denote the category of strict
$\mM$-enriched categories. Define the notion of weak equivalence in the usual way (see Section \ref{sec-globalwe}). Then $\Cat (\mM )$ has a 
tractable left proper model structure in which the generating cofibrations are obtained by free additions of generating cofibrations of $\mM$  in the
morphism space between any two objects. There is a Quillen adjunction 
$$
\Cat (\mM ) \rBotharrow \precat _{\rm proj}(\mM )
$$
and indeed the model structure on  $\Cat (\mM )$ can be used to generate the model structure on $\precat _{\rm proj}(\mM )$. 
However, $\Cat (\mM )$ is not in general cartesian. It follows that any object of $\precat _{\rm proj}(\mM )$
is equivalent to a strict $\mM$-enriched category. 
\end{theorem}
\begin{proof}
See Bergner \cite{BergnerModel} for $\mM = \mK$ and Lurie \cite{LurieTopos} for arbitrary $\mM$. The strictification result, for the case of Tamsamani $n$-groupoids,
was proven by Paoli \cite{PaoliAdvances} using techniques of $Cat ^n$-groups. 
\end{proof}

This theorem offers an alternative route to the construction of the model structure on $\precat _{\rm proj}(\mM )$,
whose proof is somewhat different from ours. The advantage is that it also gives the model structure on $\Cat (\mM )$ and hence the
strictificaton result; the disadvantage is that it doesn't give the cartesian property. The cartesian question has been treated by Rezk 
\cite{RezkCartesian} for the case of iterated Rezk categories. We leave it to the reader to explore these different points of view.

\chapter{Iterated higher categories}
\label{chap-iterate}

The conclusion of Theorem \ref{maintheorem} 
matches the hypotheses we imposed that $\mM$ be tractable, left proper and cartesian.
Therefore, we can iterate the construction to obtain various versions of model categories
for $n$-categories and similar objects. This process is inherent in the definitions of 
Tamsamani \cite{Tamsamani} and Pelissier \cite{Pelissier}. 
Rezk considered the corresponding
iteration of his definition in \cite{RezkCartesian} following Barwick, and Trimble's
definition is also iterative. 
Such an iteration is also related to Dunn's iteration of the Segal delooping machine 
\cite{Dunn}, and goes back to the well-known iterative presentation of the notion of
strict $n$-category, see Bourn \cite{Bourn} for example. 

In what follows unless otherwise indicated, the model category $\precat (\mM )$ will mean by definition the Reedy structure $\precat _{\rm Reedy}(\mM )$. 

For any $n\geq 0$ define by induction $\precat ^0(\mM ):= \mM$ and for $n\geq 1$
$$
\precat ^{n}(\mM ):= \precat (\precat ^{n-1}(\mM )).
$$
This is the model category of {\em $\mM$-enriched $n$-precategories}. Notations for objects therein will be discussed below. 

In the iterated situation, we can introduce the following definition. 

\begin{definition}
\label{fullsegal}
An $\mM$-enriched $n$-precategory $\pA\in \precat ^n(\mM )$ satisfies the {\em full Segal condition} if
it satisfies the Segal condition as an $\precat ^{n-1}(\mM )$-precategory, and 
furthermore 
inductively for any sequence of objects $x_0,\ldots , x_m\in \Ob (\pA )$ the 
$\mM$-enriched $(n-1)$-precategory
$$
\pA (x_0,\ldots , x_m)\in \precat ^{n-1}(\mM )
$$
satisfies the full Segal condition. 
\end{definition}

\begin{lemma}
If $\pA$ is a fibrant object in the (iterated Reedy) model structure on 
$\precat ^n(\mM )$, then $\pA$ satisfies the full Segal condition. 
\end{lemma}
\begin{proof}
The Segal condition for the $\precat ^{n-1}(\mM )$-precategory comes
from Proposition \ref{fibprop} (see Theorem \ref{reedyprecat}). However, if
$\pA$ is fibrant then the $\pA (x_0,\ldots , x_m)$ are fibrant in $\precat ^{n-1}(\mM )$
so by induction they also satisfy the full Segal condition. 
\end{proof}

\section{Initialization}
\label{sec-initialization}

Here are a few possible choices for $\mM$ to start with. 

If $\mM=\Sets$ is the model category of sets,
with cofibrations and fibrations being arbitrary morphisms and weak equivalences being isomorphisms, then 
$\precat ^n(\Sets )$ is the {\em model category of $n$-precats} which was considered in \cite{svk} and for which we have now
fixed up the proof. 

Let $\ast$ denote the model category with a single object and a single morphism. Then $\precat (\ast )$ is Quillen equivalent
(by a product-preserving map) to the model category $\{ \emptyset , \ast \}$ consisting of the emptyset and the one-element set,
where weak equivalences are isomorphisms. Iterating again, $\precat ^2(\ast )$ is Quillen-equivalent to $\precat (\{ \emptyset , \ast \})$.
These are both model categories of {\em graphs}, the first allowing multiple edges between nodes and the second allowing only 
zero or one edges between two nodes. The weak equivalences are defined by requiring isomorphism on the level of  $\pi _0$ i.e. the set of connected components of
a graph. These model categories of graphs are Quillen-equivalent to $\Sets$ but have the advantage that the cofibrations are monomorphisms. They are related to the notion of 
{\em setoid} in constructive type theory. 

In particular, $\precat ^{n+2}(\ast )$ is Quillen-equivalent to $\precat (\Sets )$ and should perhaps be thought of as the ``true'' model category
of $n$-categories. 

If we start with $\mM = \mK$ the Kan-Quillen model category of simplicial sets, then $PC ^n(\mK )$ is the model category of Segal $n$-precategories
introduced in \cite{descente}. 

One can imagine further constructions starting with $\mM$ as a category of diagrams or other such things. Starting with $\zz / 2$-equivariant sets
should be useful for considering $n$-categories with duals.

\section{Notations}
\label{sec-notations}

By Lemma \ref{precatdisj}, once we start iterating, the hypothesis (DISJ) of \ref{sec-interpretations} will be in vigour.
Furthermore, most of our examples of starting categories (even $\mM = \ast$) satisfy (DISJ), see Lemma \ref{startdisj}.
Whenever such is the case, 
it is reasonable to introduce an iteration of the notation $A_{n/}$ of 
Section \ref{sec-interpretations}.

One can note, on the other hand, that even based on this notation as the general framework, most of what was done in
\cite{svk} and \cite{limits} really used the notation $A(x_0,\ldots , x_n)$ at the crucial places. So, in a certain sense,
the notations we introduce here are not really the fundamental objects, nonetheless it is
convenient to have them for comparison.

In $\precat ^n(\mM )$ for any $k\leq m$ and any
multi-index $m_1,\ldots , m_k$ we can introduce the  notation 
$\pA _{m_1,\ldots , m_k/} \in \precat ^{n-k}(\mM )$ defined 
by induction on $k$. 
At the initial $k=1$ (whenever $n\geq 1$),
by noting that $\pA \in \precat (\precat ^{n-1}(\mM ))$ we can use
the notation 
$$
\pA _{m/}\in \precat ^{n-1}(\mM )
$$
considered in Section \ref{sec-interpretations}. 
Then for $k\geq 2$ define 
inductively
$$
\pA _{m_1,\ldots , m_k/} := (\pA _{m_1,\ldots , m_{k-1}/})_{m_k/} \in \precat ^{n-k}(\mM ).
$$
For $k<n$, define on the other hand 
$$
\pA _{m_1,\ldots , m_k} := \Ob (\pA _{m_1,\ldots , m_k/})\in \Sets .
$$

One can remark that the notations 
$\pA _{m_1,\ldots , m_k/}\in \precat ^{n-k}(\mM )$ 
and $\pA _{m_1,\ldots , m_k}\in \Sets$ make sense
for $k<n$ because we have seen that $\precat (\mM )$ satisfies Condition (DISJ), even if $\mM$ itself does not.

In a related direction, notice that if $\mM$ is a presheaf category then
$\precat (\mM )$ is also a presheaf category, by the discussion of Section 
\ref{sec-interpretations}. Inductively the same is true of $\precat ^n(\mM )$.
If $\mM = \presh (\Phi )$ then $\precat ^n(\mM )=\presh (\scone ^n(\Phi ))$
in the notations of Proposition \ref{upresheafcat}.

\section{The case of $n$-nerves}
\label{sec-ncatcase}

Start with $\mM = \Sets$ with the trivial model structure
(as Lurie calls it \cite{LurieTopos}), where the weak equivalences are isomorphisms and
the cofibrations and fibrations are arbitrary maps.  Iterating the construction
of Theorem \ref{maintheorem} we obtain 
the iterated Reedy model category structure 
$$
\precat ^n(\Sets ).
$$
The underlying category is the category of presheaves of sets on an iterated
version of the construction of Section \ref{sec-interpretations},
$$
\precat ^n(\Sets ) = \presh (\scone ^n(\ast )).
$$
The underlying category $\scone ^n(\ast )$ may be seen as a quotient of
$\Delta ^n=\Delta \times \cdots \times \Delta$, indeed it is the same as
the category which was denoted $\Theta ^n$ in \cite{svk} and \cite{limits}. 
If $\pA \in \precat ^n(\Sets )$ then the notation discussed in
the previous section applies, and for any  multi-index $(m_1,\ldots , m_k)$
with $k\leq n$ we get a set $\pA _{m_1,\ldots , m_k}\in \Sets$.
For $k=n$, the notation $\pA _{m_1,\ldots , m_n}:=\pA _{m_1,\ldots , m_n/}$ may be used
since the model category $\mM$ is equal to $\Sets$.  

That yields a system of notations coinciding with that of \cite{svk}, \cite{limits}, \cite{BBDSH} etc. A slight difference is that for $\pA \in  \precat ^n(\Sets )$,
and for any sequence of objects $x_0,\ldots , x_m$, what we would be denoting here by
$$
\pA (x_0,\ldots , x_m) \in \precat ^{n-1}(\Sets )
$$ 
was denoted in those preprints by $\pA _{m/}(x_0,\ldots , x_m)$.
We have dropped the subscript $(\; )_{m/}$ for brevity.

The iterated injective and Reedy model structures coincide in the case of
$\precat ^n(\Sets )$, by applying Proposition \ref{reedyinjective} inductively.  A map $\pA \rightarrt \pB$
is a cofibration if and only if, for any multiindex $m_1,\ldots ,m_k$ with $k<n$
the map $\pA _{m_1,\ldots , m_k}\rightarrt \pB_{m_1,\ldots , m_k}$ is an injection
of sets.  The monomorphism condition is not imposed at multiindices of length $k=n$, indeed
at the top level,
all maps of $\Sets$ are cofibrations for its trivial model structure. 
The reader may refer to \cite{svk} for
a fuller discussion of the notion of cofibrations.

We call $\precat ^n(\Sets )$ the category of {\em $n$-prenerves}\footnote{The
objects of $\precat ^n(\Sets )$ were also called {\em $n$-precats} in \cite{svk} \cite{limits}}; and the
objects satisfying the full Segal condition 
are the {\em $n$-nerves} of Tamsamani \cite{Tamsamani}.  A fibrant object of 
$\precat ^n(\Sets )$ is an $n$-nerve, indeed it satisfies the Segal conditions
at the last iteration (corresponding to the first element $m_1$ of a
multiindex), and furthermore the $n-1$-prenerves $\pA (x_0,\ldots , x_m)$ are themselves
fibrant in $\precat ^{n-1}(\Sets )$ so by induction they also satisfy the Segal conditions
at all of their levels. Taking the disjoint union over all sequences $x_0,\ldots , x_m$
yields $\pA _{m/}$ which is an $n-1$-nerve. 

At $n=1$, the category of $1$-prenerves is the category of simplicial sets, and
the $1$-nerves are the simplicial sets which are nerves of a $1$-category, that is
to say the category of $1$-nerves is equivalent to $\Cat$. The process $\pA \mapsto \Seg (\pA )$ is the generation of a category by generators and relations discussed in 
Section \ref{genrel1cat}.

\section{Truncation and equivalences}
\label{sec-equivalences}

The definition of weak equivalence we have adopted for $\precat (\mM )$ in general is designed for enrichment over a general model category $\mM$.
In Tamsamani's original definition of $n$-nerves, the notion of equivalence and 
the truncation operations $\tau _{\leq k}$ were
defined inductively along the way. So, in case of $\precat ^{n}(\Sets )$ 
there remains the question of equating these two definitions
of equivalences.  

For any tractable left proper cartesian model category $\mM$, define the {\em pretruncation}
$$
\tau ^p_{\leq n} : \precat ^n(\mM )\rightarrt \precat ^n(\Sets )
$$
as the functor induced by $\tau _{\leq 0}:\mM \rightarrt \Sets$. Applied to 
$\precat ^{n-k}(\mM )$ for any $0\leq k\leq n$, this gives a pretruncation functor 
$$
\tau ^p_{\leq k} : \precat ^n(\mM )\rightarrt \precat ^k(\Sets ). 
$$
If $\mM = \Sets$ and $n=k$ then it is the identity. Recall that 
for $\pA \in\precat (\mM )$ the truncation operation was defined by
$\tau _{\leq 1}(\pA ):= \tau ^p_{\leq 1}(\Seg (\pA ))$. 

\begin{remark}
It doesn't seem to be true in general that the truncation
functor from $\precat (\mM )$ to $\Cat$ 
used starting in Chapter \ref{precat1} could be expressed 
in terms of  the generators and relations operation from $1$-prenerves to $1$-nerves
as
$$
\tau _{\leq 1}(\pA ) \cong \Seg (\tau ^p_{\leq 1}(\pA )).
$$
Indeed, the operation $\Seg (\pA )$ might alter things in a way which isn't seen on the
level of $1$-truncation.
\end{remark}

One should impose the full Segal condition in order to be able to use the pretruncation. 

\begin{proposition}
If $\pA$ is an $\mM$-enriched $n$-precategory which satisfies the full Segal condition,
then for any $k\leq n$ the pretruncation $\tau ^p_{\leq k}(\pA )$ is a $k$-nerve,
and these truncations may be composed leading at $k=1$ to the usual truncation $\tau _{\leq 1}$. They are compatible with direct products. 

Suppose $\pA \rightarrt^f \pB$ is a weak equivalence in $\precat ^n(\mM )$,
and $\pA$ and $\pB$ both satisfy the full Segal condition. Then 
for any $0\leq k\leq n$, the truncation $\tau ^p_{\leq k}(f)$ is an equivalence of $k$-nerves
in the sense of \cite{Tamsamani}. 

For $\mM = \Sets$, a morphism $\pA \rightarrt^f \pB$ in $\precat ^n(\Sets )$
between $n$-nerves, is a weak equivalence if and only if it is an 
equivalence of $n$-nerves
in the sense of \cite{Tamsamani}. 
\end{proposition}
\begin{proof}
For $n=1$, see Lemma \ref{truncationproperties}. We still should show the compatibility
with composing truncation operations. 
If $\mP := \precat (\mM )$ then we have defined the truncation $\tau _{\leq 0}:\mP \rightarrt \Sets$ using the model structure of $\mP$: $\tau _{\leq 0}(\pA )$
it is the set of morphisms from $\ast$ to $\pA$ in $\Ho (\mP )$. On the other hand,
we have defined the truncation denoted also $\tau _{\leq 0}:\precat (\mM )\rightarrt  \Sets$
as sending $\pA$ to the set of isomorphism classes of $\tau _{\leq 1}(\pA )$. 
To show that they are the same, note first that
both are invariant under weak equivalences so we may assume that $\pA$ is fibrant.
Then $\Hom _{\Ho (\mP )}(\ast , \pA )$ is the set of morphisms from $\ast$ to $\pA$,
up to the relation of homotopy. The set of morphisms is just $\Ob (\pA )$, and
the relation of homotopy says that $x$ is equivalent to $y$ if and only if there
exists a map from an interval object to $\pA$ sending the endpoints to $x$ and $y$
respectively (see Lemma \ref{homotopicequiv}). If this condition holds
then looking at the image of the interval in $\tau _{\leq 1}(\pA )$ we
conclude that the points $x$ and $y$ go to the same isomorphism class.
In the other direction, if $x$ and $y$ go to isomorphic objects
in $\tau _{\leq 1}(\pA )$ then by the versality property Theorem \ref{xiinterval} plus
the contractibility of the intervals in question, Theorem \ref{xicontractible}, 
the corresponding maps $x,y:\ast \rightarrt \pA$ are homotopic. This shows that the
two versions of $\tau _{\leq 0}( \pA )$ coincide. 

Assume that $n\geq 2$
and the proposition is known for  $\precat ^{n-1}(\mM )$. For $k=1$,
the  truncation $\tau ^p_{\leq 1}(\pA )$ is the operation of
Lemma \ref{truncationproperties} which corresponds to the right truncation when
applied to $\pA$ satisfying the Segal conditions. For $k\geq 2$,
we have the truncation functor $\tau ^p_{\leq k-1}:\precat ^{n-1}(\mM )\rightarrt
\precat ^{k-1}(\Sets )$, taking weak equivalences between objects which satisfy
the full Segal condition, to weak equivalences. It follows that
if $\pA \in \precat ^n(\mM )$ satisfies the full Segal conditions,
then applying $\tau ^p_{\leq k-1}$
levelwise to $\pA$ considered as a diagram from $\Delta _{\Ob (\pA )}^o$ to 
$\precat ^{n-1}(\mM )$, it yields a diagram from $\Delta _{\Ob (\pA )}^o$ to 
$\precat ^{k-1}(\Sets )$ which again satisifes the Segal conditions, as well as the 
full Segal conditions levelwise.  But $\tau ^p_{\leq k-1}$ applied
levelwise is by definition $\tau ^p_{\leq k}$. This shows that $\tau ^p_{\leq k}(\pA )$
is a $k$-nerve. The composition of two truncation operations is again a truncation:
$$
\tau _{\leq r}^p( \tau _{\leq k}^p(\pA )) = \tau _{\leq r}^p(\pA )
$$
whenever $r\leq k$. By induction they are compatible with direct products. 

Suppose $\pA \rightarrt^f \pB$ is a weak equivalence in $\precat ^n(\mM )$,
and $\pA$ and $\pB$ both satisfy the full Segal condition. Since they
satisfy the regular Segal condition, $f$ is essentially surjective, meaning
that $\tau _{\leq 0}(\pA )\rightarrt \tau _{\leq 0}(\pB )$ is surjective, 
and induces equivalences $\pA (x,y)\rightarrt \pB (f(x),f(y))$ for all pairs
of objects $x,y\in \Ob (\pA )$. But $\pA (x,y)$ and $\pB (f(x),f(y))$ also
satisfy the full Segal condition, so by the inductive statement known for 
$\precat ^{n-1}(\mM )$,
$f$ induces equivalences of $k-1$-nerves
$$
(\tau _{\leq k}\pA )(x,y)=
\tau _{\leq k-1}\pA (x,y)\rightarrt^{\sim} \tau _{\leq k-1}(\pB (f(x),f(y)))
=(\tau _{\leq k}\pB )(f(x),f(y)).
$$
This now implies that $f$ induces an equivalence of $k$-nerves from $\tau _{\leq k}\pA$
to $\tau _{\leq k}(\pB )$. 

In the case $\mM=\Sets$, the above argument works in the other direction to show that
if $f$ is an equivalence of $n$-nerves in the sense of \cite{Tamsamani}, it is
essentially surjective and, by applying the inductive statement for 
$n-1$, it is also fully faithful. 
\end{proof}

\section{The $(n+1)$-category $nCAT$}
\label{sec-ncatcat1}

The cartesian model category structure on $\precat ^n(\mM )$ allows us to define a
structure of $\mM$-enriched $n+1$-category denoted $CAT(n;\mM )$.  In particular, starting with $\mM = \Sets$ we obtain the
$n+1$-category $nCAT = CAT(n;\Sets )$ which was originally discussed in Chapter \ref{ncats1}.

In Section \ref{sec-enr-cart}, starting with a tractable cartesian
model category $\mP$ we get the strict $\mP$-enriched category
$\Enr (\mP )$. Recall that the objects of $\Enr (\mP )$ are the cofibrant and fibrant
objects of $\mP$, and if $X,Y$ are two such objects then
the  morphism object is given by the internal $\uHom$,
$\Enr (\mP )(X,Y)=\uHom _{\mP}(X,Y)$. The identity and composition operations are 
the obvious ones.

This general discussion now applies to the model category $\mP = \precat ^n(\mM )$. 
Define
$$
nCAT (\mM ):= \Enr (\precat ^n(\mM )).
$$
Note that $nCAT (\mM )$ is a $\precat ^n(\mM )$-enriched
category, so (with the previously mentioned confusion of notation)
$$
nCAT (\mM ) \in \precat (\precat ^n(\mM ))=\precat ^{n+1}(\mM ).
$$
As $nCAT( \mM )$ is a strict $\precat ^n(\mM )$-enriched category, 
its Segal maps are isomorphisms. 
Notice that for any fibrant cofibrant objects $\pA , \pB \in \precat ^n(\mM )$,
$$
nCAT (\mM )(\pA ,\pB )= \uHom _{\precat ^n(\mM )}(\pA , \pB )
$$
is also fibrant. By the Segal isomorphisms and the fact that
fibrant and cofibrant objects are preserved by direct product, 
it follows that for any sequence
$\pA _0,\ldots  ,\pA _m$ the object
$$
nCAT (\mM )(\pA _0,\ldots  ,\pA _m ) \in \precat ^n(\mM )
$$
is fibrant. In particular, $nCAT (\mM )$ is a projectively fibrant
$\precat ^n(\mM )$-enriched precategory. 

One unfortunate consequence of the strictness on the first level is that
$nCAT (\mM )$ is not Reedy-fibrant. Therefore it isn't quite correct to write 
``$nCAT (\mM )\in \Ob ((n+1)CAT(\mM ))$''
since $nCAT(\mM )$ is not a fibrant object of $\precat ^{n+1}(\mM )$. Let $nCAT(\mM ) \rightarrow 
nCAT'(\mM )$ denote its
fibrant (and automatically cofibrant) replacement.  Then
$$
nCAT'(\mM ) \in \Ob ((n+1)CAT(\mM )).
$$
The difference between $nCAT(\mM )$ and $nCAT'(\mM )$ was one of the main obstacles which needed
to be overcome in the treatment of limits \cite{limits}.

If $\pA , \pB$ are cofibrant and fibrant $\mM$-enriched $n$-precategories,
then a morphism $f:\pA \rightarrt \pB$ corresponds to an object of 
then $\mM$-enriched $n$-precategory $\uHom (\pA ,\pB )$, or equivalently 
to a $1$-morphism in $nCAT( \mM )$. 
The morphism $f$ is a weak equivalence in $\precat ^n(\mM )$ 
if and only if it is an internal equivalence in $nCAT( \mM )$ i.e. it projects to an
isomorphism in $\tau _{\leq 1}(nCAT (\mM ))$, and we have an equivalence of categories
$$
\tau _{\leq 1}(nCAT (\mM ))\cong \Ho (\precat ^n(\mM )).
$$
This compatibility was formulated by Tamsamani in asking for $nCAT$ \cite{Tamsamani},
and may be proven using the same arguments as in the previous section. 

The above discussion applies with $\mM = \Sets$ to give the $n+1$-nerve $nCAT:= nCAT(\Sets )$ of $n$-nerves; and to $\mM = \mK$ to give the Segal $n+1$-category $nSeCAT:= nCAT (\mK )$
of Segal $n$-categories. These were used in \cite{descente} to discuss the notion
of higher stacks.


\begin{thebibliography}{MM}


\bibitem{JoyOfCats}
J. Ad\'amek, H. Herrlich, G. Strecker. 
Abstract and Concrete Categories---The Joy of Cats.
John Wiley and Sons (1990), electronic edition available at
\verb+http://katmat.math.uni-bremen.de/acc/+.


\bibitem{AdamekRosicky}
J. Ad\'amek, J. Rosicky. {\em Locally Presentable and Accessible Categories}.
London Math. Soc. Lecture Note Series {\bf 189}, Cambridge University Press (1994).


\bibitem{Adams}
J. Adams. {\em Infinite Loop Spaces}, Princeton University Press {\sc Annals
of Math. Studies} {\bf 90} (1978).

\bibitem{ArtinInventiones}
M. Artin.  Versal deformations and algebraic stacks, {\em Inventiones Math.}
{\bf 27} (1974), 165-189.

\bibitem{Badzioch}
B. Badzioch. Algebraic theories in homotopy theory. Ann. of Math. {\bf 155} (2002), 895–913.


\bibitem{BaezDolan}
J. Baez, J. Dolan. 
Higher-dimensional algebra and topological quantum
field theory. {\em Jour. Math. Phys.} {\bf 36} (1995), 6073-6105.

\bibitem{BaezDolanLetter}
J. Baez, J. Dolan. $n$-Categories, sketch of a
definition. Letter to R. Street, 29 Nov. and 3 Dec. 1995, 
{\tt http://math.ucr.edu/home/baez/ncat.def.html}.

\bibitem{BaezIntroduction}
J. Baez. An introduction to $n$-categories. 
{\em Category Theory and Computer Science (Santa Margherita Ligure 1997)}
{\sc Lect. Notes in Computer Science} {\bf 1290} (1997), 1-33.

\bibitem{BaezDolanIII}
J. Baez, J. Dolan. Higher-dimensional algebra III: $n$-categories and the
algebra of opetopes. Adv. in Math.
{\bf 135} (1998), 145-206.  


\bibitem{catBD}
J. Baez, J. Dolan. Categorification. Preprint math/9802029.


\bibitem{BaezMayTowards}
J. Baez, P. May. {\em Towards higher categories}.  
The IMA Volumes in Mathematics and its Applications Series {\bf 152},
Springer (2009).  

\bibitem{BakovicJurco}
I. Bakovi\'c, B. Jur\v{c}o. The classifying topos of a topological bicategory,
Arxiv preprint 0902.1750 (2009). 

\bibitem{BFSV}
C. Balteanu, Z. Fiedorowicz, R. Schwaenzl, and R. Vogt.
Iterated Monoidal Categories.
Preprint math.AT/9808082

\bibitem{BarwickThesis}
C. Barwick. $(\infty , n)-Cat$ as a closed model category. Doctoral
dissertation, University of Pennsylvania (2005). 

\bibitem{BarwickGottingen}
C. Barwick. $\infty$-groupoids, stacks, and Segal categories.
Seminars 2004-2005 of the Mathematical Institute, University of G\"ottingen (Y. Tschinkel, ed.). Universit\"atsverlag G\"ottingen (2005), 155-195. 



\bibitem{Barwick}
C. Barwick. On (enriched) left Bousfield localization of model categories.
Arxiv preprint arXiv:0708.2067 (2007). 

\bibitem{BarwickReedy}
C. Barwick. On Reedy  model categories. Preprint arXiv: 0708.2832 (2007). 

\bibitem{BarwickLeftRight}
C. Barwick. On left and right model categories and left and right Bousfield localizations.
To appear, {\em Homology, Homotopy and Applications}.

\bibitem{BarwickKan}
C. Barwick, D. Kan. Relative categories as another model for the homotopy theory of
homotopy theories I: the model structure. Preprint. 


\bibitem{Batanin}
M. Batanin. On the definition of weak $\omega$-category. Macquarie mathematics
report number 96/207, Macquarie University, NSW Australia.

\bibitem{Batanin2}
M. Batanin. Monoidal globular categories as a natural environment for the
theory of weak $n$-categories. {\em Adv. Math.} {\bf 136} (1998), 39-103/

\bibitem{BataninAinfty}
M. Batanin. Homotopy coherent category theory and $A_{\infty}$ structures in monoidal categories, 
{\em J. Pure and Appl. Alg.} {\bf 123} (1998), 67-103. 

\bibitem{BataninPenon}
M. Batanin. 
On the Penon method of weakening algebraic structures.
{\em J. Pure and Appl. Alg.} {\bf 172} (2002), 1-23.

\bibitem{BataninEH}
M. Batanin. The Eckmann–Hilton argument and higher operads. 
{\em Adv. in Math.} {\bf 217} (2008), 334-385.

\bibitem{BataninCisinskiWeber}
M. Batanin, D. Cisinski, M. Weber. 
Algebras of higher operads as enriched categories II. 
Preprint arXiv:0909.4715v1 (2009). 


\bibitem{Baues}
H. Baues. {\em Combinatorial homotopy and $4$-dimensional complexes.}
de Gruyter, Berlin (1991).

\bibitem{Beke}
T. Beke. Sheafifiable homotopy model categories. Math. Proc. Cambridge Phil. Soc.
{\bf 129} (2000), 447–-475.


\bibitem{Benabou}
J. B\'enabou. {\em Introduction to Bicategories}, Lect. Notes in Math. {\bf
47}, Springer-Verlag (1967).


\bibitem{BergerEckmannHilton}
C. Berger. 
Double loop spaces, braided monoidal categories and algebraic $3$-type of space.
{\em Contemp. Math.} {\bf 227} (1999), 49-65.

\bibitem{BergerCellularNerve}
C. Berger. 
A cellular nerve for higher categories.
Advances in Mathematics {\bf 169} (2002), 118-175.

\bibitem{BergerWreath}
C. Berger. 
Iterated wreath product of the simplex category and iterated loop spaces.
{\em Adv. in Math.} {\bf 213} (2007), 230-270.


\bibitem{BergerMoerdijk}
C. Berger, I. Moerdijk. On an extension of the notion of Reedy category.
Preprint arXiv: 0809.3341 (2008). 

\bibitem{BergnerModel}
J. Bergner. A model category structure on the category of simplicial categories.
Transactions of the AMS {\bf 359} (2007), 2043--2058.

\bibitem{BergnerThreeModels}
J. Bergner. Three models for the homotopy theory of homotopy theories.
{\em Topology} {\bf 46} (2007), 397-436.


\bibitem{BergnerInverting}
J. Bergner. Adding inverses to diagrams encoding algebraic structures.
{\em Homology, Homotopy Appl.} {\bf 10} (2008), 149–174.


\bibitem{BergnerSegal}
J. Bergner. A characterization of fibrant Segal categories.
{\em Proc. Amer. Math. Soc.} {\bf 135} (2007), 4031-4037.


\bibitem{BergnerRigidification}
J. Bergner. Rigidification of algebras over multi-sorted theories. Algebr. Geom. Topol. 
{\bf 6} (2006), 1925-1955.

\bibitem{BergnerSurvey}
J. Bergner. A survey of $(\infty , 1)$-categories. In 
{\em Towards higher categories} (J. Baez, P. May eds),  
The IMA Volumes in Mathematics and its Applications Series {\bf 152} (2009).  



\bibitem{BergnerSMSC}
J. Bergner. Simplicial monoids and Segal categories. Contemp. Math. {\bf 431} (2007),
59-83.

\bibitem{BergnerHoFiProd}
J. Bergner. On homotopy fiber products of homotopy theories. 
Arxiv preprint arXiv:0811.3175 (2008). 

\bibitem{Blander}
B. Blander. Local projective model structures on simplicial presheaves. 
{\em K-theory} {\bf 24}   (2001), 283-301.

\bibitem{BoardmanVogt}
J. Boardman, R. Vogt. {\em Homotopy invariant algebraic structures on
topological spaces.} Springer {\em L.N.M.} {\bf 347} (1973).

\bibitem{BondalKapranov}
A. Bondal, M. Kapranov. 
Enhanced triangulated categories. {\em Math. U.S.S.R. Sb.} {\bf 70} (1991), 93–107. 

\bibitem{BottTu}
R. Bott, L. Tu. {\em Differential Forms in Algebraic Topology}.  {\sc
Graduate Texts in Mathematics} {\bf 82} Springer, New York (1982).

\bibitem{BournDitopos}
D. Bourn. Sur les ditopos. {\em C. R. Acad. Sci. Paris} {\bf 279} (1974), 911–913.

\bibitem{Bourn}
D. Bourn. La tour de fibrations exactes des $n$-cat\'egories.
{\em Cahiers Top. G\'eom. Diff. Cat.} (1984). 

\bibitem{BournCordier}
D. Bourn and J.-M. Cordier. A general formulation of homotopy limits. {\em Jour.
Pure Appl. Alg.} {\bf 29} (1983), 129-141.

\bibitem{Bousfield-Kan}
A. Bousfield, D. Kan. {\em Homotopy limits, completions and localizations. }
Springer {\em Lecture Notes in Math.} {\bf 304} (1972).


\bibitem{BreenAsterisque}
L. Breen. On the classification of $2$-gerbs and $2$-stacks. {\em Ast\'erisque}
{\bf 225}, Soc. Math. de France (1994).

\bibitem{BreenRecent}
L. Breen. Monoidal categories and multiextensions. 
{\em Compositio Math.} {\bf 117} (1999), 295-335.



\bibitem{Breen}
L. Breen. On the classification of $2$-gerbs and $2$-stacks. {\em Ast\'erisque}
{\bf 225}, Soc. Math. de France (1994).

\bibitem{EBrown}
E. Brown, Jr. Finite computability of Postnikov complexes. {\em Ann. of Math.}
{\bf 65} (1957), 1-20.


\bibitem{KBrown}
K. Brown. Abstract homotopy theory and generalized sheaf cohomology. {\em
Trans. A.M.S.} {\bf 186} (1973), 419-458.

\bibitem{BrownGersten}
K. Brown, S. Gersten. {\em Algebraic $K$-theory as generalized sheaf cohomology}
Springer {\em Lecture Notes in Math.} {\bf 341} (1973), 266-292.

\bibitem{RBrown}
R. Brown. Groupoids and crossed objects in algebraic topology. 
{\em Homology Homotopy Appl.} {\bf 1} (1999), 1-78.

\bibitem{RBrown1}
R. Brown. Computing homotopy types using crossed $n$-cubes of
groups. {\em Adams Memorial Symposium on Algebraic Topology},
Vol 1, eds. N. Ray, G  Walker. Cambridge University Press, Cambridge
(1992) 187-210.

\bibitem{BrownGilbert}
R. Brown, N.D. Gilbert. Algebraic models of 3-types and automorphism
structures for crossed modules. {\em
Proc. London Math. Soc.} (3) {\bf 59} (1989), 51-73.

\bibitem{BrownHiggins}
R. Brown, P. Higgins. The equivalence of $\infty$-groupoids and crossed
complexes. {\em Cah. Top. Geom. Diff. Cat.} {\bf 22} (1981), 371-386.

\bibitem{BrownHiggins2}
R. Brown, P. Higgins. The classifying space of a crossed complex. {\em
Math. Proc. Camb. Phil. Soc.} {\bf 110}
(1991), 95-120.


\bibitem{BrownLoday}
R. Brown, J.-L. Loday. Van Kampen theorems for diagrams of spaces. {\em
Topology} {\bf 26} (1987), 311-335.

\bibitem{BrownLoday2}
R. Brown, J.-L. Loday. Homotopical excision, and Hurewicz theorems, for
$n$-cubes of spaces. {\em Proc. London Math. Soc.} {\bf 54} (1987), 176-192.

\bibitem{CabelloGarzon}
J. Cabello, A. Garzon. 
Closed model structures for algebraic models of $n$-types. 
Journal of Pure and Applied Algebra {\bf 103} (1995), 287-302.

\bibitem{ChengOpetopes}
The category of opetopes and the category of opetopic sets.  {\em Theory and Appl. of Categories} {\bf 11} (2003), 353-374.

\bibitem{ChengInfty}
E. Cheng. An omega-category with all duals is an omega groupoid.  {\em Appl. Cat. Struct.}
{\bf 15} (2007), 439-453.

\bibitem{ChengComparison}
E. Cheng. Comparing operadic theories of $n$-category.
Preprint arXiv:0809.2070 (2008).

\bibitem{ChengLauda}
E. Cheng, A. Lauda. {\em Higher-Dimensional Categories: an illustrated guide-book}.
\verb+http://cheng.staff.shef.ac.uk/guidebook/guidebook-new.pdf+ (2004). 

\bibitem{ChengMakkai}
E. Cheng, M. Makkai.  A note on the Penon definition of $n$-category.
Preprint arXiv:0907.3961 (2009).
 

\bibitem{CisinskiAsterisque}
D. Cisinski. {\em Les pr\'efaisceaux comme mod\`eles
des types d'homotopie.} {\sc Ast\'erisque}
{\bf 308}, S.M.F. (2006). 

\bibitem{CisinskiBatanin}
D. Cisinski. Batanin higher groupoids and homotopy types.
{\em Categories in
Algebra, Geometry and Mathematical Physics, proceedings of Streetfest}
(Batanin, Davydov, Johnson, Lack, Neeman, eds), Contemporary Math. {\bf 431} 
(2007), 171–186.

\bibitem{CisinskiDerivateurs}
D. Cisinski. 
Propri\'et\'es universelles et extensions de Kan d\'eriv\'ees.
{\em Theory and Appl. of Categories} {\bf 20} (2008), 605–649.

\bibitem{CohenLadaMay}
F. Cohen, T. Lada, J. P. May. {\em The homology of iterated loop spaces.}
Springer {\em L.N.M.} {\bf 533} (1976).


\bibitem{Cordier}
J. Cordier. Comparaison de deux cat\'egories d'homotopie de morphismes
coh\'erents. {\em Cahiers Top. G\'eom. Diff. Cat.} {\bf 30} (1989), 257-275.

\bibitem{CordierPorter86}
J. Cordier, T. Porter. Vogt's theorem on categories of homotopy coherent diagrams.
{\em Math. Proc. Camb. Phil. Soc.} {\bf 100} (1986), 65-90.

\bibitem{CordierPorter}
Cordier, Porter. Homotopy coherent category theory. 
{\em Trans. Amer. Math. Soc.} {\bf 349} (1997), 1-54.

\bibitem{Crans}
S. Crans. Quillen closed model structures for sheaves.
{\em J.  Pure Appl. Alg.} {\bf 101} (1995), 35-57. 


\bibitem{CransTensor}
S. Crans. A tensor product for Gray-categories.
{\em Theory  and Appl. of Categories} {\bf 5} (1999), 12–69.


\bibitem{CransBraidings}
S. Crans. On braidings, syllapses and symmetries. 
Cahiers de topologie et g\'eom\'etrie diff\'erentielle,
Volume 41 (2000), 2-74.


\bibitem{Curtis2}
E. Curtis. Lower central series of semisimplicial complexes. {\em Topology}
{\bf 2} (1963), 159-171.

\bibitem{Curtis}
E. Curtis. Some relations between homotopy and homology. {\em Ann. of Math.}
{\bf 82} (1965), 386-413.

\bibitem{HodgeIII}
P. Deligne. Th\'eorie de Hodge, III. {\em Publ. Math. I.H.E.S.} {\bf 44} (1974),
5-77.

\bibitem{DeligneMumford}
P. Deligne, D. Mumford. On the irreducibility of the space of curves of
a given genus. {\em Publ. Math. I.H.E.S.}  {\bf 36} (1969), 75-109.

\bibitem{lnm900}
P. Deligne, A. Ogus, J. Milne, K. Shih. {\em Hodge cycles, motives, and Shimura varieties}. {\sc Lecture Notes in Math.} {\bf 900},
Springer (1982).  

\bibitem{DrinfeldIntervalDG}
V. Drinfeld. DG quotients of DG categories. 
{\em J. of Algebra} {\bf 272} (2004), 643-691.

\bibitem{DuggerCombinatorial}
D. Dugger. Combinatorial model categories have presentations. 
{\em Adv.  in Math.} {\bf 164} (2001), 177-201. 

\bibitem{Dunn}
G. Dunn. Uniqueness of $n$-fold delooping machines. {\em J. Pure and Appl. Alg.}
{\bf 113} (1996), 159-193.

\bibitem{Duskin}
J. Duskin. Simplicial matrices and the nerves of weak $n$-categories I: nerves
of bicategories. 
{\em Theory and Applications of Categories} {\bf 9} (2002), 198–308.

\bibitem{DwyerHirschhornKan}
W. Dwyer, P. Hirschhorn, D. Kan. Model categories and  more general abstract
homotopy theory, a work in what we like to think of as progress. 

\bibitem{DwyerHirschhornKanSmith}
W. Dwyer, P. Hirschhorn, D. Kan, J. Smith. Homotopy
limit functors on model categories and homotopical categories. Mathematical Surveys and
Monographs {\bf 113} AMS, Providence,(2004).

\bibitem{DK1}
W. Dwyer, D. Kan. Simplicial localizations of categories. {\em J. Pure and Appl.
Algebra} {\bf 17} (1980), 267-284.

\bibitem{DK2}
W. Dwyer, D. Kan. Calculating simplicial localizations. {\em J. Pure and
Appl. Algebra} {\bf 18} (1980), 17-35.

\bibitem{DK3}
W. Dwyer, D. Kan. Function complexes in homotopical algebra. {\em Topology}
{\bf 19} (1980), 427-440.


\bibitem{DKS}
W. Dwyer, D. Kan, J. Smith. Homotopy commutative diagrams and their
realizations.
{\em J. Pure Appl. Algebra} {\bf 57} (1989), 5-24.

\bibitem{DwyerSpalinski}
W. Dwyer, J. Spalinski. Homotopy theories and model categories.
Handbook of Algebraic Topology, I. M. James, ed., Elsevier (1995). 

\bibitem{DyerLashoff}
E. Dyer, R. Lashoff. Homology of iterated loop spaces. {\em Amer. J. of Math.}
{\bf 84} (1962), 35-88.

\bibitem{Ellis}
G. Ellis. Spaces with finitely many nontrivial homotopy groups all of which are
finite. {\em Topology} {\bf 36} (1997), 501-504.

\bibitem{Fiedorowicz}
Z. Fiedorowicz. Classifying spaces of topological monoids and categories.
{\em Amer. J. Math.} {\bf 106} (1984), 301-350.

\bibitem{FiedorowiczVogt}
Z. Fiedorowicz, R. Vogt. 
Simplicial $n$-fold Monoidal Categories Model All $n$-fold Loop Spaces.
{\em Cah. Top. Geom. Diff. Cat.}
{\bf 44} (2003), 105-148.

\bibitem{Fukaya}
K. Fukaya. Morse homotopy, $A_{\infty}$-category and Floer homologies.
{\em Proceedings of GARC
Workshop on Geometry and Topology} (H. J. Kim, ed.), Seoul National University,
(1993).

\bibitem{Futia}
C. Futia. Weak omega categories I. Preprint arXiv:math/0404216 (2004).

\bibitem{GabrielZisman}
P. Gabriel, M. Zisman. {\em Calculus of fractions and homotopy theory.}
Springer, New York (1967).

\bibitem{Gaucher}
P. Gaucher. 
Homotopy invariants of higher dimensional categories and concurrency in computer science.
{\em Mathematical Structures in Computer Science} {\bf 10} (2000), 481-524. 

\bibitem{Giraud}
J. Giraud. {\em Cohomologie nonab\'elienne}, Grundelehren der Wissenschaften
in Einzeldarstellung {\bf 179} Springer-Verlag (1971).

\bibitem{GoerssJardine}
P. Goerss, R. Jardine. {\em Simplicial homotopy theory}.
{\sc Progress in Math.} {\bf 174}, Birkh\"{a}user (1999).  


\bibitem{GordonPowerStreet}
R. Gordon, A.J. Power, R. Street.  Coherence for tricategories {\em  Memoirs
A.M.S.} {\bf 117} (1995), 558 ff.

\bibitem{Grandis}
M. Grandis. 
Directed homotopy theory, I. The fundamental category. 
{\em Cah. Top. G\'eom. Diff. Cat.} {\bf 44} (2003), 281-316.


\bibitem{GrothendieckTohoku}
A. Grothendieck. 
Sur quelques points d'alg\`ebre homologique, I.
{\em Tohoku Math. J.}
{\bf 9} (1957), 119-221.

\bibitem{SGA1}
A. Grothendieck. {\em Revetements Etales et Groupe Fondamental (SGA I)}, 
{\sc Lecture Notes in Math.} {\bf 224}, Springer-Verlag (1971).

\bibitem{Grothendieck}
A. Grothendieck.  {\em Pursuing Stacks}, edited by G. Maltsiniotis,
to appear {\sc Documents Math\'ematiques}. See also 
\verb+http://www.grothendieckcircle.org/+.

\bibitem{GrothendieckDerivateurs}
A. Grothendieck. {\em Les D\'erivateurs}, edited by G. Maltsiniotis, available at
\verb+http://people.math.jussieu.fr/~maltsin/textes.html+.



\bibitem{Hain}
R. Hain. Completions of mapping class groups and the cycle $C-C^{-}$.
{\em Mapping class groups and moduli spaces of Riemann surfaces: proceedings of workshops held in G\"ottingen and Seattle}, {\sc Contemporary Math.} {\bf 150}, A.M.S. (1993),
75-106. 


\bibitem{HainMalcev}
R. Hain. The de rham homotopy theory of complex algebraic varieties I. 
{\em K-theory} {\bf 1} (1987),  271-324.

\bibitem{HainRelativeMalcev}
R. Hain. The Hodge de Rham theory of relative Malcev completion. 
{\em Ann. Sci. de l'E.N.S.} {\bf 31} (1998), 47-92.

\bibitem{Heller}
A. Heller. Homotopy theories. {\em Mem. Amer. Math. Soc.} {\bf  71} n. 388 (1988).

\bibitem{HermidaMakkaiPower}
C. Hermida, M. Makkai, A. Power. On weak higher-dimensional categories I.
{\em J. Pure Appl. Alg.} Part 1: {\bf 154} (2000), 221-246; Part 2: {\bf 157} (2001),
247-277.

\bibitem{HinichHAHA}
V. Hinich. Homological algebra of homotopy algebras. {\em Comm. in Algebra}
{\bf 25} (1997), 3291-3323. 

\bibitem{Hirschhorn}
P. Hirschhorn. {\em Model Categories and their Localizations}. 
{\sc Math. Surveys and Monographs} {\bf 99}, A.M.S. (2003). 

\bibitem{descente}
A. Hirschowitz, C. Simpson. Descente pour les $n$-champs. Preprint
math/9807049.

\bibitem{Hollander}
S. Hollander. 
A homotopy theory for stacks.  
{\em Israel Journal of Math.} {\bf 163} (2008) 93-124.


\bibitem{HoveyArxiv99}
M. Hovey. Monoidal model categories.  
Arxiv preprint math/9803002 (1998). 

\bibitem{Hovey}
M. Hovey. 
{\em Model categories}. 
{\sc Math. Surveys and Monographs} {\bf 63}, A.M.S. (1999). 

\bibitem{James}
I. James. Reduced product spaces. {\em Ann. of Math.} {\bf 62} (1955), 170-197.

\bibitem{Janelidze}
G. Janelidze. 
{\em Precategories and Galois theory}, Springer (1990). 


\bibitem{Jardine}
J.F. Jardine.  Simplicial presheaves, {\em J. Pure and Appl. Algebra} {\bf 47}
(1987), 35-87.

\bibitem{Johnson}
M. Johnson. The combinatorics of $n$-categorical pasting. {\em J. Pure and
Appl. Algebra} {\bf 62} (1989), 211-225.

\bibitem{JoyalLetter}
A. Joyal. Letter to A. Grothendieck (refered to in Jardine's paper).

\bibitem{JoyalQC}
A. Joyal. Quasi-categories and Kan complexes. 
{\em J. Pure Appl. Alg.}
{\bf 175} (2002), 207-222. 

\bibitem{JoyalTheta}
A. Joyal. Disks, duality and $\theta$-categories.  Preprint (1997).

\bibitem{JoyalKock}
A. Joyal, J. Kock. Weak units and homotopy $3$-types. 
{\em Categories in algebra, geometry and mathematical physics: conference and workshop in honor of Ross Street's 60th birthday}, {\sc Contemporary Math.} {\bf 431},
A.M.S. (2007), 257-276. 

\bibitem{JoyalTierney}
A. Joyal, M. Tierney. Algebraic homotopy types. Occurs as an entry in
the bibliography of \cite{BaezDolan}.

\bibitem{JoyalTierneyQCSC}
A. Joyal, M. Tierney. Quasi-categories vs Segal spaces. 
{\em Categories in algebra, geometry and mathematical physics: conference and workshop in honor of Ross Street's 60th birthday}, {\sc Contemporary Math.} {\bf 431},
A.M.S. (2007), 277-326. 


\bibitem{Kan1}
D. Kan. On c.s.s. complexes. {\em Amer. J. of Math.} {\bf 79} (1957), 449-476.

\bibitem{Kan2}
D. Kan. A combinatorial definition of homotopy groups. {\bf Ann. of Math.}
{\bf 67} (1958), 282-312.

\bibitem{Kan3}
D. Kan. On homotopy theory and c.s.s. groups. {\em Ann. of Math.} {\bf 68}
(1958), 38-53.


\bibitem{Kan4}
D. Kan. On c.s.s. categories. {\em Bol. Soc. Math. Mexicana} (1957), 82-94.

\bibitem{KapranovInv}
M. Kapranov. On the derived categories of coherent sheaves on some homogeneous spaces. 
{\em Inventiones} {\bf 92} (1988), 479-508.

\bibitem{KapranovVoevodsky}
M. Kapranov, V. Voevodsky. $\infty$-groupoid and homotopy types. {\em Cah. Top.
Geom. Diff. Cat.} {\bf 32} (1991), 29-46.

\bibitem{KaPaTo}
L. Katzarkov, T. Pantev, B. Toen. 
Algebraic and topological aspects of the schematization functor. {\em Compositio} {\bf 145} (2009),  633-686.

\bibitem{Keller}
B. Keller. Deriving DG categories. {\em Ann. Sci. E.N.S.} {\bf 27} (1994), 63-102.

\bibitem{Kelly}
G. Kelly, {\em Basic concepts of enriched category theory} London Math. Soc.
Lecture Notes {\bf 64}, Cambridge U. Press, Cambridge (1982).


\bibitem{Kock}
J. Kock. Weak identity arrows in higher categories. {\em Int. Math. Res. Papers}
(2006). 

\bibitem{Kock2}
J. Kock. 
Elementary remarks on units in monoidal categories. 
{\em Math. Proc. Cambridge Phil. Soc.} {\bf 144} (2008), 53-76.

\bibitem{KockJoyalBataninMascari}
J. Kock, A. Joyal, M. Batanin, J. Mascari. 
Polynomial functors and opetopes.
Preprint arXiv:0706.1033 (2007). 

\bibitem{Kondratieff}
G. Kondratiev. Concrete duality for strict infinity categories.
Preprint arXiv:0807.4256 (see also arXiv:math/0608436).


\bibitem{Kontsevich}
M. Kontsevich. Homological algebra of mirror symmetry.
{\em Proceedings of I. C. M.-94, Zurich} Birkh\"{a}user (1995), 120-139.


\bibitem{LaumonMB}
G. Laumon, L. Moret-Bailly.  {\em Champs Alg\'ebriques}. Springer (2000). 

\bibitem{Lawvere}
F. W. Lawvere. Functorial semantics of algebraic theories, Proc. Nat. Acad.
Sc. {\bf 50} (1963), 869-872.

\bibitem{LawvereThesis}
F. W. Lawvere. Functorial semantics of algebraic theories, Dissertation,
Columbia University 1963; Reprints in Theory Appl. Categ. {\bf 5} (2004), 23-107.

\bibitem{Leinster}
T. Leinster. A survey of definitions of $n$-category. 
{\em Theory and  Appl. of Categories} {\bf 10} (2002), 1–70.

\bibitem{LeinsterBook}
T. Leinster. {\em Higher operads, higher categories}. 
{\sc London Math. Soc. Lecture Notes} {\bf 298}, Cambridge University Press (2004).

\bibitem{Leroy}
O. Leroy.  Sur une notion de $3$-cat\'egorie adapt\'ee \`a l'homotopie.
Preprint Univ. de Montpellier 2 (1994).

\bibitem{Lewis}
L. G. Lewis. Is there a convenient category of spectra? {\em Jour. Pure and
Appl. Algebra} {\bf 73} (1991), 233-246.

\bibitem{Loday}
J.-L. Loday. Spaces with finitely many non-trivial homotopy groups. {\em J.
Pure Appl. Alg.} {\bf 24} (1982), 179-202.

\bibitem{LurieTopos}
J. Lurie. {\em Higher Topos Theory}. {\em Ann. of Math. Studies}  {\bf 170} (2009). 

\bibitem{LurieAlgebra}
J. Lurie. Derived Algebraic Geometry II--VI. Arxiv preprints (2007-2009). 

\bibitem{LurieGC}
J. Lurie. (Infinity,2)-Categories and the Goodwillie Calculus I. 
Preprint arXiv:0905.0462v2 (2009). 


\bibitem{MacDuff}
D. MacDuff. On the classifying spaces of discrete monoids. {\em Topology} {\bf
18} (1979), 313-320.

\bibitem{Mackaay}
M. Mackaay. Spherical $2$-categories and $4$-manifold invariants. {\em Adv. in Math.}
{\bf 143} (1999), 288-348.


\bibitem{MacLane}
S. MacLane. Categories for the working mathematician. Springer (1971).

\bibitem{MakkaiPare}
M. Makkai, R. Par\'e. 
{\em Accessible categories: the foundations of categorical model theory}.
{\sc Contemporary Math.} {\bf 104}, A.M.S. (1989). 


\bibitem{MaltsiniotisAsterisque}
G. Maltsiniotis. {\em La th\'eorie de l'homotopie de Grothendieck}. 
{\sc Ast\'erisque} {\bf 301}
(2005). 

\bibitem{MaltsiniotisGroGpd}
G. Maltsiniotis. Infini groupo\"'{i}des non stricts, d'apr\`es Grothendieck.  
Preprint (2007).

\bibitem{MaltsiniotisGroBat}
G. Maltsiniotis. Infini cat\'egories non strictes, une nouvelle d\'efinition. Preprint (2007).



\bibitem{Massey}
W. Massey. {\em Algebraic Topology: An Introduction}. {\sc Graduate Texts in Mathematics} {\bf 56},
Springer (1977). 

\bibitem{MaySimplicial}
J. P. May. {\em Simplicial objects in algebraic topology.} Van Nostrand (1967).

\bibitem{MayLoops}
J.P. May. {\em The geometry of iterated loop spaces} Springer {\em L.N.M.}
{\bf 271} (1972).

\bibitem{MayFibs}
J. P. May. Classifying spaces and fibrations. {\em Mem. Amer. Math. Soc.}
{\bf 155} (1975).

\bibitem{MayThomason}
J. P. May, R. Thomason. The uniqueness of infinite loop space machines. {\em
Topology} {\bf 17} (1978),  205-224.

\bibitem{MorelVoevodsky}
F. Morel, V. Voevodsky. $\aaa ^1$-homotopy theory of schemes. 
{\em Publ. Math. I.H.E.S.} {\bf 90} (1999), 45-143.

\bibitem{Moriya}
S. Moriya. 
Rational homotopy theory and differential graded category.
{\em J. of Pure and Appl. Alg.} {\bf 214} (2010), 422-439.  


\bibitem{PaoliAdvances}
S. Paoli. Weakly globular $cat^n$-groups and Tamsamani's model.
{\em  Advances in Math.} {\bf 222} (2009), 621-727.


\bibitem{Pelissier}
R. Pelissier. Cat\'egories enrichies faibles. 
Thesis, Universit\'e de Nice (2002),
\verb+http://tel.archives-ouvertes.fr/tel-00003273/fr/+ .

\bibitem{Penon}
J. Penon. Approche polygraphique des $\infty$-cat\'egories non strictes. 
{\em Cahiers Top. G\'eom. Diff. Cat.} {\bf XL-1} (1999), 31-80.

\bibitem{Power}
A. Power. Why tricategories? {\em Information and Computation} {\bf 120} (1995), 251-262.

\bibitem{Pridham}
J. Pridham. 
Pro-algebraic homotopy types.
{\em Proc. of the London Math. Soc.} {\bf 97} (2008), 273-338.

\bibitem{Quillen}
D. Quillen. {\em Homotopical algebra} Springer {\em L.N.M.} {\bf 43} (1967).

\bibitem{QuillenAnnals}
D. Quillen. Rational Homotopy Theory. {\em Ann. Math.} {\bf 90} (1969), 205-295.


\bibitem{Reedy}
C. Reedy. Homotopy theory of model categories. Preprint (1973) available from P.
Hirschhorn.

\bibitem{Rezk}
C. Rezk. A model for the homotopy theory of homotopy theory.
{\em Trans. Amer. Math. Soc.} {\bf 353} (2001), 973–1007. 

\bibitem{RezkCartesian}
C. Rezk. A cartesian presentation of weak n-categories. 
Arxiv preprint arXiv:0901.3602 (2009).

\bibitem{Riehl}
E. Riehl. On the structure of simplicial categories associated to quasi-categories.
Preprint arXiv:0912.4809 (2009). 

\bibitem{RosickyTholen}
J. Rosick\'y and W. Tholen, Left-determined model categories and universal
homotopy theories, Trans. Amer. Math. Soc. 355 (2003), 3611-3623.

\bibitem{RosickyHomotopyVarieties}
J. Rosick\'y. 
On homotopy varieties. 
{\em Adv. in Math.} {\bf 214} (2007), 525-550. 

\bibitem{RosickyCombinatorial}
J. Rosick\'y. On combinatorial model categories. 
{\em Appl. Cat. Structures} {\bf 17} (2009) 303-316. 


\bibitem{Vogt}
R. Schw\"{a}nzl, R. Vogt. Homotopy homomorphisms and the hammock localization. 
{\em Papers in honor of Jos\'e Adem}, {\em Bol. Soc. Mat. Mexicana} {\bf 37}
(1992), 431-448. 


\bibitem{SegalHspaces}
G. Segal. Homotopy everything $H$-spaces. Preprint.

\bibitem{SegalInventiones}
G. Segal. Configuration spaces and iterated loop spaces. {\em Inv. Math.} {\bf
21} (1973), 213-221.

\bibitem{Segal}
G. Segal. Categories and cohomology theories. {\em Topology} {\bf 13}
(1974), 293-312.

\bibitem{ShipleySchwede02}
B. Shipley, S. Schwede. 
Equivalences of monoidal model categories. 
{\em Algebr. Geom. Topol.} {\bf 3} (2003), 287-334. 

\bibitem{kobe}
C. Simpson. Homotopy over the complex numbers and generalized de Rham
cohomology. {\em Moduli of Vector Bundles}, M. Maruyama (Ed.) {\em Lecture
Notes in Pure and Applied Math.} {\bf 179}, Marcel Dekker (1996), 229-263.

\bibitem{flexible}
C. Simpson. Flexible sheaves. Preprint q-alg/9608025.

\bibitem{realization}
C. Simpson. The topological realization of a simplicial presheaf. Preprint
q-alg/9609004.

\bibitem{geometricN}
C. Simpson. Algebraic (geometric) $n$-stacks. Preprint
alg-geom/9609014.

\bibitem{svk}
C. Simpson. A closed model structure for $n$-categories, internal $Hom$,
$n$-stacks and generalized Seifert-Van Kampen. Preprint alg-geom/9704006.

\bibitem{limits}
C. Simpson. Limits in $n$-categories. Preprint alg-geom 9708010.

\bibitem{effective}
C. Simpson. Effective generalized Seifert-Van Kampen: how to calculate
$\Omega X$. Preprint q-alg/9710011.

\bibitem{BBDSH}
C. Simpson. On the Breen-Baez-Dolan stabilization hypothesis. Preprint 
math.CT/9810058. 

\bibitem{hty3types}
C. Simpson. 
Homotopy types of strict 3-groupoids. Preprint, math.CT/9810059.

\bibitem{SmithCombi}
J. Smith. Combinatorial model categories. Manuscript refered to in 
\cite{DuggerCombinatorial}. 

\bibitem{Stanculescu}
A. Stanculescu. 
A homotopy theory for enrichment in simplicial modules. Preprint
arXiv:0712.1319.

\bibitem{Stasheff}
J. Stasheff. Homotopy associativity of $H$-spaces, I, II. {\em Trans. Amer.
Math.
Soc.} {\bf 108} (1963), 275-292, 293-312.


\bibitem{Street}
R. Street. The algebra of oriented simplexes. {\em Jour. Pure and Appl.
Algebra} {\bf 49} (1987), 283-335.

\bibitem{Street2}
R. Street. Weak $\omega$-categories. 
{\em Diagrammatic morphisms and applications (San Francisco, 2000)}
{\em Contemporary Mathematics} {\bf 318}, AMS (2003), 207-213.

\bibitem{Tabuada}
G. Tabuada. Differential graded versus Simplicial categories. 
Preprint
arXiv:0711.3845. 

\bibitem{TabuadaSpectral}
G. Tabuada. Homotopy theory of spectral categories.
{\em Adv. in Math.} {\bf 221} (2009), 1122-1143. 

\bibitem{TamsamaniThesis}
Z. Tamsamani.  Sur des notions de $n$-categorie et $n$-groupoide non-stricte
via des ensembles multi-simpliciaux. Thesis, Universit\'e Paul Sabatier,
Toulouse (1996) available on alg-geom (95-12 and 96-07).

\bibitem{Tamsamani}
Z. Tamsamani. Sur des notions de $n$-categorie et $n$-groupoide non-stricte
via des ensembles multi-simpliciaux. {\em $K$-theory} {\bf 16} (1999), 51-99.

\bibitem{Tanre}
D. Tanre. {\em  Homotopie Rationnelle: mod\`eles de Chen, Quillen, Sullivan.}
Springer {\em Lecture Notes in Mathematics} {\bf 1025}  (1983).

\bibitem{Thomason}
R. Thomason. Uniqueness of delooping machines. {\em Duke Math. J.} {\bf 46}
(1979), 217-252.

\bibitem{ThomasonENS}
R. Thomason. 
Algebraic $K$-theory and \'etale cohomology.
{\em Ann. Sci. E.N.S.} {\bf 18} (1985), 437–552.

\bibitem{Toen}
B. Toen. Champs affines.  {\em Selecta Math.} {\bf 12} (2006), 39-134.

\bibitem{Trimble4}
T. Trimble. Notes on tetracategories,
\verb+http://math.ucr.edu/home/baez/trimble/tetracategories.html+

\bibitem{Verity}
D. Verity. Weak complicial sets I. Basic homotopy theory.
{\em Adv. in Math.} {\bf 219} (2008), 1081-1149.

\bibitem{Voevodsky}
V. Voevodsky. The Milnor conjecture. Preprint (1996). 


\bibitem{Weber}
M. Weber. 
Yoneda Structures from 2-toposes.
{\em Appl. Cat. Struct.} {\bf 15} (2007), 259-323.


\bibitem{Whitehead}
G. Whitehead. {\em Elements of Homotopy Theory}, Springer, New-York (1978).

\bibitem{WhiteheadAspheric}
J. H. C. Whitehead. On the asphericity of regions in a $3$-sphere. {\em Fund.
Math.} {\bf 32} (1939), 149-166.


\bibitem{Zawadowski}
M. Zawadowski. Lax Monoidal Fibrations. Preprint 
arXiv:0912.4464 (2009). 





\end{thebibliography}
\end{document}